\documentclass[hidelinks]{memo-l}
\usepackage{imakeidx}
\usepackage{amsmath, amssymb, amsthm, amsfonts, amstext, amsbsy, booktabs, enumitem, mathrsfs, mathtools, stmaryrd,  upgreek, verbatim}
\indexsetup{level=\section*,toclevel=section,noclearpage}
 \makeindex[name=concepts,title=Concepts]
 \makeindex[name=symbols,title=Symbols]
\usepackage[normalem]{ulem}
\usepackage[shadow]{todonotes}
\usepackage[a4paper]{geometry}
\usepackage{thmtools}
\usepackage{microtype}
\numberwithin{equation}{section}
\usetikzlibrary{arrows,backgrounds,patterns,calc,decorations.pathmorphing,positioning,fit,shapes} 

\newcommand{\markdef}[1]{\textbf{#1}} 
\newcommand{\AND}{\mathbin{\, \wedge \,}}
\newcommand{\OR}{\mathbin{\, \vee \,}}
\newcommand{\IMPLIES}{\mathbin{\, \Rightarrow \,}}
\newcommand{\IFF}{\mathbin{\, \Leftrightarrow \,}}
\renewcommand{\implies}{\Rightarrow}	
\renewcommand{\iff}{\Leftrightarrow}		
\newcommand{\FORALL}[1]{\forall {#1} \, }
\newcommand{\EXISTS}[1]{\exists {#1} \, }
\newcommand{\EXISTSONE}[1]{\exists ! {#1} \, }
\DeclareMathOperator{\Fv}{\textsc Fv} 
\DeclareMathOperator{\Subf}{Subf}
\newcommand{\LST}{\mathsf{LST}}		
\newcommand{\embeds}{%
\mathrel{\vcenter{\offinterlineskip
\hbox{$\sqsubset$}\vskip.2ex\hbox{$\sim$}}}}
\newcommand{\biembeds}{\approx}
\DeclareMathOperator{\Mod}{Mod} 
\DeclareMathOperator{\Aut}{Aut} 
\newcommand{\ZF}{\mathsf{ZF}} 		
\newcommand{\ZFC}{\mathsf{ZFC}}		
\newcommand{\AC}{\mathsf{AC}}			
\newcommand{\AComega}{{\AC}_\omega}	
\newcommand{\DC}{\mathsf{DC}}			
\newcommand{\DCR}{\DC ( \R )}			
\newcommand{\CH}{\mathsf{CH}}		
\newcommand{\GCH}{\mathsf{GCH}}	
\newcommand{\AD}{\mathsf{AD}}		
\newcommand{\ADR}{\mathsf{AD}_{\R}}	
\newcommand{\PSP}{\mathsf{PSP}}		
\newcommand{\Unif}{\mathsf{Unif}}		
\newcommand{\MA}{\mathsf{MA}}		
\newcommand{\MM}{\mathsf{MM}}		
\newcommand{\PFA}{\mathsf{PFA}}		
\newcommand{\forcing}[1]{\textbf{\upshape #1}} 	
\newcommand{\I}{\mathbf{I}}
\newcommand{\II}{\mathbf{II}}
\newcommand{\GL}{G_{\mathrm{L}}}
\newcommand{\GW}{G_{\mathrm{W}}}

\newcommand{\leqW}{\leq_{\mathrm{W}}}
\newcommand{\leql}{\leq_{\mathrm{L}}}
\newcommand{\nleql}{\nleq_{\mathrm{L}}}
\newcommand{\bSigma}{\boldsymbol{\Sigma}}
\newcommand{\bPi}{\boldsymbol{\Pi}}
\newcommand{\bGamma}{\boldsymbol{\Gamma}}
\newcommand{\bLambda}{\boldsymbol{\Lambda}}
\newcommand{\bDelta}{\boldsymbol{\Delta}}
\newcommand{\dual}{\breve}	
\newcommand{\Gdelta}{\mathbf{G}_\delta}	
\newcommand{\Fsigma}{\mathbf{F}_\sigma}	
\newcommand{\blambda}{ \boldsymbol{\lambda}}
\newcommand{\body}[1]{ [ {#1} ]} 
\DeclareMathOperator{\WO}{WO} 	
\DeclareMathOperator{\LO}{LO} 		
\DeclareMathOperator{\Tr}{Tr} 
\newcommand{\Nbhd}{{\boldsymbol N}\!} 

\newcommand{\cointerval}[2]{[{#1};{#2})}
\newcommand{\seqof}[2]{\mathopen \langle #1 \Mid #2 \mathclose \rangle} 
\newcommand{\seqofLR}[2]{\left \langle #1 \Mid #2 \right \rangle} 
\newcommand{\seq}[1]{\mathopen \langle #1 \mathclose \rangle}	
\newcommand{\seqLR}[1]{\left \langle #1 \right \rangle}	
\DeclareMathOperator{\lh}{lh}				
\newcommand{\op}[2]{{\boldsymbol \langle} {#1} , {#2} {\boldsymbol \rangle}} 	
\newcommand{\Code}[1]{\mathopen{\langle\!\langle} {#1} \mathclose{\rangle\!\rangle} }
\newcommand{\eq}[1]{{\boldsymbol [}{#1} {\boldsymbol ]}}	
\newcommand{\bdelta}{\boldsymbol{\delta}}
\newcommand{\id}{\mathrm{id}}
\newcommand{\pre}[2]{\prescript{#1}{}{#2}}				
\newcommand{\leqlex}{\leq_{\mathrm{lex}}}		
\DeclareMathOperator{\HC}{HC}			
\newcommand{\On}{\mathord{\mathrm{Ord}}}		
\newcommand{\Cn}{\mathord{\mathrm{Card}}}		
\newcommand{\Vv}{\mathord{\mathrm{V}}}		
\newcommand{\Ll}{\mathord{\mathrm{L}}}			
\newcommand{\OD}{\operatorname{OD}}			
\newcommand{\symdif}{\mathop{\triangle}}	
\newcommand{\Mid}{\boldsymbol\mid}			
\newcommand{\setof}[2]{\mathopen \{{#1}\Mid{#2} \mathclose\}} 
\newcommand{\setofLR}[2]{\left \{{#1} \Mid {#2} \right\}} 
\newcommand{\set}[1]{\mathopen \{ {#1} \mathclose \}} 
\newcommand{\setLR}[1]{\left \{ {#1} \right \}} 
\newcommand{\conc}{{}^\smallfrown}		
\newcommand{\Q}{\mathbb{Q}}	
\newcommand{\R}{\mathbb{R}}	
\newcommand{\Z}{\mathbb{Z}}		
\DeclareMathOperator{\Fin}{Fin}	
\DeclareMathOperator{\Sym}{Sym} 	
\DeclareMathOperator{\dom}{dom} 		
\DeclareMathOperator{\ran}{ran} 		
\newcommand{\into}{\rightarrowtail} 	
\newcommand{\onto}{\twoheadrightarrow}		
\DeclareMathOperator{\cof}{cof} 			
\newcommand{\cardLR}[1]{\left | #1 \right |}		
\newcommand{\card}[1]{\mathopen | #1 \mathclose |}		
\DeclareMathOperator{\ot}{ot}				

\newcommand{\V}{\mathsf{v}}
\newcommand{\MMM}{\mathfrak{M}}
\newcommand{\BBB}{\mathfrak{B}}
\newcommand{\DDD}{\mathfrak{D}}
\newcommand{\UUU}{\mathfrak{U}}
\newcommand{\QQ}{\mathbb{Q}}
\newcommand{\LL}{\mathcal{L}}
\newcommand{\bB}{\mathbf {B}}

\newcommand{\TT}{\mathbb{T}}
\newcommand{\CT}{\mathrm{CT}}
\newcommand{\OCT}{\mathrm{OCT}}
\newcommand{\RCT}{\mathrm{RCT}}
\newcommand{\GR}{\mathrm{GRAPH}}

\newcommand{\Alg}{\mathrm{Alg}}
\newcommand{\edge}{\mathrel{\mathbf{E}}}
\newcommand{\order}{\mathrel{\trianglelefteqslant}}
\DeclareMathOperator{\SUCC}{Succ}

\DeclareMathOperator{\weight}{w}
\newcommand{\appl}[2]{{#1}``{#2}}		
\newcommand{\equals}{\boldsymbol{\bumpeq}}
\newcommand{\nequals}{\not\boldsymbol{\bumpeq}}
\newcommand{\pow}{\mathscr{P}}

\DeclareMathOperator{\PROJ}{p}
\newcommand{\gr}{\operatorname{graph}}
\newcommand{\Inv}{\operatorname{Inv}}
\newcommand{\bS}{\boldsymbol{S}}
\DeclareMathOperator{\pos}{pos}
\newenvironment{enumerate-(a)}{\begin{enumerate}[label={\upshape (\alph*)}, leftmargin=2pc]}{\end{enumerate}}
\newenvironment{enumerate-(A)}{\begin{enumerate}[label={\upshape (\Alph*)}, leftmargin=2pc]}{\end{enumerate}}
\newenvironment{enumerate-(i)}{\begin{enumerate}[label={\upshape (\roman*)}, leftmargin=2pc]}{\end{enumerate}}
\newenvironment{enumerate-(I)}{\begin{enumerate}[label={\upshape (\Roman*)}, leftmargin=2pc]}{\end{enumerate}}
\newenvironment{enumerate-(1)}{\begin{enumerate}[label={\upshape (\arabic*)}, leftmargin=2pc]}{\end{enumerate}}
\theoremstyle{plain}
\newtheorem{theorem}{Theorem}[section]
\newtheorem{fact}[theorem]{Fact}

\newtheorem{proposition}[theorem]{Proposition}
\newtheorem{lemma}[theorem]{Lemma}
\newtheorem{corollary}[theorem]{Corollary}
\newtheorem{claim}{Claim}[theorem]

\theoremstyle{definition}
\newtheorem{definition}[theorem]{Definition}

\theoremstyle{remark}
\newtheorem{remark}[theorem]{Remark}
\newtheorem{remarks}[theorem]{Remarks}
\newtheorem{example}[theorem]{Example}
\newtheorem{examples}[theorem]{Examples}
\newtheorem*{notation}{Notation}
\newtheorem{convention}[theorem]{Convention}

\newtheorem*{cor:omega1underPFA}{Corollary \protect\ref{cor:omega1underPFA}}
\newtheorem*{th:Gammacompleteness}{Theorem \protect\ref{th:Gammacompleteness}}
\newtheorem*{th:mainprojective}{Theorem \protect\ref{th:mainprojective}}
\newtheorem*{th:largecardinals1}{Theorem \protect\ref{th:largecardinals1}}
\newtheorem*{th:solovaymodel2}{Theorem \protect\ref{th:solovaymodel2}}
\newtheorem*{th:parovicenko1}{Theorem \protect\ref{th:parovicenko1}}
\newtheorem*{th:parovicenko2}{Theorem \protect\ref{th:parovicenko2}}
\newtheorem*{th:QembedsintoCTalephomega}{Theorem \protect\ref{th:QembedsintoCTalephomega}}
\newtheorem*{th:Qembedsundersharps}{Theorem \protect\ref{th:Qembedsundersharps}}
\newtheorem*{th:embeddingsofprojectiveqos}{Theorem \protect\ref{th:embeddingsofprojectiveqos}}
\newtheorem*{prop:parovicenkoAC}{Proposition \protect\ref{prop:parovicenkoAC}}
\newtheorem*{prop:parovicenkoACimproved}{Proposition \protect\ref{prop:parovicenkoACimproved}}

\renewcommand{\restriction}{\mathpunct{\upharpoonright}}
\usepackage{hyperref}

\makeindex

\begin{document}

\frontmatter

\title{Souslin quasi-orders and bi-embeddability \\ of uncountable structures}

\author{Alessandro Andretta}
\address{Universit\`a di Torino \\ Dipartimento di Matematica\\ Via Carlo Alberto, 10 \\ 10123 Torino \\ Italy}
\email{alessandro.andretta@unito.it}

\author{Luca Motto Ros}
\address{Universit\`a di Torino \\ Dipartimento di Matematica\\ Via Carlo Alberto, 10 \\ 10123 Torino \\ Italy}
\email{luca.mottoros@unito.it}

\date{\today}



\begin{abstract}
We provide analogues of the results from~\cite{Friedman:2011cr,Camerlo:2012kx} (which correspond to the case \( \kappa = \omega \)) for arbitrary \( \kappa \)-Souslin quasi-orders on any Polish space, for \( \kappa \) an infinite cardinal smaller than the cardinality of \( \mathbb{R} \). 
These generalizations yield  a variety of results concerning the complexity of the embeddability relation between graphs or lattices of size \( \kappa \), the isometric embeddability relation between complete metric spaces of density character \( \kappa \), and the linear isometric embeddability relation between (real or complex) Banach spaces of density \( \kappa \).
\end{abstract}

\thanks{The authors would like to thank the anonymous referee for carefully reading the manuscript, pointing out inaccuracies, and suggesting changes that substantially improved the paper. 
All remaining errors are responsibility of the authors alone.}

\maketitle
\tableofcontents
\mainmatter

\section{Introduction} \label{sec:introduction}
\subsection{What we knew} \label{subsec:knew}

\subsubsection{Equivalence relations and classification problems} \label{subsubsec:eqrel}

The analysis of definable equivalence relations on Polish spaces (or, more generally, on standard Borel spaces) has been one of the most active areas in descriptive set theory for the last two decades --- see~\cite{Friedman:1989bs,Kechris:1999fk,Kanovei:2008qo,Gao:2009fv} for excellent surveys of the subject. 
The main goal of this research area is to classify equivalence relations by means of reductions: if \( E \) and \( F \) are equivalence relations on Polish or standard Borel spaces \( X \) and \( Y \), then \( f \colon X \to Y \) reduces \( E \) to \( F \) if and only if
\begin{equation}\label{eq:reduction}
 x \mathrel{E} x' \iff f ( x ) \mathrel{F} f ( x' ),
\end{equation}
for every \( x , x' \in X \). 
In order to obtain nontrivial results one usually imposes definability assumptions on \( E \) and \( F \) and/or on the reduction \( f \).
For example one may assume that
\begin{itemize}[leftmargin=1pc]
\item
\( E \) and \( F \) are Borel or analytic%
\footnote{As customary in descriptive set theory, analytic sets are also called \( \bSigma^{1}_{1} \).} 
and \( f \) is continuous or Borel (as it is customary when studying actions of Polish groups on standard Borel spaces~\cite{Becker:1996uq,Hjorth:2000zr,Kechris:2004jl}), or that
\item
\( E \) and \( F \) are projective and \( f \) is Borel~\cite{Harrington:1979hc}, or that 
\item \( E \), \( F \), and \( f \) belong to some inner model of determinacy, such as \( \Ll ( \R ) \)~\cite{Hjorth:1995ve,Hjorth:2000zr}.
\end{itemize} 
When the reducing function \( f \) is Borel we say that \( E \) is \emph{Borel reducible} to \( F \) (in symbols \( E \leq_{\bB} F \)) and that \( f \) is a \emph{Borel reduction} of \( E \) to \( F \).
If \( E \leq_{\bB} F \leq_{\bB} E \), then \( E \) and \( F \) are said to be \emph{Borel bi-reducible}, in symbols \( E \sim_{\bB} F \).

Borel reducibility is usually interpreted both as a topological version of the ubiquitous notion of classification of mathematical objects, and as a tool for computing cardinalities of quotient spaces in an effective way. 
If \( E \) and \( F \) are equivalence relations, then it can be argued that the statement \( E \leq_{\bB} F \) is a precise mathematical formulation of the following informal assertions:
\begin{itemize}[leftmargin=1pc]
\item
the problem of classifying the elements of \( X \) up to \( E \)-equivalence is no more complex than the problem of classifying the elements of \( Y \) up to \( F \)-equivalence;
\item
the elements of \( X \) can be classified (in a definable way) up to \( E \)-equivalence using the \( F \)-equivalence classes as invariants;
\item
the quotient space \( X / E \) has cardinality less or equal than the cardinality of \( X / F \), and this inequality can be witnessed in a concrete and definable way.
\end{itemize}

The theory of Borel reducibility has been used to gauge the complexity of many natural problems. 
One striking example is given by the classification of countable structures up to isomorphism briefly described below, whose systematic study was initiated by H.~Friedman and Stanley in~\cite{Friedman:1989bs}. 
In a pioneering work in the mid 70s~\cite{Vaught:1974kl} Vaught observed that by identifying the universe of a countable structure with \( \omega \), the collection \( \Mod^\omega_{\LL} \) of countable \( \LL \)-structures can be construed as a Polish space, and the isomorphism relation \( \cong \) on this space is an analytic equivalence relation.
In order to exploit this identification most effectively, first-order logic must be replaced by its infinitary version \( \LL_{ \omega_1 \omega} \), where countable conjunctions and disjunctions are allowed~\cite{Keisler:1971vn}: by a theorem of Lopez-Escobar the (\( \LL_{\omega_1 \omega} \)-)elementary classes, that is the sets \( \Mod^\omega_\upsigma \) of all countable models of an \( \LL_{\omega_1 \omega} \)-sentence \( \upsigma \), are exactly the Borel subsets of \( \Mod^\omega_{\LL} \) which are invariant under isomorphism. 
It follows that each elementary class is a standard Borel subspace of \( \Mod^\omega_\LL \), and that the restriction of the isomorphism relation to \( \Mod^\omega_\upsigma \), denoted in this paper either by \( \cong \restriction \Mod^\omega_\upsigma \) or by \( \cong^\omega_\upsigma \), is an analytic equivalence relation, whose complexity can then be analyzed in terms of Borel reducibility.

Other classification problems whose complexity with respect to \( \leq_\bB \) has been widely studied in the literature include e.g.\ the classification of Polish metric spaces up to isometry~\cite{Clemens:2001pu,Gao:2003qw,Clemens:2012hg,Camerlo:2015ar} and the classification of separable Banach spaces up to linear isometry~\cite{Melleray:2007di} or up to isomorphism~\cite{Ferenczi:2009fk}. 

\subsubsection{Quasi-orders and embeddability}

The concept of reduction from~\eqref{eq:reduction} can be applied to arbitrary binary relations, such as partial orders and quasi-orders (also known as preorders)~\cite{Harrington:1988ij}. 
Quasi-orders are reflexive and transitive relations, and their symmetrization gives rise to an equivalence relation.
For example, the embeddability relation \( \embeds \) between structures (graphs, combinatorial trees, lattices, quasi-orders, partial orders, \dots) is a quasi-order, and its symmetrization is the relation \( \biembeds \) of bi-embeddability.
The restriction of \( \embeds \) to the Borel set \( \Mod^ \omega _ \upsigma \), denoted by \( \embeds \restriction \Mod^\omega_\upsigma \) or by \( \embeds^\omega_\upsigma \), is an analytic quasi-order.
Similarly, the relation of bi-embeddability on \( \Mod^ \omega _ \upsigma \), denoted by \( \biembeds \restriction \Mod^\omega_\upsigma \) or \( \biembeds^\omega_\upsigma \), is an analytic equivalence relation.

In~\cite[Theorem 3.1]{Louveau:2005cq} Louveau and Rosendal proved that it is not possible to classify (in a reasonable way) all countable structures up to bi-embeddability, as the embeddability relation \( \embeds \) is as complex as possible with respect to \( \leq_\bB \) (whence also \( \biembeds \) is as complex as possible).

\begin{theorem}[Louveau-Rosendal]\label{th:LouveauRosendal}
The embeddability relation \( \embeds_{\CT}^\omega \) on countable combinatorial trees (i.e.\ connected acyclic graphs) is a \emph{\( \leq_\bB \)-complete} analytic quasi-order, that is:
\begin{enumerate-(a)}
\item\label{th:LouveauRosendal-i}
\( \embeds_{\CT}^\omega \) is an analytic quasi-order on \( \CT_ \omega \), the Polish space of countable combinatorial trees, and
\item\label{th:LouveauRosendal-ii}
every analytic quasi-order is Borel reducible to \( \embeds_{\CT}^\omega \).
\end{enumerate-(a)}
Thus also the bi-embeddability relation \( \biembeds_\CT^\omega \) on \( \CT_\omega \) is a \emph{\( \leq_\bB \)-complete} analytic equivalence relation.
\end{theorem}

\begin{remark} \label{rmk:introLR}
In an abstract setting, it is not hard to find \( \leq_\bB \)-complete quasi-orders and equivalence relations.
In fact it is easy to show that for every pointclass \( \bGamma \) closed under Borel preimages, countable unions, and projections, the collection of quasi-orders (respectively: equivalence relations) in \( \bGamma \) admits a \( \leq_\bB \)-complete element, i.e.~there is a quasi-order (respectively, an equivalence relation) \( \mathcal{U}_{ \bGamma} \in \bGamma \) such that \( R \leq_{\bB} \mathcal{U}_{ \bGamma} \), for any quasi-order (respectively: equivalence relation) \( R \in \bGamma \) --- see~\cite[Proposition 1.3]{Louveau:2005cq} and the ensuing remark.
Examples of \( \bGamma \) as above are the projective classes \( \bSigma^{1}_{n} \)'s and \( \bS ( \kappa ) \), the collection of all \( \kappa \)-Souslin sets --- see Definition~\ref{def:kappaSouslinset}.
However such \( \leq_{\bB} \)-complete \( \mathcal{U}_{ \bGamma} \)'s are usually obtained by \emph{ad hoc} constructions.
In contrast, Theorem~\ref{th:LouveauRosendal} provided the first concrete, natural example of a \( \leq_\bB \)-complete analytic equivalence relation: the relation \( \biembeds^\omega_\CT \) of bi-embeddability on combinatorial trees.%
\footnote{In~\cite{Ferenczi:2009fk} it is shown that the isomorphism relations between separable Banach spaces, and between Polish groups are other natural examples of \( \leq_\bB \)-complete analytic equivalence relations.}
\end{remark}

A \( \leq_\bB \)-complete analytic equivalence relation must contain a non-Borel equivalence class, since it reduces the equivalence relation 
\begin{equation}\label{eq:E_A}
E_A \coloneqq \setofLR{( x , y ) \in \pre{\omega}{2} \times \pre{\omega}{2} }{ x = y \vee x , y \in A} , 
\end{equation}
where \( A \subseteq \pre{\omega}{2} \) is a proper analytic set.
On the other hand, the equivalence classes of an equivalence relation \( E_G \) induced by a continuous (or Borel) action of a Polish group \( G \) are Borel~\cite[Theorem 15.14]{Kechris:1995zt}, even when \( E_G \) is \( \bSigma^{1}_{1} \) and not Borel.
Therefore there is a striking difference between the isomorphism relation \( \cong \) (which is induced by a continuous action of the group \( \Sym ( \omega ) \) of all permutations on the natural numbers) and the embeddability relation \( \embeds \): even a very simple relation like the \( E_A \) above is not Borel reducible to \( \cong \), while any \( \bSigma^{1}_{1} \) equivalence relation (in fact: any \( \bSigma^{1}_{1} \) quasi-order) is Borel reducible to \( \embeds \).
Hjorth isolated a topological property, called \emph{turbulence}, that characterizes when an equivalence relation induced by a continuous Polish group action is Borel reducible to \( \cong \)~\cite{Hjorth:2000zr}.

The next example shows that \( \embeds _{\CT}^\omega \) is also complete for partial orders of size \( \aleph_1 \), in the sense that each such partial order embeds into the quotient order of \( \embeds^\omega_\CT \).

\begin{example}\label{xmp:Parovicenko}
Let \( \preceq \) be the \( \bSigma^{0}_{2} \) quasi-order on \( \pre{ \omega }{2} \) induced by inclusion on \( \pow ( \omega ) / \Fin \), i.e. \( x \preceq y \iff \EXISTS{n}\FORALL{m \geq n } \left ( x ( m ) \leq y ( m ) \right ) \).
Then \( \preceq \) can be embedded into \( \embeds_{\CT}^\omega \) by Theorem~\ref{th:LouveauRosendal}.
By Parovi\v cenko's theorem~\cite{Parovicenko:1963fc}, under the Axiom of Choice \( \AC \) any partial order \( P \) of size \( \aleph_1 \) embeds into \( \pow ( \omega ) / \Fin \), and hence any such \( P \) can be embedded into (the quotient order of) \( \preceq \), and therefore also into (the quotient order of) \( \embeds_{\CT}^\omega \).
\end{example}

Building on~\cite{Louveau:2005cq}, in~\cite{Friedman:2011cr} S.D.~Friedman and the second author strengthened Theorem~\ref{th:LouveauRosendal} by showing that the embeddability relation on countable models is (in the terminology of~\cite{Camerlo:2012kx}) \markdef{invariantly universal},\index[concepts]{invariant universality} i.e.~that the following result holds.
 
\begin{theorem}[S.D.~Friedman-Motto Ros] \label{th:mottorosfriedman}
For every \( \bSigma^1_1 \) quasi-order \( R \) there is an \( \LL_{\omega_1 \omega} \)-sentence \( \upsigma \) such that \( R \sim_{\bB} {\embeds^\omega_\upsigma} \).
\end{theorem}

Invariant universality is a strengthening of \emph{\( \leq_\bB \)-completeness}, which just requires that every analytic quasi-order \( R \) is Borel reducible to \( \embeds \restriction \Mod^\omega_\LL \) (Theorem~\ref{th:LouveauRosendal}). 
Theorem~\ref{th:mottorosfriedman} is obtained by applying the Lopez-Escobar theorem mentioned at the end of Section~\ref{subsubsec:eqrel} to the following purely topological result --- in fact we do not know if there is a direct proof of Theorem~\ref{th:mottorosfriedman} which avoids going through Theorem~\ref{th:mottorosfriedman2}.

\begin{theorem}[S.D.~Friedman-Motto Ros] \label{th:mottorosfriedman2}
For every \( \bSigma^1_1 \) quasi-order \( R \) there is a Borel \( B \subseteq \Mod^\omega_\LL \) closed under isomorphism such that \( R \sim_{\bB} {\embeds \restriction B} \).
\end{theorem}

Various generalizations of Theorem~\ref{th:mottorosfriedman} have already appeared in the literature: for example, in~\cite{Fokina:2011dq} the authors briefly consider its computable version, while~\cite{Motto-Ros:2012ss} presents an extensive analysis of the possible interplay between the isomorphism and the embeddability relation on the same elementary class \( \Mod_{ \upsigma }^ \omega \) for \( \upsigma \) an \( \LL_{ \omega _1 \omega } \)-sentence.

\subsection{What we wanted}\label{subsec:whatwewanted}

The present paper was motivated by the quest for generalizations of Theorem~\ref{th:mottorosfriedman} in two different directions: 
\begin{enumerate-(A)}
\item\label{en:a}
considering quasi-orders belonging to more general pointclasses, such as the projective classes \( \bSigma^1_n \) (for \( n \geq 2 \)) and beyond, and 
\item\label{en:b}
using more liberal (but still definable) kind of reductions.
\end{enumerate-(A)}

Goal~\ref{en:a} is not an idle pursuit, since there are lots of equivalence relations and quasi-orders on Polish spaces naturally arising in mathematics which are not analytic. 
The following are a few examples of this sort.

\begin{example}\label{xmp:analyticquasiorders}
Consider the quasi-order \( ( \mathcal{Q} , \leq_{\bB}) \) of Borel reducibility between analytic quasi-orders.
As observed in~\cite{Louveau:2005cq}, \( \leq_{\bB} \) is a \( \boldsymbol{\Sigma}^1_3 \) relation in the codes for analytic quasi-orders, that is: there is a surjection \( q \colon \pre{\omega}{2} \onto \mathcal{Q} \) such that the relation 
\[
x \preceq_\mathcal{Q} y \IFF q ( x ) \leq_{\bB} q ( y )
\]
is \( \boldsymbol{\Sigma}^1_3 \).
By~\cite{Adams:2000gf}, the restriction of \( \preceq_\mathcal{Q} \) to (the codes for) countable Borel equivalence relations is already a \emph{proper} \( \boldsymbol{\Sigma}^1_2 \) relation.%
\footnote{To the best of our knowledge, the exact topological complexity of the full \( \preceq_{\mathcal{Q}} \) is still unknown.}
\end{example}

\begin{example}
Let \( X \coloneqq ( C ( [ 0 ; 1 ] ) )^\omega \) be the space of countable sequences of continuous, real-valued functions on the unit interval \( [ 0 ; 1] \).
Consider the following natural extension of the inclusion relation on \( X \): given \( \mathcal{F} , \mathcal{G} \in X \), we say that \( \mathcal{F} \) is \emph{essentially contained} in \( \mathcal{G} \) (in symbols \( \mathcal{F} \subseteq^{\lim} \mathcal{G} \)) if each function in \( \mathcal{F} \) can be obtained by recursively applying the pointwise limit operator to the functions in \( \mathcal{G} \). 
This relation can be equivalently described as follows. 
Denote by \( \lim ( \mathcal{F} ) \) the subset of \( C ( [ 0 ; 1 ] ) \) generated by \( \mathcal{F} \in X \) using pointwise limits, i.e.~the smallest subset of \( C ( [ 0 ; 1 ] ) \) containing \( \mathcal{F} \) and closed under the pointwise limit operation: then 
\[ 
\mathcal{F}\subseteq^{\lim} \mathcal{G} \iff \lim ( \mathcal{F} ) \subseteq \lim ( \mathcal{G} ) .
\] 
By a result of Becker (see e.g.~\cite[p.~318]{Kechris:1995zt}), it is easy to see that the relation \( \subseteq^{\lim} \) is a \emph{proper} \( \boldsymbol{\Sigma}^1_2 \) quasi-order.
\end{example}

\begin{example}
Consider the group \( \Aut ( X ) \) of all Borel automorphisms of a standard Borel space \( X \), and the conjugacy relation on it.
The elements of \( \Aut ( X ) \) can be coded as points of the Baire space, and the set \( \mathcal{A} \) of codes is \( \bPi^{1}_{1} \).
Given \( c \in \mathcal{A} \), let \( f_c \in \Aut ( X ) \) be the Borel automorphism coded by \( c \).
Then the equivalence relation \( \setof{ ( c , d ) \in \mathcal{A}^2 }{ f_c \text{ and \( f_d \) are conjugated}} \) is \( \bSigma^1_2 \)-complete~\cite{Clemens:2007yu}.
\end{example}

\begin{example}
The density of a Borel set \( A \subseteq \pre{\omega}{2} \) at a point \( x \in \pre{\omega}{2} \) is the value \( \mathscr{D}_A ( x ) \coloneqq \lim_n \mu ( A \cap \Nbhd^\omega_{ x \restriction n } ) / \mu ( \Nbhd^\omega_{ x \restriction n } ) \), where \( \mu \) is the usual measure on the Cantor space and \( \Nbhd^\omega_{ x \restriction n } \) is the basic open neighborhood determined by \( x \restriction n \) (see Section~\ref{sec:topologies}).
The limit might be undefined for some \( x \), although by the Lebesgue density theorem \( \mathscr{D}_A ( x ) \) is \( 0 \) or \( 1 \) for \( \mu \)-almost every \( x \).
Thus \( \mathscr{D}_A ( x ) \colon \pre{\omega}{2} \to [ 0 ; 1 ] \) is a partial function and we may consider the quasi-order 
\[ 
A \preccurlyeq B \iff \ran ( \mathscr{D}_A ) \subseteq \ran ( \mathscr{D}_B ) .
\]
Since the density at a point \( x \) does not change if the set \( A \) is perturbed by a null set, the quasi-order \( \preccurlyeq \) makes sense on the measure algebra as well.
It can be shown that \( \preccurlyeq \) is a \( \bPi^{1}_{2} \)-complete quasi-order on the measure algebra, and also on \( \mathbf{K} ( \pre{\omega}{2} ) \), the hyperspace of all compact subsets of \( \pre{\omega}{2} \) (see~\cite[Corollary 6.4]{Andretta:2015}).
\end{example}

\begin{example}
Given a metric space \( ( X , d_X ) \), its distance set is \( \ran ( d_X ) \coloneqq \setof{ r \in \R }{ \exists x,y \in X \, d ( x , y ) = r } \). 
A direct computation shows that the quasi-order \( \unlhd \) on the space of (codes for) Polish metric spaces \( \MMM_\omega \) (see Section~\ref{subsubsec:spacesofmetricspaces}) defined by
\[	
X \unlhd Y \iff  \ran ( d_X )  \subseteq \ran ( d_X ) 
\]
is \( \boldsymbol{\Pi}^1_2 \). 
Using~\cite[Theorem 4.5]{Camerlo:2019ke} one can easily show that \( \unlhd \) is in fact \( \boldsymbol{\Pi}^1_2 \)-complete. 
\end{example}

Goal~\ref{en:b}, that is the use of more general reductions, is motivated by the observation that in certain situations \( \leq_\bB \) is inadequate either because it could yield counterintuitive results (Examples~\ref{xmp:E_A} and~\ref{xmp:Calderoni}), or else because the class of quasi-orders is too vast, such as those belonging to some well-behaved inner model (Examples~\ref{xmp:AD} and~\ref{xmp:AC}). 

\begin{example} \label{xmp:E_A}
The equivalence relation \( E_A \) of~\eqref{eq:E_A} is not smooth with respect to Borel functions, although it is concretely classifiable.
It becomes smooth if Borel functions are replaced either by \( \sigma ( \bSigma^1_1 ) \)-measurable functions (that is: functions such that the preimage of an open set belongs to the smallest \( \sigma \)-algebra containing \( \bSigma^{1}_{1} \)), or else by absolutely \( \bDelta^{1}_{2} \)-definable functions (Definition~\ref{def:hjorth}).
\end{example}

\begin{example} \label{xmp:Calderoni}
The relations of bi-embeddability  and isomorphism on the standard Borel space of countable torsion abelian groups are \( \leq_\bB \)-incomparable, while the former is strictly weaker than the latter if \( \bDelta^{1}_{2} \)-definable functions are used, and a Ramsey cardinal is assumed~\cite{Calderoni:2018aa}.
Thus using \( \bDelta^{1}_{2} \)-definable reductions is arguably more appropriate in this situation. 
\end{example}

\begin{example} \label{xmp:AD}
Consider reductions belonging to \( \Ll ( \R ) \), the smallest inner model containing all the reals, and let \( \leq_{\Ll ( \R )} \) be the resulting reducibility relation between the quasi-orders of \( \Ll ( \R ) \).
Assume \( \AD^{\Ll ( \R )} \), i.e.~that all games with payoff set in \( \Ll ( \R ) \) are determined --- this assumption follows from strong forcing axioms (such as the Proper Forcing Axiom \( \mathsf{PFA} \) \cite{Steel:2005ab}), or from the existence of sufficiently large cardinals (e.g.~infinitely many Woodin cardinals with a measurable above them).
Then many of the results on Borel reducibility between analytic binary relations can be extended to \( \Ll ( \R ) \) by replacing \( \leq_{\bB} \) with \( \leq_{\Ll ( \R )} \), including e.g.\ the dichotomies of Silver and Glimm-Effros and the theory of turbulence --- see~\cite{Hjorth:1995ve,Hjorth:1999bh, Hjorth:2000zr}.%
\footnote{As pointed out in~\cite{Hjorth:1995ve}, the generalization of Silver's dichotomy is due to Woodin.}
In this paper we shall repeatedly use the following remarkable result~\cite[Theorem 9.18]{Hjorth:2000zr}.
\begin{theorem}[Hjorth]\label{thm:Hjorth}
If \( E_G \) is the orbit equivalence relation induced by a Polish group \( G \) acting in a turbulent way on a Polish space \( X \) (and hence \( E_G \) is \( \bSigma^1_1 \)), then 
\[
\Ll ( \R ) \models E_G \nleq_{\Ll ( \R )} {\cong \restriction \Mod^\kappa_\LL} 
\]
i.e.~\( E_G \) is not reducible in \( \Ll ( \R ) \) to the isomorphism relation \( \cong \) on the space of all \( \LL \)-structures of size \( \kappa \), for \emph{any countable language \( \LL \)} and \emph{any cardinal \( \kappa \)}.
\end{theorem}
\end{example}

\begin{example} \label{xmp:AC}
The Silver and the Glimm-Effros dichotomies have also been generalized to the \( \ZFC \)-world by considering \( \mathrm{OD}(\R) \) the inner model of real-ordinal definable sets in the Solovay's model (see~\cite{Stern:1984mw,Kanovei:1997lh,Friedman:2008oq}). 
In this framework, one compares quasi-orders in \( \mathrm{OD}(\R) \) by means of the \( \mathrm{OD}(\R) \)-reducibility \( \leq_{\mathrm{OD}(\R)} \), namely the reducibility notion obtained by considering reductions in \( \mathrm{OD}(\R) \).
\end{example}

\begin{remark} \label{rmk:ADvsAC}
Even though both Examples~\ref{xmp:AD} and~\ref{xmp:AC} mainly concern generalizations of the same dichotomies for equivalence relations, the methods used to obtain them are quite different: in the former case the extensive knowledge of models of the Axiom of Determinacy \( \AD \) is used, while in the latter forcing arguments together with absoluteness considerations are employed.
\end{remark}

Since both~\cite{Louveau:2005cq} and~\cite{Friedman:2011cr} exploit the fact that a set is \( \bSigma^1_1 \) if and only if it is \( \omega \)-Souslin, one can ask what sort of generalizations of Theorem~\ref{th:mottorosfriedman} in the direction of goal~\ref{en:a} could be attained. 
Rather than looking at the more familiar projective classes \( \bSigma^1_n \) mentioned in~\ref{en:a}, as common sense would probably suggest, working with the pointclasses \( \bS ( \kappa ) \) of \( \kappa \)-Souslin sets for \( \kappa \) an uncountable cardinal (see Section~\ref{sec:Ksouslinsets}) turns out to be the right move.%
\footnote{The idea of considering \( \kappa \)-Souslin quasi-orders is not new: in fact, \( \kappa \)-Souslin (and co-\( \kappa \)-Souslin) quasi-orders on Polish spaces have been already extensively studied in the literature, see e.g.\ \cite{Harrington:1982if, Shelah:1984kc, Kanovei:1997lh}.} 
This approach is not too restrictive, since e.g.\ under \( \AD \) one has \( \bSigma^{1}_{n} = \bS ( \kappa_n ) \) for an appropriate cardinal \( \kappa_n \).

Remark~\ref{rmk:ADvsAC} seems to suggest that if one aims at generalizing Theorem~\ref{th:mottorosfriedman} to both the \( \AC \)- and the \( \AD \)-world, then different methods should be used.
This prompts the question of which one of the two approaches is more promising for our purposes. 
On the one hand, the fact that we are going to study \( \kappa \)-Souslin quasi-orders seems to indicate that an \( \AD \)-approach similar to the one of Example~\ref{xmp:AD} should be preferred: in fact \( \AD \) imposes an extremely rich structure on the subsets of Polish spaces and provides a well-developed general theory of \( \kappa \)-Souslin sets (for quite large cardinals \( \kappa \)), while we have very little information on the structure and properties of \( \kappa \)-Souslin sets in the context of \( \AC \), where the notion of \( \kappa \)-Souslin is nontrivial only when \( \kappa \) is smaller than the continuum \( 2^{\aleph_0} \).
On the other hand, any generalization of Theorem~\ref{th:mottorosfriedman} to \( \kappa \)-Souslin quasi-orders seems to require replacing the space \( \Mod^ \omega _\LL \) with \( \Mod^ \kappa _\LL \), and the logic \( \LL_{\omega_1 \omega} \) with \( \LL_{ \kappa ^+ \kappa } \).
Unfortunately, there are two roadblocks down this path:
\begin{itemize}[leftmargin=1pc]
\item
a decent descriptive set theory on spaces like \( \Mod^\kappa_\LL \), which can be identified with the generalized Cantor space \( \pre{ \kappa}{2} \), seems to require cardinal arithmetic assumptions such as \( \card{ \pre{ < \kappa}{ \kappa}} = \kappa \) --- see e.g.\ the generalization of the Lopez-Escobar theorem in Section~\ref{subsec:LopezEscobar}, and~\cite{Vaught:1974kl, Mekler:1993kh, Friedman:2011nx} for more on these matters;
\item
the classical analysis of the logics \( \LL_{\kappa^+ \kappa} \) essentially requires the full Axiom of Choice \( \AC \).
\end{itemize}
Since both \( \AC \) and \( \card{ \pre{ < \kappa}{ \kappa}} = \kappa > \omega \) contradict the Axiom of Determinacy, any generalization of Theorem~\ref{th:mottorosfriedman} under \( \AD \) seems out of reach.
But even under \( \AC \), the cardinal condition \( \kappa^{< \kappa} = \kappa > \omega \) cannot be achieved when \( \kappa < 2^{ \aleph_0} \), and the latter is needed to guarantee that the notion of \( \kappa \)-Sousliness is not trivial. 

As there is a strong tension between the possible scenarios, any generalization of Theorem~\ref{th:mottorosfriedman} may thus seem hopeless. 

\subsection{What we did}
Despite the bleak outlook depicted in the previous section, we managed to prove some generalizations of the Louveau-Rosendal completeness result (Theorem~\ref{th:LouveauRosendal}) and of its strengthening to invariant universality (Theorem~\ref{th:mottorosfriedman}).
Below is a collection of results that will be proved in Section~\ref{sec:applications}.
We would like to emphasize that they are not proved by \emph{ad hoc} methods, but they all follow from the constructions and techniques of Section~\ref{sec:mainconstruction}--\ref{sec:alternativeapproach}.
This might serve as a partial justification for the length of this paper.

The central idea is that a quasi order \( R \) on a Polish space is reduced to the embeddability relation \( \embeds_{\CT}^{ \kappa } \) on the collection \( \CT_ \kappa \) of all combinatorial trees of size \( \kappa \).
By Theorem~\ref{th:LouveauRosendal}, if \( R \) is \( \bSigma^{1}_{1} \), that is \( \kappa \)-Souslin with \( \kappa = \omega \), then \( R \) is reducible to \( \embeds_{\CT}^{ \omega } \); our results show that the higher the complexity of \( R \), the larger the cardinal \( \kappa \) must be taken.
This is reminiscent of a well-known feature in the proofs of determinacy where a game with a complex subset of the \( \pre{\omega}{\omega} \) is transfigured in a closed game on \( \pre{ \omega }{ \kappa } \) with \( \kappa \) large.
Our reduction takes place between a Polish space and (a homeomorphic copy of) \( \pre{ \kappa }{2} \), and the complexity of the reduction will be either \( \kappa + 1 \)-Borel or \( \kappa \)-Souslin-in-the-codes.
We would like to stress that these notions are quite natural when working with \( \pre{ \kappa }{2} \).

\subsubsection*{A miscellanea of results}
For the sake of brevity, in the statements of Theorems~\ref{th:introcompleteness} and~\ref{th:introinvuniversal} the quasi-order \( R \) is tacitly assumed to be defined on some Polish or standard Borel space, and the embeddability relation on \( \Mod^\kappa_\upsigma \), the collection of models of size \( \kappa \) of the \( \LL_{\kappa^+ \kappa} \)-sentence \( \upsigma \), is denoted by \( \embeds^\kappa_\upsigma \).
Recall also that the assumption \( \AD^{\Ll ( \R )} \) used in some of the next results follows from sufficiently large cardinals or strong forcing axioms, such as \( \PFA \).

\begin{theorem}[Completeness results] \label{th:introcompleteness}
\begin{enumerate-(a)}
\item \label{th:introcompleteness-a}
Assume \( \ZFC + \PFA \).
Then \( R \leq_{\bSigma^1_2} {\embeds^{\omega_1}_{\CT}} \) and \( R \leq_{\bB}^{\omega_1} {\embeds^{\omega_1}_{\CT}} \), for every \( \boldsymbol{\Sigma}^1_2 \) quasi-order \( R \) (Corollary~\ref{cor:omega1underPFA}).
\item \label{th:introcompleteness-b}
Assume \( \ZFC + \FORALL{x \in \pre{\omega}{\omega}} ( x^\# \text{ exists} ) \).
Then \( R \leq_{\bS(\omega_2)} {\embeds^{\omega_2}_{\CT}} \) and \( R \leq_{\bB}^{\omega_2} {\embeds^{\omega_2}_{\CT}} \), for every \( \boldsymbol{\Sigma}^1_3 \) quasi-order \( R \) (Theorem~\ref{th:mainSigma13AC}).
\item \label{th:introcompleteness-c}
Assume \( \ZFC + \AD^{\Ll(\R)} \). 
Then \( R \leq_{\bS ( \omega_{ r ( n ) } ) } {\embeds^{\omega_{r(n)}}_\CT} \) and \( R \leq_{\bB}^{\omega_{ r ( n ) } } {\embeds^{\omega_{r(n)}}_\CT} \), for every \( \bSigma^1_n \) quasi-order \( R \), where \( r \colon \omega \to \omega \) is\begin{equation}\label{eq:functionrsection1}
r ( n ) = \begin{cases}
2^{k + 1 } - 2 & \text{if } n = 2 k + 1 ,
\\
2^{k + 1 } - 1 & \text{if } n = 2 k + 2 
\end{cases}
\end{equation}
(Theorem~\ref{th:embeddingsofprojectiveqos}).
\item \label{th:introcompleteness-d}
Assume \( \ZF + \DC + \AD \). 
For \( n > 0 \), let \( \kappa_n \coloneqq \blambda^1_n \) if \( n \) is odd and \( \kappa_n \coloneqq \bdelta^1_{n-1} \) if \( n \) is even.%
\footnote{See Section~\ref{subsec:projectiveordinals} for the relevant definitions.
In particular, \( \kappa_1 = \omega \), \( \kappa_2 = \omega_1 \), \( \kappa_3 = \aleph_\omega \), and \( \kappa_4 = \aleph_{\omega+1} \).}
Then \( R \leq_{\bSigma^1_n} {\embeds^{\kappa_n}_\CT} \) (and hence also \( R \leq_{\bB}^{\kappa_n} {\embeds^{\kappa_n}_\CT} \)) for every \( \bSigma^1_n \) quasi-order \( R \) (Theorem~\ref{th:mainprojective}).
\item \label{th:introcompleteness-e}
Assume \( \ZFC + \AD^{\Ll ( \R )} \). 
Then \( R \leq_{\Ll ( \R )} {\embeds^\kappa_\CT } \) for every \( \boldsymbol{\Sigma}^2_1 \) quasi-order \( R \) of \( \Ll ( \R ) \), where \( \kappa \coloneqq \boldsymbol{\delta}^2_1 \) is defined as in Section~\ref{subsec:projectiveordinals} (Theorem~\ref{th:largecardinals1}).
\end{enumerate-(a)}
\end{theorem}

\begin{remarks} \label{rmk:introdef}
\begin{enumerate-(i)}
\item \label{rmk:introdef-i}
In Theorem~\ref{th:introcompleteness}, \( \leq^ \kappa _{\bB} \) is the generalization of \( \leq_{\bB} \) to \( \kappa \) an uncountable cardinal (Definition~\ref{def:Borelreducibility}). 
The reducibilities \( \leq_{\bSigma^1_n} \) and \( \leq_{\bS ( \kappa )} \) are instead the analogue of \( \leq_{\bB} \) where the reducing functions are required to be, respectively, \( \bSigma^1_n \)-in-the-codes and \( \bS ( \kappa ) \)-in-the-codes (see Definition~\ref{def:Gammainthecodes}); with this notation, the standard Borel reducibility \( \leq_\bB \) would be denoted by \( \leq_{\bSigma^1_1} \) and \( \leq_{ \bS ( \omega )} \), respectively.
Finally, \( \leq_{\Ll ( \R )} \) is reducibility in \( \Ll(\R) \) --- see Example~\ref{xmp:AD}.
\item \label{rmk:introdef-ii}
The various statements in Theorem~\ref{th:introcompleteness} can be seen as completeness results generalizing Theorem~\ref{th:LouveauRosendal} to pointclasses \( \bGamma \) properly extending \( \bSigma^{1}_{1} \).
The move to such \( \bGamma \)'s forces us to replace the Polish space \( \CT_ \omega \) with \( \CT_ \kappa \) (for some \( \kappa > \omega \)), which is homeomorphic to (a closed subset of) the generalized Cantor space \( \pre{ \kappa }{2} \), and hence far from being Polish.
Therefore part~\ref{th:LouveauRosendal-i} of Theorem~\ref{th:LouveauRosendal} must necessarily be dropped in such generalizations (see also Remark~\ref{rmk:completeness}).
\item \label{rmk:introdef-iii}
The reader familiar with Jackson's analysis of the regular cardinals below the projective cardinals in \( \Ll ( \R ) \) will immediately see that parts~\ref{th:introcompleteness-c} and~\ref{th:introcompleteness-d} of Theorem~\ref{th:introcompleteness} are strictly related.
\end{enumerate-(i)}
\end{remarks}

The relations \( \sim^\kappa_\bB \) and \( \sim_{\Ll ( \R )} \) appearing in the next theorem are the bi-reducibility relations canonically associated to \( \leq^\kappa_\bB \) (see Remark~\ref{rmk:introdef}\ref{rmk:introdef-i}) and \( \leq_{\Ll ( \R )} \) (see Example~\ref{xmp:AD}), respectively. 

\begin{theorem}[Invariant universality results] \label{th:introinvuniversal}
\begin{enumerate-(a)}
\item \label{th:introinvuniversal-a}
Assume \( \ZF + \AC_\omega ( \R ) \). 
Then for every \( \boldsymbol{\Sigma}^1_{2} \) quasi-order \( R \) there is an \( \LL_{\omega_2 \, \omega_1} \)-sentence \( \upsigma \) such that \( R \sim^{\omega_1}_{\bB} {\embeds ^{\omega_1}_\upsigma} \) (Theorem~\ref{th:finalSigma12-b}).
\item \label{th:introinvuniversal-b}
Assume \( \ZFC + \FORALL{x \in \pre{\omega}{\omega}} ( x^\# \text{ exists} ) \).
Then for every \( \boldsymbol{\Sigma}^1_{3} \) quasi-order \( R \) there is an \( \LL_{\omega_3 \, \omega_2} \)-sentence \( \upsigma \) such that \( R \sim^{\omega_2}_{\bB} {\embeds ^{\omega_2}_\upsigma} \) (Theorem~\ref{th:mainSigma13AC}).
\item \label{th:introinvuniversal-c}
Assume \( \ZFC + \AD^{\Ll ( \R ) } \), and let \( r \colon \omega \to \omega \) be as in equation~\eqref{eq:functionrsection1}.
Then for every \( \boldsymbol{\Sigma}^1_{n} \) quasi-order \( R \) there is an \( \LL_{\omega_{ r ( n ) + 1 } \, \omega_{ r ( n ) }} \)-sentence \( \upsigma \) such that \( R \sim^{\omega_{ r ( n ) }}_{\bB} {\embeds ^{\omega_{ r ( n ) }}_\upsigma} \) (Theorem~\ref{th:embeddingsofprojectiveqos}).
\item \label{th:introinvuniversal-d}
Assume \( \ZF + \DC + \AD \). 
Let \( n \neq 0 \) be an even number.
Then for every \( \boldsymbol{\Sigma}^1_{n} \) quasi-order \( R \) there is an \( \LL_{\bdelta^1_n \, \bdelta^1_{n-1}} \)-sentence%
\footnote{See Section~\ref{subsec:projectiveordinals} for the definition of the projective ordinals \( \bdelta^1_n \), and recall that when \( n \) is even we have \( \bdelta^1_n = ( \bdelta^1_{n - 1} )^+ \) (under the assumption \( \ZF + \DC + \AD \)).}
 \( \upsigma \) such that \( R \sim^{\bdelta^1_{n - 1}}_{\bB} {\embeds ^{\bdelta^1_{ n - 1}}_\upsigma} \) (Theorem~\ref{th:mainprojective}).
\item \label{th:introinvuniversal-e}
Assume \( \ZFC + \AD^{\Ll ( \R )} \). 
Then for every \( \boldsymbol{\Sigma}^2_1 \) quasi-order \( R \) of \( \Ll ( \R ) \) there is an \( \LL_{\kappa^+ \kappa} \)-sentence \( \upsigma \) belonging to \( \Ll ( \R ) \) such that \( R \sim_{\Ll ( \R )} {\embeds^\kappa_\upsigma} \), where as in Theorem~\ref{th:introcompleteness}\ref{th:introcompleteness-e} we set \( \kappa \coloneqq \boldsymbol{\delta}^2_1 \) (Theorem~\ref{th:largecardinals1}).
\end{enumerate-(a)}
\end{theorem}

By considering the equivalence relation associated to a quasi-order, the results above can be turned into statements concerning equivalence relations on Polish or standard Borel spaces and the bi-embeddability relation \( \biembeds^\kappa_\CT \) on combinatorial trees of size \( \kappa \). 
The stark difference between bi-embeddability \( \biembeds \) and isomorphism \( \cong \) on countable models uncovered by Theorem~\ref{th:LouveauRosendal} (see the observation after Remark~\ref{rmk:introLR}) is also present in the uncountable case: assuming \( \AD^{\Ll ( \R ) } \), by Theorem ~\ref{th:introcompleteness}\ref{th:introcompleteness-e} any \( \bSigma^{2}_{1} \) equivalence relation of \( \Ll ( \R ) \) is \( \leq_{\Ll ( \R )} \)-reducible to \( \biembeds^{\bdelta^2_1}_\LL \), the bi-embeddability relation on \( \Mod^{\bdelta^2_1}_\LL \), while this badly fails if \( \biembeds \) is replaced by \( \cong \), by Hjorth's Theorem~\ref{thm:Hjorth}.

By applying Theorems~\ref{th:introcompleteness} and~\ref{th:introinvuniversal} we also get some information on the \( \bSigma^1_3 \) quasi-order \( ( \mathcal{Q} , \leq_{\bB} ) \) of Example~\ref{xmp:analyticquasiorders}. 
Such quasi-order may be seen as a (definable) embeddability relation between structures of size the continuum.
It turns out that, under suitable hypotheses, \( ( \mathcal{Q} , \leq_{\bB} ) \) can be turned into an embeddability relation between \emph{well-ordered} structures of size potentially smaller than the continuum.
Indeed, in models with choice we have:

\begin{th:Qembedsundersharps}
Assume \( \ZFC + \FORALL{x \in \pre{\omega}{\omega}} ( x^\# \text{ exists}) \).
Then the quotient order of \( ( \mathcal{Q} , \leq_{\bB} ) \) (definably) embeds into the quotient order of \( \embeds_\CT^{\aleph_ 2 } \).
Moreover, there is an \( \LL_{ \aleph _3 \, \aleph_2 } \)-sentence \( \upsigma \) such that the quotient orders of \( ( \mathcal{Q} , \leq_{\bB} ) \) and \( \embeds_ \upsigma^{\aleph_2 } \) are (definably) isomorphic.
\end{th:Qembedsundersharps}

On the other hand, in the determinacy world we get:

\begin{th:QembedsintoCTalephomega}
Assume \( \ZF + \DC + \AD \).
Then the quotient order of \( ( \mathcal{Q} , \leq_{\bB} ) \) (definably) embeds into the quotient order of \( \embeds_\CT^{\aleph_ \omega } \).
Moreover, there is an \( \LL_{\aleph_{ \omega + 1} \aleph_ \omega } \)-sentence \( \upsigma \) such that the quotient orders of \( ( \mathcal{Q} , \leq_{\bB} ) \) and \( \embeds_ \upsigma^{\aleph_ \omega } \) are (definably) isomorphic.
\end{th:QembedsintoCTalephomega}

In Section~\ref{subsec:Representingarbitrarypartialorders} we generalize the combinatorial completeness property of \( \embeds^\omega_\CT \) described in Example~\ref{xmp:Parovicenko} to the uncountable case, albeit in a weaker form, in both the \( \AC \)- and the \( \AD \)-world.

\begin{prop:parovicenkoAC} 
Assume \( \ZFC \) and let \( \omega < \kappa \leq 2^{\aleph_0} \). 
Then every partial order \( P \) of size \( \kappa \) can be embedded into the quotient order of \( \embeds_\CT^\kappa \). 
In fact, for every such \( P \) there is an \( \LL_{\kappa^+ \kappa} \)-sentence \( \upsigma \) (all of whose models are combinatorial trees) such that the quotient order of \( \embeds^{\kappa}_\upsigma \) is isomorphic to \( P \).
\end{prop:parovicenkoAC}

\begin{prop:parovicenkoACimproved}
Assume \( \ZFC \) and let \( \omega < \kappa \leq 2^{\aleph_0} \). 
Then \( {\subseteq^*_\kappa} \leq^\kappa_\bB {\embeds_\CT^\kappa} \), where \( \subseteq^*_\kappa \) is the relation on \( \pow(\kappa) \) of inclusion modulo bounded subsets. 
In particular, every linear order of size \( \aleph_{n+1} \) can be embedded into the quotient order of \( \embeds_\CT^{\aleph_n} \) (whenever \( 2^{\aleph_0} \geq \aleph_n \)).
\end{prop:parovicenkoACimproved}

\begin{th:parovicenko1}
Assume \( \ZF+ \DC + \AD \).
Let \( \kappa \) be a Souslin cardinal. 
Then every partial order \( P \) of size \( \kappa \) can be embedded into the quotient order of \( \embeds_\CT^\kappa \). 
In fact, if \( \kappa < \bdelta_{ \bS ( \kappa ) } \), then for every such \( P \) there is an \( \LL_{\kappa^+ \kappa} \)-sentence \( \upsigma \) (all of whose models are combinatorial trees) such that the quotient order of \( \embeds^{\kappa}_\upsigma \) is isomorphic to \( P \).
\end{th:parovicenko1}

Notice that by Proposition~\ref{prop:deltaSkappa}\ref{prop:deltaSkappa-c} we can apply the second part of Theorem~\ref{th:parovicenko1} to any Souslin cardinal \( \kappa \) (in a model of \( \ZF + \AD + \DC \)) which is not a regular limit of Souslin cardinals: in particular, we can take \( \kappa \) to be one of the projective ordinals \( \bdelta^1_n \).

Finally, in Section~\ref{subsec:othermodeltheoretic} we show that in all our completeness results (including the ones mentioned so far in this introduction) one may freely replace embeddings and combinatorial trees with other kinds of morphisms and structures which are relevant to graph theory and model theory, namely:
\begin{itemize}[leftmargin=1pc]
\item
we can consider full homomorphisms between graphs (or even just combinatorial trees);
\item
we may also consider embeddings between (complete) lattices, where the latter may indifferently be construed as partial orders or as bounded lattices in the algebraic sense.
\end{itemize}

\subsection{How we proved it}
In order to overcome the difficulties explained at the end of Section~\ref{subsec:whatwewanted}, we forwent the approach followed in the countable case: this would have required to first analyze the descriptive set theory of \( \Mod^ \kappa_\LL \) in order to achieve a generalization of Theorem~\ref{th:mottorosfriedman2}, and then obtain as a corollary the corresponding generalization of Theorem~\ref{th:mottorosfriedman}.
The key idea is to reverse this approach and to exploit the greater expressive power of the logic \( \LL_{\kappa^+ \kappa} \) when \( \kappa \) is uncountable and directly extend Theorem~\ref{th:mottorosfriedman} as follows: to each \( \kappa \)-Souslin quasi-order \( R \) on the Cantor space \( \pre{\omega}{2} \) we associate an \( \LL_{\kappa^+ \kappa} \)-sentence \( \upsigma \) (all of whose models are combinatorial trees) so that \( R \) is bi-reducible to \( {\embeds}^\kappa_{\upsigma} \), in symbols \( R \sim {\embeds^\kappa_\upsigma} \) (Theorem~\ref{th:main}).
The corresponding topological version is obtained as a corollary (using a restricted version of the generalized Lopez-Escobar theorem) by using formul\ae{} belonging to a sufficiently powerful fragment \( \LL_{\kappa^+ \kappa}^b \) of \( \LL_{\kappa^+ \kappa} \) (\( b \) is for \emph{bounded}).
To be more precise: given a tree \( T \) on \( 2 \times 2 \times \kappa \) such that \( R = \PROJ \body{T} \) is a quasi-order, we shall construct in \( \ZF \) 
\begin{itemize}[leftmargin=1pc]
\item
a function \( f_T \colon \pre{\omega}{2} \to \CT_ \kappa \) (see~\eqref{eq:f_T}) such that \( f_T \) reduces the quasi-order \( R \) to \( \embeds_\CT^\kappa \) and satisfies \( f_T ( x ) \cong f_T ( y ) \IFF x = y \) (Theorem~\ref{th:graphs});
\item
an \( \LL_{ \kappa ^+ \kappa } \)-sentence \( \upsigma = \upsigma_T \) and a function \( h_T \colon \Mod_{\upsigma}^ \kappa \to \pre{\omega}{2} \) such that \( \Mod_{\upsigma}^\kappa \) is the closure under isomorphism of the range of \( f_T \), and \( h_T \) reduces \( \embeds_ \CT^\kappa \) to \( R \) (Corollary~\ref{cor:inverse}).
\end{itemize}
This construction depends only on the cardinal \( \kappa > \omega \) and on the chosen witness \( T \) of the fact that \( R \) is \( \kappa \)-Souslin.
This means that the desired \( \LL_{\kappa^+ \kappa} \)-sentence \( \upsigma \) and the two reductions \( f_T , h_T \) witnessing \( R \sim {\embeds^\kappa_\upsigma} \) can be found in every inner model containing \( \kappa \) and \( T \), and that they are in fact (essentially) the same in all these models.
Moreover, the reductions involved are absolute, as they can be defined by formul\ae{} (in the language of set theory) which define reductions between \( R = \PROJ \body{T} \) and \( \embeds ^ \kappa _ \upsigma \) in every generic extension and in every inner model of the universe of sets \( \Vv \) we started with. 
Such reductions have essentially the same topological complexity of the quasi-order \( R \): for example, if \( \kappa \) is regular then they are \( \kappa + 1 \)-Borel, a natural extension of the classical notion of a Borel function (see Definition~\ref{def:weaklyBorel}). 
Therefore our results provide natural generalizations of both Theorems~\ref{th:mottorosfriedman} and~\ref{th:mottorosfriedman2} (which correspond to the basic case \( \kappa = \omega \)) to uncountable \( \kappa \)'s --- see Theorem~\ref{th:maintopology} and~\ref{th:barmaintopology}, respectively. 
All these observations lead us to consider also the more general notion of \emph{definable cardinality}, which is strictly related to the notion of \emph{definable reducibility} --- see Section~\ref{sec:definablecardinality} for a more thorough discussion on the genesis and relevance of these concepts. 

It is worth pointing out that in order to find applications in both the \( \AC \)-world and the \( \AD \)-world, the above mentioned preliminary completeness and invariant universality results for \( \kappa \)-Souslin quasi-orders (for \( \kappa \) an uncountable cardinal) \emph{are developed in \( \ZF \)}, so that they can be applied to \emph{all} situations in which \( \bS ( \kappa ) \) is a nontrivial pointclass. 
These include, among many others, the cases of the projective pointclasses \( \bSigma^1_n \), of the \( \bSigma^2_1 \) sets of \( \Ll ( \R ) \), and even of the entire \( \pow ( \R ) \) (under various determinacy or large cardinal assumptions).

Let us conclude this section with a more general comment. 
There is a common phenomenon in model theory: the study of uncountable models requires ideas and techniques quite different from the ones used in the study of countable models, and this paper is no exception.
As already observed in Section~\ref{subsec:whatwewanted} and at the beginning of this subsection, the techniques used in~\cite{Friedman:2011cr} for the countable case cannot be transferred to the uncountable case. 
Conversely, the arguments contained in this paper cannot be adapted to include the countable case, since we essentially need to be able to use a single \( \LL_{\kappa^+ \kappa} \)-sentence either to assert the existence of certain \emph{infinite} (but still small) substructures (see Section~\ref{sec:invariantlyuniversal}), or to express the \emph{well-foundedness} of a binary relation (see Section~\ref{sec:alternativeapproach}).
As it is well-known, none of these possibilities can be achieved using the logic \( \LL_{\omega_1 \omega} \).

\subsection{Classification of non-separable structures up to bi-embeddability}
As recalled at the end of Section~\ref{subsubsec:eqrel}, the problem of classifying e.g.\ Polish metric spaces up to isometry or separable Banach spaces up to linear isometry have been widely studied in the literature, but very little was known about the analogous classification problems up to bi-embeddability.
In~\cite{Louveau:2005cq} Louveau and Rosendal used their Theorem~\ref{th:LouveauRosendal} to show that many kinds of separable spaces are in fact not classifiable (in any reasonable way) with respect to the relevant notion of bi-embeddability. 
This is made precise by the following result. 
(Notice that the relations appearing in the next theorem can all be construed as analytic quasi-orders on corresponding standard Borel spaces --- see~\cite{Louveau:2005cq} for more details.)

\begin{theorem}[Louveau-Rosendal] \label{th:LouveauRosendal2}
The following relations are \emph{\( \leq_\bB \)-complete} analytic quasi-orders:
\begin{enumerate-(a)}
\item
continuous embeddability between compacta (\cite[Theorem 4.5]{Louveau:2005cq});
\item
isometric embeddability between (ultrametric or discrete) Polish metric spaces (\cite[Propositions 4.1 and 4.2]{Louveau:2005cq});
\item
linear isometric embeddability between separable Banach spaces (\cite[Theorem 4.6]{Louveau:2005cq}).
\end{enumerate-(a)}
As a consequence, the corresponding bi-embeddability relations are \emph{\( \leq_\bB \)-complete} analytic equivalence relations.
\end{theorem}

(Theorem~\ref{th:LouveauRosendal2} has been further improved in~\cite{Camerlo:2012kx}, where it is shown that it is possible to obtain analogues of Theorem~\ref{th:mottorosfriedman2} in which \( \embeds^\omega_\CT \) is replaced by e.g.\ any of the embeddability relations mentioned in Theorem~\ref{th:LouveauRosendal2}.)

Following this line of research, in Sections~\ref{subsec:isomemb} and~\ref{subsec:linearisomemb} we study the complexity of various classification problems for non-separable spaces. 
In particular, we show that in all our completeness results one may systematically replace the embeddability relation between combinatorial trees of size \( \kappa \) with e.g.\ the isometric embeddability relation between (discrete and/or ultrametric) complete metric spaces of density character \( \kappa \). 
Theorem~\ref{th:introultrametricapplications} below offers a sample of the results that may be obtained in this way (see also Theorem~\ref{th:ultrametricapplications} and Remark~\ref{rmk:variants}.)
Before stating it let us point out that \( \sqsubseteq^i \) is the relation of isometric embeddability between metric spaces, and that our use of ``completeness'' is a tad nonstandard --- see Remark~\ref{rmk:completeness}. 

\begin{theorem}\label{th:introultrametricapplications} 
\begin{enumerate-(a)}
\item
Assume \( \ZFC \). 
Then the relation \( \sqsubseteq^i \) between (discrete) complete (ultra)metric spaces of density character \( \omega_1 \) is \( \leq^{\omega_1}_\bB \)-complete for \( \bSigma^1_2 \) quasi-orders on standard Borel spaces.
If moreover we assume either \( \AD^{\Ll ( \R )} \) or \( \MA + \neg \CH + \EXISTS{ a \in \pre{\omega}{\omega}} (\omega_1^{\Ll [ a ] } = \omega_1) \), then such relation is also \( \leq_{\bSigma^1_2} \)-complete for \( \bSigma^1_2 \) quasi-orders on standard Borel spaces.
\item
Assume \( \ZFC + \FORALL{x \in \pre{\omega}{\omega} } ( x^\# \text{ exists} ) \).
The relation \( \sqsubseteq^i \) between (discrete) complete (ultra)metric spaces of density character \( \omega_2 \) is \( \leq_{\bB}^{\omega_2} \)-complete for \( \bSigma^1_3 \) quasi-orders on standard Borel spaces. 
\item
Assume \( \ZFC + \AD^{\Ll ( \R )} \). 
Then the relation \( \sqsubseteq^i \) between (discrete) complete (ultra)metric spaces of density character \( \omega_{r ( n ) } \) is \( \leq_{\bB}^{\omega_{r ( n )}} \)-complete for \( \bSigma^1_n \) quasi-orders on standard Borel spaces, where \( r \) is as in equation~\eqref{eq:functionrsection1}.
\item
Assume \( \ZF + \DC + \AD \). 
For \( 0 \neq n \in \omega \), let \( \kappa_n \) be as in Theorem~\ref{th:introcompleteness}\ref{th:introcompleteness-d}.
The relation \( \sqsubseteq^i \) between (discrete) complete (ultra)metric spaces of density character \(\kappa_n \) is both \( \leq^{\kappa_n}_{\bB} \)-complete and \( \leq_{\bSigma^1_n} \)-complete for \( \bSigma^1_n \) quasi-orders on standard Borel spaces. 
\end{enumerate-(a)}
\end{theorem}

In particular, it follows from Theorem~\ref{th:introultrametricapplications} that the problem of classifying non-separable complete metric spaces up to isometric bi-embeddability is extremely complex. 
We show that for uncountable \( \kappa \)'s, both the isometry relation and the isometric bi-embeddability relation between (discrete) complete metric spaces of density character \( \kappa \) may consistently have maximal complexity with respect to the relevant reducibility notion \( \leq^\kappa_\bB \) (Theorem~\ref{thm:completenessappendix}\ref{thm:completenessappendix-a}). 
We also show that in models of \( \AD^{\Ll(\R)} \) the relation of isometric bi-embeddability between ultrametric or discrete complete metric spaces of a given uncountable density character is way more complex (with respect to \( \Ll(\R) \)-reducibility) than the isometry relation on the same class: the former \( \leq_{\Ll(\R)} \)-reduces, among others, all \( \bSigma^{1}_{2} \) equivalence relations on a Polish or standard Borel space, while the latter cannot even \( \leq_{\Ll(\R)} \)-reduce all \( \bSigma^1_1 \) equivalence relations (see the comment after Theorem~\ref{th:upperboundsultrametric}). 
These observations show in particular that it is independent of \( \ZF + \DC \) whether e.g.\ the relation of isometric bi-embeddability between ultrametric (respectively, discrete) complete metric spaces of density character \( \omega_2 \) is \( \leq^{\omega_2}_{\bB} \)-reducible to the isometry relation on the same class of spaces (Corollary~\ref{cor:appindependence}).

Finally, in Section~\ref{subsec:linearisomemb} we show that in all previously mentioned results (including e.g.\ Theorem~\ref{th:introultrametricapplications}) one can further replace the isometric embeddability relation \( \sqsubseteq^i \) between complete metric spaces of density character \( \kappa \) with the linear isometric embeddability relation \( \sqsubseteq^{li} \) between Banach spaces of density \( \kappa \). 
Thus also the problem of classifying non-separable Banach spaces up to linear isometry or linear isometric bi-embeddability is very complex --- in fact these equivalence relations may consistently have maximal complexity as well (Corollary~\ref{cor:appindependence}).

In our opinion, these (anti-)classification results for non-separable complete metric spaces and Banach spaces constitute, together with the tight connections between \( \kappa+1 \)-Borel reducibility \( \leq^\kappa_\bB \) and Shelah's stability theory uncovered in~\cite{Friedman:2011nx}, some of the strongest motivations for pursuing the research in the currently fast-growing field of generalized descriptive set theory.

\subsection{Organization of the paper, or: How (not) to read this paper}

The aspiration of this paper, besides that of presenting new results, is to serve as a basic reference text for works dealing with generalized descriptive set theory, and at the same time to be as self-contained as possible. 
For this reason, in Sections~\ref{sec:preliminaries}--\ref{sec:Ksouslinsets} we collected all relevant definitions and surveyed all basic results involved in the rest of the work. 
These sections contain many well-known folklore facts or elementary observations which unfortunately, to the best of our knowledge, cannot be found in a unitary and organic presentation elsewhere, together with some new results which may be of independent interest. 

Due to the applications we had in mind, we developed these preliminaries trying to minimize the amount of set-theoretic assumptions or cardinal conditions required to prove them: this often forced us to find new proofs and ideas which, albeit slightly more involved than the ``classical'' ones, allowed us to work in wider or less standard contexts. 
As a byproduct, we obtain independence results concerning the generalized Cantor space \( \pre{\omega_1}{2} \) which are invisible when working under \( \ZFC + \CH \), the most popular setup in the current literature to deal with such space (see e.g.\ Remark~\ref{rem:cantorbasicproperties}\ref{rem:cantorbasicproperties-ind} or the comment after Proposition~\ref{prop:nonhomeomorphism}). 

For the reader's convenience, Sections~\ref{sec:topologies}--\ref{sec:generalizedBorelfunctions} are organized as follows: all the definitions and statements crucial for the rest of the paper are collected in its first subsections, while their proofs, a more detailed and thorough discussion on the subject under consideration, together with additional information and results, are postponed to the later subsections, which are marked with a \( * \).
\emph{Readers only interested in the results on the embeddability relation and in their applications, may safely skip these starred subsections.} 

The authors are aware of Bertrand Russell's aphorism: ``A book should have either intelligibility or correctness: to combine the two is impossible, \ldots''.
Although we tried our best to avoid writing incorrect statements, whenever a choice had to be made between intelligibility and correctness, we opted for the former.
We are also very intimidated by the second part of such aphorism: ``\ldots but to lack both is to be unworthy of such a place as Euclid has occupied in education'', hoping we did not fail so badly in this endeavor.

\subsection{Annotated content}
Here is a synopsis of the sections of the paper and their content.

\subsection*{Section~\ref{sec:preliminaries}: Preliminaries and notation.} 
Here we collect all the basic notions and facts that are taken for granted in this paper.

\subsection*{Section~\ref{sec:topologies}: The generalized Cantor space.} 
We study the basic properties of the generalized Cantor space \( \pre{\kappa}{2} \) endowed with various topologies. 
We consider both the \emph{bounded topology \( \tau_b \)} and the \emph{product topology \( \tau_p \)}, as well as many other intermediate topologies. 
We also introduce in this context the notions of Lipschitz and continuous reducibility together with the corresponding (long) reduction games, and use game theoretic arguments to prove (without any determinacy assumption) some technical properties that are used in later sections.

\subsection*{Section~\ref{sec:borelsets}: Generalized Borel sets.} 
We introduce the collection of \( \alpha \)-Borel subsets (and their effective counterpart) of a given topological space \( X \) as a generalization to higher cardinals of the classical notion of a Borel set. 
The notion of \( \alpha \)-Borel subset of \( \R \) has been thoroughly studied under determinacy assumptions in the Seventies-Eighties (see Cabal volumes), and the notion of \( \alpha \)-Borelness also turns out to be the right generalization of Borelness when \( X = \pre{\kappa}{2} \) for \( \kappa \) an uncountable cardinal. 
Using the results on Lipschitz and Wadge reducibility from the previous section, we show in particular that under \( \ZFC \) the \( \kappa+1 \)-Borel hierarchy on \( \pre{\kappa}{2} \) does not collapse for every infinite cardinal \( \kappa \).%
\footnote{This is an illuminating example of how new and more refined arguments allow us to drop unnecessary assumptions on the cardinal \( \kappa \). 
In the literature this result is usually proved (with a straightforward generalization of the classical argument using universal sets and diagonalization) only for cardinals \( \kappa \) satisfying \( \kappa^{< \kappa} = \kappa \). 
In contrast, our proof allows us to also deal with cardinals satisfying \( \kappa^{< \kappa} > \kappa \) (including the singular ones), where the required universal sets do not exist at all.}

\subsection*{Section~\ref{sec:generalizedBorelfunctions}: Generalized Borel functions.} 
We study \( \kappa+1 \)-Borel functions and \( \bGamma \)-in-the-codes functions between generalized Cantor spaces. 
These are natural generalizations in two different directions of the notion of a Borel function between Polish spaces, which would correspond to the cases \( \kappa = \omega \) and \( \bGamma = \bSigma^1_1 \) (equivalently, \( \bGamma = \bS(\omega) \)), respectively.

\subsection*{Section~\ref{sec:otherspacesandBairecategory}: The generalized Baire space and Baire category.} 
We consider the generalized Baire space \( \pre{\kappa}{\kappa} \) and some of its relevant subspaces, including e.g.\ the group \( \Sym ( \kappa ) \) of permutations of \( \kappa \), and extend to them the analysis from Section~\ref{sec:topologies}. 
We determine under which conditions \( \pre{\kappa}{\kappa} \) and \( \pre{\kappa}{2} \) are homeomorphic, notably including the cases when \( \kappa \) is singular or when we are working in models of determinacy (such cases have not been considered in the literature so far). 
Finally, we prove some Baire category results which, as usual, can be recast in terms of forcing axioms.

\subsection*{Section~\ref{sec:spacesofcodes}: Standard Borel \( \kappa \)-spaces, \( \kappa \)-analytic quasi-orders, and spaces of codes.} 
Working in \( \ZF \), we introduce a very general notion of standard Borel \( \kappa \)-space (which in the classical setting \( \ZFC + \kappa^{< \kappa} = \kappa \) coincides with the one introduced in~\cite{Motto-Ros:2011qc}), and we consider \( \kappa \)-analytic relations on such spaces. 
We also introduce various nice spaces of codes for uncountable structures or non-separable spaces, including the space \( \Mod^\kappa_\LL \) of (codes for) \( \LL \)-structures of size \( \kappa \), the space \( \MMM_\kappa \) of (codes for) complete metric spaces of density character \( \kappa \), and the space \( \BBB_\kappa \) of (codes for) Banach spaces of density \( \kappa \). 
These spaces of codes are all standard Borel \( \kappa \)-spaces, and the corresponding isomorphism/embeddability relations on them are \( \kappa \)-analytic equivalence relations/quasi-orders.

\subsection*{Section~\ref{sec:infinitarylogics}: Infinitary logics and models.} 
We recall the model theoretic notions that are used in the sequel, including the infinitary logics \( \LL_{\kappa^+ \kappa} \), and we prove various generalizations to uncountable structures of the Lopez-Escobar theorem mentioned at the end of Section~\ref{subsubsec:eqrel} (some of these generalizations were independently obtained also in~\cite{Friedman:2011nx}). 
Although a full generalization can be obtained in \( \ZFC \) only assuming \( \kappa^{< \kappa} = \kappa \), a careful analysis of the proof leads to several intermediate results. 
We also introduce the bounded logic \( \LL^b_{\kappa^+ \kappa} \), a powerful enough fragment of \( \LL_{\kappa^+ \kappa} \) which avoids the pitfalls of na\"ive generalizations of the Lopez-Escobar theorem.
This is crucial for many of our main results.

\subsection*{Section~\ref{sec:Ksouslinsets}: \( \kappa \)-Souslin sets.} 
We briefly study the pointclass \( \bS ( \kappa ) \) of \( \kappa \)-Souslin subsets of Polish spaces and the collection of Souslin cardinals. 
These notions have been extensively studied under the assumption \( \ZF + \DC + \AD \): we both review those results from this beautiful and deep area of research which are relevant to our work, and prove analogous results (when possible) under the assumption \( \ZFC \). 
In particular, we compute the exact value of the cardinal \( \bdelta_{ \bS ( \kappa )} \) associated to the pointclass \( \bS ( \kappa ) \).
To the best of our knowledge this computation has been overlooked in the literature (even in the \( \AD \)-context), so this result might be of independent interest.

\subsection*{Section~\ref{sec:mainconstruction}: The main construction.} 
In this section we show how to construct certain combinatorial trees (i.e.\ acyclic connected graphs) of size \( \kappa \) starting from a descriptive set-theoretic tree \( T \) on \( 2 \times \kappa \): this is the main technical construction that gives our completeness and invariant universality results.

\subsection*{Section~\ref{sec:embeddabilitygraphs}: Completeness.} 
Using the construction from the previous section, we show that the relation \( \embeds^\kappa_\CT \) of embeddability between combinatorial trees of size \( \kappa \) is complete for the class of \( \kappa \)-Souslin quasi-orders on Polish or standard Borel spaces.

\subsection*{Section~\ref{sec:invariantlyuniversal}: Invariant universality.} 
We improve the result from Section~\ref{sec:embeddabilitygraphs} by showing that \( \embeds^\kappa_\CT \) is indeed invariantly universal for the same class of quasi-orders.

\subsection*{Section~\ref{sec:alternativeapproach}: An alternative approach.} 
We provide a modification of the main construction and results from the previous three sections: this variant yields a generalization of Theorem~\ref{th:mottorosfriedman2} that cannot be achieved using the construction and results from Sections~\ref{sec:mainconstruction}--\ref{sec:invariantlyuniversal}.

\subsection*{Section~\ref{sec:definablecardinality}: Definable cardinality and reducibility.} 
We discuss various forms of definable cardinality and reducibility that have already been considered in the literature, and translate our main results in corresponding statements involving these concepts.

\subsection*{Section~\ref{sec:applications}: Some applications.}
We present a selection of applications of the results from Section~\ref{sec:definablecardinality} in some of the most important and well-known frameworks, leading to those natural results which, as noticed at the beginning of this introduction, were the original motivation for the entire work and for the technical developments contained in it.

\subsection*{Section~\ref{sec:furthercompletenessresults}: Further completeness results.}
We show that \( \embeds^\kappa_\CT \) is complex also from the purely combinatorial point of view by embedding in it various partial orders and quasi-orders. 
Moreover, we provide some model theoretic variants of our completeness results concerning full homomorphisms between graphs and embeddings between lattices. 
Finally, we prove some results on the complexity of the relations of isometry and isometric (bi-)embeddability between complete metric spaces of uncountable density character, as well as on the complexity of the relations of linear isometry and linear isometric (bi-)embeddability between non-separable Banach spaces.

\section{Preliminaries and notation} \label{sec:preliminaries}
\subsection{Basic notions}
Our notation is standard as in most text in set theory, such as~\cite{Moschovakis:2009fk, Kechris:1995zt, Kanamori:2003fu,Jech:2003pd}.
By \( \ZF \) we mean the Zermelo-Fr\ae nkel set theory, which is formulated in the language of set theory \( \LST \), \index[symbols]{LST@\( \LST \)} i.e.~the first-order language with \( \in \) as the only non-logical symbol; \( \subseteq \) means subset and \( \subset \) means \emph{proper} subset, \( \pow ( X ) \) is the powerset of \( X \), the disjoint union and the symmetric difference of \( X \) and \( Y \) are denoted by \( X \uplus Y \)\index[symbols]{13@\( \uplus \)} and \( X \symdif Y \)\index[symbols]{13@\( \symdif \)} respectively, and so on.
For the reader's convenience we collect here all definitions and basic facts that are used later.

\subsubsection{Ordinals and cardinals}
Ordinals are denoted by lower case Greek letters, and \( \On \) is the class of all ordinals.
A \markdef{cardinal}\index[concepts]{cardinals and cardinalities} is an ordinal not in bijection with a smaller ordinal. 
The class of all infinite cardinals is \( \Cn \), and the letters \( \kappa , \lambda , \mu\) always denote an element of \( \Cn \). 
For \( \alpha \geq \omega \) we denote by \( \alpha ^+ \) the least cardinal larger than \( \alpha \).
The \markdef{cofinality} of \( \kappa \) is the smallest cardinal \( \lambda \) such that there is a cofinal \( j \colon \lambda \to \kappa \), i.e.~\( \FORALL{\alpha < \kappa} \EXISTS{\beta < \lambda} (\alpha \leq j ( \beta ) ) \). 
A cardinal \( \kappa \) is called \markdef{regular} is \( \cof ( \kappa) = \kappa \) and \markdef{singular} otherwise.

We use
\begin{equation}\label{eq:Hessenberg}\index[symbols]{13@\( \op{ \alpha }{ \beta } \) }\index[concepts]{coding of pairs/sequences of ordinals, \( \op{ \alpha }{ \beta } \) and \( \Code{ s} \)}
 \op{\cdot }{\cdot } \colon \On \times \On \to \On , 
\end{equation}
for the inverse of the enumerating function of the well-order \( \unlhd \) on \( \On \times \On \) defined by
\[
( \alpha , \beta ) \unlhd ( \gamma , \delta ) \IFF \max \set{ \alpha , \beta } < \max \set{ \gamma , \delta } \vee \bigl [ \max \set{ \alpha , \beta } = \max \set{ \gamma , \delta } \wedge ( \alpha , \beta ) \leqlex ( \gamma , \delta ) \bigr ] ,
\]
where \( \leqlex \) is the lexicographic ordering.
This pairing function maps \( \kappa \times \kappa \) onto \( \kappa \), for \( \kappa \in \Cn \).
Let
\begin{equation}\label{eq:codingfinitesequenceofordinals}
\Code{\cdot} \colon \pre{ < \omega }{\On} \to \On \index[symbols]{13@\( \Code{s} \)}
\end{equation}
be a bijection that maps \( \pre{ < \omega }{ \kappa } \) onto \( \kappa \), for \( \kappa \in \Cn \) --- see e.g.~\cite[Lemma 2.5, p.~159]{Kunen:1980fk}.

\subsubsection{Functions and sequences} \label{subsubsec:sequences}
Unless otherwise specified, all functions \( f \colon X \to Y \) are assumed to be total.
We write \( f \colon X \into Y \) and \( f \colon X \onto Y \) to mean that \( f \) is injective and that \( f \) is surjective, respectively. 
Similarly, the symbols \( X \into Y \) and \( X \onto Y \) means that there is an injection of \( X \) into \( Y \) and that \( Y \) is a surjective image of \( X \), respectively.
The pointwise image of \( A \subseteq X \) via \( f \colon X \to Y \) is \( \setof{ f ( a )}{ a \in A} \) and it is denoted either by \( f ( A ) \) or by \( \appl{f}{A} \)\index[symbols]{13@\( \appl{f}{A} \)} --- the first notation follows the tradition in analysis and descriptive set theory, the second is common in set theory especially when \( A \) is transitive.%
\footnote{The notation \( f [ A ] \) for the pointwise image is also common in the literature, but we eschew it as square brackets are already used for equivalence classes and for the body of trees --- see Section~\ref{subsubsec:dsttrees}.}
The set of all functions from \( X \) to \( Y \) is denoted by \( \pre{X}{Y} \), while 
\begin{equation} \label{eq:injections}
\pre{X}{( Y )} \coloneqq \setofLR{ f \in \pre{ X}{Y} }{f \text{ is injective}}. \index[symbols]{13@\( \pre{X}{( Y )} \)}
\end{equation}
For \( \kappa \) an infinite cardinal let
\[ 
\forcing{Fn} ( X , Y ; \kappa ) \coloneqq \setofLR{s}{s \colon u \to Y \AND u \subseteq X \AND \card{u} < \kappa } \index[symbols]{Fn@\( \forcing{Fn} ( X , Y ; \kappa ) \)} .
\]
Let \( \pre{ < \alpha }{ Y }\coloneqq \bigcup_{ \gamma < \alpha } \pre{ \gamma }{ Y } \), and let
\[
\forcing{Fn} ( \kappa , Y ; b ) \coloneqq \pre{ < \kappa }{ Y } 
\]
be the set of all sequences with values in \( Y \) and length \( < \kappa \).
For \( \alpha , \beta \in \On \) let 
\[ 
[ \alpha ]^ \beta \coloneqq \setofLR{u \subseteq \alpha }{ \ot ( u ) = \beta } . \index[symbols]{13@\( [ \alpha ]^ \beta \), \( [ \alpha ]^{ \leq \beta} \), \( [ \alpha ]^{ < \beta} \)}
\]
An element \( u \) of \( [ \alpha ]^ \beta \) is identified with its (increasing) enumerating function \( u \colon \beta \to \alpha \).
The sets \( \pre{ \leq \alpha }{ X } \), \( [ \alpha ]^ {< \beta} \), and \( [ \alpha ]^ {\leq\beta} \) are defined similarly.
The length of \( u \in \pre{ < \alpha }{X} \) is denoted by \( \lh u \).
The concatenation of \( u \in \pre{ < \omega }{ X } \) with \( v \in \pre{ \leq \omega }{X} \) is denoted by \( u \conc v \).
When dealing with sequences of length \( 1 \) we shall often blur the distinction between the element and the sequence, and write e.g.\ \( u \conc x \) rather than \( u \conc \seqLR{ x } \).
For \( x \in X \) and \( \alpha \in \On \), the sequence of length \( \alpha \) constantly equal to \( x \) is denoted by \( x^{( \alpha ) } \)\index[symbols]{13@\( x^{( \alpha ) } \)}.
If \( \emptyset \neq u \in \pre{ < \omega }{X } \) then 
\begin{equation}\label{eq:predecessorofu}
u^\star \coloneqq u \restriction \lh ( u ) - 1 \index[symbols]{13@\( u \mapsto u^\star \)}
\end{equation}
is the finite sequence obtained by deleting the last element from \( u \).

\subsection{Choice and determinacy} \label{subsec:choicedeterminacy}
Since several results in this paper concern models of set theory where the \markdef{Axiom of Choice} (\( \AC \))\index[concepts]{Axiom of Choice, and its weak forms!full choice \( \AC \)} may or may not hold, we shall state explicitly any assumption used beyond the axioms of \( \ZF \).
Among such assumptions are the \markdef{Axiom of Countable Choices}\index[concepts]{Axiom of Choice, and its weak forms!countable choice \( \AC_\omega \), \( \AC_\omega ( \R ) \)} (\( \AC_\omega \)), the \markdef{Axiom of Dependent Choices}\index[concepts]{Axiom of Choice, and its weak forms!dependent choice \( \DC \), \( \DC ( \R ) \)} (\( \DC \)),\index[symbols]{DC@\( \DC \), \( \DCR \)} their versions restricted to the reals \( \AC_\omega ( \R ) \) and \( \DCR \), and the so-called determinacy axioms, namely \( \AD \)\index[symbols]{AD@\( \AD \), \( \ADR \)} and its stronger version \( \ADR \). 
The \markdef{Axiom of Determinacy}\index[concepts]{determinacy axioms!\( \AD \), \( \ADR \)} \( \AD \) is the statement that for each \( A \subseteq \pre{\omega}{\omega} \) the zero-sum, perfect information two-players game \( G^ \omega _A \) is \markdef{determined}.
That is to say that one of the two players \( \I \) and \( \II \) has a winning strategy in \( G^ \omega _A \) where they take turns in playing \( n_0 , n_1 , \ldots \in \omega \) and \( \I \) wins if and only if \( \seqofLR{n_i }{ i \in \omega } \in A \).
The \markdef{Axiom of Real Determinacy} \( \ADR \) asserts that every \( G^\R_A \) is determined, where \( A \subseteq \pre{ \omega }{ \R} \) and \( \I \) and \( \II \) play \( x_0 , x_1 , \ldots \in \R \).
The principle \( \AD \) implies many regularity properties, such as: every set of reals is Lebesgue measurable, has the property of Baire, and has the \markdef{Perfect Set Property}\index[concepts]{Perfect Set Property \( \PSP \)} (\( \PSP \))\index[symbols]{PSP@\( \PSP \)}, i.e.~either it is countable, or it contains a homeomorphic copy of \( \pre{\omega}{2} \).

We also occasionally consider the \markdef{Axiom of \( \kappa \)-Choices}\index[concepts]{Axiom of Choice, and its weak forms!\( \kappa \)-choices \( \AC_ \kappa \), \( \AC_ \kappa ( \R ) \)} (\( \AC_\kappa \)),\index[symbols]{ACkappa@\( \AC_ \kappa \), \( \AC_ \kappa ( \R )\)} asserting that the product of \( \kappa \)-many nonempty sets is nonempty, its restriction to the reals \( \AC_\kappa ( \R ) \), the \markdef{Continuum Hypothesis} (\( \CH \)), and various \markdef{forcing axioms} like \( \MA_{\omega_1} \), \( \PFA \), and so on.

Although \( \AD \) and \( \ADR \) forbid any well-ordering of the reals, and therefore are incompatible%
\footnote{In fact \( \AD \) contradicts \( \AC_{\omega_1} ( \R ) \), since this choice principle implies the failure of \( \PSP \).} 
with \( \AC \), they are consistent with \( \DC \). 
In fact they imply weak forms of choice.
For example, \( \ADR \) implies the \markdef{Uniformization Property}\index[concepts]{Uniformization Property \( \Unif \)} (\( \Unif \)):\index[symbols]{Unif@\( \Unif \)} every function \( f \colon \R \to \R \) has a right inverse, i.e.~a function \( g \colon \R \to \R \) such that \( f ( g ( x ) ) = x \) for all \( x \in \R \); equivalently, for every set \( A \subseteq \R^2 \) there is a function \( f \) with domain \( \setof{ x \in \R }{ \EXISTS{y \in \R} \left ( ( x , y ) \in A \right )} \) such that \( ( x , f ( x ) ) \in A \).
Uniformization does not follow from \( \AD \): in \( \Ll ( \R ) \), the canonical inner model for \( \AD \), \( \bSigma^{2}_{1} \) is the largest collection of sets that can be uniformized, and by an unpublished theorem of Woodin's, \( \ADR \) and \( \AD + \Unif \) are equivalent over \( \ZF + \DC \).\label{pag:ADdoesnotimplyuniformization}

\subsection{Cardinality} \label{subsec:cardinality}
We write \( X \asymp Y \) to say that \( X \) and \( Y \) are in bijection, and \( \eq{X}_{\asymp} \) is the collection of all sets \( Y \) of minimal rank that are in bijection with \( X \), i.e.~it is the equivalence class of \( X \) under \( \asymp \) cut-down using Scott's trick.
Working in \( \ZF \), the \markdef{cardinality}\index[concepts]{cardinals and cardinalities} of a set \( X \) is defined to be 
\[
\card{X} \coloneqq 
\begin{cases}
\text{the unique cardinal \( \kappa \asymp X \)} &\text{if \( X \) is well-orderable,}
\\
\eq{X}_{\asymp} & \text{otherwise.}
\end{cases}\index[symbols]{14@\( \card{X} \)}
\]
Thus \( \AC \) implies that every cardinality is a cardinal.
Set \( \card{X} \leq \card{Y} \) if and only if \( X \into Y \).
By the Shr\"oder-Bernstein theorem, \( \card{X} = \card{Y} \) if and only if \( \card{X} \leq \card{Y} \leq \card{X} \).

When \( \lambda \) is a cardinal \( \kappa ^ \lambda \) denotes cardinal exponentiation, i.e.~the cardinality of \( \pre{ \lambda }{ \kappa } \), while \( \kappa ^{ < \lambda } \coloneqq \card{ \pre{ < \lambda }{ \kappa } } = \sup \setofLR{ \kappa ^ \nu }{ \nu < \lambda \AND \nu \text{ a cardinal}} \).
Whenever we write \( \kappa ^ \lambda \) (for \( \lambda \) and infinite cardinal) or \( \kappa^{< \lambda} \) (for \( \lambda\) an \emph{uncountable} cardinal) the Axiom of Choice \( \AC \) is tacitly assumed. 
Given a cardinal \( \kappa \) and a set \( X \), we let \( \pow_{ \kappa } ( X ) \coloneqq \setofLR{ Y \subseteq X}{\card{Y} < \kappa } \).

The next result is straightforward under choice --- the main reason to explicitly give a proof here, is to show that it is provable in \( \ZF \).

\begin{proposition}\label{prop:sequencesofsequences}
Suppose \( X \) has at least two elements and \( \lambda \) is an infinite regular cardinal. 
Then
\[
\card{ \pre{ < \lambda }{ ( \pre{ < \lambda }{X} ) } } = \card{ \lambda \times \pre{ < \lambda }{X} } = \card{\pre{ <\lambda }{X} } .
\]
\end{proposition}

\begin{proof}
Fix distinct \( x_0 , x_1 \in X \), and for \( s \in \pre{ < \lambda }{X} \) let \( s' \coloneqq x_0 ^{( \lh s )} \conc x_1 \conc s \).
Then 
\[
 \pre{ < \lambda }{ ( \pre{ < \lambda }{X} ) } \to \pre{ < \lambda }{X} , \quad \vec{s} = \seqofLR{ s_\beta }{ \beta < \lh \vec{s}} \mapsto s'_0 \conc s'_1 \conc \dotsc \conc s'_\beta \conc \dotsc ,
\]
is injective, and \( \pre{ < \lambda }{X} \into \lambda \times \pre{ < \lambda }{X} \ \into \pre{ < \lambda }{X} \times \pre{ < \lambda }{X} \into \pre{< \lambda }{ ( \pre{ < \lambda }{X} ) } \), so we are done by the Schr\"oder-Bernstein theorem.
\end{proof}

\subsection{Algebras of sets} \label{subsec:algebras}
Let \( X \) be an arbitrary set and let \( \alpha \geq \omega \). 
A collection \( \mathcal{A} \subseteq \pow ( X ) \) is called \markdef{\( \alpha \)-algebra (on \( X \))}\index[concepts]{alphaAlgebra@\( \alpha \)-algebra} if it is closed under complements and well-ordered unions of length \( {<} \alpha \), that is \( \bigcup_{ \nu < \beta } A_\nu \in \mathcal{A} \) for every \( \beta < \alpha \) and every sequence \( \seqof{ A_\nu }{ \nu < \beta } \) of sets in \( \mathcal{A} \). 
An \( \omega \)-algebra is usually called an algebra, and an \( \alpha + 1 \)-algebra is the same as an \( \alpha ^+ \)-algebra.
In particular, an \( \omega + 1 \)-algebra (i.e. an \( \omega _1 \)-algebra) is what is usually called a \( \sigma \)-algebra.

Given a family \( \mathcal{G} \subseteq \pow ( X ) \), the \markdef{\( \alpha \)-algebra (on \( X \)) generated by \( \mathcal{G} \)} is the smallest \( \alpha \)-algebra \( \mathcal{A} \) on \( X \) such that \( \mathcal{G} \subseteq \mathcal{A} \), and is denoted by 
\[ 
\Alg ( \mathcal{G}, \alpha ) . \index[symbols]{Alg@\( \Alg ( \mathcal{G}, \alpha ) \)}
\] 
Thus the family of the Borel subsets of a topological space \( ( X , \tau ) \) is 
\[
\bB ( X ) = \Alg ( \tau , \omega + 1 ) = \Alg ( \tau , \omega_1 ) .
\]
The algebra \( \Alg ( \mathcal{G}, \alpha ) \) can be inductively generated as follows: let
\begin{align*}
\bSigma_1 ( \mathcal{G}, \alpha ) &\coloneqq \mathcal{G} 
\\
 \bPi_1 ( \mathcal{G}, \alpha ) &\coloneqq \setofLR{X \setminus A}{ A \in \mathcal{G}}
\end{align*}
and for \( \gamma > 1 \) 
\begin{align*}
\bSigma_\gamma ( \mathcal{G} , \alpha ) &\coloneqq \setofLR{\textstyle\bigcup_{ \nu < \beta} A_\nu }{ \beta < \alpha \AND \FORALL{ \nu < \beta } ( A_ \nu \in \textstyle\bigcup_{ \xi < \gamma } \bPi_{\xi} ( \mathcal{G} , \alpha ) ) } 
\\
\bPi_\gamma ( \mathcal{G} , \alpha ) &\coloneqq \setofLR{ X \setminus A }{A \in \bSigma_\gamma ( \mathcal{G} , \alpha ) } .
\end{align*}

The next result is straightforward.

\begin{lemma}\label{lem:inclusionhierarchy}
If \( 1 \leq \gamma < \delta \) and \( \delta \geq 3 \) we have that 
\[
\bSigma_{ \gamma } ( \mathcal{G} , \alpha ) \cup \bPi_{ \gamma } ( \mathcal{G} , \alpha ) \subseteq \bSigma_{ \delta } ( \mathcal{G} , \alpha ) \cap \bPi_{ \delta } ( \mathcal{G} , \alpha ) 
\]
so that if \( \mathcal{G} \subseteq \bSigma_2 ( \mathcal{G} , \alpha ) \) the inclusion holds for every \( 1 \leq \gamma < \delta \).
Moreover
\[
\Alg ( \mathcal{G} , \alpha ) = \bigcup_{ \gamma \in \On } \bSigma_{ 1 + \gamma } ( \mathcal{G} , \alpha ) = \bigcup_{\gamma \in \On} \bPi_{ 1 + \gamma } ( \mathcal{G} , \alpha ) .
\]
\end{lemma}

Notice that \( \On \) in the equation above can be replaced by any cardinal \( \lambda \) with \( \cof ( \lambda ) \geq \alpha \).

\begin{lemma}\label{lem:algebra}
Suppose \( X \) is an arbitrary set, \( \alpha \) is an infinite ordinal, \(\lambda \geq \alpha \) is a regular cardinal, and \( \mathcal{G} \subseteq \pow ( X ) \) contains at least two elements.
Then
\begin{enumerate-(a)}
\item\label{lem:algebra-a}
\( \pre{ < \lambda}{\mathcal{G}} \onto \Alg ( \mathcal{G} , \alpha ) \), and
\item\label{lem:algebra-b}
if \( \lambda = \nu^+ \) then \( \pre{ \nu }{\mathcal{G}} \onto \Alg ( \mathcal{G} , \alpha ) \).
\end{enumerate-(a)} 
In particular, if \( \AC \) holds and \( \card{\mathcal{G}} \leq \mu \) then \( \card{\Alg ( \mathcal{G} , \lambda )} \leq \mu^{ < \lambda} \) , which equals to \( \mu^{\card{ \nu}} \) if~\ref{lem:algebra-b} holds. 
\end{lemma}

\begin{proof}
\ref{lem:algebra-a} 
It is enough to construct surjections 
\[
p_ \gamma \colon \pre{ < \lambda }{\mathcal{G}} \onto \bSigma^{}_{ 1 + \gamma } ( \mathcal{G} , \alpha )
\]
for each \( \gamma < \lambda \), so that 
\[
 \lambda \times \pre{ < \lambda }{\mathcal{G}} \to \bigcup_{ \gamma < \lambda } \bSigma^{}_{ 1 + \gamma } ( \mathcal{G} , \alpha ) , \qquad ( \gamma , s ) \mapsto p_ \gamma ( s )
\]
is a surjection, and hence the result follows from Proposition~\ref{prop:sequencesofsequences}.
The construction of \( p_0 \colon \pre{ < \lambda }{\mathcal{G}} \to \bSigma^{}_{ 1 } ( \mathcal{G} , \alpha ) = \mathcal{G} \) is immediate.
Suppose \( p_ \nu \) has been defined for all \( \nu < \gamma \), and let 
\[ 
q_ \gamma \colon \gamma \times \pre{ < \lambda }{\mathcal{G}} \onto \bigcup_{ \nu < \gamma } \bSigma^{}_{ 1 + \nu } ( \mathcal{G} , \alpha ) , \qquad ( \nu , s ) \mapsto p_\nu ( s ) .
\]
Since \( \gamma \times \pre{ < \lambda }{\mathcal{G}} \subseteq \lambda \times \pre{ < \lambda }{\mathcal{G}} \), by repeated applications of Proposition~\ref{prop:sequencesofsequences} one gets \( \pre{ < \lambda }{\mathcal{G}} \onto \bigcup_{ \nu < \gamma } \bSigma^{}_{ 1 + \nu } ( \mathcal{G} , \alpha ) \), and hence \( \pre{ < \lambda }{ ( \pre{ < \lambda }{\mathcal{G}} ) } \onto \bSigma^{}_{ 1+\gamma } ( \mathcal{G} , \alpha ) \), and therefore \( \pre{ < \lambda }{\mathcal{G}} \onto \bSigma^{}_{ 1+\gamma } ( \mathcal{G} , \alpha ) \).

\smallskip

\ref{lem:algebra-b}
As in part~\ref{lem:algebra-a} one constructs surjections \( p_ \gamma \colon \pre{ \nu }{\mathcal{G}} \onto \bSigma^{}_{ 1 + \gamma } ( \mathcal{G} , \alpha ) \) for all \( \gamma < \nu^+ \), using the fact that \( \pre{ \nu }{2} \onto \nu^+ \).
\end{proof}

The \( \sigma \)-algebra of the \markdef{Borel subsets} of a topological space \( ( X , \tau ) \) is \( \bB ( X , \tau ) = \Alg ( \tau , \omega + 1 ) \).
When the topology \( \tau \) is clear from the context, it is customary write \( \bB ( X ) \), \( \bSigma^{0}_{ \alpha } ( X ) \) and \( \bPi^{0}_{ \alpha } ( X ) \) instead of \( \bB ( X , \tau ) \), \( \bSigma_{\alpha} ( \tau , \omega + 1 ) \) and \( \bPi_{\alpha} ( \tau , \omega + 1 ) \).\index[symbols]{Sigma0alpha@\( \bSigma^{0}_{ \alpha } \)}\index[symbols]{Pi0alpha@\( \bPi^{0}_{ \alpha } \)} 
Thus \( \bSigma^{0}_{1} ( X ) \) is the family \( \tau \) of all open sets of \( X \), \( \bPi^{0}_{1} ( X ) \) is the family of all closed subsets of \( X \), \( \bSigma^{0}_{2} ( X ) \) is the collection of all \( \Fsigma \) subsets of \( X \), \( \bPi^{0}_{2} \) is the collection of all \( \Gdelta \) subsets of \( X \), and so on.
Assuming \( \omega _1 \) is regular%
\footnote{It is consistent with \( \ZF \) that \( \R \) is countable union of countable sets, and hence that \emph{every} subset of \( \R \) is Borel.} 
(a fact that follows from \( \AComega ( \R ) \)), then \( \bB ( X ) = \bigcup_{ 1 \leq \alpha < \omega _1 } \bSigma^{0}_{ \alpha } ( X ) = \bigcup_{ 1 \leq \alpha < \omega _1 } \bPi^{0}_{ \alpha } ( X ) \).
If \( X \) is metrizable then every closed set is \( \Gdelta \), so \( \bSigma^{0}_{ \alpha } ( X ) \cup \bPi^{0}_{ \alpha } ( X ) \subseteq \bSigma^{0}_{ \beta } ( X ) \cap \bPi^{0}_{ \beta } ( X ) \) for every \( 1 \leq \alpha < \beta \).
From Lemma~\ref{lem:algebra} we obtain at once the following result.

\begin{corollary}\label{cor:lem:algebra}
Assume \( \AComega ( \R ) \).
If \( \tau \) is a topology on \( X \neq \emptyset \), then \( \pre{ \omega }{ \tau } \onto \bB ( X ) \).
\end{corollary}

\subsection{Descriptive set theory}
\subsubsection{Polish spaces}
A topological space is \markdef{perfect} if it has no isolated points; it is \markdef{zero-dimensional} if it has a basis of clopen sets.
A \markdef{Polish space}\index[concepts]{space!perfect, zero-dimensional}\index[concepts]{Polish space} is a separable, completely metrizable topological space; a \markdef{Polish metric space} is a complete \emph{metric} space that is separable.
If the spaces \( X_n \) (\( n \in \omega \)) are Polish metric, then so is \( \prod_n X_n \) with the product topology.
Since any countable set with the discrete topology is Polish, then the \markdef{Cantor space} \( \pre{\omega}{2} \) and the \markdef{Baire space} \( \pre{\omega}{\omega} \) are Polish: a countable dense set is given by the sequences that are eventually constant, and a complete metric for them is given by \( d ( x , y ) = 2^{-n} \) if \( n \) is least such that \( x ( n ) \neq y ( n ) \) and \( d ( x , y ) = 0 \) if there is no such \( n \).
This metric is in fact an ultrametric, so the balls are clopen and hence these spaces are zero-dimensional.
Recall that a metric \( d \) on a space \( X \) is called \markdef{ultrametric}\index[concepts]{space!ultrametric} if it satisfies the following strengthening of the triangular inequality: \( d ( x , y ) \leq \max \set{d ( x , z ) , d ( z , y ) } \) for all \( x , y , z \in X \).
The space \( \pre{\omega}{2} \) is the unique (up to homeomorphism) nonempty compact zero-dimensional perfect Polish space, and \( \pre{\omega}{\omega} \) is the unique (up to homeomorphism) nonempty zero-dimensional perfect Polish space in which all compact subsets have empty interior~\cite[Theorems 7.4 and 7.7]{Kechris:1995zt}.
Any two uncountable Polish spaces are Borel isomorphic~\cite[Theorem 15.6]{Kechris:1995zt} and since many questions in descriptive set theory are invariant under Borel isomorphism, it is customary to use \( \R \) to denote any uncountable Polish space.

\subsubsection{Pointclasses} \label{subsubsec:pointclasses}
A \markdef{general pointclass} \( \Gamma \) is an operation (i.e. a class-function) assigning to every nonempty topological space \( X \) a nonempty family \( \Gamma ( X ) \subseteq \pow ( X ) \). 
The \markdef{dual} of \( \Gamma \)\index[concepts]{pointclass!dual \( \check{\Gamma} \)} is the general pointclass defined by \( \check{\Gamma} ( X ) \coloneqq \setofLR{ X \setminus A }{ A \in \Gamma ( X ) } \).\index[symbols]{G0amma@\( \check{\Gamma} \)} 
The \markdef{ambiguous general pointclass} associated to \( \Gamma \) (or to \( \check{ \Gamma} \)) is the general pointclass \( \Delta_{\Gamma} \) defined by \( \Delta_{\Gamma} ( X ) \coloneqq \Gamma ( X ) \cap \check{\Gamma} ( X ) \).\index[symbols]{D1eltaGamma@\( \Delta_{\Gamma} \)}
A general pointclass \( \Gamma \) is \markdef{hereditary}\index[concepts]{pointclass!hereditary} if \( \Gamma ( Y ) = \setofLR{A \cap Y}{A \in \Gamma ( X )} \), for every pair of nonempty topological spaces \( Y \subseteq X \). 

A general pointclass \( \bGamma \) is said to be \markdef{boldface} if it is closed under continuous preimages, that is to say: if \( f \colon X \to Y \) in continuous and \( B \in \bGamma ( Y ) \) then \( f^{-1} ( B ) \in \bGamma ( X ) \). 
When the topological space is clear from the context, we write \( A\in \bGamma \) rather than \( A \in \bGamma ( X ) \).
General boldface pointclasses are usually denoted by Greek letters such as \( \bGamma \) and \( \bLambda \), variously decorated, and are typeset boldface, whence the name.
Examples of general boldface pointclasses are: the collection of all open sets \( \bSigma^{0}_{1} \), its dual pointclass \( \bPi^{0}_{1} \) i.e.\ the collection of all closed sets, its ambiguous pointclass \( \bDelta^{0}_{1} \) i.e.\ the collection of all clopen sets, the Borel pointclasses \( \bSigma^0_ \alpha , \bPi^0_ \alpha \) and the Borel sets \( \bB \). 

A (\markdef{boldface}) \markdef{pointclass}\index[concepts]{pointclass!boldface} is a general (boldface) pointclass restricted to \emph{Polish} spaces.
In other words, a pointclass is an operation assigning to every nonempty Polish space \( X \) a nonempty family \( \bGamma ( X ) \subseteq \pow ( X ) \), and it is boldface if \( f^{-1} ( B ) \in \bGamma ( X ) \) for every \( B \in \bGamma ( Y ) \) and every continuous function \( f \colon X \to Y \) between Polish spaces. 
If \( \bGamma ( \pre{\omega}{\omega} ) \neq \pow ( \pre{\omega}{\omega} ) \) (equivalently: \( \bGamma ( X ) \neq \pow ( X ) \) for any \emph{uncountable} Polish space \( X \)), then \( \bGamma \) is called a \markdef{proper} pointclass. 
If \( \setLR{ \emptyset , X } \subset \bGamma ( X ) \) for some space \( X \) (equivalently: \( \bGamma ( \pre{\omega}{\omega} ) \supseteq \bDelta^{0}_{1} ( \pre{\omega}{\omega} ) \)) then \( \bGamma \) is \markdef{nontrivial}. 
The boldface pointclass \( \bGamma \) is called \markdef{nonselfdual} if \( \bGamma ( \pre{\omega}{\omega} ) \neq \check{\bGamma} ( \pre{\omega}{\omega} ) \) (equivalently: \( \bGamma ( X ) \neq \check{\bGamma} ( X ) \) for any \emph{uncountable} Polish space \( X \)), and \markdef{selfdual} otherwise. 
The boldface pointclass \( \bGamma \) admits a \markdef{universal set} if there is \( \mathcal{U} \in \bGamma ( \pre{\omega}{\omega} \times \pre{\omega}{\omega} ) \) such that \( \bGamma ( \pre{\omega}{\omega} ) = \setofLR{ \mathcal{U}^{( y )} }{ y \in \pre{\omega}{\omega} } \), where \( \mathcal{U}^{( y )} \coloneqq \setofLR{ x \in \pre{\omega}{\omega}}{( x , y ) \in \mathcal{U}} \) is the \markdef{horizontal section} of \( \mathcal{U} \). 
If \( \bGamma \) admits a universal set, then it is nonselfdual: under \( \AD \) these two properties become equivalent, that is: \( \bGamma \) admits a universal set if and only if it is nonselfdual.
These definitions extend to the case of a \emph{general} boldface pointclass: \( \bGamma \) is proper (respectively, nontrivial, nonselfdual, selfdual) if its restriction to the class of Polish spaces is a proper (respectively, nontrivial, nonselfdual, selfdual) boldface pointclass. 

Besides the Borel sets \( \bB \), and the Borel pointclasses \( \bSigma^0_ \alpha , \bPi^0_ \alpha \), examples of boldface pointclasses are the projective pointclasses \( \bSigma^1_n \) (for \( n \in \omega \)), and the pointclasses associated to third-order arithmetic, such as \( \bSigma^2_1 \).
As usual, the dual of \( \bSigma^{i}_{\alpha} \) is \( \bPi^{i}_{\alpha} \) and its ambiguous part is \( \bDelta^i_\alpha \) (for \( i = 0, 1 , 2 \) and \( 0 \neq \alpha < \omega_1 \)). 
The pointclasses \( \bSigma^i_\alpha \) and \( \bPi^i_\alpha \) are hereditary and nonselfdual, the pointclass \( \bB \) is hereditary and selfdual, while the \( \bDelta^i_\alpha \)'s are selfdual and not hereditary.

The projective pointclasses are \emph{not} general boldface pointclasses, as they are not closed under continuous preimages.
To see this consider the pointclass \( \bSigma^1_1 \) of analytic sets\index[concepts]{analytic sets \( \bSigma^1_1 \)}, i.e.~for \( X \) Polish we set
\[
 \bSigma^1_1 ( X ) = \setof{ f ( \pre{\omega}{\omega} )}{ f \colon \pre{\omega}{\omega} \to X \text{ is continuous}}\cup \set{ \emptyset } .
\]
If \( Y \) is a proper \( \bPi^1_1 \) subset of \( X \), then \( Y \) is the preimage of \( X \in \bSigma^1_1 ( X ) \) under the inclusion map \( Y \hookrightarrow X \), but \( Y \) is not a continuous image of the Baire space, i.e.~\( Y \notin \bSigma^1_1(Y) \).
Closing \( \bSigma^1_1 \) under continuous preimages yields a general pointclass \( \bGamma \), but it is not clear whether \( \bGamma \) contains all Borel sets (or even the closed ones), or whether it is closed under continuous images.
Moreover we have no use for such \( \bGamma \).
This is why we stick to the usual definition of \( \bSigma^{1}_{1} \) rather than adopting the one of \( \bGamma \).

\subsubsection{The prewellordering and scale properties} \label{subsubsec:pwoscaleproperty} 
Let \( A \) be a set. 
A \markdef{(regular) norm}\index[concepts]{norm on a set} on \( A \) is a map \( \rho \colon A \to \On \) which is surjective onto some \( \alpha \in \On \). 
The ordinal \( \alpha \) is called \markdef{length} of \( \rho \). 
To each norm \( \rho \colon A \to \alpha \) we can canonically associate a \markdef{prewellordering} (i.e.~a reflexive, transitive, connected, and well-founded relation) \( \preceq_\rho \) of \( A \) by setting \( x \preceq_\rho y \iff \rho ( x ) \leq \rho ( y ) \) for all \( x , y \in A \). 
Conversely, to every prewellordering \( \preceq \) of \( A \) we can canonically associate a (unique) regular norm \( \rho \colon A \onto \alpha \) (for some \( \alpha \in \On \)) such that \( {\preceq} = {\preceq_\rho} \): in this case we say that \( \preceq \) has length \( \alpha \).

Given a boldface pointclass \( \bGamma \) and a Polish space \( X \), a norm \( \rho \) on \( A \in \bGamma ( X ) \) is a \markdef{\( \bGamma \)-norm} if there are two binary relations \( \leq_\rho^{\bGamma} \in \bGamma ( X \times X ) \) and \( \leq_\rho^{\check{\bGamma}} \in \check{\bGamma} ( X \times X ) \) on \( X \) such that for all \( y \in A \)
\[ 
\FORALL{ x \in X} \left [ ( x \in A \wedge \rho ( x ) \leq \rho ( y ) ) \IFF x \leq_\rho^{\bGamma} y \IFF x \leq_\rho^{\check{\bGamma}} y \right ].
 \] 
In other words: the initial segments of the prewellordering \( \preceq_\rho \) are uniformly \( \bDelta_{ \bGamma} \).
We say that the pointclass \( \bGamma \) is \markdef{normed}\index[concepts]{pointclass!normed} or has the \markdef{prewellordering property}\index[concepts]{prewellordering property|see{normed pointclass}} if every \( A \in \bGamma ( X ) \) admits a \( \bGamma \)-norm, for all Polish spaces \( X \).
A \( \bGamma \)-norm on a set \( A \in \bDelta_{ \bGamma} \) is automatically a \( \dual{\bGamma} \)-norm as well.
On the other hand the property of being normed does not pass to the dual: for example \( \bPi^{1}_{1} \) is normed, but \( \bSigma^{1}_{1} \) is not~\cite[Section~4.B]{Moschovakis:2009fk}.

Let \( X \) be a topological space. 
A \markdef{scale}\index[concepts]{scale} on \( A \subseteq X \) is a sequence \( \seqofLR{ \rho_n }{ n \in \omega } \) of (regular) norms on \( A \) such that for every sequence \( \seqofLR{ x_n }{ n \in \omega } \) of points from \( A \), if the \( x_n \)'s converge to \( x \in X \) and for every \( n \in \omega \) the sequence \( \seqofLR{ \rho_n ( x _i ) }{ i \in \omega } \) is eventually equal to some \( \lambda_n \in \On \), then \( x \in A \) and \( \rho_n ( x ) \leq \lambda_n \) for all \( n \in \omega \). 
When \( X = \pre{\omega}{\omega} \), the existence of a scale \( \seqofLR{ \rho _n }{ n \in \omega } \) on \( A \) is equivalent to the assertion that \( A \) is the projection of a tree on \( \omega \times \kappa \), where \( \kappa \coloneqq \sup_n \ran ( \rho _n ) \) (see Section~\ref{subsubsec:graph} for the relevant definitions). 

Given a boldface pointclass \( \bGamma \) closed under countable intersections and countable unions, a scale \( \seqofLR{ \rho_n }{ n \in \omega } \) on a set \( A \in \bGamma ( X ) \) is called \markdef{\( \bGamma \)-scale} if all the \( \rho_n \)'s are \( \bGamma \)-norms.%
\footnote{The definition of \( \bGamma \)-scale for an arbitrary boldface pointclass \( \bGamma \) requires that the norms be \emph{uniformly} in \( \bGamma \) --- see~\cite[p.~173]{Moschovakis:2009fk}.} 
The pointclass \( \bGamma \) has the \markdef{scale property}\index[concepts]{pointclass!scaled} if every \( A \in \bGamma ( X ) \) admits a \( \bGamma \)-scale for all Polish spaces \( X \).

\subsection{Trees and reductions} \label{subsec:treesandreductions}
\subsubsection{Graphs and trees} \label{subsubsec:graph}
A \markdef{graph} \( G = ( V , E ) \) consists of a nonempty set \( V \) of vertices and a set \( E \subseteq [V]^2 \) of edges.
Whenever \( V \) is clear from the context, the graph is identified with \( E \), and we write \( v_0 \mathrel{G} v_1 \) or \( v_0 \mathrel{E} v_1 \) instead of \( \setLR{ v_0, v_1 } \in E \).
A \markdef{combinatorial tree}\index[concepts]{combinatorial tree} is a connected acyclic graph. 
If a specific vertex is singled out, the resulting object is a \markdef{rooted combinatorial tree}\index[concepts]{combinatorial tree!rooted} and the chosen vertex is called \markdef{root}.
The size of a combinatorial tree is the cardinality of the set \( V \) of its vertices. 
As the nature of the elements of \( V \) is irrelevant, a combinatorial tree of size \( \kappa \) can be construed as a set of edges \( E \subseteq [ \kappa ]^2 \), and hence it can be identified with an element of \( \pre{ \kappa \times \kappa }{ 2 } \) (via its characteristic function).
Thus the space of all combinatorial trees of size \( \kappa \) is 
\begin{equation}
\CT_ \kappa \coloneqq \setofLR{f \in \pre{ \kappa \times \kappa }{ 2 } }{ f \text{ satisfies~\eqref{eq:CTdef-a}--\eqref{eq:CTdef-c}} } \label{pag:CT_kappa}
\end{equation}\index[symbols]{CT@\( \CT_\kappa \)}
with
\begin{subequations} \label{eq:CTdef}
\begin{align}
 & \FORALL{ \alpha , \beta \in \kappa } \left ( f ( \alpha , \alpha ) = 0 \wedge \left ( f ( \alpha , \beta ) = 1 \implies f ( \beta , \alpha ) = 1 \right ) \right ) \label{eq:CTdef-a} 
 \\
 & \FORALL{\alpha,\beta \in \kappa} \big [ f ( \alpha , \beta ) = 1 \OR \EXISTS{s \in \pre{< \omega}{\kappa} \setminus \{ \emptyset \}} \FORALL{i < \lh ( s ) - 1 }( f ( s ( i ) , s ( i + 1 ) ) = 1 \AND \label{eq:CTdef-b} 
 \\
& \hspace{7cm} f ( \alpha , s ( 0 ) ) = 1 \AND f ( s ( \lh ( s ) - 1 ) , \beta ) = 1 ) \big] \notag 
 \\
& \FORALL{s \in \pre{<\omega}{\kappa}} ( \lh ( s ) \geq 3 \AND \FORALL{ i < \lh ( s ) - 1 } ( f ( s ( i ) , s ( i + 1 ) ) = 1) \IMPLIES \label{eq:CTdef-c}
\\
& \hspace{7cm} f ( s ( \lh ( s ) - 1 ) , s ( 0 ) ) = 0). \notag
\end{align}
\end{subequations}
The intuition behind these formul\ae{} is that:~\eqref{eq:CTdef-a} says that \( f \) codes a graph \( G_f \),~\eqref{eq:CTdef-b} says that \( G_f \) is connected, while~\eqref{eq:CTdef-c} says that \( G_f \) is acyclic.
It is easy to check that \( \CT_ \omega \) is a \( \Gdelta \) subset of \( \pre{ \omega \times \omega }{ 2 } \) (which is homeomorphic to the Cantor space \( \pre{\omega}{2} \)), and hence \( \CT_ \omega \) is a Polish space.

\subsubsection{Descriptive set-theoretic trees} \label{subsubsec:dsttrees}
The word \emph{tree} may also refer to a different concept in descriptive set theory: a \markdef{tree on a set} \( X \) is a nonempty subset of \( \pre{ < \omega }{X} \) closed under initial segments (ordered by the prefix relation).
Such an object is often called a \markdef{descriptive set-theoretic tree},\index[concepts]{descriptive set-theoretic tree} and
\begin{equation}\label{eq:Tr}
\Tr ( X ) \index[symbols]{Tr@\( \Tr \)} 
\end{equation}
is the set of all descriptive set-theoretic trees on \( X \).
Any descriptive set-theoretic tree can be seen as a rooted combinatorial tree with root \( \emptyset \), and conversely. 
Elements of a (descriptive set-theoretic) tree \( T \) on \( X \) are called \markdef{nodes}. 
We say that \( T \) is \markdef{pruned} if every node has a proper extension, and that \( T \) is \markdef{\(< \kappa \)-branching} (for \( \kappa \) a cardinal) if every node has \( < \kappa \)-many immediate successors, that is \( \card{\setof{x \in X }{s \conc x \in T} } < \kappa \) for every \( s \in T \). 
Sometimes, \( < (\kappa +1) \)-branching trees are simply called \( \kappa \)-branching. 
The \markdef{body of} \( T \) is the set of all infinite branches of \( T \), that is the set
\[ 
\body{T} \coloneqq \setofLR{f \in \pre{ \omega }{X} }{ \FORALL{ n} ( f \restriction n \in T ) }.\index[symbols]{17@\( \body{T} \)}
\]
If \( T \) is a tree on \( X \times Y \) then the nodes are construed as pairs of sequences \( ( \seqLR{x_0 , \dots , x_n } , \seqLR{y_0 , \dots , y_n } ) \) rather than sequences of pairs \( \seqLR{ ( x_0 , y_0 ) , \dots , ( x_n , y_n) } \), and similarly the elements of \( \body{T} \) are construed as pairs \( ( f , g ) \in \pre{ \omega }{X} \times \pre{ \omega }{Y} \) such that \( ( f \restriction n , g \restriction n ) \in T \) for all \( n \). 
The \markdef{projection} (on the first coordinate) of a subset \( A \subseteq X \times Y \) of a cartesian product is 
\begin{equation}\label{eq:projectionCh2}\index[symbols]{p@\( \PROJ \)}
\PROJ A = \setof{ x \in X}{ \EXISTS{y \in Y} ( x , y ) \in A } 
\end{equation}
The projection of a tree \( T \) on \( X \times Y \) is the set
\[ 
\PROJ \body{T} \coloneqq \setofLR{f \in \pre{ \omega }{X} }{ \EXISTS{g \in \pre{ \omega }{Y} } ( f , g ) \in \body{T} }. 
\]

\subsubsection{Quasi-orders and equivalence relations} \label{subsubsec:qoander}
A binary relation \( R \) on a set \( X \) is called \markdef{quasi-order} or \markdef{preorder} if it is reflexive and transitive, and a symmetric quasi-order is called \markdef{equivalence relation}.
The set \( X \) is called \markdef{domain of \( R \)} and it is denoted by \( \dom ( R ) \).
If \( Y \subseteq X = \dom ( R ) \), we denote by \( R \restriction Y \) the restriction of \( R \) to \( Y \), i.e.~the quasi-order \( R \cap Y^2 \). 
If \( R \) is a quasi-order, then \( R^{-1} \coloneqq \setofLR{ ( y , x )}{ ( x , y ) \in R } \) is also a quasi-order, and \( E_R \coloneqq R \cap R^{-1} \) is the equivalence relation induced by \( R \).
A quasi-order \( R \) canonically induces a partial order on the quotient space \( X / E_R \), which is called \markdef{quotient order of \( R \)}.

If \( E \) is an equivalence relation on \( X \), then \( \eq{x}_E = \eq{x} \) is the \markdef{equivalence class of \( x \in X \)} and \( X / E \) is the \markdef{quotient space}.
A set \( A \subseteq X = \dom ( E ) \) is \markdef{invariant under \( E \)} if \( y \in A \) whenever \( x \mathrel{E} y \) for some \( x \in A \). 
The \markdef{\( E \)-saturation of} \( A \subseteq X \), in symbols \( \eq{A}_E \) or even just \( \eq{A} \), is the smallest invariant set containing \( A \), that is 
\[
\eq{A}_E = \eq{A} \coloneqq \bigcup_{x \in A} \eq{x} .
\]

Given a (general) boldface pointclass \( \bGamma \), we say that a quasi-order (or, more generally: a binary relation) \( R \) on \( X \) is in \( \bGamma \) if \( R \in \bGamma ( X \times X ) \). 
Notice that if \( \bGamma \) is closed under finite intersections, this implies that the associated equivalence relation \( E_R \) is in \( \bGamma \) as well.

\subsubsection{Reducibility}\label{subsubsec:reducibility}
Given two quasi-orders \( R , S \) on the sets \( X , Y \), we say that \( R \) \markdef{reduces} to \( S \) (in symbols \( R \leq S \)) if and only if there is a function \( f \colon X \to Y \) which reduces \( R \) to \( S \), i.e.~such that for all \( x , x' \in X \)
\[ 
x \mathrel{R} x' \iff f ( x ) \mathrel{S} f ( x' ) . 
\]
If \( R \leq S \) and \( S \leq R \) we say that \( R \) and \( S \) are \markdef{bi-reducible} (in symbols \( R \sim S \)).
Notice that if \( E \) is an equivalence relation and \( R \leq E \), then \( R \) is an equivalence relation as well.
If \( \mathcal{C} \) is a class of quasi-orders, we say that the quasi-order \( R \) is \markdef{complete for \( \mathcal{C} \)} if \( S \leq R \) for every \( S \in \mathcal{C} \). 
By limiting \( f \) to range in a given collection of functions we obtain a restricted form of reducibility.
For example if \( f \) is Borel or \( f \in \Ll ( \R ) \), then the notions of \markdef{Borel reducibility} \( \leq_{\bB} \) and \markdef{\( \Ll ( \R ) \)-reducibility} \( \leq_{\Ll(\R)} \) are obtained, respectively.
If \( {\leq_*} \) is a restricted form of reducibility, we say that \( R \) is \markdef{\( \leq_* \)-complete for \( \mathcal{C} \)} if \( S \leq_* R \) for every \( S \in \mathcal{C} \).
In particular \( R \) is \markdef{Borel-complete} for \( \mathcal{C} \) if for every \( S \in \mathcal{C} \) there is a Borel function that witnesses \( S \leq R \), i.e.\ \( S \leq_{\bB} R \).

\begin{remark}\label{rmk:completeness}
Our definition of ``\( R \) is (\( \leq_* \)-)complete for \( \mathcal{C} \)'' does not require that \( R \in \mathcal{C} \). 
In fact, in most applications we will have that \( \mathcal{C} \) is a collection of quasi-orders on some Polish space, while \( R \) will be the embeddability relation on some elementary class of \emph{uncountable} models. 
Notice however that \( R \) is a \( \kappa \)-analytic quasi-order in the sense of generalized descriptive set theory (Definition~\ref{def:kappa-analytic}) irrespective of the complexity of the quasi-orders in \( \mathcal{C} \).
As \( \kappa \)-analyticity is the natural counterpart of \( \bSigma^{1}_{1} \) for Polish spaces, \( R \) is arguably not more complex than the elements of \( \mathcal{C} \), as long as \( \kappa \) is related to the complexity of these elements.
This is why we use \emph{completeness} rather than \emph{hardness} in our results.
\end{remark}

The notion of reducibility between equivalence relations and quasi-orders is intimately related to the notion of cardinality of a set. 
If \( X , Y \) are arbitrary sets and \( \id ( X ) \) and \( \id ( Y ) \) denote, respectively, the equality relations on \( X \) and \( Y \), then 
\[ 
\card{X}\leq \card{Y} \IFF \id ( X ) \leq \id ( Y )
\]
and hence
\[ 
\card{X}= \card{Y} \IFF {\id ( X ) \sim \id ( Y )}.
\]
Thus reducibility between equivalence relations can be seen as a generalization of the notion of cardinality. 
Moreover, if \( E \) and \( F \) are equivalence relations on \( X \) and \( Y \), respectively, then
\begin{align*}
E \leq F & \IMPLIES \card{X / E} \leq \card{Y / F} 
\\
\shortintertext{and} 
E \sim F & \IMPLIES \card{X / E} = \card{Y / F}. 
\end{align*}
Assuming \( \AC \) the implications above can be reversed, since any surjection has a left inverse.
Under determinacy it is still possible to revert the implications when \( X \), \( Y \) are Polish spaces.
Working in \( \ZF \)
\begin{align}
 \card{X / E}\leq \card{Y / F} & \IFF \EXISTS{ R \subseteq X \times Y} \Bigl [ \FORALL{x \in X} \EXISTS{y \in Y} ( x \mathrel{R} y ) \AND \FORALL{x_1 , x_2 \in X} \FORALL{y_1 , y_2 \in Y } \notag
 \\
 &\phantom{{}\IFF \EXISTS{ R \subseteq X \times Y} \Bigl [{}} \bigl ( ( x_1 \mathrel{E} x_2 \wedge x_1\mathrel{R}y_1 \wedge x_2 \mathrel{R} y_2 \implies y_1 \mathrel{F} y_2 ) \wedge{} \label{eq:uniformizationquotients3}
 \\
 &\phantom{{}\IFF \EXISTS{ R \subseteq X \times Y} \Bigl [ \bigl({}}( y_1 \mathrel{F} y_2 \wedge x_1\mathrel{R}y_1 \wedge x_2 \mathrel{R} y_2 \implies x_1 \mathrel{E} x_2 ) \bigr )\Bigr ] ,\notag
\end{align}
so in order to get a reduction from \( E \) to \( F \) we may appeal to uniformization, as any \( g \colon X \to Y \) such that \( ( x , g ( x ) ) \in R \) for all \( x \in X \), witnesses \( E \leq F \).
The existence of such \( g \) follows form \( \ADR \) if \( E \) and \( F \) are arbitrary equivalence relations on Polish spaces, or from \( \AD + {\Vv = \Ll ( \R )} \) and \( E , F \) are \( \bDelta^2_1 \).
To see this notice that the statement on the right hand side of equivalence in~\eqref{eq:uniformizationquotients3} is \( \bSigma^{2}_{1} \) so the witness \( R \) can be taken to be \( \bSigma^{2}_{1} \) as well, and hence we can apply \( \bSigma^{2}_{1} \)-uniformization.

If a set \( A \) is the surjective image of \( X \) via some map \( f \) then, letting \( E_A \) be the equivalence relation on \( X \) defined by \( x \mathrel{E_A} x' \iff f ( x ) = f ( y ) \), the factoring map \( f \colon X / E_A \to A \), \( \eq{x}_{E_A} \mapsto f ( x ) \) is well-defined and witnesses \( \cardLR{ X / {E_A}} = \card{A} \). 
This observation can be turned into a method to compute cardinalities.
For example, in models of \( \AD \) one is mainly interested in the cardinalities of sets which are surjective images of a Polish space, and such cardinalities are called \markdef{small cardinalities}\index[concepts]{cardinals and cardinalities!small}, or \markdef{small cardinals}\label{pag:smallcardinalities} when we are dealing with well-orderable sets.
If \( A \), \( B \) are two sets whose cardinality is small, then \( {E_A \leq E_B} \implies { \card{A} \leq \card{B}} \) (again, in some special situations this implication can be reversed). 

\section{The generalized Cantor space} \label{sec:topologies} 
\subsection{Basic facts} \label{subsec:Cantorbasic}
The usual topology on the Cantor space \( \pre{\omega}{2} \) is the product topology, which is generated by the sets 
\[ 
\Nbhd_s^ \omega \coloneqq \setofLR{x \in \pre{\omega}{2} }{ s \subseteq x }
\] 
for \( s \in \pre{ < \omega }{2} \).
If \( \omega \) is replaced by some uncountable cardinal \( \kappa \), then there are at least two topologies on the \markdef{generalized Cantor space}\index[concepts]{generalized Cantor space} \( \pre{\kappa}{ 2} \) that can claim to be the natural generalization of the previous construction. 

\begin{definition} \label{def:generalizedCantor}
Let \( \kappa \) be an infinite cardinal, and let
\[ 
\Nbhd_s^ \kappa \coloneqq \setofLR{x \in \pre{\kappa}{2} }{ s \subseteq x }, \index[symbols]{Nbhd@\( \Nbhd_s^ \kappa \)}
\] 
where \( s \colon u \to 2 \) and \( u \subseteq \kappa \).
\begin{itemize}[leftmargin=1pc]
\item 
The \markdef{bounded topology}\index[concepts]{topology!bounded \( \tau_b \)}\index[symbols]{tau@\( \tau_b \), \( \tau_p \)} \( \tau_b \) on \( \pre{\kappa}{2} \) is the one generated by the collection 
\[
 \mathcal{B}_b \coloneqq \setofLR{ \Nbhd^\kappa_s }{ s \in \pre{< \kappa}{2} } . \index[symbols]{Bb@\( \mathcal{B}_b \), \( \mathcal{B}_p \)}
\] 
\item 
The \markdef{product topology}\index[concepts]{topology!product \( \tau_p \)} \( \tau_p \) on \( \pre{\kappa}{2} \) is the product of \( \kappa \) copies of \( 2 \) with the discrete topology, and hence it is generated by the collection 
\[
 \mathcal{B}_p \coloneqq \setofLR{ \Nbhd^\kappa_s }{ s \in \forcing{Fn} ( \kappa , 2 ; \omega ) }.
\]
\end{itemize}
The families \( \mathcal{B}_b \) and \( \mathcal{B}_p \) are called the \markdef{canonical basis} for \( \tau_b \) and \( \tau_p \), respectively, and their elements are called the \markdef{basic open sets}. 
\end{definition}

To simplify the notation, when \( \kappa \) is clear from the context we write \( \Nbhd_s \) instead of \( \Nbhd^\kappa_s \). 
Conversely, we write \( \tau_b ( \pre{\kappa}{2} ) \), \( \mathcal{B}_b ( \pre{\kappa}{2} ) \), \( \tau_p ( \pre{\kappa}{2} ) \), and \( \mathcal{B}_p ( \pre{\kappa}{2} ) \) instead of \( \tau_b \), \( \mathcal{B}_b \), \( \tau_p \), and \( \mathcal{B}_p \) if attention must be paid to the underlying space.
When \( \kappa = \omega \) the topologies \( \tau_b \) and \( \tau_p \) coincide, so when dealing with the Cantor space \( \pre{\omega}{2} \) we can talk about topology without further comments.
In contrast, Lemma~\ref{lem:cardinalityoftopologies} below shows that \( \tau_b \) and \( \tau_p \) are significantly different when \( \kappa \) is uncountable. 
For example 
\begin{equation}\label{eq:tau_bclosedunderintersections}
 \kappa \text{ regular} \IMPLIES \tau_b \text{ is closed under intersections of size } < \kappa 
\end{equation}
while \( \tau_p \) is never closed under countable intersections --- see Proposition~\ref{prop:topologicalproperties}\ref{prop:topologicalproperties-h} below. 

Generalizing the notion seen in Section~\ref{subsubsec:dsttrees}, a \markdef{descriptive set-theoretic tree on \( X \) of height \( \kappa \)}\index[concepts]{descriptive set-theoretic tree!of height \( \kappa \)} is a nonempty \( T \subseteq \pre{ < \kappa }{X} \) closed under initial segments.
Such \( T \) is called \markdef{pruned} if 
\[
 \FORALL{t \in T} \FORALL{ \alpha < \kappa } \EXISTS{t' \in T} [ \lh ( t' ) = \alpha \AND ( t \subseteq t' \vee t' \subseteq t ) ] .
\] 
Note that \( \emptyset \neq C \subseteq \pre{ \kappa }{2} \) is closed with respect to \( \tau_b \) if and only if it is of the form
\[
\body{T} \coloneqq \setofLR{x \in \pre{ \kappa }{2} }{ \FORALL{ \alpha < \kappa } ( x \restriction \alpha \in T ) } 
\]
with \( T \) a pruned tree on \( 2 \) of height \( \kappa \).

The literature on \( \pre{ \kappa }{2} \) with \( \kappa \) an uncountable cardinal (and on some other strictly related spaces, see Sections~\ref{subsec:Bairespace} and~\ref{subsec:spacesofmodels}) falls into two camps: most of the papers in general topology (see e.g.~\cite{Stone:1948st, Kuratowski:1966ku, Unlu:1982un, Tamano:1993tt, Chigogidze:1997cm, Chaber:1998cg, Kraszewski:2001kr, Iliadis:2012il, Chatyrko:2013ch}) deal with \( \tau_p \) or with the so-called box topology, while papers on infinitary logics (see~\cite{Vaught:1974kl, Mekler:1993kh, Halko:1996fu,Shelah:2000hs, Halko:2001kl,Shelah:2001bd, Shelah:2002lo, Shelah:2004ud, Dzamonja:2011ly, Friedman:2011nx}) use almost exclusively the bounded topology.
In this paper, instead, both the bounded and the product topology play an important role: although the statements of our main results refer to the bounded topology, several proofs make an essential use of (an homeomorphic copy of) the space \( ( \pre{\kappa}{2}, \tau_p) \) --- see Sections~\ref{sec:invariantlyuniversal} and~\ref{sec:alternativeapproach}. 
This seems to be a curious feature, and we currently do not know if the use of the product topology can be avoided at all.

\begin{remarks}\label{rmk:product}
\begin{enumerate-(i)}
\item \label{rmk:product-aa}
The name bounded topology comes from the fact that \( \tau_b \) can be equivalently defined as the topology generated by the collection of all sets of the form \( \Nbhd^\kappa_s \) with \( s \colon u \to 2 \) for some \emph{bounded} \( u \subseteq \kappa \): such a collection is an alternative basis for \( \tau_b \) which properly contains \( \mathcal{B}_b \).
 When \( \kappa \) is regular, this basis can be also described as the collection of all sets of the form \( \Nbhd^\kappa_s \) with \( s \colon u \to 2 \) for some \( u \subset \kappa \) of \emph{size \( {<} \kappa \)} (thus avoiding any reference to the ordering of \( \kappa \) in the definition of \( \tau_b \)): however, this is not true when \( \kappa \) is singular because if \( u \subseteq \kappa \) is \emph{cofinal} in \( \kappa \), then every set \( \Nbhd^\kappa_s \) for \( s \colon u \to 2 \) is a proper \( \tau_b \)-closed set.
\item\label{rmk:product-a}
The sets \( \Nbhd^\kappa_s \) with \( s \colon \setLR{ \alpha } \to 2 \) for \( \alpha \in \kappa \), form a subbasis \( \widetilde{\mathcal{B}}_p \) (which generates \( \mathcal{B}_p \)) for the product topology on \( \pre{\kappa}{2} \).
For ease of notation the elements of \( \widetilde{\mathcal{B}}_p \) are denoted by
\[
\widetilde{\Nbhd}\vphantom{\Nbhd}^\kappa_{ \alpha , i } = \widetilde{\Nbhd}_{ \alpha , i } \coloneqq \setofLR{ x \in \pre{ \kappa }{ 2 } }{ x ( \alpha ) = i } . \index[symbols]{Nbhdb@\( \widetilde{\Nbhd}_{ \alpha , i } \)}
\]
\item\label{rmk:product-b}
The definition of product topology makes sense for spaces of the form \( \pre{ A}{2} \) (or even \( \pre{ A}{X} \) for \( X \) a topological space) with \( A \) an arbitrary set, while the definition of the bounded topology requires that \( A \) be well-orderable --- see Section~\ref{subsec:spacesofmodels}.
\end{enumerate-(i)}
\end{remarks}

By Remark~\ref{rmk:product}\ref{rmk:product-aa}, if \( \kappa \) is regular then the bounded topology is generated by the sets \( \Nbhd^\kappa_s \) with \( s \in \forcing{Fn} ( \kappa , 2 ; \kappa ) \).
This suggests that when \( \kappa \) is regular \( \tau_p \) and \( \tau_b \) lie at the extrema of a spectrum of topologies on \( \pre{\kappa}{2} \).

\begin{definition}\label{def:lambdatopology} 
Let \( \lambda \leq \kappa \) be infinite cardinals.
The \markdef{\( \lambda \)-topology}\index[concepts]{topology!\( \lambda \)-topology \( \tau _ \lambda \)}\index[symbols]{taul@\( \tau_ \lambda \)} \( \tau_\lambda \) on \( \pre{ \kappa }{2} \) is generated by the basis
\[
\mathcal{B}_\lambda \coloneqq \setofLR{ \Nbhd^\kappa_s }{ s \in \forcing{Fn} ( \kappa , 2 ; \lambda ) } . \index[symbols]{Bbl@\( \mathcal{B}_ \lambda \)}
\]
In particular: \( \tau_p = \tau_ \omega \).
\end{definition}

As usual, when we need to explicitly refer to the underlying space we write \( \tau_\lambda ( \pre{\kappa}{2} ) \) and \( \mathcal{B}_\lambda ( \pre{\kappa}{2} ) \) instead of, respectively, \( \tau_\lambda \) and \( \mathcal{B}_\lambda \). 
It is easy to check that the sets in \( \mathcal{B}_p , \mathcal{B}_ \lambda , \mathcal{B}_b \) are clopen, and hence the topologies \( \tau_p , \tau _ \lambda , \tau _b \) are zero-dimensional.

\begin{lemma}\label{lem:boundedvscofinal}
Let \( \kappa \) be an uncountable cardinal, and consider the topologies \( \tau_p, \tau_ \lambda , \tau_b \) on \( \pre{ \kappa }{2} \).
\begin{enumerate-(a)}
\item\label{lem:boundedvscofinal-a}
if \( \omega \leq \lambda < \nu < \kappa \) are cardinals, then \( \mathcal{B}_ \lambda \subset \mathcal{B}_\nu \), each \( \Nbhd_s^ \kappa \in \mathcal{B}_ \nu \) is \( \tau_ \lambda \)-closed, and \( \tau_ \lambda \subset \tau_\nu \),
\item\label{lem:boundedvscofinal-b}
\( \tau_{\cof ( \kappa ) } \subseteq \tau_b \), and
\item\label{lem:boundedvscofinal-c}
 \( \tau_b = \tau_{\cof ( \kappa ) }\IFF \kappa = \cof ( \kappa ) \).
\end{enumerate-(a)} 
\end{lemma}

\begin{proof}
\ref{lem:boundedvscofinal-a}
It is clear that \( \mathcal{B}_ \lambda \subset \mathcal{B}_\nu \) and hence \( \tau_ \lambda \subseteq \tau_\nu \).
Fix \( \Nbhd_s \in \mathcal{B}_\nu \).
If \( \card{\dom s} < \lambda \), then \( \Nbhd_s \) is \( \tau_ \lambda \)-clopen, if \( \lambda \leq \card{\dom s} < \nu \), then \( \Nbhd_s = \bigcap_{ \alpha \in \dom s } \widetilde{\Nbhd}_{ \alpha , s ( \alpha ) } \) is \( \tau_ p \)-closed and hence \( \tau _ \lambda \)-closed, but it is not \( \tau_ \lambda \)-open as it does not contain any element of \( \mathcal{B} _ \lambda \).
Therefore \( \tau_ \lambda \neq \tau_ \nu \).

\smallskip

\ref{lem:boundedvscofinal-b}
To see that \( \tau_{\cof ( \kappa ) } \subseteq \tau_b \), it is enough to show that \( \mathcal{B}_{ \cof ( \kappa ) } \subseteq \tau_b \): if \( s \in \forcing{Fn} ( \kappa , 2 ; \cof ( \kappa ) ) \) then \( \Nbhd_s = \bigcup_{} \setof{ \Nbhd_t }{ t \supseteq s \AND t \in \pre{ \alpha }{2} } \), where \( \alpha = \sup \dom ( s ) < \kappa \).

\smallskip

\ref{lem:boundedvscofinal-c}
If \( \kappa \) is singular, and \( s \in \pre{ \alpha }{2} \) with \( \cof ( \kappa ) \leq \alpha < \kappa \), then \( \Nbhd_s \in \mathcal{B}_b \), but \( \Nbhd_s \) does not contain any set in \( \mathcal{B}_{ \cof ( \kappa ) }\), so it is not \( \tau_{\cof ( \kappa )} \)-open.
If \( \kappa \) is regular, then \( \mathcal{B}_b \subset \mathcal{B}_ \kappa \) so \( \tau_b \subseteq \tau_ \kappa \), and therefore \( \tau_b = \tau_ \kappa \).
\end{proof}

In Sections~\ref{sec:otherspacesandBairecategory} and~\ref{sec:infinitarylogics} we use these topologies that lie strictly between the product topology and the bounded topology, that is \( \tau _ \lambda \) with \( \omega < \lambda < \min ( \cof ( \kappa ) ^+ , \kappa ) \), while we have no use for \( \tau _ \lambda \) with \( \cof ( \kappa ) < \lambda < \kappa \).
Most of the nontrivial results on these intermediate topologies seem to require the axiom of choice. 

In order to simplify the notation, when \( \kappa \) is regular we stipulate the following 

\begin{convention}\label{cnv:basisfortau_kappa}
When \( \kappa \) is regular, the canonical basis for \( \tau_ \kappa \) is taken to be \( \mathcal{B}_b \), that is
\[
\mathcal{B}_ \kappa \coloneqq \setofLR{\Nbhd_s^ \kappa }{ s \in \pre{ < \kappa }{ 2 } } .
\]
 \end{convention}
 
\subsection{* More on \( 2^{ \kappa } \)}
\subsubsection{Lipschitz and continuous reductions} \label{subsubsec:lipschitz}
\begin{definition}\label{def:Lipschitzfunction}
Let \( \kappa \) be an infinite cardinal.
A function \( \varphi \colon \pre{ < \kappa }{2} \to \pre{ < \kappa }{2} \) is 
\begin{itemize}[leftmargin=1pc]
\item
\markdef{monotone} if \( \FORALL{s , t \in \pre{ < \kappa }{2}} ( s \subseteq t \implies \varphi ( s ) \subseteq \varphi ( t ) ) \),
\item
\markdef{Lipschitz} if it is monotone and \( \FORALL{s \in \pre{ < \kappa }{2} } ( \lh ( \varphi ( s ) ) = \lh ( s ) ) \),
\item
\markdef{continuous} if it is monotone and \( \lh ( \bigcup_{ \alpha < \kappa } \varphi ( x \restriction \alpha ) ) = \kappa \), for all \( x \in \pre{ \kappa }{2} \).
\end{itemize}
If \( \varphi \) is Lipschitz then it is continuous, and if \( \varphi \) is Lipschitz or continuous, then 
\[
f_ \varphi \colon \pre{ \kappa }{2} \to \pre{ \kappa }{2} , \quad x \mapsto \bigcup_{ \alpha < \kappa } \varphi ( x \restriction \alpha ) 
\]
is the \markdef{function induced by \( \varphi \)}.
\end{definition}

\begin{lemma}
A function \( f \colon \pre{ \kappa}{2} \to \pre{ \kappa }{2} \) is continuous with respect to \( \tau_b \) if and only if it is of the form \( f_ \varphi \) for some continuous \( \varphi \colon \pre{ < \kappa }{2} \to \pre{ < \kappa }{2} \).
\end{lemma}

\begin{proof}
If \( f \) is \( \tau_b \)-continuous, then the map \( \varphi \) defined by setting for \( s \in \pre{< \kappa}{2} \)
\[ 
\varphi ( s ) \coloneqq \text{the longest \( t \) with length \( \leq \lh ( s ) \) such that } f ( \Nbhd_s ) \subseteq \Nbhd_t 
\]
is continuous and \( f = f_ \varphi \).
Conversely if \( f = f_ \varphi \) with \( \varphi \) continuous then \( f^{-1} ( \Nbhd_t ) = \bigcup_{\varphi ( s ) \supseteq t} \Nbhd_s\), and hence \( f \) is \( \tau_b \)-continuous.
\end{proof}

We call a function \( f \colon \pre{ \kappa}{2} \to \pre{ \kappa }{2} \) \markdef{Lipschitz} if it is of the form \( f_ \varphi \) for some Lipschitz \( \varphi \) or, 
equivalently, if it is such that for all \( x,y \in \pre{\kappa}{2} \) and \( \alpha < \kappa \), \( x \restriction \alpha = y \restriction \alpha \IMPLIES f ( x ) \restriction \alpha = f ( y ) \restriction \alpha \).
(The reason for this terminology is that when \( \kappa = \omega \) then \( f \) is Lipschitz if and only if \( d ( f ( x ) , f ( y ) ) \leq d ( x , y ) \), where \( d \) is the usual metric on \( \pre{\omega}{2} \).)
It is immediate to check that the composition of Lipschitz functions is Lipschitz, and that Lipschitz functions are \( \tau_b \)-continuous.

\begin{definition} \label{def:leqL}
Let \( A , B \subseteq \pre{ \kappa }{2} \).
 We say that \( A \) is \markdef{Lipschitz reducible}\index[concepts]{Lipschitz reducibility} to \( B \), in symbols 
\[ 
 A \leql^\kappa B , \index[symbols]{22@\( \leql^\kappa \), \( \leqW^ \kappa \)}
\]
if \( A = f^{-1} ( B ) \) for some Lipschitz \( f \colon \pre{ \kappa}{2} \to \pre{ \kappa }{2} \).
Equivalently, \( A \leql^\kappa B \) if and only if Player \( \II \) has a winning strategy in \( \GL^\kappa ( A , B ) \), the \markdef{Lipschitz game of length \( \kappa \) for \( A , B \)}.
It is a zero-sum, perfect information game of length \( \kappa \), in which at each inning \( \alpha < \kappa \) the two-players \( \I \) and \( \II \) play \( x_ \alpha , y_ \alpha \in 2 \) with \( \I \) playing first and at limit levels:
\[
\begin{tikzpicture}[>=stealth]
\node (II) at (0,0) [anchor=base] { \( \II \)};
\node (I) at (0,1) [anchor=base] {\( \I \)};
\node (x_0) at (0.7,1) [anchor=base] {\( x_0 \)};
\node (x_1) at (1.2,0) [anchor=base] {\( y_0 \)};
\node (x_2) at (1.7,1) [anchor=base] {\( x_1 \)} ;
\node (x_3) at (2.2,0) [anchor=base] {\( y_1 \)};
\node at (3,1) [anchor=base] {\( \cdots \)} ;
\node at (3,0) [anchor=base] {\( \cdots \)};
\node at (4,1) [anchor=base] {\( x_ \alpha \)} ;
\node at (4.5,0) [anchor=base] {\( y_ \alpha \)};
\node at (5.5,1) [anchor=base] {\( \cdots \)};
\node at (5.5,0) [anchor=base] {\( \cdots \)};
\node (dots) at (5.5,0) [anchor=base] {\( \cdots \)};
\node (A) at (x_3.south -| II.east) {}; 
\node (B) at (I.north -| II.south east) {}; 
\node (C) at ($0.5*(A) +0.5 *(B)$){}; 
\draw (II.south east) -- (I.north -| II.south east) 
(II.north west |- C.center) -- (dots.east |- C.center); 
\begin{pgfonlayer}{background} 
\fill[gray!30,rounded corners] 
(II.south west) rectangle (dots.east |- I.north); 
\end{pgfonlayer} 
\node at (-1.5,0.5) [anchor=base] { \( \GL^\kappa ( A , B ) \)};
\end{tikzpicture} 
\]\index[symbols]{Game@ \( \GL^\kappa ( A , B ) \), \( \GW^\kappa ( A , B ) \)}
Player \( \II \) wins just in case
\[
\seqofLR{ x_ \alpha }{ \alpha < \kappa } \in A \IFF \seqofLR{ y_ \alpha }{ \alpha < \kappa } \in B .
\]
Similarly, we say that \( A \) is \markdef{continuously reducible} or \markdef{Wadge reducible}%
\footnote{W.W.\ Wadge initiated the use of games to study continuous reducibility on the Baire space in~\cite{Wadge:1983sp}.}
to \( B \), in symbols \( A \leqW^\kappa B \), \label{leqW} if there is a continuous \( f \colon \pre{ \kappa }{2} \to \pre{ \kappa }{2} \) such that \( A = f^{-1} ( B ) \).
Equivalently: \( A \leqW^\kappa B \) if and only if Player \( \II \) has a winning strategy in \( \GW^\kappa ( A , B ) \), the \markdef{Wadge game of length \( \kappa \) for \( A , B \)}.
This game is similar to \( \GL^\kappa ( A , B ) \), except that \( \II \) can pass at any round, provided that at the end a sequence of length \( \kappa \) is produced.
\end{definition}

A set \( A \subseteq \pre{ \kappa }{2} \) is \markdef{\( \leql^\kappa \)-hard for a collection} \( \mathcal{S} \subseteq \pow ( \pre{ \kappa }{2} ) \) if \( \FORALL{B \in \mathcal{S}} ( B \leql^\kappa A) \).
A set which is \( \leql^\kappa \)-hard for \( \mathcal{S} \) and moreover belongs to \( \mathcal{S} \) is said to be \markdef{\( \leql^\kappa \)-complete for} \( \mathcal{S} \).

\begin{proposition}\label{prop:hardness}
Suppose \( \mathcal{S} \subseteq\pow ( \pre{ \kappa }{2} ) \) and endow \( \pre{ \kappa }{2} \) with the bounded topology.
\begin{enumerate-(a)}
\item\label{prop:hardness-a}
If there is \( A \subseteq \pre{ \kappa }{2} \) which is \( \leql^\kappa \)-hard for \( \mathcal{S} \), then \( \mathcal{S} \neq \pow ( \pre{ \kappa }{2} ) \). 
\item\label{prop:hardness-b}
Suppose \( \mathcal{S} \) is closed under complements and continuous preimages, i.e.~\( A \in \mathcal{S} \) implies \( f^{-1} ( A ) \in \mathcal{S} \) for all continuous \( f \colon \pre{ \kappa }{2} \to \pre{ \kappa }{2} \).
Then there is no \( \leql^\kappa \)-complete set for \( \mathcal{S} \).
\end{enumerate-(a)}
\end{proposition}

\begin{proof}
\ref{prop:hardness-a}
If \( \pre{ \kappa }{2} \setminus A \nleql^\kappa A \), then \( \pre{ \kappa }{2} \setminus A \) already witnesses that \( \mathcal{S} \neq \pow ( \pre{ \kappa }{2} ) \), so we may assume that \( \pre{ \kappa }{2} \setminus A \leql^\kappa A \). 
Let \( g \colon \pre{ \kappa }{2} \to \pre{ \kappa }{2} \) be the continuous map defined by \( g ( x ) ( \alpha ) \coloneqq x ( \alpha + 1 ) \) for all \( \alpha < \kappa \), and let \( \bar{A} \coloneqq g^{-1} ( A ) \).
Player \( \I \) has a winning strategy in \( \GL^\kappa ( \bar{ A} , A ) \) by playing \( 0 \) (or \( 1 \) for that matter) at round \( 0 \) and at all limit rounds, and by following \( \sigma \) at all successor rounds, where \( \sigma \) is a winning strategy for \( \II \) in \( \GL^\kappa ( A , \pre{ \kappa }{2} \setminus A ) \).
Therefore \( \bar{A} \nleql^\kappa A \), and hence \( \bar{A} \notin \mathcal{S} \).

\smallskip

\ref{prop:hardness-b}
Towards a contradiction, suppose \( A \in \mathcal{S} \) is \( \leql^\kappa \)-complete for \( \mathcal{S} \).
Since \( \pre{ \kappa }{2} \setminus A \in \mathcal{S} \) by closure under complements, then \( \pre{ \kappa }{2} \setminus A \leql^\kappa A \). 
Then the set \( \bar{A} \) defined as in part~\ref{prop:hardness-a} contradicts the choice of \( A \): in fact, the function \( g \) witnesses \( \bar{A} \leqW^\kappa A \), so that \( \bar{A} \in \mathcal{S} \), but \( \bar{A} \nleql^\kappa A \).
\end{proof}

\subsubsection{Properties of \( \tau_p \), \( \tau_ \lambda \) and \( \tau_b \)}\label{subsubsec:Propertiesoftaus}

For the space \( \pre{\omega}{2} \) the cardinality of the topology equals the cardinality of the continuum
\[
\card{\tau_p ( \pre{\omega}{2} ) } = \card{\tau_b ( \pre{\omega}{2} )} = \card{ \pre{\omega}{2} } .
\]
The situation for \( \pre{ \kappa }{2} \) when \( \kappa > \omega \) is quite different.
First of all note that the following trivial facts:
\begin{equation}\label{eq:topologiesincluded}
\tau_p \subseteq \tau_ \lambda \AND \tau_p \subseteq \tau_b ,
\end{equation}
and
\begin{equation}\label{eq:cardinalityoftopologies}
\text{if \( \mathcal{B} \) is a basis for a topology \( \tau \), then \( \card{ \tau} \leq \card{ \pow ( \mathcal{B} ) } \),}
\end{equation}
This is an immediate consequence of the fact that \( \tau \into \pow ( \mathcal{B} ) \), \( U \mapsto \setofLR{B \in \mathcal{B} }{ B \subseteq U} \), is injective.

\begin{lemma}\label{lem:cardinalityoftopologies}
Let \( \lambda < \kappa \) be uncountable cardinals, and let \( \tau_p \coloneqq \tau_p ( \pre{\kappa}{2} ) \), \( \tau_ \lambda \coloneqq \tau_ \lambda ( \pre{ \kappa }{2} ) \), and \( \tau_b \coloneqq \tau_b ( \pre{\kappa}{2} ) \).
\begin{enumerate-(a)}
\item\label{lem:cardinalityoftopologies-a}
\( \card{\tau_p } = \card{\pow ( \kappa ) }\).
\item\label{lem:cardinalityoftopologies-b}
\( \pow ( \pre{\omega}{2} ) \into \tau_ \lambda \) and \( \pow ( \pre{\omega}{2} ) \into \tau_b \) and therefore \( \tau_ \lambda \onto \pow ( \pre{\omega}{2} ) \) and \( \tau_b \onto \pow ( \pre{\omega}{2} ) \). 
Moreover \( \card{\pow ( \pre{< \lambda}{2})} \leq \card{\tau_ \lambda } \leq \card{ \pow ( [ \kappa ]^{< \lambda } ) }\) and \( \card{\tau_b } = \card{ \pow ( \pre{ < \kappa }{2} ) }\).
\item\label{lem:cardinalityoftopologies-c}
Assume \( \AD \).
If \( \kappa \) is a surjective image of \( \R \) (i.e.~\( \kappa < \Theta \), where \( \Theta \) is as in Definition~\ref{def:Theta}), then \( \card{\tau_p } < \card{\tau_ \lambda } \) and \( \card{\tau_p } < \card{\tau_b } \).
\item\label{lem:cardinalityoftopologies-d}
Assume \( \AC \).
Then 
\begin{itemize}
\item
\( 2^ \kappa < 2^{( 2 ^{< \lambda })} \implies \card{ \tau_p } < \card{\tau_ \lambda } \),
\item
\( 2^ \kappa < 2^{(2^{< \kappa })} \iff \card{ \tau_p } < \card{\tau_b } \), 
\item
\( 2^{( \kappa ^{< \lambda })} < 2^{(2^{< \kappa })} \wedge \lambda < \min ( \cof ( \kappa )^+ , \kappa ) \implies \card{ \tau_ \lambda } < \card{\tau_b } \).
\end{itemize}
\end{enumerate-(a)}
\end{lemma}

\begin{proof}
\ref{lem:cardinalityoftopologies-a}
Since \( \mathcal{B}_p = \setofLR{\Nbhd_s}{ s \in \forcing{Fn} ( \kappa , 2 ; \omega ) } \) is a basis for \( \tau_p \) and \( \card{\forcing{Fn} ( \kappa , 2 ; \omega ) } = \kappa \), then \( \card{\tau_p} \leq \card{ \pow ( \kappa ) } \) by~\eqref{eq:cardinalityoftopologies}.
For the other inequality use the injection \( \pow ( \kappa ) \into \tau_p \), \( A \mapsto \bigcup_{ \alpha \in A} \widetilde{\Nbhd}\vphantom{\Nbhd}_{ \alpha , 1 }\).

\smallskip

\ref{lem:cardinalityoftopologies-b}
The map \( \pow ( \pre{\omega}{2} ) \to \pow ( \pre{ \kappa }{2} ) \), \( \pre{\omega}{2} \supseteq A \mapsto \bigcup_{s \in A} \Nbhd _s \), witnesses that \( \card{\pow ( \pre{\omega}{2} )} \leq \card{\tau_ \lambda} \) and \( \card{ \pow ( \pre{\omega}{2} ) } \leq \card{ \tau_b } \).

The function \( \forcing{Fn} ( \kappa , 2 ; \lambda ) \to [ \kappa ]^{< \lambda } \) defined by
\[
 s \mapsto \setof{ 2 \alpha }{ \alpha \in \dom s \wedge s ( \alpha ) = 0 } \cup \setof{ 2 \alpha + 1 }{ \alpha \in \dom s \wedge s ( \alpha ) = 1 }
\]
is an injection.
Thus \( \card{ [ \kappa ]^{< \lambda } } \geq \card{\forcing{Fn} ( \kappa , 2 ; \lambda ) } = \card{ \mathcal{B}_ \lambda } \), and hence \( \card{ \tau_ \lambda } \leq \card{ \pow ( [ \kappa ]^{< \lambda } ) } \) by~\eqref{eq:cardinalityoftopologies}.
Similarly, as \( \card{\mathcal{B}_b} = \card{ \pre{ < \kappa }{2} } \), then \( \card{\tau_b} \leq \card{\pow ( \pre{ < \kappa }{2} ) } \).
In order to prove the other inequalities, consider the map
\begin{equation} \label{eq:v}
U \colon \pre{ < \kappa }{2} \into \mathcal{B}_ b \subseteq \tau_b, \quad s \mapsto \Nbhd _{ 0^{( \lh s ) } \conc 1 \conc s} ,
\end{equation}
and notice that if \( s \in \pre{<\lambda}{2} \), then \( U(s) \in \mathcal{B}_\lambda \subseteq \tau_\lambda \). 
Then the injection \( \pow ( \pre{ < \kappa }{2} ) \into \tau_b \), \( A \mapsto \bigcup_{s \in A} U ( s ) \), witnesses \( \card{\pow ( \pre{ < \kappa }{2} ) } \leq \card{\tau_b} \), while its restriction to \( \pow ( \pre{< \lambda}{2}) \) has range contained in \( \tau_\lambda \) and witnesses \( \card{\pow ( \pre{< \lambda}{2})} \leq \card{\tau_ \lambda } \).

\smallskip

\ref{lem:cardinalityoftopologies-c}
By a standard result in determinacy (see Theorem~\ref{th:coding} below), there is a surjection of \( \pre{\omega}{2} \) onto \( \pow ( \kappa ) \), and hence \( \pre{\omega}{2} \onto \tau_p \) by~\ref{lem:cardinalityoftopologies-a}.
Let \( \tau \) be either \( \tau _ \lambda \) or \( \tau _b \), so that \( \card{\tau_p} \leq \card{\tau} \) by~\eqref{eq:topologiesincluded}.
If \( \card{\tau_p} = \card{\tau} \) then \( \pre{\omega}{2} \onto \tau \), and therefore \( \pre{\omega}{2} \onto \pow ( \pre{\omega}{2} ) \) by~\ref{lem:cardinalityoftopologies-b}, a contradiction.
Therefore \( \card{\tau_p} < \card{\tau } \). 

\smallskip

\ref{lem:cardinalityoftopologies-d} 
follows from~\ref{lem:cardinalityoftopologies-a} and~\ref{lem:cardinalityoftopologies-b}.
\end{proof} 

\begin{remark}
Assume \( \AD \) and let \( \kappa \) be such that \( \R \onto \kappa \). 
Although \( \R \) does not surject onto \( \tau_b ( \pre{\kappa}{2}) \) by Lemma~\ref{lem:cardinalityoftopologies}\ref{lem:cardinalityoftopologies-b}, at least we have that \( \R \onto \mathcal{B}_b ( \pre{\kappa}{2}) \). 
This is because \( \card{\mathcal{B}_b ( \pre{\kappa}{2}) } = \card{\pre{ < \kappa}{2}} \), \( \pre{\kappa}{2} \onto \pre{ < \kappa}{2} \), and \( \R \onto \pre{\kappa}{2} \) by Theorem~\ref{th:coding}.
As \( [ \kappa ]^{< \lambda } \into \pre{<\kappa}{2} \), this implies that \( \R \onto \mathcal{B}_ \lambda ( \pre{\kappa}{2}) \) for all \( \omega < \lambda < \min ( \cof ( \kappa ) ^+ , \kappa ) \), as well.
\end{remark}

The next result summarizes some properties of the topologies \( \tau_p \), \( \tau _ \lambda \) and \( \tau_b \) on \( \pre{ \kappa }{2} \).
We denote with \( \mathcal{C}_p \), \( \mathcal{C}_ \lambda \), \( \mathcal{C}_b \)\index[symbols]{Cb@\( \mathcal{C}_b \), \( \mathcal{C}_ \lambda \), \( \mathcal{C}_p \)} their algebras of clopen sets.

\begin{proposition}\label{prop:topologicalproperties}
Let \( \lambda < \kappa \) be uncountable cardinals, and consider the space \( \pre{ \kappa }{2} \) with one of the topologies \( \tau _p , \tau _ \lambda , \tau _b \).
\begin{enumerate-(a)}
\item\label{prop:topologicalproperties-a}
The topologies \( \tau _p , \tau _ \lambda , \tau _b \) are perfect, regular Hausdorff, and zero-dimensional.
\item\label{prop:topologicalproperties-b}
\( \tau_p \) is compact, while \( \tau_ \lambda \) and \( \tau_b \) are not.
\item\label{prop:topologicalproperties-c}
Let \( \mathcal{U}_p \), \( \mathcal{U}_ \lambda \), and \( \mathcal{U}_b \) be arbitrary open neighborhood bases of some \( x \in \pre{\kappa}{2} \) with respect to \( \tau_p \), \( \tau _ \lambda \), and \( \tau_b \).
Then \( \mathcal{U}_p \onto \kappa \), \( \mathcal{U}_b \onto \cof ( \kappa ) \), and (assuming \( \AC \)) \( \mathcal{U}_ \lambda \onto \kappa \); moreover there are \( \mathcal{U}_p \) and \( \mathcal{U}_b \) as above such that \( \card{\mathcal{U}_p} = \kappa \) and \( \card{\mathcal{U}_b } = \cof ( \kappa ) \). 
Therefore \( \tau_p \) and (assuming \( \AC \)) \( \tau_ \lambda \) are never first countable, while \( \tau_b \) is first countable if and only if \( \cof ( \kappa ) = \omega \).
\item \label{prop:topologicalproperties-d}
Let \( \mathcal{B}'_p \), \( \mathcal{B}'_ \lambda \) and \( \mathcal{B}'_b \) be arbitrary bases for the topologies \( \tau_p \), \( \tau _ \lambda \) and \( \tau_b \), respectively, on \( \pre{\kappa}{2} \).
Then \( \mathcal{B}'_p \onto \kappa \), \( \mathcal{B}'_ \lambda \onto \pre{< \lambda}{2} \) (and, assuming \( \AC \), also \( \mathcal{B}'_\lambda \onto \kappa \)), and \( \mathcal{B}'_b \onto \pre{< \kappa}{2} \). 
Therefore the topologies \( \tau_p \), \( \tau _ \lambda \) and \( \tau_b \) are never second countable.
\item\label{prop:topologicalproperties-e}
The topology \( \tau_p \) is separable if and only if \( \kappa \leq \card{\pre{\omega}{2}} \). 
Therefore under \( \AD \) the topology \( \tau_p \) is never separable, while under \( \AC \) it is separable if and only if \( \kappa \leq 2^{\aleph_0} \).
\item\label{prop:topologicalproperties-e'}
If \( D \subseteq \pre{\kappa}{2} \) is \( \tau_\lambda \)-dense, then \( D \onto \pre{< \lambda}{2} \), hence \( \tau_\lambda \) is never separable. 
Moreover there is a \( \tau_ \lambda \)-dense set of size \( \card{ [ \kappa ]^{< \lambda } } \). 
Thus assuming \( \AC \) we get that the density character of \( \tau_\lambda \) is between \( 2^{< \lambda} \) and \( \kappa^{< \lambda} \). 
Further assuming that \( \lambda \) is inaccessible, one has that \( \tau_\lambda \) has density \( \lambda \) if and only if \( \kappa \leq 2^\lambda \).
\item\label{prop:topologicalproperties-f}
If \( D \subseteq \pre{\kappa}{2} \) is \( \tau_b \)-dense, then \( D \onto \pre{< \kappa}{2} \). 
Moreover there is a \( \tau_b \)-dense set of size \( \card{\pre{< \kappa}{2}} \).
In particular \( \tau_b \) is never separable, and under \( \AC \) it has density character \( 2^{< \kappa} \).
\item\label{prop:topologicalproperties-g}
The topology \( \tau_b \) is completely metrizable if and only if it is metrizable if and only if \( \cof ( \kappa ) = \omega \); the topologies \( \tau_p \) and (assuming \( \AC \)) \( \tau_ \lambda \) are never metrizable. 
\item\label{prop:topologicalproperties-h}
The topology \( \tau_* \) with \( * \in \setLR{ p , \lambda , b } \) is closed under intersections of length \( \leq \alpha \) (for some ordinal \( \alpha \)) if and only if: \( \alpha < \omega \) when \( * = p \), \( \alpha < \lambda \) when \( * = \lambda \) (assuming \( \AC \)), \( \card{\alpha} < \cof ( \kappa ) \) when \( * = b \).
In particular~\eqref{eq:tau_bclosedunderintersections} holds.
\item\label{prop:topologicalproperties-i}
\( \mathcal{C}_p \) is a \( \omega \)-algebra, that is an algebra of sets; \( \mathcal{C}_ \lambda \) is a \( \lambda \)-algebra (assuming \( \AC \)); \( \mathcal{C}_b \) is a \( \cof ( \kappa ) \)-algebra. 
\end{enumerate-(a)}
\end{proposition}

\begin{proof}
Part~\ref{prop:topologicalproperties-a} is trivial. 

\smallskip

\ref{prop:topologicalproperties-b}
Tychonoff's theorem for \( ( \pre{ \kappa }{2} , \tau_p ) \) is provable in \( \ZF \) (see e.g.~\cite[Theorem 4]{Keremedis:2000fk}), and \( \setofLR{\Nbhd_s}{s \in \pre{\omega}{2} } \) is an infinite clopen partition of \( ( \pre{ \kappa }{2} , \tau_ \lambda ) \) and \( ( \pre{ \kappa }{2} , \tau_b ) \).
 
\smallskip

\ref{prop:topologicalproperties-c}
Fix \( x \in \pre{\kappa}{2} \), and let \( \mathcal{U} \) be one of \( \mathcal{U}_p \), \( \mathcal{U}_ \lambda \), \( \mathcal{U}_ b \).
We construct a map sending each \( U \in \mathcal{U} \) to some \( s ( U ) \subseteq x \) such that 
\begin{equation}\label{eq:prop:topologicalproperties}
\Nbhd_{s ( U )} \subseteq U, 
\end{equation}
and consider the neighborhood base \( \mathcal{U}' = \setof{ \Nbhd_{s ( U )} }{ U \in \mathcal{U} } \). 
As \( \mathcal{U} \onto \mathcal{U}' \), it is enough to show that there is a map from \( \mathcal{U}' \) onto \( \kappa \) or \( \cof ( \kappa ) \).
\begin{description}
\item[Case \( \mathcal{U} = \mathcal{U}_p \)]
choose \( v ( U ) \in [ \kappa ]^{< \omega } \) such that~\eqref{eq:prop:topologicalproperties} holds when \( s ( U ) = x \restriction v ( U ) \).
The axiom of choice is not needed here, as \( [ \kappa ]^{< \omega } \) is well-orderable.
If \( \alpha < \kappa \) and \( U \in \mathcal{U} \) is contained in \( \widetilde{\Nbhd}_{ \alpha , x ( \alpha ) } \), then \( \alpha \in v ( U ) \); this implies that \( \kappa = \bigcup_{} \setofLR{ v ( U )}{ U \in \mathcal{U}_p} \), and therefore \( \mathcal{U}' \onto \kappa \).

\item[Case \( \mathcal{U} = \mathcal{U}_ \lambda \)]
choose \( s ( U ) \in \forcing{Fn} ( \kappa , 2 ; \lambda ) \) so that~\eqref{eq:prop:topologicalproperties} holds, and let \( S = \setof{ \dom s ( U ) }{ U \in \mathcal{U} } \).
If \( \alpha < \kappa \) and \( U \in \mathcal{U} \) is contained in \( \widetilde{\Nbhd}_{ \alpha , x ( \alpha ) } \), then \( \alpha \in \dom s ( U ) \) and hence \( \kappa = \bigcup S \).
As every element of \( S \) has size \( < \lambda \), then \( \kappa \leq \lambda \cdot \card{S} \), so \( \kappa \leq \card{S} \).
Thus \( \mathcal{U}' \onto \kappa \).

\item[Case \( \mathcal{U} = \mathcal{U}_b \)]
choose \( \alpha ( U ) \in \kappa \) such that~\eqref{eq:prop:topologicalproperties} is satisfied when \( s ( U ) = x \restriction \alpha ( U ) \).
If \( \cof ( \kappa ) \nleq \card{ \mathcal{U} '} \), then \( \alpha \coloneqq \sup_{ U \in \mathcal{U}_b} \alpha ( U ) + 1 < \kappa \), so the open neighborhood \( \Nbhd_{ x \restriction \alpha } \) does not contain any element of the open neighborhood basis \( \mathcal{U}' \), a contradiction. 
\end{description}
The \( \tau_p \)-neighborhhood base \( \setof{ \Nbhd_s \in \mathcal{B}_p }{ s \subseteq x } \) has size \( \kappa \), and the \( \tau_b \)-neighborhhood base \( \setof{ \Nbhd_{ x \restriction \alpha _i } \in \mathcal{B}_b }{ i < \cof ( \kappa ) } \) has size \( \cof ( \kappa ) \), where the \( \alpha _i \) are cofinal in \( \kappa \).

The fact that \( ( \pre{ \kappa }{2} , \tau_b ) \) is first countable when \( \cof ( \kappa ) = \omega \) is immediate. 

\smallskip

\ref{prop:topologicalproperties-d} 
Fix \( x \in \pre{ \kappa }{2} \).
The set \( \mathcal{U}_p = \setof{ U \in \mathcal{B}'_p}{ x \in U } \) is a \( \tau_p \)-neighborhood base of \( x \), so \( \mathcal{U}_p \onto \kappa \) by part~\ref{prop:topologicalproperties-c}, whence \( \mathcal{B}'_p \onto \kappa \). A similar argument shows that \( \mathcal{B}'_\lambda \onto \kappa \) when assuming \( \AC \).

The map \( U \) defined in~\eqref{eq:v} has the property that \( U ( s ) \cap U ( t ) = \emptyset \) for \( s \neq t \).
Therefore the map \( \mathcal{B}'_b \to \pre{ < \kappa }{2} \) defined by
\[
B \mapsto 
 \begin{cases}
 \emptyset & \text{if } \FORALL{s\in \pre{ < \kappa }{2} } ( B \nsubseteq U ( s ) ),
 \\
 s & \text{if } B \subseteq U ( s ) ,
 \end{cases}
\]
is a well-defined surjection, and similarly for \( \mathcal{B}'_ \lambda \onto \pre{< \lambda}{2} \).

\smallskip

\ref{prop:topologicalproperties-e} 
We first show that \( \pre{I}{2} \) with the product topology is separable, when \( I \subseteq \pre{ \omega }{2} \).%
\footnote{This is a particular case of the more general statement that the product of \( \kappa \)-many separable spaces is separable if and only if \( \kappa \leq 2^{\aleph_0} \)~\cite[Exercises 3 and 4, page 86]{Kunen:1980fk}.
The reason for including the proof here is to show that it holds in \( \ZF \).}
The set
\[
E_n \coloneqq \setofLR{ \varphi \in \pre{ I}{2} }{ \FORALL{f , g \in I} \left ( f \restriction n = g \restriction n \IMPLIES \varphi ( f ) = \varphi ( g ) \right ) } 
\]
is countable since it is in bijection with the set of all functions from \( \pre{n}{2} \) to \( 2 \), and hence 
\[
E \coloneqq \bigcup_{n \in \omega } E_n
\]
is countable.
We now argue that \( E \) is dense by showing that it intersects every basic open set of \( \pre{I}{2} \) of the form
\[
V \coloneqq \setofLR{ \varphi \in \pre{I}{2} }{ \varphi ( f_1 ) = i_1 \AND \dots \AND \varphi ( f_m ) = i_m }
\]
with \( f_1 , \dots , f_m \in I \) and \( i_1 , \dots , i_m \in \set{0 , 1}\).
Let \( n \) be large enough so that \( f_1 \restriction n , \dots , f_m \restriction n \) are all distinct, and let \( \varphi \in E_n \) be such that \( \varphi ( f_j ) = i_j \) (\( j = 1, \dots , m \)).
Then \( \varphi \in E \cap V \) as required.

Conversely, suppose \( \kappa \nleq \card{ \pre{\omega}{2} } \).
Given any set \( \setofLR{ f_n }{ n \in \omega } \subseteq \pre{ \kappa }{2} \) let \( F \colon \kappa \to \pre{\omega}{2} \) be defined by
\[ 
F ( \alpha ) ( n ) \coloneqq f_n ( \alpha ) .
\]
By case assumption \( F \) cannot be injective, and hence there are \( \alpha < \beta < \kappa \) such that \( F ( \alpha ) = F ( \beta ) \), that is \( f_n ( \alpha ) = f_n ( \beta ) \) for all \( n \).
Thus \( \setofLR{ f_n }{ n \in \omega } \) is disjoint from the basic open set \( \setofLR{f \in {}^ \kappa 2}{ f ( \alpha ) = 0 \AND f ( \beta ) = 1 } \) and therefore it is not dense.

The result under \( \AD \) follows from the fact that \( \omega _1 \) does not embed into \( \pre{\omega}{2} \).

\smallskip

\ref{prop:topologicalproperties-e'}
Let \( D \subseteq \pre{\kappa}{2} \) be \( \tau_\lambda \)-dense. As the sets \( U(s) \) defined in~\eqref{eq:v} are pairwise disjoint, we obtain a surjection \( D \onto \pre{< \lambda}{2} \) by setting
\begin{equation} \label{eq:densitychar}
x\mapsto
\begin{cases}
\emptyset &\text{if \( x \notin \bigcup_{s \in \pre{ < \lambda }{2} } U ( s ) \)}
\\
s & \text{if } x \in U ( s ) .
\end{cases}
\end{equation}
Moreover, the set \( \setofLR{ \chi_u }{ u \in [ \kappa ]^{< \lambda }} \) is \( \tau _ \lambda \)-dense and has cardinality \( \card{ [ \kappa ]^{< \lambda } } \), where \( \chi_u \colon \kappa \to \setLR{0 , 1} \) is the characteristic function of \( u \). 
Finally, the additional part concerning the case when \( \lambda\) is inaccessible (assuming \( \AC \)) can be proved by replacing \( \omega\) with \( \lambda\) in the proof of part~\ref{prop:topologicalproperties-e}.

\smallskip

\ref{prop:topologicalproperties-f}
The argument is similar to that of part~\ref{prop:topologicalproperties-e'}. If \( D \) is a dense set in \( ( \pre{ \kappa }{2} , \tau_b ) \), then one gets a surjection \( D \onto \pre{< \kappa}{2} \) by replacing \( \lambda\) with
\( \kappa \) in~\eqref{eq:densitychar}. Moreover, the set \( \setofLR{ s \conc 1 \conc 0^{( \kappa )}}{ s \in \pre{< \kappa}{2}} \) is \( \tau _b \)-dense and has cardinality \( \card{\pre{<\kappa}{2}} \).

\smallskip

\ref{prop:topologicalproperties-g}
Since a metric space is first countable, by part~\ref{prop:topologicalproperties-c} it is enough to show that \( ( \pre{\kappa}{2} , \tau_b ) \) is completely metrizable when \( \cof ( \kappa ) = \omega \). 
This easily follows from the application of the Birkhoff-Kakutani theorem~\cite[\S 1.22, p.~34]{Montgomery:1955fk} to the first countable Hausdorff topological group \( ( \pre{\kappa}{2} , \tau_b ) \) with the operation \( ( x , y ) \mapsto x + y \) defined by \( ( x + y ) ( \alpha ) \coloneqq x ( \alpha ) + y ( \alpha ) \) modulo \( 2 \).
For a more direct proof, let \( \seqofLR{\lambda_n}{n \in \omega} \) be a strictly increasing sequence of ordinals cofinal in \( \kappa \), and equip \( X \coloneqq \prod_{n \in \omega} \pre{\lambda_n}{2} \) with the product of the discrete topologies on each \( \pre{\lambda_n}{2} \): then the metric \( d \) on \( X \) defined by \( d ( x , y ) = 0 \) if \( x = y \) and \( d ( x , y ) = 2^{ - n } \) with \( n \) smallest such that \( x ( n ) \neq y ( n ) \) if \( x \neq y \) is complete and compatible with the topology of \( X \). 
Therefore, \( X \) and all its closed subsets are completely metrizable. 
Since the map \( \pre{\kappa}{2} \to X\), \( x \mapsto \seqofLR{x \restriction \lambda_n }{ n \in \omega } \) is a well-defined homeomorphism between \( ( \pre{\kappa}{2}, \tau_b ) \) and a closed subset of \( X \), we get the desired result. 

\smallskip

\ref{prop:topologicalproperties-h}
Let \( \seqofLR{ U_\beta }{ \beta < \alpha} \) be a sequence sets in \( \tau_ \lambda \), with \( \alpha < \lambda \).
It is enough to show that for each \( x \in V \coloneqq \bigcap_{ \beta < \alpha } U_ \beta \) there is \( s \in \forcing{Fn} ( \kappa , 2 ; \lambda ) \) such that \( x \in \Nbhd_s \subseteq V \).
For each \( \beta < \alpha \) choose \( u_ \beta \in [ \kappa ]^{< \lambda } \) such that \( \Nbhd_{ x \restriction u_ \beta } \subseteq U_ \beta \).
Then \( v \coloneqq \bigcup_{ \beta < \alpha } u_ \beta \in [ \kappa ]^{< \lambda } \) and \( \Nbhd_{x \restriction v } \subseteq V \). 
On the other hand, the set \( \Nbhd_s \) with \( s \in \pre{ \lambda }{ 2 } \) witnesses that \( \tau_ \lambda \) is not closed under intersections of length \( \lambda \), and similarly \( \tau _p \) is not closed under infinite intersections.

The argument for \( \tau _b \) is similar.
Suppose first \( \card{ \alpha } < \cof ( \kappa ) \) and let \( \seqofLR{ U_\beta }{ \beta < \alpha} \) be a sequence of sets in \( \tau_b \), and let \( V \) be as before. 
For any \( x \in \bigcap_{\beta < \alpha} U_\beta \) we construct an \( s \in \pre{< \kappa }{ 2} \) such that \( \Nbhd _s \subseteq V \).
For each \( \beta < \alpha \) let \( f ( \beta ) \) be the smallest \( \gamma < \kappa \) such that \( \Nbhd_{ x \restriction \gamma} \subseteq U_ \beta \).
By case assumption \( \ran ( f ) \) is bounded in \( \kappa \), so \( \Nbhd_s \subseteq V \) with \( s \coloneqq x \restriction \sup \ran ( f ) \).
Suppose now \( \card{\alpha } \geq \cof ( \kappa ) \).
Let \( \seqofLR{\gamma_\beta }{ \beta < \cof ( \kappa ) } \) be increasing and cofinal in \( \kappa \), and let \( V_ \beta \coloneqq \Nbhd_{0^{( \gamma_\beta)}} \) if \( \beta < \cof ( \kappa ) \) and \( V_ \beta \coloneqq \pre{ \kappa }{2} \) if \( \cof ( \kappa ) \leq \beta < \alpha \).
Then the set \( \bigcap_{ \beta < \alpha} V_ \beta = \bigcap_{\beta < \cof ( \kappa )} V_ \beta = \set{ 0^{( \kappa )} } \) is closed and not open. 

\smallskip

\ref{prop:topologicalproperties-i} follows from~\ref{prop:topologicalproperties-h}.
\end{proof}

\begin{remarks} \label{rem:cantorbasicproperties}
\begin{enumerate-(i)}
\item \label{rem:cantorbasicproperties-i}
Even if \( ( \pre{ \kappa }{2} , \tau_b ) \) is never compact by~\ref{prop:topologicalproperties-b} above, some form of compactness is available in certain cases: for example, as shown in~\cite[Theorem 5.6]{Motto-Ros:2011qc}, \( ( \pre{\kappa}{2}, \tau_b) \) is \( \kappa \)-compact%
\footnote{In general topology, \( \kappa \)-compact spaces are also called \( \kappa \)-Lindel\"of.}
(i.e.~such that every \( \tau_b \)-open covering of \( \pre{\kappa}{2} \) has a subcovering of size \( {<} \kappa \)) if and only if \( \kappa \) is a weakly compact cardinal. 
\item \label{rem:cantorbasicproperties-ii}
By~\ref{prop:topologicalproperties-c} and~\ref{prop:topologicalproperties-d} of Proposition~\ref{prop:topologicalproperties}, the space \( ( \pre{\kappa}{2}, \tau_p ) \) with \( \omega < \kappa \leq \card{\pre{\omega}{2}} \) witnesses (in \( \ZF \)) the well-known fact that separability does not imply second countability. 
Notice that the converse implication ``a second countable space is separable'' is equivalent to \( \AC_\omega \): given nonempty sets \( A_n \), endow \( X \coloneqq \bigcup_{n} A_n \) with the topology generated by the \( A_n \)'s so that \( X \) is second countable. 
Then any function enumerating a dense subset of \( X \) yields a choice function for the \( A_n \)'s.
(The other direction of the equivalence is immediate.)

\item \label{rem:cantorbasicproperties-ind}
Part~\ref{prop:topologicalproperties-e} of Proposition~\ref{prop:topologicalproperties} shows that the separability of \( ( \pre{\omega_1}{2}, \tau_p ) \) is independent of \( \ZF + \DC \).

\item \label{rem:cantorbasicproperties-iii}
The proof of part~\ref{prop:topologicalproperties-e} of Proposition~\ref{prop:topologicalproperties} also shows that if \( A \) is arbitrary, \( X \) is a separable space, and \( \pre{ A}{X} \) is endowed with the product topology,
\[
 \pre{ A}{X} \text{ is separable} \IFF \card{A} \leq \card{ \pre{\omega}{2} }. 
\]
Thus the Perfect set Property \( \PSP \) (which follows from \( \AD \)) implies that the space \( \pre{ A}{X} \) is separable if and only if \( \card{A} \leq \omega \) or \( \card{A} = \card{ \pre{\omega}{2} } \).
\item \label{rem:cantorbasicproperties-iv}
The techniques of classical descriptive set theory heavily rely on the existence of a (complete) metric, and hence by part~\ref{prop:topologicalproperties-g} of Proposition~\ref{prop:topologicalproperties} they cannot be directly applied to \( \pre{\kappa}{2} \) with \( \tau_p \), \( \tau _ \lambda \), and \( \tau_b \) when \( \kappa > \omega \).
In fact only a handful of ``positive'' results can be generalized to e.g.\ \( ( \pre{\kappa}{2}, \tau_b ) \) (see~\cite{Friedman:2011nx, Motto-Ros:2011qc}). 
\item
If \( \cof ( \kappa ) < \kappa \) then the collection of clopen sets \( \mathcal{C}_{\cof ( \kappa )} \) and \( \mathcal{C}_b \) are distinct \( \cof ( \kappa ) \)-algebras, since \( \mathcal{B}_b \subseteq \mathcal{C}_b \) and by the proof of Lemma~\ref{lem:boundedvscofinal}\ref{lem:boundedvscofinal-c} there are sets in \( \mathcal{B}_b \) that are not \( \tau _{\cof ( \kappa )} \)-open.
\end{enumerate-(i)}
\end{remarks}

\begin{proposition}\label{cor:boundedopennonproductBorel}
Let \( \lambda < \kappa \) be uncountable cardinals, and work in \( \pre{ \kappa }{2} \).
Assume \( \AComega ( \R ) \).
There is a set which is \( \tau _ \lambda \)-open and \( \tau _b \)-open, and it is not \( \tau _p \)-Borel.
\end{proposition}

Before starting the proof, let's fix a bit of notation: if \( \nu \leq \mu \) are infinite cardinals, the \markdef{inclusion map}\index[concepts]{inclusion map between generalized Cantor spaces} \( i \colon \pre{ \nu }{2} \hookrightarrow \pre{ \mu }{2} \) is the identity if \( \nu = \mu \), and it is the map \( x \mapsto x \conc 0^{ ( \mu )} \) if \( \nu < \mu \). 

\begin{proof}
By \( \AComega ( \R ) \) there is a non-Borel set \( A \subseteq \pre{\omega}{2} \).
Then \( U \coloneqq \bigcup_{s \in A} \Nbhd^\kappa_s \) is a \( \tau_ \lambda \)-open subset of \( \pre{\kappa}{2} \).
As the inclusion map \( ( \pre{\omega}{2} , \tau_p ) \hookrightarrow ( \pre{\kappa}{2} , \tau_p ) \) is continuous and since the generalized pointclass \( \bB \) is boldface, then \( U \) cannot be Borel in \( ( \pre{\kappa}{2}, \tau_p ) \).
\end{proof}

Using the arguments above, one can prove the next result describing the continuous functions between generalized Cantor spaces.

\begin{proposition}\label{lem:inclusionmap2}
Let \( \lambda \leq \cof ( \kappa ) \) and \( \lambda ' \leq \cof ( \kappa' ) \) with \( \kappa \leq \kappa ' \) cardinals, and let \( i \colon \pre{ \kappa }{2} \hookrightarrow \pre{ \kappa '}{2} \) be the inclusion map.
Then 
\begin{enumerate-(a)}
\item\label{lem:inclusionmap2-a}
\( i \colon ( \pre{ \kappa }{2} , \tau_ \lambda ) \hookrightarrow ( \pre{ \kappa' }{2} , \tau_ {\lambda'} ) \) is continuous if and only if \( \lambda' \leq \lambda \).
Therefore \( i \) is continuous when \( \pre{ \kappa '}{2} \) is topologized with \( \tau_p \) and \( \pre{ \kappa }{2} \) is topologized with any of \( \tau_p , \tau_ \lambda , \tau_b \).
\item
If \( \pre{ \kappa '}{2} \) is given the bounded topology, then \( i \) is continuous if and only if \( \kappa = \kappa ' \).
\item
\( i \colon ( \pre{ \kappa }{2} , \tau_ b ) \hookrightarrow ( \pre{ \kappa' }{2} , \tau_ {\lambda'} ) \) is continuous if and only if \( \lambda ' \leq \cof ( \kappa ) \).
\end{enumerate-(a)}
\end{proposition}

\section{Generalized Borel sets} \label{sec:borelsets}
\subsection{Basic facts} \label{subsec:borelsetsbasic}
\subsubsection{\( \alpha \)-Borel sets}

The following definition introduces a natural generalization of the notion of Borel subset of a topological space.
It plays a central role in the analysis of models of \( \AD \), but it is also of primary interest in other areas of set theory and in general topology.

\begin{definition} \label{def:kappaBorel}
If \( X = ( X , \tau ) \) is a topological space and \( \alpha \) an ordinal, the collection of \markdef{\( \alpha \)-Borel sets (with respect to \( \tau \))}\index[concepts]{Borel! sets \( \bB \)}\index[concepts]{Borel!\( \alpha \)-Borel sets \( \bB_ \alpha \)} is
\[
\bB_\alpha ( X , \tau ) \coloneqq \Alg ( \tau , \alpha ) , \index[symbols]{Balpha@\( \bB_\alpha \), \( \bB_\infty \)}
\]
the smallest family of subsets of \( X \) containing all \( \tau \)-open sets and closed under the operations of complementation and well-ordered unions of length \( {<} \alpha \). 
A set is \markdef{\( \infty \)-Borel (with respect to \( \tau \))}\index[concepts]{Borel!\( \infty \)-Borel sets \( \bB_\infty \)} if it is \( \alpha \)-Borel (with respect to \( \tau \)) for some \( \alpha \), i.e.~if it belongs to \( \bB_\infty ( X, \tau ) \coloneqq \bigcup_{\alpha \in \On} \bB_\alpha ( X , \tau ) \). 
As usual when drop the reference to \( X \) and/or \( \tau \) when there is no danger of confusion.
\end{definition}

Using the stratification of \( \Alg ( \tau,\alpha) \) in terms of the subcollections \( \bSigma_\gamma ( \tau,\alpha ) \) described in Section~\ref{subsec:algebras}, it is easy to check that that the operation \( \bB_\alpha \) assigning to each nonempty topological space \( ( X , \tau ) \) the collection \( \bB_\alpha ( X , \tau ) \) of its \( \alpha \)-Borel subsets is a hereditary general boldface pointclass.

\begin{remarks}\label{rem:borel}
\begin{enumerate-(i)}
\item \label{rem:borel-0}
If \( \alpha \leq \alpha' \) and \( \tau \subseteq \tau' \) then \( \bB_\alpha ( X , \tau ) \subseteq \bB_{ \alpha ' } ( X , \tau' ) \).
\item\label{rem:borel-1}
Both expressions \( \bB_{ \alpha + 1} ( X, \tau) \) and \( \bB_{ \alpha ^+ } ( X, \tau ) \) denote the same collection of subsets of \( X \), but the former is often preferable when discussing models of set theory, since the term \( \alpha + 1 \) is absolute, while \( \alpha^+ \) is not.
\item\label{rem:borel-2}
\( \bB_{ \omega _1} ( X , \tau ) = \bB_{ \omega + 1} ( X , \tau) \) is the usual collection of Borel subsets of \( X = ( X , \tau ) \).
\end{enumerate-(i)}
\end{remarks}

The notion of \( \alpha \)-Borel set can be trivial under choice.
For example any subset of an Hausdorff space \( X \) is in \( \bB_ \alpha ( X ) \) if \( \alpha > \card{X} \), and therefore \( \AC \) implies that \( \bB_\infty ( X ) = \pow ( X ) \). 
On the other hand, Theorem~\ref{th:B_kappa+1nontrivial} below shows that \( \bB_{ \kappa + 1} ( \pre{ \kappa }{2} , \tau_b ) \neq \pow ( \pre{ \kappa }{2} ) \) holds in \( \ZFC \), i.e.\ that \( \bB_{ \kappa + 1} ( \pre{ \kappa }{2} , \tau_b ) \) is never a trivial class.

The situation in the \( \AD \)-world is more subtle: \( \bB_ \infty ( \pre{\omega}{2} ) = \pow ( \pre{\omega}{2} ) \) holds in every known model of \( \AD \), and in fact the received opinion is that this must always be the case --- this is one-half of the well-known conjecture that \( \AD \implies \AD^+ \) (see Section~\ref{sec:SouslinunderAD}).
Recall from Section~\ref{subsubsec:pwoscaleproperty} that the length of a prewellordering \( \preceq \) is the length of the associated regular norm; this ordinal is denoted by \( \lh ( \preceq ) \).

\begin{definition}\label{def:Theta}
\( \Theta \coloneqq \sup \setof{ \alpha \in \On }{ \text{there exists a surjection } f \colon \R \twoheadrightarrow \alpha } \).\index[symbols]{Theta@\( \Theta \)} 
Equivalently
\[ 
\Theta \coloneqq \sup \setofLR{\lh ( \preceq ) }{\preceq \text{ is a prewellordering of } \R },
\]
and \( \R \) can be replaced by any uncountable Polish space such as \( \pre{\omega}{2} \) or \( \pre{\omega}{\omega} \).
\end{definition}

It is an easy exercise to show (in \( \ZF \)) that \( \Theta \) is a cardinal, while using choice \( \Theta = ( 2^{\aleph_0})^+ \). 
In contrast, the following theorem of H.~Friedman shows that \( \Theta \) is a limit cardinal under \( \AD \) (see~\cite[Theorem 28.15]{Kanamori:2003fu} or~\cite[Exercise 7D.19]{Moschovakis:2009fk} for a proof).

\begin{theorem}[\( \AD \)]\label{th:coding}
If \( \lambda < \Theta \), then \( \R \onto \pow ( \lambda ) \). 
In particular, \( \lambda ^+ < \Theta \) for all \( \lambda < \Theta \).
\end{theorem}

In fact, \( \AD \) implies that \( \Theta \) is quite large (e.g.~larger than the first fixed point of the \( \aleph \)-sequence), it has lot of measurable cardinals below it, and so on.
For a proof of these results see e.g.~\cite{Moschovakis:2009fk} and the references contained therein. 

Under \( \AD \) it is no longer true that successor cardinals are regular --- see Section~\ref{subsec:projectiveordinals}.
The next result shows that the regular cardinals are cofinal in \( \Theta \). 

\begin{lemma}\label{lem:arbitrarilylargeregularsbelowTheta}
Assume \( \AD \).
For all \( \alpha < \Theta \) there is a regular cardinal \( \alpha \leq \lambda < \Theta \).
\end{lemma}

\begin{proof}[Sketch of the proof]
Let \( \Gamma \coloneqq \pos \Sigma^1_1 ( \preceq ) \) be defined as in~\cite[Section~7C]{Moschovakis:2009fk}, where \( \preceq \) is a prewellordering of \( \pre{\omega}{\omega} \) of length \( \alpha \).
Then \( \Gamma \) is \( \omega \)-parametrized, so \( \pre{\omega}{\omega} \onto \bGamma \) where \( \bGamma \) is the boldface version of \( \Gamma \) as defined in~\cite[Section~3H]{Moschovakis:2009fk}.
Since \( \Gamma \) satisfies the hypotheses of~\cite[Theorem 7D.8]{Moschovakis:2009fk}, then
\[ 
\lambda \coloneqq \sup \setofLR{\beta \in \On}{\beta\text{ is the length of a strict well-founded relation in } \bGamma \text{ with field in }\pre{\omega}{\omega} }
 \] 
is a regular cardinal. 
Notice that \( \lambda < \Theta \), and that \( \alpha \leq \lambda \) by our choice of \( \preceq \).
\end{proof}

\begin{proposition}[\( \AD \)] \label{prop:examples2}
If \( \kappa < \Theta \) is an infinite cardinal, then \( \R \onto \bB_{\kappa + 1} ( \pre{\omega}{2}) \), and hence \( \bB_{\kappa + 1} ( \pre{\omega}{2}) \neq \pow ( \pre{\omega}{2} ) \).
\end{proposition}

\begin{proof}
By Lemma~\ref{lem:arbitrarilylargeregularsbelowTheta} let \( \kappa < \lambda < \Theta \) be regular.
By Lemma~\ref{lem:algebra} \( \pre{< \lambda}{( \pre{\omega}{2} )} \onto \bB_{ \kappa + 1} ( \pre{\omega}{2} ) \), so it is enough to show that \( \pre{\omega}{2} \onto \pre{< \lambda}{( \pre{\omega}{2} )} \).
This follows from \( \pre{ < \lambda}{( \pre{\omega}{2} )} \into \pre{ \lambda \times \omega }{2} \into \pre{ \lambda }{2} \) and Theorem~\ref{th:coding}. 
\end{proof}

\subsubsection{Borel codes}\label{subsubsec:Borelcodes}
Let \( X = ( X , \tau ) \) be a topological space and let \( \mathcal{S} \) be a basis for \( \tau \).
An \markdef{\( \alpha \)-Borel code}\index[concepts]{Borel!\( \alpha \)-Borel code} is a pair \( C = ( T, \phi ) \) where \( T \) is a well-founded descriptive set-theoretic tree (of height \( \leq \omega \)) on some \( \beta < \alpha \) and \( \phi \colon T \to \pow ( X ) \) is a map such that
\begin{subequations}\label{eq:definitionofcode}
\begin{align}
\phi ( t ) & \in \tau \text{, if \( t \) is a terminal node of } T , \label{eq:definitionofcode-1}
\\
\phi ( t ) & = X \setminus \bigcap \setofLR{\phi ( s) }{ s \text{ is an immediate successor of \( s \) in }T } \text{, otherwise.} \label{eq:definitionofcode-2}
\end{align}
The code \( C \) canonically determines a set \( \phi ( \emptyset ) \in \bB_\alpha ( X , \tau ) \), called the \markdef{\( \alpha \)-Borel set coded by \( C \)}.
A set coded by some \( C \) is called an \markdef{effective \( \alpha \)-Borel set}\index[concepts]{Borel!effective \( \alpha \)-Borel sets \( \bB^{\mathrm{e}} _ \alpha \)} , and
\[
\bB^{\mathrm{e}} _ \alpha ( X , \tau ) \index[symbols]{Balphae@\( \bB^{\mathrm{e}} _ \alpha \)}
\]
is the collection of all these sets.
If condition~\eqref{eq:definitionofcode-1} is strengthened to 
\begin{equation}\label{eq:definitionofcode-3}
 \phi ( t ) \in \mathcal{S} \text{, if \( t \) is a terminal node of } T ,
\end{equation}
\end{subequations}
we say that the code \( ( T , \phi ) \) takes values in \( \mathcal{S} \), and write
\[
\bB^{\mathrm{e}}_ \alpha ( X , \mathcal{S} ) 
\]
for the collection of sets admitting such a code.
(As usual, we drop \( X \), \( \tau \) or \( \mathcal{S} \) from our notation when they are clear from the context.)

\begin{remarks}
\begin{enumerate-(i)}
\item
The relation \( x \in \phi ( C ) \) is absolute for transitive models of enough set theory containing \( x \) and \( C \), thus ``effective Borelness'' is a more robust notion than that of ``Borelness''. 
Clearly
\[
\bB^{\mathrm{e}} _ \alpha ( X , \tau ) \subseteq \bB _ \alpha ( X , \tau ) , 
\]
and the converse inclusion holds under \( \AC \).
On the other hand, if choice fails badly the two notions can be quite different --- \( \bB^{\mathrm{e}} _ \alpha ( \R ) \) is always a surjective image of \( \R \) if \( \alpha \) is countable,%
\footnote{Every effective \( \omega + 1 \)-Borel code for a subset of a second countable topological space is obtained from a tree on \( \omega \), and hence, using the bijection between \( \pre{< \omega}{\omega} \) and \( \omega \), its characteristic function can be identified with an element of \( \pre{\omega}{2} \).} 
but in the Feferman-L\'evy model (where \( \R \) is a countable union of countable sets) it is true that \( \bB _{ \omega + 1 } ( \R ) = \pow ( \R ) \).
\item
It is immediate that \( \bB^{\mathrm{e}}_ \alpha ( X , \mathcal{S} ) \subseteq \bB^{\mathrm{e}}_ \alpha ( X , \tau ) \), but the reverse inclusion need not hold, even if \( \AC \) is assumed.
The problem is that \( \bB^{\mathrm{e}}_ \alpha ( X , \mathcal{S} ) \subseteq \Alg ( \mathcal{S} , \alpha ) \) (and in fact equality holds under \( \AC \)), and in order to show that \( \tau \subseteq \Alg ( \mathcal{S} , \alpha ) \), one needs some further assumptions on the size of a basis of \( \tau \).
For example, if \( ( X , \tau ) \coloneqq ( \pre{ \kappa }{2} , \tau_p ) \) and \( \mathcal{S} \coloneqq \mathcal{B}_p \), then \( \ZF \) proves \( \card{\mathcal{B}_p} = \kappa \) and
\begin{equation*}
\bB ^{\mathrm{e}}_{ \kappa + 1} (\mathcal{B}_p ) = \bB^{\mathrm{e}}_{ \kappa + 1} ( \tau_p ) . 
\end{equation*}

On the other hand, if \( ( X , \tau ) \coloneqq ( \pre{ \kappa }{2} , \tau_ b ) \) and \( \mathcal{S} \coloneqq \mathcal{B}_ b \), then working in \( \ZFC + 2^{(2^{< \kappa}) } > 2^\kappa \) we have that \( \card{\tau_b} = 2^{(2^{< \kappa}) } \) by Lemma~\ref{lem:cardinalityoftopologies}\ref{lem:cardinalityoftopologies-b} and \( \card{\Alg ( \mathcal{B} _ b , \kappa + 1 )} = 2^\kappa\) by \( \card{\mathcal{B}_b} = 2^{< \kappa} \) and Lemma~\ref{lem:algebra}, whence \( \tau_ b \nsubseteq \Alg ( \mathcal{B} _ b , \kappa + 1 ) \) and 
\begin{equation*}
\bB _{ \kappa + 1} ( \tau _ b ) = \bB^{\mathrm{e}} _{ \kappa + 1} ( \tau _ b ) \neq \bB _{ \kappa + 1} ( \mathcal{B} _ b ) .
\end{equation*}
Similarly, if \( ( X , \tau ) \coloneqq ( \pre{ \kappa }{2} , \tau_ \lambda ) \) with \( \omega < \lambda \leq \cof ( \kappa ) \) and \( \mathcal{S} \coloneqq \mathcal{B}_ \lambda \), then assuming \( \ZFC + 2^{(2^{< \lambda}) } > 2^\kappa \) we have that \( \tau_ \lambda \nsubseteq \Alg ( \mathcal{B} _ \lambda , \kappa + 1 ) \) and thus \( \bB _{ \kappa + 1} ( \tau _ \lambda ) = \bB^{\mathrm{e}} _{ \kappa + 1} ( \tau _ \lambda ) \neq \bB _{ \kappa + 1} ( \mathcal{B} _ \lambda ) \).
\end{enumerate-(i)}
\end{remarks}

The notion of \( \alpha + 1 \)-Borel set can be extended to the case when the ordinal \( \alpha \) is replaced by an arbitrary set.

\begin{definition}\label{def:JBorel}
Let \( X = ( X , \tau ) \) be a topological space and let \( J \) be an arbitrary set.
\begin{enumerate-(i)}
\item
\( \bB_{J + 1} ( X , \tau ) \) is the smallest collection of subsets of \( X \), containing all open sets, and closed under complements and unions of the form \( \bigcup_{j \in J} Y_j \) with \( \setofLR{Y_j}{j \in J } \subseteq \bB_{J + 1 } ( X ) \).
\item
A code \( ( T , \phi ) \) for \( A \in \bB_{J + 1} ( X , \tau ) \) is a well-founded descriptive set-theoretic tree \( T \) (of height \( \leq \omega \)) on \( J \) together with a map \( \phi \) with domain \( T \) satisfying~\eqref{eq:definitionofcode-1}--\eqref{eq:definitionofcode-2}, and such that \( \phi ( \emptyset ) = A \).
The family of all \( A \in \bB_{J + 1} ( X , \tau ) \) which admit a code is denoted with \( \bB_{J + 1}^{\mathrm{e}} ( X , \tau ) \), and if the codes satisfy~\eqref{eq:definitionofcode-3}, we write \( \bB_{J + 1}^{\mathrm{e}} ( X , \mathcal{S} ) \). 
\end{enumerate-(i)}
\end{definition}
The family \( \bB_{J + 1} ( X ) \) is closed under unions of the form \( \bigcup_{j' \in J'} Y_{j '} \) where \( J' \) is a surjective image of \( J \); the same applies to \( \bB_{J + 1}^{\mathrm{e}} ( X ) \) provided codes for the \( Y_{j '} \) are explicitly given.
Note that 
\[
\bB_{J + 1}^{\mathrm{e}} ( X , \mathcal{S} ) \subseteq \bB_{J + 1}^{\mathrm{e}} ( X , \tau ) \subseteq \bB_{J + 1} ( X , \tau ) .
\]

\begin{lemma} \label{lem:Borel}
Let \( J \coloneqq \pre{< \kappa}{2} \).
\begin{enumerate-(a)}
\item\label{lem:Borel-a}
\( \bB^{\mathrm{e}}_{J + 1} ( \pre{ \kappa }{2} , \tau_p ) = \bB^{\mathrm{e}}_{J + 1} ( \pre{ \kappa }{2} , \tau_b ) \) and \( \bB_{J + 1} ( \pre{ \kappa }{2} , \tau_p ) = \bB_{J + 1} ( \pre{ \kappa }{2} , \tau_b ) \).
\item\label{lem:Borel-b}
\( \bB^{\mathrm{e}}_{ \kappa + 1} ( \pre{ \kappa }{2} , \tau_b ) \subseteq \bB^{\mathrm{e}}_{ J + 1} ( \pre{ \kappa }{2} , \tau_p ) \) and \( \bB_{\kappa + 1} ( \pre{\kappa}{2}, \tau_b) \subseteq \bB_{ J + 1} ( \pre{\kappa}{2}, \tau_p) \). 
\end{enumerate-(a)}
\end{lemma}

\begin{proof}
\ref{lem:Borel-a}
Any \( \tau_b \)-basic open set \( \Nbhd_s \) is the intersection of a \( {<} \kappa \)-sequence of \( \tau_p \)-basic open sets, in particular it is effective \( \kappa + 1 \)-Borel (thus also effective \( \pre{< \kappa}{2} + 1 \)-Borel, as \( \pre{< \kappa}{2} \onto \kappa \)) with respect to \( \tau_p \).
If \( U \) is \( \tau_b \)-open, then \( U = \bigcup_{s \in S} \Nbhd_s \) where \( S \coloneqq \setofLR{ s \in \pre{< \kappa }{2}} { \Nbhd_s \subseteq U } \), and from \( S \) a code witnessing that \( U \) is effectively \( \pre{< \kappa }{2} + 1 \)-Borel with respect to \( \tau _p \) can be constructed, without any appeal to choice principles.
Therefore every (effective) \( \pre{< \kappa}{2} + 1 \)-Borel set in the topology \( \tau _b \) is (effective) \( \pre{< \kappa}{2} + 1 \)-Borel in the topology \( \tau _p \).
The other inclusion follows from \( \tau _p \subseteq \tau_b \).

\smallskip

\ref{lem:Borel-b} follows from~\ref{lem:Borel-a} and \( J = \pre{<\kappa}{2} \onto \kappa \).
\end{proof}

\begin{corollary}[\( \AC \)]\label{cor:borel}
If \( 2^{< \kappa} = \kappa \), then \( \bB_{\kappa + 1} ( \pre{\kappa}{2} , \tau_p) = \bB_{\kappa + 1} ( \pre{\kappa}{2}, \tau_b) \).
\end{corollary}

Notice that Corollary~\ref{cor:borel} applies to any cardinal (not just the regular ones).
In particular in models of \( \GCH \) the \( \kappa +1 \)-Borel sets of \( \pre{ \kappa }{2} \) with respect to the topologies \( \tau_p \), \( \tau_ \lambda \), and \( \tau_b \) are the same.
This should be contrasted with the results on \( \pre{ \kappa }{ \kappa } \) --- see Remark~\ref{rmks:Bairesummarize}\ref{rmks:Bairesummarize-ii}.

\subsection{Intermezzo: the projective ordinals} \label{subsec:projectiveordinals}
Proposition~\ref{prop:examples2} shows that under \( \AD \) each general boldface pointclass \( \bB_\alpha \) is proper, as long as \( \alpha < \Theta \).
If \( \alpha \) is an odd projective ordinal then \( \bB_ \alpha ( \pre{\omega}{2} ) \) can be pinned-down in the projective hierarchy. 

\begin{definition} \label{def:projectiveordinals}
Given any boldface pointclass \( \bGamma \), let\index[concepts]{projective ordinals}
\[ 
\bdelta_{\bGamma} \coloneqq \sup \setofLR{ \lh ( \preceq ) }{ \preceq \text{ is a prewellordering of \( \R \) in } \bDelta_{\bGamma} }.
\] 
For ease of notation set
\[ 
\bdelta^1_n \coloneqq \bdelta_{\bSigma^1_n} \quad \text{and} \quad \bdelta^2_1 \coloneqq \bdelta_{ \bSigma^2_1 } . \index[symbols]{d0elta@\( \bdelta_{\bGamma} \), \( \bdelta^1_n \), \( \bdelta^2_1 \)}
\] 
The \( \bdelta^1_n \) are called \markdef{projective ordinals}.
\end{definition}
In the definition above, the set \( \R \) can be replaced by any other uncountable Polish space, such as \( \pre{\omega}{2} \) or \( \pre{\omega}{\omega} \), and the pointclass \( \bGamma \) could be replaced with its dual \( \dual{\bGamma} \), i.e.~\( \bdelta_{\bGamma} = \bdelta_{\dual{\bGamma}} \).
In particular \( \bdelta^1_n = \bdelta_{\bPi^1_n} \) and \( \bdelta^2_1 = \bdelta_{\bPi^2_1} \).
Since the initial segments of the prewellordering induced by a \( \bGamma \)-norm are in \( \bDelta_{ \bGamma} \), we have the following
\begin{fact} \label{fct:newfact}
The length of a \( \bGamma \)-norm on a set \( A \in \bGamma ( \pre{\omega}{\omega} ) \) is \( \leq \bdelta_{\bGamma} \), and if the length is \( \bdelta_{\bGamma} \) then \( A \notin \bDelta_{\bGamma} ( \pre{\omega}{\omega} ) \).
\end{fact}
In other words: the ordinal \( \bdelta_{\bGamma} \) cannot be attained by a \( \bDelta_{\bGamma} \) prewellordering of a set in \( \bDelta_{\bGamma} \).
For example: from \( \ZF + \AComega ( \R ) \) it follows that \( \bdelta^1_1 = \omega _1 \), and hence every \( \bPi^{1}_{1} \)-norm on \( \pre{\omega}{2} \) (equivalently: every \( \bDelta^{1}_{1} \) prewellordering of \( \pre{\omega}{2} \)) has countable length, but there are \( \bPi^{1}_{1} \)-norms on (proper) \( \bPi^{1}_{1} \) sets of length \( \omega _1 \). 

The ordinal \( \bdelta_{\bGamma} \) carries a lot of information on the nature of the pointclass \( \bGamma \).
In particular, the projective ordinals can be used to describe the projective pointclasses.
For example Martin and Moschovakis proved that under \( \AD+\DC \)
\begin{equation} \label{eq:odd}
\bB_{\bdelta^1_{2 n + 1}} ( \pre{\omega}{2} ) = \bDelta^1_{2 n + 1}.
\end{equation}
By the late 70's it became clear that most of the structural problems on the projective pointclasses can be reduced to questions regarding the projective ordinals, and the computation (under \( \AD \)) of the \( \bdelta^1_n \) in terms of the \( \aleph \)-function came to the fore as a crucial problem.
Recall that \( \bdelta^1_1 = \omega _1 \).
 Martin proved that \( \AD \) implies that \( \bdelta^1_2 = \omega _2 \), \( \bdelta^1_3 = \aleph_{\omega + 1} \), and \( \bdelta^1_4 = \aleph_{\omega + 2} \).
By work of Kechris, Martin and Moschovakis (see~\cite{Kechris:1978vn}), assuming \( \AD \):
\begin{enumerate-(A)}
\item\label{en-A}
all \( \bdelta^1_n \)'s are regular cardinals,%
\footnote{Thus the \( \bdelta^1_n \) are also called \markdef{projective \emph{cardinals}}.}
 in fact measurable cardinals;
\item\label{en-B}
\( \bdelta^1_{2 n + 2} = ( \bdelta^1_{2 n + 1} )^+ \);
\item\label{en-C}
\( \bdelta^1_{2 n + 1} = ( \blambda_{2 n + 1} )^+ \), where \( \blambda^1_{2 n + 1} \)\index[symbols]{lambda@\( \blambda^1_{2 n + 1} \)} is a cardinal of countable cofinality.
\end{enumerate-(A)} 
Recall that \( \kappa \) is measurable if there is a \( \kappa \)-complete non-principal ultrafilter on \( \kappa \).
The theory \( \ZF \) proves that if \( \kappa \) is measurable, then it is regular (and therefore it is a cardinal), but under \( \AD \) measurable cardinals need not to be limit cardinals, and successor cardinals need not to be regular: for example \( \omega _1 \) and \( \omega _2 \) are measurable, and for \( n \geq 3 \) the \( \omega _n \)'s are singular of cofinality \( \omega _2 \).
Working in \( \ZFC + \AD^{\Ll ( \R )} \), no \( ( \omega _n )^{\Ll ( \R )} \) can be a cardinal for \( n \geq 3 \), neither can be \( ( \aleph_ \omega )^{ \Ll ( \R ) } \), and hence 
\begin{equation}\label{eq:delta13<=aleph3}
 \bdelta^1_3 = ( \aleph_{ \omega + 1 } )^{ \Ll ( \R ) } \leq \aleph_3 .
\end{equation}

By~\ref{en-B} and~\ref{en-C} the computation of the projective cardinals in models of determinacy boils down to determine the value of \( \blambda^1_{2 n + 1} \) or, equivalently, of \( \bdelta^1_{2 n + 1} \) for \( n \geq 2 \).
This was achieved by Jackson~\cite{Jackson:1989fk} who proved the general formula
\[
 \blambda^1_{2 n + 1} = \aleph_{ \gamma ( 2 n - 1)}
\]
where
\[
\gamma ( 1 ) = \omega \quad\text{and}\quad \gamma ( n + 1 ) = \omega ^{ \gamma ( n ) } . 
\]
En route to proving these results, Jackson verified that every regular cardinal smaller than \( \bdelta_\omega \coloneqq \sup_n \bdelta_n^1 = \aleph_{ \varepsilon _0} \) is measurable, and was able to compute exactly all these cardinals: between \( \bdelta_{2 n + 1}^1 \) and \( \bdelta_{2 n + 3}^1 \) there are \( 2^{n + 1 } - 1 \) regular (in fact: measurable) cardinals. 
From this it follows that \( \bdelta^1_{2 n + 1} \) is the \( ( 2^{n + 1} - 1 ) \)-st uncountable regular cardinal.
By these results, and arguing as for~\eqref{eq:delta13<=aleph3} we obtain:

\begin{corollary}\label{cor:delta1n<aleph}
Work in \( \ZFC \) and suppose there is an inner model containing all reals and satisfying \( \AD \).
Then \( \bdelta^1_{2 n + 1} \leq \aleph_{2^{n + 1} - 1} \).
\end{corollary}

\begin{remark}
Since \( \bdelta^1_n < \Theta \leq ( 2^{\aleph_0} ) ^+ \), \eqref{eq:delta13<=aleph3} and Corollary~\ref{cor:delta1n<aleph} become trivial if the continuum is small.
On the other hand, if the continuum is large, then by~\eqref{eq:odd} large cardinals imply that every projective set is in \( \bB_{ \aleph_ \omega } ( \pre{\omega}{2} ) \).
\end{remark}

\subsection{* More on generalized Borel sets}
\subsubsection{The (generalized) Borel hierarchy}\label{subsubsec:borelhierarchy}
The collection \( \bB ( X ) \) of all Borel subsets of a topological space \( ( X , \tau ) \) is stratified in a hierarchy by \( \bSigma^0_ \alpha ( X ) \) and \( \bPi^0_ \alpha ( X ) \).
A similar result holds for generalized Borel sets: in fact \( \bB _{ \kappa + 1 } ( X ) \) is the smallest \( \kappa +1 \)-algebra on \( X \) containing \( \tau \), so with the notation of Section~\ref{subsec:algebras}, \( \bB _{ \kappa + 1 } ( X ) = \bigcup_{ \alpha \in \On } \bSigma_{\alpha } ( \tau , \kappa + 1 ) = \bigcup_{ \alpha \in \On } \bPi_{\alpha } ( \tau , \kappa + 1 ) \).
For the sake of uniformity of notation, when \( X = \pre{ \kappa }{2} \) with the bounded topology, we write \( \bSigma^{0}_{ \alpha } ( \tau_b ) \) and \( \bPi^{0}_{ \alpha } ( \tau_b ) \) instead of \( \bSigma_{\alpha } ( \tau_b , \kappa + 1 ) \) and \( \bPi_{ \alpha } ( \tau _b , \kappa + 1 ) \).

A standard result in classical descriptive set theory is that if \( X \) is an uncountable Polish space, then the Borel hierarchy \( \seqofLR{\bSigma^0_\alpha ( X ) }{1 \leq \alpha < \omega_1} \) does not collapse,%
\footnote{At least assuming \( \AComega ( \R ) \) --- see~\cite{Miller:2011fk}.}
i.e.~\( \bSigma^{0}_{ \alpha } ( X ) \subset \bSigma^{0}_{ \beta } ( X ) \) for \(1 \leq \alpha < \beta < \omega _1 \).
This follows from the fact that each \( \bSigma^{0}_{ \alpha } ( \pre{\omega}{\omega} ) \) has a universal set (Section~\ref{subsubsec:pointclasses}).
Generalizing this notion, given a boldface pointclass \( \bGamma \) and a topological space \( X \), a set \( \mathcal{U} \in \bGamma ( X \times X ) \) is \markdef{universal for \( \bGamma ( X ) \)} if \( \bGamma ( X ) = \setof{ \mathcal{U}^{( y )}}{y \in X } \) where \( \mathcal{U}^{( y )} \coloneqq \setofLR{x \in X}{ ( x , y ) \in \mathcal{U} } \).
Corollary~\ref{cor:nouniversalset} below shows that for \( \pre{ \kappa }{2} \) with the bounded topology, there may be no universal sets for \( \bSigma^0_{ \alpha } ( \tau_b ) \), but nevertheless the hierarchy does not collapse (Proposition~\ref{prop:hierarchydoesnotcollapse}).

\begin{lemma}\label{lem:hierarchydoesnotcollapse0}
In the space \( ( \pre{ \kappa }{2} , \tau_b ) \), for every \( 1 \leq \alpha < \beta \in \On \)
\[
 \bSigma^0_{ \alpha } ( \tau_b ) \cup \bPi^0_{ \alpha } ( \tau_b ) \subseteq \bSigma^0_{ \beta } ( \tau_b ) \cap \bPi^0_{ \beta } ( \tau_b ) .
 \] 
\end{lemma}

\begin{proof}
By Lemma~\ref{lem:inclusionhierarchy} it is enough to show that \( \tau_b = \bSigma^0_1 ( \tau_b ) \subseteq \bSigma^0_2 ( \tau _b ) \).
Given \( U \) open and \( \nu < \kappa \), the set \( D_\nu \coloneqq \bigcup \setof{\Nbhd_s}{ \Nbhd_s \subseteq U \wedge \lh ( s ) = \nu } \) is \( \tau_b \)-clopen, and since \( U = \bigcup_{ \nu < \kappa } D_ \nu \), then \( U \in \bSigma^0_2 ( \tau _b ) \).
\end{proof}

\begin{corollary}\label{cor:nouniversalset}
Assume \( \AC \) and that \( 2^ \kappa < 2^{(2^{< \kappa })} \).
Then for \( 1 \leq \alpha < \kappa ^+ \), neither \( \bSigma^{0}_{ \alpha } ( \tau_b ) \) nor \( \bPi^{0}_{ \alpha } ( \tau_b ) \) have a universal set.
\end{corollary}

\begin{proof}
By Lemmas~\ref{lem:cardinalityoftopologies}\ref{lem:cardinalityoftopologies-b} and~\ref{lem:hierarchydoesnotcollapse0}, \( \bSigma^{0}_{ \alpha } ( \tau _b ) \) has size \( \geq 2^{(2^{< \kappa })} \), while \( \setofLR{\mathcal{U}^{( y )}}{y \in \pre{ \kappa }{2} } \) has size \( \leq 2^ \kappa \) for all \( \mathcal{U} \subseteq \pre{ \kappa }{2} \times \pre{ \kappa }{2} \).
The case of \( \bPi^{0}_{ \alpha } ( \tau_b ) \) follows by taking complements.
\end{proof}

Note that Corollary~\ref{cor:nouniversalset} applies e.g.\ when \( \kappa = \omega _1 \) and \( \MA_{ \omega _1} \) holds.

Although there may be no universal sets for \( \bSigma^0_\alpha(\tau_b) \) and \( \bPi^0_\alpha(\tau_b) \), we always have \( \leql^\kappa \)-complete sets for such classes. 
(We use here the notation and terminology from Section~\ref{subsubsec:lipschitz}.)
For \( x , y \in \pre{ \kappa }{2} \), let \( x \oplus y \in \pre{ \kappa }{2} \) be defined by
\[
 ( x \oplus y ) ( \alpha ) \coloneqq 
	\begin{cases}
	 x ( \beta ) & \text{if } \alpha = 2 \beta 
 	\\
	 y ( \beta ) & \text{if } \alpha = 2 \beta + 1 .
 	\end{cases}
\]
For \( A \subseteq \pre{ \kappa }{2} \) and \( x \in \pre{ \kappa }{2} \), let 
\[
x \oplus A \coloneqq \setofLR{ x \oplus a }{ a \in A } \subseteq \pre{ \kappa }{2} .
\]
Using games, it is easy to check that \( A \leql^\kappa x \oplus A \) and that if \( A \) has empty interior, then \( x \oplus A \leqW^\kappa A \).

\begin{lemma}[\( \AC \)]\label{lem:completesets}
Consider the space \( ( \pre{ \kappa }{2} , \tau_b ) \).
There is a sequence \( \seqofLR{A_ \alpha }{ 1 \leq \alpha < \kappa ^+ } \) such that each \( A_ \alpha \subseteq \pre{ \kappa }{2} \) is \( \leql^\kappa \)-complete for \( \bSigma^0_{\alpha } ( \tau_b ) \). 
\end{lemma}

\begin{proof}
The set \( A_1 \coloneqq \pre{ \kappa }{2} \setminus \set{ 0^{( \kappa )}} \) is open, and if \( U \) is open then \( \II \) wins \( \GL ( U , A_1 ) \) by playing \( 0 \)'s as long as \( \I \) has not reached a position \( s \in \pre{ < \kappa }{2} \) such that \( \Nbhd_s \subseteq U \) --- if this ever happens then \( \II \) plays a \( 1 \) and then plays an arbitrary sequence.

Suppose \( \alpha > 1 \).
Choose \( \seqof{ \alpha _\nu }{ \nu < \kappa } \) such that \( 1 \leq \alpha _\nu < \alpha \), \( \alpha = \sup_{\nu < \kappa } ( \alpha_\nu + 1 ) \), and
\begin{equation*}\label{eq:alpha_nurepeated}
 \card{\setof{ \nu < \kappa }{ \alpha _\nu = \gamma }} \in \set{ 0 , \kappa }
\end{equation*}
for all \( \gamma < \alpha \).
Let \( p_ \nu \colon \pre{ \kappa }{2} \onto \pre{ \kappa }{2} \) be the map 
\( p_\nu ( x ) ( \beta ) \coloneqq x ( \op{ \nu }{ \beta } ) \), where \( \op{\cdot }{\cdot } \) is the pairing map of~\eqref{eq:Hessenberg}.
Let
\[
A_ \alpha \coloneqq \setofLR{x \in \pre{ \kappa }{2} }{ \EXISTS{\nu < \kappa } ( p_ \nu ( x ) \in \pre{ \kappa }{2} \setminus A_ { \alpha _\nu } ) } .
\]
As each \( p_ \nu \) is continuous, \( A_ \alpha \in \bSigma_{\alpha }^0 ( \tau_b ) \), so it is enough to show that \( B \leql^\kappa A_\alpha \) for all \( B \in \bSigma^0_{\alpha } ( \tau_b ) \).
Suppose \( B = \bigcup_{ \xi < \kappa } B_ \xi \) with \( B_ \xi \in \bPi^0_{ \beta _ \xi } ( \tau_b ) \) and \( \beta _ \xi < \alpha \).
Let \( g \colon \kappa \into \kappa \) be such that \( \beta _ \xi \leq \alpha _{ g ( \xi )} \) --- such a \( g \) exists by our choice of the \( \alpha_\nu \)'s.
Construct a sequence of Lipschitz functions \( \seqofLR{f_ \nu }{ \nu < \kappa } \) as follows:
\begin{subequations}
\begin{align}
& \text{if \( \nu \in \ran g \) and \( g ( \xi ) = \nu \), then \( f_ \nu \) witnesses }B_ \xi \leql^\kappa \pre{ \kappa }{2} \setminus A_{ \alpha _ \nu } ;\label{eq:f_nu1}
\\
&\text{if \( \nu \notin \ran g \), then \( f_\nu \) is constant and takes value in } A_{ \alpha _\nu} . \label{eq:f_nu2}
\end{align}
\end{subequations}
The specific bijection \( \op{\cdot}{\cdot} \) maps \( \kappa \times \kappa \) onto \( \kappa \) and guarantees that the function
\[
f \colon \pre{ \kappa }{2} \to \pre{ \kappa }{2}, \qquad f ( x ) ( \op{ \nu }{ \beta } ) \coloneqq ( f_\nu ( x ) )( \beta )
\]
is Lipschitz and since \( p_\nu ( f ( x ) ) ( \beta ) = f ( x ) ( \op { \nu }{ \beta } ) = f_\nu ( x ) ( \beta ) \), then
\[
 \FORALL{ \nu < \kappa } \left [ p_\nu ( f ( x ) ) = f_\nu ( x ) \right ] .
\]
If \( x \in B \) then \( x \in B_ \xi \) for some \( \xi < \kappa \), so \( f_\nu ( x ) \notin A_{ \alpha _\nu} \) by~\eqref{eq:f_nu1} where \( \nu \coloneqq g ( \xi ) \), and hence \( f ( x ) \in A_ \alpha \).
Conversely, if \( p_\nu ( f ( x ) ) \notin A_{ \alpha _ \nu } \) for some \( \nu < \kappa \), then \( \nu \in \ran g \) by~\eqref{eq:f_nu2}, thus \( \nu = g ( \xi ) \) for some \( \xi \), and hence \( x \in B_ \xi \) by the choice of \( f_\nu \) and therefore \( x \in B \).
This completes the proof that \( \FORALL{ x \in \pre{ \kappa }{2} } \left [ x \in B \iff f ( x ) \in A_ \alpha \right ] \).
 \end{proof}

Note that by Lemma~\ref{lem:hierarchydoesnotcollapse0}, if \( A_\alpha , A_\beta \) are as in Lemma~\ref{lem:completesets} then \( A_ \alpha \leql^\kappa A_ \beta \) for \( 1 \leq \alpha \leq \beta < \kappa^+ \). 

\begin{theorem}[\( \AC \)]\label{th:B_kappa+1nontrivial}
For every infinite cardinal \( \kappa \), \( \bB_{ \kappa + 1} ( \pre{ \kappa }{2} , \tau_ b ) \neq \pow ( \pre{ \kappa }{2} ) \). 
\end{theorem}

\begin{proof}
Let \( \seqofLR{r_ \alpha }{ \alpha < \kappa ^+ } \) be distinct elements of \( \pre{ \kappa }{2} \).
Let
\[
A \coloneqq \bigcup_{ \alpha < \kappa ^+} r_ \alpha \oplus A_ \alpha 
\]
where \( A_ \alpha \) is as in Lemma~\ref{lem:completesets}.
By Proposition~\ref{prop:hardness}\ref{prop:hardness-a} it is enough to check that \( A \) is \( \leql^\kappa \)-hard for \( \bB _{ \kappa + 1} ( \pre{ \kappa }{2} , \tau_b ) \).
If \( B \in \bB _{ \kappa + 1} ( \pre{ \kappa }{2} , \tau_b ) \) then \( B \in \bSigma^0_{ \alpha } ( \tau_b ) \) for some \( \alpha < \kappa ^+ \), and hence \( B \leql^\kappa A_ \alpha \) by Lemma~\ref{lem:completesets}.
Since \( A_ \alpha \leql^\kappa A \) we are done.
\end{proof}

\begin{proposition}[\( \AC \)] \label{prop:hierarchydoesnotcollapse}
Consider the space \( ( \pre{\kappa}{2}, \tau_b) \) and let \( 1 \leq \alpha < \beta < \kappa ^+ \). 
Then \( \bSigma^0_{ \alpha } ( \tau_b ) \neq \bPi^0_{ \alpha } ( \tau_b ) \) and \( \bSigma^0_{ \alpha } ( \tau_b ) \subset \bSigma^0_{ \beta } ( \tau_b ) \). 
\end{proposition}

\begin{proof}
For the first part, argue by contradiction and use Lemma~\ref{lem:completesets} and Proposition~\ref{prop:hardness}\ref{prop:hardness-b}. 
The second part follows from the first one and Lemma~\ref{lem:hierarchydoesnotcollapse0}. 
\end{proof}
This shows that, under \( \AC \), the \( \kappa+1 \)-Borel hierarchy on \( (\pre{\kappa}{2}, \tau_b) \) never collapses, independently of the choice of the cardinal \( \kappa \). 
However, note that \( \bSigma^0_{ 1 } ( \tau_b ) \neq \bPi^0_{ 1 } ( \tau_b ) \) already holds in \( \ZF \) since e.g.\ \( \pre{ \kappa }{2} \setminus \set{ 0^{( \kappa )} } \) is open but not closed.

As for the classical case \( \kappa = \omega \), using the fact that the classes \( \bSigma^0_\alpha(\tau_b) \) and \( \bPi^0_\alpha(\tau_b) \) are closed under continuous preimages, we get the following corollary, which is a strengthening of Proposition~\ref{prop:hardness}\ref{prop:hardness-b} for the special case \( \mathcal{S} \coloneqq \bB_{\kappa + 1 } ( \pre{\kappa}{2} , \tau_b ) \).

\begin{corollary}
There is no \markdef{\( \leqW^\kappa \)-complete set for \( \bB_{ \kappa + 1 } ( \pre{\kappa}{2} , \tau_b ) \)}, i.e.~there is no \( A \in \bB_{ \kappa + 1 } ( \pre{\kappa}{2} , \tau_b ) \) such that \( B \leqW^\kappa A \) for every \( B \in \bB_{ \kappa + 1 } ( \pre{\kappa}{2} , \tau_b ) \). 
\end{corollary}

\subsubsection{Cardinality of \( \bB_{ \kappa + 1} \).}

Cardinality considerations may be useful to show that the notion of \( \kappa+1 \)-Borelness is nontrivial.

\begin{proposition}[\( \AC \)]\label{prop:examples}
Let \( \kappa \) be an infinite cardinal.
\begin{enumerate-(a)}
\item\label{prop:examples-a}
\( \card{\bB_{\kappa + 1 } ( \pre{\omega}{2} )} = \min \set{ 2^\kappa , 2^{(2^{\aleph_0})}} \). 
Therefore \( 2^\kappa < 2^{(2^{\aleph_0})} \IFF \card{\bB_{\kappa + 1 } ( \pre{\omega}{2} )} < \card{\pow ( \pre{\omega}{2} )} \).
\item\label{prop:examples-b}
\( \card{\bB_{\kappa + 1} ( \pre{\kappa}{2}, \tau_p)} = 2^\kappa = \card{\tau_p} \). 
Therefore \( \card{\bB_{\kappa + 1} ( \pre{\kappa}{2}, \tau_p )} < 2^{(2^\kappa)} \), and hence, in particular, \( \bB_{\kappa + 1} ( \pre{\kappa}{2} , \tau_p ) \neq \pow ( \pre{\kappa}{2} ) \).
\item\label{prop:examples-c}
\( \card{\bB_{\kappa + 1} ( \pre{\kappa}{2}, \tau_b)} = 2^{(2^{<\kappa})} = \card{\tau_b} \). 
Therefore \( 2^{(2^{<\kappa})} < 2^{(2^{ \kappa })} \IFF \card{\bB_{\kappa + 1} ( \pre{\kappa}{2}, \tau_b)} < \card{ \pow (\pre{\kappa}{2})} \). 
\end{enumerate-(a)}
\end{proposition}

\begin{proof}
\ref{prop:examples-a} 
\( \card{\bB_{\kappa + 1 } ( \pre{\omega}{2} ) } \leq (2^{\aleph_0})^\kappa = 2^\kappa \) by Lemma~\ref{lem:algebra}\ref{lem:algebra-b}, and since \( \bB_{\kappa + 1 } ( \pre{\omega}{2} ) \subseteq \pow ( \pre{\omega}{2} ) \) then 
\[ 
\card{\bB_{\kappa + 1 } ( \pre{\omega}{2} ) } \leq \min \set{ 2^ \kappa , 2^{(2^{\aleph_0})} } .
\]
To prove the reverse inequality we take cases.
If \( \kappa \geq 2^{\aleph_0} \) then \( \bB_{\kappa + 1} ( \pre{\omega}{2}) = \pow ( \pre{\omega}{2} ) \), and we are done.
If \( \kappa < 2^{\aleph_0} \), let \( \seqof{ x_\alpha }{ \alpha < \kappa } \) be distinct elements of \( \pre{\omega}{2} \): then the map \( \pow ( \kappa ) \to \bB_{\kappa + 1} ( \pre{\omega}{2} ) \), \( A \mapsto \setof{x_\alpha}{\alpha \in A} \) is injective.

\smallskip

\ref{prop:examples-b}
By Lemmas~\ref{lem:cardinalityoftopologies}\ref{lem:cardinalityoftopologies-a} and~\ref{lem:algebra}, we get \( 2^\kappa = \card{\tau_p} \leq \card{\bB_{\kappa + 1} ( \pre{\kappa}{2}, \tau_p)} \leq ( 2^\kappa )^\kappa = 2^\kappa \).

\ref{prop:examples-c}
By Lemmas~\ref{lem:cardinalityoftopologies}\ref{lem:cardinalityoftopologies-b} and~\ref{lem:algebra}, and by \( \kappa \leq 2^{< \kappa} \), we get \( 2^{(2^{< \kappa})} = \card{\tau_b} \leq \card{\bB_{\kappa + 1} ( \pre{\kappa}{2}, \tau_b)} \leq (2^{(2^{<\kappa})})^\kappa = 2^{(2^{< \kappa})} \). 
\end{proof}

\begin{remark} \label{rmk:differentborel}
By parts~\ref{prop:examples-b} and~\ref{prop:examples-c} of Proposition~\ref{prop:examples}, if \( 2^\kappa < 2^{(2^{< \kappa})} \) (which follows from \( \MA_{\omega_1} \) when \( \kappa = \omega_1 \)) then \( \bB_{\kappa + 1} ( \pre{\kappa}{2}, \tau_p) \neq \bB_{\kappa + 1} ( \pre{\kappa}{2}, \tau_b) \). 
This should be contrasted with Corollary~\ref{cor:borel}.
\end{remark}
\section{Generalized Borel functions}\label{sec:generalizedBorelfunctions}
\subsection{Basic facts}\label{subsec:generalizedBorelfunctionsbasic}
\subsubsection{Borel measurable functions}\label{subsubsec:Borelmeasurablefunctions}
We now turn to \( \alpha \)-Borel measurable functions between generalized Cantor spaces. 
\emph{Unless otherwise explicitly stated} (see e.g.\ the definition of a weakly \( \mathcal{S} \)-measurable function below)\emph{, 
in this section such spaces are tacitly endowed with the bounded topology \( \tau_b \)}, so that all other derived topological notions (such as \( \alpha \)-Borelness, and so on) refer to such topology. 

\begin{definition}\label{def:weaklyBorel} 
Let \( \lambda , \mu \) be infinite cardinals, let \( f \colon \pre{ \lambda}{ 2} \to \pre{ \mu}{2} \), and let \( \mathcal{S} \subseteq \pow ( \pre{ \lambda}{ 2} ) \) be an algebra. 
\begin{itemize}[leftmargin=1pc]
\item
\( f \) is \markdef{\( \mathcal{S} \)-measurable} if and only if the preimage of every \( \tau_b \)-open subset of \( \pre{\mu}{2} \) is in \( \mathcal{S} \).
\item
\( f \) is \markdef{weakly \( \mathcal{S} \)-measurable} if and only if the preimage of every \( \tau_p \)-open subset of \( \pre{\mu}{2} \) is in \( \mathcal{S} \).
\item 
If \( \alpha \) is an infinite ordinal and \( \mathcal{S} = \bB_ \alpha ( \pre{ \lambda }{2} ) \), a (weakly) \( \mathcal{S} \)-measurable function is called a \markdef{(weakly) \( \alpha \)-Borel}\index[concepts]{measurable function}\index[concepts]{weakly \( \alpha \)-Borel function|see{measurable function}}\index[concepts]{Borel!function} function. 
Therefore \( f \) is an \(\alpha\)-Borel function if and only if \( f^{-1}(B) \in \bB_\alpha ( \pre{\lambda}{2} ) \) for every \( B \in\bB_\alpha(\pre{\mu}{2}) \).
\end{itemize}
\end{definition}

If in Definition~\ref{def:weaklyBorel} \( \alpha \)-Borel sets are replaced with their effective versions, we get the notion of \markdef{effective (weakly) \( \alpha \)-Borelness} for functions \( f \colon \pre{ \lambda}{2} \to \pre{ \mu}{2} \).
As for the case of sets, this notion coincides with its non-effective version under \( \AC \), but it may be a stronger notion in a choice-less world.
Every (effective) \( \alpha \)-Borel function is (effective) weakly \( \alpha \)-Borel because \( \tau_b \) refines \( \tau_p \). 
On the other hand, if \( f \colon \pre{\lambda}{2} \to \pre{\mu}{2} \) is weakly \( \alpha \)-Borel, then \( f \) might fail to be \( \alpha \)-Borel (see Proposition~\ref{prop:weaklyalphaBorel}\ref{prop:weaklyalphaBorel-b} below), but if \( \alpha \geq \mu \) then at least all preimages of \emph{\( \tau_b \)-basic open sets} of \( \pre{\mu}{2} \) are in \( \bB_\alpha(\pre{\lambda}{2}) \), that is
\[
 \FORALL{U \in \mathcal{B}_b ( \pre{ \mu} {2})} \left ( f^{-1} ( U ) \in \bB_\alpha(\pre{\lambda}{2}) \right ) .
\]
(Use the fact that every \( \tau_b \)-basic open set of \( \pre{ \mu}{2} \) can be written as an intersection of \( < \mu \)-many \( \tau_p \)-clopen sets). 
This also implies that when \( \alpha > \mu \), then the preimage under \( f \) of any set in the \( \mu+1 \)-algebra%
\footnote{In~\cite[Definition 15]{Friedman:2011nx} the authors call (\( \mu+1 \)-)Borel sets the elements of \( \bB_{\mu+1}(\mathcal{B}_b) \), rather than those in \( \bB_{\mu+1}(\tau_b) \), thus the functions \( f \colon \pre{\mu}{2} \to \pre{\mu}{2} \) which are weakly \( \mu+1 \)-measurable in our sense are still (generalized) Borel functions in their sense.
Note that Definition 15 of~\cite{Friedman:2011nx} coincides with our Definition~\ref{def:kappaBorel} when \( \mu^{< \mu} = \mu \), which is the setup of that paper.}
\( \bB_{\mu+1}(\mathcal{B}_b) \) on \( \pre{\mu}{2} \) generated by the canonical basis \( \mathcal{B}_b ( \pre{ \mu} {2}) \) belongs to \( \bB_\alpha( \pre{\lambda}{2} ) \).

Proposition~\ref{prop:weaklyalphaBorel}\ref{prop:weaklyalphaBorel-a} shows that under \( \AC \) the notions \( \alpha \)-Borel and weakly \( \alpha \)-Borel may also coincide if certain cardinal conditions are satisfied.

\begin{proposition}\label{prop:weaklyalphaBorel}
Let \( \lambda , \mu \) be cardinals, \( f \colon \pre{ \lambda }{2} \to \pre{ \mu}{2} \), and \( \alpha \geq \omega \).
\begin{enumerate-(a)}
\item\label{prop:weaklyalphaBorel-a}
Assume \( \AC \).
If \( 2^{< \mu} = \mu < \alpha \) and \( f \) is (effective) weakly \( \alpha \)-Borel, then \( f \) is (effective) \( \alpha \)-Borel.
In particular, if \( \mu = \omega \) (here we do not need \( \AC \)), or if \( \CH \) holds and \( \mu = \omega_1 \), then \( f \) is (effective) weakly \( \alpha \)-Borel if and only if \( f \) is (effective) \( \alpha \)-Borel.
\item\label{prop:weaklyalphaBorel-b}
If \( \lambda < \mu \) and \( \bB_\alpha ( \pre{ \lambda}{2} ) \neq \pow ( \pre{ \lambda}{2} ) \), then the inclusion map \( i \colon \pre{ \lambda }{2} \hookrightarrow \pre{ \mu }{2} \) is weakly \( \alpha \)-Borel but not \( \alpha \)-Borel.
\end{enumerate-(a)}
\end{proposition}

\begin{proof}
Part~\ref{prop:weaklyalphaBorel-a} easily follows from the fact that \( \mathcal{B}_b ( \pre{ \mu }{2} ) \) has size \( 2^{< \mu} \).

\smallskip

\ref{prop:weaklyalphaBorel-b} 
The inclusion function is not \( \alpha \)-Borel by an argument as in Proposition~\ref{cor:boundedopennonproductBorel}. 
By Lemma~\ref{lem:inclusionmap2}\ref{lem:inclusionmap2-a} we get that \( i \colon ( \pre{ \lambda}{2} , \tau_b) \hookrightarrow ( \pre{ \mu}{2}, \tau_p ) \) is continuous, and hence \( i \) is trivially weakly \(\alpha\)-Borel.
\end{proof}

Since the inclusion map of part~\ref{prop:weaklyalphaBorel-b} of Proposition~\ref{prop:weaklyalphaBorel} is arguably one of the simplest functions that one may want to consider, this suggests that weakly \( \alpha \)-Borelness is a more natural and appropriate notion of topological complexity for functions \( \pre{ \lambda}{2} \to \pre{ \mu}{2} \) when \( \lambda < \mu \).
On the other hand, it seems that in general there is no obstruction to the possibility of considering the stronger notion of \( \alpha \)-Borelness when \( \lambda \geq \mu \) (we come back to this point at the beginning of Section~\ref{subsec:topologicalcomplexity}).

As the name suggests, the notion of weakly \( \alpha \)-Borelness is quite weak.
In fact there are situations where it becomes vacuous, i.e.~every function is weakly \( \alpha \)-Borel. 
Proposition~\ref{prop:weaklyBorel} shows that the existence of a non-weakly \( \alpha \)-Borel \( f \colon \pre{ \lambda }{2} \to \pre{ \mu}{ 2} \) depends only on \( \lambda \) and \( \alpha \).

\subsubsection{$\bGamma$-in-the-codes functions} \label{subsubsec:Gammainthecodes}

When considering projective levels in models of \( \AD \), it is natural to code functions of the form \( f \colon \pre{ \omega}{2} \to \pre{ \kappa}{2} \) (for suitable \( \kappa < \Theta \)) as subsets of (products of) \( \pre{ \omega}{2} \) and \( \pre{ \omega}{\omega} \), and then require that such a code be in the projective pointclass under consideration. 
More precisely, we have the following general definition, which makes sense in arbitrary models of \( \ZF \). 
To simplify the presentation, we say that \(\rho\) is a \markdef{\( \bGamma \)-code for \( \kappa \)}\index[concepts]{Gamma@\( \bGamma \)-code for \( \kappa \)} if \( \rho \) is a \( \bGamma \)-norm of length \( \kappa \) on some \( A \in \bGamma ( \pre{ \omega}{\omega} ) \). 

\begin{definition} \label{def:Gammainthecodes}
Let \( \bGamma \) be a boldface pointclass, let \( \kappa \) be an infinite cardinal, and assume that there is a \( \bGamma \)-code \( \rho \) for \( \kappa \). 
We say that the function \( f \colon \pre{ \omega}{2} \to \pre{ \kappa}{2} \) is \markdef{\( \bGamma \)-in-the-codes}\index[concepts]{Gamma@\( \bGamma \)-in-the-codes function} (with respect to \( \rho \)) if there is \( F \in \bGamma( \pre{ \omega}{2} \times \pre{ \omega}{\omega} \times 2) \) such that for every \( \alpha < \kappa \) and \( i \in 2 \)
\[ 
f^{-1} ( \widetilde{\Nbhd}\vphantom{\Nbhd}^\kappa_{\alpha , i}) = \setofLR{ x \in \pre{ \omega}{2} }{ \EXISTS{ y \in A} ( \rho ( y ) = \alpha \wedge ( x , y , i ) \in F) },
 \] 
where \( \widetilde{\Nbhd}\vphantom{\Nbhd}^\kappa_{\alpha , i} \) is defined as in Remark~\ref{rmk:product}\ref{rmk:product-a}.
\end{definition}

\begin{remarks} \label{rmk:inthecodes}
\begin{enumerate-(i)}
\item\label{rmk:inthecodes-i}
As we shall see in Remark~\ref{rmk:rhononessential}\ref{rmk:rhononessential-b}, Definition~\ref{def:Gammainthecodes} does not really depend on the choice of \( \rho \) when \( \bGamma \) is sufficiently closed.
\item\label{rmk:inthecodes-ii}
If there is a \( \bGamma \)-code \(\rho\) for \( \kappa \), then \( \kappa \leq \bdelta_{\bGamma} \) (see Definition~\ref{def:projectiveordinals}).
Indeed, if \( \alpha < \kappa \) and \( z \in A \) is such that \( \rho ( z ) = \alpha \), then the prewellordering 
\[
 x \preceq y \IFF \bigl [ y \in A \wedge \rho ( y ) < \rho ( z ) \implies x \in A \wedge \rho ( x ) \leq \rho ( y ) \bigr ]
\]
is in \( \bDelta_{ \bGamma} \) and has length \( \alpha + 1 \). 
On the other hand, if \( \kappa < \bdelta_{\bGamma} \) then there is a \( \bGamma \)-code \(\rho\) for \( \kappa \): it is enough to let \( \rho \) be the \( \bGamma \)-norm associated to any \( \bDelta_{\bGamma} \) prewellordering of \( \pre{\omega}{\omega} \) of length \( \kappa \), which exists by \( \kappa < \bdelta_{\bGamma} \).
\item \label{rmk:inthecodes-iii}
Let \( \bGamma \subseteq \bGamma' \) be boldface pointclasses. 
If there is a \( \bGamma \)-code for \(\kappa\), then there is also a \( \bGamma' \)-code for any \( \kappa' \leq \kappa \).
\end{enumerate-(i)}
\end{remarks}

If \( \bGamma = \bSigma^1_1 \) in Definition~\ref{def:Gammainthecodes}, then \( \kappa = \bdelta^1_1 = \omega \) by Remark~\ref{rmk:inthecodes}\ref{rmk:inthecodes-ii} and~\cite[Theorem~4A.4]{Moschovakis:2009fk}, and therefore a function \( f \colon \pre{ \omega}{2} \to \pre{ \omega}{2} \) is \( \bSigma^1_1 \)-in-the-codes if and only if it is (\(\omega + 1 \)-)Borel measurable. 
More generally, if \( \bGamma \) is a nontrivial boldface pointclass closed under projections and countable intersections (e.g.\ if \( \bGamma = \bSigma^1_n \) for some \( n \geq 1 \)), then for all \( f \colon \pre{ \omega}{2} \to \pre{ \omega}{2} \)
\[
f \text{ is \( \bGamma \)-in-the-codes} \IFF f \text{ is \( \bGamma \)-measurable} \IFF f \text{ is \( \bDelta_{\bGamma} \)-measurable.}
\]
To see this, use the fact that under the hypotheses above we have that \( f \) is \( \bGamma \)-measurable if and only if the graph of \( f \) is in \( \bGamma ( \pre{ \omega}{2} \times \pre{ \omega}{2}) \).
This result can be partially generalized to uncountable \( \kappa \)'s. 
In fact, since \( \pre{ \omega}{2} \setminus f^{-1} ( \widetilde{\Nbhd}\vphantom{\Nbhd}^\kappa_{\alpha , i}) = f^{-1} ( \widetilde{\Nbhd}\vphantom{\Nbhd}^\kappa_{\alpha , 1-i}) \), if \( \bGamma \) is closed under projections and finite intersections then 
\begin{equation} \label{eq:inthecodes}
 f \colon \pre{ \omega}{2} \to \pre{ \kappa}{2} \text{ is \( \bGamma \)-in-the-codes} \IMPLIES \FORALL{\alpha < \kappa} \FORALL{i \in \setLR{0,1}} \bigl ( f^{-1} ( \widetilde{\Nbhd}\vphantom{\Nbhd}^\kappa_{\alpha , i}) \in \bDelta_{\bGamma} \bigr ) . 
\end{equation}
Lemma~\ref{lem:inthecodes} and Proposition~\ref{prop:inthecodes} below show that for certain \( \bGamma \)'s, the implication can be reversed.

\subsection{* Further results}
\begin{proposition}\label{prop:weaklyBorel}
For \( \alpha \) an infinite ordinal, the following are equivalent:
\begin{enumerate-(a)}
\item\label{prop:weaklyBorel-a}
\( \bB_\alpha ( \pre{ \lambda}{ 2} ) \neq \pow ( \pre{ \lambda}{ 2} ) \),
\item\label{prop:weaklyBorel-b}
for all cardinals \( \mu \) there is an \( f \colon \pre{ \lambda }{2} \to \pre{ \mu }{2} \) which is not weakly \( \alpha \)-Borel,
\item\label{prop:weaklyBorel-c}
there is a cardinal \( \mu \) and there is an \( f \colon \pre{ \lambda }{2} \to \pre{ \mu }{2} \) which is not weakly \( \alpha \)-Borel.
\end{enumerate-(a)}
\end{proposition}

\begin{proof}
If \( A \subseteq \pre{ \lambda}{ 2} \) is not in \( \bB_{\alpha} ( \pre{ \lambda}{ 2} ) \), then picking distinct \( y_0 , y_1 \in \pre{ \mu}{2} \) one easily sees that the function mapping the points in \( A \) to \( y_1 \) and the points in \( \pre{ \lambda}{ 2} \setminus A \) to \( y_0 \) is not weakly \( \alpha \)-Borel. 
This proves that \ref{prop:weaklyBorel-a}\( \implies \)\ref{prop:weaklyBorel-b}. 

Since \ref{prop:weaklyBorel-b}\( \implies \)\ref{prop:weaklyBorel-c} is trivial, it is enough to prove that~\ref{prop:weaklyBorel-a} follows from~\ref{prop:weaklyBorel-c}: but this is easy, as if \( f \colon \pre{ \lambda}{ 2} \to \pre{ \mu}{ 2} \) is not weakly \( \alpha \)-Borel then by definition \( f^{-1} ( U ) \notin \bB_{\alpha} ( \pre{ \lambda}{ 2} ) \) for some \( \tau_p \)-open \( U \subseteq \pre{\mu}{2} \).
\end{proof} 

We are mainly interested in (weakly) \( \kappa + 1 \)-Borel functions between the Cantor space \( \pre{ \omega}{2} \) and (a homeomorphic copy of) \( \pre{ \kappa}{2} \), with \( \kappa > \omega \) a cardinal.
From Propositions~\ref{prop:examples} and~\ref{prop:weaklyBorel} we obtain

\begin{corollary}\label{cor:prop:examples+weaklyBorel}
Assume \( \AC \).
\begin{enumerate-(a)}
\item\label{cor:prop:examples+weaklyBorel-a}
If \( \kappa \geq 2^{\aleph_0} \), then every \( f \colon \pre{ \omega}{2} \to \pre{ \kappa }{2} \) is \( \kappa + 1 \)-Borel.
\item\label{cor:prop:examples+weaklyBorel-b}
If \( 2^ \kappa < 2^{( 2^{\aleph_0})} \), then there is an \( f \colon \pre{ \omega}{2} \to \pre{ \kappa }{2} \) which is not weakly \( \kappa + 1 \)-Borel.
\item\label{cor:prop:examples+weaklyBorel-c}
If \( 2^{(2^{< \kappa})} < 2^{( 2^ \kappa) } \), then there is an \( f \colon \pre{ \kappa }{2} \to \pre{ \omega }{2} \) which is not weakly \( \kappa + 1 \)-Borel.
\end{enumerate-(a)}
\end{corollary}
Since \( \kappa+1 \)-Borel functions are in particular weakly \( \kappa+1 \)-Borel,
we could have removed the adjective ``weakly'' in parts~\ref{cor:prop:examples+weaklyBorel-b} and~\ref{cor:prop:examples+weaklyBorel-c} of Corollary~\ref{cor:prop:examples+weaklyBorel}.

From Propositions~\ref{prop:examples2} and~\ref{prop:weaklyBorel} and equation~\eqref{eq:odd} on page~\pageref{eq:odd} we obtain

\begin{corollary} \label{cor:5.6}
Assume \( \AD \).
\begin{enumerate-(a)}
\item \label{cor:5.6-a}
If \( \kappa < \Theta \), then there is an \( f \colon \pre{ \omega}{2} \to \pre{ \kappa}{2} \) which is not weakly \( \kappa + 1 \)-Borel.
\item \label{cor:5.6-b}
If \( \kappa = \blambda^1_{2 n + 1} \) as defined in condition~\ref{en-C} on page~\pageref{en-C}, then \( f \colon \pre{ \omega}{2} \to \pre{ \kappa}{ 2} \) is (weakly) \( \kappa + 1 \)-Borel if and only if it is (weakly) \( \bDelta^1_{2 n + 1} \)-measurable.
\end{enumerate-(a)}
\end{corollary}

We now turn our attention to \( \bGamma \)-in-the-codes functions.

\begin{lemma} \label{lem:inthecodes}
Let \( \bGamma \) be a nontrivial boldface pointclass closed under projections and finite intersections and unions. 
Let \( \kappa \) be an infinite cardinal, and suppose that \( \bGamma \) is closed under well-ordered unions of length \( \kappa \) and that there is a \( \bGamma \)-code \(\rho\) for \( \kappa \).
Then for every \( f \colon \pre{ \omega}{2} \to \pre{ \kappa}{2} \) the following are equivalent:
\begin{enumerate-(a)}
\item \label{lem:inthecodes-1}
\( f \) is \( \bGamma \)-in-the-codes (with respect to \( \rho \));
\item \label{lem:inthecodes-2}
\( f^{-1} ( \widetilde{\Nbhd}\vphantom{\Nbhd}^\kappa_{\alpha , i}) \in \bDelta_{\bGamma} \) for every \( \alpha < \kappa \) and \( i = 0 , 1 \);
\item \label{lem:inthecodes-3}
\( f^{-1} ( U ) \in \bDelta_{\bGamma} \) for every \( U \in \mathcal{B}_p ( \pre{ \kappa}{2} ) \).
\end{enumerate-(a)}
\end{lemma}

\begin{proof}
The equivalence between~\ref{lem:inthecodes-2} and~\ref{lem:inthecodes-3} is easy, while~\ref{lem:inthecodes-2} follows from~\ref{lem:inthecodes-1} by equation~\eqref{eq:inthecodes}.
So it is enough to show that~\ref{lem:inthecodes-2} implies~\ref{lem:inthecodes-1}.
If \( \rho \) is as in the hypotheses and \( A \in \bGamma(\pre{\omega}{\omega}) \) is its domain, it is enough to prove that the set \( \bigcup_{\alpha<\kappa , i \in 2} f^{-1} ( \widetilde{\Nbhd}\vphantom{\Nbhd}^\kappa_{\alpha , i}) \times \setofLR{ y \in A }{ \rho ( y ) = \alpha } \times \setLR{i} \) is in \( \bGamma \).
But this easily follows from the fact that by~\ref{lem:inthecodes-2} each \( f^{-1} ( \widetilde{\Nbhd}\vphantom{\Nbhd}^\kappa_{\alpha , i } ) \times \setof{ y \in A }{ \rho ( y ) = \alpha } \times \set{i } \) belongs to \( \bDelta_{\bGamma} \subseteq \bGamma \), together with the fact that \( \bGamma \) is closed under well-ordered unions of length \( \kappa \).
\end{proof}

\begin{lemma}[\( \AD + \DC \)]
Let \( \bGamma \) be a nontrivial boldface pointclass closed under projections, countable unions, and countable intersections.
Let \( \kappa \) be an infinite cardinal, and suppose that there is a \( \bGamma \)-code for \( \kappa \).
Then \( \bGamma \) is closed under well-ordered unions of length \( \kappa \).
\end{lemma}

\begin{proof}
Let \( \rho \) be a \( \bGamma \)-code for \( \kappa \) and \(A \in \bGamma ( \pre{\omega}{\omega} ) \) be its domain. 
Recall from Remark~\ref{rmk:inthecodes}\ref{rmk:inthecodes-ii} that the existence of such a \( \rho \) implies \( \kappa \leq \bdelta_{\bGamma} \).
If \( \kappa < \bdelta_{\bGamma} \) we can apply Moschovakis' Coding Lemma (see~\cite[Lemmas~7D.5 and~7D.6]{Moschovakis:2009fk}).
If \( \kappa = \bdelta_{\bGamma} \), then \( A \notin \check{\bGamma} ( \pre{ \omega}{ \omega } ) \) by Fact~\ref{fct:newfact} and by Wadge's lemma every \( B \in \bGamma ( \pre{ \omega}{ \omega } ) \) is of the form \( g^{-1} ( A ) \) for some continuous function \( g \colon \pre{ \omega}{\omega} \to \pre{ \omega}{\omega} \), and hence \( \rho \circ g \) is a \( \bGamma \)-norm on \( B \).
Therefore the pointclass \( \bGamma \) has the prewellordering property, and since \( \bGamma \) is assumed to be closed under projections and \( A \) witnesses \( \bGamma ( \pre{ \omega}{\omega} ) \neq \check{\bGamma} ( \pre{ \omega}{\omega} ) \), this implies that \( \bGamma \) is closed under well-ordered unions of arbitrary length by e.g.\ \cite[Theorem 2.16]{Jackson:2008pi}.
\end{proof}

Thus in the \( \AD \) world Lemma~\ref{lem:inthecodes} can be reformulated as follows.

\begin{proposition}[\( \AD +\DC \)] \label{prop:inthecodes}
Let \( \bGamma \) be a nontrivial boldface pointclass closed under projections, countable unions, and countable intersections.
Let \( \kappa \) be an infinite cardinal, and suppose that there is a \( \bGamma \)-code \(\rho\) for \( \kappa \).
Then for every \( f \colon \pre{ \omega}{2} \to \pre{ \kappa}{2} \) the following are equivalent:
\begin{enumerate-(a)}
\item \label{prop:inthecodes-1}
\( f \) is \( \bGamma \)-in-the-codes (with respect to \( \rho \));
\item \label{prop:inthecodes-2}
\( f^{-1} ( \widetilde{\Nbhd}\vphantom{\Nbhd}^\kappa_{\alpha , i}) \in \bDelta_{\bGamma} \) for every \( \alpha < \kappa \) and \( i = 0 , 1 \);
\item \label{prop:inthecodes-3}
\( f^{-1} ( U ) \in \bDelta_{\bGamma} \) for every \( U \in \mathcal{B}_p ( \pre{ \kappa}{2} ) \).
\end{enumerate-(a)}
\end{proposition}

If \( \bGamma \) is a nontrivial boldface pointclass closed under projections, countable unions, and countable intersections, the hypotheses of Proposition~\ref{prop:inthecodes} are fulfilled when:
\begin{itemize}[leftmargin=1pc]
\item
\( \kappa < \bdelta_{ \bGamma} \) (by definition of \( \bdelta_{ \bGamma} \)),
\item
\( \kappa = \bdelta_{ \bGamma} \) and \( \bGamma \) is a Spector pointclass closed under co-projections --- see~\cite[Exercise 4C.14]{Moschovakis:2009fk}.
\end{itemize}
In particular, Proposition~\ref{prop:inthecodes} can be applied to \( \bSigma^1_n \) and \( \kappa^+ = \bdelta^1_n \) (i.e.~\( \kappa = \blambda^1_n \) if \( n \) is odd and \( \kappa = \bdelta^1_{n-1} \) if \( n >0 \) is even), and to \( \bSigma^2_1 \) and \( \kappa = \bdelta^2_1 \). 
As we shall see in Proposition~\ref{prop:Skappainthecodesareborel} and Corollary~\ref{cor:inthecodesSouslin}, in all these (and other) interesting cases the \( \bGamma \)-in-the-codes functions \( f \colon \pre{ \omega}{2} \to \pre{ \kappa}{2} \) turn out to be just a special case of weakly \( \kappa + 1 \)-Borel functions, and for the odd levels of the projective hierarchy we further have that in fact the two notions coincide.

\begin{remarks} \label{rmk:rhononessential}
\begin{enumerate-(i)}
\item\label{rmk:rhononessential-a}
If \( \bGamma = \bSigma^{1}_{1} \) then \( \AD \) is not needed in Proposition~\ref{prop:inthecodes} and the proof goes through in \( \ZF + \AComega ( \R ) \).
\item\label{rmk:rhononessential-b}
Definition~\ref{def:Gammainthecodes} seems to depend on the particular choice of the \( \bGamma \)-norm \( \rho \): however, since Lemma~\ref{lem:inthecodes}\ref{lem:inthecodes-2} and Proposition~\ref{prop:inthecodes}\ref{prop:inthecodes-2} above do not depend on \( \rho \), Lemma~\ref{lem:inthecodes} and Proposition~\ref{prop:inthecodes} show that this is not the case for any sufficiently closed nontrivial boldface pointclass \( \bGamma \).
\end{enumerate-(i)}
\end{remarks}

\section{The generalized Baire space and Baire category} \label{sec:otherspacesandBairecategory}
\subsection{The generalized Baire space} \label{subsec:Bairespace}
So far we just considered the generalized Cantor space \( \pre{ \kappa }{ 2 } \), but similar results hold for the generalized Baire space\index[concepts]{generalized Baire space} \( \pre{ \kappa }{ \kappa } \), with \( \kappa \) an uncountable cardinal. 

\begin{definition} \label{def:generalizedBaire}
Let \( \widehat{\Nbhd}\vphantom{\Nbhd}_s^ \kappa \coloneqq \setof{x \in \pre{ \kappa }{ \kappa } }{ s \subseteq x } \),\index[symbols]{Nbhdc@\( \widehat{\Nbhd}\vphantom{\Nbhd}_s^ \kappa \)} where \( s \colon u \to \kappa \) and \( u \subseteq \kappa \).
\begin{itemize}[leftmargin=1pc]
\item 
The \markdef{bounded topology} \( \tau_b \) on \( \pre{ \kappa}{\kappa} \) is the one generated by the collection 
\[ 
\widehat{\mathcal{B}}_b \coloneqq \setof{\widehat{\Nbhd}\vphantom{\Nbhd}^\kappa_s }{ s \in \pre{ < \kappa }{ \kappa } } . \index[symbols]{Bc@\( \widehat{\mathcal{B}}_b \)}
\] 
The sets \( \widehat{\Nbhd}\vphantom{\Nbhd}^\kappa_s \) with \( s \colon u \to \kappa \) and \( u \) a \emph{bounded} subset of \( \kappa \), form a basis for \( \tau_b \) as well.
\item 
The \markdef{product topology} \( \tau_p \) on \( \pre{ \kappa }{ \kappa } \) is the product of \( \kappa \) copies of the space \( \kappa \) with the discrete topology.
A basis for \( \tau_p \) is 
\[ 
\widehat{\mathcal{B}}_p \coloneqq \setof{ \widehat{\Nbhd}\vphantom{\Nbhd}^\kappa_s }{ s \colon u \to \kappa \AND u \in [ \kappa ]^{ < \omega } } .
\]
\item
For \( \omega \leq \lambda \leq \kappa \), the \markdef{\( \lambda \)-topology} \( \tau_\lambda \) on \( \pre{ \kappa }{ \kappa } \) is the one generated by the collection 
\[ 
\widehat{\mathcal{B}}_\lambda \coloneqq \setof{\widehat{\Nbhd}\vphantom{\Nbhd}^\kappa_s }{ s \colon u \to \kappa \AND u \in [ \kappa ]^{ <  \lambda } },
\]
with the proviso as in Convention~\ref{cnv:basisfortau_kappa} that 
\begin{equation}\label{eq:conventionbasis}
 \kappa = \cof ( \kappa ) \IMPLIES \widehat{\mathcal{B}}_ \kappa \coloneqq \setof{\widehat{\Nbhd}\vphantom{\Nbhd}^\kappa_s }{ s \in  \pre{ < \kappa }{ \kappa } }  = \widehat{\mathcal{B}}_ b.
\end{equation}
\end{itemize}
\end{definition}

If \( X \) is a subspace of \( ( \pre{ \kappa}{\kappa}, \tau_* ) \) where \( * \) is one of \( b \), \( p \), or \( \lambda \), then the relative topology on \( X \) is still denoted by \( \tau_* \), so that when \( X = \pre{ \kappa}{2} \) this agrees with Definitions~\ref{def:generalizedCantor} and~\ref{def:lambdatopology}. 
Another subspace of \( \pre{ \kappa}{\kappa} \) of interest to us is 
\[ 
\Sym ( \kappa ) , \index[symbols]{Sym@\( \Sym ( \kappa ) \)}
\] 
the group of all permutations of \( \kappa \), which turns out to be an intersection of \( \kappa \)-many \( \tau_p \)-open (and hence also \( \tau_\lambda \)-open and \( \tau_b \)-open) sets, and hence a \( \bPi^0_2 ( \tau_* ) \) set, where \( * \in \{ b , p , \lambda \} \).

Most of the observations and results on \( \pre{ \kappa}{2} \) seen in the previous sections hold \emph{mutatis mutandis} for \( \pre{ \kappa}{\kappa} \).
We summarize some of them in the following proposition, in which all the topologies are understood to be on \( \pre{ \kappa}{\kappa} \).

\begin{proposition}\label{prop:summarizeonBairespace} 
Let \( \kappa \) be an uncountable cardinal, and let \( \omega \leq \lambda < \max ( \cof ( \kappa )^+ , \kappa ) \).
\begin{enumerate-(a)}
\item
If \( \lambda < \nu  < \max ( \cof ( \kappa )^+ , \kappa ) \) then  \( \tau_p = \tau_\omega \subseteq \tau_\lambda \subset \tau_\nu \subseteq \tau_{\cof ( \kappa ) }\subseteq \tau_b \), and  \( \tau_b = \tau_{ \cof ( \kappa ) }\) if and only if  \( \kappa \) is regular.
\item
\( \card{\tau_p} = \card{\pow ( \kappa )} \),  and \( \card{\tau_b} = \card{\pow ( \pre{ <\kappa}{\kappa} )} \). 
Therefore,
\begin{itemize}
\item
assuming \( \AD \), \(  \kappa < \Theta  \IMPLIES \card{\tau_p} < \card{\tau_b}  \), 
\item
assuming \( \AC  \),  \( \card{\tau_p} < \card{\tau_b} \iff 2^\kappa < 2^{( \kappa^{< \kappa})} ) \).
\end{itemize}
\item \label{prop:summarizeonBairespace-c}
The topologies  \( \tau_p \), \( \tau_\lambda \),  \( \tau_b  \) are perfect, regular Hausdorff, and zero-dimensional, and they are never \( \kappa \)-compact.
\item\label{prop:summarizeonBairespace-4}
The topology \(  \tau_* \) is not first countable (and hence neither second countable nor metrizable) and also not separable, for \( * \in \setLR{ p , \lambda } \).
\item\label{prop:summarizeonBairespace-5}
The topology \( \tau_b \) is neither second countable nor separable, and in fact assuming \( \AC \) it has density \( \kappa^{< \kappa} \). 
It is first countable if and only if it is metrizable if and only if it is completely metrizable if and only if \( \cof ( \kappa) = \omega \). 
\item
The topology \( \tau_* \) with \( * \in \setLR{ p ,  \lambda , b } \) is closed under intersections of length \( \leq \alpha \) (for some ordinal \( \alpha \)) if and only if: \( \alpha < \omega \) when \( * = p \), \( \alpha < \lambda \) when \( * = \lambda \) (assuming \( \AC \)), \( \card{\alpha} < \cof ( \kappa ) \) when \( * = b \).
Therefore the  collection of all \( \tau_* \)-clopen subsets of \( \pre{\kappa}{ \kappa } \) is a \( \omega \)-algebra if \( * = p \), is a \( \lambda \)-algebra (assuming \( \AC \)) if \( * =  \lambda  \), and  is a \( \cof ( \kappa ) \)-algebra if \( * = b \). 
\item\label{prop:summarizeonBairespace-7}
Assume \( \AC \). 
 Then \( \bB_{ \kappa + 1} ( \pre{ \kappa }{ \kappa } , \tau_ b ) \neq \pow ( \pre{ \kappa }{ \kappa } ) \).
 Moreover  \( \bSigma^0_{ \alpha } (  \pre{ \kappa }{ \kappa } ,  \tau_b ) \neq \bPi^0_{ \alpha } (  \pre{ \kappa }{ \kappa } , \tau_b  ) \) and \( \bSigma^0_{ \alpha } (  \pre{ \kappa }{ \kappa },  \tau_b  ) \subset  \bSigma^0_{ \beta } (  \pre{ \kappa }{ \kappa }, \tau_b ) \), for  \( 1 \leq \alpha < \beta < \kappa ^+ \). 
\item\label{prop:summarizeonBairespace-8}
Assume \( \AC \). 
Then \( \card{ \bB_{\kappa + 1} ( \pre{ \kappa}{\kappa}, \tau_p)} = 2^\kappa = \card{\tau_p} \) and \( \card{\bB_{\kappa + 1} ( \pre{ \kappa}{\kappa}, \tau_b ) } = 2^{(\kappa^{<\kappa})} = \card{\tau_b} \). 
Therefore, \( \card{\bB_{\kappa + 1} ( \pre{ \kappa}{\kappa}, \tau_p)} < 2^{(2^{\kappa})} \) and \( \card{\bB_{\kappa + 1} ( \pre{ \kappa}{\kappa}, \tau_b ) } < \card{\pow( \pre{ \kappa}{\kappa})} \IFF 2^{(\kappa^{<\kappa})} < 2^{(2^\kappa)} \).
In particular, \( \bB_{\kappa + 1} ( \pre{ \kappa}{\kappa}, \tau_p) \neq \pow ( \pre{ \kappa}{\kappa} ) \) and \( 2^{(\kappa^{< \kappa})} < 2^{ ( 2^\kappa )} \IMPLIES \bB_{\kappa + 1} ( \pre{ \kappa}{\kappa}, \tau_b ) \neq \pow ( \pre{ \kappa}{\kappa} ) \). 
\item
A subset of \( \pre{ \kappa}{\kappa} \) is (effective) \( \pre{ <\kappa}{\kappa} + 1 \)-Borel with respect to \( \tau_p \) if and only if it is (effective) \( \pre{ < \kappa}{\kappa} + 1 \)-Borel with respect to \( \tau_b \). 
Therefore, if we assume \( \AC \) and  \( \kappa^{< \kappa} = \kappa \) then \( \bB_{\kappa + 1} ( \pre{ \kappa}{\kappa}, \tau_p ) = \bB_{\kappa + 1} ( \pre{ \kappa}{\kappa}, \tau_b ) \).
\end{enumerate-(a)}
\end{proposition}

\begin{remarks}\label{rmks:Bairesummarize}
\begin{enumerate-(i)}
\item\label{rmks:Bairesummarize-i}
Part~\ref{prop:summarizeonBairespace-4} of Proposition~\ref{prop:summarizeonBairespace} should be contrasted with part~\ref{prop:topologicalproperties-e} of Proposition~\ref{prop:topologicalproperties}: this is one of the few differences between the generalized Cantor and Baire spaces.
\item\label{rmks:Bairesummarize-ibis}
Part~\ref{prop:summarizeonBairespace-7} of Proposition~\ref{prop:summarizeonBairespace} follows from Theorem~\ref{th:B_kappa+1nontrivial} and Proposition~\ref{prop:hierarchydoesnotcollapse} together with the fact that \(  \pre{ \kappa }{2}  \) is a closed subspace of \(  \pre{ \kappa }{ \kappa }  \) and we are dealing with hereditary general boldface pointclasses.
\item\label{rmks:Bairesummarize-ii}
If \( \kappa < \kappa^{< \kappa} \), then it may be the case that \( \card{\bB_{\kappa + 1} ( \pre{ \kappa}{\kappa}, \tau_p ) } < \card{\bB_{\kappa + 1} ( \pre{ \kappa}{\kappa}, \tau_b )} \) (see Remark~\ref{rmk:differentborel}). 
\item\label{rmks:Bairesummarize-iii}
The proof that \( ( \pre{ \kappa}{\kappa},\tau_b) \) is (completely) metrizable if and only if \( \cof(\kappa) = \omega \) (part~\ref{prop:summarizeonBairespace-5} of Proposition~\ref{prop:summarizeonBairespace}) follows the same ideas as the proof of part~\ref{prop:topologicalproperties-g} of Proposition~\ref{prop:topologicalproperties}: it is enough to identify \( \pre{ \kappa }{ \kappa } \) with \( \pre{ \kappa }{G} \) with \( G \) a group of size \( \kappa \) and then use the Birkhoff-Kakutani theorem, or replace \( \prod_{n \in \omega} \pre{ \lambda_n}{2} \) with \( \prod_{n \in \omega} \pre{ \lambda_n}{\kappa} \) in the direct proof sketched in Proposition~\ref{prop:topologicalproperties}\ref{prop:topologicalproperties-g}.
\end{enumerate-(i)}
\end{remarks}

The spaces \( \pre{ \omega}{2} \) and \( \pre{ \omega}{ \omega} \) are not homeomorphic because \( \pre{ \omega}{2} \) is compact while \( \pre{ \omega}{\omega} \) not. 
For the same reason, we have that also \( ( \pre{ \kappa}{2}, \tau_p) \) and \( ( \pre{ \kappa}{\kappa}, \tau_p) \) are never homeomorphic by Proposition~\ref{prop:topologicalproperties}\ref{prop:topologicalproperties-b} and Proposition~\ref{prop:summarizeonBairespace}\ref{prop:summarizeonBairespace-c}. 
However, we are now going to show that in models of choice the situation becomes rather different when we endow the generalized spaces \( \pre{\kappa}{2} \) and \( \pre{\kappa}{\kappa} \) with the bounded topology. 
Recall that a regular cardinal \( \kappa \) is \markdef{weakly compact} if and only if it is \emph{strong limit} (i.e.\ such that \( 2^\lambda < \kappa \) for every \( \lambda < \kappa \)), and has the \emph{tree property}, that is: for every \( T \subseteq \pre{ < \kappa}{\kappa} \), if \( 0 < \card{T \cap \pre{ \alpha}{\kappa}} < \kappa \) for all \( \alpha < \kappa \) then there is a \( \kappa \)-branch through \( T \), that is a point \( x \in \pre{ \kappa}{\kappa} \) such that \( x \restriction \alpha \in T \) for all \( \alpha < \kappa \).

\begin{lemma} \label{lem:nonhomeomorphism}
Suppose \( \kappa \) is an infinite cardinal with the tree property and such that for every cardinal \( \lambda < \kappa \) there is no injection from \( \kappa \) into \(  \pre{  \lambda }{2} \). 
\begin{enumerate-(a)}
\item\label{lem:nonhomeomorphism-a}
Every well-orderable \( \tau_b \)-closed \( C \subseteq \pre{ \kappa}{2} \) of size \( \geq \kappa \) has an accumulation point in \( C \).
\item\label{lem:nonhomeomorphism-b}
There is no continuous bijection between \( ( \pre{ \kappa}{2},\tau_b) \) and \( ( \pre{ \kappa}{\kappa}, \tau_b ) \).
\end{enumerate-(a)}
\end{lemma}

\begin{proof}
\ref{lem:nonhomeomorphism-a}
Suppose towards a contradiction that all points of \( C \) are isolated in it. 
For each \( x \in C \), let \( s_x \in \pre{ < \kappa}{2} \) be the shortest sequence with \( C \cap \Nbhd_{s_x} = \{ x \} \) (so that, in particular, \( s_x \) and \( s_{x'} \) are incomparable whenever \( x \neq x' \)), and let \( T \) be the tree generated by these \( s_x \)'s, i.e.
\[ 
T \coloneqq \setof{t \in \pre{ < \kappa}{2}} { t \subseteq s_x \text{ for some } x \in C}.
 \] 
Fix \( \alpha < \kappa \) and let \( \lambda \coloneqq \card{\alpha} < \kappa \). 
Since \( \setof{s_x}{x \in C} \) is well-orderable, then \( T \cap \pre{ \alpha}{2} \) is well-orderable too, and hence \( \card{T \cap \pre{ \alpha}{2} } = \mu \) for some cardinal \(\mu\). 
Since \( \pre{ \alpha}{2} \) and \( \pre{ \lambda}{2} \) are in bijection and there is no injection from \( \kappa \) into \( \pre{ \lambda}{2} \), we have that \( \mu < \kappa \). 
Using a similar argument, we also get \( T \cap \pre{ \alpha}{2} \neq \emptyset \) (otherwise \( s_x \in \pre{ <\alpha}{2} \) for each \( x \in C \) and we would have an injection from \( \kappa \leq \card{C} \) into \( \pre{ \lambda}{2} \)). 
Therefore by the tree property \( T \) has a \( \kappa \)-branch \( x \in \pre{ \kappa}{2} \). 
It follows from the choice of the \( s_x \)'s that \( x \) is an accumulation point of \( C \), and hence also \( x \in C \) because \( C \) is closed, a contradiction with our assumption that \( C \) consists only of isolated points.

\smallskip

\ref{lem:nonhomeomorphism-b}
Let now \( C' \coloneqq \setof{ \langle \alpha \rangle \conc 0^{(\kappa)}}{\alpha < \kappa} \subseteq \pre{ \kappa}{\kappa} \), and suppose towards a contradiction that there is a continuous bijection \( f \colon ( \pre{ \kappa}{2},\tau_b) \to ( \pre{ \kappa}{\kappa},\tau_b) \). 
Since \( C' \) is a closed well-orderable set of size \( \geq \kappa \) with no accumulation point, then so would be \( C \coloneqq f^{-1} ( C' ) \subseteq \pre{ \kappa}{2} \): but such a \( C \) cannot exists by part~\ref{lem:nonhomeomorphism-a}, and we are done.
\end{proof}

\begin{remark}
Notice that sets \( C \subseteq \pre{ \kappa}{2} \) as in part~\ref{lem:nonhomeomorphism-a} of Lemma~\ref{lem:nonhomeomorphism} do exist: the set 
\[
C \coloneqq \setof{0^{(\alpha)} \conc 1 \conc 0^{(\kappa)}}{\alpha < \kappa} \cup \{ 0^{(\kappa)} \} 
\] 
is  closed, well-orderable, and of size \( \kappa \). 
Moreover, in Lemma~\ref{lem:nonhomeomorphism}\ref{lem:nonhomeomorphism-a} we cannot drop the assumption that \( C \) be closed: the set \( \setof{0^{(\alpha)} \conc 1 \conc 0^{(\kappa)}}{\alpha < \kappa} \) is well-orderable and of size \( \kappa \), but has no accumulation point in itself.
\end{remark}

\begin{proposition}[\(\AC\)] \label{prop:homeomorphism}
Let \( \kappa \) be an infinite cardinal.
\begin{enumerate-(a)}
\item \label{prop:homeomorphism-a}
If \( \kappa \) is regular, then \( ( \pre{ \kappa}{2},\tau_b) \) and \( ( \pre{ \kappa}{\kappa},\tau_b) \) are homeomorphic if and only if \( \kappa \) is not weakly compact (equivalently, by Remark~\ref{rem:cantorbasicproperties}\ref{rem:cantorbasicproperties-i}, \( \pre{ \kappa}{2} \) is not \( \kappa \)-compact).
\item \label{prop:homeomorphism-b}
If \( \kappa \) is singular, then \( ( \pre{ \kappa}{2},\tau_b) \) and \( ( \pre{ \kappa}{\kappa},\tau_b) \) are homeomorphic if and only if \( \kappa \) is not strong limit (equivalently, 
\footnote{To see that if \( \kappa \) is singular and not strong limit then \( 2^{< \kappa} > \kappa \), let \( \cof(\kappa) \leq \lambda < \kappa \) be such that \( 2^\lambda \geq \kappa \). 
Then \( \kappa < \kappa^{\cof(\kappa)} \leq (2^\lambda)^{\cof(\kappa)} = 2^\lambda \leq 2^{< \kappa} \).}
\( 2^{< \kappa} > \kappa \)).
\end{enumerate-(a)}
\end{proposition}

Part~\ref{prop:homeomorphism-a} follows from the results in~\cite{Hung:1973fu}. 
However, for the reader's convenience we give here a simple direct proof.

\begin{proof}
\ref{prop:homeomorphism-a}
By Lemma~\ref{lem:nonhomeomorphism}\ref{lem:nonhomeomorphism-b}, \( ( \pre{ \kappa}{2}, \tau_b) \) and \( ( \pre{ \kappa}{\kappa},\tau_b) \) cannot be homeomorphic if \( \kappa \) is weakly compact. 
(Alternatively, we could use~\cite[Theorem 5.6]{Motto-Ros:2011qc} and the fact that \( \pre{ \kappa}{\kappa} \) is never \( \kappa \)-compact by Proposition~\ref{prop:summarizeonBairespace}\ref{prop:summarizeonBairespace-c}.)
Therefore we have just to show that there is a homeomorphism \( f \colon ( \pre{ \kappa}{\kappa},\tau_b) \to ( \pre{ \kappa}{2},\tau_b) \) whenever \( \kappa \) is not weakly compact.

Let us first assume that there is \( \lambda < \kappa \) with \( 2^ \lambda \geq \kappa \), so that \( \kappa^\lambda = 2^\lambda \), and fix a bijection \( g \colon \pre{ \lambda}{\kappa} \to \pre{ \lambda}{2} \).
Given \( x \in \pre{ \kappa}{\kappa} \), let \( \seqof{s^x_\alpha}{\alpha < \kappa} \) be the unique sequence of elements of \( \pre{ \lambda}{\kappa} \) such that \( x = s^x_0 \conc s^x_1 \conc \dotsc \conc s^x_\alpha \conc \dotsc \). 
The map
\[ 
f \colon \pre{ \kappa}{\kappa} \to \pre{ \kappa}{2}, \quad x \mapsto g ( s^x_0 ) \conc g ( s^x_1 ) \conc \dotsc \conc g ( s^x_\alpha ) \conc \dotsc
 \] 
is a well-defined bijection. 
Moreover, since the families 
\[ 
\setof{\widehat{\Nbhd}_s}{s \in \pre{ < \kappa}{\kappa}, \lh{s} = \lambda \cdot \alpha \text{ for some } \alpha < \kappa} \subseteq \widehat{\mathcal{B}}_b ( \pre{ \kappa}{\kappa} )
\]
and 
\[ 
\setof{\Nbhd_s}{s \in \pre{ < \kappa}{2}, \lh{s} = \lambda \cdot \alpha \text{ for some } \alpha < \kappa} \subseteq \mathcal{B}_b( \pre{ \kappa}{2})
\]
are bases for the bounded topologies on, respectively, \( \pre{ \kappa}{\kappa} \) and \( \pre{ \kappa}{2} \), we get that \( f \colon ( \pre{ \kappa}{\kappa} , \tau_b ) \to ( \pre{ \kappa}{2} , \tau_b ) \) is a homeomorphism.

Let us now assume that \( \kappa \) is strong limit, so that \( 2^{< \kappa} = \kappa \).
As \( \kappa \) is not weakly compact,  there is a tree \( T \subseteq \pre{ \kappa}{2} \) of height \( \kappa \) without a \( \kappa \)-branch. 
Let
\[ 
\partial T \coloneqq \setof{s \in \pre{ < \kappa}{2} \setminus T}{s \restriction \alpha \in T \text{ for all } \alpha < \lh{s}} \index[symbols]{41@\( \partial T \)}
 \] 
be the \markdef{boundary of  \( T \)}\index[concepts]{boundary of a tree \( \partial T \)}.
The sequences in \( \partial T \) are pairwise incomparable, and as \( T \) has no \( \kappa \)-branches, \( \setof{\Nbhd_s}{s \in \partial T} \) is a partition of \( \pre{ \kappa}{2} \). 
We claim that such partition has size \( \kappa \): towards a contradiction if \( \card{\partial T} < \kappa \), then \( \partial T \subseteq \pre{ \lambda}{2} \) for some \( \lambda < \kappa \) because \( \kappa \) is regular, and since each sequence in \( T \) can be extended to some \( s \in \partial T \), this contradicts the fact that \( T \) has height \( \kappa \).
If \( g \colon \kappa \to \partial T \) is a bijection,  the map
\[ 
f \colon \pre{ \kappa}{\kappa} \to \pre{ \kappa}{2}, \quad x \mapsto g ( x ( 0 ) ) \conc g ( x ( 1 ) ) \conc \dotsc \conc g ( x ( \alpha ) ) \conc \dotsc 
 \] 
is a homeomorphism between \( ( \pre{ \kappa}{\kappa}, \tau_b ) \) and \( ( \pre{ \kappa}{2},\tau_b ) \), as required.

\ref{prop:homeomorphism-b}
If \( \kappa \) is not strong limit, then there is \( \lambda < \kappa \) such that \( 2^\lambda \geq \kappa \). 
Arguing as in part~\ref{prop:homeomorphism-a}, we have that there is a homeomorphism \( f \colon ( \pre{ \kappa}{\kappa}, \tau_b) \to ( \pre{ \kappa}{2}, \tau_b)\). 
Conversely, assume that \( \kappa \) is strong limit, so that   \( 2^{< \kappa} = \kappa \). 
If \( \seqof{U_\alpha}{\alpha < \lambda} \) is a \( \tau_b \)-open partition of \( \pre{ \kappa}{2} \), then \( \lambda \leq \kappa \): in fact, the map assigning to each \( \alpha < \lambda \) some \( s_\alpha \in \pre{ <\kappa}{2} \) such that \( \Nbhd_{s_\alpha} \subseteq U_\alpha \) is injective and witnesses \( \lambda \leq 2^{< \kappa} = \kappa \). 
Since \( \setof{\widehat{\Nbhd}_s}{s \in \pre{ \cof(\kappa)}{\kappa}} \) is a \( \tau_b \)-open partition of \( \pre{ \kappa}{\kappa} \) of size \( > \kappa \), this shows that there is no continuous surjection \( f \colon ( \pre{ \kappa}{2},\tau_b) \to ( \pre{ \kappa}{\kappa}, \tau_b) \). 
\end{proof}

Proposition~\ref{prop:homeomorphism} shows that under \( \AC \) the spaces \( ( \pre{ \kappa}{2}, \tau_b) \) and \( ( \pre{ \kappa}{\kappa}, \tau_b) \) are in most cases homeomorphic, including e.g.\ when \( \kappa = \omega_1 \). 
The situation is quite different in models of determinacy. 
Assuming \( \AD \), the cardinal \( \omega_1 \) is measurable, so it has the tree property.
Moreover, \( \omega _1 \not\into \pre{ \omega}{2} \) by the \( \PSP \), and therefore \( ( \pre{ \omega _1 }{2},\tau_b) \) and \( ( \pre{ \omega _1 }{ \omega _1 }, \tau_b) \) are not homeomorphic by Lemma~\ref{lem:nonhomeomorphism}\ref{lem:nonhomeomorphism-b}.
This argument can be generalized to larger \emph{regular} cardinals: in fact, Steel and Woodin have shown that \( \AD + {\Vv = \Ll ( \R )} \) implies that every uncountable regular \( \kappa < \Theta \) is measurable, and hence it has the tree property, and \( \lambda^+ \not\into \pre{ \lambda}{2} \) for all \( \lambda < \Theta \) --- see~\cite{Steel:2010fk,Steel:yq}.
Woodin (unpublished) has weakened the hypothesis \( \AD + { \Vv = \Ll ( \R )} \) to \( \AD^+ + {\Vv = \Ll ( \pow ( \R ) )} \) --- see Definition~\ref{def:AD^+} for \( \AD^+ \).
Therefore by the argument above we have that:

\begin{proposition}\label{prop:nonhomeomorphism}
Assume \( \AD^+ + {\Vv = \Ll ( \pow ( \R ) )} \).
Then for every regular \( \kappa < \Theta \), the spaces \( ( \pre{ \kappa}{2},\tau_b) \) and \( ( \pre{ \kappa}{\kappa}, \tau_b) \) are not homeomorphic.
\end{proposition}

Therefore \( \AC \) is a necessary assumption for Proposition~\ref{prop:homeomorphism}. 
In fact, Propositions~\ref{prop:homeomorphism} and~\ref{prop:nonhomeomorphism} show that it is independent of \( \ZF + \DC \) whether \( ( \pre{ \omega_1}{2}, \tau_b) \) is homeomorphic to \( ( \pre{ \omega_1}{\omega_1}, \tau_b) \). 

\subsection{Baire category} \label{subsec:category}

\begin{definition}\label{def:category}
Let \( \mu \) be an infinite cardinal. 
A subset \( A \) of a topological space  \( ( X , \tau ) \) is said to be \markdef{\( \mu \)-meager} if it is the union of \(  \mu \)-many  nowhere dense sets, where a set is nowhere dense if it is disjoint from some open dense subset of \( X \). 
A set \( A \subseteq X \) is said to be \markdef{\( \mu \)-comeager} if its complement is \( \mu\)-meager.
If \( U \subseteq X \) is a nonempty open set, then \( A \subseteq X \) is \markdef{\( \mu \)-meager in \( U \)} (respectively, \markdef{\( \mu \)-comeager in \( U \)}) if \( A \cap U \) is \( \mu \)-meager (respectively, \( \mu \)-comeager) in the (topological) space \( U \) endowed with the relative topology induced by \( \tau \).
\end{definition}

Every subset of a \( \mu \)-meager set is \( \mu \)-meager as well. 
If \( \mu \leq \mu' \) are infinite cardinal and \( U \subseteq V \) are open sets of \( ( X , \tau ) \), then every set \( A \subseteq X \) which is \( \mu \)-(co)meager in \( V \) is also \( \mu' \)-(co)meager in \( U \).

\begin{definition} \label{def:muBaire}
Let \( \mu \) be an infinite cardinal.
A topological space \( X  \) is said to be \markdef{\( \mu \)-Baire} if the intersection of \( \mu \)-many open dense subsets of \( X \) is dense in \( X \).
\end{definition}

As for the classical notion of a Baire space (which corresponds to the case \( \mu = \omega \)), it is easy to check that the space \( X \) is \( \mu \)-Baire if and only if every \( \mu \)-comeager subset of \( X \) is dense, if and only if every nonempty open subset of \( X \) is not \( \mu \)-meager.
Moreover, if \( \mu \leq \mu' \) then every \( \mu' \)-Baire space is automatically \( \mu \)-Baire, and  every open subspace of a \( \mu \)-Baire space is \( \mu \)-Baire as well.

\begin{theorem}\label{th:Bairespace}
Let \( \kappa ,  \lambda \) be  cardinals such that \( \omega < \lambda  < \min ( \cof ( \kappa )^+ , \kappa ) \).
\begin{enumerate-(a)}
\item \label{th:Bairespace-1}
The space \( ( \pre{ \kappa}{\kappa}, \tau_p ) \) is \( \omega \)-Baire.
\item \label{th:Bairespace-2}
Assume \( \AC \). 
Then the space \( ( \pre{ \kappa}{\kappa}, \tau_\lambda ) \) is \( \cof ( \lambda ) \)-Baire.
\item \label{th:Bairespace-3}
Assume \( \AC \). 
Then the space \( ( \pre{ \kappa}{\kappa}, \tau_b ) \) is \( \cof ( \kappa ) \)-Baire.
\end{enumerate-(a)}
All of~\ref{th:Bairespace-1}--\ref{th:Bairespace-3} holds true when \( \pre{ \kappa}{\kappa} \) is replaced by \( \pre{ \kappa}{2} \).
\end{theorem}

Theorem~\ref{th:Bairespace} can be restated (and easily proved) in the language of forcing.
If \( \forcing{P} \) is a forcing notion (i.e.~a quasi-order), then \( \mathsf{FA}_\mu ( \forcing{P} ) \)\index[symbols]{FA@\( \mathsf{FA}_\mu ( \forcing{P} ) \)}, the \markdef{\( \mu \)-forcing axiom for \( \forcing{P} \)}, says that for any sequence \( \seqof{D_ \alpha }{ \alpha < \mu } \) of dense subsets of \( \forcing{P} \) there is a filter \( G \) intersecting all \( D_ \alpha \)'s.
Recalling that the set of partial functions is a forcing notion with reverse inclusion (see Section~\ref{subsubsec:sequences}), the preceding result asserts that certain forcing axioms hold:
\begin{enumerate}[label={\upshape (\alph**)}, leftmargin=2pc]
\item \label{?}
\( \mathsf{FA}_ \omega ( \forcing{Fn} ( \kappa , \kappa ; \omega ) ) \) and \( \mathsf{FA}_ \omega ( \forcing{Fn} ( \kappa , 2 ; \omega ) ) \).
\item \label{?}
Assume \( \AC \).
Then \( \mathsf{FA}_ {\cof ( \lambda ) } ( \forcing{Fn} ( \kappa , \kappa ; \lambda ) ) \) and \( \mathsf{FA}_ {\cof ( \lambda ) } ( \forcing{Fn} ( \kappa , 2 ; \lambda ) ) \).
\item \label{?}
Assume \( \AC \).
Then \( \mathsf{FA}_ {\cof ( \kappa ) } ( \forcing{Fn} ( \kappa , \kappa ; b ) ) \) and \( \mathsf{FA}_ {\cof ( \kappa ) } ( \forcing{Fn} ( \kappa , 2 ; b ) ) \).
\end{enumerate}

Theorem~\ref{th:Bairespace} cannot be extended to arbitrary closed subspaces of \(  \pre{ \kappa}{\kappa} \), as shown by the next example.

\begin{example} \label{xmp:notBaire}
Given a cardinal \( \kappa \) of uncountable cofinality, let \( C \subseteq \pre{ \kappa}{\kappa} \) be the collection of those \( x \in \pre{ \kappa}{2} \) such that \( \card{\setof{\alpha < \kappa}{x (\alpha) = 0}} < \omega \). 
Then \( C \) is \( \tau_{\omega_1} \)-closed, and hence also \( \tau_b \)-closed and \( \tau_\lambda \)-closed for every \( \omega_1 \leq \lambda \leq \cof(\kappa) \). 
We claim that \( C \) is not \(\omega\)-Baire with respect to any of the topologies \( \tau_p \), \( \tau_\lambda \), or \( \tau_b \). 
For each \( n \in \omega \), let \( U_n \coloneqq \setof{ x \in \pre{ \kappa}{\kappa}}{\card{\setof{\alpha < \kappa}{x(\alpha) = 0}} \geq n} \). 
Each \( U_n \) is  open and dense with respect to any of the above topologies, and \( U_n \cap C \neq \emptyset \). 
However \( \bigcap_{n \in \omega} (U_n \cap C) = \left ( \bigcap_{n \in \omega} U_n \right ) \cap C = \emptyset \), and hence \( C \) is not \(\omega\)-Baire.
\end{example}

For the results of this paper, we need to show that also \( \Sym(\kappa) \), which is a \( \bPi^0_2 \)-subset of \( \pre{ \kappa}{\kappa} \) with respect to any of our topologies, is \(\mu\)-Baire for appropriate cardinals \( \mu \). 
However, since Example~\ref{xmp:notBaire} shows that there may be even simple closed sets that are not \(\omega\)-Baire, the \( \mu \)-Baireness of \( \Sym(\kappa) \) needs to be proved independently of  Theorem~\ref{th:Bairespace}. 
To state the next result in the forcing jargon we need the following definition: for \( \kappa , \lambda \) infinite cardinals, let \index[symbols]{Inj@\( \forcing{Inj} \)}
\begin{align*}
\forcing{Inj} ( \kappa , \kappa ; \lambda ) & \coloneqq \setof{s \in \forcing{Fn} ( \kappa , \kappa ; \lambda ) }{ s \text{ is injective}} ,
\\
\forcing{Inj} ( \kappa , \kappa ; b ) & \coloneqq \setof{s \in \forcing{Fn}(\kappa,\kappa;b)  }{ s \text{ is injective}} = \pre{<\kappa}{(\kappa)}.
\end{align*}

\begin{theorem}\label{th:Bairespacesym}
Let \( \kappa ,  \lambda \) be  cardinals such that \( \omega < \lambda  < \min ( \cof ( \kappa )^+ , \kappa ) \).
\begin{enumerate-(a)}
\item \label{th:Bairespacesym-1}
The space \( ( \Sym ( \kappa ), \tau_p ) \) is \( \omega \)-Baire, i.e.~\( \mathsf{FA}_ \omega ( \forcing{Inj} ( \kappa , \kappa ; \omega ) ) \).
\item \label{th:Bairespacesym-2}
Assume \( \AC \). 
Then the space \( ( \Sym ( \kappa ), \tau_\lambda ) \) is \( \cof ( \lambda ) \)-Baire, i.e.~\( \mathsf{FA}_{\cof ( \lambda ) } ( \forcing{Inj} ( \kappa , \kappa ; \lambda ) ) \).
\item \label{th:Bairespacesym-3}
Assume \( \AC \). 
Then the space \( ( \Sym ( \kappa ), \tau_b ) \) is \( \cof ( \kappa ) \)-Baire, i.e.~\( \mathsf{FA}_{\cof ( \kappa ) } ( \forcing{Inj} ( \kappa , \kappa ; b ) )\).
\end{enumerate-(a)}
\end{theorem}

At first glance, the forcing axioms in the second part of~\ref{th:Bairespacesym-1}--\ref{th:Bairespacesym-3} may seem weaker than their counterparts for \( \mu \)-Baireness of \( \Sym ( \kappa ) \) because in all the three cases an arbitrary generic intersecting a given family of \(\mu\)-many dense subsets of the corresponding forcing poset may fail to be surjective. 
However, the arguments presented in the proof below  show that the two formulations of each point are indeed equivalent:  the existence of an arbitrary generic granted by the forcing axiom implies the existence of a surjective one, thus yielding the corresponding \( \mu \)-Baireness property for \( \Sym ( \kappa ) \).

\begin{proof}[Proof of Theorem~\ref{th:Bairespacesym}]
Let us first consider~\ref{th:Bairespacesym-2}. 
The argument is essentially the same as the one used in the classical Baire category theorem (see e.g.\ \cite[Theorem 8.4]{Kechris:1995zt}) --- we present it here for the reader's convenience. 
Let \( U \subseteq \Sym ( \kappa) \) be open and nonempty, and let \( \seqof{ U_\alpha }{ \alpha < \cof ( \lambda) } \) be a sequence of open dense subsets of \( \Sym ( \kappa) \). 
For every \( \alpha < \cof ( \lambda) \), we recursively define \( s_\alpha \in \forcing{Inj} ( \kappa  , \kappa ; \lambda) \) so that \( \widehat{\Nbhd}\vphantom{\Nbhd}_{s_\alpha} \cap \Sym ( \kappa ) \subseteq U \cap {\bigcap_{\beta<\alpha} U_\beta} \), and \( s_\beta \subseteq s_\alpha \) for \( \beta \leq \alpha < \cof ( \lambda) \). 
Let \( s_0 \in \forcing{Inj} ( \kappa , \kappa ; \lambda ) \) be such that \( \widehat{\Nbhd}\vphantom{\Nbhd}_{s_0} \cap \Sym ( \kappa ) \subseteq U \). 
Let now \( \alpha \coloneqq \gamma +1 \): since \( U_\alpha \) is open and dense in \( \Sym(\kappa) \) and \( s_\gamma \in \forcing{Inj} ( \kappa , \kappa ; \lambda ) \), \( U_\alpha \cap \widehat{\Nbhd}\vphantom{\Nbhd}_{s_\gamma} \) is nonempty and open in \( \Sym ( \kappa ) \). 
Pick \( s' \in \forcing{Inj} ( \kappa , \kappa ; \lambda ) \) so that \( \widehat{\Nbhd}\vphantom{\Nbhd}_{s'} \cap \Sym ( \kappa ) \subseteq U_\alpha \cap \widehat{\Nbhd}\vphantom{\Nbhd}_{s_\gamma} \), and notice that \( s' \) is necessarily compatible with \( s_\gamma \).
Then \( s_\alpha \coloneqq s_\gamma \cup s' \) has the required properties. 
Finally, let \( \alpha < \cof ( \lambda) \) be limit. 
Since all the \( s_\beta \)'s for \( \beta < \alpha \) are compatible and belong to \( \forcing{Inj} ( \kappa , \kappa ; \lambda )  \), \( s_\alpha \coloneqq \bigcup_{\beta < \alpha} s_\beta \in \forcing{Inj} ( \kappa , \kappa ; \lambda ) \) has the required properties. 
Let now \( t \coloneqq \bigcup_{\alpha < \cof ( \lambda)} s_\alpha \), so that \( t \) is an injective (partial) function from \( \kappa \) into itself.
Since \( \card{\dom ( t ) } = \card{\ran ( t ) } \leq \lambda < \kappa \), we can pick a bijection \( t' \) between \( \kappa \setminus \dom ( t ) \) and \( \kappa \setminus \ran ( t ) \): then \( x \coloneqq t \cup t' \in \Sym ( \kappa) \), and \( x \) witnesses \( U \cap {\bigcup_{\alpha < \cof ( \lambda)} U_\alpha} \neq \emptyset \).

\smallskip

\ref{th:Bairespacesym-1} can be proved in a similar way, using the fact that \( \kappa > \omega \) --- just notice that when dealing with the product topology we do not need the Axiom of Choice to recursively pick the sequences \( s_\alpha \) because \( \card{[ \kappa ] ^{< \omega}} = \kappa \) implies that \( \widehat{\mathcal{B}}_p ( \pre{ \kappa}{\kappa} ) \) is well-orderable. 

\smallskip

Finally, to prove~\ref{th:Bairespacesym-3} we use a sort of back-and-forth argument. 
Let \( U \subseteq \Sym ( \kappa) \) be open and nonempty, and let \( \seqof{ U_\alpha }{ \alpha < \cof ( \kappa ) } \) be a sequence of open dense subsets of \( \Sym ( \kappa) \). 
Fix an increasing sequence \( \seqof{ \lambda_\alpha }{ \alpha < \cof ( \kappa) } \) of ordinals cofinal in \( \kappa \). 
Let \( s_0 \in \pre{ < \kappa}{\kappa} \) be an injective sequence such that \( \widehat{\Nbhd}_{s_0} \cap \Sym ( \kappa) \subseteq U \), and for \( \alpha \) a limit ordinal let \( s_\alpha \coloneqq \bigcup_{\beta < \alpha} s_\beta \), which is  injective whenever all the \( s_\beta \)'s are injective. 
Let now \( \gamma \) be either \( 0 \) or a limit ordinal \( {<} \cof ( \kappa) \). 
For \( n \in \omega \) and \( i \in \setLR{ 1 , 2 } \), we recursively define \( s_{\gamma + 2n + i} \) as follows. 
If \( i = 1 \) we let \( s_{\gamma + 2 n + i} \in \pre{ < \kappa}{\kappa} \) be any injective sequence extending \( s_{\gamma+2n} \) such that \( \widehat{\Nbhd}_{s_{\gamma + 2 n + 1}} \subseteq U_{\gamma + n} \) (which exists because the \( U_\alpha \)'s are open and dense in \( \Sym ( \kappa) \)). 
If instead \( i = 2 \), let \( \eta \coloneqq \lh s_{\gamma + 2 n + 1} \), and let \( \seqof{ j_l }{ l < \delta } \) be a strictly increasing enumeration of \( \lambda_{\gamma+n} \setminus \ran ( s_{\gamma + 2 n + 1} ) \) (for a suitable \( \delta \leq \lambda_{\gamma + n} \)). 
Let \( s_{\gamma + 2n +2 } \) be the extension of \( s_{\gamma + 2 n + 1} \) of length \( \eta + \delta \) obtained by setting \( s_{\gamma+2n+2} ( \eta + l ) \coloneqq j_l \) for every \( l < \delta \), so that \( s_{\gamma + 2 n + 2 } \) is still injective and \( \lambda_{\gamma + 2 n} \subseteq \ran(s_{\gamma+2n+2}) \). 
Then it is easy to check that \( x \coloneqq \bigcup_{\alpha < \cof ( \kappa)} s_\alpha \in \Sym ( \kappa) \) and \( x \in U \cap {\bigcup_{\alpha < \cof ( \lambda)} U_\alpha} \).
\end{proof}

The following is a variant of the classical Baire property, corresponding to the case \( \mu = \omega \).

\begin{definition}
Let \( \mu \) be an infinite cardinal and \( X  \) be a topological space. 
We say that \( A \subseteq X \) has the \markdef{\( \mu \)-Baire property} if there is an open set \( U \subseteq X \) such that the symmetric difference \( A \symdif U \) is \( \mu \)-meager.
\end{definition}

As shown in~\cite{Halko:2001kl}, arguing as in the classical case one easily gets the following results (see also e.g.\ \cite[Proposition 8.22]{Kechris:1995zt} and~\cite[Proposition 8.26]{Kechris:1995zt}).

\begin{proposition} \label{prop:borelsetshaveBP}
Let \( \mu \) be an infinite cardinal and \( X  \) be a topological space. 
The class of all subsets of \( X \) having the \( \mu \)-Baire property is a \( \mu+1 \)-algebra on \( X \), and in fact it is the \( \mu+1 \)-algebra on \( X \) generated by all open sets and all \( \mu \)-meager sets. 
In particular, all sets in \( \bB_{\mu+1} ( X  ) \) have the \( \mu \)-Baire property.
\end{proposition}

\begin{proposition} \label{prop:nonmeagerarelocallycomeager}
Let \( \mu \) be an infinite cardinal and let \(X \) be a topological space. 
If \( A \subseteq X \) has the \( \mu \)-Baire property and is not \( \mu \)-meager, then there is a nonempty open set \( U \subseteq X \) such that \( A \) is \( \mu \)-comeager in \( U \).
\end{proposition}

\section{\texorpdfstring{Standard Borel \( \kappa \)-spaces, \( \kappa \)-analytic quasi-orders, and spaces of codes}{Standard Borel kappa-spaces, kappa-analytic quasi-orders, and spaces of codes}}\label{sec:spacesofcodes}
\subsection{\( \kappa \)-analytic sets} \label{subsec:kappa-analytic}
Recall that a subset \( A \) of a Polish space \( X \) is called analytic if it is a continuous image of a closed subset of the Baire space \( \pre{\omega}{\omega} \). 
Here are some reformulations of this notion, where \( \PROJ \) denotes the projection on the first coordinate, as defined in~\eqref{eq:projectionCh2}:
\begin{itemize}[leftmargin=1pc]
\item
\( A \subseteq X \) is analytic if and only if it is either empty or a continuous image of the whole Baire space \( \pre{\omega}{\omega} \);
\item
\( A \subseteq X \) is analytic if and only if it is a continuous image of a Borel subset of \( \pre{\omega}{\omega} \);
\item
\( A \subseteq X \) is analytic if and only if \( A = \PROJ F \) for some closed \( F \subseteq X \times \pre{\omega}{\omega} \);
\item
\( A \subseteq X \) is analytic if and only if \( A = \PROJ B \) for some Borel \( B \subseteq X \times \pre{\omega}{\omega} \).
\end{itemize}
It is not hard to see that the class of analytic sets contains all Borel sets and is closed under countable unions, countable intersections, and images and preimages under Borel functions.
In Section~\ref{sec:Ksouslinsets} the collection of analytic sets will be identified with \( \bS ( \omega ) \), the class of \( \omega \)-Souslin sets.

We now generalize the notion of analytic set to the uncountable context. 
To simplify the presentation, in the subsequent definition and results we endow the generalized Baire space \( \pre{\kappa}{\kappa} \) with the bounded topology, so that all the related topological notions (such as continuous functions, \( \kappa+1 \)-Borel sets, and so on) tacitly refer to \( \tau_b \). 
Of course analogous notions can be obtained by replacing, \emph{mutatis mutandis}, the topology \( \tau_b \) with any of the topologies introduced in Definition~\ref{def:generalizedBaire}. 
However, since we have no use for these variants in the rest of the paper, for the sake of simplicity we leave to the reader the burden of checking which of the properties stated below transfer to these topologies.

\begin{definition}\label{def:kappa-analytic}
A set \( A \subseteq \pre{\kappa}{\kappa} \) is called \markdef{\( \kappa \)-analytic} if it is a continuous image of a closed subset of \( \pre{\kappa}{\kappa} \).\index[concepts]{space!\( \kappa \)-analytic space}
\end{definition}

Unlike in the classical case \( \kappa = \omega \), when \( \kappa > \omega \) it is no more true in general that nonempty \( \kappa \)-analytic sets are continuous images of the whole \( \pre{\kappa}{\kappa} \) --- see~\cite{Luecke:2014ar}. 
However, all the other equivalent reformulations mentioned above remain true also in our new context. 
The key result to prove this is the following proposition.

\begin{proposition} \label{prop:borelsetsareanalytic}
Every effective \( \kappa+1 \)-Borel subset of \( \pre{\kappa}{\kappa} \) is \( \kappa \)-analytic.
\end{proposition}

\begin{proof}
We modify the proof of ~\cite[Lemma 3.9 and Proposition 3.10]{Motto-Ros:2011qc}, where the desired result is proved under%
\footnote{Under choice all \( \kappa+1 \)-Borel sets are effective.} 
 \( \AC \) and the extra cardinal assumption \( \kappa^{< \kappa} = \kappa \). 
In particular, in that proof it is argued in \( \ZF \) that given a family \( \setof{C_\alpha}{\alpha< \kappa} \) of closed subsets of \( \pre{\kappa}{\kappa} \) and a corresponding family of continuous maps \( f_\alpha \colon C_\alpha \to \pre{\kappa}{\kappa} \) with range \( A_\alpha \), one can canonically construct two closed sets \( C_\cap , C_\cup \subseteq \pre{\kappa}{\kappa} \) and continuous maps \( f_\cap \colon C_\cap \to \pre{\kappa}{\kappa} \) and \( f_\cup \colon C_\cup \to \pre{\kappa}{\kappa} \) such that \( f_\cap \) surjects onto \( \bigcap_{\alpha < \kappa} A_\alpha \) and \( f_\cup \) surjects onto \( \bigcup_{\alpha<\kappa} A_\alpha \). 
Moreover, the identity function witnesses that every closed subset of \( \pre{\kappa}{\kappa} \) is \( \kappa \)-analytic. 
We now show (in \( \ZF \), and without assuming \( \kappa^{< \kappa} = \kappa \)) that also all open sets are \( \kappa \)-analytic. 
Let \( U \subseteq \pre{\kappa}{\kappa} \) be open, and for every \( \alpha < \kappa \) set \( S_\alpha \coloneqq \setof{s \in \pre{\alpha}{\kappa}}{\widehat{\Nbhd}\vphantom{\Nbhd}^\kappa_s \subseteq U} \), so that \( U = \bigcup_{\alpha< \kappa} U_\alpha \) with \( U_\alpha \coloneqq \bigcup_{s \in S_\alpha} \widehat{\Nbhd}\vphantom{\Nbhd}^\kappa_s \). 
Notice that each \( U_\alpha \) is clopen, and therefore \( \kappa \)-analytic. 
Setting \( C_\alpha \coloneqq U_\alpha \) and \( f_\alpha \coloneqq {\id \restriction U_\alpha} \), we get that the map \( f_\cup \colon C_\cup \onto \bigcup_{\alpha< \kappa} U_\alpha \subseteq \pre{\kappa}{\kappa} \) as above witnesses that \( U \) is \( \kappa \)-analytic.

Let now \( B \subseteq \pre{\kappa}{\kappa} \) be effective \( \kappa+1 \)-Borel and let \( (T, \phi) \) be a \( \kappa+1 \)-Borel code for it. 
Using the facts mentioned in the previous paragraph, one can easily build by recursion on the rank of the nodes of the well-founded tree \( T \) a map \( g \) on \( T \) assigning to each \( t \in T \) two closed sets \( C_t , C'_t \subseteq \pre{\kappa}{\kappa} \) and two continuous functions \( f_t \colon C_t \to \pre{\kappa}{\kappa} \) and \( f'_t \colon C'_t \to \pre{\kappa}{\kappa} \) such that \( f_t \) surjects onto \( \phi ( t ) \) and \( f'_t \) surjects onto \( \pre{\kappa}{\kappa} \setminus \phi ( t ) \). 
In particular, \( C_\emptyset \) and \( f_\emptyset \) witness that \( B \) is \( \kappa \)-analytic.
\end{proof}

The proof of Proposition~\ref{prop:borelsetsareanalytic} also implies that the class of \( \kappa \)-analytic subsets of \( \pre{\kappa}{\kappa} \) is closed under finite intersections and finite unions. 
It is also closed under intersections and unions of length \( \alpha \leq \kappa \) as long as given a family of \( \kappa \)-analytic sets \( \setof{A_\beta}{\beta<\alpha} \) one can choose witnesses \( f_\beta \colon C_\beta \onto A_\beta \) of this.
Therefore if we assume \( \AC_\kappa \) the class of \( \kappa \)-analytic subsets of \( \pre{\kappa}{\kappa} \) is closed under intersections and unions of length 
\( \leq \kappa \).

\begin{corollary} \label{cor:kappaanalytic}
The following are equivalent for a set \( A \subseteq \pre{\kappa}{\kappa} \):
\begin{enumerate-(a)}
\item \label{cor:kappaanalytic-1}
\( A \) is \( \kappa \)-analytic;
\item \label{cor:kappaanalytic-2}
\( A \) is a continuous image of an effective \( \kappa+1 \)-Borel subset of \( \pre{\kappa}{\kappa} \);
\item \label{cor:kappaanalytic-3}
\( A = \PROJ F \) for some closed \( F \subseteq \pre{\kappa}{\kappa} \times \pre{\kappa}{\kappa} \);
\item \label{cor:kappaanalytic-4}
\( A = \PROJ B \) for some effective \( \kappa+1 \)-Borel \( B \subseteq \pre{\kappa}{\kappa} \times \pre{\kappa}{\kappa} \).
\end{enumerate-(a)}
\end{corollary}

\begin{proof}
\ref{cor:kappaanalytic-1} implies \ref{cor:kappaanalytic-2} and \ref{cor:kappaanalytic-3} implies \ref{cor:kappaanalytic-4} because every closed set is effective \( \kappa+1 \)-Borel. 
Moreover, \ref{cor:kappaanalytic-2} implies \ref{cor:kappaanalytic-1} by Proposition~\ref{prop:borelsetsareanalytic}, and \ref{cor:kappaanalytic-4} implies \ref{cor:kappaanalytic-2} because \( \pre{\kappa}{\kappa} \times \pre{\kappa}{\kappa} \) is homeomorphic to \( \pre{\kappa}{\kappa} \) and the projection map is continuous. 
So it is enough to show that \ref{cor:kappaanalytic-1} implies \ref{cor:kappaanalytic-3}. 
Let \( f \colon C \to \pre{\kappa}{\kappa} \) be a continuous surjection onto \( A \) with \( C \subseteq \pre{\kappa}{\kappa} \) closed. 
Since \( \pre{\kappa}{\kappa} \) is an Hausdorff space, this implies that the graph \( F \) of \( f \) is closed in \( C \times \pre{\kappa}{\kappa} \), and hence also closed in 
\( \pre{\kappa}{\kappa} \times \pre{\kappa}{\kappa} \). 
As \( A = \PROJ ( F^{-1} ) \) the result is proved.
\end{proof}

One of the main uses of Corollary~\ref{cor:kappaanalytic} is that it allows us to use (a generalization of) the Tarski-Kuratowski algorithm (see e.g.~\cite[Appendix C]{Kechris:1995zt}) to establish that a given set \( A \subseteq \pre{\kappa}{\kappa} \) is \( \kappa \)-analytic by inspecting the ``logical form'' of its definition. 

The definition of \( \kappa \)-analyticity can be extended to subsets of an arbitrary subspace \( S \) of \( \pre{\kappa}{\kappa} \): \( A \subseteq S \) is \( \kappa \)-analytic (in \( S \)) if and only if there is a \( \kappa \)-analytic subset \( A' \) of \( \pre{\kappa}{\kappa} \) such that \( A = A' \cap S \). 
Notice that if \( S \) is an effective \( \kappa+1 \)-Borel subset of \( \pre{\kappa}{\kappa} \), then by Proposition~\ref{prop:borelsetsareanalytic} and the observation following it we get that \( A \subseteq S \) is \( \kappa \)-analytic in \( S \) if and only if it is \( \kappa \)-analytic in the whole \( \pre{\kappa}{\kappa} \). 
This allows us to have a natural definition of \( \kappa \)-analytic subset for any space which is effectively \( \kappa+ 1 \)-Borel isomorphic to some \( B \in \bB^\mathrm{e}_{\kappa+1}(\pre{\kappa}{\kappa},\tau_b) \), leading to the following definitions which generalize the case \( \kappa = \omega \).

\begin{definition}
Let \( \mathscr{B} \subseteq \pow ( X ) \) be an algebra on a nonempty set \( X \). 
We say that \( ( X , \mathscr{B} ) \) is a \markdef{standard Borel \( \kappa \)-space}\index[concepts]{Borel!standard Borel \( \kappa \)-space}\index[concepts]{space!standard Borel \( \kappa \)-space}%
\footnote{Our definition of a standard Borel \( \kappa \)-space slightly differs from the one introduced in~\cite[Definition 3.6]{Motto-Ros:2011qc}; however, the two definitions essentially coincide in the setup of~\cite{Motto-Ros:2011qc}, i.e.\ when assuming \( \AC \) and \( \kappa^{< \kappa} = \kappa \).} 
if there is a topology \( \tau \) on \( X \) such that \( \mathscr{B} = \bB^{\mathrm{e}}_{ \kappa + 1 }( X , \tau ) \) and \( ( X , \tau ) \) is homeomorphic to an effective \( \kappa + 1 \)-Borel subset of \( \pre{\kappa}{\kappa} \).
If the algebra \( \mathscr{B} \) is clear from the context, we say that \( X \) is a standard Borel \( \kappa \)-space.
\end{definition}

The collection of standard Borel \( \kappa \)-spaces is closed under effective \( \kappa+1 \)-Borel subsets, that is: if \( ( X , \mathscr{B} ) \) is a standard Borel \( \kappa \)-space then for every \( B \in \mathscr{B} \) the space \( ( B , \mathscr{B} \restriction B) \) is a standard Borel \( \kappa \)-space as well, where \( \mathscr{B} \restriction B \coloneqq \setof{B' \cap B }{ B' \in \mathscr{B} } \).
In particular, for every \( B \in \bB^{\mathrm{e}}_{\kappa + 1 }(\pre{\kappa}{\kappa}, \tau_b ) \) the space 
\[ 
( B , \bB^\mathrm{e}_{\kappa + 1 }( \pre{\kappa}{\kappa} , \tau_b ) \restriction B ) = ( B , \bB^\mathrm{e}_{\kappa + 1 } ( B , \tau_b ) ) 
\] 
is a standard Borel \( \kappa \)-space. 
Moreover, the product and the disjoint union of finitely many standard Borel \( \kappa \)-spaces are again standard Borel \( \kappa \)-spaces (where a product \( X \times X' \) is equipped with the product \( \mathscr{B} \otimes \mathscr{B}' \) of the algebras \( \mathscr{B} \) and \( \mathscr{B}' \) on \( X \) and \( X' \), and a disjoint union \( X \uplus X' \) is equipped with the corresponding union algebra \( \mathscr{B} \oplus \mathscr{B}' \)). 
Assuming enough choice (\( \AC_\kappa \) suffices), the same is true for products and unions of length \( \leq \kappa \).

Notice that Definitions~\ref{def:kappaanalyticinstandardBorel} and~\ref{def:kappaanalyticqo} are independent from the choice of the witness that \( X \) is a standard Borel \( \kappa \)-space. 

\begin{definition} \label{def:kappaanalyticinstandardBorel}
Let \( X \) be a standard Borel \( \kappa \)-space. 
A set \( A \subseteq X \) is \markdef{\( \kappa \)-analytic} if for any topology \( \tau \), any set \( B \in \bB^{\mathrm{e}}_{\kappa+1}(\pre{\kappa}{\kappa},\tau_b) \), and any homeomorphism \( f \colon ( X , \tau ) \to ( B, \tau_b \restriction B ) \) witnessing that \( X \) is standard Borel, the set \( f ( A ) \) is a \( \kappa \)-analytic subset of \( \pre{\kappa}{\kappa} \) (equivalently, of \( B \)).
\end{definition}

The above definition directly implies that analogues of Proposition~\ref{prop:borelsetsareanalytic} and Corollary~\ref{cor:kappaanalytic} hold in the broader context of standard Borel \( \kappa \)-spaces.

\begin{proposition}
Let \( ( X , \mathscr{B}) \) be a standard Borel \( \kappa \)-space. 
\begin{enumerate-(a)}
\item
Every \( B \in \mathscr{B} \) is \( \kappa \)-analytic.
\item
The following are equivalent for \( A \subseteq X \):
\begin{itemize}[leftmargin=1pc]
\item
\( A \) is \( \kappa \)-analytic;
\item 
\( A \) is a continuous image of an effective \( \kappa+1 \)-Borel subset of \( \pre{\kappa}{\kappa} \) (where continuity refers to any topology \(\tau\) on \( X \) witnessing that it is a standard Borel \( \kappa \)-space);
\item
\( A = \PROJ F \) for some closed \( F \subseteq X \times \pre{\kappa}{\kappa} \) (where \( X \) is endowed with any \(\tau\) as above);
\item 
\( A = \PROJ B \) for some \( B \in \mathscr{B} \otimes \bB^\mathrm{e}_{\kappa+1}(\pre{\kappa}{\kappa}, \tau_b) \).
\end{itemize}
\end{enumerate-(a)}
\end{proposition}

\begin{definition} \label{def:kappaanalyticqo}
A \markdef{\( \kappa \)-analytic quasi-order} \( S \) on a standard Borel \( \kappa \)-space \( X \) is a quasi-order on \( X \) which is \( \kappa \)-analytic as a subset of \( X \times X \). 
If moreover \( S \) is symmetric, then it is called a \markdef{\( \kappa \)-analytic equivalence relation} (on \( X \)).
\end{definition}

Examples of natural \( \kappa \)-analytic quasi-orders and equivalence relations are given in Section~\ref{subsec:spacesofmodels}.

\subsection{Spaces of type \( \kappa \) and spaces of codes}\label{subsec:spacesofmodels}

In this section the notion of a space of type \( \kappa \) is defined --- these are spaces which are homeomorphic to \( \pre{ \kappa }{2} \) in a canonical way.
This notion is introduced to ease the study of spaces of codes for \( \LL \)-structures of size \( \kappa \) for some finite relational language \( \LL \), for complete metric spaces of density character \( \kappa \), and for Banach spaces of density \( \kappa \). 
All spaces of type \( \kappa \) and their effective \( \kappa+1 \)-Borel subsets are also standard Borel \( \kappa \)-spaces when equipped with the algebra of their effective \( \kappa+ 1 \)-Borel subsets.

\subsubsection{Spaces of type \( \kappa \)} \label{subsubsec:spacesoftypekappa}
Let \( A \) be a set of size \( \kappa \). 
Then any bijection \( f \colon \kappa \to A \) induces a bijection between \( \pre{ \kappa }{2} \) and \( \pre{ A}{2} \), so the product topology, the \( \lambda \)-topology (for \( \omega \leq \lambda < \max ( \cof ( \kappa )^+ , \kappa ) \)), and the bounded topology can be copied on \( \pre{ A}{2} \), and are denoted with \( \tau_p ( \pre{ A}{2} ) \), \(\tau_{ \lambda } ( \pre{ A}{2} ) \), and \( \tau_b ( \pre{ A}{2} ) \), respectively.
The bases for these topologies are given by
\begin{align*}
\mathcal{B}_ p ( \pre{ A } {2} ) & \coloneqq \setof{\Nbhd^A_s}{\card{s} < \omega }
\\
\mathcal{B}_ \lambda ( \pre{ A } {2} ) & \coloneqq \setof{\Nbhd^A_s}{\card{s} < \lambda }
\\
\mathcal{B}_ b ( \pre{ A } {2} ) & \coloneqq \setof{\Nbhd^A_s}{ \EXISTS{ \alpha < \kappa } \left ( \appl{ f }{ \alpha } = \dom s \right ) }
\end{align*}
where 
\[ 
\Nbhd^A_s \coloneqq \setof{ x \in \pre{ A}{2} }{ s \subseteq x }
\] 
for \( s \) a partial function from \( A \) to \( 2 \).
Note that \( \tau_p ( \pre{ A}{2} ) = \tau_ \omega ( \pre{ A}{2} ) \), since \( \mathcal{B}_p ( \pre{ A } {2} ) = \mathcal{B}_ \omega ( \pre{ A } {2} ) \).
As in Remark~\ref{rmk:product}\ref{rmk:product-a}, the collection of all 
\begin{equation} \label{eq:subbasisproduct} 
\widetilde{\Nbhd}^A_{a,i} \coloneqq \setof{x \in \pre{ A}{2}}{x ( a ) = i} \qquad (a \in A , i \in \set{ 0,1 } )
\end{equation} 
is a subbasis (generating \( \mathcal{B}_p ( \pre{ A}{2} ) \)) for \( \tau_p ( \pre{ A}{2}) \).

The definitions of \( \tau_\lambda ( \pre{ A}{2} ) \) and \( \mathcal{B}_\lambda ( \pre{ A } {2} ) \) (which includes the case of the product topology) are independent of the chosen \( f \). 
The situation for \( \tau_b ( \pre{ A}{2} ) \) is rather different --- its definition is again independent of \( f \) when \( \kappa \) is regular, but this is no more true when \( \kappa \) is singular. 
Moreover, the canonical basis \( \mathcal{B}_b ( \pre{ A } {2} ) \) always depends on \( f \), even when \( \kappa \) is regular. 
This is an unpleasant situation; however, in our applications this will not be an issue, as there will always be a \emph{canonical} bijection between \( A \) and \( \kappa \). 
To illustrate this, let us consider two representative examples.

\begin{examples}\label{xmp:canonicalbijection}
\begin{enumerate-(A)}
\item \label{xmp:canonicalbijection-1}
Consider the set \( A \coloneqq \pre{ < \omega}{2} \times \kappa \) (the topological space \( ( \pre{ A } {2} , \tau_p ) \) plays an important role in Section~\ref{sec:invariantlyuniversal}). 
Let \( \theta \colon \pre{ < \omega}{2} \to \omega \) be the unique isomorphism between \( ( \pre{ < \omega}{2} , \preceq) \) and \( ( \omega , \leq ) \), where \( \preceq \) is defined by
\[ 
u \preceq v \iff \lh u < \lh v \vee ( \lh u = \lh v \wedge u \leqlex v)
 \] 
with \( \leqlex \) the usual lexicographical order. 
Then
\begin{align*}
&f \colon \pre{ < \omega}{ 2} \times \kappa \to \kappa , & ( u , \alpha ) \mapsto \omega \cdot \alpha + \theta ( u ) 
\end{align*}
can be taken to be a standard bijection between \( A \) and \( \kappa \), as it orders \( A \) antilexicographically. 
With this bijection the notions of ``boundness'' and ``initial segment'' on \( A \) become natural and unambiguous. 

\item \label{xmp:canonicalbijection-2}
Recall from~\eqref{eq:Hessenberg} the standard pairing function for ordinals \( \op{\cdot }{\cdot } \colon \On \times \On \to \On \).
For \( n \geq 2 \) define the bijections \( f_n \colon \pre{ n}{\kappa} \to \kappa\) by setting
\begin{align*}
f_2 ( \alpha_0 , \alpha_1 ) &\coloneqq \op{ \alpha_0 }{ \alpha_1 }, 
\\
f_{n + 1} ( \alpha_0, \dotsc , \alpha_n) & \coloneqq \op{ f_n ( \alpha_0 , \dotsc , \alpha_{n - 1 } ) }{ \alpha_n }.
\end{align*} 
Each \( f_n \) can be considered as the standard bijection between \( A \coloneqq \pre{ n}{\kappa} \) and \( \kappa \).
\end{enumerate-(A)} 
\end{examples}

The next definition tries to capture the content of Examples~\ref{xmp:canonicalbijection}.

\begin{definition} \label{def:spaceoftypekappa}
Let \( \kappa \) be an infinite cardinal. 
A \markdef{space of type \( \kappa \)}\index[concepts]{space!of type \( \kappa \)} is a set of the form \( \mathcal{X} = \prod_{i \in I} \pre{ A_i}{2} \), where \( I \) is a finite set and each \( A_i \) is a set of cardinality \( \kappa \).

The bounded topology \( \tau_b ( \mathcal{X} ) \) on \( \mathcal{X} \) is the product of the bounded topologies \( \tau_b ( \pre{ A_i}{2}) \) on \( \pre{ A_i}{2} \), and is generated by the canonical basis 
\[ 
\mathcal{B}_b^{\mathcal{X}} \coloneqq \setof{ \textstyle\prod_{i \in I} \Nbhd^{A_i}_{s_i}}{\Nbhd^{A_i}_{s_i} \in \mathcal{B}_b ( \pre{ A_i}{2} ) \text{ for all } i \in I}.
\] 
Similarly, \( \tau_p ( \mathcal{X} ) \) and \( \tau_\lambda ( \mathcal{X} ) \) are defined as the product of the corresponding topologies on each factor, and their bases \( \mathcal{B}_p ( \mathcal{X} ) \) and \( \mathcal{B}_\lambda ( \mathcal{X} ) \) are defined as the products of the corresponding bases.
\end{definition}

All definitions, observations, and results concerning \( \pre{ \kappa}{2} \) considered in Sections~\ref{sec:topologies}--\ref{sec:generalizedBorelfunctions} can be applied to an arbitrary space \( \mathcal{X} \) of type \( \kappa \) as well, including the following:
\begin{itemize}[leftmargin=1pc]
\item
The collection \( \bB_\alpha ( \mathcal{X} , \tau_* ) = \bB_\alpha ( \tau_* ) \) of \emph{\( \alpha \)-Borel subsets} of \( \mathcal{X} \), where \( \tau_* \) is one of the topologies \( \tau_p ( \mathcal{X} ) \), \( \tau_\lambda ( \mathcal{X} ) \), or \( \tau_b ( \mathcal{X} ) \) (Definition~\ref{def:kappaBorel}), together with its effective counterpart \( \bB^\mathrm{e}_\alpha ( \mathcal{X},\tau_* ) \) (Section~\ref{subsubsec:Borelcodes}).
\item
The \emph{(weakly) \( \alpha \)-Borel functions} \( f \colon \mathcal{Y} \to \mathcal{Z} \), where \( \mathcal{Y} \) and \( \mathcal{Z} \) are arbitrary spaces of type \( \lambda \) and \( \mu \), respectively (Definition~\ref{def:weaklyBorel}).
\item
The \emph{\( \bGamma \)-in-the-codes function} \( f \colon \pre{ \omega}{2} \to \mathcal{X} \) for \( \bGamma \) a suitable boldface pointclass (Definition~\ref{def:Gammainthecodes}).
\item
The fact that the collection of all \( \tau_b ( \mathcal{X} ) \)-clopen subsets is a \( \cof ( \kappa ) \)-algebra on \( \mathcal{X} \) (Proposition~\ref{prop:topologicalproperties}\ref{prop:topologicalproperties-i}).
\item
The fact that assuming \( \AC \), 
if \( 2^{<\kappa} = \kappa \) then \( \bB_{\kappa + 1} ( \tau_p ( \mathcal{X} ) ) = \bB_{\kappa + 1} ( \tau_b ( \mathcal{X} ) ) \) (Corollary~\ref{cor:borel}), and therefore \( \bB_{\kappa + 1} ( \tau_\lambda ( \mathcal{X} ) ) = \bB_{\kappa + 1} ( \tau_b ( \mathcal{X} ) ) \) for every \( \omega \leq \lambda < \max ( \cof ( \kappa) ^+ , \kappa ) \).
\end{itemize}

\begin{remark} \label{rmk:partialBorelfunctions}
We often consider (weakly) \(\alpha\)-Borel functions \( f \colon A \to B \) (for suitable ordinals \(\alpha\)) where \( A \) and \( B \) are arbitrary subspaces of two spaces \( \mathcal{Y} \) and \( \mathcal{Z} \) of type \( \lambda \) and \( \mu \), respectively, each endowed with \( \tau_b \). 
By this we mean that \( f \) is (weakly) \( \mathcal{S} \)-measurable when \( B \) is endowed with the relative topologies inherited from \( \mathcal{Z} \) and \( \mathcal{S} \) is the algebra of subsets of \( A \) consisting of the traces on \( A \) of the \( \alpha \)-Borel subsets of \( \mathcal{Y} \). 
Notice that: 
\begin{enumerate-(i)}
\item
 When \( A \in \bB_{\alpha} ( \mathcal{Y} ) \), the notion of a (weakly) \( \alpha \)-Borel function \( f \colon A \to \mathcal{Z} \) is unambiguous since a set \( C \subseteq A \) is \( \alpha \)-Borel in \( A \) if and only if it is \( \alpha \)-Borel as a subset of the entire space \( \mathcal{Y} \). 
In fact, a function \( f \colon A \to \mathcal{Z} \) is (weakly) \( \alpha \)-Borel if and only if \( f = g \restriction A \) for some (weakly) \( \alpha \)-Borel function \( g \colon \mathcal{Y} \to \mathcal{Z} \).
\item
A function \( f \colon \mathcal{Y} \to B \) is (weakly) \( \alpha \)-Borel if and only if it is (weakly) \( \alpha \)-Borel as a function between \( \mathcal{Y} \) and \( \mathcal{Z} \).
\end{enumerate-(i)} 
\end{remark}

Any space \( \mathcal{X} \) of type \( \kappa \) is, by definition, homeomorphic to \( \pre{\kappa}{2} \), which is a closed (and hence effective \( \kappa+1 \)-Borel) subset of \( \pre{\kappa}{\kappa} \)\, as long as all these spaces are endowed with the same kind of topology. 
Thus \( \mathcal{X} \) equipped with the algebra of its effective \( \kappa+1 \)-Borel subsets with respect to \( \tau_b \) is also a standard Borel \( \kappa \)-space, and therefore we can consider the notions of a \( \kappa \)-analytic subset of \( \mathcal{X} \) and of a \( \kappa \)-analytic quasi-order on \( \mathcal{X} \) (or on one of its effective \( \kappa+1 \)-Borel subsets) as in Definitions~\ref{def:kappaanalyticinstandardBorel} and~\ref{def:kappaanalyticqo}.

\subsubsection{Space of codes for (\( \LL \)-)structures of size \( \kappa \)}\label{subsubsec:spacesofmodels}
For \( I \) a finite set, let 
\[ 
\LL = \setof{R_i }{ i \in I } 
\] 
be a relational%
\footnote{To simplify the presentation we retreat to relational signatures. 
This is not restrictive since any \( n \)-ary function may be identified with its graph, which is a relation of arity \( n+1 \).}
 signature.
For the sake of definiteness assume that \( I \in \omega \) and that \( n_i \) is the arity of \( R_i \), for \( i \in I \).
A \markdef{structure} or \markdef{model} for \( \LL \), sometimes called simply \( \LL \)-structure, is an object of the form \( \mathcal{A} = \seq{A ; R^{ \mathcal{A}}_i }_{i \in I} \), where \( A \) is a nonempty set and \( R^{ \mathcal{A}}_i \) is the interpretation of the symbol \( R_i \) in \( \mathcal{A} \). 
If \( \emptyset \neq B \subseteq A \) then we denote by \( \mathcal{A} \restriction B \) the restriction of \( \mathcal{A} \) to its subdomain \( B \), i.e.~\( \mathcal{A} \restriction B = \seq{B; R^{ \mathcal{A}}_i \cap \pre{ n_i}{B}}_{i \in I} \). 
When there is no danger of confusion, with abuse of notation the structure \( \mathcal{A} \) is identified with its domain \( A \).
Since the nature of the elements of \( A \) is irrelevant, a model of size \( \kappa \) is taken to have domain \( \kappa \), so that each \( R^{ \mathcal{A}}_i \) can be identified with its characteristic function \( \pre{ n_i}{ \kappa} \to 2 = \set{0 , 1} \).
Therefore any model of size \( \kappa \) can be identified, up to isomorphism, with a map 
\[ 
\textstyle\bigcup_{i \in I } \set{ i } \times \pre{ n_i}{\kappa} \to 2 , 
\] 
and hence 
\begin{equation}\label{eq:spaceofmodelsofsizekappa}
 \Mod^\kappa_{\LL} \coloneqq \textstyle\prod_{i \in I} \pre{ ( \pre{ n_i}{\kappa})}{2} \index[symbols]{Mod@\( \Mod^\kappa_{\LL} \)}
\end{equation}
can be regarded as the \markdef{space of} (codes for) \markdef{all \( \LL \)-structures of size \( \kappa \)} (up to isomorphism). 
We also set
\[ 
\Mod^{<\kappa}_\LL \coloneqq \textstyle\bigcup_{\lambda < \kappa} \Mod^\lambda_\LL \quad \text{and} \quad \Mod^\infty_{\LL} \coloneqq \textstyle\bigcup_{ \kappa \in \Cn } \Mod^ \kappa_{\LL} .
\] 
For example, if \( \LL \) is a relational language consisting of just one relational symbol of arity \( n \), then \( \Mod^\kappa_\LL = \pre{ ( \pre{ n}{\kappa} ) }{2} \), and hence by Example~\ref{xmp:canonicalbijection}\ref{xmp:canonicalbijection-2} it can be topologized by \( \tau_p \), \( \tau_\lambda \) (for \( \omega \leq \lambda < \max ( \cof ( \kappa )^ + , \kappa ) \)), or \( \tau_b \). 
For richer languages \( \LL = \setof{R_i}{ i \in I } \), the space of models \( \Mod^\kappa_\LL \) in~\eqref{eq:spaceofmodelsofsizekappa} is just a finite product of spaces of the form \( \pre{ ( \pre{ n}{\kappa} ) }{2} \): thus it is a space of type \( \kappa \), and hence a standard Borel \( \kappa \)-space.

Strictly speaking an \( X \in \Mod^ \kappa _\LL \) is a \emph{function}, but it is more convenient to think of it as an \( \LL \)-structure 
\[
X = \seq{ \kappa ; R_i^X }_{i \in I} .
\] 
The \markdef{embeddability} relation \( \embeds \)\index[symbols]{49@\( \embeds \)} on \( \Mod^ \kappa _ \LL \) is given by
\begin{multline} \label{eq:embeds}
X \embeds Y \IFF \EXISTS{g \in \pre{ \kappa }{( \kappa )} }\FORALL{i\in I} \FORALL{ \seq{ x_1 , \dotsc , x_{ n_i } } \in \pre{ n_i}{\kappa}}
\\
\bigl [ \seq{ x_1 , \dotsc , x_{ n_i } } \in R_i^X \IFF \seq{ g ( x_1 ) , \dotsc , g ( x_{ n_i } ) } \in R_i^Y \bigr ]
\end{multline}
where \( \pre{ \kappa }{( \kappa )} \) is the set of injective functions from \( \kappa \) into \( \kappa \) (see~\eqref{eq:injections}). 
The equivalence relation induced by the quasi-order \( \embeds \) is the \markdef{bi-embeddability} relation \( \biembeds \)\index[symbols]{49@\( \biembeds \)}. 
If \( \pre{ \kappa }{( \kappa )} \) in~\eqref{eq:embeds} is replaced by \( \Sym ( \kappa ) \), the group of all bijections from \( \kappa \) to \( \kappa \), the \markdef{isomorphism} relation \( \cong \) is obtained.
Thus \( \cong \) can be seen as induced by the continuous action
\begin{equation} \label{eq:action}
\Sym ( \kappa ) \times \Mod^ \kappa _\LL \to \Mod^ \kappa _\LL , \qquad ( g , X ) \mapsto g . X
\end{equation} 
where \( g . X \in \Mod^ \kappa _\LL \) is defined by
\begin{equation} \label{eq:action2}
R_i^{g . X} \coloneqq \setof{ \seq{ y_1 , \dotsc , y_{ n_i } } \in \pre{ n_i}{\kappa}}{ \seq{ g^{- 1} ( y_1 ) , \dotsc, g^{-1} ( y_{n_i } ) } \in R_i^X } \qquad ( i \in I ).
\end{equation}
The embeddability relation \( \embeds \) on \( \Mod^\kappa_\LL \) is an example of a \( \kappa \)-analytic quasi-order. 
To see this, observe that \( \pre{\kappa}{(\kappa)} \) is a closed (and hence effective \( \kappa+1 \)-Borel) subset of \( \pre{\kappa}{\kappa} \), so that \( \embeds \) is the projection on \( \Mod^\kappa_\LL \times \Mod^\kappa_\LL \) of a closed subset of \( \Mod^\kappa_\LL \times \Mod^\kappa_\LL \times \pre{\kappa}{\kappa} \) by~\eqref{eq:embeds}. 
It follows that \( \biembeds \) is a \( \kappa \)-analytic equivalence relation on the standard Borel \( \kappa \)-space \( \Mod^\kappa_\LL \).
Similarly, using the fact that \( \Sym ( \kappa ) \) is an effective \( \kappa+1 \)-Borel subset of \( \pre{\kappa}{\kappa} \) (in fact: an effective intersection of \( \kappa \)-many open sets) one sees that \( \cong \) is a \( \kappa \)-analytic equivalence relation on \( \Mod^\kappa_\LL \). 

\subsubsection{Space of codes for complete metric spaces of density character \( \kappa \)}\label{subsubsec:spacesofmetricspaces}

In the classical case \( \kappa = \omega \), there are essentially two ways for coding separable complete metric spaces, usually called \markdef{Polish metric spaces},\index[concepts]{Polish space!metric} as elements of a standard Borel space. 

The first, and more common one~\cite{Clemens:2001pu,Gao:2003qw,Louveau:2005cq}, uses the fact that the (separable) Urysohn space \( \mathbb{U} \) is universal for this class, that is to say: \( \mathbb{U} \) is itself a Polish metric space (so that all its closed subsets are Polish metric spaces as well), and every Polish metric space is isometric to some closed subset of \( \mathbb{U} \). 
Thus the space \( F(\mathbb{U}) \) of all closed subsets of \( \mathbb{U} \) endowed with its Effros-Borel structure, which is standard Borel (see e.g.\ \cite[Section 12.C]{Kechris:1995zt}), may be regarded as a space of codes for all separable complete metric spaces. 
If this approach is to be generalized to an uncountable \( \kappa \), a space which is universal for complete metric spaces of density character \( \kappa \) must be constructed. 
By~\cite{Katetov:1986tb} analogues of the Urysohn space for larger density characters may be obtained only assuming \( \AC \) and for cardinals \( \kappa \) satisfying \( \kappa^{<\kappa} = \kappa \). 
Thus this technique for coding metric spaces cannot be used here, since we want to study also choice-less models (such as models of \( \AD \)), and even in the \( \AC \) context we are interested in cardinals smaller than the continuum. 
To the best of our knowledge, there are no other kinds of universal spaces for complete metric spaces of density character \( \kappa \) in the literature, so we are forced to drop this approach.

The second way to code Polish metric spaces is to identify each of them with any of its dense subspaces, so that the original space may be recovered, up to isometry, as the completion of such a subspace (see e.g.\ \cite{Vershik:1998fk,Clemens:2012hg}). 
Fortunately, this approach does generalize in \( \ZF \) to any infinite cardinal \( \kappa \), naturally yielding to the set of codes \( \MMM_\kappa \) described below. 
Notice that the space of codes \( \MMM_\kappa \), being an effective \( \kappa+1 \)-Borel subset of a space of type \( \kappa \), carries a natural topology \( \tau_b \), and its effective \( \kappa+1 \)-Borel structure turns it into a standard Borel \( \kappa \)-space. 

Let \( \QQ^+ \) be the set of positive rational numbers, and let \( \mathcal{X} \) be the space of type \( \kappa \) defined by \( \mathcal{X} \coloneqq \pre{\kappa \times \kappa \times \QQ^+}{2} \). 
Given a complete metric space \( ( M , d ) \) of density character \( \kappa \) and a dense subset \( D = \setof{m_\alpha}{\alpha < \kappa} \) of it, we can identify \( M \) with the unique element \( x_M \in \mathcal{X}\) such that for all \( \alpha, \beta < \kappa \) and \( q \in \QQ^+ \)
\begin{equation} \label{eq:defx_M}
x_M(\alpha, \beta, q) = 1 \IFF d_M(m_\alpha,m_\beta) < q.
\end{equation}
In fact, \( M \) is isometric to the completion of the metric space \( ( \kappa , d_{x_M} ) \) where \( d_{x_M} ( \alpha , \beta ) \coloneqq \inf \setof{q \in \QQ^+}{x_M ( \alpha , \beta , q ) = 1} \) for \( \alpha,\beta < \kappa \), so that \( d_{x_M}(\alpha,\beta) = d_M(m_\alpha,m_\beta) \) for all \( \alpha, \beta < \kappa \). Consider now the space \( \MMM_\kappa \subseteq \mathcal{X} \)\index[symbols]{Mkappa@\( \MMM_\kappa \)} consisting of those \( x \in \pre{\kappa \times \kappa \times \QQ^+}{2} \) satisfying the following conditions:
\begin{align*}
\FORALL{ \alpha,\beta < \kappa} \FORALL{ q , q ' \in \QQ^+} & \left [ q \leq q' \IMPLIES x ( \alpha , \beta , q ) \leq x ( \alpha , \beta , q' ) \right ] 
\\
\FORALL{ \alpha , \beta < \kappa} \EXISTS{q \in \QQ^+} & \left [ x ( \alpha , \beta , q ) = 1\right ] 
\\
\FORALL{ \alpha < \kappa} \FORALL{q \in \QQ^+} & \left[ x ( \alpha , \alpha , q ) = 1 \right] 
\\
\FORALL{ \alpha < \beta < \kappa} \EXISTS{ q \in \QQ^+} & \left [ x ( \alpha , \beta , q ) = 0 \right ] 
\\
\FORALL{ \alpha , \beta < \kappa} \FORALL{ q \in \QQ^+} & \left [ x ( \alpha , \beta , q ) = 1 \IFF x ( \beta , \alpha , q ) = 1 \right] 
\\
\FORALL{ \alpha , \beta , \gamma < \kappa} \FORALL{q,q' \in \QQ^+} & \left[ x(\alpha,\beta,q) = 1 \AND x ( \beta , \gamma , q ' ) = 1 \IMPLIES x ( \alpha , \gamma , q + q ' ) = 1 \right ] 
\\
\FORALL{\alpha<\kappa} \EXISTS{\beta < \kappa} \EXISTS{ q \in \QQ^+} \FORALL{ \gamma < \alpha } & [ x ( \gamma , \beta , q ) = 0 ].
\end{align*}
The first six conditions are designed so that given any \( x \in \MMM_\kappa \), the (well-defined) map \( d_x \colon \kappa \times \kappa \to \R \) defined by setting 
\[ 
d_x(\alpha,\beta) \coloneqq \inf \setof{q \in \QQ^+}{x(\alpha,\beta,q) = 1} 
 \] 
is a metric on \( \kappa \); denote by \( M_x \) the completion of \( ( \kappa , d_x ) \), and notice that the last condition ensures that \( M_x \) has density character \( \kappa \). 
It is straightforward to check that the code \( x_M \) from~\eqref{eq:defx_M} of any complete metric space \( M \) of density character \( \kappa \) belongs to \( \MMM_\kappa \), and is such that \( M \) is isometric to \( M_{x_M} \); conversely, for each \( x \in \MMM_\kappa \) the space \( M_x \) is a complete metric space of density character \( \kappa \). 
Moreover, the explicit definition given above shows that \( \MMM_\kappa \in \bB^{\mathrm{e}}_{\kappa+1}( \mathcal{X}, \tau_b ) \), so that \( \MMM_\kappa \) is a standard Borel \( \kappa \)-space. 
Thus we can regard \( \MMM_\kappa \subseteq \pre{\kappa \times \kappa \times \QQ^+}{2} \), endowed with the inherited topologies and the corresponding (effective) \( \kappa +1 \)-Borel structure, as the \markdef{space of} (codes for) \markdef{all complete metric spaces of density character \( \kappa \)} (up to isometry).

The \markdef{isometric embeddability} relation \( \sqsubseteq^i \)\index[symbols]{51@\( \sqsubseteq^i \)} on \( \MMM_{\kappa} \) is given by
\[ 
x \sqsubseteq^i y \IFF \text{there is a metric-preserving map from } M_x \text{ into }M_y,
 \] 
while the \markdef{isometry} relation \( \cong^i \) on \( \MMM_\kappa \) is given by
\[ 
x \cong^i y \IFF \text{there is a metric-preserving bijection between } M_x \text{ and }M_y.
 \] 
Notice that for every \( x , y \in \MMM_\kappa \) one has 
\[ 
x \sqsubseteq^i y \IFF \text{ there is a metric-preserving map } i \colon ( \kappa , d_x ) \to M_y. 
\] 
This allows us to check, using the Tarski-Kuratowski algorithm and some standard computations, that the relation \( \sqsubseteq^i \) is a \( \kappa \)-analytic quasi-order on \( \MMM_\kappa \), and a similar observation shows that \( \cong^i \) is a \( \kappa \)-analytic equivalence relation on \( \MMM_\kappa \).
(For a blueprint of such computations, see~\cite[Lemma 4]{Clemens:2012hg}.)

We consider some natural subclasses of \( \MMM_\kappa \), such as
\[ 
\DDD_\kappa \coloneqq \setof{x \in \MMM_\kappa}{M_x \text{ is discrete}}\index[symbols]{Dkappa@\( \DDD_\kappa \)}
 \] 
or
\[ 
\UUU_\kappa \coloneqq \setof{x \in \MMM_\kappa}{M_x \text{ is ultrametric}}. \index[symbols]{Ukappa@\( \UUU_\kappa \)}
 \] 
Notice however that not all these subclasses are standard Borel: for example, \( \UUU_\kappa \) is a standard Borel \( \kappa \)-space (since it is a closed subset of \( \MMM_\kappa \)), while it can be shown that \( \DDD_\kappa \) is not.

\subsubsection{Space of codes for Banach spaces of density \( \kappa \)}\label{subsubsec:spacesofBanachspaces}
In analogy with what was done in Section~\ref{subsubsec:spacesofmetricspaces} for complete metric spaces of density character \( \kappa \), we code Banach spaces of density \( \kappa \) by identifying each of them with any of its dense subspaces closed under rational%
\footnote{To simplify the presentation, in this paper we focus on \emph{real} Banach spaces. 
However, by replacing \( \QQ \) with \( \QQ + i \QQ \) one can extend our results to complex Banach spaces of density \( \kappa \).} 
linear combinations; this gives rise to an effective \( \kappa +1 \)-Borel subset \( \BBB_\kappa \) of a space of type \( \kappa \), which can be regarded as the standard Borel \( \kappa \)-space of all Banach spaces of density \( \kappa \).

Let \( B \) be a (real) Banach space of density \( \kappa \) with norm \( \| \cdot \|_B \), and let \( D = \setof{b_\alpha}{\alpha<\kappa} \) be a dense subset of \( B \) which is also closed under rational linear combinations. 
Without loss of generality, we can assume that \( b_0 \) is the zero vector of \( B \) (so that \( b_0 \) is the unique element of \( D \) with \( B \)-norm \( 0 \)). 
Then we can identify \( B \) with an element \( x_B = ( x_B^+ , x_B^\QQ , x_B^{\| \cdot \| } ) \) of the space \( \mathcal{X} \coloneqq \pre{\kappa \times \kappa \times \kappa}{2} \times \pre{\kappa \times \QQ \times \kappa}{2} \times \pre{\kappa \times \QQ^+}{2} \) by setting for \( \alpha,\beta,\gamma < \kappa \), \( p \in \QQ \), and \( q \in \QQ^+ \)
\begin{equation} \label{eq:x_B}
\begin{aligned} 
x_B^+(\alpha,\beta,\gamma) = 1 \IFF & b_\alpha + b_\beta = b_\gamma \\
x_B^\QQ (\alpha,p,\beta) = 1 \IFF & p \cdot b_\alpha = b_\beta \\
x_B^{\|\cdot\|} (\alpha, q) = 1 \IFF & \| b_\alpha \|_B < q.
\end{aligned}
\end{equation}
The function \( x_B \) codes all the necessary informations to retrieve the normed vector space structure of \( D \), and hence of the whole \( B \).
This suggests to consider the space \( \BBB_\kappa \subseteq \mathcal{X} \)\index[symbols]{BX@\( \BBB_\kappa \)} consisting of those \( x = ( x^+, x^\QQ , x^{\|\cdot\|} ) \in \mathcal{X} \) satisfying the conditions in Table~\ref{tab:Banachspaces}.
\begin{table}
 \centering
 \begin{tabular}{@{} l @{}} %
 \toprule
\( \FORALL{\alpha , \beta < \kappa} \EXISTS{! \gamma < \kappa} [ x^+ ( \alpha , \beta , \gamma ) = 1 ] \)
\\
\( \FORALL{ \alpha < \kappa} \FORALL{p \in \QQ} \EXISTS{! \beta < \kappa} [ x^\QQ ( \alpha , p , \beta ) = 1 ] \)
\\
\( \FORALL{ \alpha , \beta , \gamma , \delta , \epsilon , \zeta < \kappa} [ x^+ ( \beta , \gamma , \delta ) = 1 \AND x^+ ( \alpha , \delta , \epsilon ) = 1 \AND x^+( \alpha , \beta , \zeta ) = 1 \IMPLIES x^+ ( \zeta , \gamma , \epsilon ) = 1] \)
\\
\( \FORALL{ \alpha , \beta , \gamma < \kappa} [x^+ ( \alpha , \beta , \gamma ) = 1 \IFF x^+ ( \beta , \alpha , \gamma ) = 1] \)
\\
\( \FORALL{\alpha < \kappa} [x^+ ( \alpha , 0 ,\alpha) = 1] \)
\\
\( \FORALL{\alpha< \kappa} \EXISTS{\beta < \kappa} [x^+(\alpha , \beta , 0) = 1] \)
\\
\( \FORALL{\alpha , \beta , \gamma < \kappa} \FORALL{p , p' \in \QQ} [x^\QQ ( \alpha , p , \beta ) = 1 \AND x^\QQ ( \beta , p' , \gamma ) = 1 \IMPLIES x^\QQ ( \alpha , pp' , \gamma ) = 1] \)
\\
\( \FORALL{\alpha < \kappa} [x^\QQ( \alpha , 1 , \alpha) = 1] \)
\\
\( \FORALL{\alpha , \beta , \gamma , \delta , \epsilon , \zeta < \kappa} \FORALL{p \in \QQ} [ x^+ ( \alpha , \beta , \gamma ) = 1 \AND x^\QQ ( \gamma , p , \delta ) = 1 \AND x^\QQ ( \alpha , p , \epsilon ) = 1 \AND x^\QQ ( \beta , p , \zeta ) = 1 \)
\\
\qquad\qquad\( {}\IMPLIES x^+ ( \epsilon , \zeta , \delta ) = 1] \)
\\
\( \FORALL{\alpha , \beta , \gamma , \delta < \kappa} \FORALL{p , p' \in \QQ^+} [x^\QQ ( \alpha , p + p' , \beta) = 1 \AND x^\QQ( \alpha , p , \gamma) = 1 \AND x^\QQ ( \alpha , p' , \delta ) = 1 \)
\\
\qquad\qquad\( {} \IMPLIES x^+( \gamma , \delta , \beta ) = 1] \) 
\\
\( \FORALL{\alpha<\kappa} \FORALL{q , q' \in \QQ^+} [ q \leq q' \IMPLIES x^{\|\cdot\|}( \alpha , q ) \leq x^{\|\cdot\|}( \alpha , q' ) ] \)
\\
\( \FORALL{\alpha< \kappa} \EXISTS{q \in \QQ^+} [ x^{\|\cdot\|}( \alpha , q) = 1] \)
\\
\( \FORALL{0 < \alpha < \kappa} \EXISTS{q \in \QQ^+} [ x^{\|\cdot\|} ( \alpha , q) =0 ] \)
\\
\( \FORALL{q \in \QQ^+} [x^{\|\cdot\|} ( 0 , q ) = 1] \)
\\
\( \FORALL{\alpha , \beta < \kappa} \FORALL{p \in \QQ} \FORALL{q \in \QQ^+} [x^\QQ ( \alpha , p , \beta ) = 1 \IMPLIES ( x^{ \| \cdot \|}( \beta , q ) = 1 \IFF x^{\|\cdot\|} ( \alpha , q / | p | ) = 1)] \) 
\\
\( \FORALL{\alpha , \beta , \gamma< \kappa}\FORALL{q , q' \in \QQ^+} [x^+( \alpha , \beta , \gamma) = 1 \AND x^{\|\cdot\|} ( \alpha , q ) = 1 \AND x^{\|\cdot\|}( \beta , q' ) = 1 \)
\\
\qquad\qquad\( {}\IMPLIES x^{\|\cdot\|}( \gamma , q + q' ) = 1] \) 
\\
\( \FORALL{\alpha< \kappa} \EXISTS{\beta < \kappa} \EXISTS{q \in \QQ^+} \FORALL{\gamma<\alpha} \FORALL{\delta , \epsilon < \kappa} [x^\QQ( \gamma , -1 , \delta ) = 1 \AND x^+ ( \beta , \delta , \epsilon ) = 1 \IMPLIES x^{\| \cdot \|}( \epsilon , q ) = 0 ] \)
\\
 \bottomrule
 \\
\end{tabular}
 \caption{The conditions defining \( \BBB_\kappa \).}
 \label{tab:Banachspaces}
\end{table}
The meaning of the conditions in this table is clear: given \( x \in \BBB_\kappa \), we can consider the normed \( \QQ \)-vector space \( D_x \) on \( \kappa \) equipped with the (well-defined) operations \( +_x \) and \( \cdot_x \) and the (well-defined) norm \( \| \cdot \|_x \) obtained by setting for \( \alpha,\beta< \kappa \) and \( p \in \QQ \)
\begin{align*}
{\alpha +_x \beta = \gamma} \IFF & {x^+(\alpha,\beta,\gamma) = 1} \\
{p \cdot_x \alpha = \beta} \IFF & {x^\QQ ( \alpha , p , \beta ) = 1} \\
\| \alpha \|_x \coloneqq & \inf \setof{q \in \QQ^+}{x^{\|\cdot\|}(\alpha,q) = 1 }.
\end{align*}
The Banach space obtained by completing the norm \( \| \cdot \|_x \) is denoted by \( B_x \). 
Then for every Banach space \( B \) of density \( \kappa \) we get that its code \( x_B \in \mathcal{X} \) defined in~\eqref{eq:x_B} belongs to \( \BBB_\kappa \) (and is such that \( B \) and \( B_{x_B} \) are linearly isometric), and, conversely, for every \( x \in \BBB_\kappa \) the space \( B_x \) is a Banach space of density \( \kappa \). 
Moreover, the explicit definition from Table~\ref{tab:Banachspaces} shows that \( \BBB_\kappa \) is an effective \( \kappa+ 1 \)-Borel subset of \( \mathcal{X} \), and hence it inherits from it the topology \( \tau_b \) and the corresponding (effective) \( \kappa+1 \)-Borel structure, which turns it into a standard Borel \( \kappa \)-space. 
Thus we can regard \( \BBB_\kappa \) as the \markdef{space of} (codes for) \markdef{all Banach spaces of density \( \kappa \)} (up to linear isometry).

The \markdef{linear isometric embeddability} relation \( \sqsubseteq^{li} \)\index[symbols]{52@\( \sqsubseteq^{li} \), \( \cong^{li} \)} on \( \BBB_{\kappa} \) is given by
\[ 
x \sqsubseteq^{li} y \IFF \text{there is a linear norm-preserving map from } B_x \text{ into }B_y,
 \] 
while the \markdef{linear isometry} relation \( \cong^{li} \) on \( \BBB_\kappa \) is given by
\[ 
x \cong^{li} y \IFF \text{there is a linear norm-preserving bijection between } B_x \text{ and \( B_y \).}
 \] 
For every \( x , y \in \BBB_\kappa \) one has that \( x \sqsubseteq^{li} y \) if and only if there is a linear norm-preserving map \( i \colon D_x \to B_y \).
Using the Tarski-Kuratowski algorithm, the relation \( \sqsubseteq^{li} \) is a \( \kappa \)-analytic quasi-order on \( \BBB_\kappa \), and \( \cong^{li} \) is a \( \kappa \)-analytic equivalence relation on \( \BBB_\kappa \).

\section{Infinitary logics and models} \label{sec:infinitarylogics}
The topic of this section is infinitary logics (as presented e.g.\ in~\cite{Barwise:1969fk}) and models of infinitary sentences. 
As in Section~\ref{subsubsec:spacesofmodels}, \( \LL = \setofLR{R_i}{ i \in I } \) with \( I \in \omega \) denotes a finite relational signature, and \( n_i \) is the arity of the symbol \( R_i \).
As usual \( \lambda , \kappa \) denote infinite cardinals.

\subsection{Infinitary logics}
\subsubsection{Syntax} \label{subsubsec:syntax}
\begin{definition}\label{def:L_kappalambda}
For \( \lambda \leq \kappa \) the set \( \LL_{\kappa \lambda} \)\index[symbols]{L@\( \LL_{\kappa \lambda} \)} of infinitary formul\ae{} for the signature \( \LL \) is defined as follows.
\begin{itemize}[leftmargin=1pc]
\item
Fix a list of objects \( \seqofLR{ \V_\alpha}{ \alpha < \lambda } \)\index[symbols]{v@\( \V_\alpha \)} called \markdef{variables}.
\item
The \markdef{atomic formul\ae{}} are finite sequences of the form \( \seqLR{ \mathord {\equals } , \V_{ \alpha _1} , \V_{ \alpha_2} } \)\index[symbols]{54@\( \equals \)} and \( \seqLR{ R_i , \V_{ \alpha _1} , \dots , \V_{ \alpha _{n_i}} } \), with \( i < I \) and \( \alpha _1 , \dots , \alpha _{n_i} < \lambda \).%
\footnote{We use the symbol \( \equals \) for the equality predicate in the infinitary logic, to distinguish it from the usual \( = \) used in the language of set theory.}
\item
\( \LL_{\kappa \lambda} \) is the smallest collection containing the atomic formul\ae{} and closed under the following operations:
\begin{description}
\item[negation]
\( \upvarphi \mapsto \seqLR{ \neg } \conc \upvarphi \);
\item[generalized conjunctions]
if \( \upvarphi_ \alpha \in \LL_{ \kappa \lambda } \) (for \( \alpha < \nu < \kappa \)), and if the total number of variables that occur free in some \( \upvarphi_ \alpha \) is \( {<} \lambda \), then \( \seqLR{\bigwedge }\conc \seqofLR{ \upvarphi_ \alpha }{ \alpha < \nu } \in \LL_{ \kappa \lambda } \);
\item[generalized existential quantification] 
\( \upvarphi \mapsto \seqLR{ \exists} \conc \seqofLR{ \V_{u ( \alpha )}}{ \alpha < \nu } \conc \upvarphi \), for an increasing \( u \colon \nu \to \lambda \) and \( \nu < \lambda \).
\end{description}
\end{itemize}
\end{definition}

\begin{remarks}
\begin{enumerate-(i)}
\item
Each formula in \( \LL_{\kappa \lambda} \) has \( {<} \lambda \) free variables occurring in it.
In fact the formal definition of \( \upvarphi \in \LL_{\kappa \lambda} \) requires the simultaneous definition of 
\[ 
\Fv ( \upvarphi ) \in [ \lambda ]^{ < \lambda }, \index[symbols]{Fv@\( \Fv ( \upvarphi ) \)}
\] 
the set of all \( \alpha \in \lambda \) such that \( \V_ \alpha \) occurs free in \( \upvarphi \).
\item
Formally, the generalized conjunction of the formul\ae{} \( \upvarphi_\alpha \in \LL_{\kappa \lambda} \) (for \( \alpha < \nu < \kappa \)) should be defined as the concatenation \( \seqLR{\bigwedge }\conc \upvarphi_0 \conc \upvarphi_1 \conc \dotsc \conc \upvarphi_\alpha \conc \dotsc \) rather than \( \seqLR{\bigwedge }\conc \seqofLR{ \upvarphi_ \alpha }{ \alpha < \nu } \), so that each formula in \( \LL_{\kappa \lambda} \) is always a(n infinite) sequence of logical symbols and symbols from \( \LL \). 
However, this other approach would then require us to prove a unique readability lemma to ensure that each of the formul\ae{} \( \upvarphi_\alpha \) can be recovered from their conjunction, something which is clear with our current definition. 
This formal presentation would be considerably more opaque, so we decided to abandon it in favor of clarity.
\item
\( \LL_{ \omega \omega } \) is ordinary first-order logic. 
\end{enumerate-(i)}
\end{remarks}

The usual syntactical notions (subformula, sentence, etc.) are defined in the obvious way, and each \( \upvarphi \) is completely described by the tree of its subformul\ae{} \( \Subf ( \upvarphi ) \), which is called the \markdef{syntactical tree} of \( \upvarphi \).

It is convenient to assume that all symbols of \( \LL_{\kappa \lambda} \) are construed as fixed elements of \( \lambda \) --- for example we can encode \( \V_ \alpha \) by the ordinal \( \op{0}{\alpha} \), the equality predicate \( \equals \) by \( \op{ 1 }{ 0 } \), the negation \( \neg \) by \( \op{1}{1} \), the generalized conjunction \( \bigwedge \) by \( \op{1}{2} \), and \( R_i \) by \( \op{ i + 2 }{ n_i } \), where \( \op{ \cdot }{ \cdot } \) is the pairing function of~\eqref{eq:Hessenberg}. 
In this way, any atomic formula can be encoded by a finite sequence of ordinals \( {<} \lambda \) (and therefore by an ordinal \( {<} \lambda \) by~\eqref{eq:codingfinitesequenceofordinals}), and an arbitrary \( \upvarphi \in \LL_{\kappa \lambda} \) can be encoded via \( \Subf ( \upvarphi ) \) with a tree on \( \kappa \). 
In fact: \( \Subf ( \upvarphi ) \) can be construed as a \emph{labelled} \( < \kappa \)-branching descriptive set-theoretic tree on \( \kappa \) of height \( \leq \omega \), with labels in \( [\lambda]^{< \lambda} \). 
Notice that also formul\ae{} in \( \LL_{\kappa^+ \lambda} \) may be encoded through their syntactical tree as (labelled) tree on \( \kappa \) --- the difference with the previous case is that now we have to consider \( \kappa \)-branching trees.
Thus:
\begin{itemize}[leftmargin=1pc]
\item
the set \( \LL_{\kappa \lambda} \) can be defined in every transitive model of \( \ZF \) containing \( \kappa \),
\item
the predicate ``\( \upvarphi \in \LL_{\kappa \lambda} \)'' is absolute for such models,
\item
\( \LL_{\kappa' \lambda'} \subseteq \LL_{\kappa \lambda} \) whenever \( \kappa' \leq \kappa \) and \( \lambda' \leq \lambda \).
\end{itemize}

We forsake the official, but awkward, notation for a language in favor of a more relaxed (if a bit inaccurate) one.
Thus:
\begin{itemize}[leftmargin=1pc]
\item
the atomic formul\ae{} are written as \( \V_{ \alpha _1} \equals \V_{ \alpha_2} \) and \( R_i ( \V_{ \alpha _1} , \dots , \V_{ \alpha _{n_i}} ) \), and \( \V_{ \alpha _1} \nequals \V_{ \alpha_2} \) is the negation of \( \V_{ \alpha _1} \equals \V_{ \alpha_2} \);
\item
the letters \( \upvarphi , \uppsi \) range over formul\ae{}, while \( \upsigma \) ranges over sentences;
\item
the negation, the generalized conjunction, and the generalized existential quantification are written as \( \neg \upvarphi \), \( \bigwedge_{ \alpha < \nu } \upvarphi _ \alpha \), and \( \EXISTS{ \seqof{ \V_ {u ( \alpha ) }}{ \alpha < \nu }} \upvarphi \) or \( \exists \V_{u ( 0 )} \exists \V_{u ( 1 )} \cdots \, \upvarphi \) or%
\footnote{Recall from Section~\ref{subsubsec:sequences} that we are identifying each \( u \in [\lambda]^{< \lambda} \) with its increasing enumerating function.}
 \( \EXISTS{ \V_u} \upvarphi \) (with \( u \in [\lambda]^{<\lambda} \)).
The generalized disjunction \( \bigvee_{ \alpha < \nu } \upvarphi _ \alpha \) and generalized universal quantification \( \FORALL{ \seqof{ \V_ {u ( \alpha ) }}{ \alpha < \nu }} \upvarphi \) or \( \forall \V_{u ( 0 )} \forall \V_{u ( 1 )} \cdots \, \upvarphi \) or \( \FORALL{ \V_u} \upvarphi \) are defined by means of the de Morgan's laws.
Ordinary conjunctions and disjunctions are obtained from generalized ones by setting \( \nu = 2 \), and from these all other connectives are defined.
Similarly, ordinary quantifications are obtained from generalized ones by setting \( \nu = 1 \);
\item
the letters \( x , y , z , \dots \) range over \( \setofLR{\V_ \alpha }{ \alpha < \lambda } \), so that \( \EXISTS{ \seqof{ x_ \alpha }{ \alpha < \nu }} \upvarphi \) is \( \EXISTS{ \seqof{ \V_ {u ( \alpha ) }}{ \alpha < \nu }} \upvarphi \) for some \( u \colon \nu \to \lambda \).
Since cardinals are additively closed, \( \EXISTS{ \seqof{x_\alpha }{ \alpha < \nu }} \EXISTS{ \seqof{ y_\beta }{ \beta < \xi } }\upvarphi \) is identified with \( \EXISTS{ \seqof{x_\alpha }{ \alpha < \nu } \conc \seqof{ y_\beta }{ \beta < \xi }} \upvarphi \),
\item
the expression \( \upvarphi ( \seqofLR{ x_\alpha }{ \alpha < \nu }) \) means that the variables that occur free in \( \upvarphi \) are among the \( \setofLR{ x_\alpha }{\alpha < \nu } \), and we will also assume that the \( x_\alpha \) are distinct and listed in an increasing order with respect to the official list. 
Therefore \( \seqofLR{ x_\alpha }{ \alpha < \nu } = \seqofLR{ \V_{u ( \alpha )} }{\alpha < \nu } \) for some unique increasing \( u \colon \nu \to \lambda \).
Such \( u \) can be identified with its range, which is an element of \( [ \lambda ] ^\nu \), so we write \( \upvarphi ( \V_u ) \) with \( u \in [ \lambda ] ^\nu \) to say that if \( \V_ \alpha \) occurs free in \( \upvarphi \), then \( \alpha \in u \).
\end{itemize}

\begin{definition}\label{def:Lbounded}
\begin{enumerate-(i)}
\item \label{def:Lbounded-a} 
The set \( \LL_{\kappa \lambda }^0 \)\index[symbols]{L0@\( \LL_{\kappa \lambda }^0 \), \( \LL_{\kappa \lambda }^b \)} of all \markdef{propositional formul\ae} consists of those \( \upvarphi \in\LL_{\kappa \lambda} \) obtained from the atomic formul\ae{} using only negation and generalized conjunctions.
\item \label{def:Lbounded-b}
The set \( \LL^b_{\kappa \lambda} \) of all \markdef{bounded formul\ae}\index[concepts]{bounded formul\ae{} (infinitary logic)} consists of those \( \upvarphi \in\LL_{\kappa \lambda} \) such that:
\[
\text{if \( \exists \seqof{ x_\alpha }{ \alpha < \nu } \uppsi \in \Subf ( \upvarphi ) \) and \( \nu \geq \omega \), then \( \uppsi \in \LL_{\kappa' \lambda }^0 \) for some \( \kappa' < \kappa \).}
\]
\end{enumerate-(i)}
\end{definition}

Thus \( \LL^b_{\kappa \lambda} \) is closed under negation, infinitary conjunctions and disjunctions (of size \( {<} \kappa \)), finite quantifications, but \emph{not} under infinitary quantifications. 
The reason for investigating this technical notion is that \( \LL^b_{ \kappa \lambda } \) avoids the counterexamples related to (the relevant direction of) the generalized Lopez-Escobar theorem (see Remark~\ref{rmk:formulaborel}).
This feature of \( \LL^b_{\kappa \lambda} \) is crucial for the results of Section~\ref{subsec:topologicalcomplexity} (see Theorems~\ref{th:maintopology} and~\ref{th:barmaintopology}).

\begin{remarks}\label{rmk:L_alphabeta}
\begin{enumerate-(i)}
\item\label{rmk:L_alphabeta-i}
We are mostly interested in logics of the form \( \LL_{\kappa^+ \kappa} \), \( \LL^b_{\kappa^+ \kappa} \) and \( \LL_{\kappa^+ \lambda} \) for \( \lambda < \kappa \). 
\item\label{rmk:L_alphabeta-ii}
As for \( \alpha \)-Borel sets, one could define the logics \( \LL_{\alpha \beta} \) for arbitrary \( \beta \leq \alpha \in \On \) in the obvious way, and with such definition one could write \( \LL_{(\kappa + 1) \kappa} \), \( \LL^b_{(\kappa + 1) \kappa} \) and \( \LL_{(\kappa + 1) \lambda} \) instead of \( \LL_{\kappa^+ \kappa} \), \( \LL^b_{\kappa^+ \kappa} \) and \( \LL_{\kappa^+ \lambda} \).
In principle, this would be preferable, as \( \kappa + 1 \) is absolute while \( \kappa^+ \) is not.
However, there are two reasons to eschew such move.
Firstly, this would run against standard notation in the literature.
Secondly, this could be source of endless minor notational quibbles regarding the use of variables in the generalized existential quantification. 
\item\label{rmk:L_alphabeta-iii}
A more substantial move would be to extend the definition of \( \LL_{ \kappa \lambda } \) to \( \LL_{A + 1 \, B + 1} \) with \( A , B \) arbitrary sets, just like \( \bB_{J + 1} \) is a generalization of \( \bB_{ \alpha } \) (Definition~\ref{def:JBorel}).
This would allow us to give a \( \ZF \)-formulation of certain results (Proposition~\ref{prop:Borelformula} parts~\ref{prop:Borelformula-2} and \ref{prop:Borelformula-3}, and Proposition~\ref{prop:Borelformula2}).
On the other hand, besides the technical problems already mentioned for \( \LL_{ \alpha \beta } \), this would render the notation quite opaque.
As we have no use for \( \LL_{ A + 1 \, B + 1} \), we decided to abandon them altogether.
\end{enumerate-(i)}
\end{remarks}

\subsubsection{Semantics}
If \( \mathcal{A} = \seqLR{A ; R_i^\mathcal{A} }_{i \in I} \) is an \( \LL \)-structure, \( \upvarphi ( \V_u ) \in \LL_{ \kappa \lambda } \) with \( u \in [ \lambda ]^{<\lambda} \), and \( s \colon u \to A \), then 
\[ 
\mathcal{A}\models \upvarphi [ s ] 
\] 
means that the formula obtained from \( \upvarphi ( \V_u ) \) by substituting each \( \V_{ \alpha } \) with \( s ( \alpha ) \) for all \( \alpha \in u \), holds true in the structure \( \mathcal{A} \).
This notion is defined recursively on the tree \( \Subf ( \upvarphi ) \) of all subformul\ae{} of \( \upvarphi \).
For example, if \( \upvarphi \) is \( \V_ \alpha \equals \V_ \beta \) or \( R_i ( \V_{ \alpha _1 } , \dots , \V_{ \alpha _n }) \), then \( \mathcal{A} \models \upvarphi [ s ] \) if and only if \( s ( \alpha ) = s ( \beta ) \) or \( \seqLR{ s ( \alpha _1 ) , \dots , s ( \alpha _n ) } \in R_i^{\mathcal{A}} \).
If instead \( \exists \V_u \uppsi \), then letting \( w \coloneqq u \cup \Fv ( \uppsi ) \)
\begin{equation}\label{eq:satisfactionofexistentialformula}
\mathcal{A} \models \upvarphi [ s ] \IFF \EXISTS{t \in \pre{ w}{A} } \left ( t \restriction \Fv ( \upvarphi ) = s \restriction \Fv ( \upvarphi ) \AND \mathcal{A} \models \uppsi [ t ] \right ) . 
\end{equation}
As in the first-order case, if \( s \restriction \Fv ( \upvarphi ) = t \restriction \Fv ( \upvarphi ) \) then \( \mathcal{A} \models \upvarphi [ s ] \IFF \mathcal{A} \models \upvarphi [ t ] \), so if \( \upvarphi \) is a sentence, i.e. \( \Fv ( \upvarphi ) = \emptyset \), then the truth \( \mathcal{A} \models \upvarphi [ s ] \) does not depend on \( s \) and we write \( \mathcal{A} \models \upvarphi \).

\begin{remarks} \label{rmk:satisfactionisnotabsolute}
\begin{enumerate-(i)}
\item
First-order formul\ae{} can be identified with specific natural numbers and \( \LL_{ \omega \omega } \) can be identified with a subset of \( \omega \).
In particular, both \( \LL_{ \omega \omega } \) and its members are well-orderable.
This fails badly for \( \LL_{ \kappa \lambda } \) when \( \kappa > \omega \).
In fact, while atomic formul\ae{} (and their negations) can be coded as ordinals \( {<} \lambda \), by taking countable conjunctions it is possible to inject \( \pre{\omega}{2} \) into \( \LL^0_{ \kappa \lambda } \).
\item
By~\eqref{eq:satisfactionofexistentialformula}, when \( \lambda \geq \omega_1 \) the existence of Skolem functions for \( \mathcal{A} \) requires \( \AC \) (or at least a well-ordering of \( \pre{< \lambda}{ \lambda} \)), even when \( A \) is well-orderable. 
\item \label{rmk:satisfactionisnotabsolute-iii}
The satisfaction relation for \( \LL_{\kappa \lambda} \) with \( \kappa \geq \omega_2 \) and \( \lambda \geq \omega_1 \) is not absolute for transitive models of \( \ZFC \). 
In fact there can be a countable structure \( \mathcal{A} \) with domain \( \omega \) and an \( \LL_{\omega_2 \omega_1} \)-sentence \( \upsigma \) such that \( \mathcal{A} \models \upsigma \) in the universe \( \Vv \), but \( \mathcal{A} \not\models \upsigma \) in a suitable forcing extension \( \Vv[G] \) of \( \Vv \). 
To see this, assume that \( \ZFC+\CH \) holds in \( \Vv \), let \( \LL = \{ P \} \) be a language consisting of just one unary relational symbol, and let \( \mathcal{A} = \seqLR{\omega ; P^{\mathcal{A}}} \in \Vv \) be such that \( P^{\mathcal{A}} \) is infinite and coinfinite.
For the sake of definiteness set \( P^{\mathcal{A}}(i) \IFF i \) is even (for \( i \in \omega \)), so that \( \mathcal{A} \) belongs to any model of \( \ZF \).
Let 
\[ 
\mathscr{D} \coloneqq \{ D \in \Vv \mid \Vv \models \text{``\( D \) is \( \forcing{Fn} ( \omega , 2 ; \omega ) \)-dense''} \},
 \] 
so that \( \mathscr{D} \) has size \( \omega_1 \) in \( \Vv \) by \( \CH \). 
The plan is to define an \( \LL_{\omega_2 \omega_1} \)-sentence \( \upsigma \) coding the existence of a \( \mathscr{D} \)-generic for the Cohen forcing \( \forcing{Fn}(\omega,2;\omega) \).
Given \( s \in \pre{<\omega}{2} \), let \( \uppsi_s = \uppsi_s ( \V_0, \dotsc, \V_{\lh ( s ) - 1 }) \) be the \( \LL_{\omega \omega} \)-formula
\[ 
{\bigwedge_{\substack{i < \lh s \\ s ( i ) = 1}} P ( \V_i )} \wedge {\bigwedge_{\substack{i < \lh s \\ s ( i ) = 0}} \neg P ( \V_i )}.
 \] 
Finally, let \( \upsigma \) be the \( \LL_{\omega_2 \omega_1} \)-sentence \( \EXISTS{ \seqof{ \V_i}{i \in \omega}} \upvarphi \), where \( \upvarphi = \upvarphi ( \seqof{ \V_i}{i \in \omega} ) \) is the formula
\[ 
 \Bigl ( \bigwedge_{ i < j < \omega} \V_i \nequals \V_j \Bigr ) \wedge \bigwedge_{D \in \mathscr{D}} \bigvee_{s \in D} \uppsi_s ( \V_0, \dotsc, \V_{\lh ( s ) - 1 }) .
 \] 
It is not hard to check that working in any \( \ZF \)-model \( W \supseteq \Vv \), if \( \seqof{a_i}{i \in \omega} \) is a sequence of elements of \( \mathcal{A} \) such that \( \mathcal{A} \models \upvarphi ( \seqof{a_i}{i \in \omega}) \), then the function \( G \colon \omega \to 2 \) defined by \( G ( i ) = 1 \IFF P^{\mathcal{A}} ( a_i ) \) is \( \mathscr{D} \)-generic for \( \forcing{Fn} ( \omega , 2 ; \omega ) \). 
Conversely, if \( G \colon \omega \to 2 \) is \( \mathscr{D} \)-generic for \( \forcing{Fn} ( \omega , 2 ; \omega ) \) then, using the fact that \( \setofLR{i \in \omega}{P^{\mathcal{A}} ( i )} \) is infinite and coinfinite, one can find a sequence \( \seqof{a_i}{i \in \omega} \) of elements of \( \mathcal{A} \) such that \( \mathcal{A} \models \upvarphi ( \seqof{a_i}{i \in \omega} ) \). 
Therefore, in any \( W \) as above it holds
\[ 
\mathcal{A} \models \upsigma \IFF \EXISTS{G \colon \omega \to 2} ( G \text{ is \( \mathscr{D} \)-generic for } \forcing{Fn}( \omega,2;\omega)).
 \] 
Therefore \( \mathcal{A} \not\models \upsigma \) in \( \Vv \), but \( \mathcal{A} \models \upsigma \) in any \( \forcing{Fn} ( \omega , 2 ; \omega ) \)-generic extension of \( \Vv \).
\end{enumerate-(i)}
\end{remarks}

As seen in Section~\ref{subsubsec:spacesofmodels}, any \( \LL \)-structure of size \( \kappa \) can be identified with an element of \( \Mod^ \kappa _\LL \) (which is a typical example of a space of type \( \kappa \)).
For \( \upvarphi ( \V_u ) \in \LL_{\nu \mu} \) and \( u \in [ \mu ]^{ < \mu } \) we set
\begin{equation} \label{eq:M_phiu}
M_{\upvarphi , u} \coloneqq \setofLR{ ( X , s ) \in \Mod^\kappa_\LL \times \pre{u}{ \kappa }} {X \models \upvarphi [ s ]} , \index[symbols]{Mphiu@ \( M_{\upvarphi , u} \)}
\end{equation}
and if \( \upsigma \) is a sentence, we set 
\begin{equation} \label{eq:modsigma}
\Mod_\upsigma^\kappa \coloneqq \setofLR{ X \in \Mod^\kappa_\LL }{ X \models \upsigma }. \index[symbols]{Modsigma@\( \Mod_\upsigma^\kappa \), \( \Mod^{ < \kappa}_\upsigma \), \( \Mod_\upsigma^\infty \)}
\end{equation}
Thus \( \Mod^\kappa_\upsigma \) is the space of all \( \LL \)-structures with domain \( \kappa \) which satisfy \( \upsigma \). 
We also let 
\[ 
\Mod^{< \kappa}_\upsigma \coloneqq \bigcup_{\lambda < \kappa} \Mod^\lambda_\upsigma \quad \text{and} \quad \Mod_\upsigma^\infty \coloneqq \bigcup_{ \kappa \in \Cn} \Mod_\upsigma^\kappa .
\]
Notice that sets of the form \( \Mod^\kappa_\upsigma \), \( \Mod^{< \kappa}_\upsigma \), and \( \Mod^\infty_\upsigma \) are always invariant under isomorphism. 

\begin{notation}\label{not:embeddability}
From now on we use \( \embeds^\kappa_\upsigma \) and \( \cong^\kappa_\upsigma \)\index[symbols]{57@\( \embeds^\kappa_\upsigma \), \( \embeds^{<\kappa}_\upsigma \), \( \embeds^\infty_\upsigma \)} for the embeddability and isomorphism relations restricted to the space of models \( \Mod^\kappa_\upsigma \). 
Similarly, we denote the restriction of the embeddability (respectively, isomorphism) relation to the spaces \( \Mod^{< \kappa}_\upsigma \) and \( \Mod^\infty_\upsigma \) by \( \embeds^{<\kappa}_\upsigma \) and \( \embeds^\infty_\upsigma \) (respectively, \( \cong^{< \kappa}_\upsigma \) and \( \cong^\infty_\upsigma \)).\index[symbols]{57@\( \cong^\kappa_\upsigma \), \( \cong^{< \kappa}_\upsigma \), \( \cong^\infty_\upsigma \)} 
A similar notation is adopted for the less frequently used bi-embeddability relation \( \biembeds \) as well: in this case we write \( \biembeds^\kappa_\upsigma \)\index[symbols]{57@\( \biembeds^\kappa_\upsigma \), \( \biembeds^{< \kappa}_\upsigma \), \( \biembeds^\infty_\upsigma \)}, \( \biembeds^{< \kappa}_\upsigma \), and \( \biembeds^\infty_\upsigma \) to denote the restriction of \( \biembeds \) to, respectively, \( \Mod^\kappa_\upsigma \), \( \Mod^{< \kappa}_\upsigma \), and \( \Mod^\infty_\upsigma \).
\end{notation} 

\begin{remark}
Fix any \( \upsigma \in \LL_{\kappa^+ \kappa} \).
As we shall see in the next section, when \( \kappa < \kappa^{< \kappa} \) (or when the assumption \( \AC \) is dropped) the set \( \Mod^\kappa_\upsigma \) may fail to be a \( \kappa+1 \)-Borel subset of \( \Mod^\kappa_\LL \), even when considering the finest topology \( \tau_b \), see Remark~\ref{rmk:formulaborel}. 
However, it is worth to point out that if \( \kappa \) is regular and either \( \upsigma \in \LL^b_{\kappa^+ \kappa} \) or \( \AC + \kappa^{<\kappa} = \kappa \) holds, then \( \Mod^\kappa_{\upsigma} \in \bB^\mathrm{e}_{\kappa+1}(\Mod^\kappa_\LL,\tau_b) \) by Corollary~\ref{cor:formulaborel1}\ref{cor:formulaborel1-ii} and Theorem~\ref{th:LopezEscobar}. 
In this case \( \Mod^\kappa_\upsigma \) is a standard Borel \( \kappa \)-space and the relations \( \embeds^\kappa_\upsigma \), \( \cong^\kappa_\upsigma \), and \( \biembeds^\kappa_\upsigma \) are, respectively, a \( \kappa \)-analytic quasi-order and two \( \kappa \)-analytic equivalence relations.
\end{remark}

The topological complexity of \( \Mod^\kappa_\upsigma \) as a subspace of the space \( \Mod^\kappa_\LL \) of type \( \kappa \) (with respect to the various natural topologies on it, see Definition~\ref{def:spaceoftypekappa}) is the subject of Section~\ref{subsec:LopezEscobar}. 
Towards this goal, we equip the spaces \( \Mod^\kappa_\LL \times \pre{u}{ \kappa } \) in~\eqref{eq:M_phiu} with the product of the topology \( \tau \) on \( \Mod^\kappa_\LL \) and \( \sigma \) on \( \pre{ u }{ \kappa } \), where \( \tau \) is either the product topology, the \( \lambda \)-topology (\( \omega \leq \lambda < \max ( \cof ( \kappa) ^+ , \kappa ) \)), or the bounded topology, and \( \sigma \) is the discrete topology on \( \pre{u}{ \kappa } \).
As before, the resulting topologies are called, respectively, product topology, \( \lambda \)-topology, and bounded topology, and are denoted by \( \tau_p \), \( \tau_\lambda \), and \( \tau_b \). 
The bijection 
\begin{equation}\label{eq:5.3}
\Mod^\kappa_\LL \times \pre{ \emptyset }{ \kappa } \to \Mod^\kappa_\LL , \qquad ( X , \emptyset ) \mapsto X 
\end{equation}
is a homeomorphism (when on both sides \( \Mod^\kappa_\LL \) is endowed with the same topology) witnessing the identification between \( M_{ \upsigma , \emptyset } \) and \( \Mod^ \kappa _ \upsigma \) for every \( \upsigma \in \LL_{\kappa \lambda} \). 

\subsection{Some generalizations of the Lopez-Escobar theorem}\label{subsec:LopezEscobar}\index[concepts]{Lopez-Escobar's theorem}
A theorem of Lopez-Escobar (see e.g.~\cite[Theorem 16.8]{Kechris:1995zt}) says that \( A \subseteq \Mod^\omega_\LL \) is Borel and invariant under isomorphism if and only if \( A = \Mod^\omega_\upsigma \) for some \( \upsigma \in \LL_{\omega_1 \omega} \).
The aim of this section is to extend this to an arbitrary infinite cardinal \( \kappa \).
Full generalizations of the Lopez-Escobar theorem have been obtained in~\cite[Theorem 4.1]{Vaught:1974kl} when \(\kappa = \kappa^{< \kappa} \), yielding the next result which follows immediately from Corollary~\ref{cor:formulaborel} with \( \lambda = \kappa \), and Proposition~\ref{prop:Borelformula2}\ref{prop:Borelformula2-c} below.

\begin{theorem}[\( \AC \)]\label{th:LopezEscobar}
If \( \kappa ^{< \kappa } = \kappa \), then a set \( A \subseteq \Mod^\kappa_\LL \) is \( \kappa+1 \)-Borel (with respect to \( \tau_b \)) and closed under isomorphism if and only if \( A = \Mod_\upsigma^{ \kappa } \) for some \( \upsigma \in \LL_{ \kappa ^+ \kappa } \).
\end{theorem}

Since \( \bB_{\kappa + 1} ( \tau_p) \subseteq \bB_{\kappa + 1} ( \tau_\lambda) \subseteq \bB_{\kappa + 1} ( \tau_b) \), Theorem~\ref{th:LopezEscobar} concerns the largest possible class of \( \kappa + 1 \)-Borel sets, but in fact under \( \kappa^{< \kappa} = \kappa \) all the above classes coincide by Corollary~\ref{cor:borel}. 
Remarks~\ref{rmk:formulaborel} and~\ref{rmk:LopezEscobar} show that without the assumption \( \kappa^{<\kappa} = \kappa \) both directions of the equivalence in Theorem~\ref{th:LopezEscobar} may fail. 
The assumption \( \kappa^{<\kappa} = \kappa \) in Theorem~\ref{th:LopezEscobar} is inconvenient for our work because, as already explained in the introduction:
\begin{itemize}[leftmargin=1pc]
\item
we need to apply a generalization of the Lopez-Escobar theorem in models where \( \AC \) fails (e.g.\ in models of determinacy);
\item
even when working in models satisfying \( \AC \), our results concern uncountable cardinals \( \kappa < \card{\pre{\omega}{2}} \), which cannot satisfy \( \kappa ^{< \kappa } = \kappa \).
\end{itemize}
However, a careful analysis of the proof of Theorem~\ref{th:LopezEscobar} reveals some new intermediate results that may be useful. 
In one direction, Corollaries~\ref{cor:formulaborel1} and~\ref{cor:formulaborel} show that if the sentence \( \upsigma \) is chosen in a suitable fragment of \( \LL_{\kappa^+ \kappa} \), then \( \Mod^\kappa_\upsigma \) may turn out to be \( \kappa+1 \)-Borel even if the assumptions \( \AC \) and \( \kappa^{< \kappa} = \kappa \) are dropped --- in fact, the bounded version \( \LL^b_{\kappa^+ \kappa} \) of the logic \( \LL_{\kappa^+ \kappa} \) has been introduced for this purpose. 
In another direction, Proposition~\ref{prop:Borelformula2} gives some interesting results even in situations when Theorem~\ref{th:LopezEscobar} cannot be applied. 
For example, part~\ref{prop:Borelformula2-a} yields in \( \ZFC \) that if \( \Mod^{\omega_1}_\LL \) is endowed with the product topology, then every \( \omega_1 + 1 \)-Borel set \( A \subseteq \Mod^{\omega_1}_\LL \) which is invariant under isomorphism is of the form \( \Mod^{\omega_1}_\upsigma \) for some \( \upsigma \in \LL_{(2^{\aleph_0})^+ \omega_1} \), independently of \( \CH \). 
More generally, when \( \kappa = \mu^+ \) and \( \Mod^\kappa_\LL \) is endowed with the product topology, then every \( \kappa + 1 \)-Borel set \( A \subseteq \Mod^\kappa_\LL \) invariant under isomorphism is always of the form \( \Mod^\kappa_\upsigma \) for some \( \LL_{(2^\mu)^+ \kappa} \)-sentence \( \upsigma \). 
Similarly, if \( 2^{\aleph_0} = \aleph_2 \) and \( 2^{\aleph_1} = \aleph_3 \) then Theorem~\ref{th:LopezEscobar} cannot be applied with \( \kappa = \omega_2 \). 
However Proposition~\ref{prop:Borelformula2}\ref{prop:Borelformula2-b} gives that if \( \Mod^{\omega_2}_\LL \) is endowed with the \( \omega_1 \)-topology, then every \( \omega_2 + 1 \)-Borel set \( A \subseteq \Mod^{\omega_2}_\LL \) invariant under isomorphism is of the form \( \Mod^{\omega_2}_\upsigma \) for some \( \upsigma \in \LL_{\omega_4 \omega_2 } \).

\subsubsection{From formul\ae{} to invariant Borel sets}
In the next results we use the sets \( M_{\upvarphi , u} \) defined in~\eqref{eq:M_phiu}.
The proof of the following lemma is a straightforward adaptation of the proof of~\cite[Proposition 16.17]{Kechris:1995zt}.

\begin{lemma} \label{lem:formulaBorel}
Let \( \lambda \leq \kappa \) and \( \upvarphi( \V_u ) \in \LL_{\kappa\lambda} \) with \( u \in [ \lambda ]^{< \lambda} \). 
The following hold:
\begin{enumerate-(a)}
\item \label{lem:formulaBorel-i}
If \( \upvarphi \) is atomic then \( M_{\upvarphi , u} \) is a basic clopen set with respect to \( \tau_p \) (and hence also with respect to \( \tau_ \lambda \) and \( \tau_b \)).
\item \label{lem:formulaBorel-ii}
If \( \upvarphi = {\neg \uppsi} \) then \( M_{\upvarphi , u} = ( \Mod^\kappa_\LL \times \pre{ u }{ \kappa } ) \setminus M_{ \uppsi , u} \).
\item \label{lem:formulaBorel-iii}
If \( \upvarphi = \bigwedge _{ \alpha < \nu } \uppsi_ \alpha \), where \( \nu < \kappa \), then \( M_{\upvarphi , u} = \bigcap_{ \alpha < \nu } M_{ \uppsi_ \alpha , u} \).
\item \label{lem:formulaBorel-iv} 
If \( \upvarphi = {\EXISTS{ \V_u} \uppsi} \) for some \( u \in [ \lambda ] ^{< \lambda} \), let \( w \coloneqq u \cup \Fv ( \uppsi ) \).
Then 
\[
M_{\upvarphi , u} = \bigcup_{t \in \pre{ w } { \kappa}} \pi ^{-1}_ t ( M_{\uppsi , w }) ,
\] 
where \( \pi_t \colon \Mod^\kappa_\LL \times \pre{u}{ \kappa} \to \Mod^\kappa_\LL \times \pre{ w}{ \kappa } \) is the continuous function 
\[ 
( x , s ) \mapsto ( x , s \restriction \Fv ( \upvarphi ) \cup t \restriction ( w \setminus \Fv ( \upvarphi ) ) ) .
\]
\end{enumerate-(a)}
\end{lemma}

In the following results of this section, the ambient space, unless otherwise indicated, is \( \Mod^ \kappa _\LL \).
Recall from Definition~\ref{def:JBorel} that for \( J \) an arbitrary set and \( \mathcal{B}_* \) the canonical basis for the topology \( \tau_* \in \set{ \tau_p , \tau_ \lambda , \tau_b } \), we denote by \( \bB^{\mathrm{e}}_{J + 1} ( \mathcal{B}_* ) \) the collection of all effective \( J + 1 \)-Borel subsets of \( \Mod^ \kappa _\LL \) with codes taking value in \( \mathcal{B}_* \).

\begin{proposition} \label{prop:formulaborel}
Let \( \lambda \leq \kappa \) and \( \upvarphi ( \V_u ) \in \LL_{\kappa^+ \lambda} \) with \( u \in [ \lambda ]^{< \lambda} \). 
Then the following hold:
\begin{enumerate-(a)}
\item \label{prop:formulaborel-i}
\( M_{\upvarphi , u} \in \bB^{\mathrm{e}}_{ J + 1} ( \tau_p ) \subseteq \bB^{\mathrm{e}}_{ J + 1} ( \tau_ \lambda ) \subseteq \bB^{\mathrm{e}}_{ J + 1} ( \tau_b ) \), where \( J \coloneqq \pre{<\lambda} {\kappa} \).
\item \label{prop:formulaborel-ii}
If \( \kappa \) is \emph{regular} and \( \upvarphi \in \LL^0_{\kappa \lambda } \), then \( M_{\upvarphi , u} \) is \( \tau_b \)-clopen.
\item\label{prop:formulaborel-iii} 
If \( \kappa \) is \emph{regular} and \( \upvarphi \in \LL^b_{\kappa^+ \lambda} \), then \( M_{\upvarphi , u} \in \bB^{\mathrm{e}}_{ \kappa + 1 } ( \tau_b ) \).
\item\label{prop:formulaborel-iv} 
If \( \lambda = \omega \), then \( M_{ \upvarphi , u} \in \bB^{\mathrm{e} }_{ \kappa + 1 } ( \tau_p ) \).
\end{enumerate-(a)}
\end{proposition}

\begin{proof}
For part~\ref{prop:formulaborel-i} it is enough to argue by induction on the complexity of the subformul\ae{} of \( \upvarphi \) and use Lemma~\ref{lem:formulaBorel} to construct an effective \( \pre{ <\lambda} { \kappa } + 1 \)-Borel code \( ( T, \phi ) \) for \( M_{\upvarphi , u} \).
We inductively define a function \( f \) from the set of the subformul\ae{} \( \uppsi \) of \( \upvarphi \) to the set of effective \( \pre{ < \lambda} { \kappa } + 1 \)-Borel codes such that \( f ( \uppsi ) \) gives a code for \( M_{\uppsi , u} \).
(The function \( f \) allows to avoid the use of \( \AC \) in the inductive steps.) 
The tree \( T \) is obtained from \( \Subf ( \upvarphi ) \), the syntactical tree of the subformul\ae{} of \( \upvarphi \), by replacing every non-splitting node given by an infinitary existential quantifier \( \exists \seqofLR{ x_ \alpha }{ \alpha < \mu } \) with a splitting node of \( T \) with (at most) \( \pre{ \mu} { \kappa } \)-many immediate successors, each corresponding to a possible witness of the existential quantification.

\smallskip

For part~\ref{prop:formulaborel-ii} argue as follows.
By inductively applying Lemma~\ref{lem:formulaBorel} to the subformul\ae{} of \( \upvarphi \) as in part~\ref{prop:formulaborel-i}, we obtain that \( M_{\upvarphi , u} \in \Alg ( \mathcal{C} , \kappa ) \) where \( \mathcal{C} \) is the collection of all basic \( \tau_p \)-clopen subsets of \( \Mod^\kappa_\LL \times \pre{u }{ \kappa } \).
Each of these sets is thus \( \tau_b \)-clopen by equation~\eqref{eq:tau_bclosedunderintersections} on page~\pageref{eq:tau_bclosedunderintersections} (see the comments following Definition~\ref{def:spaceoftypekappa}).

\smallskip

The proof of part~\ref{prop:formulaborel-iii} is similar to that of part~\ref{prop:formulaborel-i}.
In constructing an effective \( \kappa + 1 \)-Borel code for \( M_{\upvarphi , u} \), unions of size \( > \kappa \) are taken only when an infinitary existential quantification is encountered. 
By the way \( \LL^b_{\kappa^+ \lambda} \) is defined (see Definition~\ref{def:Lbounded}), the quantified subformula must be in \( \LL_{\kappa \lambda }^0 \), and hence such a union is in \( \tau_b \) by part~\ref{prop:formulaborel-ii}.
The result then follows by inductively applying Lemma~\ref{lem:formulaBorel} again.

\smallskip

Part~\ref{prop:formulaborel-iv} follows from part~\ref{prop:formulaborel-i} and the fact that \( \pre{ < \omega }{ \kappa } \) is in bijection with \( \kappa \).
\end{proof} 

By~\eqref{eq:5.3}, setting \( u = \emptyset \) in Proposition~\ref{prop:formulaborel} we get:

\begin{corollary} \label{cor:formulaborel1}
Suppose \( \lambda \leq \kappa \).
\begin{enumerate-(a)}
\item \label{cor:formulaborel1-i}
If \( \upsigma \in \LL_{\kappa^+ \lambda} \), then \( \Mod^\kappa_{\upsigma} \in \bB^{\mathrm{e}}_{ J + 1 } ( \tau_p ) \subseteq \bB^{\mathrm{e}}_{ J + 1 } ( \tau_b ) \), where \( J \coloneqq \pre{< \lambda}{\kappa} \).
\item \label{cor:formulaborel1-ii}
If \( \kappa \) is regular and \( \upsigma \in \LL^b_{\kappa^+ \kappa } \), then \( \Mod^\kappa_\upsigma \in \bB^{\mathrm{e}}_{ \kappa + 1 } ( \tau_b ) \).
\item \label{cor:formulaborel1-iii}
If \( \upsigma \in \LL_{\kappa^+ \omega} \), then \( \Mod^\kappa_\upsigma \in \bB^{\mathrm{e}}_{ \kappa + 1 } ( \tau_p ) \).
\end{enumerate-(a)}
\end{corollary}

Using Proposition~\ref{prop:formulaborel}\ref{prop:formulaborel-i} and Corollary~\ref{cor:formulaborel1}\ref{cor:formulaborel1-i} we also have:

\begin{corollary}[\( \AC \)] \label{cor:formulaborel}
Suppose \( \lambda \leq \kappa \) and \( \kappa ^{< \lambda } = \kappa \), and let \( \upvarphi ( \V_u ) \in\LL_{\kappa^+ \lambda} \) with \( u \in [ \lambda ]^{< \lambda} \). 
Then \( M_{ \upvarphi , u} \in \bB^{\mathrm{e}}_{ \kappa + 1 } ( \tau_p ) \subseteq \bB^{\mathrm{e}}_{ \kappa + 1 } ( \tau_b ) \).

Similarly, if \( \upsigma \in \LL_{\kappa^+ \lambda} \) is a sentence, then \( \Mod^\kappa_\upsigma \in \bB^{\mathrm{e}}_{ \kappa + 1 } ( \tau_p ) \subseteq \bB^{\mathrm{e}}_{ \kappa + 1 } ( \tau_b ) \).
\end{corollary}

\begin{remark}\label{rmk:formulaborel}
If \( \AC \) holds and \( \kappa ^{< \kappa } \neq \kappa \), then there may be sentences \( \upsigma \in \LL_{ \kappa ^+ \kappa } \) such that \( \Mod^\kappa_\upsigma \) is not \( \kappa + 1 \)-Borel with respect to the finest topology considered here, namely the bounded topology \( \tau_b \).
As observed in~\cite[Remark 25]{Friedman:2011nx}, by work of V\"a\"an\"anen and Shelah, if \( \lambda ^+ = \kappa < \kappa ^{ < \kappa } \) and \( \lambda ^{ < \lambda } = \lambda \) and a forcing axiom holds (and \( \omega _1^{\Ll} = \omega _1 \) if \( \lambda = \omega \)), then for some \( \upsigma \in \LL_{ \kappa \kappa } \subseteq \LL_{\kappa^+ \kappa} \) the set \( \Mod^ \kappa _{ \upsigma} \) does not have the \( \kappa \)-Baire property, so it is not \( \kappa + 1 \)-Borel with respect to \( \tau_b \) by Proposition~\ref{prop:borelsetshaveBP}.
Corollary~\ref{cor:formulaborel1}\ref{cor:formulaborel1-ii} is thus one of the main reasons to introduce the bounded logic \( \LL^b_{\kappa^+ \kappa} \).
\end{remark}

\subsubsection{From invariant Borel sets to formul\ae{}}

Let 
\[
\Inv^ \kappa _\LL \coloneqq \setofLR{A \subseteq\Mod^ \kappa _\LL}{ \FORALL{X , Y \in \Mod^ \kappa _\LL} ( X \in A \wedge Y \cong X \implies Y \in A ) }
\]
be the family of all subsets of \( \Mod^ \kappa _\LL \) which are \markdef{invariant under isomorphism}.
For every \( \upsigma \in \LL_{\nu \mu} \) and \( \mu \leq \nu' \) infinite cardinals, \( \Mod^\kappa_\upsigma \in \Inv^ \kappa _\LL \).
The following propositions provide a partial converse to Corollaries~\ref{cor:formulaborel1} and~\ref{cor:formulaborel}. 

\begin{proposition}\label{prop:Borelformula} 
Suppose \( \kappa , \lambda \) are infinite cardinals with \( \omega \leq \lambda < \max ( \cof ( \kappa )^+ , \kappa ) \).
\begin{enumerate-(a)}
\item \label{prop:Borelformula-1}
\( \bB_{\omega + 1}^{\mathrm{e}} ( \mathcal{B}_p ) \cap \Inv^ \kappa _\LL \subseteq \setofLR{ \Mod^\kappa_\upsigma }{ \upsigma \in \LL_{\kappa^+ \omega } } \).

\item \label{prop:Borelformula-2}
Assume \( \AC \). 
Then \( \bB_{\cof ( \lambda ) + 1} ( \mathcal{B}_\lambda ) \cap \Inv^ \kappa _\LL \subseteq \setofLR{ \Mod^\kappa_\upsigma }{ \upsigma \in \LL_{ ( \kappa^{< \lambda} )^+ \lambda } } \). 
\item \label{prop:Borelformula-3}
Assume \( \AC \). 
Then \( \bB_{\cof ( \kappa ) + 1} ( \mathcal{B}_b ) \cap \Inv^ \kappa _\LL \subseteq \setofLR{\Mod_\upsigma^{ \kappa } }{ \upsigma \in \LL_{( \kappa^{< \kappa} )^+ \kappa } } \).
\end{enumerate-(a)}
\end{proposition}

\begin{proposition}[\( \AC \)] \label{prop:Borelformula2}
Suppose \( \kappa \) is regular.
\begin{enumerate-(a)}
\item \label{prop:Borelformula2-a}
\( \bB_{\kappa + 1} ( \tau_p ) \cap \Inv^ \kappa _\LL \subseteq \setofLR{\Mod_\upsigma^{ \kappa } }{ \upsigma \in \LL_{( \kappa^{< \kappa} )^+ \kappa} } \).
\item \label{prop:Borelformula2-b}
Let \( \omega < \lambda < \kappa \) be such that \( \kappa^{< \lambda} = \kappa \). 
Then \( \bB_{\kappa + 1} ( \tau_\lambda ) \cap \Inv^ \kappa _\LL \subseteq \setofLR{\Mod_\upsigma^{ \kappa } }{ \upsigma \in \LL_{( \kappa^{< \kappa} )^+ \kappa} } \).
\item \label{prop:Borelformula2-c}
Assume that \( \kappa^{< \kappa} = \kappa \). 
Then \( \bB_{\kappa + 1} ( \tau_b ) \cap \Inv^ \kappa _\LL \subseteq \setofLR{\Mod_\upsigma^{ \kappa } }{ \upsigma \in \LL_{\kappa^+ \kappa}} \).
\end{enumerate-(a)}
\end{proposition}

The proofs of Propositions~\ref{prop:Borelformula} and~\ref{prop:Borelformula2} are a careful refinement of an argument already implicit in the work of Vaught~\cite{Vaught:1974kl} --- see~\cite[Proposition 16.9]{Kechris:1995zt}.

Recall from Section~\ref{subsec:Bairespace} that \( \pre{ \kappa }{ \kappa } \) and its subspaces can be equipped with several topologies, namely \( \tau_p = \tau_ \omega \), \( \tau _ \mu \) with \( \omega \leq \mu < \max ( \cof ( \kappa ) ^+ , \kappa ) \), and \( \tau_b \).
When \( \kappa \) is \emph{regular} then \( \tau_ \kappa = \tau_b \) and \( \widehat{\mathcal{B}}_ \kappa = \widehat{\mathcal{B}}_b \) as agreed in~\eqref{eq:conventionbasis}, so that the notation \( \tau_\mu \) and \( \widehat{\mathcal{B}}_ \mu \) encompass all possibilities.
If \( \kappa \) is \emph{singular}, then this is not the case and the bounded topology and its canonical basis should be treated separately.
This would cause a lot of notational inconveniences in the results that follow, so we stipulate the following:

\begin{convention}
For the rest of this section, we agree that \( \tau_ \kappa = \tau_b \) and \( \widehat{\mathcal{B}}_ \kappa = \widehat{\mathcal{B}}_b \), independently of whether \( \kappa \) is regular or not.
Thus the topologies and bases relevant for us are of the form \( \tau_\mu \) and \( \widehat{\mathcal{B}}_\mu \) with \( \mu \) in the set
\begin{equation}\label{eq:conventionbasissingular}
 \setofLR{ \nu \in \Cn}{ \omega \leq \nu \leq \cof ( \kappa ) \vee \nu = \kappa } .
\end{equation}
\end{convention}

The group \( \Sym ( \kappa ) \) inherits the relative topology (denoted again by \( \tau_ \mu \)), whose basic open sets can be written as
\[
[ s ] \coloneqq \widehat{\Nbhd}_{s^{-1}} \cap \Sym(\kappa) = \setofLR{ g \in \Sym ( \kappa ) }{ s \subseteq g^{-1} } 
\] 
with \( s \in \pre{ u}{ ( \kappa ) } \) (i.e.\ \( s \) is an injection from \( u \) to \( \kappa \), see Section~\ref{subsubsec:sequences}) and \( u \in \mathrm{D} ( \mu ) \), where 
\begin{equation}\label{eq:D(mu)}
 \mathrm{D} ( \mu ) \coloneqq \begin{cases}
[ \kappa ]^{< \mu } & \text{if } \mu \neq \kappa,
\\
\kappa & \text{if } \mu = \kappa .
\end{cases}
\end{equation}
Given a property \( \varphi \) for the elements of \( \Sym ( \kappa) \) and a nonempty \( \tau_ \mu \)-open set \( U \subseteq \Sym ( \kappa) \), we write
\[ 
\forall^*_ \mu g \in U \, \varphi ( g ) 
 \] 
to abbreviate the statement: 
\[
\setofLR{g \in \Sym ( \kappa) }{ \varphi ( g ) }\text{ is \( \cof ( \mu ) \)-comeager in } U ,
\] 
where the notion of \( \mu \)-comeagerness is defined in Section~\ref{subsec:category}.
Next we define the (local) Vaught transform of a set \( A \subseteq \Mod^\kappa_\LL \): given \( u \in \mathrm{D} ( \mu ) \), let
\[
A_ \mu ^{*u} \coloneqq \setofLR{ ( X , s ) \in \Mod_{\mathcal{L}}^ \kappa \times \pre{ u}{ ( \kappa ) } }{ \forall^*_ \mu g \in [ s ] \, \left ( g . X \in A \right ) },
\]
where \( g . X \) is as in~\eqref{eq:action2}. 

We are now going to prove Lemma~\ref{lem:LopezEscobar}, an analogue of~\cite[Proposition 16.9]{Kechris:1995zt}, from which both Propositions~\ref{prop:Borelformula} and~\ref{prop:Borelformula2} follow.
The direct adaptation of the proof of~\cite[Proposition 16.9]{Kechris:1995zt} to our context yields a formula \( \upvarphi_u \) whose variables range in \( \setof{ \V_\alpha}{\alpha < \kappa} \), while in the logic \( \LL_{ \nu \mu} \) of Lemma~\ref{lem:LopezEscobar} we can use only variables from the (possibly smaller) set \( \setof{ \V_\alpha }{ \alpha < \mu } \) (see Definition~\ref{def:L_kappalambda}). 
Such an argument would yield a \( \upvarphi_u \) which is \emph{essentially} what we require, but not quite an element of \( \LL_{ \nu \mu} \). 
To overcome this purely technical difficulty, a somewhat artificial fragment \( \LL_{ \nu \kappa}^{\downarrow } \) of \( \LL_{ \nu \kappa} \) is introduced (see the beginning of the proof of Lemma~\ref{lem:LopezEscobar}). 
In order to define such fragment, we need a few preliminary definitions.

Let \( \mathrm{D}( \mu ) \) be as in~\eqref{eq:D(mu)}, and fix \( u \in \mathrm{D} ( \mu ) \).
Set 
\[ 
u{\downarrow} \mu \coloneqq \begin{cases}
\text{the unique \( \nu < \mu \) such that } u \in [ \kappa ]^\nu & \text{if } \mu \neq \kappa, 
\\
u & \text{if } \mu = \kappa
\end{cases}
\]
so that in either case \( u {\downarrow} \mu \in [ \mu ]^{< \mu} \). 
Recalling from Section~\ref{subsubsec:sequences} that every element of \( [ \kappa ]^{ < \mu } \) is identified with its enumerating function, for any \( s \in \pre{u}{ ( \kappa ) } \) let 
\[ 
s^{\downarrow \mu} \coloneqq 
\begin{cases}
s \circ u & \text{if } \mu \neq \kappa, 
\\
s & \text{if } \mu = \kappa ,
\end{cases}
\] 
so that in either case \( s^{\downarrow \mu} \in \bigcup_{ \nu < \mu } \pre{ \nu }{ ( \kappa ) } \).
Notice that when \( u = s = \emptyset \) we have \( u{\downarrow} \mu = s^{\downarrow \mu} = \emptyset \) (for any cardinal \( \mu \) as above). 

\begin{lemma} \label{lem:LopezEscobar}
Suppose \( \lambda \leq \mu \) both belong to the set in~\eqref{eq:conventionbasissingular}, and endow the spaces \( \Sym ( \kappa ) \) and \( \Mod^\kappa_\LL \) with the topologies \( \tau_\mu \) and \( \tau_\lambda \), respectively.
Let \( A \subseteq \Mod^\kappa_\LL \) be in \( \bB^{\mathrm{e}}_{\cof ( \mu ) + 1} ( \mathcal{B}_\lambda ) \), and let \( u \in \mathrm{D}(\mu) \).
Let \( \nu \coloneqq \card{ \pre{ < \mu }{ \kappa } }^+ \), assuming \( \AC \) if \( \mu > \omega \). 
Then there is some \( \LL_{ \nu \mu} \)-formula \( \upvarphi_u ( \V_{u { \downarrow } \mu } ) \) such that for every \( X \in \Mod^ \kappa_\LL \) and \( s \in \pre{u}{( \kappa )} \)
 \[ 
 ( X , s ) \in A_\mu^{*u} \IFF ( X , s^{\downarrow \mu} ) \in M_{ \upvarphi_u , u { \downarrow} \mu } ,
\] 
where \( M_{ \upvarphi , u } \) is as in~\eqref{eq:M_phiu}.
\end{lemma}

\begin{proof}
Let \( \LL_{ \nu \kappa}^{\downarrow \mu} \) be the fragment of \( \LL_{ \nu \kappa} \) obtained by adding to Definition~\ref{def:L_kappalambda} the following restrictions:
\begin{itemize}[leftmargin=1pc]
\item
the generalized conjunction of the formul\ae{} \( \upvarphi_\alpha \) (for \( \alpha < \nu < \nu \)) may be formed only when the total number of variables occurring free in some of the \( \upvarphi_\alpha \)'s is \( < \mu \);
\item
a generalized existential quantification \( \EXISTS{\V_u} \upvarphi \) may be formed only when \( u \in [ \kappa ]^{< \mu} \).
\end{itemize}
Equivalently, \( \LL_{ \nu \kappa}^{\downarrow \mu} \) can be defined similarly to the logic \( \LL_{ \nu \mu} \) except that we may use variables from a longer list \( \seqofLR{\V_\alpha}{\alpha< \kappa} \) of length \( \kappa \) instead of using \( \mu \)-many variables. 
(When \( \mu = \kappa \) we get \( \LL_{ \nu \kappa}^{\downarrow \mu} = \LL_{ \nu \kappa} \).) 
It follows in particular that each formula in \( \LL_{ \nu \kappa}^{\downarrow \mu} \) has \( < \mu \)-many free variables occurring in it.
In order to simplify the notation, for the rest of the proof we write \( \LL_{ \nu \kappa}^{\downarrow } \), \( u {\downarrow} \), and \( s^{\downarrow} \) instead of \( \LL_{ \nu \kappa}^{\downarrow \mu} \), \( u {\downarrow} \mu \) and \( s^{\downarrow \mu } \).

By a suitable variable substitution, any \( \uppsi ' ( \V_u ) \in \LL_{ \nu \kappa}^{\downarrow } \) can be easily transformed into a corresponding \( \uppsi ( \V_{u {\downarrow} } ) \in \LL_{ \nu \mu} \) such that for all \( ( X , s ) \in \Mod^\kappa_\LL \times \pre{u}{ ( \kappa ) } \)
\[ 
( X , s ) \in M_{ \uppsi' , u } \IFF ( X , s^{ \downarrow } ) \in M_{\uppsi, u { \downarrow } }.
 \] 
Thus it is enough to find, for any given \( A \) as in its hypotheses of the lemma, an \( \LL_{ \nu \kappa}^{\downarrow } \)-formula \( \upvarphi '_u ( \V_u ) \) such that \( A^{*u}_\mu = M_{ \upvarphi ' _ u , u } \): setting \( \upvarphi ' _ u \coloneqq G ( \emptyset , u ) \) in the following construction we will be done.
 
Let \( ( T , \phi ) \) be a \( \cof ( \mu ) + 1 \)-Borel code for \( A \subseteq \Mod ^ \kappa_\LL \) such that \( \phi ( t ) \in \mathcal{B}_\lambda ( \Mod ^ \kappa_\LL ) \) for every terminal node \( t \in T \). 
We shall define a map 
\[
 G \colon T \times \mathrm{D} ( \mu ) \to \LL_{ \nu \kappa}^{\downarrow }
\] 
such that \( \Fv ( G ( t , u ) ) = u \) and
\[ 
\phi ( t )_\mu^{*u} = M_{ G ( t , u ) , u} .
\] 
The map \( G \) is defined inductively on the well-founded relation 
\[
( t , u ) \leq ( t' , u' ) \IFF t' \text{ precedes \( t \) in the ordering of } T ,
\]
using that \( ( \Sym ( \kappa ) , \tau_\mu ) \) is \( \cof ( \mu ) \)-Baire (Theorem~\ref{th:Bairespacesym}).

\smallskip

\textbf{Case 1:} \( t \) is a terminal node of \( T \), so that \( \phi ( t ) \in \mathcal{B}_\lambda ( \Mod^\kappa_\LL ) \).
Then it is easy to check that there are%
\footnote{In fact \( \uptheta \) is a Boolean combination of atomic formul\ae{} of \( \LL \) using conjunctions and disjunctions of size \( {<} \lambda \).} 
\( \uptheta ( \V_w ) \in \LL_{\nu \kappa}^{\downarrow } \) with \( w \in \mathrm{D} ( \mu ) \), and \( h \in \pre{ w }{ ( \kappa ) } \), such that 
\[
\phi ( t ) = \setofLR{ X \in \Mod ^ \kappa _ \LL }{ X \models \uptheta [ h ] } .
\] 
Then for every \( u \in \mathrm{D} ( \mu ) \) we have 
\begin{equation}\label{eq:lem:LopezEscobar}
\begin{split}
( X , s ) \in \phi ( t )_\mu^{*u} & \iff s \in \pre{u }{ (\kappa ) } \wedge \forall ^*_\mu g \in [ s ] \left ( g . X \models \uptheta [ h ] \right ) 
\\
 & \iff s \in \pre{u }{ ( \kappa ) } \wedge \forall ^*_\mu g \in [ s ] \left ( X \models \uptheta [ g ^ { - 1 } \circ h ] \right ) .
\end{split} 
\end{equation}
Let \( w' \coloneqq \setofLR{ h ( \alpha ) }{ \alpha \in w } \), and let \( \uptheta ' ( \V_{ w ' } ) \) be the formula obtained from \( \uptheta ( \V_w ) \) by substituting each (free) occurrence of \( \V_\alpha \) in \( \uptheta ( \V_w ) \) with \( \V_{ h ( \alpha ) } \) (for each \( \alpha \in w \)).
Now there are two cases:
\begin{itemize}[leftmargin=1pc]
 \item 
if \( w' \subseteq u \) then since \( g \in [ s ] \iff s \subseteq g^{-1} \) we have \( g^{ - 1 } \circ h = s \circ h = s \restriction w' \), so \( \phi ( t )_\mu^{ * u } = M_{ \upvarphi , u} \) with \( \upvarphi ( \V_u ) \in \LL_{ \nu \kappa}^{\downarrow } \) being the formula
\[ 
{\bigwedge_{\substack{\alpha , \beta \in u \\ \alpha < \beta }} \V_ \alpha \nequals \V_ \beta } \wedge \uptheta' ( \V_{w'} ) .
\]
In fact if \( ( X , s ) \in \phi ( t )_\mu^{*u} \) then by~\eqref{eq:lem:LopezEscobar} there is \( g \in [ s ] \) such that \( X \models \uptheta [ g^{-1} \circ h ] \), and hence \( X \models \uptheta [ s \circ h ] \).
From this it follows that \( X \models \upvarphi [ s ] \), i.e.\ \( (X,s) \in M_{ \upvarphi , u}\).
Conversely, if \( ( X , s ) \in M_{ \upvarphi , u} \) then \( X \models \upvarphi [ s ] \).
Therefore \( X \models \uptheta [ g^{-1} \circ h ] \) for every \( g \in [ s ] \), and hence \( ( X , s ) \in \phi ( t )_\mu^{*u} \) since \( [s] \) is trivially \( \cof(\mu) \)-comeager in itself.
\item 
if \( w' \not\subseteq u \), let \( w^{\prime\prime} \in \mathrm{D} ( \mu ) \) be smallest such that \( u \cup w' \subseteq w^{\prime\prime} \) (i.e.\ \( w^{\prime\prime} \coloneqq u \cup w' \) if \( \lambda = \mu = \kappa \) or \( \lambda \leq \mu < \kappa \), while \( w^{\prime\prime} \coloneqq \sup \setofLR{\alpha < \kappa}{\alpha \in u \cup w'} \) if \( \lambda < \mu = \kappa \)).
As every \( \cof ( \mu ) \)-comeager subset of \( [ s ] \) must intersect every \( [ r ] \) with \( r \in \pre{w^{\prime\prime}}{ ( \kappa ) } \) and \( s \subseteq r \), and since \( ( \Sym ( \kappa ) , \tau_\mu) \) is \( \cof ( \mu) \)-Baire, then
\[
\forall ^*_\mu g \in [ s ] \left ( X \vDash \uptheta [ g^{-1} \circ h] \right ) \IFF \FORALL{ r \in \pre{ w ^{\prime\prime}}{ ( \kappa ) }} ( r \supseteq s \implies X \models \uptheta' [ r \restriction w' ] ) .
\]
Arguing as in the previous case, \( \phi ( t ) _\mu^{ * u } = M_{ \upvarphi , u} \) with \( \upvarphi ( \V_u ) \in \LL_{ \nu \kappa }^{\downarrow } \) being the formula 
\[ 
\Bigl ( \bigwedge\nolimits_{\substack{ \alpha , \beta \in u \\ \alpha < \beta } } \V_ \alpha \nequals \V_ \beta \Bigr ) \AND \forall \seqofLR{\V_ \alpha }{ \alpha \in w^{\prime\prime} \setminus u } \Bigl ( \bigwedge\nolimits_{ \substack{ \alpha , \beta \in w^{\prime\prime} \\ \alpha < \beta } } \V_ \alpha \nequals \V_ \beta \implies \uptheta' ( \V_{w'} ) \Bigr ) . 
\]
\end{itemize}
In either case let \( G ( t , u ) \coloneqq \upvarphi ( \V_u ) \).

\smallskip

\textbf{Case 2:} otherwise.
Let
\[
 \SUCC ( t ) \coloneqq \setofLR{t' \in T }{ t '\text{ is an immediate successor of } t }, 
\]
which is a set of size \( \leq \cof ( \mu ) \).
By inductive hypothesis, for every \( t' \in \SUCC ( t ) \) and every \( w \in \mathrm{D} ( \mu ) \) we have that \( G ( t' , w ) \in \LL_{ \nu\kappa}^{\downarrow } \) and \( \phi ( t ' )_\mu^{* w} = M_{ G ( t' , w ) , w } \).
As \( ( \Sym ( \kappa), \tau_\mu ) \) is \( \cof ( \mu ) \)-Baire, then
\[ 
\Bigl ( \bigcap\nolimits_{t' \in \SUCC ( t ) } \phi ( t ' ) \Bigr )_\mu^{* w} = \bigcap\nolimits_{t' \in \SUCC ( t ) } \phi ( t' ) ^{*w}_\mu = \bigcap\nolimits_{t' \in \SUCC ( t ) } M_{ G ( t' , w ) , w } = M_{ \uppsi_w , w } ,
\] 
where \( \uppsi_w ( \V_w ) \) is \( \bigwedge_{t' \in \SUCC ( t ) } G ( t' , w ) \).
Given \( X \in \Mod^\kappa_\LL \), the map
\[ 
\Sym ( \kappa ) \to \Mod^\kappa_\LL, \quad g \mapsto g . X
 \] 
is continuous when both spaces are endowed with the topology \( \tau_\lambda \), and hence it is also continuous as a function between \( ( \Sym ( \kappa ) , \tau_\mu ) \) and \( ( \Mod^\kappa_\LL, \tau_\lambda ) \), as \( \tau_\mu \) refines \( \tau_ \lambda \). 
Since \( \phi(t) \in \bB^{\mathrm{e}}_{\cof ( \mu ) + 1} ( \tau_\lambda ) \subseteq \bB_{\cof ( \mu ) + 1} ( \tau_\lambda ) \) and \( \bB_{\cof ( \mu ) + 1} \) is closed under continuous preimages, we get that
\begin{equation} \label{eq:BorelsubsetofSym}
\setofLR{g \in \Sym ( \kappa ) }{ g . X \in \phi(t) } \in \bB_{\cof ( \mu ) + 1} ( \Sym ( \kappa ) , \tau_\mu ) .
 \end{equation}

Now fix an arbitrary \( u \in \mathrm{D} ( \mu ) \) in order to define \( G ( t , u ) \). 
By~\eqref{eq:BorelsubsetofSym} and Proposition~\ref{prop:borelsetshaveBP}, for every \( s \in \pre{u}{(\kappa)} \) the set \( \setofLR{g \in [ s ] }{ g . X \in \phi(t) } \) has the \( \cof ( \mu ) \)-Baire property.
Using Proposition~\ref{prop:nonmeagerarelocallycomeager}, 
\begin{align*}
 ( X , s ) \in ( \phi ( t ) )_\mu^{*u} & \IFF (X,s) \in \bigl( \Mod^\kappa_\LL \setminus \bigcap\nolimits_{t' \in \SUCC ( t ) } \phi ( t ' ) \bigr )_\mu^{*u} 
 \\
& \IFF \FORALL{u \subseteq w \in \mathrm{D}(\mu) } \FORALL{ s \subseteq r \in \pre{w} { (\kappa) }} \Bigl [ ( X , r ) \notin \bigl ( \bigcap\nolimits_{t' \in \SUCC ( t ) } \phi ( t ' ) \bigr )_\mu^{* w} \Bigr ] 
\\
& \IFF \FORALL{u \subseteq w \in \mathrm{D}(\mu) } \FORALL{ s \subseteq r \in \pre{w} { (\kappa) }} \left [ ( X , r ) \notin M_{\uppsi_w , w } \right ] , 
\end{align*}
so that \( \phi ( t )_\mu ^{*u} = M_{\upvarphi, u } \) where \( \upvarphi( \V_u ) \) is the formula 
\[ 
\Bigl (\bigwedge\nolimits_{ \substack{\alpha , \beta \in u \\ \alpha < \beta }} ( \V_ \alpha \nequals \V_ \beta ) \Bigr) \AND \bigwedge\nolimits_{u \subseteq w \in \mathrm{D}(\mu)} \forall \seqofLR{ \V_ \alpha }{ \alpha \in w \setminus u } \Bigl [ \Bigl ( \bigwedge\nolimits_{ \substack{\alpha , \beta \in w\\ \alpha < \beta }} \V_ \alpha \nequals \V_ \beta \Bigr ) \implies \neg \uppsi_w ( \V_w ) \Bigr ] .
\]
Therefore it is enough to put \( G ( t , u ) \coloneqq \upvarphi ( \V_u ) \).
\end{proof}

Varying the parameters \( \mu \) and \( \lambda \) in Lemma~\ref{lem:LopezEscobar} and taking \( u = \emptyset \), we can now prove Propositions~\ref{prop:Borelformula} and~\ref{prop:Borelformula2}.

\begin{proof}[Proof of Proposition~\ref{prop:Borelformula}] 
First notice that if \( A \subseteq \Mod^\kappa_\LL \) is invariant under isomorphism, then for every infinite cardinal \( \mu \) belonging to the set in~\eqref{eq:conventionbasissingular} and every \( X \in \Mod^\kappa_\LL \) we have \( X \in A \IFF ( X , \emptyset ) \in A_\mu^{* \emptyset} \) by Theorem~\ref{th:Bairespacesym}.
Recall also from~\eqref{eq:5.3} that for every sentence \( \upsigma \in \LL_{\card{\pre{<\mu}{\kappa}}^+ \mu} \) and every \( X \in \Mod^\kappa_\LL \) we have \( (X,\emptyset) \in M_{\upsigma,\emptyset} \IFF X \in \Mod^\kappa_\upsigma \), and that \( \emptyset^{\downarrow \mu} = \emptyset \).
Then it is enough to set \( \mu = \lambda \) and \( u = \emptyset \) in Lemma~\ref{lem:LopezEscobar}, and then further set \( \lambda = \omega \) to get~\ref{prop:Borelformula-1} and \( \lambda = \kappa \) to get~\ref{prop:Borelformula-3}.
\end{proof}

\begin{proof}[Proof of Proposition~\ref{prop:Borelformula2}]
The proof is analogous to that of Proposition~\ref{prop:Borelformula}: it is enough to set \( \mu = \kappa \) and \( u = \emptyset \) in 
Lemma~\ref{lem:LopezEscobar}, and then further set \( \lambda = \omega \) for ~\ref{prop:Borelformula2-a} and \( \lambda = \kappa \) for~\ref{prop:Borelformula2-c}. 
The assumption \( \kappa^{< \lambda} = \kappa \) in part~\ref{prop:Borelformula2-b} guarantees that \( \card{ \mathcal{B}_ \lambda } = \kappa \), so that \( \tau_ \lambda \subseteq \bB^{\mathrm{e}}_{\kappa+1} ( \mathcal{B}_\lambda ) \), and therefore \( \bB_{\kappa+1}( \tau_\lambda ) = \bB^{\mathrm{e}}_{\kappa+1} ( \mathcal{B}_\lambda ) \).
Similarly \( \kappa^{< \kappa} = \kappa \) in part~\ref{prop:Borelformula2-c} ensures that \( \bB_{\kappa+1} ( \tau_b ) = \bB^{\mathrm{e}}_{\kappa+1} ( \mathcal{B}_b ) \).
\end{proof}

\begin{remark}\label{rmk:LopezEscobar} 
As for Corollary~\ref{cor:formulaborel} (see Remark~\ref{rmk:formulaborel}), if we drop the assumption \( \kappa^{< \kappa} = \kappa \) then also Proposition~\ref{prop:Borelformula2}\ref{prop:Borelformula2-c} may fail. 
Work in \( \ZFC \).
As observed in~\cite[Theorem 4.4]{Motto-Ros:2011qc}, if \( \mu \) is regular and \( \mu < \kappa \), then there are \( 2^{(2^\mu) } \)-many \( \tau_b \)-open subsets of \( \Mod_{\LL}^ \kappa \) invariant under isomorphism, while there are \( ( \kappa ^{< \kappa } )^ \kappa = 2^ \kappa \) formul\ae{} in \( \LL_{ \kappa ^+ \kappa } \).
Thus if there is e.g.\ a regular \( \mu < \kappa \) such that \( 2^{(2^\mu) } > 2^ \kappa \) (which can happen if \( \kappa^{< \kappa} > \kappa \)), then there is also an (effective) invariant \( \kappa + 1 \)-Borel subset of \( \Mod_\LL^ \kappa \) which cannot be of the form \( \Mod^ \kappa _{ \upsigma } \) for \( \upsigma \) an \( \LL_{ \kappa ^+ \kappa } \)-sentence. 
In particular, if \( 2^{\aleph_0} = 2^{\aleph_1} \), as it is the case in models of forcing axioms like \( \MA_{\omega_1} \) or \( \PFA \), then there is an invariant \( \tau_b \)-open subset of \( \Mod^{\omega_1}_\LL \) which is not of the form \( \Mod^{\omega_1}_\upsigma \) for any \( \upsigma \in \LL_{\omega_2 \omega_1} \) (\cite[Corollary 4.6]{Motto-Ros:2011qc}). 
A similar argument shows that also Proposition~\ref{prop:Borelformula2}\ref{prop:Borelformula2-b} may fail if we do not assume that \( \kappa^{< \lambda} = \kappa \).
\end{remark}

\section{\texorpdfstring{$ \kappa $-Souslin sets}{kappa-Souslin sets}} \label{sec:Ksouslinsets}
\subsection{Basic facts}
The following definition generalizes (in a different direction from the one considered in Section~\ref{subsec:kappa-analytic}) the notion of an analytic subset of a Polish space.

\begin{definition}\label{def:kappaSouslinset}
Let \( \kappa \) be an infinite cardinal, \( \pre{\omega}{ \kappa} \) be endowed with the product of the discrete topology on \( \kappa \), and \( X \) be a Polish space.
A set \( A \subseteq X \) is called \markdef{\( \kappa \)-Souslin} if it is a continuous image of a closed subset of \( \pre{\omega}{ \kappa} \), and is called \markdef{\( \infty \)-Souslin} if there is an infinite cardinal \( \kappa \) such that \( A \) is \( \kappa \)-Souslin.
The class of all \( \kappa \)-Souslin sets is denoted by \( \bS ( \kappa ) \), and \( \bS ( \infty ) \coloneqq \bigcup_{ \kappa \in \Cn } \bS ( \kappa ) \)\index[symbols]{Skappa@\( \bS ( \kappa ) \), \( \bS ( \infty ) \)} is the collection of all \( \infty \)-Souslin sets.
\end{definition}

Thus \( \bS ( \omega ) = \bSigma^1_1 \) and \( \bSigma^1_2 \subseteq \bS ( \omega_1 ) \), but the reverse inclusion \( \bS ( \omega_1 ) \subseteq \bSigma^1_2 \) depends on the axioms we assume: it is true in models of \( \AD + \DC \), but it is false in models of \( \AC + \CH \). 
Under choice every subset of a Polish space is \( 2^{\aleph_0} \)-Souslin (see Proposition~\ref{prop:souslincardinalsunderAC}), so the notion of an \( \infty \)-Souslin set makes sense only if we work in models where \( \AC \) fails.

The class \( \bS ( \kappa ) \) is a hereditary boldface pointclass, it is closed under countable unions and countable intersections (assuming \( \AC_\omega \)), it contains all Borel sets, and it is closed under images and preimages of Borel functions (Lemma~\ref{lem:Souslin}). 
In particular, \( \bS ( \kappa ) \supseteq \bSigma^1_1 \) for every infinite \( \kappa \). 
Moreover, every \( \kappa \)-Souslin set is automatically \( \kappa' \)-Souslin for every \( \kappa' \geq \kappa \), so that \( \bS ( \kappa ) \subseteq \bS ( \kappa') \) (although it is not true in general that \( \kappa < \kappa' \IMPLIES \bS ( \kappa ) \subset \bS ( \kappa' ) \): by Corollary~\ref{cor:SouslincardinalsunderAC} this is true under \( \AC \) for all \( \kappa' \leq 2^{\aleph_0} \), but under \( \AD + \DC \) we e.g.\ have \( \bS ( \omega_1 ) = \bS ( \omega_2 ) \) --- see Section~\ref{sec:SouslinunderAD}). 
However, it is maybe worth noticing that the pointclasses \( \bS ( \kappa) \), \( \kappa \in \Cn \), may not form a well behaved hierarchy in models of \( \AC \) as it may happen that \( \check{\bS} ( \kappa ) \nsubseteq \bS ( \kappa' ) \) for \( \omega_1 \leq \kappa < \kappa' \) (see Remark~\ref{rmk:pathology}). 
This pathology is absent in models of \( \AD \), where we have
\[ 
 \bS ( \kappa) \cup \check{\bS} ( \kappa ) \subseteq \bS ( \kappa' ) \cap \check{\bS} ( \kappa' ) 
 \] 
for all Souslin cardinals \( \kappa < \kappa' \) (see Definition~\ref{def:Souslincardinal}).

\begin{remark}
Using standard arguments and the closure properties of \( \bS(\kappa) \) mentioned above, the notion of \( \kappa \)-Souslin set may be reformulated in several ways. 
For any subset \( A \) of a Polish spaces \( X \) and any infinite cardinal \( \kappa \) the following are equivalent:
\begin{itemize}[leftmargin=1pc]
\item
\( A \) is \( \kappa \)-Souslin;
\item
\( A \) is either empty or a continuous image of \( \pre{\omega}{\kappa} \);
\item
\( A \) is a continuous image of a Borel subset of \( \pre{\omega}{\kappa} \);
\item
\( A = \PROJ F \) for some closed \( F \subseteq X \times \pre{\omega}{\kappa} \);
\item
\( A = \PROJ B \) for some Borel \( B \subseteq X \times \pre{\omega}{\kappa} \),
\end{itemize}
where \( \PROJ \) is the projection map defined in~\eqref{eq:projectionCh2} on page~\pageref{eq:projectionCh2}.
\end{remark}

Since \( \bS ( \kappa ) \) is closed under countable unions, \( \kappa \)-Souslinness becomes trivial when considering a countable Polish space \( X \), as any \( A \subseteq X \) is countable, and therefore \( A \in \bS ( \omega ) \subseteq \bS ( \kappa ) \) for every infinite cardinal \( \kappa \). 
Moreover, any two uncountable Polish spaces are Borel isomorphic, the closure under images and preimages of Borel functions of \( \bS ( \kappa ) \) yields that it is enough to study \( \kappa \)-Souslin subsets of some specific uncountable Polish space.
\emph{Therefore, unless otherwise specified, from now on we confine our analysis to \( \kappa \)-Souslin subsets of (countable products of) \( \pre{ \omega}{2} \)}.
The choice of these canonical spaces is motivated by the fact that the \( \kappa \)-Souslin subsets of \( ( \pre{ \omega}{2} )^N \) (for \( 1 \leq N \leq \omega \)) admit a particularly nice representation in terms of projections of the trees introduced in the following

\begin{notation}
Recall from~\eqref{eq:Tr} on page~\pageref{eq:Tr} that \( \Tr ( Y ) \) is the set of all descriptive set-theoretic trees (of height \( \leq \omega \)) on \( Y \).
When \( Y = ( \smash{\underbrace{ 2 \times 2 \times\cdots }_{N}} ) \times \kappa \) we often write \( \Tr ( N ; 2 , \kappa ) \) instead of \( \Tr ( Y ) \).
\end{notation}

\begin{remark} \label{rmk:cardinalityS(kappa)}
As \( Y \coloneqq ( \smash{\underbrace{ 2 \times 2 \times\cdots }_{N}} ) \times \kappa \) has size \( \kappa \), the elements of \( \Tr ( Y ) \) can be identified (in \( \ZF \)) with elements of \( \pre{\kappa}{2} \).
\end{remark}

The usual representation of analytic (i.e.\ \( \omega \)-Souslin) sets (see~\cite[Proposition 25.2]{Kechris:1995zt}) extends to an arbitrary cardinal \( \kappa \), so that for every \( A \subseteq ( \pre{ \omega}{2} )^N \)
\[
A \in \bS ( \kappa ) \IFF \EXISTS{ T \in \Tr ( N ; 2 , \kappa ) } \bigl ( A = \PROJ \body{T} \bigr ) .
\]
A similar tree representation holds for \( \kappa \)-Souslin subsets of \( \pre{\omega}{\omega} \), or more generally, of any countable power \( \pre{\omega}{X} \) of a discrete countable set \( X \). 
As hinted in Section~\ref{subsubsec:pwoscaleproperty}, using such representation and the notion of leftmost branch through trees on well-orderable sets, one can easily reformulate the notion of a \( \kappa \)-Souslin set in terms of scales (see e.g.\ \cite[Lemma 2.5]{Jackson:2010ff}):
\begin{fact} \label{fct:Souslinscale}
A set \( A \subseteq \pre{\omega}{X} \) is \( \kappa \)-Souslin if and only if it admits a scale all of whose norms map into \( \kappa \).
\end{fact}

The fact that a set \( A \subseteq ( \pre{ \omega}{2} )^N \) is \( \kappa \)-Souslin gives us some structural information on it: for example, we have the following property, which is meaningful whenever \( \kappa \) is small enough compared to \( \pre{ \omega}{2} \) (see~\cite[Theorem 2C.2]{Moschovakis:2009fk}).

\begin{proposition} \label{prop:kappaPSP}
Let \( \kappa \) be an infinite cardinal and let \( A \subseteq ( \pre{ \omega}{2} ) ^N \). 
If \( A \in \bS ( \kappa ) \) then \( A \) has the \markdef{\( \kappa \)-Perfect Set Property} (\markdef{\( \kappa \)-\( \PSP \)} for short): either \( A \) has at most \( \kappa \)-many elements, or else it contains a perfect set (equivalently, a homeomorphic copy of \( \pre{\omega}{2} \)). 
\end{proposition}

This property is used in Section~\ref{sec:souslinunderAC} to show that the pointclass \( \bS ( \kappa ) \) is not trivial for small \( \kappa \). 
Another way to obtain nontrivial pointclasses is to relativize \( \kappa \)-Sousliness to some inner model, as explained in the following remark --- this approach is exploited in Section~\ref{subsubsec:reducibilityininnermodel}.

\begin{remarks}\label{rmk:Souslinininnermodel} 
\begin{enumerate-(i)}
\item
If \( W \) is an inner model (i.e. a transitive proper class model of \( \ZF \)) and \( \kappa \) is a cardinal in it, by absoluteness of existence of infinite branches, the relativization%
\footnote{To simplify the notation, we consider here only subsets of \( \pre{\omega}{2} \): the generalization of this notion to subsets of countable products \( (\pre{\omega}{2})^N \) (and to arbitrary Polish spaces) is straightforward.}
 of \( \bS ( \kappa ) \) to \( W \) is
\[
\left ( \bS ( \kappa ) \right )^W = \setofLR{A \cap (\pre{\omega}{2}) ^W }{ A \in \bS_W ( \kappa ) } \supseteq \bS_W ( \kappa ) \cap \pow ( (\pre{\omega}{2})^W ) 
\]
where%
\footnote{In the definition of \( \bS_W(\kappa) \), the tree \( T \) must belong to \( W \) but its projection \( \PROJ \body{T} \) is computed in \( \Vv \).}
\[
\bS_W ( \kappa ) \coloneqq \setofLR{\PROJ \body{T} }{ T \text{ is a tree on \( 2 \times \kappa \) belonging to } W }.
\]
In particular, if \( \pre{\omega}{2} \subseteq W \) (and hence \( \Ll ( \R ) \subseteq W \)) then
\[
\left ( \bS ( \kappa ) \right )^W = \bS_W ( \kappa ) .
\]
\item
If \( \pre{\omega}{2} \subseteq W \) (equivalently, \( \Ll ( \R ) \subseteq W \)), then \( \bS_W ( \kappa ) \) is a boldface pointclass in \( \Vv \) because any continuous function belongs to \( \Ll ( \R ) \subseteq W \). 
Moreover, \( \bS_W ( \kappa ) \) contains all Borel sets and is closed under continuous images, so that \( \bSigma^1_1 \subseteq \bS_W ( \kappa ) \subseteq \bS ( \kappa) \).
\end{enumerate-(i)}
\end{remarks}

\subsection{More on Souslin sets and Souslin cardinals}

The following lemma collects some well-known facts on \( \bS ( \kappa ) \). 

\begin{lemma}[\( \AC_\omega \)]\label{lem:Souslin}
Let \( \kappa \) be an infinite cardinal. 
\begin{enumerate-(a)}
\item\label{lem:Souslin-a}
\( \bS ( \kappa ) \) is a hereditary boldface pointclass containing all closed and open sets, closed under countable unions and countable intersections (and hence containing all Borel sets), and closed under projections (equivalently, under continuous images). 
In particular, \( \bS ( \kappa ) \) contains all analytic sets.
\item\label{lem:Souslin-b}
\( \bS ( \kappa ) \) is closed under images and preimages of (partial) functions with \( \kappa \)-Souslin graph.
\end{enumerate-(a)}
\end{lemma}

Notice that the assumption \( \AC_\omega \) can be relaxed to \( \AC_\omega ( \R ) \) if \( \AD \) is assumed and \( \kappa < \Theta \). 
Moreover, since a function between two Polish spaces is Borel if and only if its graph is \( \omega \)-Souslin (i.e.~analytic), Lemma~\ref{lem:Souslin} implies that every pointclass \( \bS ( \kappa ) \) is closed under Borel images and Borel preimages. 

\begin{proof}
The proof of~\ref{lem:Souslin-a} is standard --- see~\cite[Theorem 2B.2]{Moschovakis:2009fk}.

\smallskip

\ref{lem:Souslin-b}
Let \( X , Y \) be Polish spaces, \( A \subseteq X \) be \( \kappa \)-Souslin, and \( f \) be a partial function from \( X \) to \( Y \) with \( \kappa \)-Souslin graph.
We must show that
\[ 
B \coloneqq \setofLR{y \in Y}{\EXISTS{x \in A \cap \dom ( f ) } ( f ( x ) = y ) } \in \bS ( \kappa ).
 \] 
We can assume without loss of generality that both \( A \) and \( \gr ( f ) \) are nonempty, so let \( g \colon \pre{\omega} {\kappa} \to X \) and \( h \colon \pre{\omega}{\kappa} \to X \times Y \) be continuous surjections onto \( A \) and \( \gr ( f ) \), respectively. 
For \( i = 0 , 1 \), let \( \pi_i \colon \pre{\omega}{\kappa} \onto \pre{\omega}{\kappa} \), \( s \mapsto \seqofLR{ s ( 2 n + i ) }{ n \in \omega } \), and let \( \pi_X \) and \( \pi_Y \) be the projections of \( X \times Y \) onto the spaces \( X \) and \( Y \), respectively. 
Set
\[ 
C \coloneqq \setofLR{ s \in \pre{\omega}{\kappa}}{ g ( \pi_0 ( s ) ) = \pi_X ( h ( \pi_1 ( s ) ) ) }.
 \] 
Since all the functions involved in its definition are continuous and \( X \) is Hausdorff, \( C \subseteq \pre{\omega}{\kappa} \) is a closed set, and hence there is a continuous \( r \colon \pre{\omega}{\kappa} \to C \) such that \( r \restriction C \) is the identity map by~\cite[Proposition 2.8]{Kechris:1995zt}. 
Then \( \pi_Y \circ h \circ \pi_1 \circ r \colon \pre{\omega}{\kappa} \onto B \) is continuous, whence \( B \in \bS ( \kappa ) \).
\end{proof}

\begin{remark}\label{rmk:souslin}
If \( \kappa \) is an infinite cardinal and we assume \( \AC_\kappa \), then \( \bS ( \kappa ) \) is further closed under unions of length (at most) \( \kappa \). 
This is because if \( \seqofLR{ A_\alpha }{ \alpha < \kappa } \) is a sequence of \( \kappa \)-Souslin subsets of a Polish space \( X \), using \( \AC_\kappa \) we can chose for each \( \alpha < \kappa \) a continuous surjection \( p_\alpha \colon \pre{\omega}{\kappa} \onto A_\alpha \subseteq X \).
The surjection
\[ 
\pre{\omega}{\kappa} \onto \bigcup_{\alpha < \kappa} A_\alpha , \quad \seqLR{ \alpha _0 , \alpha _1 , \alpha _2 , \dots } \mapsto p_{ \alpha _0} ( \seqLR{ \alpha _1 , \alpha _2 , \dots } )
\] 
is continuous, whence \( \bigcup_{\alpha < \kappa} A_\alpha \in \bS ( \kappa ) \). 
In fact, assuming \( \AC \) one can further show that for every \( 0 \neq n \in \omega \) and every subset \( A \) of a Polish space \( X \), \( A \) is \( \omega_n \)-Souslin if and only if \( A \) is a union of \( \aleph_n \)-many Borel sets (see e.g.\ \cite[Proposition 13.13(f)]{Kanamori:2003fu}).
\end{remark}

The following lemma collects some classical facts relating \( \kappa \)-Sousliness to \( \kappa+1 \)-Borelness.

\begin{lemma}[Folklore]\label{lem:Souslinfolklore}
Every \( A \in \bS ( \kappa ) \) is effective \( \kappa^{+} + 1 \)-Borel (in fact, \( A \) is an effective \( \kappa^+ \)-union of effective \( \kappa + 1 \)-Borel sets), and if \( \kappa \) has uncountable cofinality then \( A \) is an effective \( \kappa \)-union of effective \( \kappa \)-Borel sets (thus it is effective \( \kappa + 1 \)-Borel).
In particular, \( \bS ( \kappa ) \subseteq \bB^{(\mathrm{e})}_{\kappa^+ + 1} \) and if \( \kappa \) has uncountable cofinality then also \( \bS ( \kappa ) \subseteq \bB^{(\mathrm{e})}_{\kappa + 1} \). 
Moreover, \( \bDelta_{\bS ( \kappa )} \subseteq \bB_{\kappa + 1} \) (independently of \( \kappa \)).
\end{lemma}

\begin{proof}
Use the characterization of \( \kappa \)-Souslin subsets of \( \pre{ \omega}{2} \) in terms of projection of trees in \( \Tr ( 2 \times \kappa) \) (see~\cite[Lemma 2.12]{Jackson:2008pi}). 
For the part concerning \( \bDelta_{\bS ( \kappa )} \), use the classical Lusin separation argument (see~\cite[Theorem 2E.2]{Moschovakis:2009fk}).
\end{proof}

\begin{proposition}\label{prop:Skappainthecodesareborel}
Suppose that there is an \( \bS ( \kappa ) \)-code for \( \kappa \), so that it makes sense to speak of a \( \bS ( \kappa ) \)-in-the-codes Borel function from \( \pre{\omega}{2} \) to \( \pre{\kappa}{2} \).
Then every \( \bS ( \kappa ) \)-in-the-codes function \( f \colon \pre{\omega}{2} \to \pre{\kappa}{2} \) is weakly \( \kappa+1 \)-Borel.
\end{proposition}

\begin{proof}
By Lemma~\ref{lem:Souslin}, \( \bS(\kappa)\) is closed under projections and finite intersections (notice that for these closure properties \( \AC_\omega \) is not needed). 
By~\eqref{eq:inthecodes} on page~\pageref{eq:inthecodes}, for every \( \alpha < \kappa \) and every \( i \in \set{ 0 , 1 } \), the set \( f^{-1} ( \widetilde{\Nbhd}\vphantom{\Nbhd}^\kappa_{\alpha , i}) \) is in \( \bDelta_{\bS(\kappa)} \). 
Since by Lemma~\ref{lem:Souslinfolklore} every set in \( \bDelta_{\bS ( \kappa )} \) is \( \kappa + 1 \)-Borel, this means that each \( f^{-1} ( \widetilde{\Nbhd}\vphantom{\Nbhd}^\kappa_{\alpha , i}) \) is \( \kappa+1 \)-Borel.
Since the sets \( \widetilde{\Nbhd}\vphantom{\Nbhd}^\kappa_{\alpha , i} \) generates (by taking \emph{finite} intersections) the canonical basis \( \mathcal{B}_p \) for the product topology \( \tau_p \) on \( \pre{\kappa}{2} \), and since \( \mathcal{B}_p \), being in bijection with \( \forcing{Fn} ( \kappa , 2 ; \omega ) \), is a well-orderable set size \( \kappa \), this implies that \( f^{-1} ( U ) \) is \( \kappa+1 \)-Borel for every \( \tau_p \)-open \( U \subseteq \pre{\kappa}{2} \), i.e.\ that \( f \) is weakly \( \kappa +1 \)-Borel. 
\end{proof}

\begin{definition}\label{def:Souslincardinal}
A cardinal \( \kappa \) is a \markdef{Souslin cardinal} if \( \bS ( \kappa ) \setminus \bigcup_{\lambda < \kappa} \bS ( \lambda) \neq \emptyset \), i.e.~if there is a \( \kappa \)-Souslin set which is not \( \lambda \)-Souslin for any \( \lambda < \kappa \).
Let 
\[
 \Xi \coloneqq \sup \setofLR{ \kappa \in \Cn }{ \kappa \text{ is a Souslin cardinal}}. \index[symbols]{Xi@\( \Xi \)}
\]
\end{definition}

By definition \( \bS ( \infty ) = \bS ( \Xi ) \), and by Lemma~\ref{lem:Souslinfolklore} \( \bS ( \infty ) \subseteq \bB_\infty \).
The converse inclusion may fail, for example \( \AD + {\Vv = \Ll ( \R )} \) implies that \( \bS ( \infty ) = \bSigma^2_1 \) and that every set of reals is \( \infty \)-Borel.
The next folklore result shows that \( \Xi \leq \Theta \).

\begin{lemma} \label{lem:Souslincardinals<Theta}
If \( \kappa \) is a Souslin cardinal, then \( \kappa < \Theta \).
In particular, \( \Theta \) is not a Souslin cardinal.
\end{lemma}

\begin{proof}
Let \( A \subseteq \pre{\omega}{\omega} \) witness that \( \kappa \) is a Souslin cardinal.
Then by Fact~\ref{fct:Souslinscale} there is a scale \( \seqofLR{ \rho_n }{ n \in \omega } \) on \( A \) such that each \( \rho_n \) maps into a cardinal \( \kappa_n \leq \kappa \). 
Let \( \lambda \coloneqq \sup_{n \in \omega } \kappa_n \leq \kappa \). 
If \( \lambda < \kappa \) then \( A \) would be \( \lambda \)-Souslin by Fact~\ref{fct:Souslinscale} again, contradicting the choice of \( A \): then \( \lambda = \kappa \), and the map \( \varphi \colon \pre{\omega}{\omega} \to \kappa \) sending each \( n \conc x \) to \( \rho_n ( x ) \) if \( x \in A \) and to \( 0 \) otherwise is surjective, whence \( \kappa < \Theta \).
\end{proof}

As \( \bSigma^1_2 \setminus \bSigma^1_1 \neq \emptyset \) and as \( \bSigma^1_1 = \bS ( \omega) \) and \( \bSigma^1_2 \subseteq \bS ( \omega_1) \), then \( \omega_1 \) is a Souslin cardinal.
The existence of a Souslin cardinal greater than \( \omega_1 \), and whether \( \Xi \) is a Souslin cardinal (and hence the largest Souslin cardinal) or, equivalently, whether \( \bS ( \Xi ) \neq \bigcup_{\kappa < \Xi } \bS ( \kappa ) \), both depend on the underlying set-theoretic assumptions --- see Sections~\ref{sec:souslinunderAC} and~\ref{sec:SouslinunderAD}.

\subsection{Souslin sets and cardinals in models with choice} \label{sec:souslinunderAC}
Under choice \( \Xi \leq 2^{\aleph_0} \) because every set of reals is trivially \( 2^{\aleph_0} \)-Souslin (Proposition~\ref{prop:souslincardinalsunderAC}). 
Therefore \( \CH \) implies that \( \Xi = \omega_1 \) (and hence \( \Xi \) is a Souslin cardinal) and that \( \bS ( \Xi) = \bS ( \omega_1) = \pow ( \pre{ \omega}{2} ) \), so the pointclass \( \bS ( \omega_1) \) and the notion of \( \omega_1 \)-Souslin set are uninteresting in this context.
Assuming instead \( \neg \CH \), the condition \( \kappa < 2^{\aleph_0} \) guarantees that \( \bS ( \kappa ) \) is a proper pointclass.

\begin{proposition}[\( \AC \)]\label{prop:Bernstein}
\( \pow ( \pre{\omega}{2} ) \setminus \bigcup_{ \kappa < 2^{\aleph_0}} \bS ( \kappa ) \neq \emptyset \), so \( 2^{\aleph_0} \) is a Souslin cardinal.
\end{proposition}

\begin{proof}
The standard construction of a Bernstein set~\cite[p.\ 38]{Kechris:1995zt} yields a set \( B \subseteq \pre{\omega}{2} \) of size \( 2^{\aleph_0} \) which does not contain any perfect set, and hence \( B \notin \bigcup_{ \kappa < 2^{\aleph_0}} \bS ( \kappa ) \) by Proposition~\ref{prop:kappaPSP}.
\end{proof}

The next result shows that in models of \( \AC \) all cardinals smaller than the continuum are Souslin cardinals, whence 
\[ 
\FORALL{ \omega \leq \kappa , \kappa' < 2^{ \aleph_0 }} \left ( \kappa \leq \kappa' \IFF \bS ( \kappa ) \subseteq \bS ( \kappa' ) \right ).
\] 
\begin{proposition} \label{prop:souslincardinalsunderAC}
Let \( \kappa \) be an infinite cardinal.
If \( A \subseteq \pre{\omega}{2} \) is an infinite set of size \( \kappa \), then \( A \in \bS ( \kappa ) \). 
If moreover \( \card{\pre{\omega}{2}} \nleq \kappa \), then \( A \in \bS ( \kappa ) \setminus \bigcup_{\lambda<\kappa} \bS ( \lambda) \). 
In particular, if \( \kappa \) is a small cardinal such that \( \card{\pre{\omega}{2}} \nleq \kappa \) and we assume \( \AC_\kappa ( \R ) \), then \( \kappa \) is a Souslin cardinal.
\end{proposition}

\begin{proof} 
Given \( A = \setofLR{x_\alpha}{\alpha < \kappa} \subseteq \pre{\omega}{2}\), let 
\[ 
T \coloneqq \setof{ ( x_\alpha \restriction n, \alpha^{ ( n ) } )}{\alpha < \kappa , n \in \omega},
\] 
where \( \alpha^{ ( n ) } \) is the sequence formed by \( n \)-many \( \alpha \)'s. 
It is straightforward to check that \( A = \PROJ \body{T} \), so that \( A \in \bS ( \kappa ) \).

Assume now that \( \card{\pre{\omega}{2}} \nleq \kappa \), so that \( A \) cannot contain a perfect set. 
If \( A \in \bS ( \lambda) \) for some infinite \( \lambda < \kappa \), then the set \( A \) would have cardinality \( \leq \lambda \) by the \( \lambda \)-\( \PSP \) (Proposition~\ref{prop:kappaPSP}), a contradiction. 
Thus \( A \in \bS ( \kappa ) \setminus \bigcup_{\lambda < \kappa} \bS ( \lambda) \).

Finally let \( \kappa \) be an infinite small cardinal, and let \( p \colon \pre{\omega}{2} \onto \kappa \) be a surjection. 
Then by \( \AC_\kappa ( \R ) \) there is a choice function for the sequence \( \seqofLR{ p^{ - 1 } ( \alpha) }{ \alpha \in \kappa } \), and the range \( A \) of such a function is a subset of \( \pre{\omega}{2} \) of cardinality \( \kappa \): the result then follows from the first part of the proposition.
\end{proof}

\begin{corollary}[\( \AC \)]\label{cor:SouslincardinalsunderAC}
All infinite \( \kappa \leq 2^{\aleph_0} \) are Souslin cardinals.
In particular, \( \Xi = 2^{\aleph_0} \).
\end{corollary}

Thus under e.g.\ \( \ZFC + \PFA \), the Souslin cardinals are exactly \( \omega \), \( \omega_1 \), and \( \omega_2 = \Xi = 2^{\aleph_0} \).
As for the cardinality of \( \bS ( \kappa ) \), notice that under choice \( \cardLR{\bS ( \kappa ) } = 2^\kappa \).
In fact \( \cardLR{\bS ( \kappa ) } \leq 2^\kappa \) by Remark~\ref{rmk:cardinalityS(kappa)}; for the other inequality use the fact that all subsets of a given \( A \subseteq \pre{\omega}{2} \) of size \( \kappa \) are \( \kappa \)-Souslin by Proposition~\ref{prop:souslincardinalsunderAC}.
Thus it may happen that \( \cardLR{\bS ( \kappa ) } < \cardLR{\bS ( \kappa' ) } \) for some \( \kappa < \kappa' \leq \Xi \).

The following proposition computes the value of \( \bdelta_{ \bS ( \kappa)} \) (for various \( \kappa \leq 2^{\aleph_0} \)) in models of \( \AC \), and should be contrasted with Proposition~\ref{prop:deltaSkappa} providing analogous computations in the \( \AD \)-context.

\begin{proposition}[\( \AC \)] \label{prop:deltaSkappaAC}
Let \( \kappa \leq 2^{\aleph_0} \) be an infinite cardinal. 
\begin{enumerate-(a)}
\item \label{prop:deltaSkappaAC-0}
\( \bdelta_{\bS ( \kappa )} \leq \kappa^+ \). 
\item \label{prop:deltaSkappaAC-a}
\( \bdelta_{ \bS ( \omega)} = \omega_1 \), \( \bdelta_{\bS ( \omega_1 )} = \omega_2 \) and \( \bdelta_{ \bS ( 2^{\aleph_0})} = (2^{\aleph_0})^+ = \Theta \). 
In particular, \( \bdelta_{ \bS ( \omega_2)} \geq \omega_2 \).
\item \label{prop:deltaSkappaAC-c}
There are models \( M \) of \( \ZFC \) in which \( 2^{\aleph_0} \) is as large as desired and \( \bdelta_{ \bS ( \kappa)} = \kappa^+ \) for all infinite \( \kappa \leq 2^{\aleph_0} \).
\end{enumerate-(a)} 
\end{proposition}

\begin{proof}
\ref{prop:deltaSkappaAC-0} 
This directly follows from the Kunen-Martin's theorem~\cite[Theorem 2G.2]{Moschovakis:2009fk}, which holds in \( \ZF \). 

\ref{prop:deltaSkappaAC-a} 
\( \bdelta_{ \bS ( \omega )} = \omega_1 \) is trivial, so let us show that \( \bdelta_{ \bS ( \omega_1 ) } = \omega_2 \). 
It is enough to show that for every \( \omega_1 \leq \alpha < \omega_2 \) there is a \( \bDelta_{ \bS ( \omega_1)} \) prewellordering \( \preceq \) of length \( \alpha\) of \( \LO \subseteq \pre{\omega \times \omega}{2} \), the set of codes for linear orderings of \( \omega \)~\cite[Section 27.C]{Kechris:1995zt}. 
For every \( \omega \leq \beta < \omega_1 \), let \( \WO_\beta \) be the set of codes for well-orders on \( \omega \) of length \( \beta \), and let \( \mathrm{NWO} \) be the set of codes for non-well-founded linear orderings of \( \omega \). 
Notice that \( \LO \) is the disjoint union of the \( \WO_\beta \)'s and \( \mathrm{NWO} \), which are all \( \bSigma^1_1 \) sets. 
Given \(\alpha\) as above, let \( i \colon \alpha \to \omega_1 \setminus \omega \) be a bijection and for every \( x , y \in \LO \) set
\[ 
x \preceq y \IFF x \in \mathrm{NWO} \OR \EXISTS{ \beta \leq \gamma < \alpha } ({x \in \WO_{i ( \beta ) }} \AND {y \in \WO_{ i ( \gamma )}}).
 \] 
The relation \( \preceq \) is a prewellordering of \( \LO \) of length \( \alpha \). 
Using that \( \WO_\beta , \mathrm{NWO} \in \bSigma^1_1 \subseteq \bS ( \omega_1 ) \) and that under \( \AC \) the pointclass \( \bS ( \omega_1 ) \) is closed under unions of size \( \omega_1 \) (Remark~\ref{rmk:souslin}), it follows that \( \preceq \) is \( \bDelta_{ \bS ( \omega_1)} \), as required. 

Finally, to show \( \bdelta_{ \bS ( 2^{\aleph_0})} = ( 2^{\aleph_0} )^+ \), fix \( 2^{\aleph_0} \leq \alpha < ( 2^{\aleph_0} )^+ \) and let \( i \colon \pre{\omega}{2} \to \alpha \) be a bijection: then the relation
\[ 
x \preceq y \IFF i ( x ) \leq i ( y )
 \] 
is a prewellordering of \( \pre{\omega}{2} \) of length \(\alpha\), and is in \( \bDelta_{ \bS ( 2^{\aleph_0})} \) because every subset of \( \pre{\omega}{2} \times \pre{\omega}{2} \) is trivially \( 2^{\aleph_0} \)-Souslin.

\ref{prop:deltaSkappaAC-c} By part~\ref{prop:deltaSkappaAC-a} it is enough to consider the case \( \omega_1 \leq \kappa < 2^{\aleph_0} \). 
By~\cite[Theorem B]{Harrington:1977aa}, there is a model \( M \) of \( \ZFC \) in which \( 2^{\aleph_0} \) is as large as desired and every set \( A \subseteq \pre{\omega}{2} \) of size \( < 2^{\aleph_0} \) is \( \bPi^1_2 \): we claim that \( M \) is as required. 
Let \( A \subseteq \pre{\omega}{2} \) be a set of size \( \kappa \). 
Fix \( \kappa \leq \alpha < \kappa^+ \) and a bijection \( i \colon A \to \alpha \). 
Then the binary relation \( \preceq \) on \( \pre{\omega}{2} \) defined by
\[ 
x \preceq y \iff x \notin A \vee (x,y \in A \wedge { i ( x ) \leq i ( y ) })
 \] 
is a prewellordering of \( \pre{\omega}{2} \) of length \(\alpha\). 
Since \( \bS ( \kappa) \) is closed under unions of size \( \kappa \), to show that \( \preceq \) is \( \bDelta_{ \bS ( \kappa)} \) it is enough to check that \( \pre{\omega}{2} \setminus A \in \bS ( \kappa) \). 
But this follows from the fact that since \( \card{A} < 2^{\aleph_0} \), our choice of \( M \) ensures \( \pre{\omega}{2} \setminus A \in \bSigma^1_2 \subseteq \bS ( \omega_1) \subseteq \bS ( \kappa) \).
\end{proof}

\begin{remark} \label{rmk:pathology}
Although it is always the case that%
\footnote{In fact, \( \check{\bS}(\omega) = \bPi^1_1 \subseteq \bSigma^1_2 \subseteq \bS(\omega_1) \).}
 \( \check{\bS}(\omega) \subseteq \bS(\omega_1) \) and \( \bS ( \kappa) \subseteq \bS ( \kappa') \) for all infinite \( \kappa \leq \kappa' \leq 2^{\aleph_0} \), there are models of \( \ZFC \) in which \( 2^{\aleph_0} \) is as large as desired but for every cardinal \( \kappa \) such that \( \kappa^+ < 2^{\aleph_0} \) there is \( A \subseteq \pre{\omega}{2} \) in \( \bS ( \omega_1) \) such that \( \pre{\omega}{2} \setminus A \notin \bS ( \kappa) \), so that in general \( \check{\bS}(\kappa) \not\subseteq \bS ( \kappa^+) \) when \( \kappa > \omega \). 
To see this let \( M \) be a model of \( \ZFC \) as in the proof of Proposition~\ref{prop:deltaSkappaAC}\ref{prop:deltaSkappaAC-c}, fix \( \kappa \) as above, and let \( A \subseteq \pre{\omega}{2} \) be such that \( \card{\pre{\omega}{2} \setminus A } = \kappa^+ \). 
Then \( A \in \bSigma^1_2 \subseteq \bS ( \omega_1) \) by our choice of \( M \), but the set \( \pre{\omega}{2} \setminus A \) does not belong to \( \bS ( \kappa) \) because it cannot satisfy the \( \kappa \)-\( \PSP \) by \( \kappa < \card{\pre{\omega}{2} \setminus A} < 2^{\aleph_0} \) (see Proposition~\ref{prop:kappaPSP}).
\end{remark}

Assuming \( \FORALL{ x \in \pre{\omega}{\omega} } ( x^\# \text{ exists}) \), the assumption \( \bdelta^1_2 = \omega _2 \) has high consistency strength.
The next result shows that this is not the case if we drop the existence of sharps.

\begin{proposition}[\( \AC \)]\label{prop:Somega1=sigma12}
Assume \( \MA + \neg \CH + \EXISTS{a \in \pre{\omega}{\omega}} (\omega_1^{\Ll [ a ] } = \omega_1) \). 
Then \( \bS(\omega_1) = \bSigma^1_2 \), and hence \( \bdelta^1_2 = \omega_2 \) by Proposition~\ref{prop:deltaSkappaAC}\ref{prop:deltaSkappaAC-a}.
\end{proposition}

\begin{proof}
As already observed, \( \bSigma^1_2 \subseteq \bS(\omega_1) \). 
For the converse, recall from Remark~\ref{rmk:souslin} that each \( A \in \bS(\omega_1) \) is a union of \( \aleph_1 \)-many Borel sets. 
Since Martin and Solovay showed in~\cite{Martin:1970ms} that under the assumptions above every union of at most \( \aleph_1 \)-many Borel sets is in \( \bSigma^1_2 \), we get \( A \in \bSigma^1_2 \) and we are done.
\end{proof}

Let us now turn our attention to \( \bS(\kappa) \)-in-the-codes functions in models with choice.
By Proposition~\ref{prop:deltaSkappaAC}\ref{prop:deltaSkappaAC-a} and Remark~\ref{rmk:inthecodes}\ref{rmk:inthecodes-ii}, in models of \( \ZFC \) it always makes sense to speak of \( \bS ( \kappa ) \)-in-the-codes functions \( f \colon \pre{\omega}{2} \to \pre{\kappa}{2} \) (see Definition~\ref{def:Gammainthecodes}) for \( \kappa = \omega \), \( \kappa = \omega_1 \), or \( \kappa = 2^{\aleph_0} \), but in some models this could well be the case for all \( \kappa \leq 2^{\aleph_0} \), see Proposition~\ref{prop:deltaSkappaAC}\ref{prop:deltaSkappaAC-c}. 
Notice also that it always makes sense to speak of \( \bSigma^1_2 \)-in-the-codes (and hence \( \bS ( \omega_1) \)-in-the-codes) functions \( f \colon \pre{\omega}{2} \to \pre{\omega_1}{2} \) \emph{even when working in \( \ZF+\DC \)}. 
This is because \( \bdelta_{\bS(\omega_1)} \geq \bdelta^1_2 \geq \omega_1 \) (by \( \bS ( \omega) = \bSigma^1_1 \subseteq \bSigma^1_2 \subseteq \bS(\omega_1) \)), and the norm \( \rho \colon \WO = \bigcup_{\omega \leq \beta < \omega_1} \WO_\beta \onto \omega_1 \) sending \( x \in \WO \) to the unique \( \alpha < \omega_1 \) such that \( x \in \WO_{\omega+\alpha} \) is easily seen to be a \( \bSigma^1_2 \)-norm on the set \( \WO \in \bPi^1_1 \subseteq \bSigma^1_2 \subseteq \bS(\omega_1) \).
Notice that the definition of \( \bS ( \kappa ) \)-in-the-codes functions does not depend on the particular choice of the \( \bS ( \kappa ) \)-norm \( \rho \) by Remark~\ref{rmk:rhononessential}\ref{rmk:rhononessential-b}. 
In fact, in models of \( \AC \) the pointclass \( \bS(\kappa) \) is closed under projections and finite unions by Lemma~\ref{lem:Souslin}, and it is closed under well-ordered unions of length \( \kappa \) by Remark~\ref{rmk:souslin}: thus it always satisfies the hypothesis of Lemma~\ref{lem:inthecodes} as soon as there is an \( \bS ( \kappa ) \)-code for \( \kappa \) (i.e.\ as soon as it makes sense to speak of \( \bS(\kappa) \)-in-the-codes functions). Using these facts, we can reformulate Lemma~\ref{lem:inthecodes} 
as follows. 

\begin{proposition}[\( \AC \)] \label{prop:inthecodesSouslinAC}
Let \( \kappa \leq 2^{\aleph_0} \) be such that there is an \( \bS ( \kappa ) \)-code for \( \kappa \). 
For every \( f \colon \pre{\omega}{2} \to \pre{\kappa}{2} \) the following are equivalent:
\begin{enumerate-(a)}
\item \label{prop:inthecodesSouslinAC-1}
\( f \) is \( \bS ( \kappa ) \)-in-the-codes; 
\item \label{prop:inthecodesSouslinAC-2}
\( f^{-1} ( \widetilde{\Nbhd}_{\alpha , i}) \in \bDelta_{\bS ( \kappa )} \) for every \( \alpha < \kappa \) and \( i = 0 , 1 \);
\item \label{prop:inthecodesSouslinAC-3}
\( f^{-1} ( U ) \in \bDelta_{\bS ( \kappa ) } \) for every \( U \in \mathcal{B}_p ( \pre{\kappa}{2} ) \).
\end{enumerate-(a)}
\end{proposition}

This in particular shows that in models of \( \ZFC \) the notion of an \( \bS(\kappa) \)-in-the-codes function is nontrivial as soon as \( \bDelta_{\bS(\kappa)} \neq \pow(\pre{\omega}{2}) \): indeed, by Proposition~\ref{prop:inthecodesSouslinAC} if \( A \subseteq \pre{\omega}{2} \) does not belong to \( \bDelta_{\bS(\kappa)} \) and \( y_0, y_1 \) are distinct points of \( \pre{\kappa}{2} \), the function \( f \colon \pre{\omega}{2} \to \pre{\kappa}{2} \) sending all points of \( A \) to \( y_0 \) and all points of \( \pre{\omega}{2} \setminus A \) to \( y_1 \) is not \( \bS(\kappa) \)-in-the-codes. 
By Proposition~\ref{prop:Bernstein}, this applies to all \( \kappa < 2^{\aleph_0} \). 
Moreover, using Remark~\ref{rmk:souslin} for small cardinals \( \kappa \) we further have:

\begin{corollary}[\( \AC \)] \label{cor:inthecodesSouslinAC}
Let \( 0 \neq n \in \omega \) be such that there is an \( \bS ( \omega_n ) \)-code for \( \omega_n \).
For every \( f \colon \pre{\omega}{2} \to \pre{\omega_n}{2} \) the following are equivalent:
\begin{enumerate-(a)}
\item \label{cor:inthecodesSouslinAC-1}
\( f \) is \( \bS ( \omega_n ) \)-in-the-codes;
\item \label{cor:inthecodesSouslinAC-2}
the set \( f^{ - 1 } ( \widetilde{\Nbhd}_{ \alpha , i } ) \) is the union of \( \aleph_n \)-many Borel sets, for every \( \alpha < \kappa \) and \( i = 0,1 \). 
\end{enumerate-(a)}
\end{corollary}

Finally, using Proposition~\ref{prop:Somega1=sigma12} and Corollary~\ref{cor:inthecodesSouslinAC} we also get the following corollary characterizing \( \bSigma^1_2 \)-in-the-codes functions in certain models of \( \ZFC \).

\begin{corollary}[\( \AC \)] \label{cor:inthecodesSouslinAC2}
Assume \( \MA + \neg \CH + \EXISTS{a \in \pre{\omega}{\omega}} ( \omega_1^{ \Ll [ a ] } = \omega_1 ) \). 
Then a function \( f \colon \pre{\omega}{2} \to \pre{\omega_1}{2} \) is \( \bSigma^1_2 \)-in-the-codes if and only if for every \( \alpha < \omega_1 \) and \( i = 0 , 1 \) the set \( f^{ - 1 } ( \widetilde{\Nbhd}_{\alpha , i}) \) is the union of \( \aleph_1 \)-many Borel sets.
\end{corollary}

Note that in Proposition~\ref{prop:inthecodesSouslinAC} and in Corollaries~\ref{cor:inthecodesSouslinAC} and~\ref{cor:inthecodesSouslinAC2} above we can replace \( \pre{\kappa}{2} \) with any space of type \( \kappa \) --- see the discussion after Definition~\ref{def:spaceoftypekappa}.

\subsection{Souslin sets and cardinals in models of determinacy} \label{sec:SouslinunderAD}
The structure of the pointclasses \( \bS ( \kappa ) \) had been extensively studied under \( \AD + \DC \).
Here we summarize the most important results, referring the reader to~\cite{Jackson:2008pi, Jackson:2010ff} for proofs and further results. 

\begin{lemma}[\( \AD \)]\label{lem:kappa<Theta=>S(kappa)noteverything}
If \( \kappa < \Theta \) then \( \bS ( \kappa ) \neq \pow ( \pre{\omega}{2} ) \).
\end{lemma}

\begin{proof}
By Remark~\ref{rmk:cardinalityS(kappa)}, \( \Tr ( 2 \times \kappa ) \into \pow ( \kappa ) \), and therefore \( \pow ( \kappa ) \onto \bS ( \kappa ) \).
By Theorem~\ref{th:coding} \( \pre{\omega}{2} \onto \pow ( \kappa ) \), so the result follows from Cantor's theorem that \( \pre{\omega}{2} \) does not surject onto \( \pow ( \pre{\omega}{2} ) \).
\end{proof}

Recall that \( \Xi \leq \Theta \) by Lemma~\ref{lem:Souslincardinals<Theta}.
By Corollary~\ref{cor:SouslincardinalsunderAC}, under \( \AC \) we have that \( \Xi = 2^{\aleph_0} < ( 2^{\aleph_0} )^+ = \Theta \), so the inequality between \( \Xi \) and \( \Theta \) is strict, but in the context of \( \AD \) the situation is different as both the cases \( \Xi = \Theta \) and \( \Xi < \Theta \) can occur. 

\begin{proposition}[\( \AD + \DC \)] \label{prop:ADR<=>Xi=Theta}
The following are equivalent:
\begin{enumerate-(a)}
\item\label{prop:ADR<=>Xi=Theta-a}
 \( \Xi = \Theta \) (i.e.~Souslin cardinals are unbounded below \( \Theta \));
\item\label{prop:ADR<=>Xi=Theta-b}
every set of reals is \( \kappa \)-Souslin for some cardinal \( \kappa \) (i.e.~\( \bS(\infty) = \bS(\Xi) = \pow(\pre{\omega}{2}) \));
\item \label{prop:ADR<=>Xi=Theta-c}
\( \Unif \) (see Section~\ref{subsec:choicedeterminacy}); 
\item \label{prop:ADR<=>Xi=Theta-d}
\( \ADR \).
\end{enumerate-(a)}
\end{proposition}

\begin{proof}
For the equivalence of~\ref{prop:ADR<=>Xi=Theta-b}, ~\ref{prop:ADR<=>Xi=Theta-c} and~\ref{prop:ADR<=>Xi=Theta-d} see~\cite[Corollary 5.12]{Ketchersid:2011gb}.
The implication \( \text{\ref{prop:ADR<=>Xi=Theta-a}}\implies \text{\ref{prop:ADR<=>Xi=Theta-b}} \) is~\cite[Remark 9.21]{Woodin:2010tn}.
For \( \neg\text{\ref{prop:ADR<=>Xi=Theta-a}}\implies \neg\text{\ref{prop:ADR<=>Xi=Theta-b}} \) apply Lemma~\ref{lem:kappa<Theta=>S(kappa)noteverything} with \( \kappa = \Xi \).
\end{proof}

On the other hand, if \( \Xi < \Theta \) then by Lemma~\ref{lem:kappa<Theta=>S(kappa)noteverything} the hereditary boldface pointclass \( \bS ( \infty ) = \bS(\Xi) \) is proper, but it can happen that \( \Xi \) is a Souslin cardinal or not.
For example, \( \AD + {\Vv = \Ll ( \R )} \) implies that \( \Xi = \bdelta^2_1 \) and that \( \bdelta^2_1 \) is a Souslin cardinal --- in fact it is the \( \bdelta^2_1 \)-th Souslin cardinal by~\eqref{eq:otofXi} below.
More generally, by results of Steel and Woodin \( \AD + \DCR \) implies that the Souslin cardinals are closed below \( \Xi \) (see~\cite{Ketchersid:2011gb}).
Woodin has isolated the following natural strengthening of \( \AD \).
\begin{definition}\label{def:AD^+}
\( \AD^+ \) is the conjunction of the following statements:\index[symbols]{AD@\( \AD^+ \)}\index[concepts]{determinacy axioms!\( \AD^+ \)}
\begin{itemize}[leftmargin=3pc]
\item
\( \DCR \),
\item
\( \bB_\infty ( \pre{\omega}{2} ) = \pow ( \pre{\omega}{2} ) \),
\item
\emph{Ordinal Determinacy}: if \( \lambda < \Theta \) and \( \pi \colon \pre{ \omega }{ \lambda } \to \pre{\omega}{\omega} \) is a continuous surjection, then \( \pi^{-1} ( A ) \) is determined, for all \( A \subseteq \pre{\omega}{\omega} \).
\end{itemize}
\end{definition} 
The axiom \( \AD^+ \) is equivalent (under \( \AD + \DCR \)) to the fact that Souslin cardinals are closed below \( \Theta \) (see~\cite[Theorem 7.2]{Ketchersid:2011gb} for a proof).
Thus assuming \( \AD^+ \), 
\[
 \Xi < \Theta \IFF \Xi \text{ is a Souslin cardinal.}
 \]
Both \( \ADR + \DC \) and \( \AD + {\Vv = \Ll ( \R )} \) imply \( \AD^+ \). 
Every model of \( \AD \) known to date does satisfy \( \AD^+ \), but it is open whether \( \AD \implies \AD^+ \).
It is known that the theory \( \AD + \neg \AD^+ \) has very high consistency strength, if consistent at all (see~\cite[Section 8]{Koellner:2010fk}).

If \( \Xi \) is a Souslin cardinal, then Souslin cardinals are unbounded below \( \Xi \) (a fact which is trivially true if \( \Xi \) is not a Souslin cardinal), \( \Xi \) is regular, and \( \bS ( \Xi ) \) is closed also under coprojections (equivalently, the dual pointclass \( \check{\bS} ( \Xi ) \) of \( \bS ( \Xi ) \) is closed under projections). 
In particular, 
\begin{equation}\label{eq:otofXi}
 \Xi \text{ is the \( \Xi \)-th Souslin cardinal.}
\end{equation}
Since Souslin cardinals are closed below \( \Xi \) they are also closed below \emph{every} Souslin cardinal \( \kappa \), and hence by a result of Kechris \( \bS ( \kappa ) \) is always a nonselfdual pointclass for such \( \kappa \)'s; thus \( \bS ( \kappa ) \) can be selfdual only if \( \kappa = \Xi \) and \( \Xi \) is not a Souslin cardinal.

A particular (and important) kind of Souslin cardinals are the ones related to the projective hierarchy:%
\footnote{As explained in~\cite[Theorem 3.28]{Jackson:2010ff}, the ``pattern'' of Souslin cardinals we are going to describe can be somehow reproduced above any \( \kappa < \Xi \) which is a limit of Souslin cardinals. 
We repeatedly use this fact in the proof of Proposition~\ref{prop:deltaSkappa}.} 
assuming \( \AD + \DC \) the Souslin cardinals below \( \bdelta_\omega \coloneqq \sup_{n \in \omega} \bdelta^1_n = \aleph_{\varepsilon_0} \) are exactly the \( \blambda^1_{2 n + 1} \)'s and the \( \bdelta^1_{2 n + 1} \)'s, and we have that \( \bS ( \blambda^1_{2 n + 1}) = \bSigma^1_{2 n + 1} \) and \( \bS ( \bdelta^1_{2 n + 1}) = \bSigma^1_{2 n + 2} \)~\cite[Theorem 2.18]{Jackson:2010ff}, and if \( \kappa \) is one of these Souslin cardinals, then \( \bdelta_{\bS ( \kappa )} = \kappa^+ \) (see Section~\ref{subsec:projectiveordinals}). 
This is generalized by the following Proposition~\ref{prop:deltaSkappa} (compare it with Proposition~\ref{prop:deltaSkappaAC}, where an analogous result is proved under \( \AC \)). 
Recall from~\cite{Jackson:2010ff} that if \( \kappa \) is a limit of Souslin cardinals and has uncountable cofinality, then we can associate to it a canonical Steel pointclass \( \bGamma_0 \) (whose existence is granted by the analysis in~\cite{Steel:2012st}), namely: \( \bGamma_0 \) is a nonselfdual boldface pointclass closed under coprojections and such that \( \bDelta_{\bGamma_0} = \bigcup_{\lambda < \kappa} \bS(\lambda) \). 
By~\cite[Lemma 3.8]{Jackson:2010ff}, \( \bS(\kappa) \) is the closure under projections of \( \bGamma_0 \), and both \( \bS(\kappa) \) and \( \bGamma_0 \) have the scale property (and hence also the prewellordering property). 
The proof of the mentioned~\cite[Lemma 3.8]{Jackson:2010ff} further shows that \( \kappa = \bdelta_{\bGamma_0} \).

\begin{proposition}[\( \AD + \DC \)] \label{prop:deltaSkappa}
Let \( \kappa \) be a Souslin cardinal. 
\begin{enumerate-(a)}
\item \label{prop:deltaSkappa-a}
Either \( \bdelta_{\bS ( \kappa )} = \kappa \) or \( \bdelta_{\bS ( \kappa )} = \kappa^+ \). 
\item \label{prop:deltaSkappa-b}
There is an \( \bS ( \kappa ) \)-code for \( \kappa \). 
\item \label{prop:deltaSkappa-c}
The following are equivalent:
\begin{enumerate-(1)}
\item \label{prop:deltaSkappa-c1}
\( \bdelta_{\bS(\kappa)} = \kappa \);
\item \label{prop:deltaSkappa-c2}
\( \bS(\kappa) \) is closed under coprojections;
\item \label{prop:deltaSkappa-c3}
either \( \kappa = \Xi \), or else \( \kappa < \Xi \) is a limit of uncountable cofinality of Souslin cardinals falling in Case III of~\cite[Theorem 3.28]{Jackson:2010ff}, i.e.\ such that its associated Steel pointclass \( \bGamma_0 \) is closed under projections;
\item \label{prop:deltaSkappa-c4}
\( \kappa \) is a regular limit of Souslin cardinals and its associated Steel pointclass \( \bGamma_0 \) is closed under projections;
\item \label{prop:deltaSkappa-c5}
\( \kappa \) is a regular limit of Souslin cardinals and \( \bS(\kappa) \) coincides with its associated Steel pointclass \( \bGamma_0 \).
\end{enumerate-(1)}
In particular, \( \bdelta_{\bS ( \kappa ) } = \kappa^+ \) whenever \( \kappa \) is not a regular limit of Souslin cardinals.
\end{enumerate-(a)} 
\end{proposition}

\begin{proof}
\ref{prop:deltaSkappa-a}
By~\cite[Theorem 2.18]{Jackson:2010ff} (see the paragraph preceding this proposition), we can assume without loss of generality that \( \kappa \) is above a limit of Souslin cardinals, so that we can use the deep and remarkable classification of the Souslin cardinals presented e.g.\ in~\cite[Theorem 3.28]{Jackson:2010ff}. 
The inequality \( \bdelta_{\bS ( \kappa )} \leq \kappa^+ \) immediately follows from the Kunen-Martin's theorem (see e.g.\ \cite[Theorem 2G.2]{Moschovakis:2009fk}), so that we just need to show \( \kappa \leq \bdelta_{\bS ( \kappa )} \). 
From~\cite[Theorem 7D.8]{Moschovakis:2009fk} we have that \( \bdelta_{\bS ( \kappa )} \) is a cardinal, and that it is of uncountable cofinality by closure of \( \bS ( \kappa ) \) under countable unions (Lemma~\ref{lem:Souslin}).

Assume first that \( \kappa \) is not a limit of Souslin cardinals itself (so that, in particular, there is a largest \( \kappa' < \kappa \) which is a limit of Souslin cardinals, which can be taken as basis to carry out the analysis provided in~\cite[Theorem 3.28]{Jackson:2010ff}). 
We consider two distinct cases, depending on whether \( \kappa \) has countable cofinality or not. 
In the former case, it follows from~\cite[Theorem 3.28]{Jackson:2010ff} that \( \kappa = \lambda_{2i+1} \) for some \( i \in \omega \) and that \( \check{\bS} ( \kappa ) = \bPi_{2i+1} \) has the scale property. 
Since the (regular) \( \check{\bS} ( \kappa ) \)-norms constituting a \( \check{\bS} ( \kappa ) \)-scale on an arbitrary set \( A \subseteq \pre{\omega}{\omega} \) in \( \check{\bS} ( \kappa ) \) are into \( \bdelta_{{\bS} ( \kappa )} \), if \( \bdelta_{\bS ( \kappa )} \leq \kappa \) then the mentioned scale would canonically yield a \( \kappa \)-Souslin representation of \( A \), and therefore we would get \( \check{\bS} ( \kappa ) \subseteq \bS ( \kappa ) \): but this contradicts the fact that \( \bS ( \kappa ) \) is nonselfdual. 
Thus we get \( \kappa < \bdelta_{\bS ( \kappa )} \), and thus also \( \bdelta_{\bS ( \kappa )} = \kappa^+ \) because \( \bdelta_{\bS ( \kappa )} \) is a cardinal. 
If instead \( \kappa \) has uncountable cofinality, then \( \bS ( \kappa ) \) has the scale property by~\cite[Lemma 3.7]{Jackson:2010ff}. 
Let \( A \subseteq \pre{\omega}{\omega} \) be such that \( A \in \bS ( \kappa ) \setminus \check{\bS} ( \kappa ) \), and consider an arbitrary \( \bS ( \kappa ) \)-scale \( \seqofLR{ \rho_n }{ n \in \omega } \) on it.
Arguing as in the proof of Lemma~\ref{lem:Souslincardinals<Theta}, by our choice of \( A \) we get
\[ 
\sup \setofLR{\alpha \in \On}{\EXISTS{n \in \omega} \EXISTS{x \in A} ( \rho_n ( x ) = \alpha ) } = \kappa.
 \] 
Moreover, since \( \kappa \) has uncountable cofinality there is \( \bar{n} \in \omega \) such that the \( \bS ( \kappa ) \)-norm \( \rho_{\bar{n}} \) on \( A \) is onto \( \kappa \). 
Then for every \( \alpha < \kappa \) we can construct a \( \bDelta_{\bS ( \kappa )} \) prewellordering \( \preceq \) of \( \pre{\omega}{\omega} \) of length \( \alpha + 2 \) by fixing \( z \in A \) such that \( \rho_{\bar{n}} ( z) = \alpha \) and then setting
\[ 
x \preceq y \iff y \notin A' \OR ( x \in A' \wedge y \in A' \wedge \rho_{\bar{n}} ( x ) \leq \rho_{\bar{n}} ( y ) ),
 \] 
where \( A' = \setofLR{x \in \pre{\omega}{\omega}}{x \in A \wedge \rho_{\bar{n}}(x) \leq \rho_{\bar{n}}(z)} \in \bDelta_{\bS(\kappa)} \).
Therefore \( \kappa \leq \bdelta_{\bS ( \kappa )} \) again.

Assume now that \( \kappa \) is a limit of Souslin cardinals. 
Then for every \( \alpha < \kappa \) there is \( \alpha < \lambda < \kappa \) such that \( \lambda \) is a Souslin cardinal and is not a limit of Souslin cardinals, so that \( \bdelta_{\bS ( \lambda)} \geq \lambda \) by the above proof. 
Since \( \lambda \leq \kappa \) implies \( \bdelta_{\bS ( \kappa )} \geq \bdelta_{\bS ( \lambda)} \), we get \( \bdelta_{\bS ( \kappa )} > \alpha \), whence \( \kappa \leq \bdelta_{\bS ( \kappa )} \).

\smallskip

\ref{prop:deltaSkappa-b}
If \( \kappa < \bdelta_{\bS ( \kappa )} \) then we are done by Remark~\ref{rmk:inthecodes}\ref{rmk:inthecodes-ii}, so let us assume \( \kappa = \bdelta_{\bS ( \kappa )} \). 
Notice that this implies that \( \kappa \) is of uncountable cofinality because, as recalled at the beginning of the proof of part~\ref{prop:deltaSkappa-a}, \( \cof ( \bdelta_{\bS ( \kappa )}) >\omega \). 
By%
\footnote{In fact we show in part~\ref{prop:deltaSkappa-c} that if \( \kappa = \bdelta_{\bS ( \kappa ) } \) then \( \kappa \) is a regular limit of Souslin cardinals, so only~\cite[Lemma 3.8]{Jackson:2010ff} needs to be applied in this case. 
However, since the proof of~\ref{prop:deltaSkappa-c} partially relies on~\ref{prop:deltaSkappa-b}, here we mentioned also~\cite[Lemma 3.7]{Jackson:2010ff} to make it evident that there is no circularity in the argument.}
~\cite[Lemmas 3.7 and 3.8]{Jackson:2010ff}, we then get that \( \bS ( \kappa ) \) has the scale property. 
Let \( A \subseteq \pre{\omega}{\omega} \) be in \( \bS ( \kappa ) \setminus \bigcup_{\lambda < \kappa} \bS ( \lambda) \), and let \( \seqofLR{ \rho_n }{ n \in \omega } \) be an \( \bS ( \kappa ) \)-scale on \( A \). 
Arguing as in the proof of Lemma~\ref{lem:Souslincardinals<Theta} and using again the fact that \( \kappa \) has uncountable cofinality, by the choice of \( A \) there is at least one \( \bar{n} \in \omega \) such that \( \rho_{\bar{n}} \) has length \( \kappa \): setting \( \rho \coloneqq \rho_{\bar{n}} \) we get the desired result.

\smallskip

\ref{prop:deltaSkappa-c} 
First we prove the equivalence between the conditions~\ref{prop:deltaSkappa-c2}--\ref{prop:deltaSkappa-c5}. 
Assume~\ref{prop:deltaSkappa-c2}, so that in particular \( \kappa > \omega \). 
Then \( \kappa \) is a regular limit of Souslin cardinals by~\cite[Lemma 3.6]{Jackson:2010ff}. 
As recalled before this proposition, in this situation \( \bS(\kappa) \) is the closure under projection of its associated Steel pointclass \( \bGamma_0 \), so that in particular \( \bS(\kappa) \supseteq \bGamma_0 \). 
If this inclusion were proper, then \( \bGamma_0 \) could not be closed under projections, and thus \( \kappa \) would fall in Case II of~\cite[Theorem 3.28]{Jackson:2010ff}: but this would contradict the closure under coprojections of \( \bS(\kappa) \). 
It follows that \( \bS(\kappa) = \bGamma_0 \), i.e.~\ref{prop:deltaSkappa-c5}. 
The implications from~\ref{prop:deltaSkappa-c5} to~\ref{prop:deltaSkappa-c4} and from~\ref{prop:deltaSkappa-c4} to~\ref{prop:deltaSkappa-c3} are obvious (recall that under our assumptions \( \bS(\kappa) \) is always closed under projections by Lemma~\ref{lem:Souslin}), so let us show that ~\ref{prop:deltaSkappa-c3} implies~\ref{prop:deltaSkappa-c2}. 
If \( \kappa = \Xi \) then the result follows from~\cite[Lemma 2.20]{Jackson:2008pi}. 
In the remaining case, \( \bS(\kappa) = \bGamma_0 \) because \( \bS(\kappa) \) is the closure under projections of \( \bGamma_0 \), and thus~\ref{prop:deltaSkappa-c2} follows from the fact that \( \bGamma_0 \) is closed under coprojections.

Since~\ref{prop:deltaSkappa-c5} easily implies~\ref{prop:deltaSkappa-c1} (because \( \kappa = \bdelta_{\bGamma_0} \)), to conclude our proof it is enough to show that~\ref{prop:deltaSkappa-c1} implies~\ref{prop:deltaSkappa-c4}. 

\begin{claim} \label{claim:regularlimitnec}
If \( \kappa = \bdelta_{\bS(\kappa)} \), then \( \kappa \) is a regular limit of Souslin cardinals.
\end{claim}

\begin{proof}[Proof of the claim]
We prove the contrapositive, i.e.\ that if \( \kappa \) is \emph{not} a regular limit of Souslin cardinals then \( \bdelta_{\bS(\kappa)} \neq \kappa \) 
(whence \( \bdelta_{\bS ( \kappa )} = \kappa^+ \) by~\ref{prop:deltaSkappa-a}). 
This is trivial for Souslin cardinals \( \kappa \) of countable cofinality: as observed in the proof of part~\ref{prop:deltaSkappa-b}, in such case \( \bdelta_{\bS ( \kappa )} \neq \kappa \) because \( \bdelta_{\bS(\kappa)} \) has uncountable cofinality by closure of \( \bS ( \kappa ) \) under countable unions. 
Now assume that \( \kappa \) has uncountable cofinality and is not a limit of Souslin cardinals, and let \( \lambda \) be the largest Souslin cardinal smaller than \( \kappa \). 
By (the proof of)~\cite[Lemma 3.7]{Jackson:2010ff}, \( \lambda \) is of countable cofinality, \( \kappa = \lambda^+ \), and \( \check{\bS} ( \lambda) \) has the prewellordering property. 
Since \( \check{\bS} ( \lambda) \), being nonselfdual, admits a universal set and is closed under coprojections, we get from e.g.\ \cite[Lemma 1.3]{Jackson:2008pi} that for every \( A \subseteq \pre{\omega}{\omega} \) such that \( A \in \check{\bS} ( \lambda) \setminus \bS ( \lambda) \) there is an \( \check{\bS} ( \lambda) \)-norm \( \rho \) on \( A \) of length \( \bdelta_{\bS ( \lambda)} \). 
Since \( \kappa \) is a Souslin cardinal, \( \bS ( \kappa ) \) is nonselfdual and \( \bS ( \lambda) \subset \bS ( \kappa ) \), whence \( \bS ( \lambda) \cup \check{\bS} ( \lambda) \subseteq \bDelta_{\bS ( \kappa )} \). 
Therefore, the prewellordering \( \preceq \) of \( \pre{\omega}{\omega} \) defined by
\begin{equation} \label{eq:pwopreceq}
x \preceq y \iff y \notin A \vee ( x \in A \wedge y \in A \wedge \rho ( x ) \leq \rho ( y ) )
 \end{equation}
is in \( \bDelta_{\bS ( \kappa )} \) and has length \( \bdelta_{\bS ( \lambda)} + 1 \). 
Since \( \bdelta_{\bS ( \lambda)} = \lambda^+ \) (because \( \lambda \), being a Souslin cardinal of countable cofinality, falls in the already considered trivial case), this implies \( \bdelta_{\bS ( \kappa )} > \lambda^+ = \kappa \).

Finally, let us assume that \( \kappa \) is a limit of Souslin cardinals such that \( \omega < \cof ( \kappa ) < \kappa \), and fix a cofinal map \( f \colon \cof ( \kappa ) \to \kappa \).
By part~\ref{prop:deltaSkappa-b}, for some \( \kappa \)-Souslin \( A \subseteq \pre{\omega}{\omega} \) there is a \( \bS ( \kappa ) \)-norm \( \rho \) on \( A \) of length \( \kappa \). 
Fix a Souslin cardinal \( \cof ( \kappa ) < \lambda < \kappa \) and some \( \bDelta_{ \bS ( \lambda )} \)-norm \( \sigma \colon \pre{\omega}{\omega} \onto \cof ( \kappa ) \) (which exists because by part~\ref{prop:deltaSkappa-a} we have \( \cof ( \kappa ) < \lambda \leq \bdelta_{\bS ( \lambda ) } \)), so that the strict well-founded relation on \( \pre{\omega}{\omega} \) associated to \( \sigma \) is in \( \bS ( \lambda ) \). 
Then apply the first Coding Lemma~\cite[Theorem 2.12]{Jackson:2010ff} to 
\[ 
R \coloneqq \setofLR{ ( z , w ) \in \pre{\omega}{\omega} \times \pre{\omega}{\omega}}{w \in A \wedge \rho ( w ) = f ( \sigma ( z ) ) } 
\]
 to obtain an \( \bS ( \lambda ) \)-set \( C \subseteq \pre{\omega}{\omega} \times \pre{\omega}{\omega} \) such that
\begin{itemize}[leftmargin=1pc]
\item
for every \( \alpha < \cof ( \kappa ) \) there is \( ( z , w ) \in C \) with \( \sigma ( z ) = \alpha \), and
\item
for every \( ( z , w ) \in C \), \( w \in A \) and \( \rho ( w ) = f ( \sigma ( z ) ) \).
\end{itemize} 
Then the prewellordering \( \preceq \) of \( \pre{\omega}{\omega} \times \pre{\omega}{\omega} \) defined by setting \( ( x , y ) \preceq ( x' , y' ) \) if and only if
\begin{align*} 
 \sigma ( x ) < \sigma ( x' ) \vee \big [ \sigma ( x ) = \sigma ( x' ) \wedge \big ( & \EXISTS { ( z , w ) \in C} ( \sigma ( z ) = \sigma ( x ) \wedge \rho ( y' ) > \rho ( w ) ) \vee \\ 
& \EXISTS { ( z , w ) \in C} ( \sigma ( z ) = \sigma ( x ) \wedge \rho ( y ) \leq \rho ( y' ) \leq \rho ( w ) ) \big ) \big] 
 \end{align*}
has length \( \sum_{\alpha < \cof ( \kappa ) } ( f ( \alpha ) + 2 ) = \kappa \). 
Moreover, by our choice of \( C \) we also have that \( ( x , y ) \preceq ( x' , y' ) \) if and only if
\begin{align*} 
 \sigma ( x ) < \sigma ( x' ) \vee \big [ \sigma ( x ) = \sigma ( x ' ) \wedge \big ( & \FORALL { ( z , w ) \in C} ( \sigma ( z ) = \sigma ( x ) \IMPLIES \rho ( y' ) > \rho ( w ) ) \vee 
 \\ 
& \FORALL { ( z , w ) \in C} ( \sigma ( z ) = \sigma ( x ) \IMPLIES \rho ( y ) \leq \rho ( y' ) \leq \rho ( w ) ) \big )\big ] , 
\end{align*}
and hence \( \preceq \) is in \( \bDelta_{\bS ( \kappa )} \) (here we are using the fact that \( \bS ( \lambda ) \cup \check{\bS}(\lambda) \subseteq \bDelta_{\bS ( \kappa )} \)). 
Since \( \pre{\omega}{\omega} \) and \( \pre{\omega}{\omega} \times \pre{\omega}{\omega} \) are homeomorphic, this shows that \( \bdelta_{\bS ( \kappa )} > \kappa \), so we are done.
\end{proof}
Now assume~\ref{prop:deltaSkappa-c1}, so that \( \kappa \) is a regular limit of Souslin cardinals by the above Claim~\ref{claim:regularlimitnec}. 
As observed in~\cite{Jackson:2008pi}, from this and~\cite[Theorem 2.1]{Steel:2012st} it follows that \( \bGamma_0 \) is closed under finite unions. 
Since \( \bGamma_0 \) is also closed under coprojections and admits a universal set (being nonselfdual), by the prewellordering property for \( \bGamma_0 \) and e.g.~\cite[Lemma 1.3]{Jackson:2008pi} we get that for any \( A \in \bGamma_0 \setminus \check{\bGamma}_0 \) there is a \( \bGamma_0 \)-norm \( \rho \) on \( A \) of length \( \bdelta_{\bGamma_0} = \kappa \). 
Assume towards a contradiction that \( \bGamma_0 \) is not closed under projections. 
Then \( \bS(\kappa) \), which is the closure under projections of \( \bGamma_0 \), would contain both \( \bGamma_0 \) and its dual \( \check{\bGamma}_0 \), so that \( \bGamma_0 \cup \check{\bGamma}_0 \subseteq \bDelta_{\bS(\kappa)} \). 
But then then prewellordering \( \preceq \) obtained from the \( \bGamma_0 \)-norm \( \rho \) as in~\eqref{eq:pwopreceq} would be in \( \bDelta_{\bS(\kappa)} \): since its length is \( \kappa+1 \), this would contradict our assumption \( \kappa = \bdelta_{\bS(\kappa)} \). 
Thus \( \bGamma_0 \) is closed under projections and~\ref{prop:deltaSkappa-c4} holds.
\end{proof}

\begin{remark} \label{rmk:deltaSkappa}
We currently do not know whether there are regular limits of Souslin cardinals which may actually fall into Case II of~\cite[Theorem 3.28]{Jackson:2010ff}, i.e.\ whether all regular limits of Souslin cardinals need to satisfy the equivalent conditions of Proposition~\ref{prop:deltaSkappa}\ref{prop:deltaSkappa-c}. 
\end{remark}

As we already observed, Souslin cardinals are always closed and unbounded below \( \Xi \): for our purposes, it is important to notice that the analysis of~\cite{Jackson:2010ff} shows that also the \emph{regular} Souslin cardinals are unbounded below \( \Xi \).

\begin{lemma}[\( \AD+\DC \)] \label{lem:regularSouslinareunbounded}
Regular Souslin cardinals are unbounded below \( \Xi \). 
In particular, for every Souslin cardinal \( \kappa \) there is a \emph{regular} Souslin cardinal \( \kappa' \geq \kappa \).
\end{lemma}

\begin{proof}
Given \( \alpha < \Xi \), let \( \bar{\kappa} \) be the smallest Souslin cardinal above \( \alpha \) and \( \kappa \) be the largest limit of Souslin cardinals \( \leq \bar{\kappa} \) (such \( \bar{\kappa} \) and \( \kappa \) exists because, as already recalled in the discussion preceding this lemma, there are club-many Souslin cardinals below \( \Xi \)). 
Applying~\cite[Theorem 3.28]{Jackson:2010ff} to such \( \kappa \), we obtain that there is \( i \in \omega \) such that \( \bar{\kappa} \leq \delta_{2i+1} \), where \( \delta_{2i+1} = (\lambda_{2i+1})^+ \) is as in any of Case I--III of~\cite[Theorem 3.28]{Jackson:2010ff}. 
Since in all these cases
\[ 
\delta_{2i+1} = \bdelta_{\bPi_{2i+1}} = \bdelta_{\bSigma_{2i+1}} = \bdelta_{\bS(\lambda_{2i+1})}, 
\]
by the Kunen-Martin's theorem (see e.g.\ \cite[Theorem 2G.2]{Moschovakis:2009fk}) the Souslin cardinal \( \delta_{2i+1} \) is the supremum of the lengths of the \( \lambda_{2i+1} \)-Souslin strict well-founded relations on \( \pre{\omega}{\omega} \), and thus it is a regular cardinal by~\cite[Lemma 2.16]{Jackson:2010ff}. 
Therefore since \( \alpha < \bar{\kappa} \leq \delta_{2i+1} < \Xi \) we are done. 
The second part of the lemma follows from the first one and the fact that if \( \Xi \) is a Souslin cardinal then it is also regular.
\end{proof}

We now consider \( \bS(\kappa) \)-in-the-codes functions in models of determinacy.
By Proposition~\ref{prop:deltaSkappa}\ref{prop:deltaSkappa-b}, if \( \kappa \) is a Souslin cardinal then it always makes sense to speak about \( \bS ( \kappa ) \)-in-the-codes functions \( f \colon \pre{\omega}{2} \to \pre{\kappa}{2} \). 
Notice also that the definition of \( \bS ( \kappa ) \)-in-the-codes functions does not depend on the particular choice of the \( \bS ( \kappa ) \)-norm \( \rho \) by Remark~\ref{rmk:rhononessential}\ref{rmk:rhononessential-b} and Lemma~\ref{lem:Souslin}. 
Using these facts, we can reformulate Proposition~\ref{prop:inthecodes} as follows. 
(When \( \kappa = \omega \) we can dispense with all determinacy assumptions in Proposition~\ref{prop:inthecodesSouslin} and Corollary~\ref{cor:inthecodesSouslin}.)

\begin{proposition}[\( \AD + \DC \)] \label{prop:inthecodesSouslin}
Let \( \kappa \) be a Souslin cardinal. 
For every \( f \colon \pre{\omega}{2} \to \pre{\kappa}{2} \) the following are equivalent:
\begin{enumerate-(a)}
\item \label{prop:inthecodesSouslin-1}
\( f \) is \( \bS ( \kappa ) \)-in-the-codes: 
\item \label{prop:inthecodesSouslin-2}
\( f^{-1} ( \widetilde{\Nbhd}_{\alpha , i}) \in \bDelta_{\bS ( \kappa )} \) for every \( \alpha < \kappa \) and \( i = 0,1 \);
\item \label{prop:inthecodesSouslin-3}
\( f^{-1} ( U ) \in \bDelta_{\bS ( \kappa ) } \) for every \( U \in \mathcal{B}_p ( \pre{\kappa}{2} ) \).
\end{enumerate-(a)}
\end{proposition}

\begin{proof}
The pointclass \( \bS ( \kappa ) \) satisfies the hypotheses of Proposition~\ref{prop:inthecodes} by Lemma~\ref{lem:Souslin}. 
\end{proof}

Notice that in Proposition~\ref{prop:inthecodesSouslin} (as well as in the subsequent Corollary~\ref{cor:inthecodesSouslin}) we can replace \( \pre{\kappa}{2} \) with any space of type \( \kappa \) --- see the discussion after Definition~\ref{def:spaceoftypekappa}. 
Moreover, arguing as in the paragraph after Proposition~\ref{prop:inthecodesSouslinAC} one may observe that the above proposition implies in particular that also in models of \( \AD + \DC \) the notion of an \( \bS(\kappa) \)-in-the-codes function \( f \colon \pre{\omega}{2} \to \pre{\kappa}{2} \) is nontrivial as soon as \( \bDelta_{\bS(\kappa)} \neq \pow ( \pre{\omega}{2} )\): by Proposition~\ref{prop:ADR<=>Xi=Theta}, this is the case for all Souslin cardinals \( \kappa \).

Recall that by Proposition~\ref{prop:Skappainthecodesareborel}, every \( \bS ( \kappa ) \)-in-the-codes function \( f \colon \pre{\omega}{2} \to \pre{\kappa}{2} \) is weakly \( \kappa + 1 \)-Borel. 
Since when \( \kappa = \blambda^1_{2 n + 1} \) for some \( n \in \omega \) we have \( \bS ( \kappa ) = \bSigma^1_{2 n + 1} \), the next corollary shows that in certain cases the two notions coincide.

\begin{corollary}[\( \AD + \DC \)] \label{cor:inthecodesSouslin}
Let \( \kappa = \blambda^1_{2 n + 1} \) for some \( n \in \omega \). 
Then a function \( f \colon \pre{\omega}{2} \to \pre{\kappa}{2} \) is \( \bSigma^1_{2 n + 1} \)-in-the-codes if and only if it is weakly \( \kappa + 1 \)-Borel.
\end{corollary}

\begin{proof} 
Since \( \bSigma^1_{2 n + 1} = \bS ( \blambda^1_{2n+1} ) \), the forward direction is just an instantiation of Proposition~\ref{prop:Skappainthecodesareborel}, while the other direction follows from Corollary~\ref{cor:5.6}\ref{cor:5.6-b} and Proposition~\ref{prop:inthecodesSouslin}.
\end{proof}

\section{The main construction}\label{sec:mainconstruction}
In~\cite{Louveau:2005cq}, the main technical construction for proving the completeness (for analytic quasi-orders) of the relation \( \embeds^\omega_\CT \) of embeddability between countable combinatorial trees is a map which given an arbitrary \( S \in \Tr ( 2 \times \omega ) \) provides a combinatorial tree \( \mathbb{G}_S \). 
The tree \( \mathbb{G}_S \) is constructed in two steps: first a fixed combinatorial tree \( \mathbb{G}_0 \) is defined, independent of \( S \), and then certain auxiliary combinatorial trees, called \emph{forks}, coding the tree \( S \) are added to \( \mathbb{G}_0 \).
In order to prove the invariant universality of the embeddability relation between countable structures, in~\cite{Friedman:2011cr} this construction is improved so that the resulting \( \mathbb{G}_S \) is rigid, i.e.\ without nontrivial automorphisms. 
As explained in that paper, there are at least two ways to ensure that the resulting structure is rigid:
\begin{enumerate-(1)}
\item \label{en:approach1}
enrich the language for graphs \( \LL = \{ \edge \} \) with an additional binary relational symbol \( \order \), and then expand \( \mathbb{G}_S \) to a so-called \emph{ordered} combinatorial tree \( \bar{\mathbb{G}}_S \) by interpreting \( \order \) as a well-order on the vertices of \( \mathbb{G}_S \) (see the proof of~\cite[Theorem 2.4]{Friedman:2011cr});
\item \label{en:approach2}
enlarge \( \mathbb{G}_0 \) to a rigid combinatorial tree \( \mathbb{G}_1 \), and then add the forks to \( \mathbb{G}_1 \) (see the proof of~\cite[Theorem 2.4]{Friedman:2011cr} or~\cite[Section 3]{Camerlo:2012kx}).
\end{enumerate-(1)}

Although the construction in~\ref{en:approach1} is simpler, such approach is slightly unnatural because it forces us to consider more complex structures, while one of the motivations in using combinatorial trees in~\cite{Louveau:2005cq} was that they are rather simple objects. 
The approach~\ref{en:approach2}, even if technically more involved, gives instead the stronger result that already \( \embeds^\omega_\CT \) is invariantly universal, and it has proven to be quite useful in the applications to infinite combinatorics, topology, analysis, and Banach space theory~\cite{Camerlo:2012kx}.

In this paper we generalize both constructions~\ref{en:approach1} and~\ref{en:approach2} to uncountable cardinals \( \kappa \).
As in the classical case, the approach using combinatorial trees (which is developed in this section and in the subsequent Sections~\ref{sec:embeddabilitygraphs}--\ref{sec:invariantlyuniversal}) is preferable as it deals with more elementary objects, and it yields full generalizations of Theorems~\ref{th:LouveauRosendal} and~\ref{th:mottorosfriedman}, as well as most of the applications, including results on non-separable (discrete, ultrametric) complete metric spaces and on non-separable Banach spaces. 
A generalization of Theorem~\ref{th:mottorosfriedman2} seems to require the approach via ordered combinatorial trees and is postponed to Section~\ref{sec:alternativeapproach}.

\medskip

Fix an infinite cardinal \( \kappa \).
We adapt the constructions from~\cite{Louveau:2005cq, Friedman:2011cr, Camerlo:2012kx} to define a map 
\[ 
\Tr ( 2 \times \kappa ) \to \CT_\kappa , \qquad S \mapsto \mathbb{G}_S , 
\]
where \( \CT_ \kappa \) is the set of (codes for) all combinatorial trees of size \( \kappa \) from~\eqref{pag:CT_kappa} on page~\pageref{pag:CT_kappa}.
Such map will then be used in Sections~\ref{sec:embeddabilitygraphs} and~\ref{sec:invariantlyuniversal} to prove that \( { \embeds^\kappa_\CT } \), the embeddability relation on \( \CT_\kappa \), is invariantly universal%
\footnote{The notion of invariant universality is given in Definition~\ref{def:invuniversal}.} 
--- and hence also complete --- for \( \kappa \)-Souslin quasi-orders. 
Since as explained in the previous paragraph we want to work just with combinatorial trees (without any additional order on their vertices), we will follow the approach~\ref{en:approach2} briefly described above, that is:
\begin{itemize}[leftmargin=1pc]
\item
construct a basic \( \mathbb{G}_0 \in \CT_\kappa \) which is independent of the given input \( S \in \Tr(2 \times \kappa) \);
\item
enlarge \( \mathbb{G}_0 \) to a suitable \( \mathbb{G}_1 \in \CT_\kappa \) (still independent of \( S \)) to get a sufficiently rigid structure;
\item
to get the final \( \mathbb{G}_S \in \CT_\kappa \), add to \( \mathbb{G}_1 \) some forks which code enough information on \( S \).
\end{itemize} 

\begin{remark} \label{rmk:RCT}
As already observed in~\cite{Louveau:2005cq}, it is easy to check that we could systematically replace combinatorial trees with \emph{rooted} combinatorial trees in all the constructions and results below --- just define the empty sequence \( \emptyset \) to be the root of the combinatorial tree \( \mathbb{G}_0 \) (and hence also of \( \mathbb{G}_1 \) and of every combinatorial tree of the form \( \mathbb{G}_S \) for \( S \in \Tr ( 2 \times \kappa) \)), and check that all proofs go through. 
\end{remark}

Notice that as for the basic case \( \kappa = \omega \) considered in~\cite{Louveau:2005cq, Friedman:2011cr, Camerlo:2012kx}, the preparatory enlargement from \( \mathbb{G}_0 \) to \( \mathbb{G}_1 \) is necessary only for the proof of the invariant universality of \( \embeds^\kappa_\CT \): a variant of the main construction in which we attach the forks directly to \( \mathbb{G}_0 \) would already enable us to prove the completeness of \( \embeds^\kappa_\CT \) (see Remark~\ref{rmk:partialembedding}).

Let us first fix some notation concerning combinatorial trees.
The language of graphs \( \LL = \setLR{ \edge} \)\label{pag:L}\index[symbols]{E@\( \edge \)} consists of one binary relational symbol, and each graph \( G = ( V , E ) \) (see Section~\ref{subsubsec:graph}) is identified with the \( \LL \)-structure \( X = \seqLR{ X ; \edge^{X} } \) with \( X \coloneqq V \) and \( \edge^{X} \coloneqq \setofLR{ (v_0,v_1) \in V^2 }{ \setLR{ v_0,v_1 } \in E } \) (so that \( \edge^{X} \) is an irreflexive and symmetric relation on \( X \)). 
Recall from page~\eqref{pag:CT_kappa} that \( \CT_\kappa \) is the collection of (codes for) all \emph{combinatorial trees} of the form \( ( \kappa , E ) \) with \( E \subseteq [ \kappa ]^2 \). 
In fact,
\[ 
\CT_\kappa = \Mod^\kappa_{\upsigma_\CT},
 \] 
where \( \Mod^\kappa_{\upsigma_\CT} \) is defined as in~\eqref{eq:modsigma} with \( \upsigma_{\CT} \) the \( \LL_{\omega_1 \omega} \)-sentence axiomatizing combinatorial trees:
\begin{multline}
\tag{$\upsigma_\CT$} \label{eq:sigmaCT}
\forall \V_0 \, \neg ( \V_0 \edge \V_0 ) \wedge \forall \V_0 \forall \V_1 \left ( \V_0 \edge \V_1 \implies {\V_1 \edge \V_0 }\right ) \wedge 
\\
 \bigwedge\nolimits_{n \in \omega} \neg \EXISTS{\seqofLR{ \V_i }{ i \leq n + 2 }} \left [ \left ( \bigwedge\nolimits_{ i < j \leq n + 2 } \V_i \nequals \V_j \right ) \wedge \left ( \bigwedge\nolimits_{i < n + 2} \V_i \edge \V_{i + 1} \right ) \wedge {\V_0 \edge \V_{n + 2}} \right ] \wedge 
 \\
\forall \V_0 \forall \V_1 \left [ \bigvee\nolimits_{n \in \omega} \EXISTS{\seqofLR{ \V_{i + 2} }{ i \leq n + 1 } } \left ( { \V_0 \equals \V_2} \wedge {\V_1 \equals \V_{n + 3}} \wedge {\bigwedge\nolimits_{i < n + 1} \V_{i + 2} \edge \V_{i + 3}} \right ) \right ].
\end{multline}\index[symbols]{sigmaCT@\( \upsigma_\CT \)}

In order to simplify the notation, we will further abbreviate the embeddability and isomorphism relations \( \embeds^\kappa_{\upsigma_\CT} \) and \( \cong^\kappa_{\upsigma_\CT} \) on \( \CT_\kappa \) (see page~\pageref{not:embeddability}) with \( \embeds^\kappa_{\CT} \) and \( \cong^\kappa_{\CT} \), respectively. 

\subsection{The combinatorial trees \( \mathbb{G}_0 \) and \( \mathbb{G}_1 \)} \label{subsec:G_0andG_1}

\begin{figure}
 \centering
\begin{tikzpicture}
\filldraw (-1, 1) circle (2pt) -- (0, 0) circle (2pt) -- (1 , 1) circle (2pt) -- (2 , 2) circle (2pt) ;
\filldraw (1,1) --(0,2) circle (2pt);
\node at (0 , 0) [label=270:\( \emptyset \)]{};
\node at (-1 , 1) [label=180:\( a \)]{};
\node (b) at (1, 1) [label= 0:\( b \)]{};
\node (c) at (0 , 2) [label=180:\( c \)]{};
\node (d) at ( 2 , 2) [label= 0:\( d \)]{};
\node (T) at (0 , -0.5) [label=270:\( T \)]{};
\filldraw (5, 2) circle (2pt) -- (5, 1) circle (2pt) -- (6, 0) circle (2pt) -- (7 , 1) circle (2pt) -- (7 , 2) circle (2pt) -- (8 , 3) circle (2pt) -- (8 , 4) circle (2pt) ;
\filldraw (7,2) --(6,3) circle (2pt) -- (6 , 4) circle (2pt) ;
\node at (6 , 0) [label=270:\( \emptyset \)]{};
\node at (5 , 1) [label=180:\( a^- \)]{};
\node at (5 , 2) [label=180:\( a^+ \)]{};
\node at (7, 1) [label= 0:\( b^- \)]{};
\node at (7, 2) [label= 0:\( b^+ \)]{};
\node at (6 , 3) [label=180:\( c^- \)]{};
\node at (6 , 4) [label=180:\( c^+ \)]{};
\node at ( 8 , 3) [label= 0:\( d^- \)]{};
\node at ( 8 , 4) [label= 0:\( d^+ \)]{};
\node at (6 , -0.5) [label=270:\( T^{\mathrm{d}} \)]{};
\end{tikzpicture}
 \caption{The doubling of a tree \( T \).}
 \label{fig:doubling}
\end{figure}
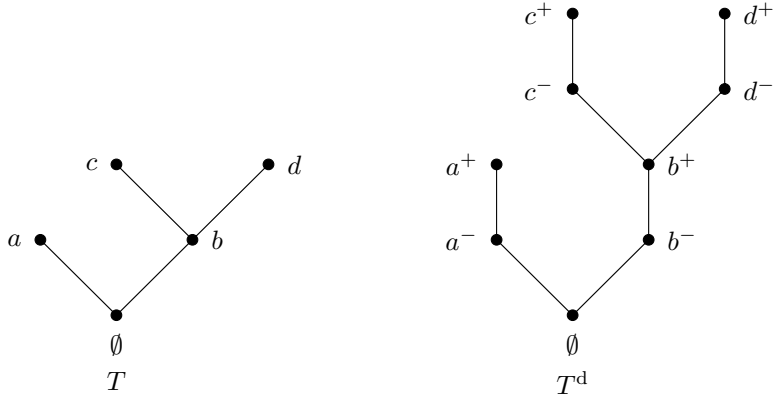

We now start our main construction.
The \markdef{doubling} of a descriptive set-theoretic tree \( T \) is the tree \( T^{\mathrm{d}} \)\index[symbols]{Td@\( T^{\mathrm{d}} \)} obtained by replacing each node \( s \) of \( T \) different from \( \emptyset \) with two nodes \( s^- < s^+ \) in \( T^{\mathrm{d}} \).
Figure~\ref{fig:doubling} shows the doubling of a finite tree.
In order to simplify the notation, the nodes \( s^+ \) of \( T^{\mathrm{d}} \) will simply be denoted by \( s \).
Moreover, \( T^{\mathrm{d}} \) will be always identified with a combinatorial tree on the set of its nodes.

Applying the doubling procedure to \( \pre{< \omega}{\kappa} \), we obtain the combinatorial tree \( \mathbb{G}_0 \).
Formally: 
\begin{definition} \label{def:G_0}
\( \mathbb{G}_0 \)\index[symbols]{Graph0@\( \mathbb{G}_0 \)} is the graph on the disjoint union
\[ 
G_0 \coloneqq \pre{< \omega}{\kappa} \uplus \setofLR{s^-}{\emptyset \neq s \in \pre{< \omega}{\kappa} } \index[symbols]{77@ \( s \mapsto s^- \)}
\] 
with edges 
\[
\edge^{\mathbb{G}_0} \coloneqq \setofLR{ \set{ s , s^-} }{ s \in \pre{< \omega}{\kappa} \setminus \set{ \emptyset} } \cup \setofLR{ \set{ s , (s \conc \alpha)^-} }{ s \in \pre{< \omega}{\kappa} , \, \alpha < \kappa} .
\]
\end{definition}

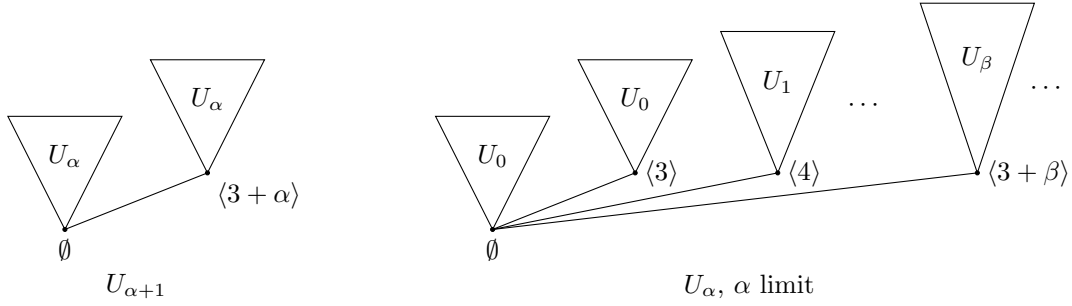
\begin{figure}
\[
\begin{tikzpicture}[scale=0.75]
 \draw(0,0) -- (-1,2)--(1,2)-- (0,0) ;
 \filldraw (0,0) circle (1pt) node[below] {$\emptyset$};
 \node at (0,1.3) {\( U_ \alpha \)};
\draw(2.5,1) -- (1.5,3)--(3.5,3)-- (2.5,1) ;
 \filldraw (2.5,1) circle (1pt) node[below right] {\( \langle 3 + \alpha \rangle \)};
 \node at (2.5,2.3) {\( U_ \alpha \)};
\draw (0,0) -- (2.5,1);
 \node at (1.25,-1) {\( U_ {\alpha + 1} \)};
 \draw (7.5,0) -- (6.5,2)--(8.5,2)-- (7.5,0) ;
 \filldraw (7.5,0) circle (1pt) node[below] {$\emptyset$};
 \node at (7.5,1.3) {\( U_ 0 \)};
 \draw (10,1) -- (9,3)--(11,3)-- (10,1) ;
 \filldraw (10,1) circle (1pt) node[right] {\( \langle 3 \rangle \)};
 \node at (10,2.3) {\( U_ 0 \)};
 \draw (12.5,1) -- (11.5,3.5)--(13.5,3.5)-- (12.5,1) ;
 \filldraw (12.5,1) circle (1pt) node[right] {\( \langle 4 \rangle \)};
 \node at (12.5,2.6) {\( U_ 1 \)};
 \node at (14,2.2) {\( \cdots \)};
 \draw (16,1) -- (15,4)--(17,4)-- (16,1) ;
 \filldraw (16,1) circle (1pt) node[right] {\( \langle 3 + \beta \rangle \)}{};
 \node at (16,3) {\( U_ \beta \)};
 \node at (17.2,2.5) {\( \cdots \)};
\draw (7.5,0) -- (10,1);
\draw (7.5,0) -- (12.5,1);
\draw (7.5,0) -- (16,1);
 \node at (12,-1) {\( U_ \alpha \), \( \alpha \) limit};
\end{tikzpicture}
\]
 \caption{The trees \( U_ \alpha \).}
 \label{fig:Ualpha}
\end{figure}

We next enlarge \( \mathbb{G}_0 \) to the new combinatorial tree \( \mathbb{G}_1 \). 
Given a descriptive set-theoretic tree \( U \subseteq \pre{ < \omega}{\kappa} \), the \markdef{width} of \( U \) is the ordinal 
\[ 
\weight ( U ) \coloneqq \sup \setofLR{\gamma + 1 }{ \seqLR{ \gamma } \in U } .\index[symbols]{w@\( \weight \)}
\] 
We construct a sequence \( \seqofLR{ U_ \alpha }{ \alpha < \kappa} \) of descriptive set-theoretic trees on \( \kappa \) as follows (see Figure~\ref{fig:Ualpha}):\index[symbols]{Ualpha@\( U_ \alpha \)} 
\begin{align*}
U_0 &= \pre{ 3}{3},
\\
U_{ \alpha + 1} & = U_ \alpha \cup \setofLR{ \seqLR{ \weight ( U_ \alpha ) } \conc s}{s \in U_\alpha }
\\
U_ \alpha & = \bigcup_{\beta < \alpha} U_\beta && \alpha \text{ limit.}
\end{align*}
By construction, for every \( \alpha < \kappa \) 
\begin{equation} \label{eq:Ualphaprop} 
\weight ( U_\alpha ) = 3 + \alpha < \kappa \quad \text{and} \quad U_ \alpha \subseteq \pre{ < \omega }{ (\weight ( U_\alpha ))} = \pre{ < \omega }{ (3+\alpha)}.
\end{equation}
Each \( U_\alpha \) is identified with the corresponding combinatorial tree (that is with the graph on the nodes of \( U_\alpha \) obtained by linking with an edge all pairs of nodes \( x , y \in U_\alpha \) such that \( x = y^\star \), for \( y^\star \) as in~\eqref{eq:predecessorofu}).
The vertex \( \emptyset \) is called the root of \( U_\alpha \). 

\begin{remark}\label{rmk:distinctstructuresinZF}
Each \( U_ \alpha \) can be construed as a combinatorial tree (with root \( \emptyset \)).
If \( \kappa = \lambda ^+ \), then using~\eqref{eq:Ualphaprop} we get that \( \card{ U_ \alpha } \leq \lambda \) for all \( \alpha < \kappa \), so each \( U_ \alpha \) is isomorphic to an element of \( \CT_ \lambda \) (at least when \( \lambda \leq \alpha \)).
Therefore by choosing a bijection between \( U_ \alpha \) and \( \lambda \) we could construct the \( U_ \alpha \)'s so that they belong to \( \CT_ \lambda \), when \( \alpha \geq \lambda \).
The appeal to \( \AC \) cannot be avoided --- if the \( \PSP \) holds and \( \kappa = \omega _1 \), then there is no \( \omega_1 \)-sequence of distinct elements of \( \CT_{ \omega } \).
In our construction we avoided this issue by relaxing the requirement that a combinatorial tree of size \( \lambda \) have domain equal to \( \lambda \). 
\end{remark}

We now define the combinatorial tree \( \mathbb{G}_1 \) by connecting a new vertex \( \hat{s} \) to each \( s \in \pre{<\omega}{\kappa} \subseteq G_0 \) and then appending a copy of \( U_{\Code{s}} \) to \( \hat{s} \) (where \( \Code{\cdot} \) is the function coding sequences by ordinals from~\eqref{eq:codingfinitesequenceofordinals}) by adding an edge between such vertex and the root of \( U_{\Code{ s}} \). 
More formally, 
\begin{definition} \label{def:G_1}
\( \mathbb{G}_1 \)\index[symbols]{Graph1@\( \mathbb{G}_1 \)} is the combinatorial tree with set of vertices
\[ 
G_1 \coloneqq G_0 \uplus \setofLR{\hat{s}}{s \in \pre{< \omega}{\kappa}} \uplus \setofLR{ ( s , t ) \in \pre{< \omega}{\kappa} \times \pre{< \omega}{\kappa}}{t \in U_{\Code{s}}}
 \] 
and edge relation \( \edge^{\mathbb{G}_1} \) defined by:
\begin{itemize}[leftmargin=1pc]
\item
if \( x , y \in G_0 \), \( x \edge^{\mathbb{G}_1} y \IFF x \edge^{\mathbb{G}_0} y \);
\item
 \( s \edge^{\mathbb{G}_1} \hat{s} \) and \( \hat{s} \edge^{\mathbb{G}_1} s \) for every \( s \in \pre{<\omega}{\kappa} \);
\item
\( \hat{s} \edge^{\mathbb{G}_1} ( s , \emptyset ) \) and \( ( s , \emptyset ) \edge^{\mathbb{G}_1} \hat{s} \) for every \( s \in \pre{< \omega}{\kappa} \);
\item
for every \( s \in \pre{<\omega}{\kappa} \) and \( t , t' \in U_{\Code{s}} \), \( ( s , t ) \edge^{\mathbb{G}_1} ( s , t' ) \IFF t \edge^{U_{\Code{ s}}} t' \);
\item
no other \( \edge^{\mathbb{G}_1} \)-relation holds.
\end{itemize} 
\end{definition}

As explained above, the purpose of moving from \( \mathbb{G}_0 \) to \( \mathbb{G}_1 \) is to obtain a more rigid structure. 
In order to prove that \( \mathbb{G}_1 \) has few nontrivial automorphisms, we first need to consider some technical properties of the \( U_\alpha \)'s which will also be useful in the proof of our invariant universality result (Section~\ref{sec:invariantlyuniversal}). 
Such technical analysis is further complicated by the fact that for our applications we need to work in \( \ZF \), and hence we have to ensure that all ``choices'' needed in the constructions can be done in a canonical way. 
However these technical details are only needed for the proof of invariantly universality: the reader who is just interested in the completeness result (Section~\ref{sec:embeddabilitygraphs}) may safely skip the rest of this section and directly jump to Section~\ref{subsec:G_T}. 

The relevant properties of the graphs \( U_\alpha \) are summarized in the next lemma.

\begin{lemma} \label{lem:propertiesL_alpha}
Let \( \alpha,\beta < \kappa \). 
\begin{enumerate-(a)}
\item \label{lem:propertiesL_alpha-a}
Except for the root \( \emptyset \) (which has degree \( 3 \) if \( \alpha = 0 \) and \( \geq 4 \) otherwise), all vertices of \( U_\alpha \) are either terminal (i.e.\ they have a unique neighbor) or have degree \( \geq 4 \).
Moreover, each vertex of \( U_ \alpha \) has degree \( \leq 1 + \weight ( U_ \alpha ) = 4 + \alpha < \kappa \).
\item \label{lem:propertiesL_alpha-b}
If \( \alpha \leq \beta \), then \( U_\alpha \subseteq U_\beta \). 
In particular, the identity map embeds \( U_\alpha \) into \( U_\beta \) and fixes the root \( \emptyset \).
\item \label{lem:propertiesL_alpha-c}
For \( \beta < \alpha < \kappa \) there is no embedding of \( U_\alpha \) into \( U_\beta \) sending \( \emptyset \) to itself. 
In particular, there is an isomorphism between \( U_\alpha \) and \( U_\beta \) fixing \( \emptyset \) if and only if \( \alpha = \beta \).
\item \label{lem:propertiesL_alpha-d}
There is no infinite path through \( U_\alpha \). 
(Equivalently, \( U_\alpha \) is well-founded when construed as a descriptive set-theoretic tree on \( \kappa \) of height \( \leq \omega \).)
\item \label{lem:propertiesL_alpha-e}
Let \( i \) be an automorphism of \( U_\alpha \). 
If \( \alpha = \gamma + 1 \) is a successor ordinal, then either \( i ( \emptyset ) = \emptyset \) or \( i ( \emptyset ) = \seqLR{ \weight ( U_\gamma ) } \) (and both the possibilities may be realized), while if \( \alpha = 0 \) or \( \alpha \) is limit then \( i ( \emptyset ) = \emptyset \).
\item \label{lem:propertiesL_alpha-f}
Let \( \mathbb{X}_ \alpha \) be the set of all trees \( X \) with domain \( \subseteq \kappa \) that are isomorphic to \( U_\alpha \).
Then there is a definable map \( \mathbb{X}_ \alpha \ni X \mapsto i_{X , \alpha} \) such that \( i_{ X, \alpha } \colon X \to U_\alpha \) is an isomorphism.
(We call such \( i_{ X, \alpha } \) the canonical isomorphism between \( X \) and \( U_ \alpha \).)
\end{enumerate-(a)}
\end{lemma}

\begin{proof}
Parts~\ref{lem:propertiesL_alpha-a},~\ref{lem:propertiesL_alpha-b} and~\ref{lem:propertiesL_alpha-d} are obvious and directly follow from the construction of the \( U_\alpha \)'s and~\eqref{eq:Ualphaprop}.

\smallskip

A combinatorial tree \( G \) with root \( r \) and \emph{without infinite paths} can be seen as a well-founded descriptive set-theoretic tree, and hence it has a rank function \( \rho_{G , r } \colon G \to \On \) defined by
\[ 
\rho_{G , r} ( a ) \coloneqq \sup \setofLR{ \rho_{ G , r } ( b ) + 1}{ b \neq a \text{ and the unique path joining \( r \) to \( b \) passes through \( a \)}}.
\]
Thus \( \rho_{G , r } ( a ) = 0 \) if and only if \( a \) has degree \( 1 \) in \( G \), and then set
\[ 
\rho ( G , r ) \coloneqq \sup \setofLR{ \rho_{ G , r } ( a ) + 1 }{a \in G, \, a \neq r} = \rho_{ G , r } ( r ).
 \] 
Notice that for every two rooted combinatorial trees \( ( G , r ) \), \( ( G ' , r ' ) \), if \( ( G , r ) \embeds ( G ' , r ' ) \) then \( \rho ( G , r ) \leq \rho ( G ' , r ' ) \) (whence in particular \( \rho ( G , r ) = \rho ( G' , r ' ) \) whenever \( ( G , r ) \cong ( G ' , r ' ) \)). 
Then~\ref{lem:propertiesL_alpha-c} easily follows from the fact that one can show by induction on \( \alpha \) that \( \rho ( U_\alpha , \emptyset ) = 3 + \alpha \).

\smallskip

We now prove~\ref{lem:propertiesL_alpha-e}.
Assume first \( \alpha = \gamma + 1 \), and let \( r \in U_\alpha \) be arbitrary. 
Then either \( \rho_{U_ \alpha , r } ( \emptyset ) = \rho_{ U_\alpha , \emptyset } ( \emptyset ) \) (if \( \seqLR{ \weight ( U_\gamma ) } \nsubseteq r \)), or else \( \rho_{U_\alpha , r } ( \weight ( U_\gamma ) ) = \rho_{ U_\alpha , \emptyset } ( \emptyset ) \) (if \( \seqLR{ \weight ( U_\gamma ) } \subseteq r \)). 
Thus if \( r \notin \{ \emptyset, \seqLR{ \weight ( U_\gamma ) } \} \), in both cases we would easily get \( \rho_{ U_\alpha , r } ( r ) > \rho_{ U_\alpha , \emptyset } ( \emptyset ) \), and hence \( \rho ( U_\alpha , r ) > \rho ( U_ \alpha , \emptyset ) \). 
Since \( i \) witnesses \( ( U_\alpha , \emptyset ) \cong ( U_\alpha , i ( \emptyset ) ) \), so that \( \rho ( U_\alpha, \emptyset ) = \rho ( U_\alpha, i ( \emptyset ) ) \), setting \( r \coloneqq i ( \emptyset ) \) in the argument above we get that \( i ( \emptyset ) \in \{ \emptyset, \seqLR{ \weight ( U_\gamma )} \} \), as required. 
To see that both the possibilities can be realized, consider the identity map and the isomorphism \( i \colon U_\alpha \to U_\alpha \) defined by
\[ 
i ( s ) \coloneqq
\begin{cases}
\weight ( U_\gamma ) \conc s & \text{if } s \in U_\gamma 
\\
t & \text{if }s = \weight ( U_\gamma ) \conc t .
\end{cases}
\] 

The case of \( \alpha = 0 \) or \( \alpha \) limit is similar: by construction, \( \rho_{ U_\alpha , i ( \emptyset )}( \emptyset ) = \rho_{U_\alpha , \emptyset } ( \emptyset ) \) independently of the value of \( i ( \emptyset ) \), and therefore if \( i ( \emptyset ) \neq \emptyset \) then \( \rho ( U_\alpha , i ( \emptyset ) ) > \rho ( U_\alpha , \emptyset ) \), contradicting the fact that \( i \) witnesses \( ( U_\alpha, \emptyset) \cong (U_\alpha, i ( \emptyset ) ) \).

\smallskip

Finally we prove~\ref{lem:propertiesL_alpha-f} by induction on \( \alpha < \kappa \). 
First a technical claim.

\begin{claim} \label{claim:rootboth}
Let \( \bar{\imath} \colon X \to U_\alpha \) be an arbitrary isomorphism.
\begin{enumerate-(i)}
\item \label{claim:rootboth-i}
Assume \( \alpha = \gamma+1 \) is a successor ordinal.
For any isomorphism \( i \colon X \to U_\alpha \)
\[
\set{ \bar{\imath}^{-1}(\emptyset) , \bar{\imath}^{-1} ( \seqLR{ \weight ( U_\gamma ) } )} = \set{ i^{-1} ( \emptyset ) , i^{-1} ( \seqLR{ \weight ( U_\gamma )} ) } .
\]
\item \label{claim:rootboth-ii}
Assume that \( \alpha \) is either \( 0 \) or a limit ordinal, and set \( \bar{\delta} \coloneqq \bar{\imath}^{-1} ( \emptyset ) \). 
Then \( i ( \bar{\delta} ) = \emptyset \) for every isomorphism \( i \colon X \to U_\alpha \).
\end{enumerate-(i)}
\end{claim}

\begin{proof}[Proof of the Claim]
\ref{claim:rootboth-i} Assume towards a contradiction that the claim fails for some isomorphism \( i \colon X \to U_\alpha \). 
Using the fact that by~\ref{lem:propertiesL_alpha-e} there is an automorphism of \( U_\alpha \) sending \( \seqLR{ \weight ( U_\gamma ) } \) to \( \emptyset \), we may assume without loss of generality that \( i^{-1}(\emptyset) \neq \bar{\imath}^{-1}(\emptyset) , \bar{\imath}^{-1} ( \seqLR{ \weight ( U_\gamma ) } ) \). 
But then \( i' \coloneqq \bar{\imath} \circ i^{-1} \) would be an automorphism of \( U_\alpha \) such that \( i'(\emptyset) \notin \{\emptyset, \seqLR{ \weight ( U_\gamma ) } \} \), contradicting part~\ref{lem:propertiesL_alpha-e}.

\ref{claim:rootboth-ii} The case \( \alpha =0 \) is trivial. 
The limit case can be treated similarly to the successor one: if \( i ( \bar{\delta} ) \neq \emptyset \) for some isomorphism \( i \), then we would also have an automorphism of \( U_\alpha \) which does not fix its root \( \emptyset \), contradicting again part~\ref{lem:propertiesL_alpha-e}. 
\end{proof}

We now come back to the proof of part~\ref{lem:propertiesL_alpha-f} of the lemma.
The case \( \alpha = 0 \) is trivial since \( U_0 \) is finite. 
Let \( \alpha = \gamma + 1 \). 
Fix an arbitrary isomorphism \( \bar{\imath} \colon X \to U_\alpha \) and let \( \bar{\delta}_0 = \min \setLR{ \bar{\imath}^{-1}(\emptyset) , \bar{\imath}^{-1} ( \seqLR{ \weight ( U_\gamma ) } ) } \) and \( \bar{\delta}_1 = \max \setLR{ \bar{\imath}^{-1}(\emptyset) , \bar{\imath}^{-1} ( \seqLR{ \weight ( U_\gamma ) } ) }\).
(The choice of \( \bar{\delta}_0 , \bar{\delta}_1 \) is uniquely determined by Claim~\ref{claim:rootboth}\ref{claim:rootboth-i}.) 
For the sake of definiteness, assume \( \bar{\delta}_0 = \bar{\imath}^{-1}(\emptyset) \) --- if this is not the case, just swap the following definitions of \( X_0 \) and \( X_1 \).
Set \( X_0 \coloneqq \bar{\imath}^{-1} ( U_\gamma ) \) and \( X_1 \coloneqq X \setminus X_0 \), so that \( \bar{\delta}_0 \in X_0 \) and \( \bar{\delta}_1 \in X_1 \). 
Notice that by Claim~\ref{claim:rootboth}\ref{claim:rootboth-i} the definition of \( X_0 , X_1 \) is independent of the choice of \( \bar{\imath} \), therefore no choice is needed here. 
Since \( \bar{\imath} \) witnesses that both \( X_0 \) and \( X_1 \) are isomorphic to \( U_\gamma \), by inductive hypothesis there are canonical isomorphisms \( i_{X_0 , \gamma} \colon X_0 \to U_\gamma \) and \( i_{X_1 , \gamma } \colon X_1 \to U_\gamma \). 
Notice that the isomorphism between \( X_j \) and \( U_\gamma \) induced by \( \bar{\imath} \restriction X_j \) sends \( \bar{\delta}_j \) to \( \emptyset \), so by Claim~\ref{claim:rootboth} (applied to \( \gamma \)) and part~\ref{lem:propertiesL_alpha-e} of the lemma we may assume without loss of generality that \( i_{X_j , \gamma } ( \bar{\delta}_j ) = \emptyset \) (for \( j = 0 , 1 \)). 
Setting for \( \delta \in X \)
\[ 
i_{X,\alpha} ( \delta ) \coloneqq 
\begin{cases}
i_{X_0 , \gamma } ( \delta ) & \text{if } \delta \in X_0 
\\
\seqLR{ \weight ( U_\gamma ) } \conc i_{X_1,\gamma} ( \delta ) & \text{if } \delta \in X_1
\end{cases} 
 \] 
we get the desired canonical isomorphism between \( X \) and \( U_\alpha \).

Finally, assume that \( \alpha \) is limit, and notice that in this case \( U_\alpha \) is the disjoint union of \( U_0 \) and all subtrees of \( U_\alpha \) with domain \( \seqLR{ 3 + \beta } \conc U_{\beta} = \setofLR{ \seqLR{ 3 + \beta } \conc t } { t \in U_{\beta} } \) for \( \beta < \alpha \).
Fix an arbitrary isomorphism \( \bar{\imath} \colon X \to U_\alpha \), and let \( \bar{\delta} \coloneqq \bar{\imath}^{-1} ( \emptyset ) \) be as in Claim~\ref{claim:rootboth}\ref{claim:rootboth-ii}. 
Set also, \( X_{\beta} \coloneqq \bar{\imath}^{-1}( \seqLR{ 3 + \beta } \conc U_{\beta} ) \) (for every \( \beta < \alpha \)) and \( X_{-1} \coloneqq X \setminus \bigcup_{\beta < \alpha} X_{\beta} \), so that 
\begin{equation} \label{eq:X_beta} 
X_{\beta} \cong U_{\beta} \quad \text{and} \quad X_{-1} \cong U_0.
\end{equation} 

\begin{claim} \label{claim:rootlimit}
For every isomorphism \( i \colon X \to U_\alpha \) and every \( \beta < \alpha \), \( i ( \bar{\delta} ) = \emptyset \), \( i ( X_{\beta} ) = \seqLR{ 3 + \beta }\conc U_{\beta} \) and \( i ( X_{-1} ) = U_0 \).
\end{claim}

\begin{proof}
Fix an arbitrary \( \beta < \alpha \). 
Since \( i ( \bar{\delta} ) = \emptyset \) by Claim~\ref{claim:rootboth}\ref{claim:rootboth-ii}, we have that \( i ( X_{ \beta } ) = \setofLR{ t \in U_\alpha}{ \seqLR{ \beta' } \subseteq t } \) for some \( \beta' < 3 + \alpha = \weight ( U_\alpha ) \).
Thus we just need to show that \( \beta' = 3 + \beta \).
Since \( X_{\beta} \) contains at least \( \card{U_0} = 13 \) points we have \( \beta' \geq 3 \), and thus \( i ( X_{\beta}) = \beta' \conc U_{\beta''} \) where \( \beta' = 3 + \beta'' \). 
Using \( i \), \( \bar{\imath} \), and~\eqref{eq:X_beta}, we get that \( U_{\beta''} \cong X_{\beta} \cong U_{\beta} \) via an isomorphism sending \( \emptyset \) to itself, so \( \beta = \beta'' \) by part~\ref{lem:propertiesL_alpha-c} of the lemma. 
Thus \( \beta' = 3 + \beta \), as required. 
The final part concerning \( X_{-1} \) and \(U_0 \) follows from the preceding one, so we are done.
\end{proof}

By Claims~\ref{claim:rootboth} and~\ref{claim:rootlimit}, the definition of \( \bar{\delta} \), of the \( X_{\beta} \)'s, and of \( X_{-1} \) is independent of the chosen \( \bar{\imath} \), so no choice is needed to define them. 
For every \( \beta < \alpha \), let \( \bar{\delta}_{\beta} \) be the unique point in \( X_{\beta} \) which is connected by an edge of \( X \) to \( \bar{\delta} \), and apply the inductive hypothesis to get canonical isomorphisms \( i_{X_{-1},0} \colon X_{-1} \to U_0 \) and \( i_{X_{\beta}, \beta} \colon X_{\beta} \to U_{\beta} \). 
Since \( \bar{\imath} \) induces an isomorphism between \( X_{\beta} \) and \( U_{\beta} \) sending \( \bar{\delta}_{\beta} \) to \( \emptyset \), by Claim~\ref{claim:rootboth} (applied to \( \beta \)) and part~\ref{lem:propertiesL_alpha-e} of the lemma we may assume without loss of generality that \( i_{X_{\beta},\beta}( \bar{\delta}_{\beta} ) = \emptyset \) as well. 
Similarly, using the fact that \( \bar{\iota} \restriction X_{-1} \) is an isomorphism between \( X_{-1} \) and \( U_0 \) sending \( \bar{\delta} \) to \( \emptyset \), we also have that \( i_{X_{-1},0}(\bar{\delta}) = \emptyset \). 
Therefore setting
\[ 
i_{X,\alpha}(\delta) \coloneqq
\begin{cases}
i_{X_{-1},0}(\delta) & \text{if } \delta \in X_{-1} 
\\
\seqLR{ 3 + \beta } \conc i_{X_{\beta},\beta}(\delta) & \text{if } \delta \in X_{\beta} \text{ for some } \beta < \alpha,
\end{cases}
 \] 
we get the desired canonical isomorphism \( i_{X,\alpha} \colon X \to U_\alpha \).
\end{proof}

We next use some of the properties of the \( U_\alpha \)'s described in Lemma~\ref{lem:propertiesL_alpha} to obtain a rigidity property of \( \mathbb{G}_1 \) which will be crucial in the proof of the invariant universality of \( \embeds^\kappa_\CT \).

\begin{lemma}\label{lem:quasirigid}
Every automorphism of \( \mathbb{G}_1 \) is the identity on \( G_0 \).
\end{lemma}

\begin{proof}
Let \( j \) be an automorphism of \( \mathbb{G}_1 \), and for every \( s \in \pre{< \omega}{\kappa} \) let \( U_s \subseteq G_1 \) be the copy of \( U_{\Code{s}} \) in \( \mathbb{G}_1 \), that is the substructure of \( \mathbb{G}_1 \) with domain \( \setofLR{(s,t) \in \pre{< \omega}{\kappa} \times \pre{< \omega}{\kappa}}{t \in U_{\Code{s}}} \) (see Definition~\ref{def:G_1}).
It suffices to prove that \( j ( s ) = s \) for every \( s \in \pre{<\omega}{\kappa} \subseteq G_0 \subseteq G_1 \). 
First notice that \( j ( \pre{<\omega}{\kappa}) = \pre{<\omega}{\kappa} \) because \( \pre{<\omega}{\kappa} \subseteq G_1 \) is the set of all elements of \( \mathbb{G}_1 \) of degree \( \kappa \) (by Lemma~\ref{lem:propertiesL_alpha}\ref{lem:propertiesL_alpha-a} and the fact that all vertices of the form \( s^- \) or \( \hat{s} \) have degree \( 2 \), see also the subsequent Lemma~\ref{lem:G_T}).
In particular, \( j \) maps \( G_1 \setminus \pre{<\omega}{\kappa} \) in itself.
Consider the point \( ( s , \emptyset ) \in G_1 \setminus \pre{< \omega}{\kappa} \), which is the root of \( U_s \): since it has degree \( \geq 4 \) and distance \( 2 \) from \( s \), we must have \( j ( s , \emptyset ) = ( j ( s ) , \emptyset ) \in U_{ j ( s )} \) (recall that necessarily \( j (G_1 \setminus \pre{< \omega}{\kappa}) \subseteq G_1 \setminus \pre{< \omega}{\kappa} \)). 
It follows that \( j \restriction U_s \) is an isomorphism between \( U_s \) and \( U_{j ( s ) } \) which sends the root of \( U_s \) to the root of \( U_{ j ( s ) } \). Therefore \( \Code{s} = \Code{j ( s ) } \) by Lemma~\ref{lem:propertiesL_alpha}\ref{lem:propertiesL_alpha-c}, and thus \( j ( s ) = s \) by injectivity of \( \Code{\cdot} \).
\end{proof}

\subsection{The combinatorial trees \( \mathbb{G}_S \)} \label{subsec:G_T}

As described at the beginning of this section, the combinatorial tree \( \mathbb{G}_S \)\index[symbols]{Graph2@\( \mathbb{G}_S \)} associated to some \( S \in \Tr ( 2 \times \kappa ) \) is obtained by adding some forks to \( \mathbb{G}_1 \).

\begin{figure} 
 \centering
\begin{tikzpicture}
\filldraw (0, 0) circle (2pt) -- (1 , 0) circle (2pt) -- (2 , 0) circle (2pt) -- (3 , 0) circle (2pt) -- ( 4, 0 ) 
circle (2pt) -- (5 , 0) circle (2pt)-- (6 , 0) circle (2pt)-- (7 , 0) circle (2pt) ;
\filldraw[loosely dotted] (7,0)--(9,0);
\filldraw (3,0) --(3,1) circle (2pt);
\node at (3 , 1) [label=90:\( 0^{(3)}1 \)]{};
\node at (0 , 0) [label=270:\( \emptyset \)]{};
\node at (1 , 0) [label=270:\( 0^{(1)} \)]{};
\node at (2 , 0) [label=270:\( 0^{(2)} \)]{};
\node at (3 , 0) [label=270:\( 0^{(3)} \)]{};
\node at (4 , 0) [label=270:\( 0^{(4)} \)]{};
\node at (5 , 0) [label=270:\( 0^{(5)} \)]{};
\node at (6 , 0) [label=270:\( 0^{(6)} \)]{};
\node at (7 , 0) [label=270:\( 0^{(7)} \)]{};
\end{tikzpicture}
 \caption{The fork \( F_u \), when \( \theta ( u ) = 3 \) (i.e.~\( u = \emptyset \))}
 \label{fig:T_u}
\end{figure}
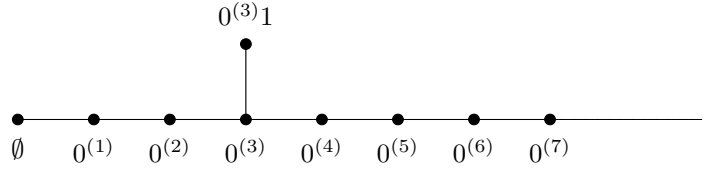

Fix any injection \( \theta \colon \pre{< \omega}{ 2} \to \omega \) such that 
\begin{subequations}
\begin{align}
& \theta ( u ) \text{ is odd for all } u \in \pre{< \omega}{2} , \label{eq:thetaoddvalues} 
\\ 
& \theta ( \emptyset ) = 3 \text{ and } \theta ( u ) > \theta ( \emptyset ) \text{ for all } u \in \pre{<\omega}{2} \setminus \{ \emptyset \} , \label{eq:theta0=3} 
\\ 
& \card{\theta(u) - \theta ( v ) } > 4 \cdot \max \setLR{ \lh u , \lh v } \text{ for all distinct } u , v \in \pre{<\omega}{2} . \label{eq:thetadistances}
\end{align} 
\end{subequations}
Such an injection can be easily constructed by induction on the length of \( u \in \pre{< \omega}{2} \) by e.g.\ setting \( \theta(\emptyset) \coloneqq 3 \) and by letting \( \theta \restriction \pre{n+1}{2} \) be an arbitrary injection into the set 
\[ 
\setofLR{ \max \setofLR{\theta(u)}{u \in \pre{n}{2}} + 5 k ( n + 1 )}{1 \leq k \leq 2^{n+1}} .
\]
For \( u \in \pre{< \omega}{2} \) we define the \markdef{fork (coding \( u \))} to be the graph \( F_u \) on 
\[ 
\setofLR{ w \in \pre{< \omega}{ 2} }{ w \subseteq 0^{( \omega )} \vee w = 0^{( \theta ( u ) )} \conc 1 }
\]
connecting each sequence \( w \neq \{ \emptyset \} \) to its immediate predecessor \( w^\star = w \restriction ( \lh w ) - 1 \) (Figure~\ref{fig:T_u}).
Then \( F_u \) is a combinatorial tree consisting of three disjoint branches departing from \( 0^{( \theta ( u ) )} \), which is the unique vertex of degree \( 3 \); one branch is infinite, one has length \( 1 \), and one has odd length \( \theta ( u ) > 1 \).
Since any embedding respects degrees,
\begin{subequations}
\begin{align}
& u \neq v \IMPLIES \text{there is no embedding of \( F_u \) into \( F_v \) fixing } \emptyset, \label{eq:fixempty}
\\
& F_u \text{ is rigid, i.e.~its unique automorphism is the identity,} \label{eq:F_urigid}
\\
& u = v \IFF F_u \cong F_v . \label{eq:=iffiso}
\end{align}
\end{subequations}
Using the fact that the domain of \( F_u \) is an infinite subset of \( \pre{ < \omega }{2} \), the graph \( F_u \) can be then identified with a graph on \( \omega \) via a canonical bijection 
\begin{equation}\label{eq:e_u}
 e_u \colon \dom ( F_u ) \to \omega ,
\end{equation}
sending \( \emptyset \) to \( 0 \).

\begin{figure}
 \centering
\begin{tikzpicture}[scale=1.5]
 \draw[->] (0,0) -- (3,-1) node[right] {$x$} ;
 \draw[->] (0,0) -- (0,2) node[above] {$z$};
\draw[->] (0,0) -- (2.4,1.2) node[above] {$y$} ;
\filldraw[fill=gray!30,draw=black] (2.2,0.9) -- (1,0) -- (3,-0.3) --cycle node[right] {$G_1$};
\fill[black] (1.7,0) circle (1pt) node[below] {$s$};
\draw[thick] (1.7,0) -- (1.6,2) node[left] {$F_{u,s}$};
\draw[thick] (1.6,0.25) -- (1.68,0.25) ;
\draw[thick] (1.7,0) -- (1.8,2) node[right] {$F_{v,s}$};
\draw[thick] (1.8,0.5) -- (1.725,0.5);
\end{tikzpicture}
 \caption{The graph \( \mathbb{G}_S \).}
 \label{fig:G_T}
\end{figure}
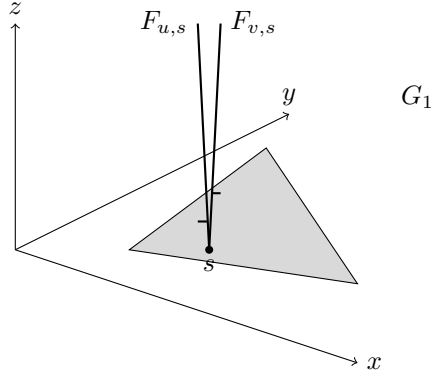

For each \( ( u , s ) \in \pre{<\omega}{(2 \times \omega)} \), we fix an isomorphic copy \( F_{u , s} \) of \( F_u \), so that the sets of vertices of \( F_{u , s} \) and \( F_{v , t} \) are disjoint whenever \( (u,s) \neq (v,t) \), and such that each \( F_{u,s} \) is also disjoint from the domain \( G_1 \) of \( \mathbb{G}_1 \).
More precisely, we let 
\begin{equation}
F_{u , s}= \setLR{( u , s)} \times F_u \label{eq:F_{u,s}} 
\end{equation}
be the graph on \( \setof{ ( u , s , w ) }{ w \subseteq 0^{( \omega )} \vee w = 0^{( \theta ( u ) )} \conc 1 } \) with set of edges defined by
\[ 
( u , s , w ) \mathrel{F_{u , s}} ( u , s , z ) \IFF w \mathrel{ F_u} z . 
\]
Following~\cite[Theorem 3.1]{Louveau:2005cq} and~\cite[Theorem 3.9]{Friedman:2011cr}, to each \( S \in \Tr ( 2 \times \kappa ) \) we now associate the combinatorial tree \( \mathbb{G}_S \) obtained by joining \( \mathbb{G}_1 \) and the \( F_{u , s} \) for \( ( u , s ) \in S \) via the identification of any vertex \( s \in \pre{<\omega}{\kappa} \subseteq G_1 \) with each vertex of the form \( ( u , s , \emptyset ) \).
Thus the domain \( G_S \) of \( \mathbb{G}_S \) can be identified with 
\begin{multline*}
\pre{ < \omega }{ \kappa } \uplus \setofLR{ s^-}{ \emptyset \neq s \in\pre{ < \omega }{ \kappa } } \uplus \setofLR{ \hat{s}}{ \emptyset \neq s \in\pre{ < \omega }{ \kappa } } \uplus \setofLR{ ( s , t ) \in \pre{<\omega}{\kappa} \times \pre{ < \omega}{\kappa}}{t \in U_{\Code{ s}}} 
 \\
\uplus \setofLR{( u , s , w ) }{ ( u , s ) \in S \AND ( u , s , w ) \in F_{ u , s} \AND w \neq \emptyset} .
\end{multline*} 
(Figure~\ref{fig:G_T} is a clumsy attempt to visualize \( \mathbb{G}_S \) in the three-dimensional space: \( \mathbb{G}_1 \) is the grey area in the \( x y \)-plane, \( s \in \pre{ < \omega }{ \kappa } \) is a vertex of \( \mathbb{G}_1 \) in the set \( \pre{< \omega}{\kappa} \subseteq G_0 \subseteq G_1 \), and the forks \( F_{u , s} , F_{ v , s} \) are growing vertically.)

To simplify the presentation, we also introduce the following notation for the relevant subsets of the domain \( G_S \) of \( \mathbb{G}_S \) (when necessary, each of these sets is also identified with the corresponding substructure of \( \mathbb{G}_S \)):
\begin{subequations} \label{eq:substructuresofG_T}
\begin{align}
 \mathrm{Seq} ( \mathbb{G}_S ) &\coloneqq \setofLR{ s \in G_S}{s \in \pre{ < \omega}{ \kappa} } ; \label{eq:substructuresofG_T-a}
\\
 \mathrm{Seq}^- ( \mathbb{G}_S ) &\coloneqq \setofLR{ s^- \in G_S}{s \in \pre{ < \omega}{ \kappa} } ; \label{eq:substructuresofG_T-b}
\\
 \widehat{\mathrm{Seq}}( \mathbb{G}_S ) & \coloneqq \setofLR{\hat{s} \in G_S}{s \in \pre{ < \omega}{\kappa}} ; \label{eq:substructuresofG_T-c}
\\
 \mathrm{U}_s ( \mathbb{G}_S ) & \coloneqq \setofLR{ ( s , t ) \in \pre{<\omega}{\kappa} \times \pre{< \omega}{\kappa}}{ t \in U_{\Code{ s} }} \text{ for every } s \in \pre{< \omega}{\kappa} ; \label{eq:substructuresofG_T-d}
\\
 \mathrm{U} ( \mathbb{G}_S ) & \coloneqq \bigcup \setofLR{ \mathrm{U}_s ( \mathbb{G}_S ) }{s \in \pre{< \omega }{\kappa}} ; \label{eq:substructuresofG_T-e}
\\ 
 \mathrm{F}'_{ u , s } ( \mathbb{G}_S ) & \coloneqq \setofLR{ ( u , s , w ) \in F_{u , s } }{w \neq \emptyset } \text{ for every } ( u , s ) \in S ; \label{eq:substructuresofG_T-f}
\\
 \mathrm{F}(\mathbb{G}_S) &\coloneqq \bigcup \setofLR{ \mathrm{F}'_{ u , s } ( \mathbb{G}_S ) }{(u,s) \in S} . \label{eq:substructuresofG_T-g}
\end{align}
\end{subequations}

\begin{remark} \label{rmk:substructuresofG_T}
The subsets and the corresponding substructures of \( \mathbb{G}_S \) defined in~\eqref{eq:substructuresofG_T-a}--\eqref{eq:substructuresofG_T-f} do not depend at all on the specific \( S \in \Tr ( 2 \times \kappa ) \) under consideration, but only on the parameters \( s \) and \( u \) (when they appear in the notation). 
Therefore, to further simplify the notation we can safely set
\[ 
\mathrm{Seq} \coloneqq \mathrm{Seq} ( \mathbb{G}_S ) \qquad (\text{for some/any } S \in \Tr ( 2 \times \kappa ) ) ,
 \] 
and define in a similar way the structures \( \mathrm{Seq}^- \), \( \widehat{\mathrm{Seq}} \), \( \mathrm{U}_s \), \( \mathrm{U} \), \( \mathrm{F}'_{ u , s } \), and \( \mathrm{F} \). 
\end{remark}

In particular, we have \( G_0 = \mathrm{Seq} \uplus \mathrm{Seq}^- \), \( G_1 = G_0 \uplus \widehat{\mathrm{Seq}} \uplus \mathrm{U} \), and \( G_S = G_1 \uplus \mathrm{F} \). 
Finally, the graph \( \mathbb{G}_S \) is then further identified in a canonical way with its copy on \( \kappa \) (which is thus an element of \( \CT_\kappa \)), using the bijections \( \Code{\cdot} \), \( \op{ \cdot}{\cdot } \), and \( e_u \) in the obvious way. 

We end this section by listing some properties of the vertices of \( \mathbb{G}_S \) in term of distances and degrees which will be useful in the subsequent proofs (Sections~\ref{sec:embeddabilitygraphs} and~\ref{sec:invariantlyuniversal}), leaving to the reader to check their validity (using when necessary Lemma~\ref{lem:propertiesL_alpha}).
Given a combinatorial tree \( G \in \CT_\kappa \), a \markdef{\( \Z \)-chain of \( G \)} is a sequence \( \seqofLR{x_z}{z \in \Z} \) of pairwise distinct elements of \( G \) such that \( x_z \edge^{G} x_{z + 1} \) for every \( z \in \Z \).

\begin{lemma} \label{lem:G_T}
Let \( S \in \Tr (2 \times \kappa ) \).
\begin{enumerate-(a)}
\item \label{lem:G_T-a}
The vertices in \( \mathrm{Seq} \) are the unique vertices of \( \mathbb{G}_S \) having degree \( \kappa \); they are also the unique vertices of \( \mathbb{G}_S \) with degree \( \geq 4 \) all of whose neighbors have degree \( = 2 \). 
The distance between two elements of \( \mathrm{Seq} \) is always even.
\item \label{lem:G_T-b}
The vertices in \( \mathrm{Seq}^- \) are the unique vertices of \( \mathbb{G}_S \) with degree \( 2 \) and all of whose neighbors are in \( \mathrm{Seq} \).
\item \label{lem:G_T-c}
The vertices of the form \( ( u , s , 0^{( \theta ( u ) )}) \) (for some \( ( u , s ) \in S \)) are the unique vertices of \( \mathbb{G}_S \) having degree (at least) \( 3 \), odd distance from the vertices in \( \mathrm{Seq} \), and belonging to a \( \Z \)-chain of \( \mathbb{G}_S \). 
All other vertices in \( \mathrm{F}'_{u,s} \) have degree at most \( 2 \).
\item \label{lem:G_T-c'}
The vertices in \( \mathrm{F} \) are the unique vertices in \( \mathbb{G}_S \) with the following three properties:
\begin{itemize}[leftmargin=1pc]
\item
they have degree \( < 4 \);
\item
all their neighbors are either in \( \mathrm{Seq} \) or have degree \( < 4 \) as well;
\item
at least one of their neighbors has degree \( < 4 \).
\end{itemize} 
\item \label{lem:G_T-d}
The vertices in \( \mathrm{U} \) are the unique vertices of \( \mathbb{G}_S \) with the following two properties:
\begin{itemize}[leftmargin=1pc]
\item
they have degree \( \neq 2 \); 
\item
they are adjacent to a vertex of degree \( \geq 4 \) which in turn has another neighbor (distinct from the described one) with degree \( \geq 4 \) as well.
\end{itemize}
They do not belong to any \( \Z \)-chain of \( \mathbb{G}_S \), and they have degree strictly smaller than \( \kappa \).
\item \label{lem:G_T-e}
The vertices in \( \widehat{\mathrm{Seq}} \) are the unique vertices in \( \mathbb{G}_S \) with the following three properties:
\begin{itemize}[leftmargin=1pc]
\item
they have degree \( 2 \);
\item
at least one of their neighbors has degree \( \geq 4 \) and belongs to \( \mathrm{U} \);
\item
all their neighbors have degree \( \geq 4 \).
\end{itemize}
\end{enumerate-(a)} 
\end{lemma} 

\begin{remark}
In many cases, the list of conditions in Lemma~\ref{lem:G_T} for characterizing the vertices belonging to the various substructures of \( \mathbb{G}_S \) is overkill. 
For example, by suppressing the condition ``all their neighbors have degree \( \geq 4 \)'' in~\ref{lem:G_T-e} we would still get a correct characterization of the vertices in \( \widehat{\mathrm{Seq}} \); however, such extra requirement is what makes the conjunction of the conditions from~\ref{lem:G_T-e} incompatible with the conjunction of the conditions from~\ref{lem:G_T-c'}. 
This ``pairwise incompatibility'' of the characterizations of the various substructures of \( \mathbb{G}_S \) will become a very useful feature when we will have to render them with corresponding \( \LL_{\kappa^+ \kappa} \)-formul\ae{} in order to prove our invariant universality result --- see Remark~\ref{rmk:disj}.
\end{remark}

\section{Completeness}\label{sec:embeddabilitygraphs}
In this section we will show that from any tree \( T \) on \( 2 \times 2 \times \kappa \) witnessing that \( R = \PROJ \body{T} \) is a \( \kappa \)-Souslin quasi-order, one can build a function \( f_T \colon \pre{\omega}{2} \to \CT_ \kappa \) that reduces \( R \) to \( \embeds ^ \kappa _ \CT \).
The function \( f_T \) is constructed in three steps:
\begin{itemize}[leftmargin=1pc]
\item
firstly the tree \( T \) is replaced by a better tree \( \tilde{T} \) called \emph{faithful representation of \( R \)} (Lemma~\ref{lem:normalform});
\item
a function \( \Sigma _T \colon \pre{\omega}{2} \to \Tr ( 2 \times \kappa ) \) is constructed, so that \( \Sigma _T \) reduces \( R \) to \( \leq_{\max}^ \kappa \), where the latter is a quasi-order on \( \Tr (2 \times \kappa) \) which is independent of \( T \) (Lemma~\ref{def:Sigma_T});
\item
finally, we compose \( \Sigma_T \) with the map \( \Tr ( 2 \times \kappa ) \to \CT_ \kappa \), \( S \mapsto \mathbb{G}_S \) defined in Section~\ref{sec:mainconstruction}, and then check that the resulting function \( f_T ( x ) = \mathbb{G}_{ \Sigma _T ( x ) } \) is the required reduction (Theorem~\ref{th:graphs}).
\end{itemize}

\subsection{Faithful representations of $\kappa$-Souslin quasi-orders} \label{subsec:tildeT}

Let \( \op{\cdot }{\cdot } \colon \On \times \On \to \On \) be the pairing function as in~\eqref{eq:Hessenberg}, and let 
\[
\boldsymbol{\varrho} ( \alpha , \beta ) \coloneqq 2 + \op{ \alpha }{\beta } . \index[symbols]{rho@\( \boldsymbol{\varrho} \), \( \bar{\boldsymbol{\varrho}} \)}
\]
Then \( \boldsymbol{\varrho} \colon \On \times \On \to \On \setminus \setLR{ 0,1 } \) is a bijection that maps \( \kappa \times \kappa \) onto \( \kappa \setminus \setLR{ 0 , 1 } \) for all cardinals \( \kappa \geq \omega \), and it satisfies \( \boldsymbol{\varrho} ( n , m ) > n , m \) for all \( n , m \in \omega \).
Thus the map 
\begin{equation} \label{eq:rho}
 \bar{\boldsymbol{\varrho} }\colon \pre{\leq \omega}{ ( \On \times \On )} \to \pre{ \leq \omega}{ ( \On \setminus \setLR{ 0 , 1} ) }
\end{equation}
defined using \( \boldsymbol{\varrho} \) coordinate-wise
\[
\bar{\boldsymbol{\varrho}} \bigl ( \seqLR{ ( \alpha _0 , \beta _0 ) , ( \alpha _1 , \beta _1), \dots } \bigr ) \coloneqq \seqLR{ \boldsymbol{\varrho} ( \alpha _0 , \beta _0 ) , \boldsymbol{\varrho} ( \alpha _1 , \beta _1), \dots }
\]
is a bijection.

\begin{definition}
For \( \kappa \) an infinite cardinal, let 
\[
 \TT_\kappa \coloneqq \Tr ( 2 \times 2 \times \kappa ) \index[symbols]{TT@\( \TT_\kappa \), \( \TT \)}
\]
and \( \TT \coloneqq \bigcup_{\kappa \in \Cn } \TT_\kappa \). 
\end{definition}

Given a \( \kappa \)-Souslin quasi-order \( R \) on \( \pre{\omega}{2} \), we would like to have a witness \( \tilde{T} \in \TT_\kappa \) of the fact that \( R \) is \( \kappa \)-Souslin which also witnesses transitivity and reflexivity of \( R \) at all finite levels in a uniform way. 
To be more specific, we want reflexivity to be witnessed by almost all elements in \( \pre{ \omega }{ \kappa } \) and that witnesses of transitivity are given by \( \bar{\boldsymbol{\varrho}} \).

\begin{definition} \label{def:faithfulrepresentation}
Let \( \kappa \) be an infinite cardinal, and \( R \) be a \( \kappa \)-Souslin quasi-order on \( \pre{\omega}{2} \). 
A tree \( \tilde{T} \in \TT_\kappa \) is called \markdef{faithful representation of \( R \)} if the following conditions hold:
\begin{enumerate-(1)}
\item \label{def:faithfulrepresentation-i}
 \( R = \PROJ \body{\tilde{T}} \);
\item \label{def:faithfulrepresentation-ii}
\( \FORALL{ u , n , s } \left ( \lh u = \lh s + 1 \IMPLIES ( u , u , n \conc s ) \in \tilde{T} \right ) \) (reflexivity);
\item \label{def:faithfulrepresentation-iii}
\( \FORALL{ u , v , w , s , t } \left ( ( u , v , s ) , ( v , w , t ) \in \tilde{T} \IMPLIES ( u , w , \bar{\boldsymbol{\varrho}} ( s , t ) ) \in \tilde{T} \right ) \) (transitivity);
 \item\label{def:faithfulrepresentation-iv}
 \( \FORALL{u , v} ( ( u , v , 0^{( \lh u )}) \in \tilde{T} \implies u = v ) \).
\end{enumerate-(1)}
\end{definition}

\begin{remark}\label{rmk:faithfulrepresentation}
Notice that the statement ``\( \tilde{T} \in \TT_\kappa \) is a faithful representation (of \( R = \PROJ \body{\tilde{T}} \))'' is absolute for transitive models \( M \) of \( \ZF \) containing \( \tilde{T} \) and such that \( \kappa \in \Cn^M \). 
Moreover, a faithful representation \( \tilde{T} \) not only ``combinatorially'' reflects the properties of the quasi-order \( R = \PROJ \body{\tilde{T}} \), but it makes the \( \kappa \)-Sousliness of \( R \) more robust because in every \( \ZF \)-model \( M \) as above the tree \( \tilde{T} \) continues to define a canonical \( \kappa \)-Souslin quasi-order \( R^M_{\tilde{T}} \coloneqq ( \PROJ \body{\tilde{T}} )^M \) which is coherent with \( R \) in the following sense:
\begin{equation} \label{eq:coherenceforquasiorders}
\FORALL{ x , y \in \pre{\omega}{2} \cap M} ( x \mathrel{R} y \IFF x \mathrel{R^M_{\tilde{T}}} y ) .
 \end{equation}
In fact, by Definition~\ref{def:faithfulrepresentation}\ref{def:faithfulrepresentation-ii}--\ref{def:faithfulrepresentation-iii} we get that for all \( x , y , z \in ( \pre{\omega}{2} )^M \):
\begin{enumerate-(i)}
\item
every \( \xi \in ( \pre{\omega}{\kappa} )^M \) with \( \xi ( 0 ) \in \omega \) witnesses \( ( x , x ) \in (\PROJ \body{ \tilde{T} })^M \);
\item
if \( \xi_0 , \xi_1 \in (\pre{\omega}{\kappa})^M \) witness, respectively, \( ( x , y ) \in (\PROJ \body{ \tilde{T} })^M \) and \( ( y , z ) \in (\PROJ \body{ \tilde{T} })^M \), then \( \bar{\boldsymbol{\varrho}} ( \xi_0 , \xi_1 ) \) witnesses \( ( x , z ) \in (\PROJ \body{ \tilde{T} })^M \).
\end{enumerate-(i)}
Therefore \( R^M_{\tilde{T}} = ( \PROJ \body{\tilde{T}})^M \) is a quasi-order. 
The coherence condition~\eqref{eq:coherenceforquasiorders} easily follows by absoluteness of the existence of infinite branches through \( \tilde{T} \).

These absoluteness properties of faithful representations \( \tilde{T} \) will be exploited in Section~\ref{subsec:absolute}.
\end{remark} 

We are now going to show that every \( \kappa \)-Souslin quasi-order \( R \) on \(\pre{\omega}{2} \) admits a faithful representation \( \tilde{T} \). 
Indeed, the following variation of the construction given in the proof of~\cite[Theorem 2.4]{Louveau:2005cq} shows how to construct such a \( \tilde{T} \) starting from an arbitrary \( T \in \TT_\kappa \) with \( R = \PROJ \body{T} \). 

For \( T \in \TT_\kappa \) let 
\[
\hat{T} \coloneqq T \cup \setofLR{ ( u , u , s ) \in \pre{< \omega}{2} \times \pre{< \omega}{2} \times \pre{< \omega}{\kappa} }{ \lh u = \lh s } .
\] 
Then \( \hat{T} \in \TT_\kappa \), and \( \PROJ \body{T} = \PROJ \body{\hat{T}} \) whenever \( \PROJ \body{T} \) is a reflexive relation on \( \pre{\omega}{2} \).
Recall from~\eqref{eq:predecessorofu} that if \( u \neq \emptyset \) then \( u^\star = u \restriction ( \lh u - 1 ) \).
Inductively define:
\begin{align*}
 \tilde{T}_0 & \coloneqq \setLR{( \emptyset , \emptyset , \emptyset )} \cup \setof{ ( u , v , 0 \conc s) }{ ( u^\star , v^\star , s ) \in \hat{T} }
\\
\tilde{T}_{n + 1} & \coloneqq \setLR{( \emptyset , \emptyset , \emptyset )} \cup \setof{ ( u , v , ( n + 1 ) \conc s ) }{ ( u , v , n \conc s ) \in \tilde{T}_n } \cup {}
\\
 & \hphantom{{}={}}{}\cup \setof{ ( u , w , ( n + 1 ) \conc \bar{\boldsymbol{\varrho}} ( s , t ) )}{ \exists v \left ( ( u , v , n \conc s ) , ( v , w , n \conc t ) \in \tilde{T}_n \right ) }.
\end{align*}
It is immediate to check that for all \( n \in \omega \)
\begin{itemize}[leftmargin=1pc]
\item
\( s \neq \emptyset \wedge ( u , v , s ) \in \tilde{T}_n \IMPLIES s ( 0 ) = n \),
\item
\( \tilde{T}_n \in \TT_\kappa \),
\item
\( s \neq \emptyset \wedge ( u , v , s ) \in \hat{T} \IMPLIES ( u , v , n \conc s^\star ) \in \tilde{T}_n \), and in particular \( \PROJ \body{\hat{T} } \subseteq \PROJ \body{\tilde{T}_n} \).
\end{itemize}
Then
\begin{equation}\label{eq:S_T}
 \tilde{T} \coloneqq \big ( \bigcup\nolimits_{n \in \omega } \tilde{T}_n \big) \setminus \setof{( u , v , 0^{( k )})}{ u , v \in \pre{ k }{2} \wedge u \neq v \wedge k > 0 } \in \TT_\kappa.
\end{equation}
By construction, if \( ( u , v , s ) \in \tilde{T} \) and \( s \neq \emptyset \) then \( s ( 0 ) \in \omega \) and \( ( u , v , s ) \in \tilde{T}_{s ( 0 )} \setminus \bigcup_{j < s ( 0 ) } \tilde{T}_j \), so
\[
\tilde{T} =\big ( \bigcup \nolimits_{n \geq 1 } \tilde{T}_n \big ) \cup ( \tilde{T}_0 \setminus \setof{( u , v , 0^{( k )})}{ u , v \in \pre{k }{2} \wedge u \neq v \wedge k > 0 } ) .
\]
 
\begin{lemma} \label{lem:normalform}
Let \( \kappa \) be an infinite cardinal, and let \( T \in \TT_\kappa \). 
If \( R = \PROJ \body{ T } \) is a quasi-order, then the tree \( \tilde{T} \in \TT_\kappa \) defined in~\eqref{eq:S_T} is a faithful representation of \( R \).
\end{lemma}

\begin{proof}
We have to show that \( \tilde{T} \) satisfies the four conditions~\ref{def:faithfulrepresentation-i}--\ref{def:faithfulrepresentation-iv} of Definition~\ref{def:faithfulrepresentation}.

The tree \( \tilde{T} \) satisfies~\ref{def:faithfulrepresentation-iv}, and it satisfies also ~\ref{def:faithfulrepresentation-ii} by definition of \( \hat{T} \). 
To prove~\ref{def:faithfulrepresentation-iii}, we can assume that \( s , t \neq \emptyset \) and set \( n \coloneqq s ( 0 ) \) and \( m \coloneqq t ( 0 ) \), so that \( ( u , v , s ) \in \tilde{T}_n \) and \( ( v , w , t ) \in \tilde{T}_m \).
We also let \( \bar{s} \coloneqq \seqofLR{ s ( i ) }{ 1 \leq i < \lh s } \) and \( \bar{t} \coloneqq \seqofLR{ t ( i ) }{ 1 \leq i < \lh t } \).
Then \( ( u , v , k \conc \bar{s} ) , ( v , w , k \conc \bar{t} ) \in \tilde{T}_k \) for all \( k \geq \max \{ n , m \} \), and in particular for \( k \coloneqq \boldsymbol{\varrho} ( n , m ) - 1 \).
Therefore
\[ 
( u , w , ( k + 1 ) \conc \bar{\boldsymbol{\varrho}} ( \bar{s} , \bar{t} ) ) = ( u , w , \boldsymbol{\varrho} ( n , m ) \conc \bar{\boldsymbol{\varrho}} ( \bar{s} , \bar{t} ) ) = ( u , w , \bar{\boldsymbol{\varrho}} ( s , t ) ) \in \tilde{T}_{\boldsymbol{\varrho} ( n , m ) } \subseteq \bigcup\nolimits_{n \in \omega} \tilde{T}_n.
\] 
Since \( \boldsymbol{\varrho} \) never takes value \( 0 \) then \( \bar{\boldsymbol{\varrho}} ( s , t ) \neq 0^{ ( \lh u ) } \), thus \( ( u , w , \bar{\boldsymbol{\varrho}} ( s , t ) ) \in \tilde{T} \), as required.

It remains to prove~\ref{def:faithfulrepresentation-i}. 
One direction is easy: \( R = \PROJ \body{ \hat{T} } \) (since \( R \) is reflexive) and \( \PROJ \body{\hat{T} } = \PROJ \body{ \tilde{T}_0 } \subseteq \PROJ \body{\tilde{T}_1} \subseteq \PROJ \body{\tilde{T}} \).
Since \( \PROJ \body{\tilde{T}} \subseteq \PROJ \body{ \bigcup_{n} \tilde{T}_n } \) it is enough to prove that \( \PROJ \body{ \bigcup_{n} \tilde{T}_n } \subseteq R \), so we may assume that \( ( x , y ) \in \PROJ \body{ \bigcup_{n} \tilde{T}_n } \) and let \( \xi \in \pre{\omega}{ \kappa} \) be a witness of this.
Then \( \forall k ( x \restriction k , y \restriction k , \xi \restriction k ) \in \tilde{T}_{\xi ( 0 )} \), so that \( ( x , y ) \in \PROJ \body{ \tilde{T}_{ \xi ( 0 )} } \). 
Therefore it is enough to prove by induction on \( n \) that 
\begin{claim}
For every \( n \in \omega \), 
\( \PROJ \body{ \tilde{T}_n } \subseteq R \), and in fact \( \PROJ \body{ \tilde{T}_n } = R \) since \( R \subseteq \PROJ \body{ \tilde{T}_0 } \subseteq \PROJ \body{ \tilde{T}_n } \).
\end{claim}

\begin{proof}[Proof of Claim]
The proof is by induction on \( n \in \omega \).
The case \( n = 0 \) is obvious as \( \PROJ \body{\tilde{T}_0 } = \PROJ \body{ \hat{T} } = R \) (by reflexivity of \( R \)), so assume \( \PROJ \body{ \tilde{T}_n } \subseteq R \), choose an arbitrary \( ( x , y ) \in \PROJ \body{ \tilde{T}_{ n + 1 } } \) and let \( \bar{\xi} \in \pre{\omega}{ \kappa} \) be such that \( ( x , y , ( n + 1 ) \conc \bar{\xi} ) \in \body{ \tilde{T}_{n + 1} } \). 
Because of the definition of \( \tilde{T}_{n + 1} \) we have to distinguish two cases:
\begin{description}
\item[Case 1]
\( \EXISTS{{} ^\infty k}[ ( x \restriction k, y \restriction k , n \conc ( \bar{\xi} \restriction ( k - 1 ) ) ) \in \tilde{T}_n ] \).
Then \( ( x , y , n \conc \bar{\xi} ) \in \body{ \tilde{T}_n } \), so that \( ( x , y ) \in \PROJ \body{ \tilde{T}_n } \subseteq R \) (by inductive hypothesis). 
\item[Case 2]
\( \FORALL{{}^\infty k} \EXISTS{ v_k} [ ( x \restriction k , v_k , n \conc \bar{\xi}_0 \restriction ( k - 1 ) ) , ( v_k , y \restriction k , n \conc \bar{\xi} _1 \restriction ( k - 1) ) \in \tilde{T}_n ] \), where \( \bar{\xi}_0 , \bar{\xi}_1 \in \pre{\omega}{ \kappa} \) are the unique elements such that \( \bar{\xi} =\bar{\boldsymbol{\varrho}} ( \bar{\xi}_0 , \bar{\xi}_1 ) \).
Now notice that the collection of all possible \( v_k \)'s as above form an infinite finitely-branching tree (there are infinitely many \( v_k \)'s because such witnesses must be distinct for different \( k > 0 \) as \( \lh v_k = k \)), so by K\"onig's lemma there is an infinite branch \( z \in \pre{ \omega}{2} \) through it which has the property that \( ( x \restriction k , z \restriction k , n \conc \bar{\xi}_0 \restriction ( k - 1 ) ) , ( z \restriction k , y \restriction k , n \conc \bar{\xi}_1 \restriction ( k - 1 ) ) \in \tilde{T}_n \) for every \( k > 0 \). 
Therefore \( \xi_0 \coloneqq n \conc \bar{\xi}_0 \) and \(\xi_1 \coloneqq n \conc \bar{\xi}_1 \) witness \( ( x , z ) , ( z , y ) \in \PROJ \body{ \tilde{T}_n } \subseteq R \), so that \( ( x , y ) \in R \) by the transitivity of \( R \). \qedhere
\end{description} 
\end{proof}
This concludes the proof of the lemma.
\end{proof}

\subsection{The quasi-order \( \leq_{\max} \) and the reduction \( \Sigma _T \)}

\begin{definition}\label{def:max}
Given trees \( S , S' \in \Tr ( 2 \times \kappa ) \), let
\[ 
S \leq^\kappa_{\max} S' \index[symbols]{87@\( \leq^\kappa_{\max} \), \( \preceq^\kappa_{\max} \)}
\] 
if and only if there is a Lipschitz (i.e.~a monotone and length-preserving) function \( \varphi \colon \pre{< \omega}{ \kappa} \to \pre{< \omega}{\kappa} \) such that for all \( u \in \pre{< \omega}{ 2} \) and \( s \in \pre{< \omega}{ \kappa} \) of the same length 
\[
 ( u , s ) \in S \implies ( u, \varphi ( s ) ) \in S'. 
\]
If \( \varphi \) can be taken to be injective, we write \( S \preceq^\kappa_{\max} S' \).
\end{definition}

As observed in~\cite{Louveau:2005cq}, if we restrict our attention to \emph{normal} trees on \( 2 \times \kappa \) (that is trees \( S \) such that \( ( u , t ) \in S \) whenever there is \( s \in \pre{ \lh t }{\kappa} \) such that \( t \) is pointwise bigger then \( s \) and \( ( u , s ) \in S \)) then the quasi-orders \( \leq^\kappa_{\max} \) and \( \preceq^\kappa_{\max} \) coincide. 
However, unlike the case \( \kappa = \omega \) considered in~\cite{Louveau:2005cq}, in the uncountable case we cannot require the tree \( \tilde{T} \) defined in~\eqref{eq:S_T} to be normal, so the two quasi-orders \( \leq^\kappa_{\max} \) and \( \preceq^\kappa_{\max} \) must be dealt with separately.

\begin{definition}\label{def:Sigma_T}
For \( T \in \TT_\kappa \), let \( \Sigma _T \colon \pre{\omega}{2} \to \Tr ( 2 \times \kappa ) \) be defined as
\[ 
 \Sigma _T ( x ) \coloneqq \setofLR{ ( u , s ) }{ ( u , x \restriction \lh u , s ) \in \tilde{T} } , \index[symbols]{SigmaT@\( \Sigma _T \)}
\]
where \( \tilde{T} \) is as in~\eqref{eq:S_T}. 
\end{definition}

Recall that \( \Code{\cdot} \colon \pre{ < \omega }{\On } \to \On \) is the bijection of~\eqref{eq:codingfinitesequenceofordinals} and that it maps \( \pre{<\omega}{\kappa} \) onto \( \kappa \).
\begin{lemma} \label{lem:max}
Let \( T \in \TT_\kappa \) be such that \( R = \PROJ \body{ T } \) is a quasi-order. 
\begin{enumerate-(a)}
\item \label{lem:max-i}
\( \Sigma _T \) simultaneously reduces \( R \) to \( \leq^\kappa_{\max} \) and \( \preceq^\kappa_{\max} \).
In particular, both \( \leq^\kappa_{\max} \) and \( \preceq^\kappa_{\max} \) are complete for \( \kappa \)-Souslin quasi-orders.
\item \label{lem:max-ii}
If \( x , y \in \pre{\omega}{2} \) are such that \( x \mathrel{R} y \), then there is a witness \( \varphi \) of \( \Sigma_T ( x ) \preceq^\kappa_{\max} \Sigma_T ( y ) \) such that \( \Code{s} \leq \Code{ \varphi ( s )} \) for every \( s \in \pre{< \omega}{\kappa} \).
\item \label{lem:max-iii}
\( \Sigma_T \) is injective.
\end{enumerate-(a)}
\end{lemma}

\begin{proof}
\ref{lem:max-i}
Since \( \preceq^\kappa_{\max} \) refines \( \leq^\kappa_{\max} \), it is enough to show that if \( R = \PROJ \body{ T } \) is a quasi-order (so that \( \tilde{T} \) is a faithful representation of \( R \) by Lemma~\ref{lem:normalform}) then for every \( x , y \in \pre{ \omega}{2} \), 
\[ 
{\Sigma_T ( x ) \leq^\kappa_{\max} \Sigma_T ( y )} \IMPLIES {x \mathrel{R}y} \IMPLIES {\Sigma_T ( x ) \preceq^\kappa_{\max} \Sigma_T ( y )}. 
\]
The proof is identical to the one of~\cite[Theorem 2.5]{Louveau:2005cq}.
Suppose first that \( \varphi \) witnesses \( \Sigma_T ( x ) \leq^\kappa_{\max} \Sigma_T ( y ) \). 
Let \( \xi \coloneqq \bigcup_{k \in \omega} \varphi ( 0^{ ( k ) }) \). 
By reflexivity of \( \tilde{T} \) (Definition~\ref{def:faithfulrepresentation}\ref{def:faithfulrepresentation-ii}), \( ( x \restriction k, 0^{( k )}) \in \Sigma_T ( x ) \), and hence \( ( x \restriction k , \varphi ( 0^{ ( k )} ) ) \in \Sigma_T ( y ) \) for all \( k \).
But this means that \( ( x , y , \xi ) \in \body{ \tilde{T} } \), and hence \( ( x , y ) \in R = \PROJ \body{\tilde{T}} \) (where for the last equality we use Definition~\ref{def:faithfulrepresentation}\ref{def:faithfulrepresentation-i}).

Assume now that \( \xi \in \pre{\omega}{ \kappa} \) witnesses \( ( x , y ) \in \PROJ \body{ \tilde{T} } = R \).
Then \( \xi ( 0 ) \in \omega \).
For \( s \in \pre{< \omega}{\kappa} \), let 
\begin{equation} \label{eq:fhi}
\varphi ( s ) \coloneqq \bar{\boldsymbol{\varrho}} ( s , \xi \restriction \lh s ) .
\end{equation} 
Since the function \( \boldsymbol{\varrho} \) used to define \( \bar{\boldsymbol{\varrho}} \) is injective, then \( \varphi \) is injective as well. 
Suppose now that \( ( u , s ) \in \Sigma_T ( x ) \) and let \( k \coloneqq \lh u = \lh s \).
Then \( s ( 0 ) \in \omega \) and \( ( u , x \restriction k , s ) \in \tilde{T} \). 
On the other hand \( ( x \restriction k, y \restriction k , \xi \restriction k ) \in \tilde{T} \), therefore \( ( u , y \restriction k , \bar{\boldsymbol{\varrho}} ( s , \xi \restriction k ) ) \in \tilde{T} \) by transitivity of \( \tilde{T} \) (Definition~\ref{def:faithfulrepresentation}\ref{def:faithfulrepresentation-iii}), so \( ( u , \varphi ( s ) ) \in \Sigma_T ( y ) \).

\smallskip

\ref{lem:max-ii} 
The map \( \varphi \) defined in~\eqref{eq:fhi} will do.

\smallskip

\ref{lem:max-iii} This follows from the fact that \( \tilde{T} \) satisfies Definition~\ref{def:faithfulrepresentation}\ref{def:faithfulrepresentation-iv}, as if \( x \neq y \) and \( k \in \omega \) is such that \( x \restriction k \neq y \restriction k \), then \( ( x \restriction k , 0^{(k)} ) \in \Sigma_T ( x ) \setminus \Sigma_T( y ) \).
\end{proof}

\subsection{Reducing \( \leq_{\max}^ \kappa \) to \( \embeds^ \kappa_\CT \)}

Now we are ready to prove the main result of this section.
Given an infinite cardinal \( \kappa \) and \( T \in \TT_\kappa \), define the function 
\begin{equation}\label{eq:f_T}
f_T \colon \pre{ \omega}{2} \to \CT_{\kappa} , \qquad x \mapsto \mathbb{G}_{ \Sigma _T ( x )}, \index[symbols]{fT@\( f_T \)}
\end{equation}
where \( \Sigma_T \) is as in Definition~\ref{def:Sigma_T} and the combinatorial tree \( \mathbb{G}_{\Sigma_T ( x ) } \) associated to \( \Sigma_T ( x ) \in \Tr ( 2 \times \kappa ) \) is defined as in Section~\ref{subsec:G_T}. 

\begin{theorem} \label{th:graphs}
Let \( \kappa \) be an infinite cardinal and \( T \in \TT_\kappa \).
If \( R = \PROJ \body{ T } \) is a quasi-order, then the map \( f_T \) defined in~\eqref{eq:f_T} is such that:
\begin{enumerate-(a)}
\item \label{th:graphs-a}
\( f_T \) reduces \( R \) to the embeddability relation \( \embeds_\CT^\kappa \);
\item \label{th:graphs-b}
\( f_T \) reduces \( = \) on \( \pre{\omega}{2} \) to the isomorphism relation \( \cong_\CT^\kappa \).	
\end{enumerate-(a)} 
In particular, \( \embeds_\CT^\kappa \) is complete for \( \kappa \)-Souslin quasi-orders on \( \pre{\omega}{2} \), i.e.~every \( \kappa \)-Souslin quasi-order on \( \pre{\omega}{ 2} \) is reducible to the embeddability relation on \( \CT_\kappa \). 
\end{theorem}

\begin{proof}
The proof of \ref{th:graphs-a} is similar to the ones of~\cite[Theorem 3.1]{Louveau:2005cq} and~\cite[Theorem 3.9]{Friedman:2011cr}, but it is simplified a little bit by our different choice of the map \( \theta \). 
Let \( x , y \in \pre{\omega}{2} \), and assume first that \( x \mathrel{R} y \). 
By Lemma~\ref{lem:max}\ref{lem:max-ii}, there is a \( \varphi \colon \pre{<\omega}{\kappa} \to \pre{<\omega}{\kappa} \) witnessing \( \Sigma_T ( x ) \preceq^\kappa_{\max} \Sigma_T ( y ) \) such that \( \Code{s} \leq \Code{ \varphi ( s )} \) for every \( s \in \pre{ < \omega}{\kappa} \). 
Set
\begin{itemize}[leftmargin=1pc]
\item
\( i ( s ) \coloneqq \varphi ( s ) \), \( i ( s^- ) \coloneqq ( \varphi ( s ) )^- \), and \( i ( \hat{s} ) \coloneqq \widehat{\varphi ( s )} \) for every \( s \in \pre{< \omega}{\kappa} \);
\item
\( i ( s , t ) \coloneqq ( \varphi ( s ) , t ) \) for every \( s \in \pre{ < \omega}{\kappa} \) and \( t \in U_{\Code{s} } \) (this definition is well given by \( \Code{s} \leq \Code{\varphi ( s )} \) and Proposition~\ref{lem:propertiesL_alpha}\ref{lem:propertiesL_alpha-b});
\item
\( i ( u , s , w ) \coloneqq ( u, \varphi ( s ) , w ) \) for every \( ( u , s ) \in \Sigma_T ( x ) \), \( ( u , s , w ) \in \mathrm{F}_{ u , s } \), and \( w \neq \emptyset \). 
(This is well-defined because \( ( u , s ) \in \Sigma_T ( x ) \IMPLIES ( u , \varphi ( s ) ) \in \Sigma_T ( y ) \) by our choice of \( \upvarphi \)).
\end{itemize}
It is easy to check that \( i \) is the desired embedding of \( f_T ( x ) \) into \( f_T ( y ) \).

Conversely, let \( j \) be an embedding of \( f_T ( x ) = \mathbb{G}_{\Sigma_T ( x )} \) into \( f_T ( y ) = \mathbb{G}_{\Sigma_T ( y )} \).
By Lemma~\ref{lem:max} it is enough to show that \( \Sigma_T ( x ) \preceq^\kappa_{\max} \Sigma_T ( y ) \). 
Since embeddings cannot decrease degrees, by Lemma~\ref{lem:G_T}\ref{lem:G_T-a} we get \( j (\mathrm{Seq}) \subseteq \mathrm{Seq} \). 
Moreover, Lemma~\ref{lem:G_T} implies that each vertex of the form \( ( u , s , 0^{ ( \theta ( u ) ) } ) \) is sent into a vertex of the same form because the properties characterizing these vertices listed in Lemma~\ref{lem:G_T}\ref{lem:G_T-c} are preserved by embeddings (and \( j (\mathrm{Seq}) \subseteq \mathrm{Seq} \)). 
In particular it follows that \( j ( G_0 \cup \mathrm{F}(\mathbb{G}_{\Sigma_T ( x )})) \subseteq G_0 \cup \mathrm{F}(\mathbb{G}_{\Sigma_T ( y )}) \).
Since \( ( \emptyset , \emptyset , 0^{ ( 3 ) } ) \) has distance \( 3 \) (which is the minimal value attained by \(\theta\)) from \( \emptyset \in \mathrm{Seq} \) by~\eqref{eq:theta0=3}, we get \( j ( \emptyset , \emptyset , 0^{ ( 3 ) } ) = ( \emptyset , \emptyset , 0^{ ( 3 ) } ) \), whence \( j ( \emptyset ) = \emptyset \). 
Arguing by induction on \( \lh s \) (and using injectivity of \( j \) and \( j ( \mathrm{Seq} ) \subseteq \mathrm{Seq} \)) we then get that \( \varphi \coloneqq j \restriction \mathrm{Seq} \colon \pre{ < \omega}{\kappa} \to \pre{<\omega}{\kappa} \) is an injective Lipschitz map. 
\begin{claim}\label{claim:theta}
For each \( ( u , s ) \in \Sigma_T ( x ) \), \( j ( u , s , 0^{ ( \theta ( u ) ) } ) = ( u , \varphi ( s ) , 0^{( \theta ( u ) )} ) \) (recall that \( \varphi(s) \coloneqq j(s) \)).
\end{claim}

\begin{proof}[Proof of Claim]
Recall that \( j \) must send \( ( u , s , 0^{ ( \theta ( u ) ) } ) \) into a vertex of the same form, so let \( ( v , t ) \in \Sigma_T ( y ) \) be such that \( j ( u , s , 0^{ ( \theta ( u ) ) }) = ( v , t , 0^{ ( \theta ( v ) ) } ) \). 
 Note that \( \theta ( u ) \) is the distance in \( f_T ( x ) \) between \( s \) and \( ( u , s , 0^{ ( \theta ( u ) ) } ) \) and \( \theta ( v ) \) is the distance in \( f_T ( y ) \) between \( t \) and \( ( v , t , 0^{( \theta ( v ) ) }) \). 
 Moreover, the path in \( f_T ( x ) \) between the nodes \( s \) and \( ( u , s , 0^{ ( \theta ( u ) ) } ) \) is mapped by \( j \) to the (unique) path in \( f_T ( x ) \) between the vertices \( \varphi ( s ) \coloneqq j ( s ) \) and \( ( v , t , 0^{ ( \theta ( t ) ) } ) \), which necessarily passes through \( t \). 
 This implies that \( \theta ( u ) - \theta ( v ) \) is the distance in \( f_T ( y ) \) between \( \varphi ( s ) \) and \( t \): but such distance is \( \leq 4 \cdot \max \setLR{ \lh \varphi ( s ) , \lh t } \) (because the latter is an upper bound for the length of the path in \( f_T(y) \) which goes from \( \varphi(s) \) to \( \emptyset \) and then back to \( t \)), and since \( \lh u = \lh s = \lh \varphi ( s ) \) and \( \lh v = \lh t \), this implies \( u = v \) by~\eqref{eq:thetadistances}, whence also \( \varphi ( s ) = t \) because such vertices have then distance \( \theta(u) - \theta(v) = 0 \).
\end{proof}

Claim~\ref{claim:theta} easily implies that \( ( u , s ) \in \Sigma_T ( x ) \IMPLIES ( u , \varphi ( s ) ) \in \Sigma_T ( y ) \), i.e.\ that \( \varphi \) witnesses \( \Sigma_T ( x ) \preceq^\kappa_{\max} \Sigma_T ( y ) \).

\smallskip

\ref{th:graphs-b}
Fix an isomorphism \( j \) between \( f_T ( x ) = \mathbb{G}_{\Sigma_T ( x )} \) and \( f_T ( y ) = \mathbb{G}_{\Sigma_T ( y )} \).
It follows from parts~\ref{lem:G_T-a},~\ref{lem:G_T-c}, and~\ref{lem:G_T-d} of Lemma~\ref{lem:G_T} that \( j ( \mathrm{Seq} ) = \mathrm{Seq} \), \( j ( G_0 ) = G_0 \), and \( j ( G_1 ) = G_1 \). 
In particular, \( j \restriction G_1 \) is an automorphism of \( \mathbb{G}_1 \).
Therefore by Lemma~\ref{lem:quasirigid}, \( j \restriction G_0 \) is the identity, and so is \( \varphi \coloneqq j \restriction \mathrm{Seq} \).
Thus the second part of the proof of~\ref{th:graphs-a} shows that \( \Sigma_T ( x ) \subseteq \Sigma_T ( y ) \) (because \( \varphi \) is now the identity map). 
Replacing \( j \) with \( j^{-1} \) in this argument, we obtain \( \Sigma_T ( y ) \subseteq \Sigma_T ( x ) \), so that \( \Sigma_T ( x ) = \Sigma_T ( y ) \). 
Since \( x \mapsto \Sigma_T ( x ) \) is injective by Lemma~\ref{lem:max}\ref{lem:max-iii}, this implies that \( f_T \) reduces \( = \) to \( \cong \), as required.
\end{proof}

\begin{remark}\label{rmk:partialembedding} 
Notice that the second half of the proof of part~\ref{th:graphs-a} shows that \( x \mathrel{R} y \) whenever there is an embedding of the subgraph of \( f_T(x) \) with domain \( G_0 \cup \mathrm{F}(\mathbb{G}_{\Sigma_T ( x )}) \) into the subgraph of \( f_T(y) \) with domain \( G_0 \cup \mathrm{F}(\mathbb{G}_{\Sigma_T ( y )}) \). 
This feature will be used in Sections~\ref{sec:alternativeapproach} and~\ref{subsubsec:changingmorphism} --- see the proofs of Theorems~\ref{th:OCT} and~\ref{thm:mainhomomorphism}.
\end{remark}

\subsection{Some absoluteness results} \label{subsec:someabsoluteness}

By closely inspecting the constructions provided in Sections~\ref{sec:mainconstruction} and~\ref{sec:embeddabilitygraphs}, one easily sees that the definition of the map \( f_T \) from~\eqref{eq:f_T} only requires the knowledge of the parameters \( T \) and \( \kappa \), and that such definition is uniform in those parameters and independent of the transitive model of \( \ZF \) we are working in. 
In fact, since the tree \( \tilde{T} \) in~\eqref{eq:S_T} is definable from \( T \) and \( \kappa \), then the function sending \( (\kappa,T) \) to the map \( f_T \) is definable (without parameters) via an \( \LST \)-formula, which moreover is absolute for transitive models of \( \ZF \). 
To be more precise, let \( M \) be an arbitrary transitive model of \( \ZF \), \( \kappa \) be a cardinal in \( M \), and \( T \in (\TT_\kappa)^M \): then, working inside \( M \), we can define the function \( f_T^M \coloneqq f_T \) as in~\eqref{eq:f_T}, which continues to be a reduction of \( R^M \coloneqq ( \PROJ \body{T})^M \) to \( (\embeds_\CT^\kappa)^M \) as long as \( R^M \) is a quasi-order in \( M \) (because Theorem~\ref{th:graphs}, which is proved in \( \ZF \), holds in \( M \)). 
With this notation, we then get the definability and absoluteness results briefly discussed below, which will be used in Section~\ref{subsec:absolute}.

\begin{fact} \label{prop:f_Tdefinable}
There is an \( \LST \)-formula \( {\sf\Psi}_{f_T} ( x_0 , x_1 , z_0 , z_1 ) \)\index[symbols]{Psi@\( {\sf \Psi}_{f_T} \)} with the following properties.
\begin{enumerate-(a)}
\item \label{prop:f_Tdefinable-a}
For every transitive model \( M \) of \( \ZF \), \( \kappa \in \Cn^{M} \) and \( T \in (\TT_\kappa)^{M} \), the formula \( {\sf\Psi}_{f_T}( x_0 , x_1 , \kappa , T ) \) defines in \( M \) (the graph of) \( f_T^M \), that is: for every \( x \in (\pre{\omega}{2})^M \) and \( X \in (\CT_\kappa)^M \)
\[ 
f_T^M ( x ) = X \IFF {M \models {\sf \Psi}_{f_T} [ x , X , \kappa , T ]}.
\] 
\item \label{prop:f_Tdefinable-b}
Let \( N \) be another transitive model of \( \ZF \) with \( \kappa \in \Cn^{N} \) and \( T \in (\TT_\kappa)^{N} \), and let \( f_T^{N} \) be the function defined in \( N \) by the formula \( {\sf\Psi}_{f_T} ( x_0 , x_1 , \kappa , T) \) --- see part~\ref{prop:f_Tdefinable-a}. 
Then
\[ 
\FORALL{x \in (\pre{\omega}{2})^{M} \cap (\pre{\omega}{2})^{N}} (f_T^{M} ( x ) = f_T^{N} ( x ) ),
 \] 
i.e.\ \( f_T^{M} \) and \( f_T^{N} \) coincide on the common part of their domain.
\end{enumerate-(a)}
\end{fact}

Indeed, \( {\sf\Psi}_{f_T} ( x_0 , x_1 , z_0 , z_1 ) \) is the formalization in the language of set theory of the construction of the combinatorial tree \( x_1 \coloneqq f_T ( x_0 ) = \mathbb{G}_{\Sigma_T ( x_0 )} \in \CT_\kappa \) starting from the parameters \( x_0 \in \pre{\omega}{2} \), \( z_0 \coloneqq \kappa \), and \( z_1 \coloneqq T \).
 We leave to the reader to check that such formalization is indeed possible. 
For part~\ref{prop:f_Tdefinable-b}, notice that for every \( x \in ( \pre{\omega}{2} )^{M} \cap (\pre{\omega}{2} )^{N} \) the two combinatorial trees \( f_T^{M} ( x ) = ( \mathbb{G}_{\Sigma_T ( x ) })^M \) and \( f_T^{N} ( x ) = (\mathbb{G}_{\Sigma_T ( x ) })^{N} \) must coincide because they are explicitly computed in \( \ZF \) using just \( x \), \( \kappa \) and \( T \) as parameters and all the bijections \( \Code{\cdot} \), \( \op{ \cdot}{\cdot } \), and \( e_u \) involved in their coding as structures on \( \kappa \) are absolute between transitive models of \( \ZF \).

\begin{lemma} \label{lem:canonicalembedding}
Given an infinite cardinal \( \kappa \) and a tree \( T \in \TT_\kappa \) such that \( R = \PROJ \body{ T } \) is a quasi-order, let \( f_T \) be the map defined in~\eqref{eq:f_T}. 
Then for arbitrary \( x , y \in \pre{\omega}{2} \), if \( f_T ( x ) \embeds f_T ( y ) \) then there is a (canonical) witness of this fact which is explicitly \( \LST \)-definable \( \ZF \) using only \( x \), \( y \), \( \kappa \), and \( T \) as parameters.
\end{lemma}

\begin{proof}
If \( f_T ( x ) \embeds f_T ( y ) \), then \( x \mathrel{R} y \) by Theorem~\ref{th:graphs}\ref{th:graphs-a}. 
Pick the leftmost branch \( b_{x , y} \) such that \( ( x , y , b_{x , y}) \in \body{ \tilde{ T }} \) (where \( \tilde{T} \) is the faithful representation of \( R \) constructed from \( T \) as in~\eqref{eq:S_T}), and apply~\eqref{eq:fhi} to \( \xi \coloneqq b_{x , y} \) to get a (canonical) \( \varphi \colon \pre{<\omega}{\kappa} \to \pre{< \omega}{\kappa} \) witnessing \( \Sigma_T ( x ) \preceq^\kappa_{\max} \Sigma_T ( y ) \).
Then use this \(\varphi\) to define the desired canonical embedding \( i \) of \( f_T ( x ) \) into \( f_T ( y ) \) as in the first part of the proof of Theorem~\ref{th:graphs}\ref{th:graphs-a}.
\end{proof}

Given an infinite cardinal \( \kappa \), the embeddability relation on \( \Mod^\kappa_\LL \), being defined by a \( \Sigma_1 \) \( \LST \)-formula, is upward absolute but in general not downward absolute: for example, for suitable \( X , Y \in \CT_\kappa \) one may be able to add by forcing an embedding of \( X \) into \( Y \) (so that \( X \embeds Y \) holds in a suitable generic extension) even if \( X \) does not embed into \( Y \) in the ground model.
In contrast, we are now going to show that the embeddability relation on the range of an \( f_T \) as in~\eqref{eq:f_T} is always absolute for transitive models of \( \ZF \) containing \( \kappa \) and \( T \).

\begin{proposition} \label{prop:embeddabilityabsolute}
Let \( M_0 \), \( M_1 \) be transitive models of \( \ZF \), and let \( \kappa \) and \( T \) be such that \( \kappa \in \Cn^{M_i} \), \( T \in ( \TT_\kappa )^{M_i} \), and \( R^{M_i} \coloneqq ( \PROJ\body{T})^{M_i} \) is a quasi-order in \( M_i \) (for \( i = 0 , 1 \)). 
Recall that \( f_T^{M_0} ( z ) = f_T^{M_1} ( z ) \) for all \( z \in (\pre{\omega}{2})^{M_0} \cap (\pre{\omega}{2})^{M_1} \) by Fact~\ref{prop:f_Tdefinable}\ref{prop:f_Tdefinable-b}, so that we can unambiguously set \( f_T ( z ) \coloneqq f_T^{M_0} ( z ) = f_T^{M_1} ( z ) \) for all such \( z \).
Then for every \( x , y \in (\pre{\omega}{2})^{M_0} \cap (\pre{\omega}{2})^{M_1} \)
\[ 
M_0 \models f_T ( x ) \embeds f_T ( y ) \IFF M_1 \models f_T ( x ) \embeds f_T ( y ).
 \] 
\end{proposition}

\begin{proof}
By Theorem~\ref{th:graphs}\ref{th:graphs-a} (which holds both in \( M_0 \) and \( M_1 \)) and absoluteness of existence of infinite branches through descriptive set-theoretic trees on \( \kappa \), we have that 
\begin{align*}
M_0 \models f_T ( x ) \embeds f_T ( y ) & \IFF M_0 \models \EXISTS{\xi \in \pre{\omega}{\kappa}} ( ( x , y , \xi ) \in \body{ T } ) 
\\
& \IFF M_1 \models \EXISTS{\xi \in \pre{\omega}{\kappa} } ( ( x , y , \xi ) \in \body{T} ) 
\\
& \IFF M_1 \models f_T ( x ) \embeds f_T ( y ). 
\qedhere
\end{align*} 
\end{proof}
\section{Invariant universality} \label{sec:invariantlyuniversal}
As mentioned in the introduction, the isomorphism relation \( {\cong}^\omega_\upsigma \) on the set \( \Mod^\omega_\upsigma \) of countable models of an \( \LL_{\omega_1 \omega} \)-sentence \( \upsigma \) is an analytic equivalence relation, and these equivalence relations have been extensively studied in the literature (see~\cite{Vaught:1974kl,Becker:1996uq} and the references therein).
In~\cite{Louveau:2005cq} the analytic quasi-order \( {\embeds }^\omega_\upsigma \) of embeddability between countable models of \( \upsigma \) (where e.g.\ \( \upsigma \coloneqq \upsigma_\CT \) is the \( \LL_{\omega_1 \omega} \)-sentence axiomatizing combinatorial trees from page~\pageref{eq:sigmaCT}) has been shown to be \( \leq_B \)-complete for the class of analytic quasi-orders.
In~\cite{Friedman:2011cr} a strengthening of Borel-completeness (called invariant universality in~\cite{Camerlo:2012kx}) which bears on both \( \cong \) and \( \embeds \) has been introduced.
Here this notion is generalized to arbitrary infinite cardinals.

\begin{definition} \label{def:invuniversal}
Let \( \mathcal{C} \) be a class of quasi-orders, \( \LL \) be a finite relational language, and \( \kappa \) be an infinite cardinal.
The embeddability relation \( {\embeds} ^\kappa_\LL \) is \markdef{invariantly universal for \( \mathcal{C} \)}\index[concepts]{invariant universality} if for every \( R \in \mathcal{C} \) there is an \( \LL_{\kappa^+ \kappa} \)-sentence \( \upsigma \) such that \( R \sim {\embeds} ^\kappa_\upsigma \).
\end{definition} 

A localized version of invariant universality can be defined in a similar way.

\begin{definition} \label{def:invuniversallocal}
Let \( \mathcal{C}\), \( \LL \), and \( \kappa \) be as in Definition~\ref{def:invuniversal}. 
Given an \( \LL_{\kappa^+ \kappa} \)-sentence \( \uptau \), the embeddability relation \( {\embeds} ^\kappa_{\uptau} \) is \markdef{invariantly universal for \( \mathcal{C} \)} if for every \( R \in \mathcal{C} \) there is an \( \LL_{\kappa^+ \kappa} \)-sentence \( \upsigma \) such that \( \Mod_ \upsigma^ \kappa \subseteq \Mod_ \uptau ^ \kappa \) and \( R \sim {\embeds} ^\kappa_\upsigma \).
\end{definition}

As for the case of (\( \leq_* \)-)completeness, when in Definitions~\ref{def:invuniversal} and~\ref{def:invuniversallocal} the reducibility \( \leq \) is replaced by one of its restricted forms \( \leq_* \) we speak of \markdef{\( \leq_* \)-invariant universality}.

\begin{figure}
\[
\begin{tikzpicture}[scale=1.5]
\node at (0,4) [above] {\( \pre{ \kappa }{ 2 } \)};
\draw (2,0) rectangle (6,4);
\node at (4,4) [above] {\( \CT_ \kappa \)};
\draw[rounded corners=9pt] (2.5,0.5) rectangle (5.5,3.5);
\draw[->,>=stealth] (4,2) .. controls (2, 1.1) and (1, 1.6) .. node[below]{\( h_T \)}(10pt,2);
\filldraw[fill=white,very thick] (4,2) circle (1.2);
\filldraw[fill=gray!30] (4,2) circle (0.7);
\node at (4,2) {\( \ran (f_T) \)};
\draw[->,>=stealth] (0,2) .. controls (1, 2.8) and (2, 2.4) .. node[above]{\( f_T \)} (3.3,2);
\filldraw[draw=black, fill=gray!30] (0,2) ellipse [x radius=10pt, y radius=2cm] ;
\draw[rounded corners=9pt] (8,0) rectangle (9,4);
\node at (8.5,4) [above] {\( \pre{ ( \pre{ < \omega }{2} \times \kappa ) }{ 2 } \)};
\draw[very thin] (3,1.366) -- (3,-0.5);
\draw[very thin] (5,0.5) -- (5,-0.5);
\node at (3,-0.5) [below]{\( \Mod^ \kappa _{ \upsigma_T} \)};
\node at (5,-0.5) [below]{\( \Mod ^ \kappa _\Uppsi \)};
\draw[->,>=stealth] (5.5,2) .. controls (6.25, 2.2) and (7.5, 2.2) .. node[above] {\( g \)}(8,2);
\end{tikzpicture}
\]
\caption{The reductions used in Section~\ref{sec:invariantlyuniversal}.
The set \( \Mod^ \kappa _{ \upsigma_T} \) is the saturation of \( \ran (f_T) \), where \( f_T \) is the map defined in~\eqref{eq:f_T}.}\label{fig:reductions}
\end{figure}
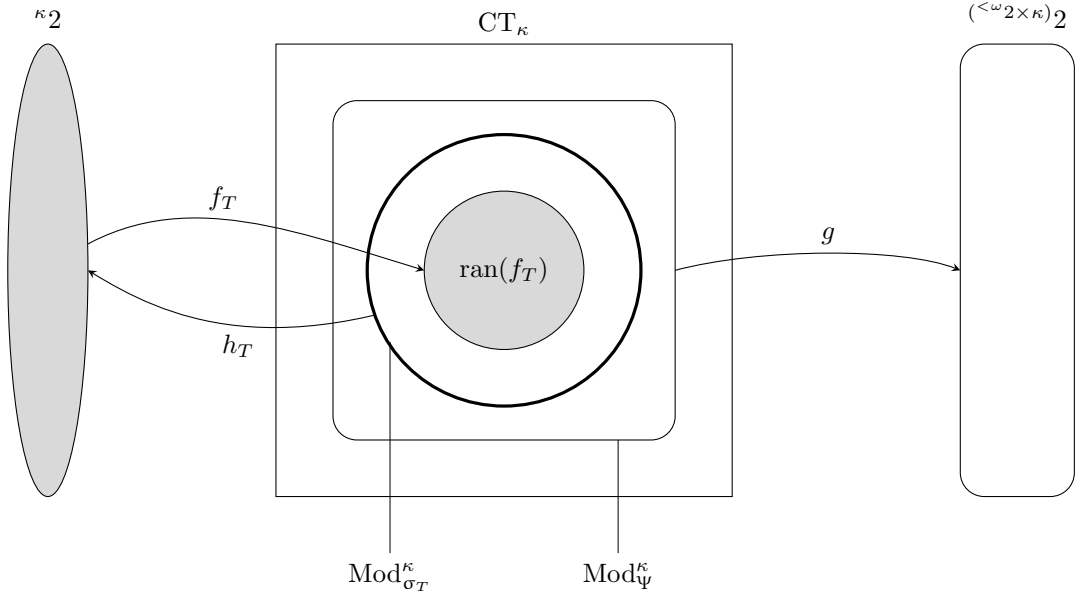

In order to establish the invariant universality of \( \embeds^\kappa_\CT \) for \( \kappa \)-Souslin quasi-orders on \( \pre{\omega}{2} \), in this section we will have to construct some infinitary sentences, denoted by \( \Uppsi \) and \( \upsigma_T \), and some maps between their sets of models of size \( \kappa \): the reader is advised to refer to Figure~\ref{fig:reductions} in order to keep track of these functions.

\subsection{An \( \LL_{\kappa^+ \kappa} \)-sentence \( \Uppsi \) describing the structures \( \mathbb{G}_S \).} \label{subsec:Uppsi}
\emph{Henceforth we fix an uncountable cardinal \( \kappa \).} 
We will now begin the definition of a sequence of nine \( \LL_{\kappa^+ \kappa} \)-sentences \( \Upphi_0, \dotsc, \Upphi_8 \) which will be crucial for the proof of our main result. 
These sentences try to describe with some accuracy the common properties of the structures of the form \( \mathbb{G}_S \) for \( S \in \Tr(2 \times \kappa) \) defined in Section~\ref{sec:mainconstruction}.
To help in understanding the intended meaning of such sentences, we will freely use the following two conventions (besides the ones already explained in Section~\ref{subsubsec:syntax}):
\begin{itemize}[leftmargin=1pc]
\item
we will use metavariables \( x , y , z , x_\alpha , y_\alpha , z_\alpha \) (possibly with various decorations or different subscripts) instead of the \( \V_\alpha \)'s;
\item
we will consider (infinitary) conjunctions and disjunctions over \emph{sets} of formul\ae{} of \emph{size \( \leq \kappa \)} (instead of conjunctions and disjunctions over \emph{sequences of length \( < \kappa^+ \)}, as Definition~\ref{def:L_kappalambda} would officially require), as long as it is clear that such sets can be well-ordered in a canonical way in \( \ZF \) (usually by means of the coding functions \( \Code{\cdot} \), \( \op {\cdot}{\cdot} \), and \( e_u \) for \( u \in \pre{<\omega}{2} \) from~\eqref{eq:codingfinitesequenceofordinals},~\eqref{eq:Hessenberg} and~\eqref{eq:e_u}, respectively, which are absolute for transitive models of \( \ZF \)). 
\end{itemize}
This means that, formally, the \( \LL_{\kappa^+ \kappa} \)-sentence \( \Upphi_i \) is obtained from the displayed one by substituting in the natural way each metavariable with a corresponding variable in the official list \( \seqofLR{ \V_\alpha }{ \alpha < \kappa } \), and by well-ordering in a canonical way all the sets of subformul\ae{} to which an (infinitary) conjunction or disjunction is applied. 
We will explicitly perform this formalization just for the first few formul\ae{}, leaving to the reader all other cases.

Given a variable \( x \) and \( 0 \neq n \in \omega \), let \( \mathsf{d}_{<n} ( x ) \) denote the \( \LL_{\omega \omega} \)-formula
\begin{equation}\tag{$\mathsf{d}_{ < n } ( x ) $}
\FORALL{ \seqof{ x_i }{ i < n }} \Bigl ( \bigwedge_{i < n} x_i \edge x \IMPLIES \bigvee_{i < j < n} x_i \equals x_j \Bigr )
\end{equation}
(where if \( n = 1 \) we agree that \( \bigvee_{i < j < n} x_i \equals x_j \) is any inconsistent sentence), and abbreviate by \( \mathsf{d}_{\geq n} ( x ) \), \( \mathsf{d}_{=n} ( x ) \), and \( \mathsf{d}_{\neq n} ( x ) \) the \( \LL_{\omega \omega} \)-formul\ae{} \( \neg ( \mathsf{d}_{<n} ( x ) ) \), \( \mathsf{d}_{< n+1} ( x ) \wedge \mathsf{d}_{\geq n} ( x ) \), and \( \neg (\mathsf{d}_{=n} ( x ) ) \), respectively. 
If \( a \) is a vertex of a graph \( X \), then \( X \models \mathsf{d}_{<n} [ a ] \) if and only if \(a \) has degree \( < n \) in \( X \).
To completely formalize the formula \( \mathsf{d}_{ < n } ( x ) \), one should first fix \( \alpha < \kappa \) and enumeration \( \seqofLR{ c_m }{ m < \frac{n(n-1)}{2} } \) of all pairs \( ( i , j ) \in \omega^2 \) with \( i < j < n \), and then define \( \mathsf{d}_{<n} ( \V_\alpha ) \) as
\[
\FORALL{ \seqof{ \V_{\alpha + i + 1} }{ i < n } } \Bigl ( \bigwedge_{i < n} \V_{\alpha+i+1} \edge \V_\alpha \IMPLIES \bigvee_{m < \frac{ n ( n - 1 ) }{2}} \V_{\alpha + i_m + 1 } \equals \V_{\alpha + j_m + 1} \Bigr ),
\]
where \( ( i_m , j_m ) \in \omega^2 \) is the unique pair such that \( c_m = ( i_m , j_m ) \). 
A similar formalization may easily be obtained for the following auxiliary \( \LL_{\omega \omega} \)-formul\ae{}:
\begin{equation}\tag{$\mathsf{Seq} ( x ) $}
\mathsf{d}_{\geq 4} ( x ) \wedge \FORALL{y} ( x \edge y \IMPLIES \mathsf{d}_{=2} ( x ) ) ;
\end{equation}

\begin{equation}\tag{$\mathsf{Seq}^- ( x )$}
\mathsf{d}_{=2} ( x ) \wedge \FORALL{y} ( x \edge y \IMPLIES \mathsf{Seq} ( y ) ) ;
\end{equation}

\begin{equation}\tag{$\mathsf{F} ( x )$}
\mathsf{d}_{<4} ( x ) \wedge \FORALL{y} [ x \edge y \IMPLIES ({\mathsf{d}_{< 4} ( y ) } \vee {\mathsf{Seq} ( y ) } ) ] \wedge \EXISTS{z} [ {x \edge z} \wedge {\mathsf{d}_{< 4} ( z ) }] .
\end{equation}

\begin{equation}\tag{$\mathsf{U} ( x )$}
\mathsf{d}_{\neq 2} ( x ) \wedge \EXISTS{y,z} [{z \nequals x} \wedge {x \edge y } \wedge {y \edge z} \wedge { \mathsf{d}_{\geq 4} ( y ) } \wedge {\mathsf{d}_{\geq 4}(z)}] ;
\end{equation}

\begin{equation}\tag{$\widehat{\mathsf{Seq}} ( x )$}
\mathsf{d}_{=2} ( x ) \wedge \EXISTS{y} [ \mathsf{U} ( y ) \wedge {x \edge y} \wedge {\mathsf{d}_{\geq 4} ( y ) } ] \wedge \FORALL{z}[ {z \edge x} \IMPLIES \mathsf{d}_{\geq 4}(z)] ;
\end{equation}

\begin{remark} \label{rmk:substructures}
It is not hard to check that for every \( S \in \Tr ( 2 \times \kappa ) \) and every vertex \( a \in \mathbb{G}_S \) (where \( \mathbb{G}_S \) is defined as in Section~\ref{subsec:G_T}) we have
\[ 
\mathbb{G}_S \models \mathsf{Seq} [ a ] \IFF a \in \mathrm{Seq}(\mathbb{G}_S),
 \] 
and analogous results hold for the other formul\ae{} \( \mathsf{Seq}^- ( x ) \), \( \mathsf{U} ( x ) \), \( \widehat{\mathsf{Seq}} ( x ) \), and \( \mathsf{F} ( x ) \) and the corresponding substructures \( \mathrm{Seq}^-( \mathbb{G}_S ) \), \( \mathrm{U} ( \mathbb{G}_S ) \), \( \widehat{\mathrm{Seq}}( \mathbb{G}_S ) \), and \( \mathrm{F} ( \mathbb{G}_S ) \) of \( \mathbb{G}_S \) defined in~\eqref{eq:substructuresofG_T-a}--\eqref{eq:substructuresofG_T-g}.
\end{remark}

We are now ready to introduce the first \( \LL_{\kappa^+ \kappa} \)-sentence \( \Upphi_0 \) (which is actually an \( \LL_{\omega_1 \omega} \)-sentence):
\begin{equation}\tag{$\Upphi_0$}
\upvarphi_\CT \wedge \FORALL{x} [\mathsf{Seq} ( x ) \vee \mathsf{Seq}^- ( x ) \vee \widehat{\mathsf{Seq}} ( x ) \vee \mathsf{U} ( x ) \vee \mathsf{F} ( x ) ] .
\end{equation}

\begin{remark} \label{rmk:disj}
Notice that the formul\ae{} appearing in the disjunction of \( \Upphi_0 \) are mutually exclusive: each element of an arbitrary \( \LL \)-structure \( X \) satisfying \( \Upphi_0 \) can realize at most one of \( \mathsf{Seq} ( x ) \), \( \mathsf{Seq}^- ( x ) \), \( \widehat{\mathsf{Seq}} ( x ) \), \( \mathsf{U} ( x ) \), or \( \mathsf{F} ( x )\).
\end{remark}

In any combinatorial tree of the form \( \mathbb{G}_S \) (for \( S \in \Tr(2 \times \kappa) \)) one has that for every vertex \( a \) in \( \mathrm{U} ( \mathbb{G}_S ) \) there is a unique \( b \in \widehat{\mathrm{Seq}}( \mathbb{G}_S ) \) such that \( a \) is connected to \( b \) by a (finite) chain all of whose intermediate points are in \( \mathrm{U}( \mathbb{G}_S ) \) as well. 
This property of \( a \) and \( b \) is rendered by the following \( \LL_{\omega_1 \omega} \)-formula:
\begin{equation}\tag{$\mathsf{root} ( y , x ) $}
\mathsf{U} ( y ) \wedge \widehat{\mathsf{Seq}} ( x ) \wedge \bigvee_{n < \omega} \EXISTS{ \seqofLR{ x_i }{ i \leq n } } \Bigl [ {\bigwedge_{0 < i \leq n } \mathsf{U} ( x_ i )} \wedge {x \equals x_0} \wedge {x_n \edge y} \wedge {\bigwedge_{ i < n } x_i \edge x_{ i + 1 }} \Bigr ] .
\end{equation}
This allows us to write the \( \LL_{\omega_1 \omega} \)-formula \( \Upphi_1 \), which expresses the above mentioned property of the structures \( \mathbb{G}_S \) (the symbol \( \exists ! \) denotes the quantifier ``there exists a unique'', which may be expressed using the other quantifiers and connectives in the usual way):
\begin{equation}\tag{$\Upphi_1$}
\FORALL{y} [ \mathsf{U} ( y ) \IMPLIES \EXISTSONE{ x} (\mathsf{root} ( y , x ) ) ] .
\end{equation}

Let now \( X \) be an arbitrary \( \LL \)-structure \emph{ of size \( < \kappa \)}, and \( i \colon X \to \kappa \) be any injection. 
We denote by 
\[ 
\uptau^i_{\mathsf{qf}} ( X ) ( \seqofLR{ \V_\alpha }{ \alpha \in \ran ( i ) } )
 \] 
the \markdef{quantifier-free type of \( X \) (induced by \( i \))}, i.e.\ the \( \LL^0_{\kappa \kappa} \)-formula
\[ 
\bigwedge_{\substack{x,y \in X \\ x \neq y}} ( \V_{ i ( x ) } \nequals \V_{ i ( y ) } ) \wedge \bigwedge_{\substack{ x , y \in X \\ x \edge^X y} } ( \V_{ i ( x ) } \edge \V_{ i ( y ) } ) \wedge \bigwedge_{\substack{ x , y \in X \\ \neg ( x \edge^X y ) } } \neg ( \V_{ i ( x ) } \edge \V_{ i ( y ) } ) .
 \] 
To completely formalize this sentence, one of course need to well-order the infinitary conjunctions above using the given injection \( i \). 
Notice that if \( Y \) is an \( \LL \)-structure and \( \seqofLR{ a_\alpha }{ \alpha \in \ran ( i ) }\), \( \seqofLR{ b_\alpha }{ \alpha \in \ran ( i ) } \) are two sequences of elements of \( Y \) such that both \( Y \models \uptau^i_{\mathsf{qf}}( X ) [ \seqofLR{ a_\alpha }{ \alpha \in \ran ( i ) } ] \) and \( Y \models \uptau^i_{\mathsf{qf}} ( X ) [ \seqofLR{ b_\alpha }{ \alpha \in \ran ( i ) } ] \), then \( Y \restriction \setofLR{ a_\alpha }{ \alpha \in \ran ( i ) } \) and \( Y \restriction \setofLR{ b_\alpha }{ \alpha \in \ran(i) } \) are isomorphic via the map \( a_\alpha \mapsto b_\alpha \) (in fact, they are both isomorphic to \( X \)). 
In this paper, the previous procedure will be applied only to structures \( X \) which are \emph{canonically well-orderable} in \( \ZF \) using the coding maps \( \Code{\cdot} \), \( \op {\cdot}{\cdot} \), and \( e_u \) for \( u \in \pre{<\omega}{2} \) in the obvious way --- in fact, the domains of these structures will in general be subsets of size \( < \kappa \) of \( \pre{< \omega}{\kappa} \). 
In all such cases we will thus have a canonical injection \( i = i_X \colon X \to \kappa\) (namely, the one induced by the coding map \( \Code{\cdot} \)): in order to simplify the notation, we will then safely drop the reference to such \( i \), replace variables with metavariables, and call the resulting expression \markdef{qf-type of \( X \)}. 
In the formul\ae{} below we will denote the qf-type of such an \( \LL \)-structure \( X \) simply by
\[ 
\uptau_{\mathsf{qf}} ( X ) ( \seqofLR{ x_i }{ i \in X } ) .
 \] 

The next \( \LL_{\kappa^+ \kappa} \)-sentences \( \Upphi_2 \), \( \Upphi_3 \), and \( \Upphi_4 \) will complete the description of the substructure \( \mathbb{G}_1 \) of any \( \mathbb{G}_S \) (see Lemma~\ref{lem:descriptionofG_1}). 
To this aim, for each \( s \in \pre{<\omega}{\kappa} \) we first introduce the following auxiliary \( \LL_{\kappa^+ \kappa} \)-formul\ae{}:
\begin{multline}\tag{$\widehat{\mathsf{Seq}}_s ( x ) $}
\widehat{\mathsf{Seq}} ( x ) \AND \EXISTS{ \seqof{ x_i }{ i \in U_{\Code{s}} } } \biggl [ \bigwedge_{i \in U_{ \Code{s} }} \mathsf{root} ( x_i , x ) \AND \FORALL{y} \bigl ( \mathsf{root} ( y , x ) \IMPLIES \bigvee_{i \in U_{ \Code{s} }} y \equals x_i \bigr ) 
\\ 
 \AND {x \edge x_\emptyset } \AND \bigl ( {\bigwedge_{\emptyset \neq i \in U_{\Code{s}}} \neg ( x \edge x_i)} \bigr ) \AND \uptau_{\mathsf{qf}} ( U_{\Code{s}} ) ( \seqofLR{ x_i }{ i \in U_{\Code{s}} } ) \biggr ] ;
\end{multline}

\begin{equation}\tag{$\mathsf{Seq}_s ( x )$}
\mathsf{Seq} ( x ) \wedge \EXISTS{y} ( \widehat{\mathsf{Seq}}_s ( y ) \wedge {x \edge y} ) ;
\end{equation}

\begin{equation}\tag{$\mathsf{Seq}^-_s ( x )$}
\mathsf{Seq}^- ( x ) \wedge \EXISTS{ y , z} ( \mathsf{Seq}_{s^\star} ( y ) \wedge \mathsf{Seq}_s ( z ) \wedge {y \edge x} \wedge { x \edge z} ) .
\end{equation}
Recall from~\eqref{eq:predecessorofu} that for \( \emptyset \neq s \in \pre{< \omega}{\kappa}\) we set \( s^\star = s \restriction (\lh s -1) \), so \( \mathsf{Seq}^-_s ( x ) \) is defined only for \( s \neq \emptyset \). 

Notice that if \( S \in \Tr(2 \times \kappa) \) and \( a \in \mathbb{G}_S \), then 
\[ 
\mathbb{G}_S \models \widehat{\mathsf{Seq}}_s [ a ] \IFF a = \widehat{s} .
\]
Similarly,
\[ 
\mathbb{G}_S \models \mathsf{Seq}_s [ a ] \IFF a = s \qquad \text{and} \qquad \mathbb{G}_S \models \mathsf{Seq}^-_s [ a ] \IFF a = s^- .
\]
Now let
\begin{equation}\tag{$\Upphi_2$}
\FORALL{x} \Bigl [ \widehat{\mathsf{Seq}} ( x ) \IMPLIES \bigvee_{s \in \pre{<\omega}{\kappa}} \widehat{\mathsf{Seq}}_s ( x ) \Bigr ] \wedge \bigwedge_{s \in \pre{<\omega}{\kappa}} \EXISTSONE{x} (\widehat{\mathsf{Seq}}_s ( x ) ) ;
\end{equation}

\begin{equation}\tag{$\Upphi_3$}
\FORALL{x} [\mathsf{Seq} ( x ) \IMPLIES \EXISTSONE{y} (\widehat{\mathsf{Seq}} ( y ) \wedge { x \edge y} )] \wedge \bigwedge_{s \in \pre{<\omega}{\kappa}} \EXISTSONE{x} \mathsf{Seq}_s ( x ) ;
\end{equation}

\begin{equation}\tag{$\Upphi_4$}
\FORALL{x} \Bigl [ \mathsf{Seq}^- ( x ) \IMPLIES \bigvee_{\emptyset \neq s \in \pre{<\omega}{\kappa}} \mathsf{Seq}^-_s ( x ) \Bigr ] \wedge \bigwedge_{\emptyset \neq s \in \pre{<\omega}{\kappa}} \EXISTSONE{x} (\mathsf{Seq}^-_s ( x ) ) .
\end{equation}

It is not hard to see that \( \mathbb{G}_S \models \bigwedge_{i \leq 4} \Upphi_i \) for every \( S \in \Tr(2 \times \omega) \). 
We are now going to show in Lemma~\ref{lem:descriptionofG_1} that every model of \( \bigwedge_{i \leq 4} \Upphi_i \) contains a substructure (canonically) isomorphic to \( \mathbb{G}_1 \). 
Let us first fix some notation. 
For \( X \in \Mod^\kappa_\LL \), set
\begin{subequations}
\begin{align}
\mathrm{Seq} ( X ) &\coloneqq \setofLR{a \in X}{X \models \mathsf{Seq} [ a ]}
 \\
\mathrm{Seq}^- ( X ) &\coloneqq \setofLR{a \in X}{X \models \mathsf{Seq}^- [ a ]} 
\\
\widehat{\mathrm{Seq}} ( X ) & \coloneqq \setofLR{a \in X}{X \models \widehat{\mathsf{Seq}} [ a ]} 
\\
\mathrm{U} ( X ) & \coloneqq \setofLR{a \in X}{X \models \mathsf{U} [ a ]} 
\\
\mathrm{F} ( X ) & \coloneqq \setofLR{a \in X}{X \models \mathsf{F} [ a ]} 
\\
\mathbb{G}_0 ( X ) & \coloneqq \mathrm{Seq} ( X ) \cup \mathrm{Seq}^- ( X ) 
\\
\mathbb{G}_1 ( X ) & \coloneqq \mathbb{G}_0 ( X ) \cup \widehat{\mathrm{Seq}} ( X ) \cup \mathrm{U} ( x ) . \label{eq:substructuresofX-g}
\end{align} 
\end{subequations}
(Notice that when \( X = \mathbb{G}_S \) for some \( S \in \Tr(2 \times \kappa) \) this notation is coherent with the one established in~\eqref{eq:substructuresofG_T-a}--\eqref{eq:substructuresofG_T-c}, \eqref{eq:substructuresofG_T-e} and~\eqref{eq:substructuresofG_T-g}, and that the unions in the definition of \( \mathbb{G}_0 ( X ) \) and \( \mathbb{G}_1 ( X ) \) are necessarily disjoint by Remark~\ref{rmk:disj}.)

\begin{lemma} \label{lem:descriptionofG_1}
For every \( X \in \CT_\kappa \), if \( X \models \bigwedge_{i \leq 4} \Upphi_i \), then \( \mathbb{G}_1 ( X ) \cong \mathbb{G}_1 \).
Moreover this is witnessed by a specific \( \iota_X \colon \mathbb{G}_1 ( X ) \to \mathbb{G}_1 \), which we call the canonical isomorphism.
\end{lemma}

\begin{proof}
Recall the \( \LL \)-structures \( \mathrm{Seq} \), \( \mathrm{Seq}^- \), \( \widehat{\mathrm{Seq}} \), \( \mathrm{U}_s \), and \( \mathrm{U} \) defined in Remark~\ref{rmk:substructuresofG_T}. 
We will canonically define some partial isomorphisms 
\begin{align*} 
\iota_{\widehat{\mathrm{Seq}}} \colon & \widehat{\mathrm{Seq}} ( X ) \to \widehat{\mathrm{Seq}} 
\\
\iota_{\mathrm{U}} \colon & \mathrm{U} ( X ) \to \mathrm{U} 
\\
\iota_{\mathrm{Seq}} \colon & \mathrm{Seq} ( X ) \to \mathrm{Seq} 
\\
\iota_{\mathrm{Seq}^-} \colon & \mathrm{Seq}^- ( X ) \to \mathrm{Seq}^- .
\end{align*}
and then show that the union \( \iota_X \) of these maps is the desired canonical isomorphism.

By \( X \models \Upphi_2 \), we get that there is a bijection \( j \colon \widehat{\mathrm{Seq}} ( X ) \to \pre{< \omega}{\kappa} \) (namely, the map sending \( a \in \widehat{\mathrm{Seq}} ( X ) \) to the unique \( s \in \pre{< \omega}{\kappa} \) such that \( X \models \widehat{\mathsf{Seq}}_s [a] \)), so we can define the bijection \( \iota_{\widehat{\mathrm{Seq}}} \colon \widehat{\mathrm{Seq}} ( X ) \to \widehat{\mathrm{Seq}} \colon a \mapsto \widehat{j ( a )} \). 

For each \( s \in \pre{<\omega}{\kappa} \), let 
\[
 \mathrm{U}_s ( X ) \coloneqq \setofLR{ a \in X}{X \models \mathsf{root} [ a , j^{-1} ( s ) ] } .
\] 
(Notice that this notation is again coherent with~\eqref{eq:substructuresofG_T-d} when \( X = \mathbb{G}_S \) for some \( S \in \Tr(2 \times \kappa) \).)
Then \( \setofLR{ \mathrm{U}_s ( X ) }{s \in \pre{< \omega}{\kappa}} \) is a partition of \( \mathrm{U} ( X ) \) by \( X \models \Upphi_1 \). 
Moreover, using \( X \models \widehat{\mathsf{Seq}}_s [ j^{-1} ( s ) ] \) we have that each \( \mathrm{U}_s ( X ) \) is isomorphic to \( \mathrm{U}_s \) via some canonical \( \iota_s \) which maps the unique vertex in \( \mathrm{U}_s ( X ) \) adjacent to \( j^{-1} ( s ) \) to the point \( ( s ,\emptyset ) \). 
(To see that we do not need any choice to pick the isomorphisms \( \iota_s \), use Lemma~\ref{lem:propertiesL_alpha}\ref{lem:propertiesL_alpha-f} to first get a canonical isomorphism \( i_{\mathrm{U}_s ( X ) , \Code{s}} \colon \mathrm{U}_s ( X ) \to U_{\Code{s}} \), notice that by Lemma~\ref{lem:propertiesL_alpha}\ref{lem:propertiesL_alpha-e} we can always assume that the map \( i_{\mathrm{U}_s ( X ) , \Code{s}} \) sends the unique vertex in \( \mathrm{U}_s ( X ) \) adjacent to \( j^{-1} ( s ) \) to \( \emptyset \), and then set \( \iota_s ( a ) \coloneqq (s, i_{\mathrm{U}_s ( X ) , \Code{s} }(a) ) \) for every \( a \in \mathrm{U}_s ( X ) \).) 
Set \( \iota_{\mathrm{U}} \coloneqq \bigcup_{s \in \pre{< \omega}{\kappa}} \iota_s \), and notice that it is a well-defined bijection between \( \mathrm{U} ( X ) \) and \( \mathrm{U} \) because the \( \mathrm{U}_s ( X ) \) are pairwise disjoint.

\begin{claim} \label{claim:G_1minusG_0}
\( \iota' \coloneqq \iota_{\widehat{\mathrm{Seq}}} \cup \iota_{\mathrm{U}} \) is a partial isomorphism between \( \mathbb{G}_1 ( X ) \setminus \mathbb{G}_0 ( X ) = \widehat{\mathrm{Seq}} ( X ) \cup \mathrm{U} ( X ) \) and \( \mathbb{G}_1 \setminus \mathbb{G}_0 = \widehat{\mathrm{Seq}} \cup \mathrm{U} \).
\end{claim}

\begin{proof}[Proof of the Claim]
The fact that \( \iota' \) is a bijection between \( \mathbb{G}_1 ( X ) \setminus \mathbb{G}_0 ( X ) \) and \( \mathbb{G}_1 \setminus \mathbb{G}_0 \) is obvious, so we just need to show that it preserves the edge relation.
Let \( a , b \in \widehat{\mathrm{Seq}} ( X ) \cup \mathrm{U} ( X ) \). 
If \( X \models \widehat{\mathsf{Seq}} [ a ] \), then all neighbors of \( a \) must have degree \( \geq 4 \), so if \( a \edge^X b \), then \( X \not\models \widehat{\mathsf{Seq} } [ b ] \): this shows that any two points in \( \widehat{\mathrm{Seq}} ( X ) \) are not connected by en edge. 

Suppose now \( a , b \in \mathrm{U} ( X ) \). 
If \( a \in \mathrm{U}_s ( X ) \) and \( b \in \mathrm{U}_t ( X ) \) for distinct \( s , t \in \pre{<\omega}{\kappa} \), then \( \neg ( a \edge^X b) \), because otherwise \( X \models \mathsf{root} [ a , j^{-1} ( s ) ] \wedge \mathsf{root} [ a , j^{-1} ( t ) ] \), contradicting \( X \models \Upphi_1 \). 
If instead \( a,b \in \mathrm{U}_s ( X ) \) for the same \( s \in \pre{<\omega}{X} \), then \( \iota' ( a ) = \iota_s ( a ) \in \mathrm{U}_s \) and \( \iota' ( b ) = \iota_s ( b ) \in \mathrm{U}_s \), so 
\[ 
{a \edge^X b} \IFF { \iota' ( a ) \edge^{\mathrm{U}_{s}} \iota' ( b ) } \IFF {\iota' ( a ) \edge^{\mathbb{G}_1} \iota' ( b ) } 
\] 
by the choice of \( \iota_s \) and \( \mathbb{G}_1 \restriction \mathrm{U}_s = \mathrm{U}_s \).

Finally, assume that \( a \in \widehat{\mathrm{Seq}} ( X ) \) and \( b \in \mathrm{U} ( X ) \), and let \( s \in \pre{<\omega}{\kappa} \) be such that \( j ( a ) = s \) (so that \( \iota' ( a ) = \iota_{\widehat{\mathrm{Seq}}}(a) = \hat{s} \)). 
If \( b \in \mathrm{U}_t ( X ) \) for some \( s \neq t \in \pre{<\omega}{\kappa} \), then \( \neg ( a \edge^X b ) \), as otherwise \( X \models \mathsf{root}[ b , j^{-1} ( s ) ] \wedge \mathsf{root}[ b , j^{-1} ( t ) ] \), contradicting \( X \models \Upphi_1 \) again. 
If instead \( b \in \mathrm{U}_s ( X ) \), then \( a \edge^X b \IFF \iota' ( b ) = \iota_s ( b ) = ( s , \emptyset ) \) by our choice of \( \iota_s \), so that \( {a \edge^X b} \IFF {\iota' ( a ) \edge^{\mathbb{G}_1} \iota' ( b )} \) by Definition~\ref{def:G_1}.

By checking the definition of \( \edge^{\mathbb{G}_1} \) in Definition~\ref{def:G_1}, the previous observations suffice to show the desired result. 
\end{proof}

Let us now consider an arbitrary point \( a \in \mathrm{Seq} ( X ) \). 
Since \( X \models \Upphi_3 \), there is a unique \( b \in \widehat{\mathrm{Seq}} ( X ) \) with \( b \edge^X a \): therefore, we can unambiguously set \( \iota_{\mathrm{Seq}} ( a ) \coloneqq j ( b ) \), and check that by \( X \models \Upphi_3 \) the map \( \iota_{\mathrm{Seq} ( X ) } \colon \mathrm{Seq} ( X ) \to \mathrm{Seq} \) is a bijection. 
Moreover, for every \( a \in X \) and \( s \in \pre{< \omega}{\kappa} \) we have \( X \models \mathrm{Seq}_s [ a ] \iff \iota_{\mathrm{Seq}} ( a ) = s \). 

Finally, if \( a \in \mathrm{Seq}^- ( X ) \), then by \( X \models \Upphi_4 \) there is a unique \( \emptyset \neq s \in \pre{<\omega}{\kappa} \) such that \( X \models \mathsf{Seq}^-_s [ a ] \), so that the map \( \iota_{\mathrm{Seq}^-} \) sending \( a \) to \( \iota_{\mathrm{Seq}^-}(a) \coloneqq s^- \) is a bijection between \( \mathrm{Seq}^- ( X ) \) and \( \mathrm{Seq}^- \).

Consider now the canonical bijection \( \iota_X \colon \mathbb{G}_1 ( X ) \to \mathbb{G}_1 \) defined by
\begin{equation}\label{eq:i_X} 
\iota_X \coloneqq \iota_{\mathrm{Seq}} \cup \iota_{\mathrm{Seq}^-} \cup \iota_{\widehat{\mathrm{Seq}}} \cup \iota_U,
\end{equation}
which is well-defined since the functions appearing in the union have pairwise disjoint domains. 
We claim that \( \iota_X \) is an isomorphism between \( \mathbb{G}_1 ( X ) \) and \( \mathbb{G}_1 \), so let us fix arbitrary \( a,b \in \mathbb{G}_1 ( X ) \). 
Since by Claim~\ref{claim:G_1minusG_0} we already know that \( \iota_X \restriction (\mathbb{G}_1 ( X ) \setminus \mathbb{G}_0 ( X ) ) = \iota' \) is an isomorphism between \( \mathbb{G}_1 ( X ) \setminus \mathbb{G}_0 ( X ) \) and \( \mathbb{G}_1 \setminus \mathbb{G}_0 \), we may assume without loss of generality that \( a \in \mathbb{G}_0 ( X ) \). 
Suppose first that \( b \in \mathbb{G}_1 ( X ) \setminus \mathbb{G}_0 ( X ) \). 
If \( a \in \mathrm{Seq}^- ( X ) \), then \( \neg (a \edge^X b) \) because \( a \) has only two neighbors each of which must be in \( \mathrm{Seq} ( X ) \) (by \( X \models \mathsf{Seq}^- [ a ] \)), while \( b \notin \mathbb{G}_0 ( X ) \supseteq \mathrm{Seq} ( X ) \) by case assumption. 
If instead \( a \in \mathrm{Seq} ( X ) \) and \( b \in \widehat{\mathrm{Seq}} ( X ) \), then \( a \edge^X b \iff \iota_X ( b ) = \widehat{\iota_X ( a )} \) by definition of \( \iota_{\mathrm{Seq}} \) and \( \iota_{\widehat{\mathrm{Seq}}} \). 
Finally, if \( a \in \mathrm{Seq} ( X ) \) and \( b \in \mathrm{U} ( x ) \), then \( \neg(a \edge^X b) \) because \( a \) must have only neighbors of degree \( 2 \) by \( X \models \mathsf{Seq} [ a ] \), while \( b \) has degree \( \neq 2 \) in \( X \) by \( X \models \mathrm{U} [ b ] \).

Suppose now that \( a,b \in \mathbb{G}_0 ( X ) \), and let us assume that in fact \( a \in \mathrm{Seq}^- ( X ) \). 
Then \( X \models \mathsf{Seq}^-_s [ a ] \) for some \( \emptyset \neq s \in \pre{<\omega}{\kappa} \), and \( \iota_X ( a ) = \iota_{\mathrm{Seq}^-}(a) = s^- \) by definition of \( \iota_{\mathrm{Seq}^-} \). 
This implies that \( a \) has only two neighbors \( c_0, c_1 \) in \( X \), and they are such that \( X \models \mathsf{Seq}_{s^\star}[c_0] \wedge \mathsf{Seq}_s [ c_1 ] \): therefore \( \iota_X ( c_0 ) = \iota_{\mathrm{Seq}} ( c_0 ) = s^\star \) and \( \iota_X ( c_1 ) = \iota_{\mathrm{Seq}}(c_1) = s \). 
It follows that \( b \) is connected by an edge to \( a \) if and only if \( b = c_0 \vee b = c_1 \) if and only if \( \iota_X ( b ) = s^\star \vee \iota_X ( b ) = s \). 
The same argument (with \( a \) and \( b \) switched) takes care of the case \( b \in \mathrm{Seq}^- ( X ) \), so we just need to consider the case \( a , b \in \mathrm{Seq} ( X ) \). 
But then \( \neg (a \edge^X b ) \) because \( b \) has degree \( \geq 4 \) by \( X \models \mathsf{Seq} [ b ] \), while all neighbors of \( a \) must have degree \( 2 \) by \( X \models \mathsf{Seq} [ a ] \).

Checking the definition of \( \edge^{\mathbb{G}_1} \) in Definition~\ref{def:G_1} again, it is now easy to check that the above observations suffice to show that \( \iota_X \colon \mathbb{G}_1 ( X ) \to \mathbb{G}_1 \) is an isomorphism. 
\end{proof}

Given \( u \in \pre{<\omega}{2} \) and a sequence \( \seqofLR{x_i}{i \in F_u} \) of variables, let \( \mathsf{F}_u( \seqofLR{x_i}{i \in F_u} ) \) abbreviate the following \( \LL_{\omega_1 \omega} \)-formula:
\begin{equation}\tag{$\mathsf{F}_u( \seqofLR{x_i}{i \in F_u} )$}
\Bigl ( \bigwedge_{\emptyset \neq i \in F_u} \mathsf{F} ( x_i ) \Bigr ) \wedge \mathsf{Seq} ( x_\emptyset ) \wedge \uptau_{\mathsf{qf}} ( F_u ) ( \seqofLR{ x_i }{ i \in F_u } ) .
\end{equation}
Given \( X \in \Mod^\kappa_\LL \) and \( u \in \pre{<\omega}{2} \), we call \markdef{\( X \)-fork (coding \( u \))} any substructure of \( X \) determined by a sequence of points \( \seqofLR{a_i \in X}{i \in F_u} \) such that \( X \models \mathsf{F}_u [ \seqofLR{a_i}{i \in F_u} ] \), and the point \( a_\emptyset \) is called \markdef{root} of such an \( X \)-fork. 

We now provide \( \LL_{\omega_1 \omega_1} \)-sentences \( \Upphi_5, \dotsc, \Upphi_8 \) which, together with the previous ones \( \Upphi_0, \dotsc, \Upphi_4 \), complete the description of an \( \LL \)-structure of the from \( \mathbb{G}_S \) (for any \( S \in \Tr ( 2 \times \omega ) \)): 
\begin{equation*} \tag{$\Upphi_5$}
\FORALL{x} \biggl [ \mathsf{F} ( x ) \IMPLIES \bigvee_{u \in \pre{<\omega}{2}} \EXISTS{\seqofLR{ x_i }{ i \in F_u }} \Bigl ( \mathsf{F}_u ( \seqofLR{ x_i }{ i \in F_u } ) \wedge {\bigvee_{\emptyset \neq i \in F_u} x \equals x_i } \Bigr ) \biggr ];
\end{equation*}

\begin{multline*}\tag{$\Upphi_6$}
 \forall x \bigwedge_{u \in \pre{< \omega}{ 2}} \forall \seqofLR{ y_i }{ i \in F_u } \biggl [ \vphantom{\bigwedge_{\emptyset \neq i \in F_u}} \mathsf{F}_u ( \seqofLR{ y_i }{ i \in F_u } ) \wedge {\bigwedge_{i \in F_u} x \nequals y_i} 
 \\ 
{} \IMPLIES {\bigwedge_{\emptyset \neq i \in F_u} \Big ( \neg ( x \edge y_i ) \wedge \neg ( y_i \edge x ) \Big )} \biggr ] ;
\end{multline*}

\begin{multline*} {\tag{$\Upphi_7$}}
\bigwedge_{\substack{u \in \pre{< \omega}{ 2} \\ v \in \pre{< \omega}{ 2 }}} \FORALL{ \seqofLR{ x_i }{ i \in F_u }} \FORALL{ \seqofLR{ y_j }{ j \in F_v }} \biggl [ \mathsf{F}_u ( \seqofLR{ x_i }{ i \in F_u } ) \wedge \mathsf{F}_v ( \seqofLR{ y_{j} }{ j \in F_{v}} ) 
\\ 
 \IMPLIES {\bigwedge_{\substack{\emptyset \neq i \in F_u \\ \emptyset \neq j \in F_v}} ( x_i \nequals y_j )} \vee \Bigl ( \bigl ( {\bigwedge_{i \in F_u} \bigvee_{j \in F_v} ( x_i \equals y_j )} \bigr ) \wedge \bigl ( {\bigwedge_{j \in F_v} \bigvee_{i \in F_u} ( y_j \equals x_i )} \bigr ) \Bigr ) \biggr ] ;
\end{multline*}

\begin{multline*}\tag{$\Upphi_8$}
\bigwedge_{u \in \pre{ < \omega}{2}} \FORALL{ \seqofLR{ x_i}{i \in F_u }} \FORALL{ \seqofLR{ y_j}{ j \in F_u } } \biggl [ \vphantom{\bigwedge_{ i , j \in F_u }} \mathsf{F}_u ( \seqofLR{ x_i }{ i \in F_u } ) \wedge \mathsf{F}_u ( \seqofLR{ y_i }{ i \in F_u } ) 
\\
 \IMPLIES \bigl ( {\bigwedge_{ i , j \in F_u } ( x_i \nequals y_j )} \bigr ) \vee \bigl ( \vphantom{\bigwedge_{ i , j \in F_u }} {\bigwedge_{ i \in F_u} ( x_i \equals y_i ) } \bigr ) \biggr ] .
\end{multline*}
To understand the meaning of the above sentences, observe that for every \( X \in \Mod^\kappa_\LL \) we have:
\begin{itemize}[leftmargin=1pc]
\item
\( X \models \Upphi_5 \) if and only if every \( a \in \mathrm{F} ( X ) \) belongs to an \( X \)-fork (coding some \( u \in \pre{<\omega}{2} \));
\item
\( X \models \Upphi_6 \) if and only if given an arbitrary \( X \)-fork, if a point \( a \in X \) does not belong to that fork then it can be connected by an edge only to its root;
\item
\( X \models \Upphi_7 \) if and only if any two distinct \( X \)-forks may share only their roots;
\item
\( X \models \Upphi_8 \) if and only if any two distinct \( X \)-forks coding \emph{the same \( u \in \pre{< \omega}{2} \)} must be disjoint.
\end{itemize}

\smallskip

Finally, let \( \Uppsi \) be the \( \LL_{\kappa^+ \kappa} \)-sentence 
\begin{equation}\label{eq:Psi}
 \bigwedge_{0 \leq i \leq 8 }\Upphi_i . \index[symbols]{Psi@\( \Uppsi \)}
\end{equation}

The following lemma is straightforward (see also Remark~\ref{rmk:substructures}).

\begin{lemma} \label{lem:range}
Let \( T \in \TT_\kappa \), \( f_T \) be as in~\eqref{eq:f_T}, and \( \Uppsi \) be the \( \LL_{\kappa^+ \kappa} \)-sentence in~\eqref{eq:Psi}.
Then \( f_T(x) \models \Uppsi \) for every \( x \in \pre{\omega}{2} \), whence \( \ran (f_T) \subseteq \Mod^\kappa_{\Uppsi} \) (see Figure~\ref{fig:reductions}).
\end{lemma} 

\subsection{A classification of the structures in \( \Mod^\kappa_\Uppsi \) up to isomorphism}

Define the map
\begin{equation} \label{eq:defg}
 g \colon \Mod^\kappa_{\Uppsi} \to \pre{\pre{< \omega}{ 2} \times \kappa}{2} 
\end{equation}
(see Figure~\ref{fig:reductions}) as follows: given \( X \in \Mod^\kappa_{\Uppsi} \) and \( (u, \alpha) \in \pre{< \omega}{2} \times \kappa \), set \( g ( X ) ( u,\alpha ) = 1 \) if and only if there is an \( X \)-fork coding \( u \) whose root \( a_\emptyset \in X \) is such that \( X \models \mathsf{Seq}_s [ a_\emptyset ] \) for the unique \( s \coloneqq s_\alpha \in \pre{< \omega}{\kappa} \) with \( \Code{s} = \alpha \), that is: 
\begin{equation} \label{eq:defgwithformula}
g ( X ) ( u , \alpha ) = 1 \IFF X \models \EXISTS{\seqofLR{x_i}{i \in F_u}} [ \mathsf{F}_u ( \seqofLR{x_i}{i \in F_u} ) \wedge \mathsf{Seq}_{s_\alpha} ( x_\emptyset ) ] .
 \end{equation}

\begin{proposition} \label{prop:classification}
The map \( g \) from~\eqref{eq:defg} reduces \( \cong \) to \( = \).
\end{proposition}

\begin{proof}
Let \( X , Y \in \Mod^\kappa_{\Uppsi} \) be isomorphic via some map \( \iota \), fix \( (u, \alpha) \in \pre{< \omega}{ 2} \times \kappa \), and let \( s = s_\alpha \in \pre{<\omega}{\kappa} \) be such that \( \Code{s} = \alpha \). 
Then for every sequence \( \seqofLR{a_i}{i \in F_u} \) 
\[ 
X \models \mathsf{F}_u [ \seqofLR{a_i}{i \in F_u} ] \wedge \mathsf{Seq}_s [ a_\emptyset ] \iff Y \models \mathsf{F}_u [ \seqofLR{ \iota ( a_i ) }{i \in F_u} ] \wedge \mathsf{Seq}_s [ \iota ( a_\emptyset ) ] . 
\]
It follows that \( g ( X ) = g ( Y ) \) by~\eqref{eq:defgwithformula}.

Conversely, let \( X , Y \in \Mod^\kappa_{\Uppsi} \), so that they both satisfy \( \Upphi_i \) for \( 0 \leq i \leq 8 \), and assume that \( g ( X ) = g ( Y ) \).
By Lemma~\ref{lem:descriptionofG_1}, there are canonical isomorphisms \( \iota_X \colon \mathbb{G}_1 ( X ) \to \mathbb{G}_1 \) and \( \iota_Y \colon \mathbb{G}_1 ( Y ) \to \mathbb{G}_1 \), where \( \mathbb{G}_1 ( X ) \) and \( \mathbb{G}_1 ( Y ) \) are defined as in~\eqref{eq:substructuresofX-g} and \( \iota_X \), \( \iota_Y \) are the maps from~\eqref{eq:i_X}.
We will now extend the isomorphism \( \iota^{-1}_Y \circ \iota_X \colon \mathbb{G}_1 ( X ) \to \mathbb{G}_1 ( Y ) \) to an isomorphism 
\[
\iota_{ X , Y } \colon X \to Y .
\] 
Let \( a \in X \setminus \mathbb{G}_1 ( X ) \). 
By \( X \models \Upphi_0 \) we have \( X \models \mathsf{F} [ a ] \), and hence by \( X \models \Upphi_5 \) there are \( u \in \pre{ < \omega}{2} \) and \( \seqofLR{ a_i }{i \in F_u} \) such that \( X \models \mathsf{F}_u [ \seqofLR{ a_i }{ i \in F_u } ] \) and \( a = a_{\bar{\imath}} \) for some \( \emptyset \neq \bar{\imath} \in F_u \). 
In particular, \( X \models \mathsf{Seq} [ a_\emptyset ] \), so that by \( X \models \Upphi_2 \wedge \Upphi_3 \) there is (a unique) \( s \in \pre{<\omega}{\kappa} \) such that \( X \models \mathsf{Seq}_s [ a_\emptyset ] \) (see also the definition of \( \iota_{\mathrm{Seq}} \) in the proof of Lemma~\ref{lem:descriptionofG_1}). 
By definition of \( g \) we have \( g ( X ) ( u , \Code{s} ) = 1 \), so \( g ( Y ) ( u , \Code{s} ) = 1 \) by \( g ( X ) = g ( Y ) \). 
Let \( \seqofLR{b_i}{i \in F_u} \) be a sequence of elements of \( Y \) such that \( Y \models \mathsf{F}_u [ \seqofLR{b_i}{i \in F_u} ] \wedge \mathsf{Seq}_s [ b_\emptyset ] \), and set \( \iota_{ X , Y } ( a ) \coloneqq b_{\bar{\imath}} \). 
The definition of \( \iota_{ X , Y } ( a ) \) seems to depend on the choice of the sequence \( \seqofLR{b_i}{i \in F_u} \), but the next claim shows that this is not the case. 

\begin{claim}\label{claim:welldefined}
Suppose \( s \in \pre{<\omega}{\kappa} \), \( u \in \pre{< \omega}{ 2} \), and \( \seqofLR{ a_i }{ i \in F_u } \) is a sequence of elements of \( X \) such that \( X \models \mathsf{F}_u [ \seqofLR{ a_i }{ i \in F_u } ] \wedge \mathsf{Seq}_s [ a_\emptyset ] \).
Then there is a unique sequence \( \seqofLR{ b_i }{ i \in F_u } \) of elements of \( Y \) such that \( Y \models \mathsf{F}_u [ \seqofLR{ b_i }{ i \in F_u } ] \wedge \mathsf{Seq}_s [ b_\emptyset ] \).
Therefore \( \iota_{ X , Y } ( a_i ) = b_i \) for every \( i \in F_u \).
\end{claim}

\begin{proof}[Proof of the Claim]
Given two sequences \( \seqofLR{b_i }{ i \in F_u } \) and \( \seqofLR{ b'_i }{ i \in F_u } \) as above, then \( b_\emptyset = b'_\emptyset \) by \( Y \models \mathsf{Seq}_s [ b_\emptyset ] \wedge \mathsf{Seq}_s [ b'_\emptyset ] \) and \( Y \models \Upphi_3 \). 
This, together with \( Y \models \mathsf{F}_u [ \seqofLR{ b_i }{ i \in F_u } ] \) and \( Y \models \mathsf{F}_u [ \seqofLR{ b'_i }{ i \in F_u } ] \), implies \( b_i = b'_i \) for all \( i \in F_u \) by \( Y \models \Upphi_8 \).
\end{proof}

This shows that \( \iota_{ X , Y } \colon X \to Y \) is a well-defined map. 
We next check that it is a bijection. 

\begin{claim}
 \( \iota_{ X , Y } \) is injective.
\end{claim}

\begin{proof}[Proof of the Claim]
Let \( a,a' \in X \) be distinct.
If at least one of \( a , a' \) belongs to \( \mathbb{G}_1 ( X ) \), then \( \iota_{ X , Y } ( a) \neq \iota_{ X , Y } ( a') \) because \( \iota_{ X , Y } \restriction \mathbb{G}_1 ( X ) = \iota^{-1}_Y \circ \iota_X \) is a bijection between \( \mathbb{G}_1 ( X ) \) and \( \mathbb{G}_1 ( Y ) \), and \( \iota_{ X , Y } ( X \setminus \mathbb{G}_1 ( X ) ) \subseteq Y \setminus \mathbb{G}_1 ( Y ) \) by construction. 

Thus we can assume \( a , a' \notin \mathbb{G}_1 ( X ) \), so that \( X \models \mathsf{F} [ a ] \) and \( X \models \mathsf{F} [ a' ] \) by \( X \models \Upphi_0 \). 
Assume towards a contradiction that \( \iota_{ X , Y } ( a ) = \iota_{ X , Y } ( a' ) \). 
Since \( X \models \Upphi_5 \), there are \( u,v \in \pre{<\omega}{2} \) and sequences \( \seqofLR{a_i}{i \in F_u} \) and \( \seqofLR{a'_i}{j \in F_v} \) such that \( X \models \mathsf{F}_u [ \seqofLR{a_i}{i \in F_u} ] \), \( X \models \mathsf{F}_v [ \seqofLR{a'_i}{j \in F_v} ] \), and \( a = a_{\bar{\imath}} \), \( a' = a'_{\bar{\jmath}} \) for suitable \( \emptyset \neq \bar{\imath} \in F_u \) and \( \emptyset \neq \bar{\jmath} \in F_v \). 
Moreover, by \( X \models \Upphi_2 \wedge \Upphi_3 \) there are \( s , t \in \pre{<\omega}{\kappa} \) such that \( X \models \mathsf{Seq}_s [ a_\emptyset ] \) and \( X \models \mathsf{Seq}_t [ a'_\emptyset ] \). 
Since \( g ( X ) = g ( Y ) \) there are two sequences \( \seqofLR{b_i }{ i \in F_u } \) and \( \seqofLR{b'_j }{ j \in F_v } \) such that \( Y \models \mathsf{F}_u [ \seqofLR{ b_i }{ i \in F_u } ] \wedge \mathsf{Seq}_s [ b_\emptyset ] \) and \( Y \models \mathsf{F}_v [ \seqofLR{ b'_j }{ j \in F_v } ] \wedge \mathsf{Seq}_t [ b'_\emptyset ] \) (and moreover such sequences are unique by Claim~\ref{claim:welldefined}).
Then by definition of \( \iota_{ X , Y } \) we get \( \iota_{ X , Y } ( a_i ) = b_i \) and \( \iota_{ X , Y } ( a'_j ) = b'_j \) for every \( i \in F_u \) and \( j \in F_v \), so that in particular \( \iota_{ X , Y } ( a ) = b_{\bar{\imath}} \) and \( \iota_{ X , Y } ( a' ) = b'_{\bar{\jmath}} \). 
Since \( Y \models \Upphi_7 \) and 
\begin{equation}\label{eq:b_i}
 b_{\bar{\imath}} = \iota_{ X , Y } ( a ) = \iota_{ X , Y } ( a' ) = b'_{\bar{\jmath}} ,
\end{equation}
it follows that 
\begin{equation}\label{eq:b_n}
\setofLR{ b_i }{i \in F_u} = \setofLR{ b'_j }{j \in F_v},
\end{equation}
Since \( Y \models \uptau_{\mathsf{qf}} ( F_u ) [ \seqofLR{b_i}{i \in F_u} ] \wedge \uptau_{\mathsf{qf}} ( F_v )[ \seqofLR{b'_j}{j \in F_v} ] \), the substructure of \( Y \) with domain \( \setofLR{ b_i }{i \in F_u} = \setofLR{ b'_j }{j \in F_v} \) is isomorphic to \( F_u \) via the map \( b_i \mapsto i \) and to \( F_v \) via the map \( b'_j \mapsto j \), so that \( F_u \cong F_v \). 
Thus \( u = v \) by~\eqref{eq:=iffiso}, and by~\eqref{eq:b_i} and \( Y \models \Upphi_8 \) we get that \( b_i = b'_i \) for all \( i \in F_u = F_v \). 
Notice that this fact, together with~\eqref{eq:b_i} and the fact that all the \( b_i \)'s are necessarily distinct, also implies \( \bar{\imath} = \bar{\jmath} \). 
Moreover, \( b_\emptyset = b'_\emptyset \).
Since \( a_\emptyset,a'_\emptyset \in \mathbb{G}_1 ( X ) \), \( \iota_{ X , Y }(a_\emptyset) = b_\emptyset \), \( \iota_{ X , Y }(a'_\emptyset) = b'_\emptyset \), and \( \iota_{ X , Y } \restriction \mathbb{G}_1 ( X ) = \iota_Y^{-1} \circ \iota_X \) is a bijection, we get \( a_\emptyset = a'_\emptyset \). 
Since \( X \models \Upphi_8 \) and \( u = v \), we then get \( a_i = a'_i \) for all \( i \in F_u = F_v \), and recalling that we showed \( \bar{\imath} = \bar{\jmath} \) we finally get \( a = a_{\bar{\imath}} = a'_{\bar{\imath}} = a'_{\bar{\jmath}} = a' \), a contradiction. 
\end{proof}

\begin{claim}
\( \iota_{ X , Y } \) is surjective. 
\end{claim}

\begin{proof}[Proof of the Claim]
Given \( b \in Y \) we want to find \( a \in X \) with \( \iota_{ X , Y } ( a ) = b \). 
If \( b \in \mathbb{G}_1 ( Y ) \), then this follows from the fact that \( \iota_{ X , Y } \restriction \mathbb{G}_1 ( X ) \) is a bijection between \( \mathbb{G}_1 ( X ) \) and \( \mathbb{G}_1 ( Y ) \), so we may assume \( b \in Y \setminus \mathbb{G}_1 ( Y ) \). 
Since \( Y \models \Upphi_0 \wedge \Upphi_5 \), this implies that there is \( u \in \pre{<\omega}{2} \) and a sequence \( \seqofLR{b_i}{i \in F_u} \) of elements of \( Y \) such that \( Y \models \mathsf{F}_u [ \seqofLR{b_i}{i \in F_u} ] \) and \( b = b_{\bar{\imath}} \) for some \( \emptyset \neq \bar{\imath} \in F_u \). 
Moreover, by \( Y \models \Upphi_2 \wedge \Upphi_3 \) we know that there is also \( s \in \pre{<\omega}{\kappa} \) such that \( Y \models \mathsf{Seq}_s [ b_\emptyset ] \), so \( g ( Y ) ( u , \Code{s} ) = 1 \). 
Since we are assuming \( g ( X ) = g ( Y ) \), this means that there is a sequence \( \seqofLR{a_i}{i \in F_u } \) of elements of \( X \) such that \( X \models \mathsf{F}_u [ \seqofLR{a_i}{i \in F_u} ] \wedge \mathsf{Seq}_s [ a_\emptyset ] \). 
By our definition of \( \iota_{X,Y} \) and Claim~\ref{claim:welldefined} we then have \( \iota_{ X , Y } ( a_{\bar{\imath}} ) = b_{\bar{\imath}} = b \).
\end{proof}

It remains to show that \( \iota_{ X , Y } \) is also an isomorphism, i.e.\ that it preserves the edge relation.
Fix \( a , a ' \in X \). 
Since \( \iota_{ X , Y } \restriction \mathbb{G}_1 ( X ) = \iota^{-1}_Y \circ \iota_X \) is already an isomorphism between \( \mathbb{G}_1 ( X ) \) and \( \mathbb{G}_1 ( Y ) \), we may assume without loss of generality that \( a \notin \mathbb{G}_1 ( X ) \). 
By the usual argument repeatedly used above, we then get from \( g ( X ) = g ( Y ) \) that there are \( u \in \pre{<\omega}{2} \), \( s \in \pre{<\omega}{\kappa} \), a sequence \( \seqofLR{a_i}{i \in F_u } \) of points of \( X \), and \( \bar{\imath} \in F_u \) such that \( X \models \mathsf{F}_u [ \seqofLR{a_i}{i \in F_u} ] \wedge \mathsf{Seq}_s [ a_\emptyset ] \) and \( a = a_{\bar{\imath}} \), together with a (unique) sequence of points \( \seqofLR{b_i}{i \in F_u} \) of \( Y \) such that \( Y \models \mathsf{F}_u [ \seqofLR{b_i}{i \in F_u} ] \wedge \mathsf{Seq}_s [ b_\emptyset ] \), so that \( \iota_{ X , Y } ( a_i ) = b_i \) for every \( i \in F_u \) by construction (in particular, \( \iota_{ X , Y } ( a ) = b_{\bar{\imath}} \)). 
We distinguish two cases. 
If \( a' \neq a_i \) for every \( i \in F_u \), then \( \iota_{ X , Y } ( a' ) \neq b_i \) for every \( i \in F_u \) by injectivity of \( \iota_{ X , Y } \) and the fact that \( \iota_{ X , Y }( a_i ) = b_i \) for all \( i \in F_u \). 
But since both \( X \) and \( Y \) satisfy \( \Upphi_6 \), this implies that \( a = a_{\bar{\imath}} \) is not \( \edge^X \)-related to \( a' \) and \( \iota_{X,Y} ( a ) = b_{\bar{\imath}} \) is not \( \edge^Y \)-related to \( \iota_{ X , Y } ( a' ) \). 
If instead \( a' = a_{\bar{\jmath}} \) for some \( \bar{\jmath} \in F_u \), then \( \iota_{ X , Y }(a') = b_{\bar{\jmath}} \), and hence 
\[ 
a = a_{\bar{\imath}} \edge^X a_{\bar{\jmath}} = a' \IFF \bar{\imath} \edge^{F_u} \bar{\jmath} \IFF \iota_{ X , Y } ( a ) = b_{\bar{\imath}} \edge^Y b_{\bar{\jmath}} = \iota_{ X , Y } ( a' ) . \qedhere
 \] 
\end{proof}

Recall from~\eqref{eq:defgwithformula} that for every \( X \in \Mod^\kappa_{\Uppsi} \) and \( ( u , \alpha ) \in \pre{<\omega}{2} \times \kappa \), we have set \( g ( X ) ( u , \alpha ) = 1 \) if and only if \( X \) satisfies the \( \LL_{\kappa^+ \kappa} \)-sentence
\[ 
\EXISTS{\seqofLR{x_i}{i \in F_u}} [ \mathsf{F}_u ( \seqofLR{x_i}{i \in F_u} ) \wedge \mathsf{Seq}_s ( x_\emptyset ) ]
 \] 
where \( s = s_\alpha \in \pre{< \omega}{\kappa} \) is such that \( \Code{s} = \alpha \).
For technical reasons we need to replace such a sentence with one belonging to the bounded version \( \LL^b_{\kappa^+ \kappa} \) of \( \LL_{\kappa^+ \kappa} \) (see Definition~\ref{def:Lbounded}\ref{def:Lbounded-b}) --- this will be crucial for the results in Section~\ref{subsec:topologicalcomplexity}.

\begin{lemma} \label{lem:bettersentence}
Let \( X \in \Mod^\kappa_{\Uppsi} \) and \( ( u , \alpha ) \in \pre{<\omega}{2} \times \kappa \). 
Then \( g ( X ) ( u , \alpha ) = 1 \IFF X \models \upsigma_{u,\alpha} \), where \( \upsigma_{u , \alpha} \) is the \( \LL^b_{\kappa^+ \kappa} \)-sentence
\begin{multline}\tag{$\upsigma_{u , \alpha }$}
\EXISTS{x} \EXISTS{y} \left [ \mathsf{Seq} ( x ) \AND \widehat{\mathsf{Seq}} ( y ) \AND {x \edge y} \AND \EXISTS{\seqofLR{z_i}{i \in U_\alpha}} \bigg({y \edge z_\emptyset} \AND \uptau_{\mathsf{qf}}(U_\alpha)( \seqofLR{z_i}{i \in U_\alpha}) \bigg ) \right . 
 \\
\wedge \neg \EXISTS{\seqofLR{w_j}{j \in U_{\alpha + 1}}} \bigg ( {y \edge w_\emptyset} \AND \uptau_{\mathsf{qf}}(U_{\alpha+1})(\seqofLR{w_j}{j \in U_{\alpha+1}}) \bigg ) 
\\
\left . \vphantom{\widehat{\mathsf{Seq}} ( y )} \AND \EXISTS{\seqofLR{l_k}{k \in F_u^{\theta(u)} }} \bigg ( l_\emptyset \equals x \AND \left( \bigwedge\nolimits_{\emptyset \neq k \in F_u^{\theta(u)} } \mathsf{F} ( l_k ) \right) \AND \uptau_{\mathsf{qf}}(F_u^{\theta(u)} ) \left (\seqofLR{l_k}{k \in F_u^{\theta(u)} } \right) \bigg ) \right ] ,
\end{multline}
where we set \( F_u^{\theta(u)} \coloneqq F_u \cap \pre{\theta(u)+1}{2} \).
\end{lemma}

\begin{proof}
Using Proposition~\ref{lem:propertiesL_alpha}\ref{lem:propertiesL_alpha-c}, it is not hard to see that if \( g ( X ) (u,\alpha) = 1 \) then \( X \models \upsigma_{u,\alpha} \). 
For the other direction, assume that \( X \models \upsigma_{u,\alpha} \), and let \( a \), \( \hat{a} \), \( \seqofLR{b_i}{i \in U_\alpha} \), and \( \seqofLR{c_k }{ k \in F_u^{\theta(u)} } \) be (sequences of) elements of \( X \) such that 
\begin{multline*}
X \models \mathsf{Seq} [ a ] \AND \widehat{\mathsf{Seq}} [ \hat{a} ] \AND \bigg ({a \edge \hat{a}} \AND {\hat{a} \edge b_\emptyset} \AND \uptau_{\mathsf{qf}}( U_\alpha) [ \seqofLR{b_i}{i \in U_\alpha} ] \bigg ) 
\\
\AND \neg \EXISTS{\seqofLR{w_j}{j \in U_{\alpha+1}}} \bigg ({\hat{a} \edge w_\emptyset} \AND \uptau_{\mathsf{qf}}(U_{\alpha+1})(\seqofLR{w_j}{j \in U_{\alpha+1}}) \bigg ) 
\\
\AND \bigg ( {c_\emptyset \equals a} \AND \left( \bigwedge\nolimits_{\emptyset \neq k \in F_u^{\theta(u)} } \mathsf{F} [ c_k ] \right) 
\AND \uptau_{\mathsf{qf}}(F_u^{\theta(u)} )\left [\seqofLR{c_k}{k \in F_u^{\theta(u)} }\right] \bigg ) .
 \end{multline*}
It follows from \( X \models \Upphi_2 \wedge \Upphi_3 \) and \( X \models \mathsf{Seq} [ a ] \wedge \widehat{\mathsf{Seq}} [ \hat{a} ] \wedge {a \edge \hat{a}} \) that \( X \models \mathsf{Seq}_s [ a ] \) and \( X \models \widehat{\mathsf{Seq}}_s[\hat{a}] \) for one and the same \( s \in \pre{< \omega}{\kappa} \): we claim that \( \Code{s} = \alpha \). 
Indeed, from \( X \models \uptau_{\mathsf{qf}}(U_\alpha)(\seqofLR{b_i}{i \in U_\alpha}) \) and \( \hat{a} \edge^X b_\emptyset \) it easily follows that \( X \models \mathsf{U} [ b_i ] \wedge \mathsf{root} [ b_i , \hat{a} ] \) for every \( i \in U_\alpha \). 
Arguing by contradiction, one sees that Lemma~\ref{lem:propertiesL_alpha}\ref{lem:propertiesL_alpha-c} together with \( X \models \widehat{\mathsf{Seq}}_s [ \hat{a} ] \) and \( X \models \uptau_{\mathsf{qf}}( U_\alpha)(\seqofLR{b_i}{i \in U_\alpha}) \) imply \( \Code{s} \geq \alpha \). 
Moreover, by Lemma~\ref{lem:propertiesL_alpha}\ref{lem:propertiesL_alpha-b} we get that \( \Code{s} > \alpha \) would contradict 
\[ 
X \models \neg \EXISTS{\seqofLR{w_j}{j \in U_{\alpha+1}}} \bigg ( {\hat{a} \edge w_\emptyset} \wedge \uptau_{\mathsf{qf}}(U_{\alpha+1})(\seqofLR{w_j}{j \in U_{\alpha+1}}) \bigg ) .
\] 
Therefore, \( \Code{s} = \alpha \).
Since \( X \models c_\emptyset \equals a \), we have also showed that \( X \models \mathsf{Seq}_s[c_\emptyset] \) for the unique \( s \in \pre{<\omega}{\kappa} \) such that \( \Code{s} = \alpha \). 

Fix any \( \emptyset \neq k \in F_u^{\theta(u)} \). 
Since \( X \models \mathsf{F}[c_k] \), by \( X \models \Upphi_5 \) we get that there are \( {v_k} \in \pre{< \omega}{2} \) and a sequence \( \seqofLR{d^k_i}{i \in F_{v_k}} \) of elements of \( X \) such that \( X \models \mathsf{F}_{v_k}[\seqofLR{d^k_i}{i \in F_{v_k}}] \) and \( X \models c_k \equals d^k_{i_k} \) for some \( \emptyset \neq i_k \in F_{v_k} \). 
In particular, \( X \models \mathsf{Seq}[d^k_\emptyset] \): we claim that \( d^k_\emptyset = c_\emptyset \). 
If not, then by \( X \models \bigwedge_{i \leq 4} \Upphi_i \) and Lemma~\ref{lem:descriptionofG_1} there would be a path of length \( \geq 2 \) connecting \( c_\emptyset \) to \( d^k_\emptyset \) which is totally contained in \( \mathbb{G}_0(X) \subseteq \mathbb{G}_1(X) \). 
On the other hand, the set of vertices \( \setofLR{c_{k'}}{k' \subseteq k} \cup \setofLR{d^k_i}{i \subseteq i_k} \) would contain a path of length \( \geq 2 \) connecting \( c_\emptyset \) to \( d^k_\emptyset \) which, except for its extreme points, is totally contained in \( \mathrm{F}(X) \). 
Since \( \mathbb{G}_1(X) \cap \mathrm{F}(X) = \emptyset \) by Remark~\ref{rmk:disj}, the two paths would then be distinct, contradicting the fact that \(X \) is acyclic by \( X \models \Upphi_0 \). 
A similar argument shows that by acyclicity of \( X \) the two paths \( \seqofLR{c_{k'}}{k' \subseteq k} \) and \( \seqofLR{d^k_i}{i \subseteq i_k} \) joining \( c_\emptyset = d^k_\emptyset \) to \( c_k = d^k_{i_k} \) must coincide. 
Thus, in particular, \( c_{\seq{0}} = d^k_{\seq{0}} \) for every \( \emptyset \neq k \in F_u^{\theta(u)} \). 
Fix now \( k , k' \in F_u^{\theta(u)} \setminus \{ \emptyset \} \).
Since \( d^k_{\seq{0}} = c_{\seq{0}} = d^{k'}_{\seq{0}} \), by \( X \models \Upphi_7 \) it follows that 
\[ 
\setofLR{d^k_i}{i \in F_{v_k}} = \setofLR{d^{k'}_i}{i \in F_{v_{k'}}}.
 \] 
 Since the maps \( i \mapsto d^k_i \) and \( i \mapsto d^{k'}_i \) witness that both \( F_{v_k} \) and \( F_{v_{k'}} \) are isomorphic to the same structure \( X \restriction \setofLR{d^k_i}{i \in F_{v_k}} = X \restriction \setofLR{d^{k'}_i}{i \in F_{v_{k'}}} \), we have \( F_{v_k} \cong F_{v_{k'}} \), whence \( v_k = v_{k'} \) by~\eqref{eq:=iffiso}. 
It then follows from \( X \models \Upphi_8 \) and \( d^k_{\seq{0}} = d^{k'}_{\seq{0}} \) that \( d^k_i = d^{k'}_i \) for all \( i \in F_{v_k} = F_{v_{k'}} \). 
Set \( v \coloneqq v_k \) and \( d_i \coloneqq d^k_i \) for some/any \( \emptyset \neq k \in F_u^{\theta(u)} \) and all \( i \in F_v \). 
Since \( \setofLR{c_k}{k \in F_u^{\theta(u)}} \subseteq \setofLR{d_i}{i \in F_v} \) and \( c_{0^{\theta(u)}} \) has three distinct neighbors among the \( c_k \)'s, necessarily \( c_{0^{\theta(u)}} = d_{0^{\theta(v)}} \). 
From this and \( c_\emptyset = d_\emptyset \) it also follows that \( \setofLR{c_k}{k \subseteq 0^{\theta(u)}} = \setofLR{d_i}{i \subseteq 0^{\theta(v)}} \) (here we use again the fact that \( X \) is acyclic), whence \( \theta(u) = \theta(v) \), and hence also \( u = v \) by injectivity of \( \theta \). 
Since \( d_\emptyset = c_\emptyset \) and \( X \models \mathsf{Seq}_s[c_\emptyset] \), it thus follows that \( \seqofLR{d_i}{i \in F_u} \) witnesses that
\[ 
X \models \EXISTS{\seqofLR{x_i}{i \in F_u}}[\mathsf{F}_u(\seqofLR{x_i}{i \in F_u}) \AND \mathsf{Seq}_s(x_\emptyset)],
 \] 
and since \( \Code{s}= \alpha \) we have \( g(X)(u,\alpha) = 1 \), as desired.
\end{proof}

\subsection{The invariant universality of \( \embeds^\kappa_\CT \)}

Endow \( \pre{\pre{< \omega}{ 2} \times \kappa} {2} \) with the \emph{product} topology \( \tau_p \), so that \( \left (\pre{\pre{<\omega}{2} \times \kappa}{2}, \tau_p \right ) \) is homeomorphic to \( ( \pre{\kappa}{2}, \tau_p ) \) (see Section~\ref{subsubsec:spacesoftypekappa} and, in particular, Example~\ref{xmp:canonicalbijection}\ref{xmp:canonicalbijection-1}). 

\begin{lemma} \label{lem:open}
For every (\( \tau_p\)-)open set \( V \subseteq \pre{\pre{< \omega}{ 2} \times \kappa}{2} \) there is an \( \LL^b_{ \kappa^+ \kappa } \)-sentence \( \upsigma_V \) such that \( g^{ - 1 } ( V ) = \Mod^\kappa_{ \upsigma_V } \cap \Mod^\kappa_{ \Uppsi} \) (equivalently, \( g^{-1}(V) = \Mod^\kappa_{\upsigma_V \wedge \Uppsi} ) \).
\end{lemma}

The fact that the sentence \( \upsigma_V \) belongs to the fragment \( \LL^b_{\kappa^+ \kappa} \) (and not just to \( \LL_{\kappa^+ \kappa} \)) will be crucially used in Section~\ref{subsec:topologicalcomplexity} to prove that certain maps are (effective) \( \kappa+ 1\)-Borel --- see Lemma~\ref{lem:inverse} and Theorem~\ref{th:maintopology}.

\begin{proof}
Recall from~\eqref{eq:subbasisproduct} that \( \mathcal{S} = \setofLR{ \widetilde{\Nbhd}^A_{ ( u , \alpha ) , i } }{ {( u , \alpha ) \in \pre{< \omega}{ 2} \times \kappa} \AND {i = 0,1} } \) is a subbasis for the product topology on \( \pre{A}{2} \), where \( A \coloneqq \pre{< \omega}{ 2} \times \kappa \).
The map \( V \mapsto \upsigma_ V \) will be defined first on \( \mathcal{S} \), then on the canonical basis \( \mathcal{B}_p \) (which is generated by \( \mathcal{S} \) by taking finite intersections), and finally on all \( \tau_p \)-open sets.

When \( V = \widetilde{\Nbhd}^A_{ ( u ,\alpha ) , i} \) let \( \upsigma_V \) be \( \upsigma_{ u, \alpha } \) if \( i = 1 \) or \( \neg \upsigma_{u,\alpha} \) otherwise. 
By Lemma~\ref{lem:bettersentence} we get \( g^{-1} ( V ) = \Mod^\kappa_{\upsigma_V} \cap \Mod^\kappa_\Uppsi \).

Pick \( V \in \mathcal{B}_p \) and let \( U_1 ,\dots , U_n \in \mathcal{S} \) be such that \( V = U_1 \cap \dots \cap U_n \).
Then \( g^{-1} ( V ) = \Mod^\kappa_{\upsigma_V } \cap \Mod^\kappa_\Uppsi \), where
\( \upsigma_V \) is \( \upsigma_{U_1} \wedge \dots \wedge\upsigma_{U_n} \).
 
Given now an arbitrary \( \tau_p \)-open set \( V \subseteq \pre{\pre{< \omega}{ 2} \times \kappa}{2} \), let \( B_V \) be the collection of those \( U \in \mathcal{B}_p \) which are contained in \( V \): since the cardinality of \( B_V \) is at most \( \card{[\kappa]^{< \omega}} = \kappa \) (which is the cardinality of \( \mathcal{B}_p \)), we get the desired result letting \( \upsigma_V \) be \( \bigvee_{U \in B_V} \upsigma_U \).
\end{proof}

\begin{lemma} \label{lem:continuous} 
Let \( T \in \TT_\kappa \) and \( f_T \) be the map from~\eqref{eq:f_T}. 
Then the map \( g \circ f_T \colon \pre{\omega}{ 2} \to \pre{\pre{< \omega}{ 2} \times \kappa}{2} \) is continuous.
\end{lemma}

\begin{proof}
Fix \( (u, \alpha) \in A = \pre{< \omega}{ 2} \times \kappa \), and let \( s \in \pre{<\omega}{\kappa} \) be such that \( \Code{s} = \alpha \). 
Let 
\[
 D_1 \coloneqq \setofLR{ v \in \pre{ \lh u }{2} }{ ( u , v , s ) \in \tilde{T} } ,
\]
where \( \tilde{T} \) is the tree obtained from \( T \) as in~\eqref{eq:S_T}, and let \( D_0 \coloneqq \pre{ \lh u}{2} \setminus D_1 \). 
Then by definition of \( f_T \) and \( g \) we get \( (g \circ f_T ) ^{-1} ( \widetilde{\Nbhd}^A_{ ( u , \alpha ) , i }) = \bigcup_{v \in D_i } \Nbhd_v \).
\end{proof}

\begin{corollary} \label{cor:closed}
Let \( T \in \TT_\kappa \) and \( f_T \) be as in~\eqref{eq:f_T}. 
For every closed set \( C \subseteq \pre{\omega}{ 2} \), \( (g \circ f_T)(C) \) is closed. 
\end{corollary}

\begin{proof}
The set \( C \) is compact since it is a closed subset of the compact space \( \pre{ \omega}{2} \). 
Therefore, since \( g \circ f_T \) is continuous by Lemma~\ref{lem:continuous}, \( (g \circ f_T)(C) \) is compact as well: but then it is also closed because \( \pre{\pre{< \omega}{2 \times \kappa}}{2} \) is a Hausdorff space.
\end{proof}

\begin{corollary} \label{cor:saturation}
Let \( T \in \TT_\kappa \) and \( f_T \) be as in~\eqref{eq:f_T}. 
Then there is an \( \LL_{\kappa^+ \kappa} \)-sentence \( \upsigma_T \) such that the closure under isomorphism of \( \ran ( f_T ) \) is \( \Mod^\kappa_{\upsigma_T} \), i.e.
\[ 
\Mod^\kappa_{\upsigma_T} = \setofLR{X \in \Mod^\kappa_\LL}{X \cong f_T ( x ) \text{ for some } x \in \pre{\omega}{2} } .
\]
\end{corollary}

Notice that since \( \ran (f_T) \subseteq \Mod^\kappa_\Uppsi \) by Lemma~\ref{lem:range}, we thus have in particular (see Figure~\ref{fig:reductions})
\[ 
\ran (f_T) \subseteq \Mod^\kappa_{\upsigma_T} \subseteq \Mod^\kappa_{\Uppsi}.
 \] 

\begin{proof}
Let \( V \) be the complement of \( \ran ( g \circ f_T ) = ( g \circ f_T ) ( \pre{\omega}{ 2} ) \), which is open by Corollary~\ref{cor:closed}. 
Use Lemma~\ref{lem:open} to find an \( \LL^b_{\kappa^+ \kappa} \)-sentence \( \hat{\upsigma} \coloneqq \upsigma_V \) such that \( g^{-1} ( V ) = \Mod^\kappa_{\hat{\upsigma}} \cap \Mod^\kappa_\Uppsi \). 
Finally, let \( \upsigma_T \) be \( \Uppsi \wedge \neg \hat{\upsigma} \), so that 
\begin{equation}\label{eq:rangf}
g^{-1} ( \ran ( g \circ f_T ) ) = \Mod^\kappa_{\upsigma_T}.
\end{equation}

If \( X \in \Mod^\kappa_\LL \) is such that \( X \cong f_T ( x ) \) for some \( x \in \pre{ \omega}{2} \), then \( X \in \Mod^\kappa_{\Uppsi} \) (use the fact that \( f_T(X) \in \Mod^\kappa_{\Uppsi} \) by Lemma~\ref{lem:range} and that \( \Mod^\kappa_{\Uppsi} \) is invariant under isomorphism): therefore \( g ( X ) \) is defined and equals \( ( g \circ f_T ) ( x ) \) by Proposition~\ref{prop:classification}, so that \( X \in g^{-1} ( \ran(g \circ f_T)) \) and hence \( X \models \upsigma_T \) by~\eqref{eq:rangf}. 
Conversely, if \( X \in \Mod^\kappa_{\upsigma_T} = \Mod^\kappa_{\Uppsi \wedge \neg \hat{\upsigma}} \subseteq \Mod^\kappa_{\Uppsi} \) then by~\eqref{eq:rangf} there must be \( x \in \pre{ \omega}{2} \) such that \( g ( X ) = ( g \circ f_T ) ( x ) \): but then \( X \cong f_T ( x ) \) by Proposition~\ref{prop:classification} again.
\end{proof}

We now define the last reduction \( h_T \colon \Mod^\kappa_{\upsigma_T} \to \pre{\omega}{2} \) of Figure~\ref{fig:reductions}.

\begin{definition}\label{def:defh_T}
Let \( T \in \TT_\kappa \) be such that \( R = \PROJ{\body{T}} \) is a quasi-order, \( f_T \) be as in~\eqref{eq:f_T}, and \( \upsigma_T \) as in Corollary~\ref{cor:saturation}. 
Then we define the map \( h_T \colon \Mod^\kappa_{\upsigma_T} \to \pre{\omega}{2} \) by setting for all \( X \in \Mod^\kappa_{\upsigma_T} \)
\begin{equation}\label{eq:defh_T}
h_T ( X ) = x \iff f_T ( x ) \cong X.
\end{equation}
\end{definition}

Notice that the function \( h_T \) is well-defined: by Corollary~\ref{cor:saturation}, if \( X \models \upsigma_T \) then there is at least one \( x \) satisfying~\eqref{eq:defh_T}, and by Theorem~\ref{th:graphs}\ref{th:graphs-b} such \( x \) is unique because if \( y \in \pre{\omega}{2} \) is distinct from \( x \), then \( f_T ( x ) \not\cong f_T ( y ) \), whence \( f_T ( y ) \not\cong X \).

\begin{corollary} \label{cor:inverse}
Let \( T \in \TT_\kappa \), \( f_T \) be as in~\eqref{eq:f_T}, and \( \upsigma_T \) be as in Corollary~\ref{cor:saturation}. 
If \( R = \PROJ{\body{T}} \) is a quasi-order, then the map \( h_T \colon \Mod^\kappa_{\upsigma_T} \to \pre{\omega}{2} \) from Definition~\ref{def:defh_T} simultaneously reduces \( \embeds \) to \( R \) and \( \cong \) to \( = \).
\end{corollary}

\begin{proof}
Let \( X , Y \in \Mod^\kappa_{\upsigma_T} \) and assume that \( R = \PROJ \body{ T } \) is a quasi-order.
Let \( x \coloneqq h_T ( X ) \) and \( y \coloneqq h_T ( Y ) \).
Since \( f_T ( x ) \cong X \) and \( f_T ( y ) \cong Y \) by~\eqref{eq:defh_T}, by Theorem~\ref{th:graphs} we have
\[ 
{X \embeds Y} \iff { f_T ( x ) \embeds f_T ( y ) } \iff {x \mathrel{R}y}, 
\]
and
\[ 
{X \cong Y} \iff {f_T ( x ) \cong f_T ( Y ) } \iff x = y. \qedhere 
\]
\end{proof}

Summing up what we obtained in this section, we have the following theorem (see also Figure~\ref{fig:reductions}). 

\begin{theorem} \label{th:main}
Let \( \kappa \) be an infinite cardinal and \( T \in \TT_\kappa \).
If \( R = \PROJ \body{ T } \) is a quasi-order, then there are \( \upsigma_T , f_T , h_T \) such that:
\begin{enumerate-(a)}
\item \label{th:main-a}
\( \upsigma_T \) is an \( \LL_{\kappa^+ \kappa} \) sentence;
\item \label{th:main-b}
\( f_T \) reduces \( R \) to \( \embeds ^\kappa_{\upsigma_T} \) and \( = \) to \( \cong ^\kappa_{\upsigma_T} \);
\item \label{th:main-c}
\( h_T \) reduces \( \embeds ^\kappa_{\upsigma_T} \) to \( R \) and \( \cong ^\kappa_{\upsigma_T} \) to \( = \);
\item \label{th:main-d}
\( h_T \circ f_T = \id \) and \( (f_T \circ h_T) ( X ) \cong X \) for every \( X \in \Mod^\kappa_{\upsigma_T} \).
\end{enumerate-(a)}
In particular, \( \embeds_\CT^\kappa \) is invariantly universal for \( \kappa \)-Souslin quasi-orders on \( \pre{\omega}{2} \), that is to say: for every \( \kappa \)-Souslin quasi-order \( R \) on \( \pre{ \omega}{2} \) there is an \( \LL_{\kappa^+ \kappa} \) sentence \( \upsigma \) all of whose models are combinatorial trees such that \( R \sim {\embeds ^\kappa_\upsigma} \).
\end{theorem}

\begin{proof}
For part~\ref{th:main-d}, observe that for every \( x \in \pre{\omega}{2} \) we must have \( f_T ((h_T \circ f_T)(x)) \cong f_T(x) \) and \( (h_T \circ f_T ) ( x ) = x \) by~\eqref{eq:defh_T}.
The other condition on \( f_T \circ h_T \) is just a rewriting of~\eqref{eq:defh_T}, so we are done.
\end{proof}

Since \( R \sim {\embeds ^\kappa_{\upsigma_T}} \) by~\ref{th:main-b} and~\ref{th:main-c}, it follows that the quotient order of \( R \) is bi-embeddable with the quotient order of \( {\embeds^\kappa_{\upsigma_T}} \). 
However,~\ref{th:main-d} implies the following stronger result.

\begin{corollary}\label{cor:maintheorem}
Let \( \kappa \) be an infinite cardinal and \( T \in \TT_\kappa \).
If \( R = \PROJ \body{ T } \) is a quasi-order, then there is an \( \LL_{\kappa^+ \kappa} \)-sentence \( \upsigma_T \) (all of whose models are combinatorial trees) such that the quotient orders of \( R \) and \( {\embeds^\kappa_{\upsigma_T}} \) are \emph{isomorphic}.
\end{corollary}

Notice also that using Corollary~\ref{cor:maintheorem} all results on invariant universality of Sections~\ref{sec:definablecardinality},~\ref{sec:applications}, and~\ref{subsubsec:changingmorphism} could be reformulated in terms of (definable) isomorphisms between the induced quotient orders.

\begin{remarks} \label{rem:maintheorem}
\begin{enumerate-(i)}
\item \label{rem:maintheorem-1}
It is not hard to check that all models of the \( \LL_{\kappa^+ \kappa} \)-sentence \( \Uppsi \), and hence also the models of the sentence \( \upsigma_T \) from Corollary~\ref{cor:saturation}, are necessarily of size \( \kappa \). 
Therefore \( \Uppsi \) \emph{characterizes} the cardinal \( \kappa \) in the sense of~\cite[p.~59]{Knight:1977mb}, once in the original definition we replace the logic \( \LL_{\omega_1 \omega} \) with the more powerful \( \LL_{\kappa^+ \kappa} \). 
This also implies that \( \Mod^\infty_{\upsigma_T} = \Mod^\kappa_{\upsigma_T} \), so that the conclusion of Theorem~\ref{th:main} could also be reformulated as follows: 
\begin{quote}
For every \( \kappa \)-Souslin quasi-order \( R \) on \( \pre{\omega}{2} \) there is an \( \LL_{\kappa^+ \kappa} \)-sentence \( \upsigma \) all of whose models are combinatorial trees such that \( R \sim {\embeds^\infty_\upsigma} \). 
\end{quote}
Notice that all results on invariant universality of Sections~\ref{sec:definablecardinality},~\ref{sec:applications}, and~\ref{subsubsec:changingmorphism} could be reformulated analogously.
\item \label{rem:maintheorem-2}
Every model \( X \) of \( \Uppsi \) (and hence every model of \( \upsigma_T \) for \( T \in \TT_\kappa \)) admits an \markdef{\( \LL_{\kappa^+ \kappa} \)-Scott sentence}, i.e.\ an \( \LL_{\kappa^+ \kappa} \)-sentence \( \upsigma^X \) such that for every \( Y \in \Mod^\kappa_\LL \)
\[ 
Y \models \upsigma^X \iff Y \cong X .
 \] 
The sentence \( \upsigma^X \) is obtained by applying an argument similar to that of Corollary \ref{cor:saturation}: Let \( V \coloneqq \pre{\pre{<\omega}{2} \times \kappa}{2} \setminus \setLR{ g ( X ) } \) (which is open since \( \pre{\pre{<\omega}{2} \times \kappa}{2} \) is a Hausdorff space), and let \( \upsigma^X \) be the \( \LL_{\kappa^+ \kappa} \)-sentence \( \Uppsi \wedge \neg \upsigma_V \), where \( \upsigma_V \) is as in Lemma~\ref{lem:open}.
\end{enumerate-(i)}
\end{remarks}

It is a well-known classical result due to D.~Scott that when \( \kappa = \omega \) then Remark~\ref{rem:maintheorem}\ref{rem:maintheorem-2} is true \emph{for every} \( X \in \Mod^\omega_\LL \), see e.g.\ 
\cite[Corollary~16.10]{Kechris:1995zt}. 
However, it is an old result, probably due to the Finnish school, that when \( \kappa \) is uncountable such condition may fail for some \( X \in \Mod^\kappa_{\LL} \).
The next example from~\cite{Vaananen:2011tg} shows this.
(We thank S.~D.~Friedman and the referee for having brought it to our attention.)

\begin{example} \label{xmp:Scottsentence}
Let \( \LL \) be the order language consisting of one binary relational symbol, set \( \kappa = \omega_1 \), and consider the collection of all \( \omega_1 \)-like dense linear orders on \( \omega_1 \) without a minimum (briefly, \( \omega_1 \)-DLO), where a linear order is called \( \omega_1 \)-like if all its proper initial segments are countable~\cite[Definition 9.6]{Vaananen:2011tg}. 
We build two \( \omega _1 \)-DLOs \( X , Y \) using the construction from~\cite[Definition 9.8]{Vaananen:2011tg}. 
Equip \( \Q_0 \coloneqq \Q \cap ( 0 ; 1 ) \), \( \Q_1 \coloneqq \Q \cap \cointerval{0}{1} \), and \( \omega_1 \) with the usual orders, and let 
\[ 
X \coloneqq \omega_1 \times \Q_0 
 \] 
and
\[ 
Y \coloneqq \setLR{ 0 } \times \Q_0 \cup (\omega_1 \setminus \setLR{ 0 }) \times \Q_1
 \] 
be endowed with the lexicographical ordering. 
Then, after identifying them with corresponding structures on \( \omega_1 \), \( X \) and \( Y \) are easily seen to be \( \omega_1 \)-DLO. 
Moreover, they are not isomorphic because \( Y \) contains a closed unbounded (with respect to its ordering) set of order type \( \omega_1 \), while \( X \) does not contain such a substructure (see~\cite[Lemma 9.9]{Vaananen:2011tg}, where \( X \) and \( Y \) are denoted by \( \Phi ( \emptyset ) \) and \( \Phi ( \omega _1 ) \)). 
Player \( \II \) has a winning strategy in the Ehrenfeuch-Fra\"iss\'e game \( \mathrm{EF}_\omega^{\omega_1} ( X , Y ) \) by~\cite[Lemma 9.10]{Vaananen:2011tg}, and thus \( X \) and \( Y \) are \emph{\( \LL_{\omega_2 \omega_1} \)-equivalent} (i.e.\ they satisfy the same \( \LL_{\omega_2 \omega_1} \)-sentences) by~\cite[Theorem 9.26]{Vaananen:2011tg}. 
Since \( X \) and \( Y \) are non-isomorphic but \( \LL_{\omega_2 \omega_1} \)-equivalent, none of \( X \) and \( Y \) admits an \( \LL_{\omega_2 \omega_1} \)-Scott sentence. 
\end{example}

Working in \( \ZFC \) and assuming \( \kappa^{< \kappa} = \kappa > \omega \), easier counterexamples may be isolated. 
However, the construction in Example~\ref{xmp:Scottsentence} has the further merit of being carried out in \( \ZF \) and without further assumptions on \( \kappa = \omega_1 \) --- it is easy to check that when applying to our setup the (proofs of the) relevant results from~\cite{Vaananen:2011tg} no choice is needed, as \( X , Y \) have well-ordered domains.
In particular, Example~\ref{xmp:Scottsentence} works also in models of \( \ZFC + \neg \CH \).

\subsection{More absoluteness results} \label{subsec:moreabsoluteness}

Continuing the work in Section~\ref{subsec:someabsoluteness} on absoluteness of the definition of the map \( f_T \) from~\eqref{eq:f_T}, we are now going to observe that also the \( \LL_{\kappa^+ \kappa} \)-sentence \( \upsigma_T \) from Corollary~\ref{cor:saturation} and the map \( h_T \) from Definition~\ref{def:defh_T} are (essentially) absolute between transitive models of \( \ZF \) containing all the relevant parameters. 
These results, together with those from Section~\ref{subsec:someabsoluteness}, will be crucially used in Section~\ref{subsec:absolute}.

Given an arbitrary transitive model \( M \) of \( \ZF \), \( \kappa \in \Cn^{M} \), and \( T \in (\TT_\kappa)^M \), carry out the constructions from Section~\ref{sec:invariantlyuniversal} inside \( M \) (this is possible by our choice of \( M \)), i.e.\ let 
\[ 
g^M \coloneqq (g)^M \colon (\Mod^\kappa_\Uppsi)^M \to ( \pre{\pre{< \omega}{ 2} \times \kappa}{2})^M 
\] 
be defined as in~\eqref{eq:defg}--\eqref{eq:defgwithformula}, and let 
\[ 
\upsigma^M_T \coloneqq (\upsigma_T)^M \in (\LL_{\kappa^+ \kappa})^M 
\] 
be the sentence obtained in the proof of Corollary~\ref{cor:saturation} (where all objects involved are now computed in \( M \)). 
Note that by Corollary~\ref{cor:saturation}, which holds in \( M \), we have that \( ( \Mod^\kappa_{\upsigma^M} )^M \) is the saturation of the range of the map \( f^M_T \) from Section~\ref{subsec:someabsoluteness} (as computed in \( M \)), and that \( M \models \text{``}\Mod^\kappa_{\upsigma^M} = ( g^M )^{-1} ( \ran ( g^M \circ f_T^M ) ) \text{''} \). 
With this notation, we then get the following absoluteness result for the sentence \( \upsigma_T \).
 
\begin{proposition} \label{prop:upsigma_Tisabsolute}
Let \( M_0 \), \( M_1 \) be transitive models of \( \ZF \), and let \( \kappa \) and \( T \) be such that \( \kappa \in \Cn^{M_i} \) and \( T \in (\TT_\kappa)^{M_i} \) (for \( i = 0,1 \)).
Then \( \upsigma^{M_0}_T = \upsigma^{M_1}_T \).
\end{proposition}

\begin{proof} 
First notice that the set \( \forcing{Fn} ( \pre{< \omega}{ 2} \times \kappa , 2 ; \omega ) \) is absolute for transitive models of \( \ZF \) containing \( \kappa \). 
Since the elements of this set determine (both in \( M_0 \) and in \( M_1 \)) the canonical basis \( \mathcal{B}_p = \mathcal{B}_p ( \pre{\pre{< \omega}{ 2} \times \kappa}{2} ) \) for the product topology \( \tau_p \) on \( \pre{\pre{< \omega}{ 2} \times \kappa}{2} \), by the way we defined \( \upsigma_T^{M_i} \) it is enough to show that for every \( s \in \forcing{Fn} ( \pre{< \omega}{ 2} \times \kappa , 2 ; \omega ) \)
\begin{equation} \label{eq:iffabsolute} 
M_0 \models \text{``}\ran(g^{M_0} \circ f_T^{M_0}) \cap \Nbhd^A_s \neq \emptyset \text{''} \qquad \text{if and only if} \qquad M_1 \models \text{``} \ran(g^{M_1} \circ f_T^{M_1}) \cap \Nbhd^A_s \neq \emptyset \text{''},
\end{equation}
where \( A \coloneqq \pre{< \omega}{ 2} \times \kappa \) (see Section~\ref{subsubsec:spacesoftypekappa}). 
In fact, in the proof of Corollary~\ref{cor:saturation} we have set \( \upsigma_T^{M_i} \coloneqq \Uppsi \wedge \neg \upsigma^{M_i}_{V_i} \) with \( V_i = (\pre{A}{2})^{M_i} \setminus \ran(g^{M_i} \circ f^{M_i}_T) \), where \( \upsigma^{M_i}_{V_i} \) is computed in \( M_i \) according to the proof of Lemma~\ref{lem:open}. 
By inspecting such proof, it is easy to check that the \( (\LL_{\kappa^+ \kappa})^{M_i} \)-sentence \( \upsigma^{M_i}_{V_i} \) only depends on the set of \( s \in \forcing{Fn} ( \pre{< \omega}{ 2} \times \kappa , 2 ; \omega ) \) for which \( M_i \models \Nbhd^A_s \subseteq V_i \): thus if the equivalence in~\eqref{eq:iffabsolute} is true, this set is the same when computed in \( M_0 \) or \( M_1 \), whence 
\[ 
\upsigma_T^{M_0} = \Uppsi \wedge \neg \upsigma^{M_0}_{V_0} = \Uppsi \wedge \neg \upsigma^{M_1}_{V_1} = \upsigma_T^{M_1} ,
\]
as desired.

We now prove~\eqref{eq:iffabsolute}, starting with the implication from left to right. 
By Lemma~\ref{lem:continuous} (which holds in \( M_0 \)), the function \( g^{M_0} \circ f^{M_0}_T \) is continuous in \( M_0 \). 
Therefore, if \( M_0 \models \text{``} \ran(g^{M_0} \circ f^{M_0}_T) \cap \Nbhd^A_s \neq \emptyset \text{''} \), then there is \( t \in \pre{< \omega}{2} \) such that \( M_0 \models \text{``} \Nbhd^\omega_t \subseteq (g^{M_0} \circ f^{M_0}_T)^{-1}(\Nbhd^A_s) \text{''} \), in particular \( M_0 \models \text{``}(g^{M_0} \circ f^{M_0}_T)(\bar{x}) \in \Nbhd^A_s \text{''} \) with \( \bar{x} \coloneqq t {}^\smallfrown{} 0^{(\omega)} \). 
Since \( (\pre{<\omega}{2})^{M_0} = (\pre{< \omega}{2})^{M_1} \), it follows that \( \bar{x} \in (\pre{\omega}{2})^{M_0} \cap (\pre{\omega}{2})^{M_1} \) and hence \( f^{M_0}_T(\bar{x}) = f_T^{M_1}(\bar{x}) \) by Fact~\ref{prop:f_Tdefinable}\ref{prop:f_Tdefinable-b}. 
By Lemma~\ref{lem:range} (which holds in both \( M_0 \) and \( M_1 \)) and the definition of the map \( g \) given in~\eqref{eq:defg}--\eqref{eq:defgwithformula}, both \( g^{M_0} \) and \( g^{M_1} \) are defined on \( X \coloneqq f^{M_0}_T(\bar{x}) = f_T^{M_1}(\bar{x}) \) and \( g^{M_0}(X) = g^{M_1}(X) \). 
Therefore \( M_1 \models \text{``}(g^{M_1} \circ f_T^{M_1})(\bar{x}) \in \Nbhd^A_s \text{''} \) as well, witnessing \( M_1 \models \text{``} \ran(g^{M_1} \circ f^{M_1}_T) \cap \Nbhd^A_s \neq \emptyset \text{''} \).

The reverse implication is proved in a similar way by switching the role of \( M_0 \) and \( M_1 \). 
\end{proof}

Let \( M \), \( \kappa \), \( T \), and \( \upsigma_T^M \) be as in the paragraph preceding Proposition~\ref{prop:upsigma_Tisabsolute}. 
Applying Definition~\ref{def:defh_T} in \( M \), we also get the map 
 \[ 
 h_T^M \coloneqq (h_T)^M \colon (\Mod^\kappa_{\upsigma^M_T})^M \to (\pre{\omega}{2})^M .
 \] 
Notice that since Corollary~\ref{cor:inverse} holds in \( M \), the map \( h^M_T \) reduces, in the sense of \( M \), the embeddability relation \( (\embeds^\kappa_{\upsigma^M})^M \) to \( R^M \coloneqq ( \PROJ \body{T})^M \) (as long as \( R^M \) is a quasi-order in \( M \)). 

\begin{proposition} \label{prop:h_Tisabsolute} 
\begin{enumerate-(a)}
\item \label{prop:h_Tisabsolute-a}
There is an \( \LST \)-formula \( {\sf\Psi}_{h_T} ( x_0 , x_1 , z_0 , z_1 ) \) such that for every transitive model \( M \) of \( \ZF \) with \( \kappa \in \Cn^M \) and \( T \in (\TT_\kappa)^M \), the graph of \( h^M_T \) is defined in \( M \) by \( {\sf \Psi}_{h_T}(x_0,x_1,\kappa,T) \), that is: for every \( X \in (\Mod^\kappa_{\upsigma^M_T})^M \) and \( x \in (\pre{\omega}{2})^M \)
\[ 
h^M_T(X) = x \iff M \models {\sf \Psi}_{h_T}[X,x,\kappa,T].
 \] 
\item \label{prop:h_Tisabsolute-b}
Let \( M_0 \), \( M_1 \) be transitive models of \( \ZF \), and let \( \kappa \) and \( T \) be such that \( \kappa \in \Cn^{M_i} \) and \( T \in (\TT_\kappa)^{M_i} \) (for \( i = 0,1 \)), so that \( \upsigma^{M_0}_T = \upsigma^{M_1}_T \) by Proposition~\ref{prop:upsigma_Tisabsolute}. 
If 
\begin{equation}\label{eq:prop:h_Tisabsolute}
 ( \pre{\kappa}{\kappa})^{M_0} \subseteq (\pre{\kappa}{\kappa})^{M_1} ,
\end{equation}
then setting \( \upsigma \coloneqq \upsigma^{M_0}_T = \upsigma^{M_1}_T \) we have that for every \( X \in (\Mod^\kappa_\upsigma)^{M_0} \cap (\Mod^\kappa_\upsigma)^{M_1} \)
\begin{equation} \label{eq:coherenceh_T}
h^{M_0}_T(X) = h_T^{M_1}(X).
\end{equation}
\end{enumerate-(a)}
\end{proposition}
 
 \begin{proof}
For part~\ref{prop:h_Tisabsolute-a} it is enough to observe that equation~\eqref{eq:defh_T} is rendered by the \( \LST \)-formula \( {\sf\Psi}_{h_T} ( x_0 , x_1 , z_0 , z_1 ) \)
\[ 
{x_0 \in \Mod^{z_0}_\LL} \wedge {x_1 \in \pre{\omega}{2}} \wedge \EXISTS{ i \in \Sym(z_0)} \FORALL{ x_2 } ( {\sf\Psi}_{f_T} ( x_1 , x_2 , z_0 , z_1 ) \implies i \text{ witnesses } x_2 \cong x_0 ) ,
 \] 
where \( {\sf\Psi}_{f_T} ( x_1 , x_2 , z_0 , z_1 ) \) is as in Fact~\ref{prop:f_Tdefinable} and ``\( i \) witnesses \( x_2 \cong x_0 \)'' stands for
\begin{equation} \label{eq:isom}
\forall \alpha , \beta < \kappa \, \bigl ( \alpha \edge^{x_2} \beta \iff i ( \alpha ) \edge^{x_0} i ( \beta )\bigr ) .
\end{equation}

We now prove part~\ref{prop:h_Tisabsolute-b}. 
Let \( X \in (\Mod^\kappa_\upsigma)^{M_0} \cap (\Mod^\kappa_\upsigma)^{M_1} \) and set \( x \coloneqq h^{M_1}_T(X) \in (\pre{\omega}{2})^{M_1} \) and \( x' \coloneqq h^{M_0}_T(X) \in (\pre{\omega}{2})^{M_0} \subseteq (\pre{\omega}{2})^{M_1} \) (the latter inclusion follows from our assumption \( (\pre{\kappa}{\kappa})^{M_0} \subseteq (\pre{\kappa}{\kappa})^{M_1} \)). 
Then there is \( i \in (\Sym(\kappa))^{M_0} \) such that \( M_0 \models \text{``} i \text{ witnesses } f^{M_0}_T(x') \cong X \text{''} \). 
Notice that since \( (\Sym ( \kappa ) )^{M_0} \subseteq ( \Sym( \kappa ) ) ^{M_1} \) (again by our assumption \( (\pre{\kappa}{\kappa})^{M_0} \subseteq (\pre{\kappa}{\kappa})^{M_1} \)), then \( i \in ( \Sym( \kappa ) ) ^{M_1} \) as well.
Since \( f_T^{M_0}(x') = f_T^{M_1} ( x' ) \) by Fact~\ref{prop:f_Tdefinable}\ref{prop:f_Tdefinable-b}, and since \( f_T^{M_1} ( x' ) \in M_1 \), the sentence ``\( i \) is an isomorphism between the structures \( f^{M_0}_T(x') \) and \( X \in \Mod^\kappa_\LL \)'' can be formalized by the \( \Delta_0 \) \( \LST \)-formula \emph{with parameters \( i \), \( f^{M_0}_T ( x' ) \), \( X \), and \( \kappa \) in \( M_0 \cap M_1 \)} provided in~\eqref{eq:isom} (where we replace \( x_2 \) and \( x_0 \) by, respectively, \( f^{M_0}_T ( x' ) \) and \( X \)), and is therefore absolute between \( M_0 \) and \( M_1 \). 
Therefore \( M_1 \models \text{``} f_T^{M_1} ( x' ) \cong X \text{''}\). 
By our choice of \( x \) we also have \( M_1 \models \text{``} f_T^{M_1} ( x ) \cong X \text{''} \), and therefore it follows from Theorem~\ref{th:graphs}\ref{th:graphs-b}, which holds in the \( \ZF \)-model \( M_1 \), that \( x' = x \). 
Thus~\eqref{eq:coherenceh_T} is satisfied and we are done. 
\end{proof}

\begin{remark} \label{rmk:newrmkabsolute}
Condition~\eqref{eq:prop:h_Tisabsolute} is artificial and can be removed in many cases.
For example, if \( M_0 , M_1 \) satisfy all hypotheses of Proposition~\ref{prop:h_Tisabsolute}\ref{prop:h_Tisabsolute-b} except for~\eqref{eq:prop:h_Tisabsolute}, but there is a transitive model \( N \) of \( \ZF \) containing both \( M_0 , M_1\), then the conclusion still holds.
In fact, since \( ( \pre{ \kappa }{ \kappa } )^{M_i} \subseteq ( \pre{ \kappa }{ \kappa } )^{N} \) then \( h_T^{M_i} ( X ) = h_T^{N} ( X )\) by Proposition~\ref{prop:h_Tisabsolute}\ref{prop:h_Tisabsolute-b}, whence \( h_T^{M_0} ( X ) = h_T^{M_1} ( X ) \).
In particular this applies when 
\begin{itemize}
\item
\( M_0 , M_1 \subseteq \Vv \) are transitive models of \( \ZF \) --- take \( N = \Vv \);
\item
\( M_0 , M_1 \) are generic extensions of \( \Vv \) --- if \( M_i = \Vv [ G_i ] \) with \( G_i \subseteq \forcing{P}_i \) generic over \( \Vv \), then take \( N = \Vv [ G ] \) with \( G \subseteq \forcing{P}_0 \times \forcing{P}_1 \) generic over \( \Vv \).
\end{itemize}
\end{remark}

\section{An alternative approach}\label{sec:alternativeapproach}
As explained at the beginning of Section~\ref{sec:mainconstruction}, two different approaches have been employed in the literature to show that the embeddability relation on countable structures is invariantly universal, namely the approaches~\ref{en:approach1} and~\ref{en:approach2} briefly described on page~\pageref{en:approach1}. 
In Sections~\ref{sec:mainconstruction}--\ref{sec:invariantlyuniversal} we successfully followed approach~\ref{en:approach2}, and this will provide generalizations to uncountable cardinals \( \kappa \) of Theorems~\ref{th:LouveauRosendal} and~\ref{th:mottorosfriedman} (see Sections~\ref{sec:definablecardinality} and~\ref{sec:applications}). 
However, our proof of Theorem~\ref{th:main} does not yield a generalization of Theorem~\ref{th:mottorosfriedman2}: the reason is that the \( \LL_{\kappa^+ \kappa} \)-sentence \( \upsigma_T \) from Corollary~\ref{cor:saturation} is quite complex, and hence \( \Mod^\kappa_{\upsigma_T} \) is usually far from being a \( \kappa+1 \)-Borel subset of \( \Mod^\kappa_\LL \) (unless we assume \( \AC + \kappa^{< \kappa} = \kappa \), see Section~\ref{subsec:LopezEscobar} and, in particular, Remark~\ref{rmk:formulaborel}). 

In this section we are going to show that although approach~\ref{en:approach1} forces us to consider a less natural kind of structures, it allows us to further obtain some sort of generalization of Theorem~\ref{th:mottorosfriedman2} to uncountable \( \kappa \)'s: this is essentially because with this alternative approach we will be able to associate to each \( \kappa \)-Souslin quasi order \( R = \PROJ{\body{T}} \) a sentence \( \bar{\upsigma}_T \) in the \emph{bounded} logic \( \bar{\LL}^b_{\kappa^+ \kappa} \) still having the property that \( R \sim {\embeds}^\kappa_{\bar{\upsigma}_T} \), thus further obtaining that when \( \kappa \) is regular \( \Mod^\kappa_{\bar{\upsigma}_T} \) is a \( \kappa+1 \)-Borel subset of \( \Mod^\kappa_{\bar{\LL}} \) by Corollary~\ref{cor:formulaborel1}\ref{cor:formulaborel1-ii}.

The alternative construction we are going to provide consists of the following steps:
\begin{itemize}[leftmargin=1pc]
\item
we expand the \( \LL \)-structure \( \mathbb{G}_0 \in \CT_\kappa \) defined in Section~\ref{subsec:G_0andG_1} to a so-called \emph{ordered} combinatorial tree \( \bar{\mathbb{G}}_0 \) by interpreting the new symbol \( \order \) in the extended language \( \bar{\LL} = \{ \edge, \order \} \)\index[symbols]{106@\( \order \)} as a well-founded order on (the nodes of) \( \mathbb{G}_0 \);
\item
given \( S \in \Tr(2 \times \kappa) \), we directly add to \( \bar{\mathbb{G}}_0 \) some forks as in Section~\ref{subsec:G_T} (without first enlarging \( \bar{\mathbb{G}}_0 \) to an analogue of \( \mathbb{G}_1 \)), leaving the interpretation of \( \order \) unchanged.
\end{itemize} 

Let us now fix some further notation concerning the \( \bar{\LL} \)-structures that we are going to consider. 
Following~\cite{Friedman:2011cr}, a combinatorial tree with a supplementary transitive relation defined on a subset of its set of vertices will be called an \markdef{ordered combinatorial tree}.
This is a bit of a misnomer, since this extra relation need not be an ordering, although in what follows this will always be the case.
To formalize this concept, consider the expansion \( \bar{\LL} = \setLR{ \edge , \order } \) of the graph language \( \LL \) (see page~\pageref{pag:L}) with \( \order \) a binary relational symbol.
Then \( \OCT_\kappa \subseteq \Mod^\kappa_{\bar{\LL}} \) is the collection of all ordered combinatorial trees of size \( \kappa \) (up to isomorphism), i.e.~the set of all \( X = \seqLR{\kappa ; {\edge^X} , {\order^X} } \) such that \( \seqLR{ \kappa ; {\edge^X} } \) is a combinatorial tree and \( \order^X \) is a transitive relation. 
Formally, \( \OCT_\kappa \coloneqq \Mod^\kappa_{\upsigma_\OCT} \),\index[symbols]{OCT@\( \OCT_\kappa \)} where \( \upsigma_{\OCT} \) is the following \( \bar{\LL}_{\omega_1 \omega} \)-sentence axiomatizing ordered combinatorial trees:
\begin{equation}
\tag{$\upsigma_\OCT$} \upsigma_\CT \wedge \forall \V_0 \forall \V_1 \forall \V_2 \left ( {\V_0 \trianglelefteqslant \V_1} \wedge {\V_1 \trianglelefteqslant \V_2} \implies {\V_0 \trianglelefteqslant \V_2}\right ). \index[symbols]{sigmaOCT@\( \upsigma_{\OCT} \)}
\end{equation}

Similarly to the case of combinatorial trees, in order to simplify the notation we will abbreviate the embeddability relation \( \embeds \restriction \Mod^\kappa_{\upsigma_{\OCT}} \) and the isomorphism relation \( \cong \restriction \Mod^\kappa_{\upsigma_\OCT} \) (see page~\pageref{not:embeddability}) with \( \embeds_\OCT^\kappa \) and \( \cong_\OCT^\kappa \), respectively.

\subsection{Completeness}

To begin with, we slightly modify the construction presented in Section~\ref{sec:mainconstruction} by first adding a well-founded order on the vertices of the combinatorial tree \( \mathbb{G}_0 \) from Definition~\ref{def:G_0}. 
Such order is essentially obtained using in the obvious way the bijection \( \Code{\cdot} \colon \pre{< \omega}{\On} \to \On \) from~\eqref{eq:codingfinitesequenceofordinals}. 
More precisely, define 
\begin{equation}\label{eq:pound}
\# \colon \pre{< \omega}{\kappa} \uplus \setofLR{s^-}{\emptyset \neq s \in \pre{< \omega}{\kappa} } \to \kappa \index[symbols]{106@\( \# \)}
\end{equation}
by letting 
\[
\begin{split}
\# ( s ) &\coloneqq 2 \Code{ s }
\\
 \# ( s^- ) &\coloneqq 
 \begin{cases}
2 \Code{ s } - 1 & \text{if } \Code { s } \in \omega ,
\\
2 \Code { s } + 1&\text{otherwise.}
\end{cases}
\end{split}
\]
Letting \( s_\alpha \coloneqq \#^{-1} ( \alpha) \), one gets that \( s_\alpha \in \pre{< \omega}{ \kappa} \) (i.e.~\( s_\alpha \) is not of the form \( s^- \)) if and only if \( \alpha \) is even. 
The ordered combinatorial tree \( \bar{\mathbb{G}}_0 \) is then the expansion of the structure \( \mathbb{G}_0 = \seqLR{G_0 ; \edge^{G_0}} \) obtained by interpreting \( \order \) as the well-order induced on \( G_0 \) by the above map \( \# \). 
Formally:

\begin{definition}\label{def:barG_0}
\( \bar{\mathbb{G}}_0 \) is the \( \bar{\LL} \)-structure on 
\[ 
\bar{G}_0 \coloneqq G_0 = \pre{< \omega}{\kappa} \uplus \setofLR{s^-}{\emptyset \neq s \in \pre{< \omega}{\kappa} }
 \] 
with edge relation \( \edge^{\bar{\mathbb{G}}_0} \) defined by
\[ 
a \edge^{\bar{\mathbb{G}}_0} b \IFF a \edge^{\mathbb{G}_0} b
\] 
and order relation \( \order^{\bar{\mathbb{G}}_0} \) defined by
\[ 
a \order^{\bar{\mathbb{G}}_0} b \IFF \# ( a ) \leq \# ( b ) ,
\] 
where \( a,b \in \bar{G}_0 \), \( \edge^{\mathbb{G}_0} \) is as in Definition~\ref{def:G_0}, and \( \# \) is the map from~\eqref{eq:pound}.
\end{definition}

Notice that \( \order^{\bar{\mathbb{G}}_0} \) is isomorphic to \( \kappa \), so \( \bar{\mathbb{G}}_0 \) is rigid.
We next obtain the structures \( \bar{\mathbb{G}}_S \) (for \( S \in \Tr ( 2 \times \kappa ) \)) by joining \( \bar{\mathbb{G}}_0 \) and the forks \( F_{ u , s } \) from~\eqref{eq:F_{u,s}} (for every \( ( u , s ) \in S \)) via the identification of \( s \in \bar{G}_0 \) with the vertex \( ( u , s , \emptyset ) \) of \( F_{ u , s } \) --- note that this construction is exactly the same described in Section~\ref{subsec:G_T} except that now we are avoiding the enlargement to the structure \( \mathbb{G}_1 \). 

\begin{definition} \label{def:barG_T}
Let \( S \in \Tr ( 2 \times \kappa ) \). 
Then \( \bar{\mathbb{G}}_S \)\index[symbols]{Graph3@\( \bar{\mathbb{G}}_S \)} is the \( \bar{\LL} \)-structure on 
\[ 
\bar{G}_S \coloneqq \bar{G}_0 \uplus \setofLR{ ( u , s , w ) }{ ( u , s ) \in S \wedge ( u , s , w ) \in F_{ u , s } \wedge w \neq \emptyset}
 \] 
with edge relation \( \edge^{\bar{\mathbb{G}}_S} \) defined by
\[ 
a \edge^{\bar{\mathbb{G}}_S} b \IFF a \edge^{\mathbb{G}_S} b
 \] 
and order relation \( \order^{\bar{\mathbb{G}}_S} \) defined by
\[ 
a \order^{\bar{\mathbb{G}}_S} b \iff a , b \in \bar{G}_0 \wedge a \order^{\bar{\mathbb{G}}_0} b,
 \] 
where \( a , b \in \bar{G}_S \) and \( \edge^{\mathbb{G}_S} \) is the edge relation on \( \mathbb{G}_S \) defined in Section~\ref{subsec:G_T}.
\end{definition}

Notice that the above is well-defined because \( \bar{G}_S \subseteq G_S \): in fact, using the notation introduced in~\eqref{eq:substructuresofG_T-a}--\eqref{eq:substructuresofG_T-g} and Remark~\ref{rmk:substructuresofG_T}
\[ 
\bar{G}_S = \mathrm{Seq} \uplus \mathrm{Seq}^- \uplus \bigcup \setofLR{\mathrm{F}'_{u,s}}{(u,s) \in S} = G_S \setminus \widehat{\mathrm{Seq}}.
 \] 
Each structure \( \mathbb{G}_S \) is then further identified in a canonical way with its copy on \( \kappa \) (which is thus an element of \( \OCT_\kappa \)), using the bijections \( \Code{\cdot} \), \( \op{ \cdot}{\cdot } \), and \( e_u \) from~\eqref{eq:codingfinitesequenceofordinals}, \eqref{eq:Hessenberg}, and~\eqref{eq:e_u} in the obvious way. 

Finally, given an infinite cardinal \( \kappa \) and \( T \in \TT_\kappa \) we let \( \bar{f}_T \) be the variant of the function \( f_T \) from~\eqref{eq:f_T} obtained by composing the map \( \Sigma_T \colon \pre{\omega}{2} \to \Tr ( 2 \times \kappa ) \) from Definition~\ref{def:Sigma_T} with the modified map \( S \mapsto \bar{\mathbb{G}}_S \), that is: 
\begin{equation}\label{eq:defbarf_T}
\bar{f}_T \colon \pre{ \omega}{2} \to \OCT_\kappa , \qquad x \mapsto \bar{\mathbb{G}}_{\Sigma_T ( x )}. \index[symbols]{fT@\( \bar{f}_T \)}
\end{equation}

\begin{theorem} \label{th:OCT}
Let \( \kappa \) be an infinite cardinal and \( T \in \TT_\kappa \). 
If \( R = \PROJ \body{ T } \) is a quasi-order, then the map \( \bar{f}_T \) defined in~\eqref{eq:defbarf_T} is such that:
\begin{enumerate-(a)}
\item \label{th:OCT-1}
\( \bar{f}_T \) reduces \( R \) to the embeddability relation \( \embeds_\OCT^\kappa \);
\item \label{th:OCT-2}
\( \bar{f}_T \) reduces \( = \) on \( \pre{ \omega}{2} \) to the isomorphism relation \( \cong_\OCT^\kappa \).
\end{enumerate-(a)}
In particular, \( \embeds^\kappa_\OCT \) is complete for \( \kappa \)-Souslin quasi-orders.
\end{theorem}

\begin{proof}
For part~\ref{th:OCT-1} just check that the same proof of Theorem~\ref{th:graphs}\ref{th:graphs-a} goes through also with the new construction with the following minor modifications:
\begin{itemize}[leftmargin=1pc]
\item
in the forward direction, to get the embedding between \( \bar{f}_T(x) = \bar{\mathbb{G}}_{\Sigma_T(x)} \) and \( \bar{f}_T(y) = \bar{\mathbb{G}}_{\Sigma_T(y)} \) when \( x \mathrel{R} y \) we of course just consider the restriction of \( i \) to \( \bar{G}_{\Sigma_T(x)} \subseteq G_{\Sigma_T(x)} \) instead of the full embedding \( i \colon \mathbb{G}_{\Sigma_T(x)} \to \mathbb{G}_{\Sigma_T(y)} \) described in the original proof;
\item
for the backward direction, we use Remark~\ref{rmk:partialembedding} and the fact that \( \bar{\mathbb{G}}_{\Sigma_T(x)} \) and \( \bar{\mathbb{G}}_{\Sigma_T(y)} \) are just \( \bar{\LL} \)-expansions of \( \mathbb{G}_{\Sigma_T(x)} \restriction (G_0 \cup \mathrm{F}(\mathbb{G}_{\Sigma_T ( x )})) \) and \( \mathbb{G}_{\Sigma_T(y)} \restriction (G_0 \cup \mathrm{F}(\mathbb{G}_{\Sigma_T ( y )})) \), respectively.
\end{itemize}

For part~\ref{th:OCT-2}, fix an isomorphism \( j \) between \( \bar{f}_T ( x ) = \bar{\mathbb{G}}_{\Sigma_T ( x )} \) and \( \bar{f}_T ( y ) = \bar{\mathbb{G}}_{\Sigma_T(y)} \). 
Then \( j ( \bar{G}_0 ) = \bar{G}_0 \) because the vertices in \( \bar{G}_0 \subseteq \bar{G}_{\Sigma_T(x)}, \bar{G}_{\Sigma_T ( y )} \) are the unique ones which are in the domain of the orders \( \order^{\bar{\mathbb{G}}_{\Sigma_T(x)}} \) and \( \order^{\bar{\mathbb{G}}_{\Sigma_T(y)}} \), respectively. 
Since the relational symbol \( \order \) is interpreted in both \( \bar{\mathbb{G}}_{\Sigma_T(x)} \) and \( \bar{\mathbb{G}}_{\Sigma_T(y)} \) as the same well-order on \( \bar{G}_0 \) (of the same order type \( \kappa \)), we get that \( j \restriction \bar{G}_0 \) is the identity function. 
Arguing as in Theorem~\ref{th:graphs}\ref{th:graphs-b}, we then get that this implies \( \Sigma_T(x) = \Sigma_T(y) \), and hence that \( x=y \) by injectivity of \( \Sigma_T \) (Lemma~\ref{lem:max}\ref{lem:max-iii}).
\end{proof}

\subsection{Invariant universality} 

Henceforth we again fix an \emph{uncountable} cardinal \( \kappa \), and show that \( \embeds^\kappa_{\OCT} \) is invariantly universal for \( \kappa \)-Souslin quasi-orders. 
Following Section~\ref{sec:invariantlyuniversal}, we will first provide an \( \bar{\LL}_{\kappa^+ \kappa} \)-formula \( \bar{\Uppsi} \) describing the common parts of the \( \bar{\mathbb{G}}_S \) (for \( S \in \Tr(2 \times \kappa) \)), and then classify the structures in \( \Mod^\kappa_{\bar{\Uppsi}} \) up to isomorphism using the elements of \( \pre{\pre{<\omega}{2} \times \kappa}{2} \) as invariants. 
The main improvement of the new construction is that \( \bar{\Uppsi} \) will be a sentence in the bounded logic \( \bar{\LL}^b_{\kappa^+ \kappa} \), so that when \( \kappa \) is regular \( \Mod^\kappa_{\bar{\Uppsi}} \) is an effective \( \kappa+1 \)-Borel set by Corollary~\ref{cor:formulaborel1}\ref{cor:formulaborel1-ii}. 
As in Section~\ref{sec:invariantlyuniversal}, to simplify the presentation we will freely use metavariables and consider (infinitary) conjunctions and disjunctions over (canonically) well-orderable sets of formul\ae{} of size \( \leq \kappa \). 

Given an ordinal \( \alpha < \kappa \) and a sequence of variables \( \seqofLR{ x_\beta }{ \beta < \alpha } \), let \( \bar\uprho ( \seqofLR{ x_\beta }{ \beta < \alpha } ) \) be the \( \bar{\LL}^0_{ \card{\alpha}^+ \omega} \)-formula 
\begin{equation}\tag{$\bar\uprho$}
\bigwedge_{\beta < \gamma < \alpha} \left ( x_\beta \nequals x_\gamma \AND x_\beta \order x_\gamma \AND \neg(x_\gamma \order x_\beta) \right ) .
\end{equation} 
Notice that for any \( \bar{\LL} \)-structure \( X \) with domain \( \kappa \) and any assignment \( s \in \pre{\alpha}{ \kappa} \), \( X \models \bar\uprho [ s ] \) if and only if \( s \) is injective and \( \ran ( s ) \) is well-ordered by \( \order ^X \) in order type \( \alpha \).
Now let \( \bar\Upphi_1 \) be the \( \bar{\LL}^b_{\kappa^+ \kappa} \)-sentence given by the conjunction of the following:
\begin{enumerate-(i)}
 \item \label{en:Phi0-i}
\( \forall x , y \bigl [ \exists z \left ( x\order z \vee z \order x \right ) \wedge \exists z \left ( y\order z \vee z \order y \right) \IMPLIES x\order y \vee y\order x \bigr ] \);
 \item \label{en:Phi0-ii}
 \( \forall x , y , z ( x\order y \wedge y\order z \implies x\order z ) \);
 \item \label{en:Phi0-iii}
 \( \forall x , y ( x\order y \wedge y \order x \implies x \equals y ) \);
 \item \label{en:Phi0-iv}
 \( \neg \exists \seqofLR{ x_n}{ n < \omega } [ \bigwedge_{n<m<\omega} ( x_m \nequals x_n \wedge x_m\order x_n ) ] \);
\item \label{en:Phi0-v}
\( \bigwedge_{\alpha < \kappa} \exists \seqofLR{ x_\beta }{ \beta < \alpha } \bar\uprho ( \seqofLR{ x_\beta }{ \beta < \alpha }) \);
\item \label{en:Phi0-vi}
\( \neg \exists x \bigwedge_{\alpha < \kappa} \exists \seqofLR{ y_\beta }{ \beta < \alpha } \bar\uprho ( \seqofLR{ y_\beta }{ \beta < \alpha } \conc x ) \).
\end{enumerate-(i)}
If \( X \) is a structure as above and satisfies~\ref{en:Phi0-i}--\ref{en:Phi0-iv} then \( \order ^X \) is a well-order.%
\footnote{We don't need to assume \( \DC \) to have the equivalence between ``well-foundness'' and ``absence of descending chains'' because in our situation the domain of the structure \( X \) is always assumed to be the cardinal \( \kappa \), which carries a natural well-order.}
If moreover it satisfies~\ref{en:Phi0-v} the length of \( \order ^X \) is \( \geq \kappa \), and if it satisfies~\ref{en:Phi0-vi} the length is \( {<} \kappa + 1 \).
Therefore if \( X \models \bar\Upphi_1 \) then \( \order ^X \) is a well-order of length \( \kappa \), possibly defined just on a subset of the domain of \( X \).

For every ordinal \( \alpha < \kappa \) let \( \bar\uprho_ \alpha ( x ) \) be the \( \bar{\LL}^b_{\kappa^+ \kappa} \)-formula 
\begin{equation}\tag{$\bar\uprho_\alpha$}
 \exists \seqofLR{x_ \beta }{ \beta < \alpha } \bar\uprho ( \seqofLR{x_ \beta }{ \beta < \alpha } \conc x ) \wedge \neg \exists \seqofLR{x_ \beta }{ \beta < \alpha + 1 } \bar\uprho ( \seqofLR{x_ \beta }{ \beta < \alpha + 1} \conc x ) .
\end{equation}
stating that \( x \) is the \( \alpha \)-th element in the order \( \order \), at least whenever \( \order \) is (interpreted in) a well-order of length greater than \( \alpha \).
Given ordinals \( \alpha , \beta < \kappa \) and two variables \( x_0 , x_1 \), let \( \bar\uppsi_{\alpha , \beta } ( x_0 , x_1 ) \) be the formula \( x_0 \edge x_1 \) if \( s_\alpha \edge^{\bar{\mathbb{G}}_0} s_{\beta} \) or \( \neg ( x_0 \edge x_1 ) \) otherwise, where \( s_\alpha \coloneqq \#^{-1}(\alpha) \) and \( \# \) is as in~\eqref{eq:pound}. 
Let \( \bar\Upphi_2 \) be the \( \bar{\LL}^b_{ \kappa^+ \kappa} \)-sentence:
\begin{equation}\tag{\( \bar\Upphi_2 \)}
\bigwedge_{\alpha < \beta < \kappa} \forall x , y \left (\bar\uprho_ \alpha ( x ) \wedge \bar\uprho_ { \beta} ( y ) \implies \bar\uppsi_{ \alpha , \beta } ( x , y ) \wedge \bar\uppsi_{ \beta , \alpha } ( y , x ) \right ) .
\end{equation} 
Notice that if an \( \bar{\LL} \)-structure \( X \) satisfies \( \bar\Upphi_1 \wedge \bar\Upphi_2 \), then its restriction to the field of \( \order ^X \) is isomorphic to \( \bar{\mathbb{G}}_0 \), and moreover such an isomorphism is unique by the rigidity of \( \bar{\mathbb{G}}_0 \).

For every \( u \in \pre{< \omega}{2} \), let \( \bar{F}_u \) be the \( \bar{\LL} \)-expansion of \( F_u \) in which \( \emptyset \order^{\bar{F}_u} \emptyset \) and no other \( \order^{\bar{F}_u} \)-relation holds. 
We call \( \emptyset \) the \markdef{root} of \( \bar{F}_u \). 
Notice that \( \bar{F}_u \) still satisfies analogues of~\eqref{eq:fixempty}--\eqref{eq:=iffiso}. 
In particular, for any \( \bar{\LL} \)-structure \( X \cong \bar{F}_u \) the isomorphism between \( X \) and \( \bar{F}_u \) is unique, and we can thus unambiguously call root of \( X \) the unique element of \( X \) which is mapped by such isomorphism to the root \( \emptyset \) of \( \bar{F}_u \).

Given \( n , m \in \omega \), \( u \in \pre{< \omega}{ 2} \), and two variables \( x_0,x_1 \), let \( \bar{\upchi}^u_{ n , m } ( x_0 , x_1 ) \) be the \( \bar{\LL} \)-formula \( x_0 \edge x_1 \) if \( e_u^{-1}(n) \edge^{F_u} e_u^{-1}(m) \) and \( \neg ( x_0 \mathrel{\mathbf{E}} x_1 ) \) otherwise, where \( e_u \) is the bijection of~\eqref{eq:e_u}. 
Then \( \bar\upchi_u ( \seqofLR{ x_n }{ n \in \omega }) \) denotes the \( \bar{\LL}^0_{\omega_1 \omega} \)-formula
\begin{multline} \tag{$\bar{\upchi}_u$}
 \bigwedge_{ n < m < \omega} ( x_n \nequals x_m ) \wedge \bigwedge_{ n , m < \omega} \bar{\upchi}^u_{ n , m } ( x_n , x_m ) 
 \\ 
 \wedge ( x_0 \order x_0 ) \wedge \bigwedge_{\substack{n , m < \omega \\ n > 0}} ( \neg ( x_n \order x_m ) \wedge \neg ( x_m \order x_n ) ).
\end{multline}
 Notice that for any \( \bar{\LL} \)-structure \( X \) and any \( \seqofLR{a_n }{ n \in \omega } \in \pre{\omega}{ X} \), 
\begin{equation}\label{eq:barchi_u}
 \begin{split}
X \models \bar\upchi_u [ \seqofLR{a_n }{ n \in \omega } ] & \iff \text{the substructure of \( X \) with domain } \setofLR{a_n}{n \in \omega } 
 \\
 &\hphantom{{}\iff{}}\text{is isomorphic to } \bar{F}_u \text{ (via the map \( a_n \mapsto e_u^{-1}(n) \))}.
\end{split}
\end{equation}
Moreover such isomorphism must again be unique by~\eqref{eq:F_urigid}.
Let now \( \bar\Upphi_3 \) be the following \( \mathcal{L}^b_{\omega_1 \omega_1} \)-sentence:
\begin{equation} \tag{$\bar\Upphi_3$}
\FORALL{ x} \bigg[ \neg ( x \order x ) \implies \bigvee_{u \in \pre{< \omega}{ 2}} \exists \seqofLR{ y_n }{ n < \omega } \bigg (\bar\upchi_u(\seqofLR{y_n}{n < \omega}) \AND \bigvee_{0 \neq i < \omega} x \equals y_i \bigg ) \bigg ] . 
 \end{equation} 
This means that if an \( \bar{\LL} \)-structure \( X \) satisfies \( \bar\Upphi_1 \wedge \bar\Upphi_3 \), then each element of \( X \) which is not in the field of \( \order ^X \) belongs to some substructure which is isomorphic to \( \bar{F}_u \) (for some \( u \in \pre{< \omega}{ 2} \)).

Let now \( \bar\Upphi_4 \) be the following \( \bar{\LL}^b_{\omega_1 \omega_1} \)-sentence:
\begin{multline}\tag{$\bar\Upphi_4$}
 \FORALL{ x} \FORALL{ \seqofLR{ y_n }{ n < \omega }} \bigg[ \bigvee_{u \in \pre{< \omega}{ 2}} \bar\upchi_u ( \seqofLR{ y_n }{ n < \omega } ) \wedge \bigwedge_{ n < \omega} ( x \nequals y_n ) 
 \\
 \implies \bigwedge_{0 \neq n < \omega} ( \neg ( x \edge y_n ) \wedge \neg ( x \order y_n ) \wedge \neg ( y_n \order x ) )\bigg ] . 
\end{multline}
Note that if \( X \) is an \( \bar{\LL} \)-structure such that \( X \models \upsigma_{\OCT} \wedge \bar\Upphi_4 \), then any of its substructures which is isomorphic to some \( \bar{F}_u \) is ``isolated'' from the rest of \( X \), meaning that each element of such substructure (except for its root) is neither \( \edge^X \)-related nor \( \order ^X \)-related to any other element of \( X \) which does not belong to said substructure.

Let \( \bar\Upphi_5 \) be the following \( \bar{\LL}^b_{\omega_1 \omega_1} \)-sentence:
\begin{multline} \tag{$\bar\Upphi_5$}
\forall \seqofLR{ x_n }{ n < \omega } \forall \seqofLR{ y_m }{ m < \omega } \bigg [ \bigvee_{u \in \pre{< \omega}{ 2}} \bar\upchi_u ( \seqofLR{ x_n }{ n < \omega } ) \wedge \bigvee_{v \in \pre{< \omega}{ 2 }} \bar\upchi_v ( \seqofLR{ y_{m} }{ m < \omega } ) 
\\
\implies \bigwedge_{n < \omega} ( x_n \equals y_n ) \vee \bigwedge_{\substack{n , m < \omega \\ ( n , m ) \neq ( 0 , 0)}} ( x_n \nequals y_m ) \bigg ] .
\end{multline}
This means that if an \( \bar{\LL} \)-structure \( X \) satisfies \( \bar\Upphi_5 \), then any two of its substructures which are isomorphic to, say, \( \bar{F}_u \) and \( \bar{F}_v \), either coincide or else share at most the same root.

Finally, let \( \bar\Upphi_6 \) be the following \( \bar{\LL}^b_{\omega_1 \omega_1} \)-sentence:
\begin{multline}\tag{$\bar\Upphi_6$}
\forall \seqofLR{ x_n}{n < \omega } \forall \seqofLR{ y_m}{ m < \omega } \bigg[ \bigvee_{u \in \pre{<\omega}{2}} \left ( \bar\upchi_u ( \seqofLR{ x_n }{ n \in \omega } ) \wedge \bar\upchi_u (\seqofLR{ y_m }{ m \in \omega } ) \right )
\\
{} \implies \bigwedge_{n \in \omega} ( x_n \equals y_n ) \vee \bigwedge_{n,m \in \omega} ( x_n \nequals y_m ) \bigg].
\end{multline}
This means that if \( X \models \bar\Upphi_6 \), where \( X \) is an \( \bar{\LL} \)-structure, then any two of its substructures which are isomorphic to \emph{the same \( \bar{F}_u \)}, either they coincide or else they are completely disjoint.

Let finally \( \bar\Uppsi \) be the \( \bar{\LL}^b_{\kappa^+ \kappa} \)-sentence 
\begin{equation}\label{eq:barPsi}
 \upsigma_{\OCT} \wedge \bar\Upphi_1 \wedge \bar\Upphi_2 \wedge \bar\Upphi_3 \wedge \bar\Upphi_4 \wedge \bar\Upphi_5 \wedge \bar\Upphi_6 . \index[symbols]{Psib@\( \bar\Uppsi \)}
\end{equation}

It is easy to check that:

\begin{lemma} \label{lem:barrange}
Let \( T \in \TT_\kappa \) and \( \bar{f}_T \) be as in~\eqref{eq:defbarf_T}.
Then \( \ran(\bar{f}_T) \subseteq \Mod^\kappa_{\bar\Uppsi} \).
\end{lemma} 

We now classify again the structures in \( \Mod^\kappa_{\bar{\Uppsi}} \) using the elements of \( \pre{\pre{<\omega}{2} \times \kappa}{2} \) as invariant. 
For \( X \in \Mod^\kappa_{\bar\Uppsi} \), define \( \bar{g} ( X ) \colon \pre{< \omega}{ 2} \times \kappa \to 2 \) by letting \( g ( X ) ( u , \alpha ) = 1 \) if and only if there is a substructure of \( X \) isomorphic to \( \bar{F}_u \) whose root is the \( \alpha \)-th element with respect to the order \( \order ^X \). 
Formally,
\begin{equation} \label{eq:defbarg}
\bar{ g } \colon \Mod^\kappa_{\bar\Uppsi} \to \pre{\pre{< \omega}{ 2} \times \kappa}{2} , \qquad \bar{g} ( X ) ( u , \alpha ) = 1 \IFF X \models \bar\upsigma_{u , \alpha} , 
\end{equation}
where \( \bar\upsigma_{u , \alpha} \) is the \( \bar{\LL}^b_{\kappa^+ \kappa} \)-sentence 
\begin{equation}\tag{$\bar\upsigma_{u , \alpha}$}
 \EXISTS{ x} \left [ \bar\uprho_{ \alpha } ( x ) \wedge \EXISTS{ \seqofLR{ z_n }{ n < \omega } } (\bar\upchi_u ( \seqofLR{ z_n }{ n < \omega } ) \AND x \equals z_0) \right ].
\end{equation} 

The following proposition is the analogue in the new context of Proposition~\ref{prop:classification}, and its proof is quite similar to (but simpler than) the original one. 
However, for the reader's convenience we fully reprove it here because such classification result lies at the core of the proof of the invariant universality of \( \embeds^\kappa_\OCT \).

\begin{proposition} \label{prop:barclassification}
The map \( \bar{g} \) defined in~\eqref{eq:defbarg} reduces \( \cong \) to \( = \).
\end{proposition}

\begin{proof}
Let \( X , Y \in \Mod^\kappa_{\bar\Uppsi} \) be isomorphic.
Then for every \( (u, \alpha) \in \pre{< \omega}{ 2} \times \kappa \), \( X \models \bar\upsigma_{u , \alpha} \iff Y \models \bar\upsigma_{u , \alpha} \), whence \( \bar{g} ( X ) = \bar{g} ( Y ) \) by~\eqref{eq:defbarg}.

Conversely, assume \( \bar{g} ( X ) = \bar{g} ( Y ) \) for some \( X , Y \in \Mod^\kappa_{\bar\Uppsi} \). 
Since \( X \) and \( Y \) are models of \( \bar\Uppsi \) they satisfy \( \upsigma_{\OCT} \) and \( \bar\Upphi _1 , \dotsc , \bar\Upphi _6 \).
Let \( X' \) and \( Y' \) be the substructures of \( X \) and \( Y \) whose domains are the fields of the orderings \( \order ^X \) and \( \order ^Y \), respectively.
As \( X \) and \( Y \) satisfy \( \bar\Upphi_1 \wedge \bar\Upphi_2 \), then \( X' \cong Y' \) via a (unique) isomorphism \( \iota_{ X' , Y' } \).
We will now extend \( \iota_{X',Y'} \) to an isomorphism 
\[
\iota_{ X , Y } \colon X \to Y .
\] 
Let \( a \in X \setminus X' \). 
Since \( X \models \bar\Upphi_3 \), there are \( \seqofLR{ a_n}{n < \omega} \in \pre{\omega}{ X } \) and \( 0 \neq \bar{\imath} < \omega \) such that \( a = a_{\bar{\imath}} \) and \( X \models \bar\upchi_u [ \seqofLR{ a_n }{ n < \omega } ] \) for some \( u \in \pre{< \omega}{ 2} \). 
In particular, \( a_0 \) is in the field \( X' \) of \( \order ^X \), which is a well-order of length \( \kappa \) since \( X \models \bar\Upphi_1 \). 
Let \( \alpha < \kappa \) be such that \( a_0 \) is the \( \alpha \)-th element in this order: then by definition of \( \bar{g} \) we have \( \bar{g} ( X ) ( u , \alpha ) = 1 \), so \( \bar{g} ( Y ) ( u , \alpha ) = 1 \) by case assumption. 
Let \( b_0 \in Y' \) be the \( \alpha \)-th element with respect to the order \( \order ^Y \) --- \( b_0 \) is well-defined because \( Y \models \bar\Upphi_1 \).
Choose \( \seq{ b_1 , b_2 , \dotsc} \in \pre{\omega}{Y} \) such that \( Y \models \bar\upchi_u [ \seqofLR{ b_n }{ n < \omega } ] \), and set \( \iota_{ X , Y } ( a ) \coloneqq b_{\bar{\imath}} \).
The definition of \( \iota_{ X , Y } ( a ) \) seems to depend on the choice of the \( b_n \)'s, but the next claim shows that this is not the case.
\begin{claim}\label{claim:barwelldefined}
Suppose \( \alpha < \kappa \), \( u \in \pre{< \omega}{ 2} \), and \( \seqofLR{ a_n }{ n < \omega } \) is a sequence of elements of \( X \) such that \( a_0 \) is the \( \alpha \)-th element in \( \order^X \) and \( X \models \bar\upchi_u [ \seqofLR{ a_n }{ n < \omega } ] \).
Then there is a unique sequence \( \seqofLR{ b_n }{ n < \omega } \) of elements of \( Y \) such that \( b_0 \) is the \( \alpha \)-th element in \( \order ^Y \) and \( Y \models \bar\upchi_u [ \seqofLR{ b_n }{ n < \omega } ] \).
Therefore \( \iota_{ X , Y } ( a_n) = b_n \) for every \( n < \omega \).
\end{claim}

\begin{proof}[Proof of the Claim]
Given two sequences \( \seqofLR{b_n }{ n< \omega } \) and \( \seqofLR{ b'_n }{ n < \omega } \) as in the conclusion of the claim, we have that \( b_0 = b'_0 \) and both \( Y \models \bar\upchi_u [ \seqofLR{ b_n }{ n < \omega } ] \) and \( Y \models \bar\upchi_u [\seqofLR{ b'_n }{ n < \omega } ] \).
Since \( Y \models \bar\Upphi_6 \), it follows that \( b_n = b'_n \) for all \( n < \omega \).
\end{proof}

\begin{claim}
 \( \iota_{ X , Y } \) is injective.
\end{claim}

\begin{proof}[Proof of the Claim]
Let \( a , a' \in X \) be distinct.
If either \( a \) or \( a' \) belongs to \( X' \), then \( \iota_{ X , Y } ( a) \neq \iota_{ X , Y } ( a') \) because \( \iota_{ X , Y } \restriction X' = \iota_{X' , Y' } \) is a bijection between \( X' \) and \( Y' \), and \( \iota_{ X , Y } ( X \setminus X' ) \subseteq Y \setminus Y' \) by construction. 
So we can assume \( a , a' \notin X' \).

Assume towards a contradiction that \( \iota_{ X , Y } ( a ) = \iota_{ X , Y } ( a') \). 
Since \( X \models \bar\Upphi_1 \wedge \bar\Upphi_3 \), let \( \alpha , \beta < \kappa \), \( u , v \in \pre{< \omega}{ 2} \), \( \seqofLR{a_n }{ n < \omega } , \seqofLR{ a'_n }{ n < \omega } \in \pre{ \omega }{X} \), and \( 0 < \bar{\imath} , \bar{\jmath} < \omega \) be such that \( a_0 \) and \( a'_0 \) are, respectively, the \( \alpha \)-th and \( \beta \)-th elements of \( \order^X \), \( X \models \bar\upchi_u [ \seqofLR{ a_n }{ n < \omega } ] \), \( X \models \bar\upchi_v [ \seqofLR{ a'_n }{ n < \omega } ] \), \( a = a_{\bar{\imath}} \), and \( a' = a'_{\bar{\jmath}} \). 
By Claim~\ref{claim:barwelldefined} let \( \seqofLR{b_n }{ n < \omega } \) and \( \seqofLR{b'_n }{ n < \omega } \) be the unique sequences such that \( b_0 \) and \( b'_0 \) are, respectively, the \( \alpha \)-th and \( \beta \)-th elements of \( \order^Y \), \( Y \models \bar\upchi_u [ \seqofLR{ b_n }{ n < \omega } ] \), and \( Y \models \bar\upchi_v [ \seqofLR{ b'_n }{ n < \omega } ] \).
Then \( b_n = \iota_{ X , Y } ( a_n) \) and \( b'_n = \iota_{ X , Y } ( a'_n ) \) for every \( n \), so that \( \iota_{ X , Y } ( a ) = b_{\bar{\imath}} \) and \( \iota_{ X , Y } ( a' ) = b'_{\bar{\jmath}} \). 
Since \( Y \models \bar\Upphi_5 \) and 
\begin{equation}\label{eq:barb_i}
 b_{\bar{\imath}} = \iota_{ X , Y } ( a ) = \iota_{ X , Y } ( a') = b'_{\bar{\jmath}} ,
\end{equation}
it follows that 
\begin{equation}\label{eq:barb_n}
 \forall n \left ( b_n = b'_n\right ).
\end{equation}
Let \( Z \) be the substructure of \( Y \) with domain \( \setofLR{b_n}{ n \in \omega } \).
Since \( Y \models (\bar\upchi_u \wedge \bar\upchi_v ) [ \seqofLR{ b_n }{ n < \omega } ] \), then \( \bar{F}_u \cong Z \cong \bar{F}_v \) by~\eqref{eq:barchi_u}.
Thus \( u = v \) by (the analogue of)~\eqref{eq:=iffiso}. 

Notice that \( a_0 , a'_0 \in X' \) and \( b_0,b'_0 \in Y' \), thus \( b_0 = b'_0 \) by~\eqref{eq:barb_n}, and hence \( a_0 = a'_0 \) as \( \iota_{ X , Y } \restriction X' = \iota_{X',Y'} \) is a bijection between \( X' \) and \( Y' \). 
Moreover, \( Y \models \bar\upchi_u [ \seqofLR{ b_n }{ n < \omega } ] \) implies that the \( b_n \)'s are all distinct. 
Therefore \( \bar\imath = \bar\jmath \) because \( b_{\bar\imath} = b_{\bar\jmath} \) by~\eqref{eq:barb_i} and~\eqref{eq:barb_n}. 
Since \( u = v \), it follows that \( X \models \bar\upchi_u [ \seqofLR{ a_n }{ n < \omega } ] \) and \( X \models \bar\upchi_u [ \seqofLR{ a'_n }{ n < \omega } ] \).
Therefore, \( a_0 = a'_0 \) and \( X \models \bar\Upphi_6 \) imply \( \forall n (a_n = a'_n) \), whence \( a = a_{\bar\imath} = a'_{\bar\imath} = a'_{\bar\jmath} = a' \), a contradiction.
\end{proof}

\begin{claim}
\( \iota_{ X , Y } \) is surjective and \( \iota_{Y,X} = \iota^{-1}_{ X , Y } \). 
\end{claim}

\begin{proof}[Proof of the Claim]
First notice that \( \iota_{Y,X} \restriction Y' = \iota_{Y',X'} \) is an isomorphism between \( Y' \) and \( X' \) and is the inverse of \( \iota_{X',Y'} = \iota_{ X , Y } \restriction X' \) (by the uniqueness of the isomorphism between \( X' \) and \( Y' \)). 
Let now \( b \in Y \setminus Y' \). 
Since \( Y \models \bar\Upphi_1 \wedge \bar\Upphi_3 \), there are \( \alpha < \kappa \), \( u \in \pre{< \omega}{ 2} \), \( \seqofLR{b_n }{ n \in \omega } \in \pre{ \omega }{ Y} \), and \( 0 \neq \bar\imath <\omega \) such that \( b_0 \) is the \( \alpha \)-th element in the order \( \order ^Y \), \( Y \models \bar\upchi_u [ \seqofLR{ b_n }{ n < \omega } ] \), and \( b = b_{\bar\imath} \). 
Applying Claim~\ref{claim:barwelldefined} with the role of \( X \) and \( Y \) interchanged, there is a unique sequence \( \seqofLR{a_n }{ n \in \omega } \in \pre{ \omega }{X} \) such that \( a_0 \) is the \( \alpha \)-th element in \( \order^X \) and \( X \models \bar\upchi_u [\seqofLR{ a_n }{ n < \omega } ] \).
Then \( \iota_{Y,X} ( b ) = a_{\bar\imath} \) by definition. 
Using Claim~\ref{claim:barwelldefined}, it is easy to check that \( \iota_{ X , Y } ( a_{\bar\imath} ) = b \), so we are done.
\end{proof}

We have shown that \( \iota_{ X , Y } \colon X \to Y \) is a bijection. 
It remains to be proved that it is also an isomorphism. 
Since \( \iota^{-1}_{ X , Y } = \iota_{Y,X} \), it is enough to show just one direction of each equivalence, i.e.~that for every \( a , a' \in X \) 
 \begin{align*}
 a \edge^X a' & \IMPLIES \iota_{ X , Y } ( a ) \edge^Y \iota_{ X , Y } ( a' )
 \\
 a \order ^X a' & \IMPLIES \iota_{ X , Y } ( a ) \order ^Y \iota_{ X , Y } ( a' ) .
 \end{align*}
If \( a , a' \in X' \) the result follows from the fact that \( \iota_{ X , Y } \restriction X' = \iota _{ X' , Y'} \) is an isomorphism.
In particular, since \( a \order^X a' \) implies \( a , a' \in X' \), the second implication holds.
Thus it is enough to focus on the first implication and assume, without loss of generality, that \( a \notin X' \) (both \( \edge^X \) and \( \edge^Y \) are symmetric relations by \( X \models \upsigma_{\OCT} \) and \( Y \models \upsigma_{\OCT} \)).
Since \( X \models \bar\Upphi_1 \wedge \bar\Upphi_3 \), there are \( \alpha < \kappa \), \( u \in \pre{< \omega}{ 2} \), \( \seqofLR{a_n }{ n \in \omega } \in \pre{ \omega }{ X} \), and \( 0 \neq \bar\imath < \omega \) such that \( a_0 \) is the \( \alpha \)-th element in \( \order^X \), \( X \models \bar\upchi_u [ \seqofLR{ a_n }{ n < \omega } ] \), and \( a = a_{\bar\imath} \). 
If \( a \edge^X a' \), then \( X \models \upsigma_{\OCT} \wedge \bar\Upphi_4 \) implies \( a' = a_{\bar\jmath} \) for some \( \bar\jmath \neq \bar\imath \). 
By Claim~\ref{claim:barwelldefined}, let \( \seqofLR{b_n }{ n \in \omega } \) be the unique sequence such that \( b_0 \) is the \( \alpha \)-th element in \( \order^Y \) and \( Y \models \bar\upchi_u [ \seqofLR{ b_n }{ n < \omega } ] \), so that \( b_n = \iota_{ X , Y } ( a_n ) \) for all \( n \in \omega \) by definition.
Since \( a_{\bar\imath} = a \edge^X a' = a_{\bar\jmath} \) and \( X \models \bar{\upchi}_u[\seqofLR{a_n}{n < \omega}] \) imply that the subformula \( \bar{\upchi}^u_{ \bar\imath , \bar\jmath } ( x_{\bar\imath} , x_{\bar\jmath} ) \) of \( \bar\upchi_u ( \seqofLR{ x_n }{ n < \omega } ) \) is \( x_{\bar\imath} \edge x_{\bar\jmath} \), 
and since \( Y \models \bar\upchi_u [ \seqofLR{ b_n }{ n < \omega } ] \),
we get \( \iota_{ X , Y } ( a ) = b_{\bar\imath} \edge^Y b_{\bar\jmath} = \iota_{ X , Y } ( a' ) \).
\end{proof}

Endow \( \pre{\pre{< \omega}{2} \times \kappa}{2} \) with the product topology \( \tau_p \) (see Section~\ref{subsubsec:spacesoftypekappa} and, in particular, Example~\ref{xmp:canonicalbijection}\ref{xmp:canonicalbijection-1}). 
The following is the analogue of Lemma~\ref{lem:open} and can be proven in the same way (just replace the \( \LL^b_{\kappa^+ \kappa} \)-sentence \( \upsigma_{u,\alpha} \) with the \( \bar{\LL}^b_{\kappa^+ \kappa} \)-sentence \( \bar\upsigma_{u,\alpha} \)).

\begin{lemma} \label{lem:baropen}
For every (\( \tau_p \)-)open set \( V \subseteq \pre{\pre{< \omega}{ 2} \times \kappa}{2} \) there is an \( \bar{\LL}^b_{ \kappa^+ \kappa } \)-sentence \( \bar \upsigma_V \) such that \( g^{ - 1 } ( V ) = \Mod^\kappa_{ \bar\upsigma_V \wedge \bar\Uppsi} \).
\end{lemma}

The following lemma, which corresponds to Lemma~\ref{lem:continuous}, requires instead a slightly different argument due to the different coding we used.

\begin{lemma} \label{lem:barcontinuous}
Let \( T \in \TT_\kappa \) and \( \bar{f}_T \) be the map defined in~\eqref{eq:defbarf_T}. 
Then the map \( \bar{g} \circ \bar{f}_T \colon \pre{\omega}{ 2} \to \pre{\pre{< \omega}{ 2} \times \kappa}{2} \) is continuous.
\end{lemma}

\begin{proof}
Arguing as in the proof of Lemma~\ref{lem:continuous}, since \( \mathcal{S} = \{ \widetilde{\Nbhd}^A_{ ( u , \alpha ) , i } \Mid ( u , \alpha ) \in \pre{< \omega}{ 2} \times \kappa \wedge {i = 0,1} \} \) is a subbasis for the product topology on \( \pre{A}{2} \) (where \( A \coloneqq \pre{< \omega}{ 2} \times \kappa \)), it is enough to check that \( ( \bar{g} \circ \bar{f}_T ) ^{-1} ( \widetilde{\Nbhd}^A_{ ( u , \alpha ) , i } ) \) is open in \( \pre{\omega}{2} \) for every \( (u, \alpha) \in \pre{< \omega}{ 2} \times \kappa \) and \( i = 0,1 \). 
Fix \( (u, \alpha) \in \pre{< \omega}{ 2} \times \kappa \). 
If \( \alpha \) is odd then \( ( \bar{g} \circ \bar{f}_T ) ^{-1} ( \widetilde{\Nbhd}^A_{ ( u , \alpha ) , 1 } ) = \emptyset \) and \( ( \bar{g} \circ \bar{f}_T ) ^{-1} ( \widetilde{\Nbhd}_{ ( u , \alpha ) , 0 }) = \pre{\omega}{2} \).
If \( \alpha \) is even let 
\[
 A_1 \coloneqq \setofLR{ v \in \pre{ \lh u }{2} }{ ( u , v , s_\alpha ) \in \tilde{T} } ,
\]
where \( \tilde{T} \) is the tree obtained from \( T \) as in~\eqref{eq:S_T} and \( s_\alpha = \#^{-1}(\alpha) \), and let \( A_0 \coloneqq \pre{ \lh u}{2} \setminus A_1 \). 
Then by definition of \( \bar{f}_T \) and \( \bar{g} \) we get \( (\bar{g} \circ \bar{f}_T ) ^{-1} ( \widetilde{\Nbhd}_{ ( u , \alpha ) , i }) = \bigcup_{v \in A_i } \Nbhd_v \).
\end{proof}

Finally, the next corollaries are the counterparts of Corollaries~\ref{cor:closed} and~\ref{cor:saturation}: their proofs are obtained from the original ones by systematically replacing Lemmas~\ref{lem:range}, ~\ref{lem:open}, and~\ref{lem:continuous}, Proposition~\ref{prop:classification}, and Corollary~\ref{cor:closed}, with Lemmas~\ref{lem:barrange}, ~\ref{lem:baropen}, and~\ref{lem:barcontinuous}, Proposition~\ref{prop:barclassification}, and Corollary~\ref{cor:barclosed}, respectively.

\begin{corollary} \label{cor:barclosed}
Let \( T \in \TT_\kappa \) and \( \bar{f}_T \) be as in~\eqref{eq:defbarf_T}. 
For every closed set \( C \subseteq \pre{\omega}{ 2} \), \( ( \bar{g} \circ \bar{f}_T ) ( C ) \) is closed. 
\end{corollary}

\begin{corollary} \label{cor:barsaturation}
Let \( T \in \TT_\kappa \) and \( \bar{f}_T \) be as in~\eqref{eq:defbarf_T}. 
There is an \( \bar{\LL}^b_{\kappa^+ \kappa} \)-sentence \( \bar\upsigma_T \) such that the closure under isomorphism of \( \ran ( \bar{f}_T ) \) is \( \Mod^\kappa_{\bar\upsigma_T} \).
\end{corollary}

We finally defined the ``inverse'' map of \( \bar{f}_T \) analogously to Definition~\ref{def:defh_T}.

\begin{definition}\label{def:defbarh_T}
For \( T \in \TT_\kappa \), \( \bar{f}_T \) as in~\eqref{eq:defbarf_T} and \( \bar\upsigma_T \) as in Corollary~\ref{cor:barsaturation} , let \( \bar{h}_T \colon \Mod^\kappa_{\bar\upsigma_T} \to \pre{\omega}{2} \) be defined by
\begin{equation}\label{eq:defbarh_T}
\bar{h}_T ( X ) = x \iff \bar{f}_T ( x ) \cong X.
\end{equation}
\end{definition}

Note that the function \( \bar{h}_T \) is well-defined by Corollary~\ref{cor:barsaturation} and Theorem~\ref{th:OCT}. 
The following corollary is analogous to Corollary~\ref{cor:inverse}, and it completes the proof of Theorem~\ref{th:barmain}: it can be proved with an argument similar to the original one (but replacing Theorem~\ref{th:graphs} with Theorem~\ref{th:OCT}).

\begin{corollary} \label{cor:barinverse}
Let \( T \in \TT_\kappa \), \( \bar{f}_T \) be as in~\eqref{eq:defbarf_T}, \( \bar\upsigma_T \) as in Corollary~\ref{cor:barsaturation}, and \( \bar{h}_T \) as in Definition~\ref{def:defbarh_T}.
If \( R = \PROJ \body{ T } \) is a quasi-order then \( \bar{h}_T \) simultaneously reduces \( \embeds \) to \( R \) and \( \cong \) to \( = \).
\end{corollary}

The following is the analogue of Theorem~\ref{th:main}.

\begin{theorem} \label{th:barmain}
Let \( \kappa \) be an infinite cardinal and \( T \in \TT_\kappa \).
If \( R = \PROJ \body{ T } \) is a quasi-order, then there are \( \bar\upsigma_T , \bar{f}_T , \bar{h}_T \) such that:
\begin{enumerate-(a)}
\item \label{th:barmain-a}
\( \bar\upsigma_T \) is an \( \bar{\LL}^b_{\kappa^+ \kappa} \) sentence;
\item \label{th:barmain-b}
\( \bar{f}_T \) reduces \( R \) to \( \embeds^\kappa_{\bar\upsigma_T} \) and \( = \) to \( \cong^\kappa_{\bar\upsigma_T} \);
\item \label{th:barmain-c}
\( \bar{h}_T \) reduces \( \embeds^\kappa_{\bar\upsigma_T} \) to \( R \) and \( \cong^\kappa_{\bar\upsigma_T} \) to \( = \);
\item \label{th:barmain-d}
\( \bar{h}_T \circ \bar{f}_T = \id \) and \( (\bar{f}_T \circ \bar{h}_T) ( X ) \cong X \) for every \( X \in \Mod^\kappa_{\bar\upsigma_T} \).
\end{enumerate-(a)}
In particular, \( \embeds_\OCT^\kappa \) is invariantly universal for \( \kappa \)-Souslin quasi-orders on \( \pre{\omega}{2} \) (equivalently, on any uncountable Polish or standard Borel space).
\end{theorem}

By Theorem~\ref{th:barmain}\ref{th:barmain-d} we further have:
\begin{corollary}\label{cor:barmaintheorem}
Let \( \kappa \) be an infinite cardinal and \( T \in \TT_\kappa \).
If \( R = \PROJ \body{ T } \) is a quasi-order, then there is an \( \bar{\LL}^b_{\kappa^+ \kappa} \)-sentence \( \bar\upsigma_T \) such that the quotient orders of \( R \) and \( {\embeds^\kappa_{\bar\upsigma_T}} \) are \emph{isomorphic}.
\end{corollary}

Theorem~\ref{th:barmain}\ref{th:barmain-a} should be contrasted with Theorem~\ref{th:main}\ref{th:main-a}, in that it involves a formula belonging to the fragment \( \bar{\LL}^b_{\kappa^+ \kappa} \) of \( \bar{\LL}_{\kappa^+ \kappa} \): this yields to the fact that by Corollary~\ref{cor:formulaborel1}\ref{cor:formulaborel1-ii} 
\begin{equation} \label{eq:barupsigmaisBorel} 
\text{if \( \kappa \) is \emph{regular}, } \Mod^\kappa_{\bar\upsigma_T} \text{ is an effective \( \kappa+1 \)-Borel subset of } \Mod^\kappa_{\bar\LL} \text{ (with respect to \( \tau_b \))} .
\end{equation} 
This important observation will allow us to get generalizations of Theorem~\ref{th:mottorosfriedman2} to all uncountable regular \( \kappa \)'s (see Theorem~\ref{th:barmaintopology}).

Moreover, Remark~\ref{rem:maintheorem}\ref{rem:maintheorem-1} still holds after replacing the corresponding \( \LL_{\kappa^+ \kappa} \)-sentence \( \Uppsi \) with the \( \bar{\LL}^b_{\kappa^+ \kappa} \)-sentence \( \bar\Uppsi \). 
Remark~\ref{rem:maintheorem}\ref{rem:maintheorem-2} can be improved to the following Remark~\ref{rem:barmaintheorem}\ref{rem:barmaintheorem-1}, and since only bonded formul\ae{} are now involved in it, we further get Remark~\ref{rem:barmaintheorem}\ref{rem:barmaintheorem-2}

\begin{remarks} \label{rem:barmaintheorem}
\begin{enumerate-(i)}
\item \label{rem:barmaintheorem-1}
Every model \( X \) of \( \bar\Uppsi \) (and hence every model of \( \bar\upsigma_T \) for \( T \in \TT_\kappa \)) admits an \( \bar{\LL}^b_{\kappa^+ \kappa} \)-Scott sentence \( \bar{\upsigma}^X \).
The sentence \( \bar\upsigma^X \) is again obtained by applying the argument of Corollary~\ref{cor:barsaturation} (see also Corollary~\ref{cor:saturation}) to the singleton of \( \bar{h}_T ( X ) \).

\item \label{rem:barmaintheorem-2}
By the previous observation and Corollary~\ref{cor:formulaborel1}\ref{cor:formulaborel1-ii}, it follows that if \( \kappa \) is regular then the isomorphism type \( \setofLR{ Y \in \Mod^\kappa_{\bar{\LL}} }{ Y \cong X } \) of every \( X \in \Mod^\kappa_{\bar\Uppsi} \) is an effective \( \kappa + 1 \)-Borel set with respect to the bounded topology.
\end{enumerate-(i)}
\end{remarks}

Notice that if \( \kappa^{<\kappa} = \kappa \), then by the generalized Lopez-Escobar theorem (i.e.\ Theorem~\ref{th:LopezEscobar}) a structure \( X \in \Mod^\kappa_\LL \) admits an \( \LL_{\kappa^+ \kappa} \)-Scott sentence (for \( \LL \) and arbitrary finite relational language) if and only if its isomorphism type is \( \kappa + 1 \)-Borel with respect to the bounded (equivalently, to the product, or to the \( \lambda \)-)topology, i.e.~\ref{rem:barmaintheorem-1} and~\ref{rem:barmaintheorem-2} are essentially equivalent once we drop the restriction \( \bar{\upsigma}^X \in \bar{\LL}^b_{\kappa^+ \kappa} \). 
However, recall that in most of the applications considered in this paper, either we will work in models of \( \AD \) (where \( \AC \) fails), or we will deal with cardinals \( \kappa \) which are strictly smaller than \( 2^{\aleph_0} \) (see the discussions in Sections~\ref{sec:introduction} and~\ref{sec:Ksouslinsets}), so in both cases we lack the crucial condition \( \kappa^{< \kappa} = \kappa \), even when e.g.\ \( \kappa = \omega_1 \). 
This is why we preferred to consider separately the two statements~\ref{rem:barmaintheorem-1} and~\ref{rem:barmaintheorem-2} of Remark~\ref{rem:barmaintheorem}.

As in the case of Remark~\ref{rem:barmaintheorem}\ref{rem:barmaintheorem-1} (see the discussion after Remarks~\ref{rem:maintheorem}), Remark~\ref{rem:barmaintheorem}\ref{rem:barmaintheorem-2} is again well-known when \( \kappa = \omega \) (for any countable language \( \LL \) and any countable structure \( X \in \Mod^\omega_{\LL} \), see e.g.\ 
\cite[Theorem~16.6]{Kechris:1995zt}). 
However, when \( \kappa \) is uncountable then such condition may fail again for some \( X \in \Mod^\kappa_{\mathcal{L}} \): when \( \kappa^{< \kappa} = \kappa \) and \( \AC \) is assumed this already follows from Example~\ref{xmp:Scottsentence} and Theorem~\ref{th:LopezEscobar}, but since as discussed above we will often work in a different context we point out the following independent counterexample due to S.~D.~Friedman.

\begin{example} \label{xmp:isomorphismtype}
Let \( \mathcal{L} \) be the order language, \( \kappa = \omega_1 \), and endow both \( \pre{\omega_1}{2} \) and \( \Mod^{\omega_1}_\mathcal{L} \) with the bounded topology. 
Consider the structure 
\[ 
X \coloneqq \omega_1 \times \Q_1,
 \] 
where \( \Q_1 \) is defined as in Example~\ref{xmp:Scottsentence} and \( X \) is endowed with the lexicographical ordering: we claim that the isomorphism type of \( X \) is not \( \kappa + 1 \)-Borel. 
Otherwise, it would be \( \kappa + 1 \)-Borel also the club filter \( \mathrm{CUB} \subseteq \pre{\omega_1}{2} \) defined by
\[ 
\mathrm{CUB} \coloneqq \setofLR{ x \in \pre{\omega_1}{2} }{ \setofLR{ \alpha < \omega_1 }{ x ( \alpha ) = 1 } \text{ contains a club of } \omega_1 }
 \] 
(because there is a continuous map \( f \colon \pre{\omega_1}{2} \to \Mod^{\omega_1}_\mathcal{L} \) such that \( f^{-1} ( \setofLR{ Y \in \Mod^{\omega_1}_\mathcal{L} }{ Y \cong X }) = \mathrm{CUB} \), see below).
But this is impossible, because every \( \kappa + 1 \)-Borel subset \( B \subseteq \pre{\omega_1}{2} \) has the \( \kappa \)-Baire property by Proposition~\ref{prop:borelsetshaveBP}, while \( \mathrm{CUB} \) does not have the \( \kappa \)-Baire property (see~\cite{Halko:2001kl}): in fact, for every intersection \( C \subseteq \pre{\omega_1}{2} \) of at most \( \omega_1 \)-many open dense subsets of \( \pre{\omega_1}{2} \) and every \( s \in \pre{<\omega_1}{2} \) there are \( x , y \in C \) such that \( x \in \Nbhd_s \setminus \mathrm{CUB} \) (which implies that \( U \triangle \mathrm{CUB} \) is non-meager for every nonempty open set \( U \)) and \( y \in \mathrm{CUB} \) (which implies that \( \mathrm{CUB} \) itself is non-meager), so Proposition~\ref{prop:nonmeagerarelocallycomeager} gives the desired result. 
The reduction \( f \) of \( \mathrm{CUB} \) to the isomorphism class of \( X \) is defined as follows: given \( x \in \pre{\omega_1}{2} \), let \( f ( x ) \) be (an isomorphic copy on \( \omega_1 \) of) the structure \( Y_x \) obtained by ordering lexicographically the set%
\footnote{We define the initial segment \( \setLR{ 0 } \times \Q_1 \) of \( Y_x \) independently from \( x \) because we want to guarantee that \( Y_x \) always has a minimal element, as \( X \) does.}
\[ 
\setLR{ 0 } \times \Q_1 \cup \left ( \setofLR{\alpha }{ x ( \alpha ) = 1 } \times \Q_1 \right ) \cup \left ( \setofLR{ \alpha }{ x ( \alpha ) = 0 } \times \Q_0 \right ),
 \] 
where \( \Q_0 , \Q_1 \) are again defined as in Example~\ref{xmp:Scottsentence} (it is left to the reader to verify that we can identify each \( Y_x \) with an isomorphic copy on \( \omega_1 \) in such a way that the resulting \( f \) is continuous). 
Notice that, in particular, when \( x \in \pre{\omega_1}{2} \) is constantly equal to \( 1 \) then \( Y_x = X \). 
If \( x \in \mathrm{CUB} \), let \( C \subseteq \omega_1 \) be a club of \( \omega_1 \) such that \( 0 \in C \) and \( x(\alpha) = 1 \) for every \( 0 \neq \alpha \in C \), and fix an order preserving bijection \( j \colon C \to \omega_1 \). 
Then the map \( i \) sending \( (\alpha , 0) \) to \( (j(\alpha),0) \) (for \( \alpha \in C \)) can be extended to an isomorphism between \( Y_x \) and \( X \) (using the fact that any two countable dense linear orders without maximum and minimum are isomorphic to \( \Q \), and hence to each other). 
Conversely, if \( Y_x \cong X \) then \( Y_x \) contains a closed unbounded (with respect to its ordering) set of order type \( \omega_1 \), whence
\( x \in \mathrm{CUB} \). 
Therefore \( x \in \mathrm{CUB} \iff f ( x ) \cong X \), as required.
\end{example}
\section{Definable cardinality and reducibility} \label{sec:definablecardinality}
As recalled in Section \ref{subsec:cardinality}, under \( \AC \) the cardinality \( \card{X} \) is defined to be the unique cardinal \( \kappa \) in bijection with \( X \).
The resulting theory, though, is somewhat unsatisfactory, since there are no universally accepted methods to settle simple questions like e.g.~the continuum problem: by work of G\"odel~\cite{Godel:1938dz} and Cohen~\cite{Cohen:1963ys} the theory \( \ZFC \) does not settle the statement \( 2^{\aleph_0} = \aleph_1 \). 
Moreover, knowing that \( \card{\R } = \kappa \) for some cardinal \( \kappa \), gives little information on such \( \kappa \) or on the possible bijections between \( \R \) and \( \kappa \): in fact, there is no \emph{definable} or \emph{natural} way to explicitly well-order the reals in \( \ZFC \) (see e.g.~\cite{Shelah:1990sw}). 
This should be contrasted with the fact that, in practice, when \( \card{I} \leq \card{J} \) one would like to have an explicit witness for this fact, namely a reasonably defined injection from \( I \) to \( J \). 
These observations yield the notion of \emph{definable cardinality}.
Several definitions of this concept have been considered in the literature (usually depending on the kind of problem one is dealing with), and they all amount to restricting the objects under consideration and the functions used to compute cardinalities to some reasonably simple class, e.g.~to the class of functions belonging to some inner model such as \( \Ll ( \R ) \), \( \OD ( \R ) \), and so on. 
It is quite remarkable that when replacing the notion of cardinality with its definable version, all obstacles (i.e.~the independence results on the size of simple sets) simply disappear (see Section~\ref{subsubsec:reducibilityininnermodel}).

As for cardinalities, also the notion of reducibility is not completely satisfactory unless we impose some sort of definability condition on (i.e.~we restrict the class of) the reductions that can be used. 
For example, let us consider the Vitali equivalence relation \( E_0 \) on the real line \( \R \) defined by 
\begin{equation}\label{eq:E0}
r \mathrel{E_0} r' \IFF r - r' \in \Q .
\end{equation}
Under \( \AC \) we have that \( E_0 \leq \id ( \R ) \), but on the other hand there is no Baire-measurable (in particular, no Borel) witness for this reduction. 
This kind of phenomenon suggests that it could be more interesting to consider what could be called \emph{definable reducibility}.
Again, this is a vague notion and we will need to specify which definability conditions we are interested in. 

By the observations contained in Section~\ref{subsubsec:reducibility}, the notions of (definable) cardinality and of (definable) reducibility are strictly related, and all results about (definable) reducibility can be immediately translated into results about (definable) cardinality.
For example, we can e.g.~derive from the subsequent Theorem~\ref{th:pag:smallcardinalities} and Remark~\ref{rmk:smallcardinalities} a nice characterization of small cardinalities (see Section~\ref{subsubsec:reducibility}) in terms of infinitary sentences, or, to be more precise, in terms of the bi-embeddability relation on the models of such sentences: in fact, under \( \ADR \) the cardinality of any set \( A \) is small if and only if there is a (regular) cardinal \( \kappa < \Theta \) and an \( \LL_{\kappa^+ \kappa} \)-sentence \( \upsigma \) such that \( \card{A} = \cardLR{\Mod^\kappa_\upsigma/{\biembeds}} \), where \( \biembeds \) denotes as usual the bi-embeddability relation. 
In other words, under \( \ADR \) the set of small cardinalities can be written as
\[ 
\setofLR{ \cardLR{\Mod^\kappa_\upsigma/{\biembeds}} } { \kappa < \Theta, \upsigma \in \LL_{\kappa^+ \kappa}}. 
\]

In this section we will consider many different incarnations of definable cardinality and definable reducibility which have appeared in the literature, and using these notions (or naturally adapting them to our new setup) we will show that the completeness and invariant universality results obtained in Sections~\ref{sec:embeddabilitygraphs}--\ref{sec:alternativeapproach} can be interpreted as results on definable reducibility and hence, by the discussion above, also on definable cardinality. 

\subsection{Topological complexity} \label{subsec:topologicalcomplexity} 

The first case of definable reducibility we will analyze is the Borel reducibility \( \leq_{\bB} \): as discussed in the introduction, in this case the objects are \emph{analytic} relations on \( \pre{ \omega}{2} \) or, equivalently, on any uncountable Polish or standard Borel space, and the reductions are \emph{Borel functions} between such spaces.
This is nowadays the standard reducibility notion between \( \bSigma^1_1 \) quasi-orders and equivalence relations --- see e.g.~\cite{Becker:1996uq, Hjorth:2000zr, Kechris:2004jl, Kanovei:2008qo, Gao:2009fv}.
We will generalize this notion to our context by considering \( \kappa + 1 \)-Borel functions (instead of just Borel functions) as reductions between quasi-orders defined on (subsets of) spaces of type \( \lambda \), with \( \lambda \) an infinite cardinal. 
\emph{From this point onward, unless otherwise explicitly stated, we will conform to the standard practice in the existing literature of endowing all such spaces with the bounded topology \( \tau_b \),} and hence continuity and \( \kappa + 1 \)-Borelness (of both functions and sets) will always tacitly refer to this topology. 

When replacing Borel functions with \( \kappa + 1 \)-Borel functions, we have to decide whether we want to use the strong or the weak formulation of such class of functions (see Definition~\ref{def:weaklyBorel}). 
Of course the smaller is the class of reducing functions, the stronger are the results asserting that a certain quasi-order is reducible to another: so in general we should prefer \( \kappa + 1 \)-Borel reductions to the weakly \( \kappa + 1 \)-Borel ones. 
However, as discussed after Proposition~\ref{prop:weaklyalphaBorel}, when considering functions \( f \colon \pre{\lambda}{2} \to \pre{\mu}{2} \) with \( \lambda < \mu \) it is often more meaningful to consider the weaker notion of \( \kappa + 1 \)-Borelness --- even the inclusion map \( f \) defined in Proposition~\ref{prop:weaklyalphaBorel}\ref{prop:weaklyalphaBorel-b} can fail to be (non-weakly) \(\kappa + 1 \)-Borel. 
The next example shows that a similar phenomenon appears when considering reducibilities between quasi-orders.

\begin{example} \label{xmp:weaklynonBorelreducibilities}
Let \( \kappa \) and \( \lambda < \mu \) be infinite cardinals, and assume that \( \bB_{\kappa + 1} ( \pre{\lambda}{2} , \tau_b) \neq \pow ( \pre{\lambda}{2} ) \). 
Let \( E \) be the identity relation on \( \pre{\lambda}{2} \) (i.e.~\( x \mathrel{E} y \IFF x = y \)), and let \( F \) be its ``extension'' to the space \( \pre{\mu}{2} \), that is set for \( x , y \in \pre{\mu}{2} \)
\[ 
x \mathrel{F} y \iff x \restriction \lambda = y \restriction \lambda.
 \] 
The map \( \pre{\mu}{2} \to \pre{\lambda}{2} \) sending \( x \in \pre{ \mu}{2} \) to \( x \restriction \lambda \) is continuous (and hence also \( \kappa + 1 \)-Borel), and reduces \( F \) to \( E \). 
However, there is no \( \kappa + 1 \)-Borel reduction of \( E \) to \( F \). 
Indeed, let \( f \colon \pre{\lambda}{2} \to \pre{\mu}{2} \) be a reduction of \( E \) to \( F \), \( A \subseteq \pre{\lambda}{2} \) be a set not in \( \bB_{\kappa + 1} ( \pre{\lambda}{2} , \tau_b ) \), and set 
\[ 
A' \coloneqq \setofLR{ f ( x ) \restriction \lambda}{x \in A} .
 \] 
Then \( U \coloneqq \bigcup_{s \in A'} \Nbhd^\mu_s \) is \( \tau_b \)-open, but \( f^{ - 1 } ( U ) = A \) since \( f \) reduces \( E \) to \( F \), whence \( f^{ - 1 } ( U ) \notin \bB_{\kappa + 1} ( \pre{\lambda}{2} , \tau_b ) \).
\end{example}

Example~\ref{xmp:weaklynonBorelreducibilities} shows that if we decide to consider only (non-weakly) \( \kappa + 1 \)-Borel functions as reductions, then the identity relation \( E \) on \( \pre{\lambda}{2} \) would be \emph{strictly} more complex than its ``extension'' \( F \) to \( \pre{\mu}{2} \): this contradicts our intuition, which suggests that \( E \) and \(F \) should have the same complexity with respect to any reasonable notion of reducibility. 
All this discussion justifies the forthcoming definition of \( \kappa+ 1 \)-Borel reducibility, and our choice of the reducing functions should then appear quite natural.

\begin{definition} \label{def:Borelreducibility}
Let \( \lambda , \mu, \kappa \) be infinite cardinals, and let \( \mathcal{X} , \mathcal{Y} \) be spaces of type \( \lambda \) and \( \mu \), respectively. 
Let \( R \) and \( S \) be quasi-orders defined on (subsets of) \( \mathcal{X} \) and \( \mathcal{Y} \), respectively.
We say that \( R \) is \markdef{\( \kappa + 1 \)-Borel reducible} to \( S \) (in symbols \( R \leq^\kappa_{\bB} S \))\index[symbols]{116@\( \leq^\kappa_{\bB} \), \( \leq^\kappa_{\mathrm{B}} \), \( \sim^\kappa_{\bB} \), \( \sim^\kappa_{\mathrm{B}} \)} if and only if:%
\footnote{See Remark~\ref{rmk:partialBorelfunctions} for more on \emph{partial} (weakly) \( \kappa+1 \)-Borel functions.}
\begin{description}
\item[case \( \lambda \geq \mu \)]
there is a \( \kappa + 1 \)-Borel reduction \( f \colon \dom ( R ) \to \dom ( S ) \) of \( R \) to \( S \);

\item[case \( \lambda < \mu \)]
there is a \emph{weakly} \( \kappa + 1 \)-Borel function \( f \colon \dom ( R ) \to \dom ( S ) \) reducing \( R \) to \( S \).
\end{description} 

The quasi-orders \( R \) and \( S \) are \markdef{\( \kappa + 1 \)-Borel bi-reducible} (in symbols \( R \sim^\kappa_{\bB} S \)) if both \( R \leq^\kappa_{\bB} S \) and \( S \leq^\kappa_{\bB} R \).

Finally, when replacing (weakly) \( \kappa + 1 \)-Borel functions with their effective counterparts in the previous definitions, we get the notions of \markdef{effective \( \kappa + 1 \)-Borel (bi-)reducibility}, denoted by \( \leq^\kappa_{\mathrm{B}} \) and \( \sim^\kappa_{\mathrm{B}} \), respectively.
\end{definition}

Notice that when \( \lambda = \mu = \kappa = \omega \), the quasi-order \( \leq^\kappa_{\bB} \) coincide{\color{blue}s} with the usual Borel reducibility \( \leq_{\bB} \) mentioned in the introduction and at the beginning of this subsection.
Therefore the next completeness result naturally generalizes Theorem~\ref{th:LouveauRosendal} to uncountable \( \kappa \)'s.

\begin{theorem} \label{th:definable=borel}
Let \( \kappa \) be an infinite cardinal and \( R \) be a \( \kappa \)-Souslin quasi-order on \( \pre{ \omega}{2} \). 
Then \( R \leq^\kappa_{\bB} {\embeds_\CT^\kappa} \). 
In particular, \( \embeds_\CT^\kappa \) is \emph{\( \leq^\kappa_{\bB} \)-complete} for the class of \( \kappa \)-Souslin quasi-orders on \( \pre{\omega}{2} \), and the same is true when replacing \( \leq^\kappa_{\bB} \) with \( \leq^\kappa_{\mathrm{B}} \). 
\end{theorem}

\begin{proof}
The function \( f_T \colon \pre{\omega}{2} \to \CT_\kappa \) defined in~\eqref{eq:f_T} is easily seen to be continuous when both spaces are endowed with the \emph{product} topology \( \tau_p \). 
Since the topologies \( \tau_p \) and \( \tau_b \) coincide on \( \pre{\omega}{2} \), this shows that \( f_T \) is (effective) weakly \( \kappa + 1 \)-Borel (indeed, it is even effective weakly \( \alpha \)-Borel for any \( \alpha \geq \omega \)).
Since by Theorem~\ref{th:graphs}\ref{th:graphs-a} the function \( f_T \) reduces \( R \) to \( \embeds_\CT^\kappa \), the result follows. 
\end{proof}

As observed in Section~\ref{subsubsec:Gammainthecodes}, there is another natural generalization of the notion of Borel functions, namely \( \bGamma \)-in-the-codes functions (for suitable boldface pointclasses \( \bGamma \)): this yields a corresponding generalization of the Borel reducibility \( \leq_{\bB} \).

\begin{definition} \label{def:S(kappa)-reducibility}
Let \( \bGamma \) be a boldface pointclass, \( \kappa\) be an infinite cardinal for which there is a 
\( \bGamma \)-code (so that 
\(\kappa \leq \bdelta_{\bGamma} \) by Remark~\ref{rmk:inthecodes}\ref{rmk:inthecodes-ii}), and \( \mathcal{X} \) be a space of type \( \kappa \).
Given two quasi-orders \( R \) and \( S \) on, respectively, \( \pre{\omega}{2} \) and (a subset of) \( \mathcal{X} \), we say that \( R \) is \markdef{\( \bGamma \)-reducible} to \( S \) (in symbols, \( R \leq_{\bGamma} S \)) if there is a \( \bGamma \)-in-the-codes function \( f \colon \pre{\omega}{2} \to \dom(S) \subseteq \mathcal{X} \) which reduces \( R \) to \( S \).
\end{definition}

\begin{remark} \label{rmk:thereisacode}
Since one can speak of \( \bGamma \)-in-the-codes functions only if there is a \( \bGamma \)-code for \( \kappa \), \emph{in all the results mentioning the reducibility \( \leq_{\bGamma} \) we will tacitly assume that such a code exists}. When assuming \( \AD \) this is actually granted by Proposition~\ref{prop:deltaSkappa}\ref{prop:deltaSkappa-b}. Similarly, in the \( \AC \) world this implicit assumption will in most cases follow from the other hypotheses (see e.g.~Theorem~\ref{th:embeddingsofprojectiveqos} and the ensuing remark), although there are some cases, such as Theorem~\ref{th:mainSigma13AC}, in which the situation is less clear.
\end{remark}

Thus, in particular, \( \leq_{\bSigma^1_1 } \) coincides with \( \leq_{\bB} \) (notice that \( \bGamma = \bSigma^1_1 \) automatically implies \( \kappa = \omega \) in Definition~\ref{def:S(kappa)-reducibility}, so that in this case \( \mathcal{X} \) is homeomorphic to the classical Cantor space \( \pre{\omega}{2} \)). 
The next theorem can thus be seen as another natural generalization to uncountable \( \kappa \)'s of Theorem~\ref{th:LouveauRosendal}.

\begin{theorem} \label{th:Gammacompleteness}
\begin{enumerate-(a)}
\item \label{th:Gammacompleteness-a}
(\( \AC \)) Let \( \kappa \leq 2^{\aleph_0} \) be such that there is a 
\( \bS(\kappa) \)-code for it. 
Then the embeddability relation \( \embeds_\CT^\kappa \) is \emph{\( \leq_{\bS ( \kappa )} \)-complete} for the class of \( \kappa \)-Souslin quasi-orders on \( \pre{\omega}{2} \), i.e.~\( R \leq_{\bS ( \kappa )} {\embeds_\CT^\kappa} \) for every \( \kappa \)-Souslin quasi-order \( R \) on \( \pre{\omega}{2} \).
\item \label{th:Gammacompleteness-b}
(\( \AD + \DC \))
Let \( \kappa \) be a Souslin cardinal. 
Then the embeddability relation \( \embeds_\CT^\kappa \) is \emph{\( \leq_{\bS ( \kappa )} \)-complete} for the class of \( \kappa \)-Souslin quasi-orders on \( \pre{\omega}{2} \).
\end{enumerate-(a)} 
\end{theorem}

Observe that in Theorem~\ref{th:Gammacompleteness} it always makes sense to consider the notion of \( \leq_{\bS ( \kappa )} \)-reducibility (for part~\ref{th:Gammacompleteness-b} use Proposition~\ref{prop:deltaSkappa}\ref{prop:deltaSkappa-b}). Moreover, recall that 
the condition on \( \kappa \) in part~\ref{th:Gammacompleteness-a} is automatically satisfied if \( \kappa = \omega_1 \), \( \kappa = \omega_2 \), or \( \kappa = 2^{\aleph_0} \) by Proposition~\ref{prop:deltaSkappaAC}\ref{prop:deltaSkappaAC-a}.

\begin{proof}
Let \( R = \PROJ{\body{T}} \) be a \( \kappa \)-Souslin quasi-order on \( \pre{\omega}{2} \). 
As noticed at the beginning of the proof of Theorem~\ref{th:definable=borel}, the function \( f_T \colon \pre{\omega}{2} \to \CT_\kappa \subseteq \Mod^\kappa_\LL \) defined in~\eqref{eq:f_T} is easily seen to be continuous when both spaces are endowed with the product topology \( \tau_p \), so for every \( U \) in the canonical basis \( \mathcal{B}_p(\Mod^\kappa_\LL) \) for the product topology we have
\[ 
f_T^{-1} ( U ) \in \bSigma^0_1( \pre{\omega}{2} ) \subseteq \bDelta^1_1 ( \pre{\omega}{2} ) \subseteq \bDelta_{\bS ( \kappa )} ( \pre{\omega}{2} ) .
\] 
By Propositions~\ref{prop:inthecodesSouslinAC} and~\ref{prop:inthecodesSouslin}, under our assumptions \( f_T \) is then \( \bS ( \kappa ) \)-in-the-codes. 
Since \( f_T \) reduces \( R \) to \( \embeds_\CT^\kappa \) by Theorem~\ref{th:graphs}\ref{th:graphs-a}, we get \( R \leq_{\bS ( \kappa )} {\embeds_\CT^\kappa} \), as required.
\end{proof}

We now consider some results concerning invariant universality which generalize Theorems~\ref{th:mottorosfriedman} and~\ref{th:mottorosfriedman2} to \emph{regular} uncountable \( \kappa \)'s.
Given an \( \LL_{\kappa^+ \kappa} \)-sentence \( \upsigma \), say that a map \( h \colon \Mod^\kappa_\upsigma \to \pre{ \omega}{2} \) is \markdef{\( \LL_{\kappa^+ \kappa} \)-measurable} if for every open set \( U \subseteq \pre{ \omega}{2} \) (equivalently, for every basic open set \( U = \Nbhd^\omega_s \) with \( s \in \pre{<\omega}{2} \)) there is an \( \LL_{\kappa^+ \kappa} \)-sentence \( \upsigma_U \) such that \( h^{-1} ( U ) = \Mod^\kappa_{\upsigma_U} \cap \Mod^\kappa_{\upsigma} = \Mod^\kappa_{\upsigma_U \wedge \upsigma} \). 
Similarly, we can define the notion of \markdef{\( \LL^b_{\kappa^+ \kappa} \)-measurability}. 

Notice that by Corollary~\ref{cor:formulaborel1}\ref{cor:formulaborel1-ii}, if \( \kappa \) is regular and \( h \colon \Mod^\kappa_\upsigma \to \pre{ \omega}{2} \) is \( \LL^b_{\kappa^+ \kappa} \)-measurable, then \( h \) is also (effective) \( \kappa+1 \)-Borel measurable.

\begin{lemma} \label{lem:inverse}
The map \( h_T \colon \Mod^\kappa_{\upsigma_T} \to \pre{ \omega}{2} \) from Definition~\ref{def:defh_T} is \( \LL^b_{\kappa^+ \kappa} \)-measurable, and also (effective) \( \kappa + 1 \)-Borel measurable when \( \kappa \) is regular,.
\end{lemma}

\begin{proof}
Let \( f_T \) be the function defined in~\eqref{eq:f_T}, and \( \Uppsi \) be the \( \LL_{\kappa^+ \kappa} \)-sentence defined in~\eqref{eq:Psi}. 
Then for every open set \( U \subseteq \pre{\omega}{2} \) 
\[
h_T^{-1} ( U ) = \setofLR{ X \in \Mod^\kappa_{\Uppsi} }{ \exists x \in U ( X \cong f_T ( x ) ) } = g^{-1} ( ( g \circ f_T ) ( U ) ) .
\] 
Since \( C \coloneqq \pre{\omega}{2} \setminus U \) is closed, \( ( g \circ f_T ) ( C ) \) is closed as well by Corollary~\ref{cor:closed}, so let \( \upsigma_V \) be the \( \LL^b_{\kappa^+ \kappa} \)-sentence given by Lemma~\ref{lem:open} applied to \( V \coloneqq \pre{\pre{<\omega}{2} \times \kappa}{2} \setminus ( g \circ f_T ) ( C ) \). 
Then \( h_T^{-1} ( U ) = g^{-1} ( ( g \circ f_T ) ( U ) ) = g^{-1} ( V ) \cap \Mod^\kappa_{\upsigma_T} = \Mod^\kappa_{\upsigma_V} \cap \Mod^\kappa_{\upsigma_T} \).
\end{proof}

The next theorem generalizes Theorem~\ref{th:mottorosfriedman} to \emph{regular} uncountable \( \kappa \)'s.

\begin{theorem} \label{th:maintopology}
Let \( \kappa \) be an infinite \emph{regular} cardinal. 
Then \( \embeds_\CT^\kappa \) is \emph{\( \leq^\kappa_{\bB} \)-invariantly universal} for \( \kappa \)-Souslin quasi-orders on \( \pre{\omega}{2} \), that is: for every \( \kappa \)-Souslin quasi-order \( R \) on \( \pre{ \omega}{2} \) there is an \( \LL_{\kappa^+ \kappa} \)-sentence \( \upsigma \) all of whose models are combinatorial trees such that \( R \sim^\kappa_{\bB} {\embeds^\kappa_\upsigma} \), and in fact even \( R \sim^\kappa_{\mathrm{B}} {\embeds^\kappa_\upsigma} \).
\end{theorem}

\begin{proof} 
Let \( T \in \TT_\kappa \) be such that \( R = \PROJ \body{ T } \), let \( \upsigma_T \), \( f_T \) and \( h_T \) be as in Theorem~\ref{th:main}, and set \( \upsigma \coloneqq \upsigma_T \). 
Since by Theorem~\ref{th:main} we already know that \( f_T \) and \( h_T \) reduce \( R \) and \( \embeds^\kappa_\upsigma \) to each other, it remains to show that they are (weakly) \( \kappa + 1 \)-Borel functions: arguing as in the proof of Theorem~\ref{th:definable=borel}, we get that the function \( f_T \) is (effective) weakly \( \kappa + 1 \)-Borel (and in fact also effective weakly \( \alpha \)-Borel for any \( \alpha \geq \omega \)), while the function \( h_T \) is (effective) \( \kappa + 1 \)-Borel by regularity of \( \kappa \) and Lemma~\ref{lem:inverse}. 
\end{proof}

Notice that the \( \LL_{\kappa^+ \kappa} \)-sentence \( \upsigma \) we provided in the proof of Theorem~\ref{th:maintopology} (which is the \( \LL_{\kappa^+ \kappa} \)-sentence \( \upsigma_T \) from Theorem~\ref{th:main}, see also Corollary~\ref{cor:saturation}) does not belong to the fragment \( \LL^b_{\kappa^+ \kappa} \),%
\footnote{More precisely: the sentence \( \upsigma_T \) provided in the proof of Corollary~\ref{cor:saturation} is the conjunction of an \( \LL^b_{\kappa^+ \kappa} \)-sentence with the \( \Uppsi \in \LL_{\kappa^+ \kappa} \) from~\eqref{eq:Psi}. 
Since the latter does not belong to the bounded logic \( \LL^b_{\kappa^+ \kappa} \), the same applies to \( \upsigma_T \).} 
and hence we cannot in general guarantee that \( \Mod^\kappa_\upsigma \) be a \( \kappa+1 \)-Borel subset of \( \Mod^\kappa_\LL \) (see Section~\ref{subsec:LopezEscobar}, and in particular Remark~\ref{rmk:formulaborel}). 
Therefore, if one aims at generalizing Theorem~\ref{th:mottorosfriedman2} to (regular) uncountable \( \kappa \)'s, then Theorem~\ref{th:main} needs to be replaced with Theorem~\ref{th:barmain} (and \( \CT_\kappa \) with \( \OCT_\kappa \)). 
First notice that by replacing in the proof of Lemma~\ref{lem:inverse} the map \( f_T \) with the function \( \bar{f}_T \) from~\eqref{eq:defbarf_T}, \( \Uppsi \) with the \( \bar{\LL}^b_{\kappa^+ \kappa} \)-sentence \( \bar{\Uppsi} \) from~\eqref{eq:barPsi}, Corollary~\ref{cor:closed} with Corollary~\ref{cor:barclosed}, and Lemma~\ref{lem:open} with Lemma~\ref{lem:baropen}, we get the following variant of it.

\begin{lemma} \label{lem:barinverse}
The map \( \bar{h}_T \colon \Mod^\kappa_{\bar\upsigma_T} \to \pre{ \omega}{2} \) from Definition~\ref{def:defbarh_T} is \( \bar\LL^b_{\kappa^+ \kappa} \)-measurable, and hence, if \( \kappa \) is regular, also (effective) \( \kappa + 1 \)-Borel measurable.
\end{lemma}

Then arguing as in the proof of Theorem~\ref{th:maintopology} we get the following generalization of Theorem~\ref{th:mottorosfriedman2}.

\begin{theorem} \label{th:barmaintopology}
Let \( \kappa \) be an infinite \emph{regular} cardinal. 
Then for every \( \kappa \)-Souslin quasi-order \( R \) on \( \pre{ \omega}{2} \) there is an effective \( \kappa+1 \)-Borel set \( B \subseteq \Mod^\kappa_{\bar{\LL}} \) closed under isomorphism (all of whose elements are ordered combinatorial trees) such that \( R \sim^\kappa_{\bB} {\embeds \restriction B} \), and in fact even \( R \sim^\kappa_{\mathrm{B}} {\embeds \restriction B} \).
\end{theorem}

\begin{proof} 
Let \( T \in \TT_\kappa \) be such that \( R = \PROJ \body{ T } \), let \( \bar\upsigma_T \), \( \bar{f}_T \) and \( \bar{h}_T \) be as in Theorem~\ref{th:barmain}, and set \( \bar\upsigma \coloneqq \bar\upsigma_T \). 
Then \( B \coloneqq \Mod^\kappa_{\bar\upsigma_T} \subseteq \Mod^\kappa_{\bar\LL}\) is effective \( \kappa+1 \)-Borel by~\eqref{eq:barupsigmaisBorel}.
Since by Theorem~\ref{th:barmain} we know that \( \bar{f}_T \) and \( \bar{h}_T \) reduce \( R \) and \( \embeds^\kappa_{\bar\upsigma} \) to each other, it remains to show that they are both (weakly) \( \kappa + 1 \)-Borel functions. 
For the function \( \bar{f}_T \) we argue as in the proof of Theorem~\ref{th:definable=borel}: \( \bar{f}_T \colon \pre{\omega}{2} \to \Mod^\kappa_{\bar\LL} \) is easily seen to be continuous when both spaces are endowed with the product topology \( \tau_p \), and hence it is (effective) weakly \( \kappa +1 \)-Borel (in fact, even effective weakly \( \alpha\)-Borel for any \( \alpha \geq \omega \)). 
The function \( \bar{h}_T \) is effective \( \kappa + 1 \)-Borel by regularity of \( \kappa \) and Lemma~\ref{lem:barinverse}, so we are done. 
\end{proof}

\subsection{Absolutely definable reducibilities}\label{subsec:absolute} 

Even if Borel reducibility (i.e.~\( \bSigma^1_1 \)-reducibility) is probably the most natural way to compare the complexity of \( \bSigma^1_1 \) equivalence relations and quasi-orders, one sometimes needs to generalize this notion to the projective levels: for example, in~\cite{Hjorth:1995qf, Hjorth:2000zr} the so-called \emph{absolutely \( \bDelta^1_2 \)-reducibility} \( \leq_{\mathrm{a}\bDelta^1_2} \), i.e.~the reducibility notion obtained using absolutely \( \bDelta^1_2 \) maps as reductions, has proven to be useful in this area. 
Here we recall the definition of such functions as presented in~\cite{Hjorth:2000zr}.

\begin{definition}[{\cite[Definition 9.1]{Hjorth:2000zr}}] \label{def:hjorth}
Let \( \HC \) be the collection of hereditarily countable sets. 
A function \( f \colon \pre{\omega}{\omega} \to \HC \) is \emph{absolutely \( \bDelta^1_2 \)} if there is some parameter \( p \in \pre{\omega}{\omega} \) and an \( \LST \)-formula \( {\sf\Psi} ( x_0 , x_1 , z_0 ) \) such that:
\begin{enumerate-(1)}
\item \label{def:hjorth-1}
for all \( x \in \pre{\omega}{\omega} \) and \( y \in \HC \), \( f ( x ) = y \) if and only if there is some countable transitive set \( \mathcal{M} \) with \( x , y , p \in \mathcal{M}\)  such that \( ( \mathcal{M}, \in ) \models {\sf \Psi} [ x , y , p ] \);
\item \label{def:hjorth-2}
\( {\sf\Psi} ( x_0 , x_1 , p ) \) \emph{absolutely} defines a function, in the sense that in all generic extensions of the universe \( \Vv \) it continues to be the case that for all \( x \in \pre{\omega}{\omega} \) there exists\footnote{Such \( y \) is unique by \( \Sigma^1_2 \)-absoluteness.} \( y \in \HC \) and a countable transitive \( \mathcal{M} \) with \( (\mathcal{M}, \in ) \models {\sf\Psi} [ x , y , p ] \).
\end{enumerate-(1)} 
\end{definition}
Such definition is then naturally extended to cover all functions of interest in~\cite[Section 9]{Hjorth:2000zr}, including functions between arbitrary Polish spaces, functions between \( \omega_1 \) and \( \HC \), and so on. 
\begin{remark}\label{rmk:coherent}
The reference to the countable transitive \( \mathcal{M} \) in Definition~\ref{def:hjorth} is added to have a \( \bDelta^1_2 \) definition of \( f \), but if we are only interested in the \emph{absolute definability} of \( f \) (without specifying the complexity of such a definition) then it is natural to only require  that the \( \LST \)-formula \( {\sf\Psi}(x_0,x_1,p) \) defines (in \( \Vv \)) the graph of \( f \colon \pre{\omega}{\omega} \to \HC \), and it continues to define a function between \( \pre{\omega}{\omega} \) and \( \HC \) in all generic extensions of the universe \( \Vv \) in a coherent way, that is, if \( x \in \Vv [ G_0 ] \cap \Vv [ G_1 ] \), then \( \Vv [ G_0 ] \models {\sf\Psi} [ x , y ,p ] \) and \( \Vv [ G_1 ] \models {\sf\Psi} [ x , y ,p ] \) for the same \( y \).
\end{remark}

In the realm of standard Borel spaces, (\( \omega+1 \)-)Borel functions are absolutely \( \bDelta^{1}_{2} \)-definable, so \( \leq_{\mathrm{a}\bDelta^1_2} \) is used only whenever suitable \( \bDelta^{1}_{1} \)-reductions are not available.
In contrast, when \( \kappa \) is uncountable \( \kappa + 1 \)-Borel functions need not be absolute, so one might wonder whether the reductions appearing in our main results are indeed absolute.
In this section we generalize Hjorth's approach and show that, essentially, the reductions obtained in the proof of Theorem~\ref{th:main} cannot be destroyed by passing to set-forcing extensions or to inner models. 
Since these are the two main techniques for proving the independence of a given mathematical assertion from the chosen set-theoretic axiomatization, Theorems~\ref{th:absolutedefin} and~\ref{th:absolutedefinold} essentially show that the invariant universality of the embeddability relation \( \embeds^\kappa_\CT \) is absolute for transitive models of \( \ZF \) containing all relevant parameters.

Before stating the main results of this section (Theorems~\ref{th:absolutedefin} and~\ref{th:absolutedefinold}), we first need to adapt Definition~\ref{def:hjorth} to our context.
As done by Hjorth, also Definition~\ref{def:absolute} is given just for the specific functions which are relevant to the results of this paper. 
However, with some extra work such a definition could be easily adapted to a general definition of an absolutely definable function \( f \) between (definable subsets) of spaces \( \mathcal{X}, \mathcal{Y} \) of type, respectively, \( \lambda,\mu \in \Cn \).

\begin{definition} \label{def:absolute}
Let \( \kappa \) be an infinite cardinal and \( \upsigma \in \LL_{\kappa^+ \kappa} \). 
A function \( f \colon \pre{\omega}{2} \to \Mod^\kappa_\upsigma \) is \markdef{absolutely (\( p \)-)definable} if there is a parameter \( p \in \pre{\kappa}{\kappa} \) and an \( \LST \)-formula \( {\sf \Psi}_f ( x_0 , x_1, z_0, z_{1}) \) such that for all forcing notions \( \forcing{P}_0, \forcing{P}_1 \in \Vv \) and all \( \Vv \)-generic \( G_i \subseteq \forcing{P}_i \) such that \( \kappa \in \Cn ^{\Vv [ G_i ] } \), the following conditions hold:
\begin{enumerate-(1)}
\item 
for all \( x \in \pre{\omega}{2} \) and \( y \in \Mod^\kappa_\upsigma \), \( f ( x ) = y \) if and only if \( \Vv \models {\sf \Psi}_f [ x , y , \kappa , p ] \);
\item
\( \forcing{P}_i \) forces that \( {\sf \Psi}_f ( x_0 , x_1 , \kappa , p ) \) defines (the graph of) a function \( f^{\Vv [ G_i ]} \colon ( \pre{\omega}{2} )^{\Vv [ G_i ]} \to ( \Mod^\kappa_\upsigma )^{\Vv [ G_i ]} \), that is 
\[
\Vv [ G_i ] \models \forall x \in \pre{\omega}{2} \exists ! X \in \Mod^\kappa_\upsigma {\sf \Psi}_f [ x , X , \kappa , p ] 
\]
\item
\( f^{\Vv [ G_0 ]} \) and \( f^{\Vv [ G_1 ]} \) are coherent, that is for all \( x \in ( \pre{\omega}{2} )^{\Vv [ G_0 ]} \cap ( \pre{\omega}{2} )^{\Vv [ G_1 ]} \)
\begin{equation}\label{eq:coherence}
f^{\Vv [ G_0 ]} ( x ) = f^{\Vv [ G_1 ]} ( x ) .
\end{equation}
\end{enumerate-(1)}%
Absolutely (\( p \)-)definable functions \( \Mod^ \kappa _{\upsigma}\to \pre{\omega}{2} \) are defined similarly.
\end{definition}

\begin{remarks} 
\begin{enumerate-(i)}
\item
The coherence condition~\eqref{eq:coherence} does not explicitly appear in Hjorth's Definition~\ref{def:hjorth} since in that case it is automatically satisfied---see Remark~\ref{rmk:coherent}.
\item
The restriction of \( f^{\Vv [ G ]} \) to \( (\pre{\omega}{2})^{\Vv [ G ]} \setminus (\pre{\omega}{2})^{\Vv} \) depends in general on the chosen \( p \) and \( {\sf \Psi}_f \) (but \( f^{\Vv [ G ]} \restriction (\pre{\omega}{2})^{\Vv} \) does not by~\eqref{eq:coherence}). 
Similar considerations hold for an absolutely definable \( f \colon \Mod^\kappa_\upsigma \to \pre{\omega}{2} \).
\end{enumerate-(i)}
\end{remarks}

Using absolutely definable functions as reductions (in \( \Vv \)) between a quasi-order \( R \) on \( \pre{\omega}{2} \) and an embeddability relation \( \embeds^\kappa_\upsigma \) (for some \( \LL_{\kappa^+ \kappa} \)-sentence \( \upsigma \)) would yield an analogue of Hjorth's absolutely \( \bDelta^1_2 \) reducibility, which may be dubbed \emph{absolutely definable reducibility}. 
However, when \( R \) is a \( \kappa \)-Souslin quasi-order and \( T \in \TT_\kappa \) is a \emph{faithful representation} of it (see Definition~\ref{def:faithfulrepresentation}), then by Remark~\ref{rmk:faithfulrepresentation} in all forcing extensions \( \Vv [ G ] \) of \( \Vv \) we have a canonical extension \( R^{\Vv [ G ]}_T \coloneqq (\PROJ \body{T})^{\Vv [ G ]} \) of \( R \) which is still a \( \kappa \)-Souslin quasi-order and is coherent with \( R \), i.e.~\( R^{\Vv [ G ]}_T \) coincides with \( R \) on their common domain \( (\pre{\omega}{2})^{\Vv} \). 
Therefore when comparing such an \( R \) with an embeddability relation \( \embeds^\kappa_\upsigma \), it is natural to require that the absolutely definable functions involved continue to be reductions between \( R^{\Vv [ G ]}_T \) and \( (\embeds^\kappa_\upsigma)^{\Vv[G]} \) in all forcing extensions as in Definition~\ref{def:absolute} (and not just in \( \Vv \)). This leads us to the following stronger definition. 

\begin{definition} \label{def:absolutelydefinablereducibility}
Let \( \kappa \) be an infinite cardinal, \( T \in \TT_\kappa \) be a faithful representation of a \( \kappa \)-Souslin quasi-order \( R \), and \( \upsigma \) be an \( \LL_{\kappa^+ \kappa} \)-sentence. 
We say that \( R \) is \markdef{absolutely (\( p\)-)definably reducible} to \( \embeds^\kappa_\upsigma \), in symbols
\[
 R \leq_{\mathrm{aD}} {\embeds^\kappa_\upsigma} \index[symbols]{121@\( \leq_{\mathrm{aD}} \), \( \sim_{\mathrm{aD}} \)},
\]
if there is a reduction \( f \colon \pre{\omega}{2} \to \Mod^\kappa_\upsigma \) of \( R \) to \( \embeds^\kappa_\upsigma \) such that for some parameter \( p \in \pre{\kappa}{\kappa} \) and some \( \LST \)-formula \( {\sf\Psi}_f ( x_0 , x_1 , z_0 , z_1 ) \) the following hold:
\begin{enumerate-(1)}
\item 
\( {\sf\Psi}_f ( x_0 , x_1 , \kappa , p ) \) absolutely defines \( f \);
\item
for all forcing notions \( \forcing{P} \) and all \( \Vv \)-generic \( G \subseteq \forcing{P} \) such that \( \kappa \in \Cn ^{\Vv [ G ] } \)
\[ 
\Vv[G] \models \text{``} f^{\Vv[G]} \text{ reduces } R^{\Vv[G]}_T \text{ to } \embeds^\kappa_\upsigma \text{''},
 \] 
where as in Definition~\ref{def:absolute} we denote by \( f^{\Vv[G]} \) the map defined in \( \Vv[G] \) by \( {\sf\Psi}_f ( x_0 , x_1 , \kappa , p ) \). 
\end{enumerate-(1)}
The notions of absolutely definable reducibility of \( \embeds^\kappa_\upsigma \) to \( R \) and of \markdef{absolutely definable bi-reducibility} \( \sim_{\mathrm{aD}} \) are defined similarly.
\end{definition}

We are now ready to reformulate our main invariant universality result Theorem~\ref{th:main} in terms of absolutely definable reducibility.

\begin{theorem} \label{th:absolutedefin}
Let \( \kappa \) be an infinite cardinal. 
Then \( \embeds^\kappa_\CT \) is \emph{\( \leq_{\mathrm{aD}} \)-invariantly universal} for \( \kappa \)-Souslin quasi-orders on \( \pre{\omega}{2} \).

More precisely: for every faithful representation \( T \in \TT_\kappa \) of a \( \kappa \)-Souslin quasi-order \( R \) on \( \pre{\omega}{2} \) there is an \( \LL_{\kappa^+ \kappa} \)-sentence \( \upsigma \) and two \( \LST \)-formul\ae{} \( {\sf\Psi}_{f_T} ( x_0 , x_1 , z_0 , z_1 ) \) and \( {\sf\Psi}_{h_T} ( x_0 , x_1 , z_0 , z_1 ) \) such that:
\begin{enumerate-(a)}
\item \label{th:absolutedefin-b}
\( {\sf\Psi}_{f_T} ( x_0 , x_1 , \kappa , T ) \) absolutely defines a function \( f_T \colon \pre{\omega}{2} \to \Mod^\kappa_\upsigma \) witnessing \( R \leq_{\mathrm{aD}} {\embeds^\kappa_\upsigma} \);%
\footnote{Here \( T \) is identified with a parameter in \( \pre{\kappa}{\kappa} \) via its characteristic functions and the bijection \( \Code{\cdot} \) from~\eqref{eq:codingfinitesequenceofordinals}. 
The same identification is tacitly applied every time that \( T \) is used as a parameter in an \( \LST \)-formula absolutely defining a function.}
\item \label{th:absolutedefin-c}
\( {\sf\Psi}_{h_T}(x_0,x_1,\kappa,T) \) absolutely defines a function \( h_T \colon \Mod^\kappa_\upsigma \to \pre{\omega}{2} \) witnessing \( {\embeds^\kappa_\upsigma} \leq_{\mathrm{aD}} R \);
\item \label{th:absolutedefin-d}
for all forcing notions \( \forcing{P} \) and all \( \Vv \)-generic \( G \subseteq \forcing{P} \) such that \( \kappa \in \Cn ^{\Vv [ G ] } \), we have
\( h_T^{\Vv [ G ] } \circ f_T^{\Vv [ G ] } = \id^{\Vv [ G ] } \) and \( {\Vv [ G ] } \models ( f^{\Vv [ G ] }_T \circ h^{\Vv [ G ] }_T ) ( X ) \cong X \) for every \( X \in ( \Mod^\kappa_\upsigma )^{\Vv [ G ] } \).
\end{enumerate-(a)}
\end{theorem}

\begin{proof}
Let \( \upsigma \) be the \( \LL_{\kappa^+ \kappa} \)-sentence obtained in the proof of Corollary~\ref{cor:saturation}, \( {\sf\Psi}_{f_T} ( x_0 , x_1 , \kappa , T ) \) be the \( \LST \)-formula from Fact~\ref{prop:f_Tdefinable}, and \( {\sf \Psi}_{h_T} ( x_0 , x_1 , z_0 , z_1 ) \) be the \( \LST \)-formula from Proposition~\ref{prop:h_Tisabsolute}\ref{prop:h_Tisabsolute-a}: we claim that \( \upsigma \), \( {\sf\Psi}_{f_T} \), and \( {\sf \Psi}_{h_T} \) are as required. 

Fix any forcing notion \( \forcing{P} \) and a \( \Vv \)-generic \( G \subseteq \forcing{P} \) such that \( \kappa \in \Cn ^{\Vv [ G ] } \).
Observe that if \( \upsigma^{\Vv[G]}_T \coloneqq ( \upsigma_T )^{\Vv[G]} \) is the \( ( \LL_{\kappa^+ \kappa} )^{\Vv[G]} \)-sentence coming from (the proof of) Corollary~\ref{cor:saturation} when applied in \( \Vv[G] \) (see Section~\ref{subsec:moreabsoluteness}),
then 
\begin{equation} \label{eq:upsigmaisthesame} 
\upsigma = \upsigma^{\Vv[G]}_T 
\end{equation} 
by Proposition~\ref{prop:upsigma_Tisabsolute}.
By Fact~\ref{prop:f_Tdefinable}, the \( \LST \)-formula \( {\sf\Psi}_{f_T} ( x_0 , x_1 , \kappa , T ) \) absolutely defines \( f_T \) in the strong sense that: the function \( f_T^{\Vv[G]} \) defined by \( {\sf\Psi}_{f_T} ( x_0 , x_1 , \kappa , T ) \) in \( \Vv[G] \) is exactly the function \( (f_T)^{\Vv[G]} \) obtained as in~\eqref{eq:f_T} once all the construction is carried out in \( \Vv[G] \). 
In particular, since Theorem~\ref{th:graphs}\ref{th:graphs-a} holds in \( \Vv[G] \), condition~\ref{th:absolutedefin-b} of the present theorem is satisfied. 
By Proposition~\ref{prop:h_Tisabsolute}, the \( \LST \)-formula \( {\sf\Psi}_{h_T}(x_0,x_1, \kappa , T) \) absolutely defines \( h_T \) (the coherence condition~\eqref{eq:coherence} in Definition~\ref{def:absolute} is guaranteed by Proposition~\ref{prop:h_Tisabsolute}\ref{prop:h_Tisabsolute-b}, whose extra condition can be removed in our setup by Remark~\ref{rmk:newrmkabsolute}). More precisely, when evaluated in \( \Vv[G] \) the formula \( {\sf\Psi}_{h_T}(x_0,x_1, \kappa , T) \) defines the function \( h_T^{\Vv[G]} = ( h_T )^{\Vv[G]} \) computed (in \( \Vv[G] \)) according to Definition~\ref{def:defh_T}, whose domain is \( ( \Mod^\kappa_{\upsigma^{\Vv[G]}_T} )^{\Vv[G]} = ( \Mod^\kappa_\upsigma )^{\Vv[G]} \) (the latter equality follows from~\eqref{eq:upsigmaisthesame}). 
Since Theorem~\ref{th:main} holds in \( \Vv[G] \), both conditions~\ref{th:absolutedefin-c} and~\ref{th:absolutedefin-d} of the present theorem are satisfied.
\end{proof}

\begin{remarks} \label{rmk:elementaryclassrelativizedVG}
\begin{enumerate-(i)}
\item \label{rmk:elementaryclassrelativizedVG-i}
Recall that when \( \kappa > \omega \) neither the embeddability relation \( \embeds \restriction \Mod^\kappa_\LL \) nor the statement ``\( X \models \upsigma \)'' (for \( X \in \Mod^\kappa_\LL \) and \( \upsigma \) an arbitrary \( \LL_{\kappa^+ \kappa} \)-sentence) are in general absolute for transitive models of \( \ZF \) (see the observation before Proposition~\ref{prop:embeddabilityabsolute} and Remark~\ref{rmk:satisfactionisnotabsolute}\ref{rmk:satisfactionisnotabsolute-iii}). 
However, for the specific \( \LL_{\kappa^+ \kappa} \)-sentences considered in Theorem~\ref{th:absolutedefin} (that is, for those sentences constructed as in the proof of Corollary~\ref{cor:saturation}) one can show, using Proposition~\ref{prop:embeddabilityabsolute}, that the embeddability relations \( (\embeds^\kappa_\upsigma)^{\Vv} \) and \( (\embeds^\kappa_\upsigma)^{\Vv[G]} \) are coherent in the following strong sense:
\begin{enumerate-(1)}
\item \label{th:absolutedefinVG-e}
\( ( \Mod^\kappa_\upsigma )^{\Vv} = ( \Mod^\kappa_\upsigma )^{\Vv[G]} \cap \Vv \);
\item \label{th:absolutedefinVG-f}
\( ( \embeds^\kappa_\upsigma )^{\Vv[G]} \) and \( ( \embeds^\kappa_\upsigma )^{\Vv} \) coincide on their common domain \( (\Mod^\kappa_\upsigma)^{\Vv} \), that is: for every \( X,Y \in \Mod^\kappa_\upsigma \)
\[ 
\Vv \models {X \embeds Y} \IFF \Vv[G] \models {X \embeds Y}. 
 \] 
\end{enumerate-(1)}
\item \label{rmk:elementaryclassrelativizedVG-ii}
It follows from the absoluteness results above that if \( \upsigma \) is an \( \LL_{\kappa^+ \kappa} \)-sentence obtained as in Corollary~\ref{cor:saturation} and \( \forcing{P} \) does not add new reals, then every \( X \in (\Mod^\kappa_{\upsigma_T})^{\Vv[G]} \) has an isomorphic copy in the ground model \( \Vv \). 
\end{enumerate-(i)}
\end{remarks}

So far we just considered generic (or upward) absoluteness, that is absoluteness with respect to forcing extensions. However, similar results hold when considering e.g.\ absoluteness with respect to transitive \( \ZF \)-models \( M \subseteq \Vv \) (thus including inner models as a particular case). In this new setup, Definition~\ref{def:absolute} might be reformulated as follows.

\begin{definition} \label{def:absoluteinner}
Let \( \kappa \) be an infinite cardinal and \( \upsigma \in \LL_{\kappa^+ \kappa} \). 
A function \( f \colon \pre{\omega}{2} \to \Mod^\kappa_\upsigma \) is \markdef{(downward) absolutely (\( p \)-)definable} if there is a parameter \( p \in \pre{\kappa}{\kappa} \) and an \( \LST \)-formula \( {\sf \Psi}_f ( x_0 , x_1 , z_0 , z_{1} ) \) such that for every transitive \( \ZF \)-model \( M \subseteq \Vv \) with \( \kappa,p,\upsigma \in M \), the following hold:
\begin{enumerate-(1)}
\item
\( {\sf \Psi}_f(x_0,x_1,\kappa,p) \) defines in \( M \) (the graph of) a function \( f^M \colon (\pre{\omega}{2})^M \to (\Mod^\kappa_\upsigma)^M \), that is
\[ 
(M, \in) \models \FORALL{ x \in \pre{\omega}{2} } \EXISTSONE{ X \in \Mod^\kappa_\upsigma } {\sf\Psi}_f[x,X,\kappa,p] ;
\]
\item \label{def:absoluteinner-2}
for every \( x \in \pre{\omega}{2} \cap M \) and every \( X \in \Mod^\kappa_\upsigma \cap M \), 
\[ 
f(x) = X \IFF (M, \in) \models {\sf\Psi}_f [ x , X , \kappa , p ] .
 \] 
\end{enumerate-(1)}
\end{definition}
Notice that the fact that \( {\sf\Psi} \) defines \( f \) in \( \Vv \) is included in part~\ref{def:absoluteinner-2} by taking \( M = \Vv \). 
 The same condition also yields the analogue of the coherence condition~\eqref{eq:coherence} in Definition~\ref{def:absolute}: indeed, each \( f^M \) is required to be the restriction of \( f \) to \( (\pre{\omega}{2})^M \). 
 
Using Definition~\ref{def:absoluteinner}, one could then introduce a notion of \markdef{(downward) absolute definable reducibility} \( \leq^{\mathrm{dw}}_{\mathrm{aD}} \) mirroring in the natural way Definition~\ref{def:absolutelydefinablereducibility} --- the unique thing to notice here is that in this case \( R^M_T \) will be the restriction of the quasi-order \( R \) to \( (\pre{\omega}{2})^M \). 
Theorem~\ref{th:absolutedefin} would then turn into the following ``downward absoluteness'' result, which can be proved in the same vein. 
(Part~\ref{th:absolutedefinold-a}, which must explicitly be added here because \( M \) is now a subclass of \( \Vv \) which in principle might not contain the desired \( \upsigma \), follows from the analogue of~\eqref{eq:upsigmaisthesame}.)

\begin{theorem} \label{th:absolutedefinold}
Let \( \kappa \) be an infinite cardinal. 
Then \( \embeds^\kappa_\CT \) is \emph{\( \leq^{\mathrm{dw}}_{\mathrm{aD}} \)-invariantly universal} for \( \kappa \)-Souslin quasi-orders on \( \pre{\omega}{2} \).

More precisely: for every faithful representation \( T \in \TT_\kappa \) of a \( \kappa \)-Souslin quasi-order \( R \) on \( \pre{\omega}{2} \) there is an \( \LL_{\kappa^+ \kappa} \)-sentence \( \upsigma \) and two \( \LST \)-formul\ae{} \( {\sf\Psi}_{f_T} ( x_0 , x_1 , z_0 , z_1 ) \) and \( {\sf\Psi}_{h_T} ( x_0 , x_1 , z_0 , z_1 ) \) such that for every transitive \( \ZF \)-model \( M \subseteq \Vv \) with \( T \in M \) (and hence \( \kappa \in M \)) the following conditions hold:
\begin{enumerate-(a)}
\item \label{th:absolutedefinold-a}
\( \upsigma \in (\LL_{\kappa^+ \kappa})^M \) (equivalently, by absoluteness of the statement ``\( \upsigma \in \LL_{\kappa^+ \kappa} \)'', \( \upsigma \in M \));
\item \label{th:absolutedefinold-b}
\( {\sf\Psi}_{f_T} ( x_0 , x_1 , \kappa , T ) \) (downward) absolutely defines a function \( f_T \colon \pre{\omega}{2} \to \Mod^\kappa_\upsigma \) witnessing \( R \leq^{\mathrm{dw}}_{\mathrm{aD}} {\embeds^\kappa_\upsigma} \);
\item \label{th:absolutedefinold-c}
\( {\sf\Psi}_{h_T}(x_0,x_1,\kappa,T) \) (downward) absolutely defines a function \( h_T \colon \Mod^\kappa_\upsigma \to \pre{\omega}{2} \) witnessing \( {\embeds^\kappa_\upsigma} \leq^{\mathrm{dw}}_{\mathrm{aD}} R \);
\item \label{th:absolutedefinold-d}
\( h_T^M \circ f_T^M = \id^M \) and \( M \models ( f^M_T \circ h^M_T ) ( X ) \cong X \) for every \( X \in ( \Mod^\kappa_\upsigma )^M \).
\end{enumerate-(a)}
\end{theorem}

\begin{remarks} \label{rmk:elementaryclassrelativized}
Considerations similar to those in Remark~\ref{rmk:elementaryclassrelativizedVG} can be made also in this new context. If \( M \subseteq \Vv \) is any transitive \( \ZF \)-model and \( \upsigma_T \) is an \( \LL_{\kappa^+ \kappa} \)-sentence obtained as in Corollary~\ref{cor:saturation} starting from some \( T \in \TT_\kappa \cap M \), then
\begin{enumerate-(i)}
\item \label{rmk:elementaryclassrelativized-i}
\( ( \Mod^\kappa_\upsigma )^M = ( \Mod^\kappa_\upsigma )^{\Vv} \cap M \) and for every \( X,Y \in ( \Mod^\kappa_\upsigma )^{M} \subseteq \Mod^\kappa_\upsigma \)
\[ 
\Vv \models {X \embeds Y} \IFF M \models {X \embeds Y}. 
 \] 
\item \label{rmk:elementaryclassrelativized-ii}
if \( M \) contains all the reals of \( \Vv \), then every \( X \in \Mod^\kappa_{\upsigma_T} \) has an isomorphic copy in \( M \), that is, \( (\Mod^\kappa_{\upsigma_T})^M \) meets all isomorphism classes of \( \Mod^\kappa_{\upsigma_T} \). 
In particular, \( \Mod^\kappa_{\upsigma_T} \) is always faithfully represented in every inner model containing \( \Ll ( \R ) \).
\end{enumerate-(i)} 
\end{remarks}

We conclude this section by pointing out that we just considered absoluteness with respect to generic extensions or transitive \( \ZF \)-submodels of \( \Vv \) only for the sake of simplicity. 
However, it is easily seen that our absoluteness results naturally extend to much wider contexts, such as the generic multiverse investigated by several people (including J.D.~Hamkins, S.D.~Friedman, W.H.~Woodin). 
Unfortunately, it is difficult (if not impossible) to coherently define a unique setup encompassing all these possibilities, so we explicitly considered only two of the most relevant setups.

\subsection{Reducibilities in an inner model} \label{subsubsec:reducibilityininnermodel} 
As a further generalization of absolutely \( \bDelta^1_2 \)-reducibility, one could consider reducibility under determinacy within \( \Ll ( \R ) \); in other words, assuming \( \AD^{\Ll ( \R ) } \) and working inside \( \Ll ( \R ) \) we consider arbitrary reductions between quasi-orders. 
The rationale for this is that \( \Ll ( \R ) \) includes anything one could consider reasonably definable, like e.g.~all Polish spaces up to homeomorphism, all separable Banach spaces up to linear isometry, all Borel and projective sets and functions, and so forth. 
Moreover \( \AD^{\Ll ( \R )} \) yields a very detailed picture of \( \Ll ( \R ) \), much akin to what the axiom of constructibility \( \Vv = \Ll \) does for \( \Ll \). 
This approach to definable reducibility has been explicitly considered e.g.~in~\cite{Hjorth:1995ve, Hjorth:1999bh} and~\cite[Chapter 9]{Hjorth:2000zr}.

\subsubsection{Cardinality and reducibility in $\Ll ( \R )$.}
Because of the special nature of \( \Ll ( \R) \) (see the subsequent Remark~\ref{rmk:crucialL(R)}\ref{rmk:crucialL(R)-i}), many results about \markdef{\( \Ll ( \R ) \)-reducibility} can be recast in terms of \markdef{\( \Ll ( \R ) \)-cardinality} and conversely (see also Section~\ref{subsubsec:reducibility}, the introduction to Section~\ref{sec:definablecardinality}, and point~\ref{LR-4} below).

\begin{definition} \label{def:LRcardinality}
For \( A,B \in \Ll ( \R ) \) we let
\[
\begin{split}
\cardLR{A}_{\Ll ( \R )} \leq \cardLR{B}_{\Ll ( \R )} & \IFF \EXISTS{f \in \Ll ( \R )} ( f \colon A \into B )
\\
 & \IFF \Ll ( \R ) \models \card{A} \leq \card{B}.
\end{split}
\]
\end{definition}

The following are a sample of results in this context which are relevant for our discussion.
By Definition~\ref{def:LRcardinality}, comparing \( \Ll ( \R ) \)-cardinalities amounts to work inside \( \Ll ( \R ) \), so for ease of notation \emph{we will now step into this model, assume \( \AD \), and drop the subscripts, writing \( \card{A} \) instead of \( \card{A}_{\Ll ( \R )} \)}.

\begin{enumerate-(A)}
\item \label{LR-1}
Every subset of \( \R \) is either countable or else it contains a copy of \( \pre{\omega}{2} \) (so there is no intermediate cardinality between \( \omega \) and \( \card{\R} \)), and since \( \R \) is not well-orderable, then \( \card{\R} \) is incomparable with (the cardinality of) any uncountable ordinal.
This, in a sense, solves the continuum problem under \( \AD \).
\item \label{LR-2}
More generally, Woodin showed that \( \R \) is essentially the unique obstruction for a set \( A \) to be well-orderable: either \( \card{A} \leq \card{\alpha} \) for some ordinal \( \alpha \), or \( \card{\R } \leq \card{A} \) (see e.g.~\cite[Theorem 2.8 and Corollary 2.9]{Hjorth:1999bh}). 
Note that this generalizes Silver's dichotomy~\cite{Silver:1980fh} to arbitrary sets.
\item \label{LR-3}
Similarly, Hjorth~\cite[Theorem 2.6]{Hjorth:1995ve} extended the Glimm-Effros dichotomy \cite{Harrington:1990qa} to arbitrary sets \( A \) by showing that either \( \card{A} \leq \card{\pow ( \alpha)} \) for some ordinal \( \alpha \), or else \( \cardLR{\R / E_0 } \leq \cardLR{A} \), where \( E_0 \) is as in~\eqref{eq:E0}.
\item \label{LR-4}
Furthermore, \emph{every} set is of the form \( \bigcup_{\alpha < \kappa} A_\alpha \), where \( \kappa \in \Cn \) and each \( A_\alpha \) is of small cardinality, i.e.~for every \( \alpha < \kappa \) there is an equivalence relation \( E_ \alpha \) on \( \R \) such that \( \card{A_ \alpha } = \card{\R / E_ \alpha } \) (see~\cite[proof of Theorem 2.6]{Hjorth:1995ve} and~\cite[Lemma 2.13]{Hjorth:1999bh}).
Therefore all cardinalities in \( \Ll ( \R ) \) can be analyzed in terms of ordinals and reducibility between quasi-orders on \( \R \) (or, equivalently, on \( \pre{\omega}{2} \)).
\end{enumerate-(A)}

\begin{remarks} \label{rmk:crucialL(R)}
\begin{enumerate-(i)}
\item \label{rmk:crucialL(R)-i}
In~\ref{LR-2}--\ref{LR-4} above, it is crucial that we work in \( \Ll ( \R ) \), so that, in particular, every set is definable using only reals and ordinals as parameters. 
However, under \( \AD \) alone one still has that~\ref{LR-2}--\ref{LR-4} are true e.g.~for sets of small cardinality, i.e.~for sets \( A \) such that \( \card{A} = \card{\R / E} \) for some equivalence relation \( E \) on \( \R \) (without definability conditions on \( E \)), or even in the wider context of real-ordinal definable sets.
\item \label{rmk:crucialL(R)-ii}
Dichotomy~\ref{LR-2} shows that under \( \AD \) if \( A \) is arbitrary, \( X \) is a separable space, and \( \pre{ A}{X} \) is endowed with the product topology,
\[
 \pre{ A}{X} \text{ is separable} \IFF \card{A} \leq \omega \OR \card{A} = \card{\R} . 
\]
\end{enumerate-(i)}
\end{remarks}

Stepping-back into the universe of sets \( \Vv \) (where we may assume that \( \AC \) holds), all the results above can be restated as assertions about \( \Ll ( \R ) \)-cardinalities (as introduced in Definition~\ref{def:LRcardinality}) under the assumption \( \AD^{\Ll ( \R )} \).
For example~\ref{LR-2} reads as follows: assuming \( \AD^{\Ll ( \R )} \), for all \( A\in \Ll ( \R ) \) either \( \card{A}_{\Ll ( \R )} \leq \card{ \alpha }_{\Ll ( \R )} \) or else \( \card{\R}_{\Ll ( \R )} \leq \card{ A }_{\Ll ( \R )} \).

\subsubsection{Cardinality and reducibility in an inner model}

The approach above can be generalized by considering any inner model \( W \) of \( \ZF \) containing all the reals, so that in particular \( W \supseteq \Ll ( \R ) \).
If \( W \) is constructed in a canonical, explicit way, such as \( \Ll ( \R ) \) or \( \OD ( \R ) \), then it is reasonable to describe the objects in \( W \) as \emph{definable}.
Assuming e.g.~sufficiently large cardinal axioms or determinacy assumptions, or working in some special model of set theory (like the Solovay models), this approach gives in general a nice definable reducibility and cardinality theory. 
Therefore in what follows we will construe ``reasonably definable'' as ``belonging to the inner model \( W \supseteq \Ll ( \R ) \)'' under consideration. 
The resulting relations of \markdef{\( W \)-reducibility} \( \leq_W \) and \markdef{\( W \)-bi-reducibility} \( \sim_W \) are defined in the obvious way, namely:

\begin{definition} \label{def:Wreducibility}
Let \( W \supseteq \R \) be an inner model and \( R, S \in W \) be quasi-orders. 
Then 
\[
\begin{split}
R \leq_W S & \IFF \EXISTS{f \in W} (f \text{ reduces } R \text{ to } S) 
\\
& \IFF W \models R \leq S,
\end{split}
\]
and \( {R \sim_W S} \iff {{R \leq_W S} \wedge {S \leq_W R}} \).%
\footnote{The symbol \( \leq_W \), which is a reducibility for quasi-orders, should not be confused with \( \leqW^ \kappa \), the Wadge reducibility between subsets of \( \pre{ 2}{ \kappa } \) from Sections~\ref{sec:topologies} and~\ref{sec:borelsets}.}
\end{definition}

Using the fact that they hold in the \( \ZF \)-model \( W \), Theorems~\ref{th:graphs} and~\ref{th:main} can be reformulated as follows:

\begin{theorem} \label{th:graphsinnermodel}
Let \( W \supseteq \R \) be an inner model and suppose that \( \kappa \in \Cn^W \) and \( R \in \bS_W ( \kappa ) \).%
\footnote{See Remark~\ref{rmk:Souslinininnermodel} for the definition of \( \bS_W(\kappa) \).}
Let \( S \coloneqq ( \embeds^\kappa_{\CT})^W \) be the embeddability relation between combinatorial trees of size \( \kappa \) \emph{in \( W \)}, i.e.~\( S \) is the unique quasi-order (in \( \Vv \)) such that \( W \models \text{``}{S = {\embeds^\kappa_{\CT}}}\text{''} \). 
Then \( R \leq_W S \).
Therefore \( S \) is \emph{\( \leq_W \)-complete} for quasi-orders in \( \bS_W ( \kappa ) \). 
\end{theorem}

\begin{theorem} \label{th:maininnermodel}
Let \( W \supseteq \R \) be an inner model and let \( \kappa \in \Cn^W \). 
Then for every \( R \in \bS_W ( \kappa ) \) there is \( \upsigma \in ( \LL_{\kappa^+ \kappa})^W \subseteq \LL_{\kappa^+ \kappa} \) such that \( R \sim_W {\embeds ^\kappa_\upsigma} \).
\end{theorem}

\begin{remark}\label{rmk:maininnermodel}
As in~\cite[Chapter 9]{Hjorth:2000zr}, all terms in the statement of Theorem \ref{th:maininnermodel} must be relativized to \( W \). 
Thus e.g.~\( R \leq_W {\embeds ^\kappa_\upsigma} \) is construed as: there is \( f \colon \R \to ( \Mod^\kappa_\upsigma )^W \) such that for every \( x , y \in \R \)
\[
x \mathrel{R} y \IFF W \models f ( x ) \embeds f ( y ) . 
\] 
Unfortunately, the fact that we are forced to use \( ( \Mod^\kappa_\upsigma)^W \) instead of \( \Mod^\kappa_\upsigma \) (and that, in general, \( ( \Mod^\kappa_\upsigma )^W \subset \Mod^\kappa_\upsigma \)) forbids to formally restate Theorem~\ref{th:maininnermodel} in terms of \( \leq_W \)-invariant universality as introduced in Definitions~\ref{def:invuniversal} and~\ref{def:invuniversallocal}.
Notice also that in (the statement of) Theorem~\ref{th:maininnermodel} we have \( ( \Mod^\kappa_\upsigma)^W = \Mod^\kappa_\upsigma \cap W \) by Remark~\ref{rmk:elementaryclassrelativized}\ref{rmk:elementaryclassrelativized-i}: this is because in this case the formula \( \upsigma \) provided by Theorem~\ref{th:main} is obtained as in Corollary~\ref{cor:saturation}, that is it is of the form \( \upsigma_T \) for some \( T \in ( \TT_\kappa)^W \). 
Moreover, by Remark~\ref{rmk:elementaryclassrelativized}\ref{rmk:elementaryclassrelativized-ii} we also get that the \( \LL_{\kappa^+ \kappa} \)-sentence \( \upsigma \coloneqq \upsigma_T \) under consideration has the remarkable property that every model (in \( \Vv \)) of \( \upsigma \) has an isomorphic copy belonging to the inner model \( W \), so that even though they may fail to belong to \( M \), the \( \Vv \)-relations \( \cong^\kappa_\upsigma \) and \( \embeds^\kappa_\upsigma \) are at least faithfully represented in \( W \).
\end{remark}

The preorder \( \leq_W \) from Definition~\ref{def:Wreducibility} can be in principle extended to compare quasi-orders which are not necessarily in \( W \), namely for \emph{arbitrary} quasi-orders \( R , S \) of \( \Vv \) we can set
\[ 
R \leq'_W S \IFF \EXISTS{f \in W} (f \text{ reduces } R \text{ to } S).
 \] 
However, by definition \( R \leq'_W S \) can hold only if \( \dom(R) = \dom(f) \in W \) for some/any \( f \) witnessing \( R \leq'_W S \).
For this reason, we can restate in this more general context only the completeness result (Theorem~\ref{th:graphs}) and not the invariant universality one (Theorem~\ref{th:main}), as in the latter case \( \Mod^\kappa_{\upsigma_T} \) need not to belong to \( W \) --- as recalled in Remark~\ref{rmk:maininnermodel}, we are just guaranteed that \( \Mod^\kappa_{\upsigma_T} \cap W = ( \Mod^\kappa_{\upsigma_T})^W \in W \).

\begin{theorem} \label{th:W'}
Let \( W \supseteq \R \) be an inner model, \( \kappa \) be a cardinal, and \( R \in \bS_W ( \kappa ) \). 
Then \( R \leq'_W {\embeds^\kappa_{\CT}} \).
Therefore \( \embeds^\kappa_{\CT} \) is \emph{\( \leq'_W \)-complete} for quasi-orders in \( \bS_W ( \kappa ) \).
\end{theorem}

\begin{proof}
Let \( T \in ( \TT_\kappa)^W \) be such that \( R = \PROJ \body{ T } \). 
Since \( \kappa \in \Cn^W \) (because \( W \) is an inner model) and \( T \in W \), by (the proof of) Theorem~\ref{th:absolutedefin} the function \( f_T \) from~\eqref{eq:f_T} is absolutely definable (using only \( T \in W \) as a parameter) and reduces \( R \) to \( \embeds^\kappa_{\CT} \). 
Since \( \R \subseteq W \), by~\eqref{eq:coherence} we get \( f_T^W = f_T \), whence \( f_T \in W \). 
Therefore, \( f_T \) witnesses \( R \leq'_W {\embeds^\kappa_{\CT}} \), as required.
\end{proof}

\section{Some applications} \label{sec:applications}
As we will see in this section, the results obtained in Sections~\ref{sec:embeddabilitygraphs}--\ref{sec:definablecardinality} yield several corollaries: inside any model where the \( \kappa \)-Souslin quasi-orders on \( \pre{\omega}{2} \) (or, equivalently, on any Polish or standard Borel space) form an interesting class, we get a \emph{completeness} and \emph{invariant universality} result for the embeddability relation on combinatorial trees of size \( \kappa \). 
Since there are lot of situations of this kind, in what follows we will just explicitly state a few of them which correspond to some of the most relevant cases. 
All these results admit several variants which will not be explicitly mentioned but that could be of interest on their own, namely:
\begin{enumerate-(1)}
\item
the relations \( \leq^\kappa_{\bB} \) and \( \sim^\kappa_{\bB} \) could systematically be replaced by their effective counterparts \( \leq^\kappa_{B} \) and \( \sim^\kappa_{B} \);
\item
in each statement, we could equivalently consider the wider class of \( \kappa \)-Souslin quasi-orders defined on \emph{arbitrary Polish spaces} (or even on \emph{arbitrary standard Borel spaces}) instead of its restriction to the \( \kappa \)-Souslin quasi-orders on \( \pre{\omega}{2} \) --- this is because every two uncountable Polish or standard Borel spaces are Borel isomorphic and \( \bS ( \kappa) \) is closed under Borel (pre)images by Lemma~\ref{lem:Souslin};
\item
elementary classes of the form \( \Mod^\kappa_\upsigma \) for some \( \LL_{\kappa^+ \kappa} \)-sentence \( \upsigma \) obtained as in Corollary~\ref{cor:saturation} could always be equivalently substituted by \( \Mod^\infty_\upsigma \) by Remark~\ref{rem:maintheorem}\ref{rem:maintheorem-1};
\item
by Corollary~\ref{cor:maintheorem}, in all the subsequent results concerning the bi-reducibility between two quasi-orders \( R \) and \( S \) we could further add that the quotient orders of \( R \) and \( S \) are in fact \emph{isomorphic}.
\end{enumerate-(1)}

\subsection{$\bSigma^1_2$ quasi-orders}

Recall that since the \( \leq_{\bSigma^1_1} \)-reducibility (see Definition~\ref{def:S(kappa)-reducibility}), which is the same as the \( \leq_{ \bS ( \omega)} \)-reducibility, coincides with the classical Borel reducibility \( \leq_\bB^\omega \), by Theorems~\ref{th:LouveauRosendal} and~\ref{th:mottorosfriedman} we have that:
\begin{description}[leftmargin=1.5pc]
\item[Completeness]
\( \embeds_\CT^\omega \) is \emph{\( \leq_{\bSigma^1_1} \)-complete} (equivalently, \emph{\( \leq_{ \bS ( \omega)} \)-complete}) for \( \bSigma^1_1 \) quasi-orders on \( \pre{\omega}{2} \), i.e.\ \( R \leq_{\bSigma^1_1 } {\embeds_\CT^\omega} \) (equivalently, \( R \leq_{ \bS ( \omega)} {\embeds_\CT^\omega} \)) for every \( \bSigma^1_1 \) quasi-order \( R \) on \( \pre{\omega}{2} \);
\item[Invariant universality]
\( \embeds_\CT^\omega \) is also \emph{\( \leq_\bB^\omega \)-invariantly universal} (and hence also \( \leq_{\bB}^\omega \)-complete) for \( \bSigma^1_1 \) quasi-orders on \( \pre{\omega}{2} \), i.e.\ for every such \( R \) there is an \( \LL_{\omega_1 \omega} \)-sentence \( \upsigma \) such that \( R \sim_\bB^{\omega} {\embeds^\omega_\upsigma} \).
\end{description}
The following theorems generalize the above results to the next level of the projective hierarchy. 
(It is easy to check that in all the situations below it always makes sense to consider \( \bSigma^1_2 \)-in the-codes and \( \bS ( \omega_1) \)-in-the-codes functions \( f \colon \pre{\omega}{2} \to \Mod^{\omega_1}_{\LL} \) by the observations following Proposition~\ref{prop:Somega1=sigma12} and Lemma~\ref{lem:regularSouslinareunbounded}.)

\begin{theorem}\label{th:finalSigma12-a}
Assume either \( \AC \) or \( \AD + \DC \). 
Then the relation \( \embeds_\CT^{\omega_1} \) is \( \leq_{ \bS ( \omega_1)} \)-complete%
\footnote{By Proposition~\ref{prop:Skappainthecodesareborel}, we could replace \( \leq_{ \bS ( \omega_1)} \)-completeness with \( \leq^{\omega_1}_\bB \)-completeness. 
However the latter is a weaker notion, and in fact the corresponding completeness result can be proved already in \( \ZF + \AC_{\omega}(\R) \) --- see Theorem~\ref{th:finalSigma12-b}.}
for \( \bSigma^1_2 \) quasi-orders on \( \pre{\omega}{2} \).
\end{theorem}

\begin{proof}
Apply Theorem~\ref{th:Gammacompleteness}, using the fact that \( \bSigma^1_2 \subseteq \bS ( \omega_1) \). 
\end{proof}

Notice that in Theorem~\ref{th:finalSigma12-a} when assuming \( \AD + \DC \) we could replace the somewhat artificial notion of \( \leq_{ \bS ( \omega_1)} \)-completeness with the more natural notion of \( \leq_{\bSigma^1_2} \)-completeness because in this case \( \bS ( \omega_1) = \bSigma^1_2 \) (see the paragraph before Proposition~\ref{prop:deltaSkappa}). 
A similar strengthening can be obtained in models of \( \AC \) as well under some additional set-theoretic assumptions. 

\begin{theorem}[\( \AC \)] \label{th:sigma12improved}
Assume either \( \AD^{\Ll ( \R )} \) or \( \MA + \neg \CH + \EXISTS{a \in \pre{\omega}{\omega}} (\omega_1^{\Ll [ a ] } = \omega_1) \). 
Then \( \embeds_\CT^{\omega_1} \) is \( \leq_{\bSigma^1_2} \)-complete for \( \bSigma^1_2 \) quasi-orders on \( \pre{\omega}{2} \).
\end{theorem}

\begin{proof}
Under \( \ZFC + \MA + \neg \CH + \EXISTS{a \in \pre{\omega}{\omega}} (\omega_1^{\Ll[a]} = \omega_1) \) we get \( \bS ( \omega_1) = \bSigma^1_2 \) by Proposition~\ref{prop:Somega1=sigma12}, so the result follows from Theorem~\ref{th:finalSigma12-a}.

Let us now assume \( \AD^{ \Ll ( \R ) } \), and recall that \( \omega_1^{ \Ll ( \R ) } = \omega_1 \) and \( ( \pre{\omega}{2} )^{ \Ll ( \R ) } = \pre{\omega}{2} \). 
Let \( R \) be a \( \bSigma^1_2 \) quasi-order on \( \pre{\omega}{2} \), so that \( R \) is \( \bSigma^1_2 \) also in \( \Ll ( \R ) \) and hence \( R \in \bS_{ \Ll ( \R ) }(\omega_1) \). 
Let \( T \in (\TT_{\omega_1})^{ \Ll ( \R ) } \) be such that \( R = \PROJ \body{T} \).
Then by (the proof of) Theorem~\ref{th:W'} the function \( f_T \) from~\eqref{eq:f_T} is a reduction of \( R \) to \( \embeds^{\omega_1}_\CT \) and belongs to \( \Ll ( \R ) \). 
Under our assumptions, \( \ZF + \AD + \DC \) holds in \( \Ll ( \R ) \) (whence also \( \Ll ( \R ) \models \text{``} \bS ( \omega_1) = \bSigma^1_2 \text{''}\)), so applying Theorem~\ref{th:Gammacompleteness}\ref{th:Gammacompleteness-b} in \( \Ll ( \R ) \) we have that 
\[ 
 \Ll ( \R ) \models f_T \text{ is \( \bSigma^1_2 \)-in-the-codes}.
\]
Thus \( f_T \) is actually \( \bSigma^1_2 \)-in-the-codes (in \( \Vv \)) by Shoenfield's \( \bSigma^1_2 \)-absoluteness.
\end{proof}

 As for invariant universality, we get the following result in \( \ZF + \AC_\omega ( \R ) \).

\begin{theorem}[\( \AC_\omega(\R) \)] \label{th:finalSigma12-b}
The relation \( \embeds_\CT^{\omega_1} \) is \( \leq_{\bB}^{\omega_1} \)-invariantly universal (and hence also \( \leq_{\bB}^{\omega_1} \)-complete) for \( \bSigma^1_2 \) quasi-orders on \( \pre{\omega}{2} \).
\end{theorem}

\begin{proof}
Since \( \bSigma^1_2 \subseteq \bS ( \omega_1) \), it is enough to use Theorem~\ref{th:maintopology}, which can be applied because \( \AC_\omega(\R) \) implies that \( \omega_1 \) is a regular cardinal.
\end{proof}

Theorems~\ref{th:finalSigma12-a}--\ref{th:finalSigma12-b} are more interesting in all cases in which the topological notions involved (such as \( \bS ( \omega_1) \)-in-the-codes functions, \( \omega_1 + 1 \)-Borel sets and functions, and so on) are nontrivial, i.e.\ when \( \omega_1 \) is ``small enough'' with respect to the cardinality of the continuum. 
This is the case if we work in models of \( \AD \). 
If instead we work in models of \( \ZFC \), then as observed in Sections~\ref{sec:topologies}--\ref{sec:otherspacesandBairecategory} and~\ref{sec:Ksouslinsets} all the relevant topological notions trivialize under \( \CH \). 
Thus we are naturally lead to work in models of \( \ZFC + \neg \CH \). 
By the observation following Proposition~\ref{prop:inthecodesSouslinAC}, this already gives that the notion of a \( \bS ( \omega_1) \)-in-the-code function \( f \colon \pre{\omega}{2} \to \Mod^{\omega_1}_{\LL} \) is nontrivial. 
Moreover, if we further assume that the inequality \( 2^{\aleph_1} < 2^{(2^{\aleph_0})} \) holds (which is e.g.\ the case in models of forcing axioms like \( \MA_{\omega_1} \), \( \PFA \), and so on), then also the notions of a (weakly) \( \omega_1 + 1 \)-Borel function \( f \colon \pre{\omega}{2} \to \Mod^{\omega_1}_\LL \) becomes interesting (see Corollary~\ref{cor:prop:examples+weaklyBorel}\ref{cor:prop:examples+weaklyBorel-b}). 
This discussion shows that the following instantiation of Theorems~\ref{th:finalSigma12-a} and~\ref{th:sigma12improved} is nontrivial.

\begin{corollary}[\( \AC+\PFA \)] \label{cor:omega1underPFA}
The relation \( \embeds_{\CT}^{\omega_1} \) is \( \leq_{\bSigma^1_2} \)-complete (and hence also \( \leq_{ \bS ( \omega_1)} \)- and \( \leq^{\omega_1}_\bB \)-complete) for \( \bSigma^1_2 \) quasi-orders on \( \pre{\omega}{2} \).
\end{corollary}

\begin{proof}
\( \AD^{ \Ll ( \R ) } \) follows from \( \PFA \) by~\cite{Steel:2005ab}, so Theorem~\ref{th:sigma12improved} can be applied.
\end{proof}

If besides \( 2^{\aleph_1} < 2^{(2^{\aleph_0})} \) we further assume%
\footnote{Since \( 2^{\aleph_0} = 2^{\aleph_1} \) implies \( 2^{(2^{< \aleph_1})} = 2^{(2^{\aleph_1})} \), \( \MA_{\omega_1} \) (and even the stronger \( \PFA \), \( \MM \), and so on) are not sufficient to ensure the extra cardinal condition under discussion. 
However, it could still be the case that the notion of an \( \omega_1 +1 \)-Borel function \( f \colon \Mod^{\omega_1}_{\LL} \to \pre{\omega}{2} \) is nontrivial also in models of forcing axioms for reasons different from the cardinality considerations of Corollary~\ref{cor:prop:examples+weaklyBorel}.} 
that \( 2^{(2^{< \aleph_1})} < 2^{(2^{\aleph_1})} \), then also the notion of an \( \omega_1 +1 \)-Borel function \( f \colon \Mod^{\omega_1}_{\LL} \to \pre{\omega}{2} \) becomes nontrivial (see Corollary~\ref{cor:prop:examples+weaklyBorel}\ref{cor:prop:examples+weaklyBorel-c}). 
Therefore, a situation of interest in which Theorem~\ref{th:finalSigma12-b} can be applied is \( 2^\kappa = \kappa^{++} \) for every \( \kappa \leq \aleph_3 \).

\subsection{Projective quasi-orders}
In this section we will generalize Theorems~\ref{th:LouveauRosendal} and~\ref{th:mottorosfriedman} to even larger projective levels (under suitable assumptions).

\subsubsection{Models of \( \AC \)}

Martin showed in~\cite{Martin:2012ma} that \( \ZFC + \FORALL{ x \in \pre{\omega}{\omega}} (x^\# \text{ exists}) \) implies that all \( \bSigma^1_3 \) sets are \( \omega_2 \)-Souslin. 
Using this fact and applying Theorems~\ref{th:Gammacompleteness}\ref{th:Gammacompleteness-a} and~\ref{th:maintopology}, we get the following further generalizations of Theorems~\ref{th:LouveauRosendal} and~\ref{th:mottorosfriedman}. 
(Recall that in models with choice it always makes sense to speak of \( \bS ( \omega_2) \)-in-the-codes functions \( f \colon \pre{\omega}{2} \to \Mod^{\omega_2}_{\LL} \) by the observation following Proposition~\ref{prop:Somega1=sigma12}, and that every such \( f \) is automatically weakly \( \omega_2 +1 \)-Borel by Proposition~\ref{prop:Skappainthecodesareborel}.)

\begin{theorem}[\( \AC \)] \label{th:mainSigma13AC}
Assume that \( x^\# \) exists for all \( x \in \pre{\omega}{\omega} \).
\begin{enumerate-(a)}
\item \label{th:mainSigma13AC-a}
The relation \( \embeds_\CT^{\omega_2} \) is \( \leq_{ \bS ( \omega_2)} \)-complete for \( \bSigma^1_3 \) quasi-orders on \( \pre{\omega}{2} \). 
\item \label{th:mainSigma13AC-b}
The relation \( \embeds_\CT^{\omega_2} \) is \( \leq_{\bB}^{\omega_2} \)-invariantly universal (and hence also \( \leq_{\bB}^{\omega_2} \)-complete) for \( \bSigma^1_3 \) quasi-orders on \( \pre{\omega}{2} \).
\end{enumerate-(a)}
\end{theorem}

An interesting application of this theorem is when considering the quasi-order \( ( \mathcal{Q} , \leq_{\bB} ) \) of Example~\ref{xmp:analyticquasiorders}, that is the relation of Borel reducibility between analytic quasi-orders. 
As recalled in the introduction, such quasi-order may be seen as a (definable) embeddability relation between structures of size the continuum \( 2^{\aleph_0} \): the next result shows that \( ( \mathcal{Q} , \leq_{\bB} ) \) can be turned into an embeddability relation on structures of size \( \aleph_2 \), \emph{independently of the actual value of \( 2^{\aleph_0} \)}. 

\begin{theorem}[\( \AC \)] \label{th:Qembedsundersharps}
Assume that \( x^\# \) exists for all \( x \in \pre{\omega}{\omega} \).
Then the quotient order of \( ( \mathcal{Q} , \leq_{\bB} ) \) (definably) embeds into the quotient order of \( \embeds _\CT^{\aleph_ 2 } \).

Moreover there is an \( \LL_{ \aleph _3 \, \aleph_2 } \)-sentence \( \upsigma \) such that the quotient orders of \( ( \mathcal{Q} , \leq_{\bB} ) \) and \( \embeds_ \upsigma^{\aleph_2 } \) are even (definably) isomorphic.
\end{theorem}

\begin{proof}
For the additional part about the \( \LL_{\aleph_3 \, \aleph_2} \)-sentence \( \upsigma \), we use Corollary~\ref{cor:maintheorem}.
\end{proof}

Recall that the assumption \( \FORALL{ x \in \pre{\omega}{\omega}} (x^\# \text{ exists}) \) is equivalent over \( \ZFC \) to \( \bSigma^1_1 \)-determinacy (in fact, even to \( {<} \omega^2\text{-} \bPi^1_1 \)-determinacy) by results of Harrington and Martin (see~\cite{Harrington:1978ht}). 
Assuming more determinacy, we can extend Theorem~\ref{th:mainSigma13AC} to all projective levels (recall that the axiom \( \AD^{\Ll ( \R )} \) used in Theorem~\ref{th:embeddingsofprojectiveqos} follows both from the existence of infinitely many Woodin cardinals with a measurable above, and from strong forcing axioms such as \( \PFA \)).

\begin{theorem}[\( \AC \)]\label{th:embeddingsofprojectiveqos}
Assume \( \AD^{\Ll ( \R )} \). 
Then there is a monotone function \( r \colon \omega \to \omega \) such that: 
\begin{enumerate-(a)}
\item \label{th:embeddingsofprojectiveqos-a}
the relation \( \embeds_\CT^{\omega_{r(n)}} \) is \( \leq_{ \bS ( \omega_{r(n)})} \)-complete for \( \bSigma^1_n \) quasi-orders on \( \pre{\omega}{2} \);
\item \label{th:embeddingsofprojectiveqos-b}
the relation \( \embeds_\CT^{\omega_{r(n)}} \) is \( \leq_{\bB}^{\omega_{r(n)}} \)-invariantly universal (and hence also \( \leq_{\bB}^{\omega_{r(n)}} \)-complete) for \( \bSigma^1_n \) quasi-orders on \( \pre{\omega}{2} \).
\end{enumerate-(a)}
An upper bound for \( r \) is given by:
\[
 r ( n ) \leq \begin{cases}
2 ^{k + 1} - 2 & \text{ if } n = 2 k + 1 ,
\\
2^{k + 1} - 1 & \text{ if } n = 2 k + 2 .
 \end{cases}
 \]
\end{theorem}

The proof below will in particular show that, under the hypotheses of the theorem, there is an \( \bS(\omega_{r(n)}) \)-code for \( \omega_{r(n)} \) for each \( n \), so that in part~\ref{th:embeddingsofprojectiveqos-a} it makes sense to speak of \( \bS ( \omega_{r(n)}) \)-in-the-codes functions.

\begin{proof}
As every element of \( \Ll ( \R ) \) is definable reals and ordinals, and since any \( \omega \)-sequence of reals belongs to \( \Ll ( \R ) \), then \( \AC \) implies \( \DC^{\Ll ( \R )} \).
Therefore under our assumptions \( \Ll ( \R ) \models \AD + \DC \). 
Recall that projective sets are absolute between \( \Vv \) and \( \Ll ( \R ) \), in particular \( \bdelta^1_n = (\bdelta^1_n)^{ \Ll ( \R ) } \), but while \( \bdelta^1_n \) is always a cardinal in \( \Ll ( \R ) \) by \( \Ll ( \R ) \models \AD + \DC \), it may be just an ordinal in \( \Vv \)). 
Let \( \kappa_n \coloneqq \bdelta^1_{n-1} = (\bdelta^1_{n-1})^{ \Ll ( \R ) } \) if \( n \) is even and \( \kappa_n \coloneqq (\blambda^1_n)^{ \Ll ( \R ) } \) if \( n \) is odd, and let \( r(n) \) be such that \( \omega_{r(n)} = \card{\kappa_n} \) (where the cardinality of \( \kappa_n \) is computed in \( \Vv \)). 
By Corollary~\ref{cor:delta1n<aleph} we get the above upper bound for the function \( r \colon \omega \to \omega \) (in the odd case we use the fact that \( (\blambda^1_{n})^{ \Ll ( \R ) } \) has countable cofinality in \( \Ll ( \R ) \), and hence its \( \Vv \)-cardinality is collapsed at least to the \( \Vv \)-cardinality of the largest regular cardinal of \( \Ll ( \R ) \) below it). 
Since \( \Ll ( \R ) \models \AD + \DC \), by~\cite[Theorem 2.18]{Jackson:2010ff} (see the observation before Proposition~\ref{prop:deltaSkappa}) we get
\[ 
\bSigma^1_n = ( \bS ( \kappa_n ) )^{ \Ll ( \R ) } \subseteq \bS_{ \Ll ( \R ) }( \kappa_n ) \subseteq \bS ( \omega_{r ( n ) }).
\]
Notice that since \( ( \bS ( \kappa_n ) )^{ \Ll ( \R ) } \subseteq \bS ( \omega_{ r ( n ) }) \) and \( \kappa_n \geq \omega_{ r ( n ) } \), it follows from Proposition~\ref{prop:deltaSkappa}\ref{prop:deltaSkappa-b} (applied in the model \( \Ll(\R) \) with \( \kappa \coloneqq \kappa_n \)) and Remark~\ref{rmk:inthecodes}\ref{rmk:inthecodes-iii} that it always makes sense to speak of \( \bS ( \omega_{r(n)}) \)-in-the-codes functions \( f \colon \pre{\omega}{2} \to \Mod^{\omega_{ r ( n ) }}_{\LL} \) and that, in particular, the hypotheses of Theorems~\ref{th:Gammacompleteness}\ref{th:Gammacompleteness-a} are satisfied for \( \kappa \coloneqq \omega_{r(n)} \), although it may happen that \( \omega_{ r ( n ) } = 2^{\aleph_0} \) if the continuum is smaller than \( \aleph_\omega \). 
Moreover, the \( \omega_{r ( n ) } \) are always successor cardinals, and thus regular in model of \( \AC \), so that the hypotheses of Theorem~\ref{th:maintopology} are satisfied for \( \kappa \coloneqq \omega_{r(n)} \).
Thus, to obtain the desired results it is enough to apply Theorems~\ref{th:Gammacompleteness}\ref{th:Gammacompleteness-a} and~\ref{th:maintopology}.
\end{proof} 

\subsubsection{Models of \( \AD \)} \label{subsubsec:modelsofAD}

By~\cite[Theorem 2.18]{Jackson:2010ff}, assuming \( \ZF + \AD + \DC \) we get that \( \bSigma^1_n = \bS ( \kappa_n) \), where \( \kappa_n \) is such that \( \bdelta^1_n = \kappa_n^+ \), that is \( \kappa_n \coloneqq \blambda^1_n \) if \( n \) is odd and \( \kappa_n \coloneqq \bdelta^1_{n-1} \) otherwise. 
Therefore, applying Theorems~\ref{th:Gammacompleteness}\ref{th:Gammacompleteness-b} and~\ref{th:maintopology} we get the following \( \AD \)-analogue of Theorem~\ref{th:embeddingsofprojectiveqos}.

\begin{theorem}[\( \AD + \DC \)] \label{th:mainprojective}
Let \( 0 \neq n \in \omega \). 
\begin{enumerate-(a)}
\item \label{th:mainprojective-a}
The relation \( \embeds_\CT^{\kappa_n} \) is \( \leq_{\bSigma^1_n} \)-complete for \( \bSigma^1_n \) quasi-orders on \( \pre{\omega}{2} \).
\item \label{th:mainprojective-b}
Let \( n \) be an \emph{even} number. 
Then the relation \( \embeds_\CT^{\bdelta^1_{n-1}} \) is \( \leq^{\bdelta^1_{n-1}}_{\bB} \)-invariantly universal (and hence also \( \leq^{\bdelta^1_{n-1}}_{\bB} \)-complete) for \( \bSigma^1_n \) quasi-orders on \( \pre{\omega}{2} \).
\end{enumerate-(a)}
\end{theorem}

\begin{remarks} \label{rmk:proj}
\begin{enumerate-(i)}
\item \label{rmk:proj-i}
In part~\ref{th:mainprojective-a} we can speak of \( \leq_{\bSigma^1_n} \)-completeness because in the \( \AD \)-world \( \bS ( \kappa_n) = \bSigma^1_n \).
Notice that \( \leq_{\bSigma^1_n} \)-completeness implies \( \leq^{\kappa_n}_\bB \)-completeness by Proposition~\ref{prop:Skappainthecodesareborel}, and that the two notions coincide if \( n \) is odd by Corollary~\ref{cor:inthecodesSouslin}.
\item \label{rmk:proj-ii}
The restriction in part~\ref{th:mainprojective-b} comes from the fact that \( \kappa_n \) is regular if and only if \( n \) is even, and therefore Theorem~\ref{th:maintopology} can be applied only to the even levels of the projective hierarchy. 
However, if we drop the requirement that the reductions be (weakly) \( \kappa_n + 1 \)-Borel, then Theorem~\ref{th:mainprojective}\ref{th:mainprojective-b} would be true also for the odd levels by Theorem~\ref{th:main}.
\end{enumerate-(i)}
\end{remarks}

As for Theorem~\ref{th:mainSigma13AC}, also Theorem~\ref{th:mainprojective} can be applied with \( n = 3 \) to the quasi-order \( ( \mathcal{Q} , \leq_{\bB} ) \) of Borel reducibility between analytic quasi-orders: in this case, such definable embeddability is turned into the embeddability relation between structures of size \( \kappa_3 = \aleph_\omega \). 
Since \( \aleph_\omega \) is a singular cardinal, to get the second part of the following theorem use the observation in Remark~\ref{rmk:proj}\ref{rmk:proj-ii} together with Corollary~\ref{cor:maintheorem}, rather than Theorem~\ref{th:mainprojective}\ref{th:mainprojective-b}.

\begin{theorem}[\( \AD + \DC \)]\label{th:QembedsintoCTalephomega}
The quotient order of \( ( \mathcal{Q} , \leq_{\bB} ) \) (definably) embeds into the quotient order of \( \embeds _\CT^{\aleph_ \omega } \).

Moreover there is an \( \LL_{ \aleph _{\omega+1} \, \aleph_\omega } \)-sentence \( \upsigma \) such that the quotient orders of \( ( \mathcal{Q} , \leq_{\bB} ) \) and \( \embeds_ \upsigma^{\aleph_\omega } \) are even (definably) isomorphic.
\end{theorem}

\subsection{More complex quasi-orders in models of determinacy} \label{subsec:modelsofAD+}
Assuming \( \ZF + \AD + \DC \), it is natural to consider the largest boldface pointclass to which our results can be applied, namely the collection \( \bS ( \Xi ) = \bS ( \infty ) \) of all \( \infty \)-Souslin sets, where \( \Xi \) is the supremum of all Suslin cardinals (Definition~\ref{def:Souslincardinal}). 
Assuming \( \AD \) in \( \Ll ( \R ) \), we have that \( \Xi = \bdelta^2_1 \) and \( \bS ( \Xi ) = \bSigma^2_1 \). 
In this case \( \bdelta^2_1 \) is a Souslin cardinal, so applying again Theorems~\ref{th:Gammacompleteness} and~\ref{th:maintopology} we get the following completeness and invariant universality results.

\begin{theorem}[\( \AD + \Vv = \Ll ( \R ) \)] \label{th:mainLR}
\begin{enumerate-(a)}
\item \label{th:mainLR-a}
The embeddability relation \( \embeds_\CT^{\bdelta^2_1} \) is \( \leq_{\bSigma^2_1} \)-complete for \( \bSigma^2_1 \) quasi-orders (equivalently, \( \infty \)-Souslin quasi-orders) on \( \pre{\omega}{2} \).
\item \label{th:mainLR-b}
The relation \( \embeds_\CT^{\bdelta^2_1}\) is \( \leq^{\bdelta^2_1}_{\bB} \)-invariantly universal (and hence also \( \leq^{\bdelta^2_1}_{\bB} \)-complete) for \( \bSigma^2_1 \) quasi-orders (equivalently, \( \infty \)-Souslin quasi-orders) on \( \pre{\omega}{2} \).
\end{enumerate-(a)}
\end{theorem}

\begin{remarks}
\begin{enumerate-(i)}
\item
In Theorem~\ref{th:mainLR} we need not explicitly assume \( \DC \) as this follows from \( \AD \) in \( \Ll ( \R ) \) by~\cite{Kechris:1984il}.
\item
Theorem~\ref{th:maintopology} can be applied to get part~\ref{th:mainLR-b} as \( \bSigma^2_1 \) is also closed under coprojections and therefore the cardinal \( \bdelta^2_1 \) is regular by~\cite[Theorem 7D.8]{Moschovakis:2009fk} (in fact, it is a weakly inaccessible cardinal). 
\end{enumerate-(i)}
\end{remarks}

Theorem~\ref{th:mainLR} can be analogously reformulated in every model of \( \ZF + \AD + \DC \) in which \( \Xi \) is a Souslin cardinal, i.e.~in every model of \( \ZF + \AD^+ + \DC \) (see Section~\ref{sec:SouslinunderAD}): in fact, in this case \( \Xi \) is always a (limit) regular cardinal (see e.g.\ \cite[Lemma 2.20]{Jackson:2008pi}).

\begin{theorem}[\( \AD^+ + \DC \)] \label{th:levelbylevel}
\begin{enumerate-(a)}
\item 
The embeddability relation \( \embeds_\CT^{\Xi} \) is \( \leq_{ \bS ( \Xi)} \)-complete
 for \( \infty \)-Souslin quasi-orders on \( \pre{\omega}{2} \).
\item 
The relation \( \embeds_\CT^{\Xi}\) is \( \leq^{\Xi}_{\bB} \)-invariantly universal (and hence also \( \leq^{\Xi}_{\bB} \)-complete) for \( \infty \)-Souslin quasi-orders on \( \pre{\omega}{2} \).
\end{enumerate-(a)}
\end{theorem}

\begin{remark} \label{rmk:smallcardinalitiesanticipated}
Applying Theorem~\ref{th:Gammacompleteness} and~\ref{th:maintopology}, a level-by-level version of Theorem~\ref{th:levelbylevel} is obtained. 
Assuming \( \AD + \DC \) it is possible to assign to every quasi-order \( R \) in \( \bDelta_{ \bS ( \Xi)} \) a regular Souslin cardinal \( \kappa_R < \Xi \) such that:
\begin{enumerate-(a)}
\item 
\( R \leq_{ \bS ( \kappa_R)} {\embeds_\CT^{\kappa_R}} \);
\item 
 \( R \sim^{\kappa_R}_\bB {\embeds^{\kappa_R}_{\upsigma}} \) for some \( \LL_{\kappa^+_R \, \kappa_R} \)-sentence \( \upsigma \) all of 
whose models are combinatorial trees.
\end{enumerate-(a)}
To see this, it is enough to let \( \kappa_R < \Xi \) be the smallest \emph{regular} Souslin cardinal such that \( R \in \bS ( \kappa_R ) \), whose existence is granted by Lemma~\ref{lem:regularSouslinareunbounded}.
\end{remark}

Further assuming \( \AD_\R \), we get global completeness and invariant universality results. 
Notice that in this case it does not make much sense to consider \( \leq_{\bS ( \Theta ) } \)-reducibility, since under \( \ZF + \ADR + \DC \) every subset of a Polish space is in \( \bS ( \infty ) = \bS ( \Theta) \) by Proposition~\ref{prop:ADR<=>Xi=Theta}, and hence \( \leq_{\bS ( \Theta ) } \) would coincide with the reducibility \( \leq \) (without any definability condition on the reductions).

\begin{theorem}[\( \AD_\R + \DC \)] \label{th:pag:smallcardinalities}
\begin{enumerate-(a)}
\item
The embeddability relation \( \embeds_\CT^ {\Theta} \) is \( \leq^{\Theta}_{\bB} \)-complete for \emph{arbitrary} quasi-orders on \( \pre{\omega}{2} \).
\item 
Assume that \( \Theta \) is regular. 
Then the relation \( \embeds_\CT^ {\Theta} \) is also \( \leq^{\Theta}_{\bB} \)-invariantly universal for \emph{arbitrary} quasi-orders on \( \pre{\omega}{2} \).
\end{enumerate-(a)}
\end{theorem}

\begin{proof}
By Proposition~\ref{prop:ADR<=>Xi=Theta}, under \( \AD_\R + \DC \) every subset of \( \pre{\omega}{2} \) is \( \infty \)-Souslin and \( \Xi = \Theta \). 
So it is enough to apply Theorems~\ref{th:definable=borel} and~\ref{th:maintopology}.
\end{proof}

\begin{remark} \label{rmk:smallcardinalities}
As for Theorem~\ref{th:levelbylevel}, also in this case one can formulate and prove level-by-level versions of Theorem~\ref{th:pag:smallcardinalities}, namely: under \( \AD_\R + \DC \) (which implies \( \Xi = \Theta \) by Proposition~\ref{prop:ADR<=>Xi=Theta}), one can assign to \emph{every} quasi-order \( R \) on \( \pre{\omega}{2} \) a regular Souslin cardinal \( \kappa_R < \Theta \) such that:
\begin{enumerate-(a)}
\item 
\( R \leq_{ \bS ( \kappa_R)} {\embeds_\CT^{\kappa_R}} \);
\item 
 \( R \sim^{\kappa_R}_\bB {\embeds^{\kappa_R}_{\upsigma}} \) for some \( \LL_{\kappa^+_R \, \kappa_R} \)-sentence \( \upsigma \) all of 
whose models are combinatorial trees.
\end{enumerate-(a)}
(Notice that in this case we necessarily have \( \kappa_R < \Theta \) by Lemma~\ref{lem:Souslincardinals<Theta}.)
\end{remark}

\subsection{\( \Ll ( \R ) \)-reducibility}

In this section we present some results concerning the \( \Ll ( \R ) \)-reducibility considered in~\cite{Hjorth:1995ve,Hjorth:1999bh, Hjorth:2000zr} --- see Example~\ref{xmp:AD} and Section~\ref{subsubsec:reducibilityininnermodel}. 
We recall once more that the axiom \( \AD^{ \Ll ( \R ) } \) used in Theorem~\ref{th:largecardinals1} follows from both large cardinals and forcing axioms. 
Notice also that the pointclass \( \bGamma^2_1 \coloneqq ( \bSigma^2_1 )^{\Ll ( \R )} \) is closed (in \( \Vv \)) under preimages and images of continuous functions because all such functions belong to \( \Ll ( \R ) \).

\begin{theorem}[\( \AC \)] \label{th:largecardinals1}
Assume \( \AD^{ \Ll ( \R ) } \).
\begin{enumerate-(a)}
\item \label{th:largecardinals1-a}
The relation \( \embeds^{\bdelta^2_1}_\CT \) is \( \leq_{ \Ll ( \R ) } \)-complete (and also \( \leq'_{ \Ll ( \R ) } \)-complete) for quasi-orders on \( \pre{\omega}{2} \) belonging to \( \bGamma^2_1 \coloneqq ( \bSigma^2_1 )^{\Ll ( \R )} \).
\item \label{th:largecardinals1-b}
For every \( \bGamma^2_1 \) quasi-order \( R \) there is an \( \LL_{( \bdelta^2_{1} )^+ \, \bdelta^2_{1}} \)-sentence \( \upsigma \) (in \( \Ll ( \R ) \)) such that \( R \sim_{\Ll ( \R )} {\embeds ^{\bdelta^2_1}_\upsigma} \).
\end{enumerate-(a)} 
\end{theorem}

\begin{proof}
Use Theorems~\ref{th:graphsinnermodel},~\ref{th:W'} and~\ref{th:maininnermodel}.
\end{proof}

\begin{remarks}
\begin{enumerate-(i)}
\item
When dealing with \( \Ll ( \R ) \)-reducibility, the relations \( \embeds_\CT^\kappa \) and \( \embeds ^\kappa_\upsigma \) (for \( \kappa \in ( \Cn )^{ \Ll ( \R ) } \) and \( \upsigma \in (\LL_{\kappa^+ \kappa})^{ \Ll ( \R ) } \)) must be construed as \( ( \embeds_\CT^\kappa )^{ \Ll ( \R ) } \) and \( ( \embeds ^\kappa_\upsigma )^{ \Ll ( \R ) } \), respectively --- see Remark~\ref{rmk:maininnermodel}.
\item
As already noticed in the introduction after Theorem~\ref{th:introinvuniversal}, Theorem~\ref{th:largecardinals1} implies that every equivalence relation in \( \bGamma^2_1 \), a quite large boldface pointclass in \( \Vv \) which includes e.g.~all projective levels, is \( \Ll ( \R ) \)-reducible to a bi-embeddability relation \( \biembeds^\kappa_{\LL} \) (on structures of an appropriate uncountable size \( \kappa \)). 
This should be strongly contrasted with the case of the isomorphism relation: by Example~\ref{xmp:AD} (see~\cite[Theorem 9.18]{Hjorth:2000zr}), there are even \( \bSigma^1_1 \) equivalence relations \( E \) such that \( E \nleq_{\Ll ( \R )} {\cong ^\kappa_\LL} \) for \emph{every} \( \kappa \in ( \Cn )^{ \Ll ( \R ) } \) (and hence also for every cardinal in \( \Vv \)).
\item
In the statement of Theorem~\ref{th:largecardinals1} we could replace \( \Ll ( \R ) \) with \emph{any} inner model \( W \) containing all the reals of the universe.
\item
Working inside inner models satisfying stronger forms of determinacy, one can extend Theorem~\ref{th:largecardinals1} to more complex quasi-orders.
For example if \( W \) is an inner model of \( \AD_\R \) containing all the reals of the universe, then Theorem~\ref{th:largecardinals1} holds for all quasi-orders on \( \pre{\omega}{2} \) belonging to \( W \), with \( \bdelta^2_1 \) and \( \leq_{\Ll ( \R ) } \) replaced by \( \Theta^W \) and \( \leq_W \) respectively.
In this case we also get a characterization of the quasi-orders in \( W \) in terms of \( \leq_W \)-reducibility, namely: a quasi-order \( R \) on \( \pre{ \omega }{2} \) belongs to \( W \) if and only if \( R \leq_W {\embeds_\CT^{\Theta^W}} \).
\end{enumerate-(i)}
\end{remarks}

\section{Further completeness results} \label{sec:furthercompletenessresults}
\subsection{Representing arbitrary partial orders as embeddability relations}\label{subsec:Representingarbitrarypartialorders}
The methods developed in this paper yield also some purely combinatorial results, showing that embeddability relations can be quite complex.
Recall from Example~\ref{xmp:Parovicenko} that in \( \ZFC \), every partial order of size \( \kappa = \aleph_1 \) can be embedded into (the quotient order of) \( \embeds_\CT^\omega \): this is shown by proving that the relation \( ( \pow ( \omega ), \subseteq^* ) \) of inclusion modulo finite subsets on \( \pow ( \omega) \) is Borel reducible to \( \embeds_\CT^\omega \), and then using Parovicenko's theorem to embed any partial order \( P \) of size \( \aleph_1 \) into the inclusion relation on \( \pow ( \omega ) / \Fin \).
If we assume enough choice, weak forms of the above fact can be obtained for all uncountable \emph{small cardinals} --- for larger cardinals see~\cite{Mildenberger:2012ke}.

\begin{proposition} \label{prop:parovicenkoAC}
Let \( \omega < \kappa \leq 2^{\aleph_0} \) and assume \( \AC_\kappa ( \R) \). 
Then every partial order \( P \) of size \( \kappa \) can be embedded into the quotient order of \( \embeds_\CT^\kappa \). 
In fact, for every such \( P \) there is an \( \LL_{\kappa^+ \kappa} \)-sentence \( \upsigma \) (all of whose models are combinatorial trees) such that the quotient order of \( \embeds^{\kappa}_\upsigma \) is isomorphic to \( P \).
\end{proposition}

\begin{proof} 
Without loss of generality we may assume that \( P \) is \( ( \kappa , \unlhd ) \) . 
We will now assign to each \( \alpha < \kappa \) a combinatorial tree \( G_\alpha \) of size \( \kappa \) such that for all \( \alpha , \beta < \kappa \)
\begin{equation} \label{eq:P}
\alpha \unlhd \beta \IFF G_\alpha \embeds G_\beta.
 \end{equation}
The map \( \alpha \mapsto G_\alpha \) yields the desired embedding of \( P \) into the quotient order of \( \embeds^\kappa_\CT \).

Since by Theorem~\ref{th:LouveauRosendal} the relation of equality on \( \R \) is (Borel) reducible to \( \embeds^\omega_\CT \) and \( \kappa \leq 2^{\aleph_0} \), using \( \AC_\kappa ( \R) \) we can pick a sequence \( \seqof{S_\delta}{\delta < \kappa} \) of countable combinatorial trees (with disjoint domains) such that \( S_\delta \not\embeds S_{\delta'} \) for all distinct \( \delta , \delta' < \kappa \), and fix an arbitrary element \( x_\delta \in S_\delta \). 
Given \( \alpha < \kappa \), set 
\[
C_\alpha \coloneqq \setofLR{ 2\cdot \gamma }{ \gamma < \kappa } \cup \setofLR{ 2 \cdot \gamma +1 }{ \gamma < \kappa , \gamma \unlhd \alpha }. 
\]
The combinatorial tree \( G_\alpha \) is then (an isomorphic copy with domain \( \kappa \) of the graph) defined on the disjoint union 
\[
\set{ r_\alpha } \uplus {\bigcup \setofLR{S_\delta }{\delta \in C_\alpha} }
\] 
by connecting the vertex \( r_\alpha \) with each \( x_\delta \in S_\delta \) (for \( \delta \in C_\alpha \)). 
Using the fact that embeddings cannot decrease degrees of vertices, that \( r_\alpha \) has always degree \( \kappa \) and all other vertices of \( G_\alpha \) have degree \( \leq \omega < \kappa \), and that for all \( \alpha,\beta < \kappa \)
\[ 
\alpha \unlhd \beta \IFF \setofLR{\gamma < \kappa}{\gamma \unlhd \alpha} \subseteq \setofLR{\gamma < \kappa}{\gamma \unlhd \alpha} \IFF C_\alpha \subseteq C_\beta,
 \] 
we easily get that~\eqref{eq:P} is satisfied.

For the second part, one first check that the structure \( G_\alpha \) admits a Scott sentence \( \upsigma_\alpha \in \LL_{\kappa^+ \kappa} \) (see Remark~\ref{rem:maintheorem}\ref{rem:maintheorem-2} for the definition). 
Such a Scott sentence is obtained by formalizing as in Section~\ref{subsec:Uppsi} the conjunction of the following statements (here we use that \( \kappa \geq \omega_1 \) and that each \( S_\delta \), being countable, can be described by its quantifier-free type, which is an \( \LL^0_{\omega_1 \omega} \)-formula):
\begin{enumerate-(1)}
\item
the structure is a combinatorial tree;
\item
there is a unique vertex of degree \( > \omega \);
\item \label{en:description}
all other vertices belong to a unique countable substructure of \( G_\alpha \) which is isomorphic to some \( S_\delta \), and all these substructures are ``disjoint'' (vertices of different substructures are never connected by an edge);
\item
the vertex with degree \( > \omega \) is connected to a unique vertex in each of the substructures described in~\ref{en:description}, and such vertex ``corresponds'' to \( r_\delta \in S_\delta \);
\item
for each \( \delta < \kappa \), there is at most one substructure as in~\ref{en:description} which is isomorphic to \( S_\delta \);
\item
for each \( \delta < \kappa \), there is a substructure as in~\ref{en:description} which is isomorphic to \( S_\delta \) just in case \( \delta \in C_\alpha \).
\end{enumerate-(1)}
Then letting \( \upsigma \) be the \( \LL_{\kappa^+ \kappa} \) sentence \( \bigvee_{\alpha < \kappa} \upsigma_\alpha \), we get that \( \Mod^\kappa_\upsigma \) is the closure under isomorphism (inside \( \CT_\kappa \)) of the family \( \setofLR{G_\alpha}{\alpha < \kappa} \). 
Therefore \( P \) and the quotient order of \( \embeds^\kappa_\upsigma \) are isomorphic by~\eqref{eq:P}.
\end{proof}

The first part of Proposition~\ref{prop:parovicenkoAC} may be improved to the following result, in which we denote by \( \subseteq^*_\kappa \) the relation \( ( \pow ( \kappa), \subseteq^*) \) of inclusion modulo bounded subsets on \( \pow( \kappa) \).

\begin{proposition} \label{prop:parovicenkoACimproved}
Let \( \omega < \kappa \leq 2^{\aleph_0} \) and assume \( \AC_\kappa ( \R) \). 
Then \( {\subseteq^*_\kappa} \leq^\kappa_{\bB} {\embeds^\kappa_\CT} \).
\end{proposition}

\begin{proof}
Let \( T_0 \), \( T_1 \), and \( \seqofLR{S_\delta}{\delta<\kappa} \) be countable combinatorial trees such that:
\begin{itemize}[leftmargin=1pc]
\item
\( S_\delta \not\embeds S_{\delta'} \) for all distinct \( \delta,\delta' < \kappa \);
\item
\( S_\delta \not \sqsubseteq T_i \) and \( T_i \not \sqsubseteq S_\delta \) for every \( \delta < \kappa \) and \( i = 0,1 \);
\item
\( T_0 \sqsubseteq T_1 \) but \( T_1 \not \sqsubseteq T_0 \).
\end{itemize}
Combinatorial trees as above can easily be obtained by applying Theorem~\ref{th:LouveauRosendal} to the Borel quasi-order on \( \R \times \set{ 0 , 1 } \) defined by setting \( ( r , i ) \mathrel{R} ( s , j ) \IFF r = s \AND i \leq j \), and then using \( \AC_\kappa ( \R ) \).
For each \( \delta < \kappa \) and \( i = 0 , 1 \) fix arbitrary elements \( x_\delta \in S_\delta \) and \( y_i \in T_i \).

To each \( X \subseteq \kappa \) we now associate a combinatorial tree \( G_X \) of size \( \kappa \) as follows. 
Fix distinct vertices \( r \), \( r_\alpha \), and \( p_{ \alpha , \delta } \) for \( \alpha , \delta < \kappa \). 
For each \( \alpha , \delta < \kappa \) set
\[ 
i^X_{\alpha , \delta } \coloneqq
\begin{cases}
1 & \text{if } \delta< \alpha \text{ or } \delta \in X 
\\
0 & \text{otherwise.}
\end{cases}
 \] 
Append distinct copies of \( S_\delta \) and \( T_{ i^X_{ \alpha , \delta }} \) to the vertex \( p_{ \alpha , \delta } \) by connecting it with distinct edges to the (copies of the) distinguished vertices \( x_\delta \) and \( y_{i^X_{ \alpha , \delta } } \) of \( S_\delta \) and \( T_{ i^X_{ \alpha , \delta } } \), respectively. 
Then append all the combinatorial trees obtained in this way (for a fixed \( \alpha \) and arbitrary \( \delta < \kappa \)) to the vertex \( r_\alpha \) by connecting to it with an edge the vertices \( p_{ \alpha , \delta } \). 
Finally, add an edge between \( r \) and each \( r_\alpha \): the resulting combinatorial tree is \( G_X \).

The map associating (a suitable copy on \( \kappa \) of) the combinatorial tree \( G_X \) to each \( X \in \pow ( \kappa) \) is \( \kappa + 1 \)-Borel (in fact: continuous): we claim that it also reduces \( \subseteq^*_\kappa \) to \( \embeds \). 
Fix \( X,Y \subseteq \kappa \) and assume first that there is \( \beta < \kappa \) such that \( \delta \in X \IMPLIES \delta \in Y \) for all \( \beta \leq \delta < \kappa \). 
Consider the partial map \( e \) between \( G_X \) and \( G_Y \) sending \( r \) to itself, each \( r_\alpha \) to \( r_{\beta+\alpha} \), and, accordingly, each \( p_{\alpha,\delta} \) to \( p_{\beta+\alpha,\delta} \). 
Since the choice of \( \beta \) ensures that \( i^X_{\alpha,\delta} \leq i^Y_{\beta+\alpha,\delta} \) for all \( \delta < \kappa \), the map \( e \) can be completed to an embedding of \( G_X \) into \( G_Y \). 
Conversely, let \( e \) be an embedding of \( G_X \) into \( G_Y \). 
Since \( r \) is the unique vertex in both \( G_X \) and \( G_Y \) having \( \kappa \)-many neighbors of degree \( \kappa \), the embedding \( e \) must send \( r \) to itself and, consequently, \( r_0 \) to \( r_\beta \) for some \( \beta < \kappa \). 
Our choice of the \( S_\delta \) and \( T_i \) then implies that \( e ( p_{0,\delta}) = p_{\beta,\delta} \) and that \( e \) embeds \( T_{i^X_{0,\delta}} \) into \( T_{i^Y_{\beta,\delta}} \) (for all \( \delta < \kappa \)). 
Since \( T_i \sqsubseteq T_j \IFF i \leq j \), this shows that \( i^X_{0,\delta} \leq i^Y_{\beta,\delta} \) for all \( \delta < \kappa \), i.e.\ that \( \delta \in X \IMPLIES \delta \in Y \) for all \( \beta \leq \delta <\kappa \): thus \( X \subseteq^* Y \), as desired. 
\end{proof}

From Proposition~\ref{prop:parovicenkoACimproved} we in particular obtain that, under its assumptions, any partial order that can be embedded into (the quotient order of) \( ( \pow ( \kappa), \subseteq^*) \) can also be embedded into (the quotient order of) \( \sqsubseteq^\kappa_\CT \). 
This = applies to all partial orders of size \( \kappa \), but also to many other interesting cases. 
For example, P.\ Schlicht and K.\ Thompson have recently verified (personal communication) that a straightforward adaptation to the uncountable case of Parovicenko's proof shows that under \( \AC \) every linear order of size \( \aleph_{n+1} \) can be embedded into the quotient order of \( ( \pow ( \aleph_n ) , \subseteq^*) \).

In Theorem~\ref{th:parovicenko1} below we provide a counterpart to Proposition~\ref{prop:parovicenkoAC} in the \( \AD \)-world. 
The proof relies on the methods and results developed in this paper, and currently we do not see any simpler argument.

\begin{lemma}[\( \AD+\DC \)] \label{lem:parovicenko1}
Let \( \kappa \) be a Souslin cardinal, and \( P \) be an arbitrary partial order of size \( \kappa \). 
Then there is a \( \kappa \)-Souslin quasi-order \( R \) on \( \pre{\omega}{2} \) such that \( P \) embeds into the quotient order of \( R \). 
If \( \kappa < \bdelta_{\bS ( \kappa )} \) then \( R \) can be chosen so that its quotient order is isomorphic to \( P \).
\end{lemma}

\begin{proof}
We can assume that \( P \) is of the form \( ( \kappa , \unlhd ) \) itself. 
Since \( \pre{\omega}{2} \) and \( \pre{\omega}{\omega} \) are Borel isomorphic and \( \bS ( \kappa ) \) is closed under Borel preimages by Lemma~\ref{lem:Souslin}, it is enough to find a quasi-order \( R \) as in the statement, but defined on \( \pre{\omega}{\omega} \) instead of \( \pre{\omega}{2} \). 
Let \( \rho \) be an \( \bS ( \kappa ) \)-norm of length \( \kappa \) defined on a set \( A \subseteq \pre{\omega}{\omega} \) belonging to \( \bS ( \kappa ) \) (which exists by Proposition~\ref{prop:deltaSkappa}\ref{prop:deltaSkappa-b}), and let \( < \) be it strict part, i.e.~for \( x , y \in \pre{\omega}{\omega} \) set
\[ 
x < y \IFF x , y \in A \wedge \rho ( x ) < \rho ( y ).
 \] 
Then \( {<} \in \bS ( \kappa ) \), and its rank function (which is just the \( \bS ( \kappa ) \)-norm \( \rho \)) is onto \( \kappa \). 
Now consider the function \( f \colon \kappa \times \kappa \to \pow ( \pre{\omega}{2} ) \) defined by
\[ 
f ( \alpha , \beta ) \coloneqq
\begin{cases}
\pre{\omega}{2} & \text{if } \alpha \unlhd \beta , 
\\
\emptyset & \text{otherwise}.
\end{cases}
 \] 
By Moschovakis' first coding lemma~\cite[Lemma 7D.5]{Moschovakis:2009fk}, there is an \( \bS ( \kappa ) \) set \( C \subseteq \pre{\omega}{\omega} \times \pre{\omega}{\omega} \times \pre{\omega}{2} \) which is a \emph{choice set} for \( f \), i.e.~such that for all \( x_1 , x_2 \in \pre{\omega}{\omega} \), \( y \in \pre{\omega}{2} \), and \( \alpha , \beta < \kappa \) the following hold:
\begin{enumerate-(1)}
\item 
\( ( x_1 , x_2 , y ) \in C \IMPLIES x_1, x_2 \in A \wedge y \in f ( \rho ( x_1 ) , \rho ( x_2 ) ) \);
\item
\( f ( \alpha , \beta ) \neq \emptyset \IMPLIES \EXISTS{x_1,x_2 \in A} \EXISTS{y \in \pre{\omega}{2}} [ \rho ( x_1 ) = \alpha \wedge \rho ( x_2 ) = \beta \wedge ( x_1 , x_2 , y ) \in C] \).
\end{enumerate-(1)}
Let \( R \) be the quasi-order on \( \pre{\omega}{\omega} \) defined by 
\begin{multline*}
z_1 \mathrel{R} z_1 \IFF z_1 = z_2 \OR [z_1, z_2 \in A \AND \EXISTS{x_1,x_2 \in A} \EXISTS{y \in \pre{\omega}{2}}( \rho ( x_1 ) = \rho ( z_1 ) \AND \\ \rho ( x_2 ) = \rho ( z_2 ) \AND ( x_1 , x_2 , y ) \in C)] .
 \end{multline*}
Then \( R \) is \( \kappa \)-Souslin by the closure properties of \( \bS ( \kappa) \) (see Lemma~\ref{lem:Souslin}), and it is straightforward to check that, by definition of \( C \) and the fact that it is a choice set for \( f \), the map \( \alpha \mapsto \setofLR{z \in A}{\rho ( z ) = \alpha} \) is an isomorphism between \( P \) and the quotient order of \( R \restriction (A \times A) \), and hence it is an embedding of \( P \) into the quotient order of the whole \( R \).

The additional part concerning those \( \kappa \) which are smaller than \( \bdelta_{\bS ( \kappa )} \) follows from the fact that in this case in the argument above we can we can let \( \rho \) be the rank function of any \( \bDelta_{\bS ( \kappa )} \) prewellordering of the whole \( \pre{\omega}{\omega} \) of length \( \kappa \), so that \( A = \pre{\omega}{\omega} \).
\end{proof}

\begin{theorem}[\( \AD + \DC \)] \label{th:parovicenko1}
Let \( \kappa \) be a Souslin cardinal. 
Then every partial order \( P \) of size \( \kappa \) can be embedded into the quotient order of \( \embeds_\CT^\kappa \). 
In fact, if \( \kappa < \bdelta_{ \bS ( \kappa ) } \), then for every such \( P \) there is an \( \LL_{\kappa^+ \kappa} \)-sentence \( \upsigma \) (all of whose models are combinatorial trees) such that the quotient order of \( \embeds^{\kappa}_\upsigma \) is isomorphic to \( P \).
\end{theorem}

Notice that by Proposition~\ref{prop:deltaSkappa}\ref{prop:deltaSkappa-c} the second part of Theorem~\ref{th:parovicenko1} can be applied exactly when \( \kappa \) is such that \( \bS ( \kappa) \) is not closed under coprojections, and thus, in particular, when \( \kappa \) is one of the projective cardinals \( \bdelta^1_n \) or, more generally, when \( \kappa \) is not a regular limit of Souslin cardinals.

\begin{proof}
Let \( R \) be as in Lemma~\ref{lem:parovicenko1}, and let \( T \in \TT_\kappa \) be such that \( R = \PROJ \body{ T } \). 
By Theorem~\ref{th:graphs}, the function \( f_T \) defined in~\eqref{eq:f_T} reduces \( R \) to \( \embeds_\CT^\kappa \), and thus its quotient map 
\[
 \hat{f} _T \colon \pre{\omega}{\omega} / R \to \CT_\kappa / {\biembeds_\CT^\kappa} , \quad \eq{x}_{E_R} \mapsto \eq{ f_T ( x ) }_{\biembeds_\CT^\kappa} 
\] 
where \( {\biembeds_\CT^\kappa} \) is the bi-embeddability relation on \( \CT_\kappa \), embeds the quotient order \( \pre{\omega}{\omega} / R \) into the quotient order \( \CT_\kappa / {\biembeds_\CT^\kappa} \). 
Therefore, composing the embedding of \( P \) into \( \pre{\omega}{\omega} / R \) with \( \hat{f}_T \) gives the desired embedding.

If \( \kappa < \bdelta_{\bS ( \kappa)} \), then by Lemma~\ref{lem:parovicenko1} there is a \( \kappa \)-Souslin quasi-order \( R \) on \( \pre{\omega}{2} \) whose quotient order is isomorphic to \( P \). 
By Corollary~\ref{cor:maintheorem}, there is an \( \LL_{\kappa^+ \kappa} \)-sentence \( \upsigma \) (all of whose models are combinatorial trees) such that the quotient order of \( \embeds^\kappa_\upsigma \) is isomorphic to \( \pre{\omega}{\omega} / R \), and hence also to \( P \). 
\end{proof}

Currently we do not know if it is possible to obtain an analogue of Proposition~\ref{prop:parovicenkoACimproved} in the \( \AD \)-context.

\subsection{Other model theoretic examples} \label{subsec:othermodeltheoretic}
In this section we consider some model theoretic variants of our results on the embeddability relation between combinatorial trees. 
We consider two kinds of variations: in the first one, we change the morphism between combinatorial trees under consideration by replacing embeddings with full homomorphisms, while in the second one we consider the embeddability relation again but we change the nature of the underlying structures. 
Of course these are only a few examples of the many interesting model theoretic variants of the problem considered in this paper (namely, the descriptive set-theoretic complexity of natural relations between uncountable structures) which would deserve some attention in the near future.

\subsubsection{Changing the morphism: full homomorphisms between combinatorial trees} \label{subsubsec:changingmorphism}
Recall that \( \LL = \setLR{ \edge } \) is the language of graphs consisting of a single binary relational symbol.

\begin{definition} \label{def:fullhomomorphism}
Let \( X = \seq{X; \edge^X} \) and \( Y = \seq{Y; \edge^Y} \) be two \( \LL \)-structures. 
A map \( f \colon X \to Y \) is called \markdef{full homomorphism} between \( X \) and \( Y \) if for all \( a,b \in X \)
\[ 
a \edge^X b \IFF f(a) \edge^Y f(b).
 \] 
\end{definition}

Thus embeddings are just \emph{injective} full homomorphisms. 
Moreover, since the identity function on an \( \LL \)-structure is a full homomorphism, and since the composition of two full homomorphisms is still a full homomorphism, the binary relation defined below is reflexive and transitive, i.e.\ a quasi-order.

\begin{notation} \label{notation:fullhomomorphism}
Given two \( \LL \)-structures \( X = \seq{X; \edge^X} \) and \( Y = \seq{Y; \edge^Y} \), we set \( X \embeds^h Y \) if and only if there exists a full homomorphism between \( X \) and \( Y \).
\end{notation}

Replacing \( \embeds \) with \( \embeds^h \) in Definitions~\ref{def:invuniversal} and~\ref{def:invuniversallocal} we get a corresponding notion of invariant universality for \( \embeds^h \).

\begin{definition}\label{def:invunivesalhom}
Let \( \mathcal{C} \) be a class of quasi-orders, \( \LL \) be the graph language, and \( \kappa \) be an infinite cardinal. 
Given an \( \LL_{\kappa^+ \kappa} \)-sentence \( \uptau \), the relation \( {\embeds}^h \restriction \Mod^\kappa_{\uptau} \) is \markdef{invariantly universal for \( \mathcal{C} \)} if for every \( R \in \mathcal{C} \) there is an \( \LL_{\kappa^+ \kappa} \)-sentence \( \upsigma \) such that \( \Mod_ \upsigma^ \kappa \subseteq \Mod_ \uptau ^ \kappa \) and \( R \sim {\embeds}^h \restriction \Mod^\kappa_\upsigma \).

As usual, when the reducibility \( \leq \) is replaced by one of its restricted form \( \leq_* \) we speak of \markdef{\( \leq_* \)-invariant universality}.
\end{definition}

We are now going to show that for every infinite cardinal \( \kappa \) and every \( T \in \TT_\kappa \), the relations \( \embeds \) and \( \embeds^h \) coincide on \( \Mod^\kappa_{\upsigma_T} \), where \( \upsigma_T \) is as in Corollary~\ref{cor:saturation}. 
To prove this, we will use the following variant of~\cite[Proposition 3.10]{Friedman:2011cr}, which can be proved in the same way.

\begin{lemma} \label{lem:basichomomorphism}
Suppose that \( G , G' \) are combinatorial trees and \( f \colon G \to G' \) is a full homomorphism. 
If \( a , b \in G \) are distinct and \( f ( a ) = f ( b ) \), then \( a \) and \( b \) have both degree \( 1 \) in \( G \) and have geodesic distance \( 2 \) from each other (i.e.\ they share their unique adjacent vertex).
\end{lemma}

This essentially shows that a full homomorphism between combinatorial trees is already almost injective, and allows us to prove the following.

\begin{theorem} \label{thm:mainhomomorphism}
Let \( \kappa \) be an uncountable cardinal, \( T \in \TT_\kappa \) be such that \( R = \PROJ \body{ T } \) is a quasi-order, and \( \upsigma_T \) be as in Corollary~\ref{cor:saturation}. 
Then the relation \( {\embeds}^h \restriction \Mod^\kappa_{\upsigma_T} \) coincides with \( \embeds^\kappa_{\upsigma_T} \) (that is, with the embeddability relation \( {\embeds} \restriction \Mod^\kappa_{\upsigma_T} \)).
\end{theorem}

\begin{proof}
First recall that by Corollary~\ref{cor:saturation} the set \( \Mod^\kappa_{\upsigma_T} \) is the closure under isomorphism of \( \ran ( f_T ) \), where \( f_T \) is as in~\eqref{eq:f_T}. 
Since embeddings are in particular full homomorphisms, it is enough to show that for every \( x,y \in \pre{\omega}{2} \)
\[ 
f_T ( x ) \embeds^h f_T ( y ) \IMPLIES f_T ( x ) \embeds f_T ( y ) .
 \] 
Let \( j \) be a full homomorphism between \( f_T ( x ) = \mathbb{G}_{\Sigma_T ( x )} \) and \( f_T ( y ) = \mathbb{G}_{\Sigma_T ( y )} \). 
By Lemma~\ref{lem:basichomomorphism}, its restriction to \( G_0 \cup \mathrm{F} ( \mathbb{G}_{\Sigma_T ( x )} ) \) (see Definition~\ref{def:G_0} and the notation introduced in~\eqref{eq:substructuresofG_T-a}--\eqref{eq:substructuresofG_T-g} of Section~\ref{subsec:G_T}) is injective, and hence an embedding. 
Arguing as in the second half of the proof of Theorem~\ref{th:graphs}\ref{th:graphs-a}, we thus get that \( j \restriction ( G_0 \cup \mathrm{F} ( \mathbb{G}_{\Sigma_T ( x )} ) ) \) is an embedding of the subgraph of \( f_T ( x ) \) with domain \( G_0 \cup \mathrm{F} ( \mathbb{G}_{\Sigma_T ( x )}) \) into the subgraph of \( f_T ( y ) \) with domain \( G_0 \cup \mathrm{F} ( \mathbb{G}_{\Sigma_T ( y )} ) \). 
By the proof of Theorem~\ref{th:graphs}\ref{th:graphs-a} (see Remark~\ref{rmk:partialembedding}), this implies that \( x \mathrel{R} y \), whence \( f_T ( x ) \embeds f_T ( y ) \) by Theorem~\ref{th:graphs}\ref{th:graphs-a} again.
\end{proof}

\begin{remark}\label{rmk:variantmainhomomorphism}
By adapting in the obvious way Definitions~\ref{def:fullhomomorphism} and~\ref{def:invunivesalhom} to the language \( \bar{\LL} \) of \emph{ordered} combinatorial trees (see Section~\ref{sec:alternativeapproach}), one can easily check that in Theorem~\ref{thm:mainhomomorphism} we can also replace \( \upsigma_T \) with the sentence \( \bar\upsigma_T \) from Corollary~\ref{cor:barsaturation}. 
In fact, the ordered combinatorial trees \( \bar{f}_T ( x ) \) and \( \bar{f}_T ( y ) \) (where \( \bar{f}_T \) is as in~\eqref{eq:defbarf_T}) already have domain \( G_0 \cup \mathrm{F} ( \mathbb{G}_{\Sigma_T ( x )}) \) and \( G_0 \cup \mathrm{F}( \mathbb{G}_{\Sigma_T ( y )}) \), respectively. 
Thus any full homomorphism between \( \bar{f}_T ( x ) \) and \( \bar{f}_T ( y ) \) is an embedding by Lemma~\ref{lem:basichomomorphism} (which holds for ordered combinatorial trees as well, since these are \( \bar{\LL} \)-expansions of combinatorial trees).
\end{remark}

Theorem~\ref{thm:mainhomomorphism} and Remark~\ref{rmk:variantmainhomomorphism} show that in all the completeness and invariant universality results from Sections~\ref{sec:embeddabilitygraphs}--\ref{sec:applications} and~\ref{subsec:Representingarbitrarypartialorders} we could systematically replace the embeddability relation \( \embeds \) with the one induced by full homomorphisms \( \embeds^h \).

\subsubsection{Changing the structures: partial orders and lattices}
At the end of~\cite[Section 3.1]{Louveau:2005cq}, Louveau and Rosendal provided a continuous map from countable connected graphs to lattices (viewed as partial orders) which reduces the embeddability relation to itself. 
Such a construction can straightforwardly be adapted to the uncountable case.
For technical reasons which will be clear shortly we have to modify the construction a bit.

\begin{definition} \label{def:lattices}
Let \( \kappa \) be any infinite cardinal. 
To any graph \( G \) on \( \kappa \) we associate the following lattice \( L_G = ( \kappa , \preceq_G ) \):
\begin{itemize}[leftmargin=1pc]
\item
\( 0 \preceq_G \alpha \preceq_G 1 \) for every \( \alpha \in \kappa \)
\item
if in \( G \) there is an edge between \( \alpha, \beta \in \kappa \), then both \( 2+(2 \cdot \alpha) \preceq_G 2 + (2\cdot\op{\alpha}{\beta}+1) \) and \( 2 + ( 2 \cdot \beta ) \preceq_G 2 + ( 2 \cdot \op{\alpha}{\beta} + 1 ) \) (where \( \op{\cdot}{\cdot} \) is the pairing function from~\eqref{eq:Hessenberg});
\item
no other \( \preceq_G \)-relation holds.
\end{itemize}
\end{definition}

Notice that the \( L_G \)'s are actually \emph{complete} lattices.

\begin{remark} \label{rmk:lattices}
The map \( G \mapsto L_G \) is continuous when both spaces of models (namely, the space of graphs of size \( \kappa \) and the space of lattices of size \( \kappa \)) are endowed with the same topology \( \tau_p \) or \( \tau_b \).
Our modification of the original Louveau-Rosendal construction is required to obtain continuity with respect to the topology \( \tau_p \).
\end{remark}

\begin{theorem} \label{thm:lattices}
The map \( G \mapsto L_G \) from Definition~\ref{def:lattices} reduces the embeddability (respectively, the isomorphism) relation between connected graphs to the embeddability (respectively, the isomorphism) relation between lattices.
\end{theorem}

\begin{proof}
This is just a minor variation of the proof of~\cite[Theorem 3.3]{Louveau:2005cq}. 
If \( j \) is an embedding between the connected graphs \( G \) and \( G' \), then the map defined by 
\begin{align*}
0 & \mapsto 0 
\\
 1 & \mapsto 1 
 \\ 
2 + ( 2 \cdot \alpha ) & \mapsto 2 + ( 2 \cdot j ( \alpha ) ) 
\\
 2 + ( 2 \cdot \op{\alpha}{\beta} + 1 ) & \mapsto 2 + ( 2 \cdot \op{ j ( \alpha ) }{ j ( \beta ) } + 1 ) 
\end{align*}
is an embedding of \( L_G \) into \( L_{G'} \). 

Conversely, let \( j \) be an embedding of \( L_G \) into \( L_{G'} \). 
Notice that in both \( L_G \) and \( L_{G'} \) we have that an ordinal \( 1 \neq \alpha \in \kappa \) is an immediate predecessor of the maximum \( 1 \) if and only if it is odd (here we use the fact that \( G \) and \( G' \) are connected). 
Since \( j \) is an embedding, we thus get that for every \( \alpha \in \kappa \) there is \( \gamma_\alpha \in \kappa \) such that \( j ( 2 + ( 2 \cdot \alpha)) = 2 + ( 2 \cdot \gamma_\alpha) \). 
Since \( \alpha,\beta \in \kappa \) are connected by an edge in \( G \) if and only if \( 2 + (2 \cdot \alpha) \) and \( 2 + (2 \cdot \beta) \) share a \( \preceq_G \)-successor distinct from the maximum \( 1 \) (and the same is true for \( G' \)), it follows that the map \( \alpha \mapsto \gamma_\alpha \) is an embedding of \( G \) into \( G' \).

The same proof works for the case of the isomorphism relation, so we are done.
\end{proof}

Theorem~\ref{thm:lattices} and Remark~\ref{rmk:lattices} show that in all the completeness results from Sections~\ref{sec:embeddabilitygraphs}--\ref{sec:applications} and~\ref{subsec:Representingarbitrarypartialorders} we could systematically replace combinatorial trees with (complete) lattices, and hence, in particular, with partial orders. 
Moreover, both Theorem~\ref{thm:lattices} and Remark~\ref{rmk:lattices} remain true if the lattice \( L_G \) associated to the graph \( G \) is construed as a bounded lattice in the algebraic sense, that is as a structure in the language consisting of two binary function symbols (the join and the meet operations \( \vee \) and \( \wedge \)) and two constant symbols (the minimum and the maximum \( 0 \) and \( 1 \)). 
Thus we could also further replace combinatorial trees with such algebraic structures.

\subsection{Isometry and isometric embeddability between complete metric spaces of density character \( \kappa \)} \label{subsec:isomemb}

Recall from Section~\ref{subsubsec:spacesofmetricspaces} the standard Borel \( \kappa \)-space \( \MMM_\kappa \) of (codes for) complete metric spaces of density character \( \kappa \), and its subspaces \( \DDD_\kappa \) and \( \UUU_\kappa \) consisting of, respectively, discrete and ultrametric spaces. 
In this section we will provide some informations on the complexity of the isometry relation \( \cong^i \), and the isometric embeddability relation \( \sqsubseteq^i \), which are both \( \kappa \)-analytic quasi-orders on \( \MMM_\kappa \), and on their restrictions to \( \DDD_\kappa \) and \( \UUU_\kappa \). Some of the results concerning \( \cong^i \) already appeared in~\cite{mottoros:2016bu}, but we include them in our presentation as well for the reader's convenience.

\subsubsection{The discrete case} \label{subsubsec:isomembdiscrete}

We first consider the case of \emph{discrete} metric spaces of density character (equivalently, of size) \( \kappa \). 
Fix strictly positive \( r_0, r_1 \in \R \) such that \( r_0 < r_1 \leq 2r_0 \). 
To each graph \( G \) on \( \kappa \) we associate the discrete metric space \( D_G \) on \( \kappa \) with distance \( d_G \) defined by
\[ 
d_G ( \alpha,\beta) \coloneqq
\begin{cases}
0 & \text{if } \alpha = \beta 
\\
r_0 & \text{if }\alpha \neq \beta \text{ and \( \alpha \) and \( \beta \) are adjacent in } G 
\\
r_1 & \text{if } \alpha \neq \beta \text{ and \( \alpha \) and \( \beta \) are not adjacent in } G. 
\end{cases}
 \] 
Our choice of \( r_0 \) and \( r_1 \) guarantees that the triangular inequality is satisfied by \( d_G \). 
Since the space \( D_G \) is already defined on \( \kappa \), it can canonically be identified with its code \( x_G \in \DDD_\kappa \) obtained by setting, as in~\eqref{eq:defx_M}, \( x_G ( \alpha,\beta,q) = 1 \IFF d_G ( \alpha,\beta) < q \) for every \( \alpha, \beta < \kappa \) and \( q \in \QQ^+ \) .

\begin{remark} \label{rmk:continuitydiscrete}
Let \( \Mod^\kappa_\GR \) be the space of all graphs on \( \kappa \). 
Then the map 
\begin{equation} \label{eq:deftheta}
\theta_D \colon \Mod^\kappa_\GR \to \DDD_\kappa \subseteq \MMM_\kappa, \qquad G \mapsto x_G
\end{equation}
is continuous when both \( \Mod^\kappa_\GR \) and \( \MMM_\kappa \) are endowed with the same topology \( \tau_p \) or \( \tau_b \).
\end{remark}

\begin{lemma} \label{lem:reductiontodiscrete}
The map \( \theta_D \) from~\eqref{eq:deftheta} simultaneously reduces \( \cong \) to \( \cong^i \) and \( \embeds \) to \( \sqsubseteq^i \).
\end{lemma}

\begin{proof}
It is easy to check that given any \( G,G' \in \Mod^\kappa_\GR \) and any \( \varphi \colon \kappa \to \kappa \), the map \( \varphi \) is an isomorphism (respectively, an embedding) between \( G \) and \( G' \) if and only if it is an isometry (respectively, an isometric embedding) between \( D_G \) and \( D_{G'} \).
\end{proof}

This yields the following lower bounds for the complexity of \( \cong^i \restriction \DDD_\kappa \) and \( \sqsubseteq^i \restriction \DDD_\kappa \). 
As usual, we endow both \( \Mod^\kappa_\GR \) and \( \MMM_\kappa \) with the bounded topology and \( \DDD_\kappa \) with the induced relative topology. 
Moreover, to simplify the notation we denote by \( \cong^\kappa_\GR \) and \( \embeds^\kappa_\GR \) the isomorphism and the embeddability relation on \( \Mod^\kappa_\GR \), respectively.

\begin{theorem} \label{th:reductiontodiscrete}
Let \( \kappa \) be any infinite cardinal.
\begin{enumerate-(a)}
\item \label{th:reductiontodiscrete-a}
\( {\cong^\kappa_\GR} \leq^\kappa_\bB {\cong^i \restriction \DDD_\kappa} \) and \( {\embeds^\kappa_\GR} \leq^\kappa_\bB {\sqsubseteq^i \restriction \DDD_\kappa} \).
\item \label{th:reductiontodiscrete-b}
The relation \( \sqsubseteq^i \restriction \DDD_\kappa \) is \emph{\( \leq^\kappa_\bB \)-complete} for the class of \( \kappa \)-Souslin quasi-orders on Polish or standard Borel spaces.
\item \label{th:reductiontodiscrete-c}
(\( \AC \)) 
If \( \kappa \leq 2^{\aleph_0} \) and there is an \( \bS(\kappa) \)-code for \( \kappa \) (which is always the case for \( \kappa = \omega \), \( \kappa = \omega_1 \), and \( \kappa = 2^{\aleph_0} \)), then \( \sqsubseteq^i \restriction \DDD_\kappa \) is \emph{\( \leq_{\bS(\kappa)} \)-complete} for \( \kappa \)-Souslin quasi-orders on Polish or standard Borel spaces.
\item \label{th:reductiontodiscrete-d}
(\( \AD + \DC \)) 
If \( \kappa \) is a Souslin cardinal, then \( \sqsubseteq^i \restriction \DDD_\kappa \) is \emph{\( \leq_{\bS(\kappa)} \)-complete} for \( \kappa \)-Souslin quasi-orders on Polish or standard Borel spaces.
\item \label{th:reductiontodiscrete-e}
Let \( W \supseteq \R \) be an inner model and let \( \kappa \in \Cn^W \). 
Then \( \sqsubseteq^i \restriction \DDD_\kappa \) is \emph{\( \leq_W \)-complete} for quasi-orders in \( \bS_W(\kappa) \).
\end{enumerate-(a)}
\end{theorem}

\begin{proof}
\ref{th:reductiontodiscrete-a} 
By Lemma~\ref{lem:reductiontodiscrete}, this is witnessed by the map from~\eqref{eq:deftheta}, which is \( \tau_b \)-continuous (and hence \( \kappa+1 \)-Borel) by Remark~\ref{rmk:continuitydiscrete}.

\ref{th:reductiontodiscrete-b} 
We argue as in the proof of Theorem~\ref{th:definable=borel}. 
Let \( R = \PROJ{\body{T}} \) be an arbitrary \( \kappa \)-Souslin quasi-order on \( \pre{\omega}{2} \), and consider the map \( \theta_D \circ f_T \colon \pre{\omega}{2} \to \DDD_\kappa \subseteq \MMM_\kappa \), where \( f_T \) is as in~\eqref{eq:f_T} and \( \theta_D \) is as in~\eqref{eq:deftheta}. 
By Remark~\ref{rmk:continuitydiscrete}, such map is continuous when both \( \pre{\omega}{2} \) and \( \MMM_\kappa \) are endowed with the product topology, and hence it is (effective) weakly \( \kappa+1 \)-Borel. 
Moreover, it reduces \( R \) to \( \sqsubseteq^i \) by Theorem~\ref{th:graphs} and Lemma~\ref{lem:reductiontodiscrete}, so we are done.

\ref{th:reductiontodiscrete-c}--\ref{th:reductiontodiscrete-d} 
Argue as in the proof of Theorem~\ref{th:Gammacompleteness}, using again the fact that \( \theta_D \circ f_T \colon \pre{\omega}{2} \to \DDD_\kappa \) is continuous when both spaces are endowed with the product topology and that it reduces the \( \kappa \)-Souslin quasi-order \( R = \PROJ{\body{T}} \) to \( \sqsubseteq^i \).

\ref{th:reductiontodiscrete-e} 
Use Theorem~\ref{th:graphsinnermodel} together with the fact that the map \( \theta_D \) from~\eqref{eq:deftheta} is definable in \( W \) and that Lemma~\ref{lem:reductiontodiscrete}, which is a theorem of \( \ZF \), is true in \( W \) as well.
\end{proof}

By Theorem~\ref{th:reductiontodiscrete}, in all the completeness results from Sections~\ref{sec:definablecardinality}, \ref{sec:applications}, and~\ref{subsec:Representingarbitrarypartialorders} we may systematically replace \( \embeds^\kappa_\CT \) with \( \sqsubseteq^i \restriction \DDD_\kappa \); a sample of (stronger forms) of these variants will be presented in Section~\ref{subsubsec:ultrametric} (see Theorems~\ref{th:ultrametricapplications}, \ref{th:ultrametricombinatorial} and Remark~\ref{rmk:variants}). 

\subsubsection{The ultrametric case} \label{subsubsec:ultrametric}

In this section we consider the case of complete \emph{ultrametric} spaces of density character \( \kappa \) and prove some strengthenings of the results from the previous section. 
Most of the definitions and results are direct generalizations to the uncountable context of the material from~\cite{Camerlo:2015ar} --- see also~\cite{mottoros:2016bu}. 

Let \( U = (U,d_U) \) be an ultrametric space. 
Then by~\cite[Lemma 2.20]{mottoros:2016bu} for every non-negative \( r \in \R \) and every dense \( D \subseteq U \)
\begin{equation} \label{eq:realizeonadense}
\EXISTS{ x , y \in U} ( d_U ( x , y ) = r ) \IFF \EXISTS{ x , y \in D} ( d_U ( x , y ) = r ).
\end{equation}
It follows that if \( U \) has density character \( \kappa \), then its set of (realized) distances
\[ 
D ( U ) \coloneqq \setofLR{r \in \R}{\EXISTS{ x , y \in U} ( d_U ( x , y ) = r)}
 \] 
has size at most \( \kappa \). 
Conversely, given any set \( A \subseteq \R \) of non-negative reals containing \( 0 \) and of size \( \leq \kappa \), we can form a complete%
\footnote{Completeness follows automatically from the fact that \( U(A) \) is discrete.} 
ultrametric space \( U ( A ) = ( U_A , d_A ) \) of density character \( \kappa \) such that \( D ( U ( A ) ) = A \) by setting \( U_A \coloneqq \setofLR{ ( r , \alpha ) }{ r \in A, \alpha < \kappa} \) and for every \( r , r' \in A \) and \( \alpha , \beta < \kappa \)
\begin{equation} \label{eq:canonicalultrametric}
 d_A ( ( r , \alpha ) , ( r' , \beta ) ) \coloneqq
\begin{cases}
0 & \text{if } r = r' \text{ and } \alpha = \beta 
\\
\max \setLR{ r , r' } & \text{otherwise.}
\end{cases}
 \end{equation}
Thus we have that a set \( A \subseteq \R \) is the set of (realized) distances of a complete ultrametric space of density character \( \kappa \) if and only if 
\( 0 \in A \subseteq \cointerval{0}{+ \infty} \) and \( A \) has size \( \leq \kappa \). 
This suggests to consider the following subclasses of \( \UUU_\kappa \) (where for simplicity we denote by \( U_x \) the complete space coded by \( x \in \UUU_\kappa \) --- see~Section~\ref{subsubsec:spacesofmetricspaces}).

\begin{definition} \label{def:subclassofultrametric}
Let \( \kappa \) be an infinite cardinal and \( A \subseteq \R \) be such that \( 0 \in A \subseteq\cointerval{0}{+\infty } \) and \( \card{A} \leq \kappa \). 
Then we set
\[ 
\UUU_\kappa(A) \coloneqq \setofLR{x \in \UUU_\kappa}{D ( U_x ) \subseteq A}
\] 
and
\[ 
\UUU^\star_\kappa(A) \coloneqq \setofLR{x \in \UUU_\kappa}{D ( U_x ) = A}.
\] 
\end{definition}
Using~\eqref{eq:realizeonadense}, it is easy to see that both \( \UUU_\kappa ( A ) \) and \( \UUU^\star_\kappa ( A ) \) are effective \( \kappa+1 \)-Borel subsets of \( \UUU_\kappa \): thus since the latter is a standard Borel \( \kappa \)-space, \( \UUU_\kappa ( A ) \) and \( \UUU^\star_\kappa(A) \) are standard Borel \( \kappa \)-spaces as well. 
Notice also that the spaces \( \UUU^{(\star)}_\kappa(A) \) are always nonempty because \( U(A) \in \UUU^\star_\kappa(A) \subseteq \UUU_\kappa(A) \). 
However, the sets \( A \) that can be considered in these constructions may vary from a model of set theory to another: for example, in models of \( \AD \) only countable \( A \)'s can satisfy the conditions in Definition~\ref{def:subclassofultrametric} because there are no uncountable well-orderable subsets of \( \R \) by the (\( \omega \)-)\( \PSP \). 

We are now going to show that when a set \( A \) as in Definition~\ref{def:subclassofultrametric} is ill-founded with respect to the usual ordering of \( \R \), then the restriction of the isometry relation \( \cong^i \) and of the isometric embeddability relation \( \sqsubseteq^i \) to \( \UUU^\star_\kappa(A) \) are quite complex. 
This in particular yields strengthenings of the results from Section~\ref{subsubsec:isomembdiscrete} because when \( 0 \) is not an accumulation point of \( A \) we have 
\[ 
\UUU^\star_\kappa ( A ) \subseteq \UUU_\kappa ( A ) \subseteq \UUU_\kappa \cap \DDD_\kappa.
\]

The following construction essentially corresponds to the special case \(\alpha = \omega \) in~\cite[Section 3]{mottoros:2016bu}, where it is shown that it yields a reduction from isomorphism to isometry; here we are going to show that the same construction yields also a reduction of embeddability to isometric embeddability.
Fix a strictly decreasing sequence \( \vec{r} = \seqofLR{r_n}{n \in \omega} \) of positive real numbers. 
Given a rooted combinatorial tree \( G \) on \( \kappa \), define a partial order \( \preceq_G \) on \( \kappa \) by setting 
\[ 
\alpha \preceq_G \beta \IFF \text{ the unique path in \( G \) connecting the root \( r \) of \( G \) to \( \beta \) passes through \( \alpha \)}. 
\]
Define a metric \( d^{\vec{r}}_G \) on \( \kappa \) by setting \( d^{\vec{r}}_G ( \alpha , \beta ) \coloneqq r_{ n ( \alpha , \beta )} \), where 
\[ 
n ( \alpha , \beta ) \coloneqq \max \setof{ l_G ( \gamma ) }{\gamma \in \kappa, \gamma \preceq_G \alpha , \beta} 
\] 
and \( l_G ( \gamma ) \) is the length of the path connecting \( r \) to \( \gamma \), and then consider the completion \( U^{\vec{r}}_G = (U_G,d^{\vec{r}}_G) \) of the space \( (\kappa,d^{\vec{r}}_G) \). 
It is clear from the construction that \( U^{\vec{r}}_G \) is a complete ultrametric space of density character \( \kappa \), and that its canonical code \( u^{\vec{r}}_G \in \UUU_\kappa \) defined by setting for every \( \alpha ,\beta < \kappa \) and \( q \in \Q^+ \)
\[ 
u^{\vec{r}}_G ( \alpha , \beta, q ) = 1 \IFF d^{\vec{r}}_G ( \alpha , \beta ) < q
 \] 
belongs to \( \UUU_\kappa(A) \) for any \( A \) as in Definition~\ref{def:subclassofultrametric} containing all the \( r_n \)'s.

\begin{remark} \label{rmk:continuityultrametric}
Let \( \RCT_\kappa \) be the space of (codes for) all rooted combinatorial trees on \( \kappa \), \( A \) be as in Defintion~\ref{def:subclassofultrametric} and ill-founded, and \( \vec{r} \) be a strictly decreasing sequence of elements of \( A \). 
Then the map 
\begin{equation} \label{eq:defthetaU}
\theta^{\vec{r}}_U \colon \RCT_\kappa \to \UUU_\kappa(A) \subseteq \MMM_\kappa, \qquad G \mapsto u^{\vec{r}}_G
\end{equation}
is continuous when both \( \RCT_\kappa \) and \( \UUU_\kappa(A) \subseteq \MMM_\kappa \) are endowed with the same topology \( \tau_p \) or \( \tau_b \).
\end{remark}

\begin{lemma} \label{lem:reductiontoultrametric}
The map \( \theta^{\vec{r}}_U \) from~\eqref{eq:defthetaU} simultaneously reduces \( \cong \) to \( \cong^i \) and \( \embeds \) to \( \sqsubseteq^i \).
\end{lemma}

\begin{proof}
If \( \inf \vec{r} > 0 \), then \( U^{\vec{r}}_G = (\kappa, d^{\vec{r}}_G) \) and one can simply use the straightforward adaptation to the uncountable context of the proof of~\cite[Theorem 5.2]{Camerlo:2015ar} --- the construction provided in this paper coincide with the one from~\cite[Section 5.1]{Camerlo:2015ar} once we identify each element of \( \RCT_\kappa \) with any of its isomorphic copies having domain a subset of \( \pre{<\omega}{\kappa} \) closed under subsequences and root the empty sequence \( \emptyset \).

If instead \( \inf \vec{r} = 0 \), the construction used in this paper is slightly different from the ones used in~\cite[Theorem 4.4]{Gao:2003qw} and~\cite[Proposition 4.2]{Louveau:2005cq}. 
In fact, in this case \( U^{\vec{r}}_G \) may be identified with a space on \( G \uplus \body{ G } \), where 
\[ 
\body{ G } \coloneqq \setofLR{b \in \pre{\omega}{\kappa}}{ b ( 0 ) = r \wedge \FORALL{n \in \omega} \bigl ( b ( n ) \mathrel{G} b ( n + 1 ) \bigr ) } 
\] 
is the set of all infinite \( \omega \)-branches through \( G \) starting from its root \( r \), while in ~\cite{Gao:2003qw} and~\cite{Louveau:2005cq} the authors (essentially) considered only \( G \in \RCT_\kappa \) without terminal vertices and associated to each of them only the subspace \( [G] \) of \( U^{\vec{r}}_G \). 
Although the two constructions are very close to each other (and in fact essentially equivalent), for the reader's convenience we give a sketch of the proof of the desired result based only on our new construction.

Let \( G , G' \in \RCT_\kappa \). We first deal with isomorphism and isometry, giving just the main ideas and referring the reader to~\cite[Section 3]{mottoros:2016bu} for more details.
By construction, any isomorphism (respectively, embedding) between \( G \) and \( G' \) naturally extends to an isometry (respectively, isometric embedding) between \( U^{\vec{r}}_G = G \uplus \body{ G } \) and \( U^{ \vec{r} }_{G'} = G' \uplus \body{ G' } \). 
For the other direction, first notice that the (adaptation to the uncountable context of the) proof of~\cite[Theorem 5.2]{Camerlo:2015ar} shows that if there exists an isometry (respectively, an isometric embedding) \( \psi \) between \( (\kappa, d^{\vec{r}}_G) \) and \( ( \kappa , d^{\vec{r}}_{G'} ) \), then \( G \) is isomorphic to (respectively, embeds into) \( G' \). 
It is thus enough to show that if there is an isometry (respectively, an isometric embedding) \( \varphi \colon U^{\vec{r}}_G \to U^{\vec{r}}_{G'} \), then there is also an isometry (respectively, an isometric embedding) \( \psi \colon (\kappa, d^{\vec{r}}_G) \to (\kappa, d^{\vec{r}}_{G'}) \).
Let \( \varphi \colon U^{\vec{r}}_G \to U^{\vec{r}}_{G'} \) be a metric-preserving map, and assume first that it is surjective, and therefore hence an isometry. 
Since by construction a point is isolated in \( U^{\vec{r}}_G \) if and only if it belongs to \( G \) (and the same is true when replacing \( G \) with \( G' \)), then \( \varphi(G) = G' \): thus \( \psi \coloneqq \varphi \restriction G \) is an isometry between \( (\kappa, d^{\vec{r}}_G) \) and \( (\kappa,d^{\vec{r}}_{G'}) \) and we are done. 

Assume now that \( \varphi \) is just an isometric embedding (i.e.\ not necessarily surjective). 
Then for some \( \alpha \in G \) we may have \( \varphi(\alpha) = b_\alpha \in \body{ G' } \). 
However, since \( \alpha \) is isolated in \( U^{\vec{r}}_G \), then \( b_\alpha \) must be isolated in the range of \( \varphi \): this implies that there is \( n_\alpha \in \omega \) such that \( \ran \varphi \) does not contain any \( \gamma \in G' \) with \( b_\alpha(n_\alpha) \preceq_{G'} \gamma \) nor any \( b' \in \body{ G' } \) such that \( b_\alpha ( n_\alpha ) = b' ( n_\alpha ) \).
By our choice of \( n_\alpha \), it easily follows that for every \( \beta \in G \) distinct from \( \alpha \) we have \( d^{\vec{r}}_{G'} ( b_\alpha , \varphi ( \beta ) ) = d^{\vec{r}}_{G'} ( b_\alpha ( n_\alpha ) , \varphi ( \beta ) ) \). 
Thus the map \( \psi \colon G \to G' \) defined by setting for every \( \alpha \in G \)
\[ 
\psi(\alpha) \coloneqq
\begin{cases}
\varphi(\alpha) & \text{if } \varphi ( \alpha) \in G' 
\\
b_\alpha(n_\alpha) & \text{if } \varphi(\alpha) \in \body{ G' }
\end{cases}
 \] 
is a well-defined isometric embedding between \( (\kappa, d^{\vec{r}}_G) \) and \( (\kappa,d^{\vec{r}}_{G'}) \), so we are done.
\end{proof}

Similarly to the discrete case, Remark~\ref{rmk:continuityultrametric} and Lemma~\ref{lem:reductiontoultrametric} yields the following lower bounds for the complexity of \( \cong^i \restriction \UUU_\kappa(A) \) and \( \sqsubseteq^i \restriction \UUU_\kappa(A) \). 
As usual, we endow both \( \RCT_\kappa \) and \( \MMM_\kappa \) with the bounded topology and \( \UUU_\kappa(A) \) with the induced relative topology. 
Moreover, to simplify the notation we denote by \( \cong^\kappa_\RCT \) and \( \embeds^\kappa_\RCT \) the isomorphism and the embeddability relation on \( \RCT_\kappa \), respectively.

\begin{theorem} \label{th:reductiontoultrametric}
Let \( \kappa \) be an infinite cardinal and \( A \) be an ill-founded subset of \( \R \) satisfying the conditions of Definition~\ref{def:subclassofultrametric}. \begin{enumerate-(a)}
\item \label{th:reductiontoultrametric-a}
\( {\cong^\kappa_{\RCT}} \leq^\kappa_\bB {\cong^i \restriction \UUU_\kappa(A)} \) and \( {\embeds^\kappa_{\RCT}} \leq^\kappa_\bB {\sqsubseteq^i \restriction \UUU_\kappa(A)} \).
\item \label{th:reductiontoultrametric-b}
The relation \( \sqsubseteq^i \restriction \UUU_\kappa(A) \) is \emph{\( \leq^\kappa_\bB \)-complete} for the class of \( \kappa \)-Souslin quasi-orders on Polish or standard Borel spaces.
\item \label{th:reductiontoultrametric-c}
(\( \AC \)) 
If \( \kappa \leq 2^{\aleph_0} \) and there is an \( \bS(\kappa) \)-code for \( \kappa \) (which is always the case for \( \kappa = \omega \), \( \kappa = \omega_1 \), and \( \kappa = 2^{\aleph_0} \)), then \( \sqsubseteq^i \restriction \UUU_\kappa(A) \) is \emph{\( \leq_{\bS(\kappa)} \)-complete} for \( \kappa \)-Souslin quasi-orders on Polish or standard Borel spaces.
\item \label{th:reductiontoultrametric-d}
(\( \AD + \DC \)) If \( \kappa \) is a Souslin cardinal, then \( \sqsubseteq^i \restriction \UUU_\kappa(A) \) is \emph{\( \leq_{\bS(\kappa)} \)-complete} for \( \kappa \)-Souslin quasi-orders on Polish or standard Borel spaces.
\item \label{th:reductiontoultrametric-e}
Let \( W \supseteq \R \) be an inner model and let \( \kappa \in \Cn^W \). 
Then \( \sqsubseteq^i \restriction \UUU_\kappa(A) \) is \emph{\( \leq_W \)-complete} for quasi-orders in \( \bS_W(\kappa) \).
\end{enumerate-(a)}
\end{theorem}

\begin{proof}
As observed in Remark~\ref{rmk:RCT}, in all the constructions and results from Sections~\ref{sec:mainconstruction}--\ref{sec:invariantlyuniversal} and~\ref{sec:definablecardinality} we could systematically replace combinatorial trees with rooted combinatorial trees. 
Then it is enough to argue as in the proof of Theorem~\ref{th:reductiontodiscrete} but replacing Remark~\ref{rmk:continuitydiscrete} and Lemma~\ref{lem:reductiontodiscrete} with Remark~\ref{rmk:continuityultrametric} and Lemma~\ref{lem:reductiontoultrametric}, respectively. 
\end{proof}

By Theorem~\ref{th:reductiontoultrametric}, in all the completeness results from Sections~\ref{sec:definablecardinality}, \ref{sec:applications}, and~\ref{subsec:Representingarbitrarypartialorders} we may systematically replace \( \embeds^\kappa_\CT \) with \( \sqsubseteq^i \restriction \UUU_\kappa \): here is a sample of the statements that one may obtain in this way. 

\begin{theorem}\label{th:ultrametricapplications} 
\begin{enumerate-(a)}
\item
(\( \AC \)) 
The relation \( \sqsubseteq^i \restriction \UUU_{\omega_1} \) is \( \leq^{\omega_1}_\bB \)-complete for \( \bSigma^1_2 \) quasi-orders on Polish or standard Borel spaces.
If moreover we assume either \( \AD^{\Ll(\R)} \) or \( \MA + \neg \CH + \EXISTS{a \in \pre{\omega}{\omega}} ( \omega_1^{\Ll [ a ] } = \omega_1 ) \), then \( \sqsubseteq^i \restriction \UUU_{\omega_1} \) is also \( \leq_{\bSigma^1_2} \)-complete for \( \bSigma^1_2 \) quasi-orders on Polish or standard Borel spaces.
\item
(\( \AC \)) 
Assume that \( x^\# \) exists for all \( x \in \pre{\omega}{\omega} \). 
Then the relation \( \sqsubseteq^i \restriction \UUU_{\omega_2} \) is \( \leq_{\bB}^{\omega_2} \)-complete for \( \bSigma^1_3 \) quasi-orders on Polish or standard Borel spaces. 
In particular, the quotient order of \( ( \mathcal{Q} , \leq_{\bB} ) \) (definably) embeds into the quotient order of \( \sqsubseteq^i \restriction \UUU_{\aleph_ 2 } \).
\item
(\( \AC \)) 
Assume \( \AD^{\Ll ( \R )} \). 
Then the relation \( \sqsubseteq^i \restriction \UUU_{\omega_{r(n)}} \) is \( \leq_{\bB}^{\omega_{r(n)}} \)-complete for \( \bSigma^1_n \) quasi-orders on Polish or standard Borel spaces, where \( r \colon \omega \to \omega \) is as in Theorem~\ref{th:embeddingsofprojectiveqos}.
\item
(\( \AD + \DC \)) 
For \( 0 \neq n \in \omega \), let \( \kappa_n \) be such that \( \bdelta^1_n = \kappa_n^+ \) (see Section~\ref{subsubsec:modelsofAD}).
Then the relation \( \sqsubseteq^i \restriction \UUU_{\kappa_n} \) is both \( \leq^{\kappa_n}_{\bB} \)-complete and \( \leq_{\bSigma^1_n} \)-complete for \( \bSigma^1_n \) quasi-orders on Polish or standard Borel spaces. 
In particular, the quotient order of \( ( \mathcal{Q} , \leq_{\bB} ) \) (definably) embeds into the quotient order of \( \sqsubseteq^i \restriction\UUU_{\aleph_ \omega } \).
\item \label{th:ultrametricapplications-h}
(\( \AC \)) Assume \( \AD^{\Ll(\R)} \). 
Then the relation \( \sqsubseteq \restriction \UUU_{\bdelta^2_1} \) is \( \leq_{\Ll(\R)} \)-complete for quasi-orders on \( \pre{\omega}{2} \) (or on arbitrary Polish or standard Borel spaces in \( \Ll(\R) \)) belonging to \( \bGamma^2_1 \coloneqq ( \bSigma^2_1 )^{\Ll ( \R )} \).
\item
(\( \AC \)) 
Let \( W \) be a transitive inner model containing \( \R \) and satisfying \( \ADR \), and let \( \kappa \coloneqq \Theta^W \). 
(The existence of such a \( W \) follows from a Woodin cardinal which is limit of Woodin cardinals, or even less.)
Then the relation \( \sqsubseteq^i \restriction \UUU_{\kappa} \) is \( \leq_W \)-complete for quasi-orders on \( \pre{\omega}{2} \) (or on arbitrary Polish or standard Borel spaces in \( W \)) belonging to \( W \). 
In fact, for every quasi-order \( R \) on \( \pre{\omega}{2} \) 
\[ 
R \in W \IFF R \leq_W {\sqsubseteq^i \restriction \UUU_\kappa} .
\]
\end{enumerate-(a)}
\end{theorem}

Moreover, combining Theorem~\ref{th:reductiontoultrametric} with the results from Section~\ref{subsec:Representingarbitrarypartialorders} we get that the quasi-order \( {\sqsubseteq^i} \restriction \UUU_\kappa \) is complex also from the combinatorial point of view.

\begin{theorem} \label{th:ultrametricombinatorial}
\begin{enumerate-(a)}
\item
Let \( \omega < \kappa \leq 2^{\aleph_0} \) and assume \( \AC_\kappa(\R) \). 
Then the relation \( \subseteq^*_\kappa \) on \( \pow(\kappa) \) of inclusion modulo bounded subsets is \( \kappa+1 \)-Borel reducible to \( {\sqsubseteq^i} \restriction \UUU_\kappa \). 
In particular, every partial order \( P \) of size \( \kappa \) can be embedded into the quotient order of \( {\sqsubseteq^i} \restriction \UUU_\kappa \). 
Further assuming \( \AC \) and \( 2^{\aleph_0} \geq \aleph_n \) (for some \( n \in \omega \)), we also get that every linear order of size \( \aleph_{n+1} \) can be embedded into the quotient order of \( {\sqsubseteq^i} \restriction \UUU_{\aleph_n} \).
\item
Assume \( \AD + \DC \) and let \( \kappa \) be a Souslin cardinal. 
Then every partial order \( P \) of size \( \kappa \) can be embedded into the quotient order of \( \sqsubseteq^i \restriction \UUU_\kappa \).
\end{enumerate-(a)}
\end{theorem}

\begin{proof}
Use again the fact that in Proposition~\ref{prop:parovicenkoACimproved} and Theorem~\ref{th:parovicenko1} we could have replaced combinatorial trees with rooted combinatorial trees, and then apply Theorem~\ref{th:reductiontoultrametric}\ref{th:reductiontoultrametric-a}.
\end{proof}

\begin{remark} \label{rmk:variants}
In Theorems~\ref{th:ultrametricapplications} and~\ref{th:ultrametricombinatorial} we may also consider just isometric embeddability between ultrametric spaces in \( \UUU_\kappa ( A ) \) for any ill-founded set of distances \( A \) (provided that \( A \) is as in Definition~\ref{def:subclassofultrametric}). 
As already observed, when \( A \) is bounded away from \( 0 \) then all spaces in \( \UUU_\kappa ( A ) \) are discrete (and hence of size \( \kappa \)). 
This implies that in Theorems~\ref{th:ultrametricapplications} and~\ref{th:ultrametricombinatorial} we may replace the collection of complete ultrametric spaces of density character \( \kappa \) with any subclass of \( \MMM_\kappa \) containing one of these \( \UUU_\kappa ( A ) \), including the following notable examples (all spaces below are intended to be complete metric and of density character \( \kappa \)):
\begin{itemize}[leftmargin=1pc]
\item
discrete (ultrametric) spaces;
\item
(ultrametric) spaces of size \( \kappa \);
\item
locally compact (ultrametric) spaces;
\item
zero-dimensional spaces.
\end{itemize}
\end{remark}

We conclude this section by observing that in Theorem~\ref{th:reductiontoultrametric} we could also systematically replace all occurrences of \( \UUU_\kappa(A) \) with the smaller \( \UUU^\star_\kappa(A) \). 
This is because one can find (for any ill-founded \( A \) satisfying the conditions of Definition~\ref{def:subclassofultrametric}) a modification \( \theta^{\star}_{U,A} \colon \RCT_\kappa \to \UUU^\star_\kappa(A) \) of the map \( \theta^{\vec{r}}_U \) from~\eqref{eq:defthetaU} such that:
\begin{enumerate-(1)}
\item \label{thetastar-a}
the map \( \theta^{\star}_{U,A} \) is continuous when both \( \RCT_\kappa \) and \( \UUU_\kappa^\star(A) \) are endowed with the same topology \( \tau_p \) or \( \tau_b \);
\item \label{thetastar-b}
the map \( \theta^{\star}_{U,A} \) simultaneously reduces isomorphism to isometry and embeddability to isometric embeddability.
\end{enumerate-(1)}

In fact, given such an \( A \subseteq \R \) fix a strictly decreasing sequence \( \vec{r} = \seqofLR{r_n}{n \in \omega} \) of points from \( A \) such that there is \( \bar{r} \in A \) with \( \bar{r} > r_0 \). 
Given \( G \in \RCT_\kappa \), consider the space \( U^{\vec{r}}_G \) coded by \( \theta^{\vec{r}}_U ( G ) = u^{\vec{r}}_G \in \UUU_\kappa ( A ) \), and notice that since \( G \) has size \( \kappa \) (in particular, it has vertices distinct from its root) each point of \( U^{\vec{r}}_G \) realizes the distance \( r_0 \). 
Then let \( U^A_G \) be the disjoint union of \( U^{\vec{r}}_G \) and \( U ( A \setminus \set{ r_0 } ) \) (where \( U ( A \setminus \set{ r_0 } ) = ( U_{ A \setminus \set{ r_0 } }, d_{A \setminus \set{ r_0 } }) \) is defined as in~\eqref{eq:canonicalultrametric}), and endow it with the metric \( d^A_G \) extending both \( d^{\vec{r}}_G \) and \( d_{A \setminus \set{ r_0 } } \) obtained by setting for every \( x \in U^{\vec{r}}_G \), \( r \in A \setminus \set{ r_0 } \), and \( \alpha < \kappa \)
\[ 
d^A_G ( x , ( r , \alpha ) ) \coloneqq \max \set { \bar{r} , r }.
 \] 
A code \( u^A_G \in \UUU^\star ( A ) \) for \( U^A_G \) can be obtained in a continuous-in-\( G \) way by e.g.\ copying the code \( u^{\vec{r}}_G \) on the odd ordinal numbers and (a code for) the space \( U ( A \setminus \set{ r_0 } ) \) on the even ordinal numbers --- we leave to the reader to carry out the details of such coding procedure. 
Consider the map \( \theta^\star_{U,A} \colon \RCT_\kappa \to \UUU^\star_\kappa(A) \subseteq \MMM_\kappa \) defined by setting \( \theta^{\star}_{U,A}(G) \coloneqq u^A_G \) for every \( G \in \RCT_\kappa \). 
By the coding procedure briefly sketched above, condition~\ref{thetastar-a} is satisfied. 
To show that also~\ref{thetastar-b} is satisfied, first notice that by Lemma~\ref{lem:reductiontoultrametric} we only need to check that for all \( G,G' \in \RCT_\kappa \)
\[ 
{\theta^{\vec{r}}_U ( G ) \cong^i \theta^{\vec{r}}_U ( G' ) } \IFF { \theta^{\star}_{ U , A } ( G ) \cong^i \theta^{\star}_{ U , A } ( G' ) } \quad \text{ and } \quad {\theta^{\vec{r}}_U ( G ) \sqsubseteq^i \theta^{\vec{r}}_U ( G ' ) } \IFF { \theta^{\star}_{ U , A } ( G ) \sqsubseteq^i \theta^{\star}_{ U , A } ( G' ) }.
 \] 
The forward direction is obvious, so assume that \( \theta^{\star}_{U,A} ( G ) \cong^i \theta^{\star}_{ U , A } ( G' ) \) (respectively, \( \theta^{\star}_{ U , A } ( G ) \sqsubseteq^i \theta^{\star}_{U , A } ( G' ) \)) and let \( \varphi \) be an isometry (respectively, an isometric embedding) between the spaces \( U^A_G \) and 
\( U^A_{G'} \) coded by \( \theta^{\star}_{ U , A } ( G ) \) and \( \theta^{\star}_{ U , A } ( G' ) \). 
Then since by construction the points in \( U^{\vec{r}}_G \subseteq U^A_G \) are the unique realizing the distance \( r_0 \) in \( U^A_G \) (and the same is true when replacing \( G \) with \( G' \)), then \( \varphi \restriction U^{\vec{r}}_G \) is an isometry (respectively, an isometric embedding) between \( U^{\vec{r}}_G \) and \( U^{\vec{r}}_{G'} \), so that \( \theta^{\vec{r}}_U ( G ) \cong^i \theta^{\vec{r}}_U ( G' ) \) (respectively, \( \theta^{\vec{r}}_U ( G ) \sqsubseteq^i \theta^{\vec{r}}_U ( G ' ) \)).

\subsection{Linear isometry and linear isometric embeddability between Banach spaces of density \( \kappa \)} \label{subsec:linearisomemb}

Let \( \boldsymbol{c_0} = (c_0, \|\cdot\|_\infty) \) be the separable (real) Banach space of vanishing \( \omega \)-sequences of reals endowed with the usual pointwise 
operations and the sup norm \( \| \cdot \|_\infty \).
Building on previous work by Louveau and Rosendal~\cite{Louveau:2005cq}, in~\cite[Section 5.6]{Camerlo:2012kx} it was defined a Borel map simultaneously reducing isomorphism and embeddability between countable graphs to, respectively, linear isometry and linear isometric embeddability between separable Banach spaces isomorphic to \( \boldsymbol{c_0} \). 
We are now going to adapt such construction to the uncountable context in order to study the complexity of the relations of linear isometry \( \cong^{li} \) and linear isometric embeddability \( \sqsubseteq^{li} \) between non-separable Banach spaces. 
To this aim we first have to define the right analogue of the space \( \boldsymbol{c_0} \).

\begin{definition} \label{def:c_0^kappa}
Given an infinite cardinal \( \kappa \), let \( c_0^\kappa \subseteq \pre{\kappa}{\R} \) be the collection of all \( \kappa \)-sequences \( x = \seqof{ x_\alpha }{ \alpha < \kappa} \) of reals such that for all \( \varepsilon \in \R^+ \) we have \( \card{ x_\alpha } < \varepsilon \) for all but finitely many \( \alpha < \kappa \).

The Banach space \( \boldsymbol{c^\kappa_0} = ( c^\kappa_0 , \| \cdot \|_\infty ) \) is then obtained by endowing \( c^\kappa_0 \) with the usual pointwise operations and the sup norm \( \| \cdot \|_\infty \).
\end{definition}

In particular, \( \boldsymbol{c^\omega_0} = \boldsymbol{c_0} \).
A basis \( \seqofLR{e_\alpha}{\alpha < \kappa} \) for \( \boldsymbol{c^\kappa_0} \) is obtained by letting \( e_\alpha \) be the \( \kappa \)-sequence with value \( 1 \) on coordinate \( \alpha \) and \( 0 \) elsewhere; then each element \( x \in c^\kappa_0 \) may be uniquely written as \( \sum_\alpha x_\alpha e_\alpha \).
Notice that \( \boldsymbol{c^\kappa_0} \) has density \( \kappa \) (a dense subset is given by the collection of all \( x \in c^\kappa_0 \cap \pre{\kappa}{\Q} \) having only finitely many non-null coordinates), and that each \( x \in c^\kappa_0 \) has countable support, that is it has at most countably many non-null coordinates. 

Given any graph \( G \) on \( \kappa \), let \( \boldsymbol{X_G} = (c^\kappa_0, \| \cdot \|_G) \) be the (real) Banach space on \( c^\kappa_0 \) equipped with the pointwise operations and the norm defined by 
\[ 
\left \| \sum\nolimits_\alpha x_\alpha e_\alpha \right \|_G \coloneqq \sup \setofLR{ \card{ x_i } + \textstyle\frac{ \card{ x_j } }{3 - \chi_G ( i , j ) } }{ i \neq j \in \kappa },
 \] 
where \( \chi_G \colon \kappa \times \kappa \to \set{ 0 , 1 } \) is the characteristic function of the graph relation of \( G \). 
It is easy to verify that \( \| \cdot \|_G \) is equivalent to \( \| \cdot \|_\infty \), as \( \| \sum_\alpha x_\alpha e_\alpha \|_\infty \leq \| \sum_\alpha x_\alpha e_\alpha \|_G \leq \frac{3}{2} \| \sum_\alpha x_\alpha e_\alpha \|_\infty \). 
Moreover, \( \boldsymbol{X_G} \) has density \( \kappa \) and thus it can be coded as in~\eqref{eq:x_B} by an element \( x_{\boldsymbol{X_G}} = ( x_{ \boldsymbol{ X_G}}^+ , x_{ \boldsymbol{ X_G} }^\Q , x_{\boldsymbol{X_G}}^{\| \cdot \|}) \) of the standard Borel \( \kappa \)-space \( \BBB_\kappa \) introduced in Section~\ref{subsubsec:spacesofBanachspaces}.

\begin{remark} \label{rmk:continuityBanach}
The map 
\begin{equation} \label{eq:defthetaB}
\theta_B \colon \Mod^\kappa_\GR \to \BBB_\kappa, \qquad G \mapsto x_{\boldsymbol{X_G}}
\end{equation}
is continuous when both \( \Mod^\kappa_\GR \) and \( \BBB_\kappa \) are endowed with the same topology \( \tau_p \) or \( \tau_b \).
\end{remark}

The proof of the next lemma is identical%
\footnote{It is enough to check that all the properties of the norm \( \| \cdot \|_G \) used in the original proof are maintained when passing to the uncountable context, including e.g.\ the fact that the sup in the definition of \( \|\cdot\|_G \) is attained, or the fact that when \( \| \sum_\alpha x_\alpha e_\alpha\|_G \geq \varepsilon \) for some \( \varepsilon \in \R^+ \), then there are only finitely many coordinates \( i,j \in \kappa \) for which \( | x_i | + \frac{ | x_j | }{ 3 - \chi_G ( i , j ) } > \varepsilon \).}
to the one of~\cite[Lemma 5.20]{Camerlo:2012kx}, so we omit it here.

\begin{lemma} \label{lem:reductiontoBanach}
The map \( \theta_B \) from~\eqref{eq:defthetaB} simultaneously reduces \( \cong \) to \( \cong^{li} \) and \( \embeds \) to \( \sqsubseteq^{li} \).
\end{lemma}

Arguing as in the case of metric spaces one can straightforwardly check that Remark~\ref{rmk:continuityBanach} and Lemma~\ref{lem:reductiontoBanach} yield the following lower bounds for the complexity of \( \cong^{li} \restriction \BBB_\kappa \) and \( \sqsubseteq^{li} \restriction \BBB_\kappa \). 
(As usual, in what follows both \( \Mod^\kappa_\GR \) and \( \BBB_\kappa \) are endowed with the bounded topology.)

\begin{theorem} \label{th:reductiontoBanach}
Let \( \kappa \) be an infinite cardinal. 
\begin{enumerate-(a)}
\item \label{th:reductiontoBanach-a}
\( {\cong^\kappa_{\GR}} \leq^\kappa_\bB {\cong^{li} \restriction \BBB_\kappa} \) and \( {\embeds^\kappa_{\GR}} \leq^\kappa_\bB {\sqsubseteq^{li} \restriction \BBB_\kappa} \).
\item \label{th:reductiontoBanach-b}
The relation \( \sqsubseteq^{li} \restriction \BBB_\kappa \) is \emph{\( \leq^\kappa_\bB \)-complete} for the class of \( \kappa \)-Souslin quasi-orders on Polish or standard Borel spaces.
\item \label{th:reductiontoBanach-c}
(\( \AC \)) 
If \( \kappa \leq 2^{\aleph_0} \) and there is an \( \bS(\kappa) \)-code for \( \kappa \) (which is always the case for \( \kappa = \omega \), \( \kappa = \omega_1 \), and \( \kappa = 2^{\aleph_0} \)), then \( \sqsubseteq^{li} \restriction \BBB_\kappa \) is \emph{\( \leq_{\bS ( \kappa )} \)-complete} for \( \kappa \)-Souslin quasi-orders on Polish or standard Borel spaces.
\item \label{th:reductiontoBanach-d}
(\( \AD + \DC \)) 
If \( \kappa \) is a Souslin cardinal, then \( \sqsubseteq^{li} \restriction \BBB_\kappa \) is \emph{\( \leq_{\bS ( \kappa)} \)-complete} for \( \kappa \)-Souslin quasi-orders on Polish or standard Borel spaces.
\item \label{th:reductiontoBanach-e}
Let \( W \supseteq \R \) be an inner model and let \( \kappa \in \Cn^W \). 
Then \( \sqsubseteq^{li} \restriction \BBB_\kappa \) is \emph{\( \leq_W \)-complete} for quasi-orders in \( \bS_W ( \kappa) \).
\end{enumerate-(a)}
\end{theorem}

Theorem~\ref{th:reductiontoBanach} shows in particular that in all the completeness results from Sections~\ref{sec:definablecardinality}, \ref{sec:applications} and ~\ref{subsec:Representingarbitrarypartialorders} we may further replace \( \embeds^\kappa_\CT \) with \( \sqsubseteq^{li} \restriction \BBB_\kappa \). 
To see a sample of the statements that one may obtain in this way just systematically replace \( \sqsubseteq^i \restriction \UUU_\kappa \) with \( \sqsubseteq^{li} \restriction \BBB_\kappa \) in Theorems~\ref{th:ultrametricapplications} and~\ref{th:ultrametricombinatorial}.

\subsection{* Further results on the classification of nonseparable metric and Banach spaces}
We collect some facts that follow easily by combining theorems from the literature with the constructions developed insofar.
Although some of the results are not strictly related to the main topic of this paper, they do fit naturally into our framework and they do not require much extra work.

First we notice that many of the relations considered in the previous section, including the isometry relation \( \cong^i \), the isometric embeddability relation \( \sqsubseteq^i \), and their analogue for Banach spaces \( \cong^{li} \) and \( \sqsubseteq^{li} \) are consistently as complex as possible, even when restricted to some specific subclasses. 

\begin{theorem} \label{thm:completenessappendix}
\begin{enumerate-(a)}
\item \label{thm:completenessappendix-a}
Assume \( \Vv = \Ll \), and let \( \kappa = \lambda^+ \) be such that \( \lambda^\omega = \lambda \) (equivalently: \( \lambda \) is either a successor cardinal or a limit cardinal of uncountable cofinality). 
Then the following relations are \emph{\( \leq^\kappa_{\bB} \)-complete} for the collection of all \( \kappa \)-analytic equivalence relations on standard Borel \( \kappa \)-spaces:
\begin{itemize}[leftmargin=1pc]
\item
\( \cong^i \restriction \DDD_\kappa \);
\item
\( \cong^i \restriction \UUU_\kappa(A) \), where \( A \subseteq \R \) satisfies the conditions of Definition~\ref{def:subclassofultrametric} and contains a strictly decreasing chain of length \( 2 \cdot \omega +1 \) (with respect to the usual ordering of \( \R \));
\item
\( \cong^{li} \restriction \BBB_\kappa \).
\end{itemize}
\item \label{thm:completenessappendix-b}
Assume \( \AC \) and let \( \kappa \) be a weakly compact%
\footnote{\cite{Mildenberger:2012ke} shows that the same is true for any cardinal \( \kappa \) satisfying the equality \( \kappa^{< \kappa} = \kappa \).}
cardinal. 
Then the following relations are \emph{\( \leq^\kappa_{\bB} \)-complete} for the collection of all \( \kappa \)-analytic quasi-orders on standard Borel \( \kappa \)-spaces:
\begin{itemize}[leftmargin=1pc]
\item
\( \sqsubseteq^i \restriction \DDD_\kappa \);
\item
\( \sqsubseteq^{li} \restriction \BBB_\kappa \).
\end{itemize}
\end{enumerate-(a)}
\end{theorem}

\begin{proof}
\ref{thm:completenessappendix-a}
In~\cite[Corollary 2.15]{Hyttinen:2012fj} it is proved that under our assumptions there is a first-order theory \( T \) such that the isomorphism relation on the collection \( \Mod^\kappa_T \) of models of size \( \kappa \) of \( T \) is \( \leq^\kappa_\bB \)-complete for \( \kappa \)-analytic equivalence relations on \( \pre{\kappa}{2} \).
This can be extended to \( \kappa \)-analytic equivalence relations on arbitrary standard Borel \( \kappa \)-spaces by adapting the proof of~\cite[Lemma 6.8]{Motto-Ros:2011qc} to our context, using the fact that any standard Borel \( \kappa \)-space is \( \kappa+ 1 \)-Borel isomorphic to some \( B \in \bB_{\kappa+1}(\pre{\kappa}{2}, \tau_b) \). 
 By~\cite[Theorem 5.5.1]{Hodges:1993mi}, \( T \) can be interpreted in the theory of graphs: more precisely, there is a \( \tau_b \)-continuous map \( \Mod^\kappa_T \to \Mod^\kappa_\GR \) which simultaneously reduces \( \cong \restriction \Mod^\kappa_T \) to \( \cong^\kappa_\GR \) and \( \embeds \restriction \Mod^\kappa_T \) to \( \embeds^\kappa_\GR \). 
In particular, this shows that under our assumptions \( \cong^\kappa_\GR \) is \( \leq^\kappa_\bB \)-complete for \( \kappa \)-analytic equivalence relations on \( \pre{\kappa}{2} \) as well. 
Thus the result concerning discrete metric spaces follows from Theorem~\ref{th:reductiontodiscrete}\ref{th:reductiontodiscrete-a}, while the one for Banach spaces follows from Theorem~\ref{th:reductiontoBanach}\ref{th:reductiontoBanach-a}. 
The case of ultrametric spaces follows the same reasoning lines but is slightly more complex: a full proof can be found in~\cite[Sections 3 and 4]{mottoros:2016bu}.

\ref{thm:completenessappendix-b}
The argument is similar to that of part~\ref{thm:completenessappendix-a}. 
In~\cite[Corollary 9.5]{Motto-Ros:2011qc} it is proved that under our assumptions the embeddability relation on \emph{generalized trees}%
\footnote{We call \markdef{generalized tree} any partial order \( T = (T, \leq_T) \) that is a tree in the model theoretic sense, that is: for every \( t \in T \), the set \( \setofLR{t' \in T}{t' \leq_T t} \) is linearly ordered by \( \leq_T \).}
of size \( \kappa \) is \( \leq^\kappa_\bB \)-complete for \( \kappa \)-analytic quasi-orders on standard Borel \( \kappa \)-spaces. 
By the above argument, this implies that the same is true for \( \embeds^\kappa_\GR \). 
Thus the result follows again from Theorem~\ref{th:reductiontodiscrete}\ref{th:reductiontodiscrete-a} and Theorem~\ref{th:reductiontoBanach}\ref{th:reductiontoBanach-a}.
\end{proof}

\begin{remark}
Since in Theorem~\ref{thm:completenessappendix}\ref{thm:completenessappendix-a} the set \( A \) can be chosen to be bounded away from \( 0 \), Remark~\ref{rmk:variants} applies here as well: as a consequence, the isometry relation on any of the classes of metric spaces mentioned in that remark is consistently as complex as possible. 
Moreover, arguing as at the end of Section~\ref{subsec:isomemb} one sees that we may as well replace \( \UUU_\kappa ( A ) \) with \( \UUU_\kappa^\star ( A ) \) in Theorem~\ref{thm:completenessappendix}\ref{thm:completenessappendix-a}.
\end{remark}

Whether \( \sqsubseteq^i \restriction \UUU_\kappa \) may consistently have maximal complexity (i.e.\ whether it can be added to Theorem~\ref{thm:completenessappendix}\ref{thm:completenessappendix-b}) seems to be open: the main problem is that we do not know e.g.\ whether \( \embeds^\kappa_\RCT \) can be \( \leq^\kappa_\bB \)-complete for \( \kappa \)-analytic quasi-orders on \( \pre{\kappa}{2} \) (equivalently, whether it has the same complexity as \( \embeds^\kappa_\GR \)).
However, we can at least observe that \( \cong^\kappa_\GR \) and \( \embeds^\kappa_\GR \) are always upper bounds for the complexity of \( \cong^i \restriction \UUU_\kappa \) and \( \sqsubseteq^i \restriction \UUU_\kappa \), respectively. 

\begin{theorem} \label{th:upperboundsultrametric}
Let \( \kappa \) be any infinite cardinal. 
Then \( {\cong^i \restriction \UUU_\kappa} \leq^\kappa_\bB {\cong^\kappa_\GR} \) and \( {\sqsubseteq^i \restriction \UUU_\kappa} \leq^\kappa_\bB {\embeds^\kappa_\GR} \).
\end{theorem}

\begin{proof} 
Argue as in the proof of~\cite[Theorem 4.4]{Gao:2003qw} and use the fact that \( {\cong \restriction \Mod^\kappa_{\hat{\LL}}} \leq^\kappa_\bB {\cong^\kappa_\GR} \) and \( {\embeds \restriction \Mod^\kappa_{\hat{\LL}}} \leq^\kappa_\bB {\embeds^\kappa_\GR} \) for any countable language \( \hat{\LL} \) by~\cite[Theorem 5.5.1]{Hodges:1993mi} (see also the proof of Theorem~\ref{thm:completenessappendix}\ref{thm:completenessappendix-a}).
\end{proof}

In contrast, for discrete metric spaces we can complement Theorem~\ref{th:reductiontodiscrete}\ref{th:reductiontodiscrete-a} and determine the exact complexity with respect to \( \leq^\kappa_\bB \) of both \( \cong^i \restriction \DDD_\kappa \) and \( \sqsubseteq \restriction \DDD_\kappa \). 

\begin{theorem} \label{th:upperboundsdiscrete}
Let \( \kappa \) be any infinite cardinal. 
Then \( {\cong^i \restriction \DDD_\kappa} \sim^\kappa_\bB {\cong^\kappa_\GR} \) and \( {\sqsubseteq^i \restriction \DDD_\kappa} \sim^\kappa_\bB {\embeds^\kappa_\GR} \).
\end{theorem}

\begin{proof}
By Theorem~\ref{th:reductiontodiscrete}\ref{th:reductiontodiscrete-a}, we only need to show that \( {\cong^i \restriction \DDD_\kappa} \leq^\kappa_\bB {\cong^\kappa_\GR} \) and \( {\sqsubseteq^i \restriction \DDD_\kappa} \leq^\kappa_\bB {\embeds^\kappa_\GR} \). 
Consider the countable language \( \hat{\LL} \coloneqq \setofLR{P_q}{q \in \QQ^+} \), where each \( P_q \) is a binary relational symbol. 
To each \( x \in \DDD_\kappa \) associate the \( \hat{\LL} \)-structure \( \mathcal{A}_x \) on \( \kappa \) defined by setting for each \( q \in \QQ^+ \) and \( \alpha, \beta < \kappa \)
\[ 
P^{\mathcal{A}_x}_q ( \alpha , \beta ) \IFF x ( \alpha , \beta , q ) = 1.
 \] 
Then the map \( x \mapsto \mathcal{A}_x \) is continuous when both \( \DDD_\kappa \) and \( \Mod^\kappa_{\hat{\LL}} \) are endowed with the bounded topology, and moreover for every \( x,y \in \DDD_\kappa \)
\[ 
{x \cong^i y} \IFF {\mathcal{A}_x \cong \mathcal{A}_y} \quad \text{and} \quad {x \sqsubseteq^i y} \IFF {\mathcal{A}_x \embeds \mathcal{A}_y}.
 \] 
Indeed, using the notation from Section~\ref{subsubsec:spacesofmetricspaces}, a map \( f \colon \kappa \to \kappa \) is an isometry (respectively, an isometric embedding) between \( M_x = (\kappa,d_x) \) and \( M_y = (\kappa,d_y) \) if and only if it is an isomorphism (respectively, an embedding) between \( \mathcal{A}_x \) and \( \mathcal{A}_y \).
This shows that \( {\cong^i \restriction \DDD_\kappa} \leq^\kappa_\bB {\cong \restriction \Mod^\kappa_{\hat{\LL}}} \) and \( {\sqsubseteq^i \restriction \DDD_\kappa} \leq^\kappa_\bB {\embeds \restriction \Mod^\kappa_{\hat{\LL}}} \). 
Since as discussed in the proof of Theorem~\ref{thm:completenessappendix}\ref{thm:completenessappendix-a} we have \( {\cong \restriction \Mod^\kappa_{\hat{\LL}}} \leq^\kappa_\bB {\cong^\kappa_\GR} \) and \( {\embeds \restriction \Mod^\kappa_{\hat{\LL}}} \leq^\kappa_\bB {\cong^\kappa_\GR} \) by~\cite[Theorem 5.5.1]{Hodges:1993mi}, we are done.
\end{proof}

Since Theorems~\ref{th:upperboundsultrametric} and~\ref{th:upperboundsdiscrete} are proved in \( \ZF \), it follows that their conclusions are true in \( \Ll ( \R) \). 
This shows that in models of \( \AD^{\Ll ( \R)} \) the relation of isometric bi-embeddability between ultrametric or discrete complete metric spaces of uncountable density character is way more complex than the isometry relation on the same class.
In fact, by Theorem~\ref{th:ultrametricapplications}\ref{th:ultrametricapplications-h} every equivalence relation in \( \bGamma^2_1 \) (a quite large boldface pointclass in \( \Vv \) which includes e.g.\ all projective levels) is \( \Ll ( \R ) \)-reducible to the isometric bi-embeddability relation relation on \( \UUU_\kappa \cap \DDD_\kappa \) (for a suitable cardinal \( \kappa \)), while by Example~\ref{xmp:AD} (see~\cite[Theorem 9.18]{Hjorth:2000zr}) and Theorems~\ref{th:upperboundsdiscrete} and~\ref{th:upperboundsultrametric} there are \( \bSigma^1_1 \) equivalence relations \( E \) such that \( E \nleq_{\Ll ( \R )} {\cong^i \restriction \DDD_\kappa } \) and \( E \nleq_{\Ll ( \R )} {\cong^i \restriction \UUU_\kappa } \) for \emph{every} \( \kappa \in ( \Cn )^{ \Ll ( \R ) } \), and therefore for \emph{every} \( \kappa \in \Cn \).

Combining this observation with Theorem~\ref{thm:completenessappendix} we further get the following interesting independence result.

\begin{corollary} \label{cor:appindependence}
It is independent of \( \ZF + \DC \) whether for \( \kappa = \omega_2 \) (or for any \( \kappa \) which is a successor of some \( \lambda \) satisfying \( \lambda^\omega = \lambda \)) the relation of isometric bi-embeddability between ultrametric (respectively, discrete) complete metric spaces of density character \( \kappa \) is \( \leq^\kappa_{\bB} \)-reducible to the isometry relation on the same class of spaces.
\end{corollary}

\begin{proof}
Since the relation of isometric bi-embeddability is a \( \kappa \)-analytic equivalence relation on the standard Borel \( \kappa \)-space \( \MMM_\kappa \), under 
\( \Vv = \Ll \) it is \( \leq^\kappa_\bB \)-reducible to \( \cong^i \restriction ( \UUU_\kappa \cap \DDD_\kappa) \) by Theorem~\ref{thm:completenessappendix}\ref{thm:completenessappendix-a} (and the ensuing remark). 
On the other hand, by the observation preceding this corollary we get that under \( \AD + {\Vv = \Ll ( \R ) } \) there are \( \bSigma^1_1 \) equivalence relations on \( \R \) which are not reducible (even without definability conditions on the reductions that may be used) to \( \cong^i \restriction \DDD_\kappa \) or to \( \cong^i \restriction \UUU_\kappa \), while all such relations are \( \leq^\kappa_\bB \)-reducible to the isometric bi-embeddability relation on \( \UUU_\kappa \cap \DDD_\kappa \) by Theorem~\ref{th:reductiontoultrametric}\ref{th:reductiontoultrametric-b} (and Remark~\ref{rmk:variants}).
Thus in this case the isometric bi-embeddability relation is not \( \leq^\kappa_\bB \)-reducible to \( \cong^i \restriction \DDD_\kappa \) or \( \cong^i \restriction \UUU_\kappa \).
\end{proof}

\backmatter
\chapter*{Indexes}
Here you will find two indexes, one for the \textbf{concepts} and one for the \textbf{symbols}.
The index of symbols is essentially divided in two parts: in the first part we listed in the order they appear in the text all  those symbols (such as \( \equals \)) that cannot easily be placed in  alphabetical order, and in the second part we list  lexicographically  all the other symbols (such as \(  \Sigma _T \)).
\printindex[concepts]
\printindex[symbols]

\bibliographystyle{amsalpha}
\bibliography{BibliographyAM}

\newcommand{\etalchar}[1]{$^{#1}$}
\providecommand{\bysame}{\leavevmode\hbox to3em{\hrulefill}\thinspace}
\providecommand{\MR}{\relax\ifhmode\unskip\space\fi MR }
\providecommand{\MRhref}[2]{%
  \href{http://www.ams.org/mathscinet-getitem?mr=#1}{#2}
}
\providecommand{\href}[2]{#2}
\begin{thebibliography}{CMMR18b}

\bibitem[AC18]{Andretta:2015}
A.~Andretta and R.~Camerlo, \emph{Analytic sets of reals and the density
  function in the {C}antor space}, European Journal of Mathematics (2018).

\bibitem[AK00]{Adams:2000gf}
S.~Adams and A.~S. Kechris, \emph{Linear algebraic groups and countable {B}orel
  equivalence relations}, J. Amer. Math. Soc. \textbf{13} (2000), no.~4,
  909--943 (electronic).

\bibitem[Bar69]{Barwise:1969fk}
J.~Barwise, \emph{Infinitary logic and admissible sets}, Journal of Symbolic
  Logic \textbf{34} (1969), no.~2, 226--252.

\bibitem[BK96]{Becker:1996uq}
H.~Becker and A.~S. Kechris, \emph{The descriptive set theory of {P}olish group
  actions}, London Mathematical Society Lecture Note Series, no. 232, Cambridge
  University Press, Cambridge, 1996.

\bibitem[CGK01]{Clemens:2001pu}
J.~Clemens, S.~Gao, and A.~S. Kechris, \emph{Polish metric spaces: their
  classification and isometry groups}, Bull. Symbolic Logic \textbf{7} (2001),
  no.~3, 361--375.

\bibitem[CGP98]{Chaber:1998cg}
J.~Chaber, G.~Gruenhage, and R.~Pol, \emph{On {S}ouslin sets and embeddings in
  integer-valued function spaces on {$\omega\sb 1$}}, Topology Appl.
  \textbf{82} (1998), no.~1-3, 71--104, Special volume in memory of Kiiti
  Morita.

\bibitem[CK13]{Chatyrko:2013ch}
V.~A. Chatyrko and A.~Karassev, \emph{The (dis)connectedness of products in the
  box topology}, Questions Answers Gen. Topology \textbf{31} (2013), no.~1,
  11--21.

\bibitem[Cle07]{Clemens:2007yu}
J.~Clemens, \emph{Classifying {B}orel automorphisms}, J. Symbolic Logic
  \textbf{72} (2007), no.~4, 1081--1092.

\bibitem[Cle12]{Clemens:2012hg}
\bysame, \emph{Isometry of {P}olish metric spaces}, Ann. Pure Appl. Logic
  \textbf{163} (2012), no.~9, 1196--1209.

\bibitem[CM97]{Chigogidze:1997cm}
A.~Chigogidze and J.~Martin, \emph{Fixed point sets of autohomeomorphisms of
  uncountable products}, Topology Appl. \textbf{80} (1997), no.~1-2, 63--71.

\bibitem[CMMR13]{Camerlo:2012kx}
R.~Camerlo, A.~Marcone, and L.~Motto~Ros, \emph{Invariantly universal analytic
  quasi-orders}, Transactions of the American Mathematical Society \textbf{365}
  (2013), no.~4, 1901--1931.

\bibitem[CMMR18a]{Mildenberger:2012ke}
F.~Calderoni, H.~Mildenberger, and L.~Motto~Ros, \emph{Uncountable structures
  are not classifiable up to bi-embeddability}, preprint submitted (2018).

\bibitem[CMMR18b]{Camerlo:2015ar}
R.~Camerlo, A.~Marcone, and L.~Motto~Ros, \emph{On isometry and isometric
  embeddability between ultrametric {P}olish spaces}, Adv. Math. \textbf{329}
  (2018), 1231--1284. \MR{3783437}

\bibitem[CMMR19]{Camerlo:2019ke}
\bysame, \emph{Polish metric spaces with fixed distance set}, preprint
  submitted (2019).

\bibitem[Coh63]{Cohen:1963ys}
P.~Cohen, \emph{The independence of the continuum hypothesis}, Proceedings of
  the National Academy of Sciences of the United States of America \textbf{50}
  (1963), 1143--1148.

\bibitem[CT19]{Calderoni:2018aa}
F.~Calderoni and S.~Thomas, \emph{The bi-embeddability relation for countable
  abelian groups}, Trans. Amer. Math. Soc. \textbf{371} (2019), no.~3,
  2237--2254. \MR{3894051}

\bibitem[DV11]{Dzamonja:2011ly}
M.~D{\v{z}}amonja and J.~V{\"a}{\"a}n{\"a}nen, \emph{Chain models, trees of
  singular cardinality and dynamic {EF}-games}, J. Math. Log. \textbf{11}
  (2011), no.~1, 61--85.

\bibitem[FFH{\etalchar{+}}12]{Fokina:2011dq}
E.~Fokina, S.~D. Friedman, V.~Harizanov, J.~F. Knight, C.~McCoy, and
  A.~Montalb\'an, \emph{Isomorphism relations on computable structures},
  Journal of Symbolic Logic \textbf{77} (2012), no.~1, 122--132.

\bibitem[FHK14]{Friedman:2011nx}
S.~D. Friedman, T.~Hyttinen, and V.~Kulikov, \emph{Generalized descriptive set
  theory and classification theory}, Mem. Amer. Math. Soc. \textbf{230} (2014),
  no.~1081, vi+80. \MR{3235820}

\bibitem[FK08]{Friedman:2008oq}
S.~D. Friedman and V.~Kanovei, \emph{Some natural equivalence relations in the
  {S}olovay model}, Abh.\ Math.\ Semin.\ Univ.\ Hambg. \textbf{78} (2008),
  no.~1.

\bibitem[FLR09]{Ferenczi:2009fk}
V.~Ferenczi, A.~Louveau, and C.~Rosendal, \emph{The complexity of classifying
  separable {B}anach spaces up to isomorphism}, J. Lond. Math. Soc. (2)
  \textbf{79} (2009), no.~2, 323--345. \MR{2496517 (2010c:46017)}

\bibitem[FMR11]{Friedman:2011cr}
S.~D. Friedman and L.~Motto~Ros, \emph{Analytic equivalence relations and
  bi-embeddability}, Journal of Symbolic Logic \textbf{76} (2011), no.~1,
  243--266.

\bibitem[FS89]{Friedman:1989bs}
H.~Friedman and L.~Stanley, \emph{A {B}orel reducibility theory for classes of
  countable structures}, J. Symbolic Logic \textbf{54} (1989), no.~3, 894--914.

\bibitem[Gao09]{Gao:2009fv}
S.~Gao, \emph{Invariant {D}escriptive {S}et {T}heory}, Monographs and Textbooks
  in Pure and Applied Mathematics, no. 293, CRC Press, Taylor\&Francis Group,
  2009.

\bibitem[GK03]{Gao:2003qw}
S.~Gao and A.~S. Kechris, \emph{On the classification of {P}olish metric spaces
  up to isometry}, Mem. Amer. Math. Soc. \textbf{161} (2003), no.~766, viii+78.

\bibitem[G{\"o}d38]{Godel:1938dz}
K.~G{\"o}del, \emph{The consistency of the axiom of choice and the generalized
  continuum hypothesis}, Proceedings of the National Academy of Sciences of the
  United Staes of America \textbf{24} (1938), 556--557.

\bibitem[Hal96]{Halko:1996fu}
A.~Halko, \emph{Negligible subsets of the generalized {B}aire space {$\omega\sp
  {\omega\sb 1}\sb 1$}}, Ann. Acad. Sci. Fenn. Math. Diss. (1996), no.~107, 38.

\bibitem[Har77]{Harrington:1977aa}
L.~Harrington, \emph{Long projective wellorderings}, Annals of Mathematical
  Logic \textbf{12} (1977), 1--24.

\bibitem[Har78]{Harrington:1978ht}
\bysame, \emph{Analytic determinacy and {$0\sp{\sharp }$}}, J. Symbolic Logic
  \textbf{43} (1978), no.~4, 685--693.

\bibitem[Hjo95]{Hjorth:1995ve}
G.~Hjorth, \emph{A dichotomy for the definable universe}, Journal of Symbolic
  Logic \textbf{60} (1995), no.~4, 1199--1207.

\bibitem[Hjo99]{Hjorth:1999bh}
\bysame, \emph{Orbit cardinals: on the effective cardinalities arising as
  quotient spaces of the form {$X/G$} where {$G$} acts on a polish space
  {$X$}}, Israel Journal of Mathematics \textbf{111} (1999), 221--261.

\bibitem[Hjo00]{Hjorth:2000zr}
\bysame, \emph{Classification and orbit equivalence relations}, Mathematical
  Surveys and Monographs, no.~75, American Mathematical Society, 2000.

\bibitem[HK95]{Hjorth:1995qf}
G.~Hjorth and A.~S. Kechris, \emph{Analytic equivalence relations and
  {U}lm-type classifications.}, Journal of Symbolic Logic \textbf{60} (1995),
  no.~4, 1273--1300.

\bibitem[HK15]{Hyttinen:2012fj}
T.~Hyttinen and V.~Kulikov, \emph{On $\mathbf{\Sigma}^1_1$-complete equivalence
  relations on the generalized {B}aire space}, Mathematical Logic Quarterly
  \textbf{61} (2015), 66--81.

\bibitem[HKL90]{Harrington:1990qa}
L.~Harrington, A.~S. Kechris, and A.~Louveau, \emph{A {G}limm-{E}ffros
  dichotomy for {B}orel equivalence relations}, Journal of the American
  Mathematical Society \textbf{3} (1990), 903--928.

\bibitem[HMS88]{Harrington:1988ij}
L.~Harrington, D.~Marker, and S.~Shelah, \emph{Borel orderings}, Trans. Amer.
  Math. Soc. \textbf{310} (1988), no.~1, 293--302.

\bibitem[HN73]{Hung:1973fu}
H.~Hung and S.~Negrepontis, \emph{Spaces homeomorphic to $(2^\alpha)_\alpha$},
  Bulletin of the American Mathematical Society \textbf{79} (1973), no.~1,
  143--146.

\bibitem[Hod93]{Hodges:1993mi}
W.~Hodges, \emph{Model theory}, Encyclopedia of Mathematics and Its
  Applications, no.~42, Cambridge University Press, 1993.

\bibitem[HS79]{Harrington:1979hc}
L.~Harrington and R.~Sami, \emph{Equivalence relations, projective and beyond},
  Logic {C}olloquium '78 ({M}ons, 1978), Stud. Logic Foundations Math.,
  vol.~97, North-Holland, Amsterdam, 1979, pp.~247--264. \MR{567673
  (82d:03080)}

\bibitem[HS82]{Harrington:1982if}
L.~Harrington and S.~Shelah, \emph{Counting equivalence classes for
  co-$\kappa$-{S}ouslin equivalence relations.}, van Dalen, D. (eds.) et al.,
  Proc.\ Conf.\ in Prague (LC'80), North Holland, Amsterdam. 147--152., 1982.

\bibitem[HS01]{Halko:2001kl}
A.~Halko and S.~Shelah, \emph{On strong measure zero subsets of {$\sp
  \kappa2$}}, Fund. Math. \textbf{170} (2001), no.~3, 219--229.

\bibitem[Ili12]{Iliadis:2012il}
S.~D. Iliadis, \emph{A note on hierarchies of {B}orel type sets}, Topology
  Appl. \textbf{159} (2012), no.~7, 1702--1704.

\bibitem[Jac89]{Jackson:1989fk}
S.~Jackson, \emph{A{D} and the very fine structure of {$L({\bf R})$}}, Bull.
  Amer. Math. Soc. (N.S.) \textbf{21} (1989), no.~1, 77--81.

\bibitem[Jac08]{Jackson:2008pi}
\bysame, \emph{Suslin cardinals, partition properties, homogeneity.
  {I}ntroduction to {P}art {II}.}, Kechris, Alexander S. (ed.) et al., Games,
  scales, and Suslin cardinals. The Cabal Seminar, Vol. I. Reprints of papers
  and new material based on the Los Angeles Caltech-UCLA Logic Cabal Seminar
  1976--1985. Cambridge: Cambridge University Press; Urbana, IL: Association
  for Symbolic Logic (ASL). Lecture Notes in Logic 31, 273--313., 2008.

\bibitem[Jac10]{Jackson:2010ff}
\bysame, \emph{Structural consequences of {AD}.}, Foreman, Matthew (ed.) et
  al., Handbook of set theory. In 3 volumes. Dordrecht: Springer. 1753--1876,
  2010.

\bibitem[Jec03]{Jech:2003pd}
T.~Jech, \emph{Set theory}, Springer-Verlag, Berlin, 2003. \MR{1940513
  (2004g:03071)}

\bibitem[Kan97]{Kanovei:1997lh}
V.~Kanovei, \emph{An {U}lm-type classification theorem for equivalence
  relations in {S}olovay model}, Journal of Symbolic Logic \textbf{62} (1997),
  no.~4, 1333--1351.

\bibitem[Kan03]{Kanamori:2003fu}
A.~Kanamori, \emph{The higher infinite}, second ed., Springer Monographs in
  Mathematics, Springer-Verlag, Berlin, 2003.

\bibitem[Kan08]{Kanovei:2008qo}
V.~Kanovei, \emph{Borel equivalence relations}, University Lecture Series,
  vol.~44, American Mathematical Society, Providence, RI, 2008.

\bibitem[Kat88]{Katetov:1986tb}
M.~Kat{\v{e}}tov, \emph{On universal metric spaces}, General topology and its
  relations to modern analysis and algebra, {VI} ({P}rague, 1986), Res. Exp.
  Math., vol.~16, Heldermann, Berlin, 1988, pp.~323--330. \MR{952617
  (89k:54066)}

\bibitem[Kec78]{Kechris:1978vn}
A.~S. Kechris, \emph{{${\rm AD}$} and projective ordinals}, Wadge degrees and
  projective ordinals. {T}he {C}abal {S}eminar. {V}olume {II} (A.~S. Kechris,
  B.~L{\"o}we, and J.~R. Steel, eds.), Lecture Notes in Logic, vol.~37,
  Association for Symbolic Logic, Berlin, 1978, pp.~91--132.

\bibitem[Kec84]{Kechris:1984il}
\bysame, \emph{The axiom of determinacy implies dependent choices in {$L({\bf
  R})$}}, J. Symbolic Logic \textbf{49} (1984), no.~1, 161--173.

\bibitem[Kec95]{Kechris:1995zt}
\bysame, \emph{Classical descriptive set theory}, Graduate Text in Mathematics,
  no. 156, Springer-Verlag, Heidelberg, New York, 1995.

\bibitem[Kec99]{Kechris:1999fk}
\bysame, \emph{New directions in descriptive set theory}, Bull. Symbolic Logic
  \textbf{5} (1999), no.~2, 161--174.

\bibitem[Kei71]{Keisler:1971vn}
H.~J. Keisler, \emph{Model theory for infinitary logic. {L}ogic with countable
  conjunctions and finite quantifiers}, North-Holland Publishing Co.,
  Amsterdam, 1971, Studies in Logic and the Foundations of Mathematics, Vol.
  62.

\bibitem[Ker00]{Keremedis:2000fk}
K.~Keremedis, \emph{The compactness of {$2^{\bf R}$} and the axiom of choice},
  Math. Log. Q. \textbf{46} (2000), no.~4, 569--571.

\bibitem[Ket11]{Ketchersid:2011gb}
R.~Ketchersid, \emph{More structural consequences of {${\rm AD}$}}, Set theory
  and its applications, Contemp. Math., vol. 533, Amer. Math. Soc., Providence,
  RI, 2011, pp.~71--105.

\bibitem[KM04]{Kechris:2004jl}
A.~S. Kechris and B.~D. Miller, \emph{Topics in orbit equivalence}, Lecture
  Notes in Mathematics, no. 1852, Springer-Verlag, Berlin, Heidelberg, 2004.

\bibitem[Kni77]{Knight:1977mb}
J.~F. Knight, \emph{A complete ${L}_{\omega_1 \omega}$-sentence characterizing
  $\aleph_1$}, Journal of Symbolic Logic \textbf{42} (1977), no.~1, 59--62.

\bibitem[Kra01]{Kraszewski:2001kr}
J.~Kraszewski, \emph{Properties of ideals on the generalized {C}antor spaces},
  J. Symbolic Logic \textbf{66} (2001), no.~3, 1303--1320.

\bibitem[Kun80]{Kunen:1980fk}
K.~Kunen, \emph{Set theory}, Studies in Logic and the Foundations of
  Mathematics, vol. 102, North-Holland Publishing Co., Amsterdam, 1980.

\bibitem[Kur66]{Kuratowski:1966ku}
K.~Kuratowski, \emph{Topology. {V}ol. {I}}, New edition, revised and augmented.
  Translated from the French by J. Jaworowski, Academic Press, New York-London;
  Pa\'nstwowe Wydawnictwo Naukowe, Warsaw, 1966.

\bibitem[KW10]{Koellner:2010fk}
P.~Koellner and W.~H. Woodin, \emph{Large cardinals from determinacy}, Handbook
  of set theory. {V}ols. 1, 2, 3, Springer, Dordrecht, 2010, pp.~1951--2119.
  \MR{2768702}

\bibitem[LR05]{Louveau:2005cq}
A.~Louveau and C.~Rosendal, \emph{Complete analytic equivalence relations},
  Trans. Amer. Math. Soc. \textbf{357} (2005), no.~12, 4839--4866 (electronic).

\bibitem[LS15]{Luecke:2014ar}
P.~L{\"u}cke and P.~Schlicht, \emph{Continuous images of closed sets in
  generalized {B}aire spaces}, Israel Journal of Mathematics \textbf{209}
  (2015), no.~1, 421--461.

\bibitem[Mar12]{Martin:2012ma}
D.~A. Martin, \emph{Projective sets and cardinal numbers: some questions
  related to the continuum problem}, Wadge degrees and projective ordinals.
  {T}he {C}abal {S}eminar. {V}olume {II}, Lect. Notes Log., vol.~37, Assoc.
  Symbol. Logic, La Jolla, CA, 2012, pp.~484--508.

\bibitem[Mel07]{Melleray:2007di}
J.~Melleray, \emph{Computing the complexity of the relation of isometry between
  separable {B}anach spaces}, MLQ Math. Log. Q. \textbf{53} (2007), no.~2,
  128--131.

\bibitem[Mil11]{Miller:2011fk}
A.~W. Miller, \emph{A {D}edekind finite {B}orel set}, Arch. Math. Logic
  \textbf{50} (2011), no.~1-2, 1--17.

\bibitem[Mos09]{Moschovakis:2009fk}
Y.~N. Moschovakis, \emph{Descriptive set theory}, second ed., Mathematical
  Surveys and Monographs, vol. 155, American Mathematical Society, Providence,
  RI, 2009.

\bibitem[MR12]{Motto-Ros:2012ss}
L.~Motto~Ros, \emph{On the complexity of the relations of isomorphism and
  bi-embeddability}, Proceedings of the American Mathematical Society
  \textbf{140} (2012), no.~1, 309--323.

\bibitem[MR13]{Motto-Ros:2011qc}
\bysame, \emph{The descriptive set-theoretical complexity of the embeddability
  relation on models of large size}, Annals of Pure and Applied Logic
  \textbf{164} (2013), no.~12, 1454--1492.

\bibitem[MR17]{mottoros:2016bu}
\bysame, \emph{Can we classify complete metric spaces up to isometry?}, Boll.
  Unione Mat. Ital. \textbf{10} (2017), no.~3, 369--410. \MR{3691805}

\bibitem[MS70]{Martin:1970ms}
D.~A. Martin and R.~M. Solovay, \emph{Internal cohen extensions}, Annals of
  Mathematical Logic \textbf{2} (1970), 143--178.

\bibitem[MV93]{Mekler:1993kh}
A.~Mekler and J.~V\"a\"an\"anen, \emph{Trees and \( \boldsymbol{\Pi}^1_1
  \)-subsets of \( {}^{\omega_1} \omega_1 \)}, Journal of Symbolic Logic
  \textbf{58} (1993), no.~3, 1052--1070.

\bibitem[MZ55]{Montgomery:1955fk}
D.~Montgomery and L.~Zippin, \emph{Topological transformation groups},
  Interscience Publishers, New York-London, 1955.

\bibitem[Par63]{Parovicenko:1963fc}
I.~I. Parovi{\v c}enko, \emph{A universal bicompact of weight $\aleph_1$},
  Dokl. Akad. Nauk SSSR \textbf{150} (1963), 36--39.

\bibitem[She84]{Shelah:1984kc}
S.~Shelah, \emph{On co-$\kappa$-{S}ouslin relations}, Israel Journal of
  Mathematics \textbf{47} (1984), 139--153.

\bibitem[She01]{Shelah:2001bd}
\bysame, \emph{Strong dichotomy of cardinality}, Results Math. \textbf{39}
  (2001), no.~1-2, 131--154.

\bibitem[She04]{Shelah:2004ud}
\bysame, \emph{On nice equivalence relations on {${}\sp \lambda 2$}}, Arch.
  Math. Logic \textbf{43} (2004), no.~1, 31--64.

\bibitem[Sil80]{Silver:1980fh}
J.~H. Silver, \emph{Counting the number of equivalence classes of {B}orel and
  coanalytic equivalence relations}, Annals of Mathematical Logic \textbf{18}
  (1980), 1--28.

\bibitem[Ste84]{Stern:1984mw}
J.~Stern, \emph{On {L}usin's restricted continuum problem}, Annals of
  Mathematics \textbf{120} (1984), 7--37.

\bibitem[Ste05]{Steel:2005ab}
J.~R. Steel, \emph{{\sf P{FA}} implies {${\sf AD}^{\rm L( \mathbb R ) }$}},
  Journal of Symbolic Logic \textbf{70} (2005), no.~4, 1255--1296.

\bibitem[Ste10]{Steel:2010fk}
J.~Steel, \emph{An outline of inner model theory}, Foreman, Matthew (ed.) et
  al., Handbook of set theory. In 3 volumes. Dordrecht: Springer. 1595--1684,
  2010.

\bibitem[Ste12]{Steel:2012st}
J.~R. Steel, \emph{Closure properties of pointclasses}, Wadge degrees and
  projective ordinals. {T}he {C}abal {S}eminar. {V}olume {II}, Lect. Notes
  Log., vol.~37, Assoc. Symbol. Logic, La Jolla, CA, 2012, pp.~102--117.

\bibitem[Sto48]{Stone:1948st}
A.~H. Stone, \emph{Paracompactness and product spaces}, Bull. Amer. Math. Soc.
  \textbf{54} (1948), 977--982.

\bibitem[SV00]{Shelah:2000hs}
S.~Shelah and J.~V{\"a}{\"a}n{\"a}nen, \emph{Stationary sets and infinitary
  logic}, J. Symbolic Logic \textbf{65} (2000), no.~3, 1311--1320.

\bibitem[SV02]{Shelah:2002lo}
S.~Shelah and P.~V{\"a}is{\"a}nen, \emph{On equivalence relations second order
  definable over {$H(\kappa)$}}, Fund. Math. \textbf{174} (2002), no.~1, 1--21.

\bibitem[SW90]{Shelah:1990sw}
S.~Shelah and H.~Woodin, \emph{Large cardinals imply that every reasonably
  definable set of reals is {L}ebesgue measurable}, Israel Journal of
  Mathematics \textbf{70} (1990), 381--393.

\bibitem[SWar]{Steel:yq}
J.~R. Steel and W.~Hugh Woodin, \emph{{HOD} is a core model}, Kechris,
  Alexander S. (ed.) et al., Ordinal definability and recursion theory. The
  Cabal Seminar, Vol. III. Reprints of papers and new material based on the Los
  Angeles Caltech-UCLA Logic Cabal Seminar 1976--1985, (to appear).

\bibitem[TT93]{Tamano:1993tt}
Ken-ichi Tamano and Hui Teng, \emph{Normality and covering properties of open
  sets of uncountable products}, Topology Proc. \textbf{18} (1993), 313--322.

\bibitem[{\"U}nl82]{Unlu:1982un}
Yusuf {\"U}nl{\"u}, \emph{Spaces for which the generalized {C}antor space
  {$2\sp{J}$} is a remainder}, Proc. Amer. Math. Soc. \textbf{86} (1982),
  no.~4, 673--678.

\bibitem[V{\"a}{\"a}11]{Vaananen:2011tg}
J.~V{\"a}{\"a}n{\"a}nen, \emph{Models and games}, Cambridge Studies in Advanced
  Mathematics, vol. 132, Cambridge University Press, Cambridge, 2011.

\bibitem[Vau75]{Vaught:1974kl}
R.~Vaught, \emph{Invariant sets in topology and logic}, Fundamenta
  Mathematic\ae{} \textbf{82} (1974/75), 269--294.

\bibitem[Ver98]{Vershik:1998fk}
A.~M. Vershik, \emph{The universal {U}ryson space, {G}romov's metric triples,
  and random metrics on the series of natural numbers}, Uspekhi Mat. Nauk
  \textbf{53} (1998), no.~5(323), 57--64.

\bibitem[Wad83]{Wadge:1983sp}
William~W. Wadge, \emph{Reducibility and determinateness on the baire space},
  Ph.D. thesis, 1983.

\bibitem[Woo10]{Woodin:2010tn}
W.~H. Woodin, \emph{The axiom of determinacy, forcing axioms, and the
  nonstationary ideal}, second ed., de Gruyter Series in Logic and its
  Applications, vol.~1, Walter de Gruyter \& Co., Berlin, 2010.

\end{thebibliography}

\printindex

\end{document}